\documentclass[a4paper]{book}
\pdfoutput=1 

\usepackage{blindtext}
\usepackage[T1]{fontenc}
\usepackage{graphicx}
\usepackage{chngcntr}

\usepackage{tabularx} 
\usepackage{url}
\usepackage{color}            
\usepackage{amssymb, amsmath}
\usepackage{mathtools}
\usepackage{dsfont}

\usepackage{setspace}
\newcommand{\fullname}{Luiz Frederic Wagner}
\newcommand{\titel}{Pseudo-Holomorphic Hamiltonian Systems and Complex Coadjoint Orbits}

\usepackage{hyperref}
\hypersetup{
pdftitle=\titel,
pdfauthor=\fullname,
pdfproducer={pdfeTex 3.14159-1.30.6-2.2},
colorlinks=false,
pdfborder=0 0 0	
}

\usepackage{amsthm}
\theoremstyle{definition}
\newtheorem{definition}{Definition}[section]
\newtheorem{theorem}[definition]{Theorem}
\newtheorem{corollary}[definition]{Corollary} 
\newtheorem{lemma}[definition]{Lemma}
\newtheorem{example}[definition]{Example}
\newtheorem{remark}[definition]{Remark}
\newtheorem{proposition}[definition]{Proposition}

\newtheorem{conjecture}[definition]{Conjecture}
\newtheorem{auxiliarylemma}[definition]{Auxiliary lemma}

\newtheorem*{definition*}{Definition}
\newtheorem*{theorem*}{Theorem}
\newtheorem*{corollary*}{Corollary} 
\newtheorem*{lemma*}{Lemma}
\newtheorem*{example*}{Example}
\newtheorem*{remark*}{Remark}
\newtheorem*{proposition*}{Proposition}
\newtheorem*{notation*}{Notation}
\newtheorem*{conjecture*}{Conjecture}
\newtheorem*{auxiliarylemma*}{Auxiliary lemma}

\newcommand{\mH}{\mathcal{H}}

\newcommand{\actdiski}{\mathcal{A}^{D^{z_0}_R}_{\mH; 1}}
\newcommand{\actdiskii}{\mathcal{A}^{D^{z_0}_R}_{\mH; 2}}
\newcommand{\acti}{\mathcal{A}^{D}_{\mH; 1}}
\newcommand{\actii}{\mathcal{A}^{D}_{\mH; 2}}
\newcommand{\nab}[1]{\prescript{#1}{}{\nabla}}
\newcommand{\tor}[1]{\prescript{#1}{}{T}}
\newcommand{\hvec}[1]{\Gamma (T^{(1,0)} #1)}

\newcommand{\End}{\operatorname{End}}
\newcommand{\Spec}{\operatorname{Spec}}
\newcommand{\Rec}{\operatorname{Rec}}
\newcommand{\id}{\operatorname{id}}
\newcommand{\im}{\operatorname{im}}
\newcommand{\can}{\operatorname{can}}
\newcommand{\pa}[1]{\partial_{ #1 }}
\newcommand{\pra}[2]{\frac{\partial #1}{\partial #2}}

\newcommand{\R}{\mathbb{R}}
\newcommand{\C}{\mathbb{C}}
\newcommand{\Qua}{\mathbb{H}}
\newcommand{\Z}{\mathbb{Z}}

\newcommand{\Sym}{\operatorname{Sym}}
\newcommand{\Span}{\operatorname{Span}}
\newcommand{\U}{\operatorname{U}}
\newcommand{\SU}{\operatorname{SU}}
\newcommand{\Or}{\operatorname{O}}
\newcommand{\GL}{\operatorname{GL}}
\newcommand{\Sp}{\operatorname{Sp}}
\newcommand{\Fix}{\operatorname{Fix}}
\newcommand{\Ad}{\operatorname{Ad}}
\newcommand{\Aut}{\operatorname{Aut}}
\newcommand{\ad}{\operatorname{ad}}
\newcommand{\tr}{\operatorname{tr}}
\newcommand{\pr}{\operatorname{pr}}
\newcommand{\sk}[2]{\left\langle #1 , #2 \right\rangle}
\newcommand{\skcdot}{\left\langle \cdot , \cdot \right\rangle}
\newcommand{\kks}{\omega_{\text{KKS}}}
\newcommand{\KKS}{\Omega_{\text{KKS}}}

\newcommand{\ver}{\operatorname{vert}}
\newcommand{\geo}{\operatorname{geo}}
\newcommand{\vol}{\operatorname{vol}}
\newcommand{\ev}{\operatorname{ev}}
\newcommand{\sgn}{\operatorname{sgn}}
\newcommand{\rlc}{r.l.c. }
\newcommand{\clc}{c.l.c. }

\usepackage{enumitem}
\setlist[enumerate,1]{label={(\roman*)}}

\usepackage{tikz}
\usepackage{tikz-cd}
\usetikzlibrary{arrows.meta}
\usetikzlibrary{angles}
\usepackage{caption}
\usepackage{subcaption}
\usepackage{todonotes}
\usepackage{wrapfig}
\usepackage{graphicx, array, blindtext}

\makeatletter
\let\oldtheequation\theequation
\renewcommand\tagform@[1]{\maketag@@@{\ignorespaces#1\unskip\@@italiccorr}}
\renewcommand\theequation{(\oldtheequation)}
\makeatother

\title{\titel}
\author{\fullname}
\date{\today}

\begin{document}

%

\begin{titlepage}
\thispagestyle{empty}
\begin{center}
	\Large{\textbf{Pseudo-Holomorphic Hamiltonian Systems and\\ Kähler Duality of Complex Coadjoint Orbits}}\\
	\doublespacing 
	\vspace{1cm}
	\large{\textbf{Dissertation}}\\
	\vspace{1cm}	
	\large
	zur Erlangung des akademischen Grades\\
	Dr.rer.nat.\\	
	\vspace{0.4cm}
	\begin{figure}[h]
	\centering
	\includegraphics[width=0.4\textwidth]{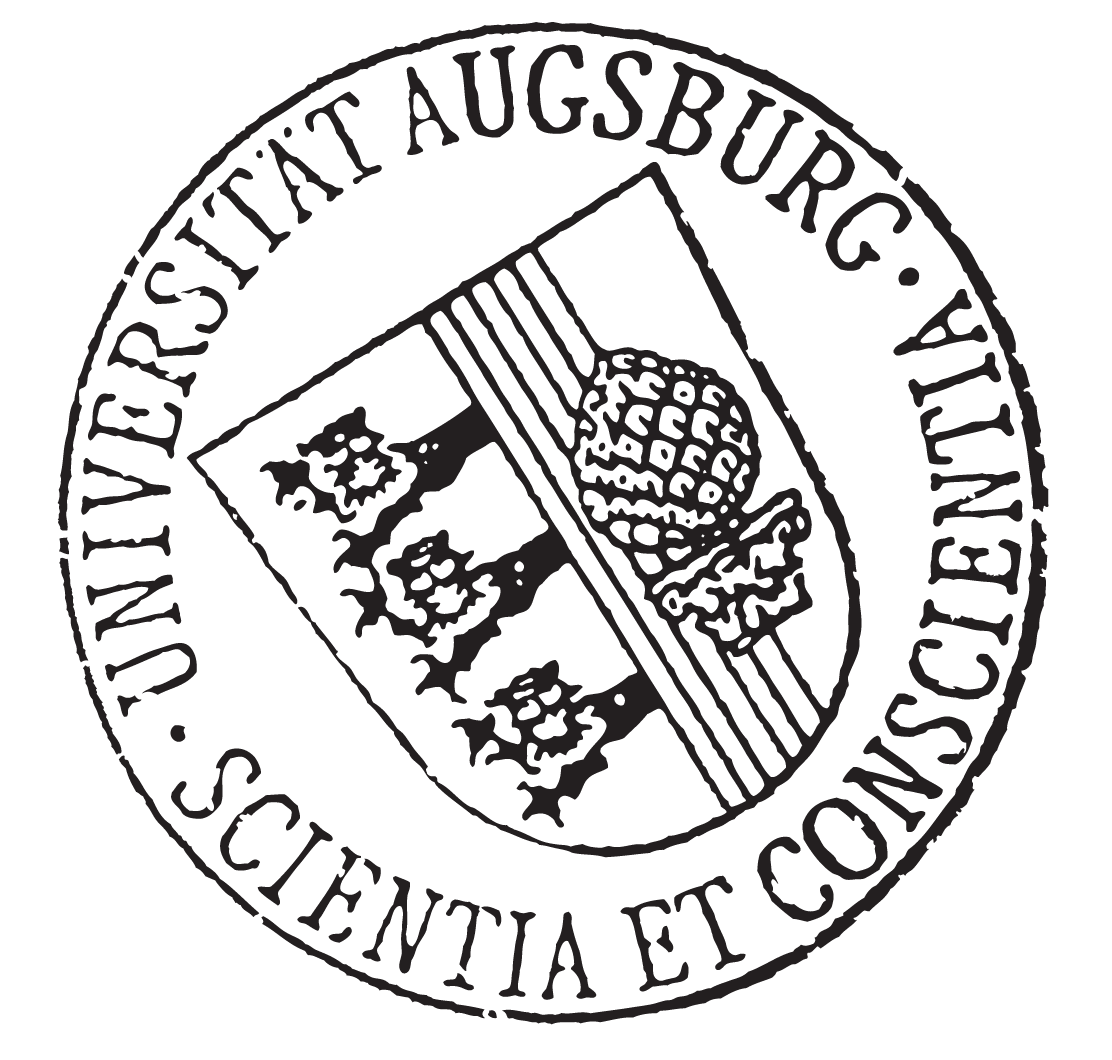}
	\end{figure}
	eingereicht an der\\
	Mathematisch-Naturwissenschaftlich-Technischen Fakultät\\
	der Universität Augsburg\\
	\vspace{1cm}
	von\\
	\textbf{Luiz Frederic Wagner}\\
	\vspace{1cm}
	Augsburg, August 2024
	\vfill
\end{center}

\newpage
\thispagestyle{empty}
\vfill

\begin{center}
 \begin{tabular}{rl}
  Erstgutachter: & Prof. Dr. Kai Cieliebak, Universität Augsburg\\
  Zweitgutachter: & Prof. Dr. Marc Nieper-Wißkirchen, Universität Augsburg\\
  Drittgutachter: & Dr. Stefan Nemirovski, Ruhr-Universität Bochum\\
  Tag der mündlichen Prüfung: & 04.11.2024
 \end{tabular}
\end{center}

\end{titlepage}

\pagenumbering{roman}

\tableofcontents

\newpage
\thispagestyle{empty}\phantom{This page is intentionally left blank!}
\newpage

\pagenumbering{arabic}

\chapter*{Acknowledgments}
\addcontentsline{toc}{chapter}{Acknowledgments}
\markboth{}{Acknowledgments}
First and foremost, I would like to thank Prof. Dr. Kai Cieliebak for his supervision of my doctorial studies and my thesis. I cannot recall a conversation with him where I did not leave with a better understanding of mathematics. Furthermore, I would like to express my gratitude toward Prof. Dr. Marc Nieper-Wißkirchen and Dr. Stefan Nemirovski who agreed to write the examiner's reports. On top of that, I am very grateful to the symplectic group in Augsburg for insightful discussions and their support over the years, in particular Prof. Dr. Urs Frauenfelder who has been a great teacher, PD Dr. Lei Zhao who motivated me to write down \autoref{chap:PHHS} of this thesis, and my ``proofreaders'' Dr. Airi Takeuchi, Hanna Häußler, Marián Poppr, Martin Konrad, Milan Zerbin, and Shuaipeng Liu. Of course, I also express my thanks toward the other members of the group as well as former ones and guests, including (but not limited to) Dr. Kevin Ruck, Dr. Evgeny Volkov, Dr. Yannis Bähni, Dr. Cengiz Aydin, Dr. Amanda Hirschi, Dr. Yuchen Wang, Zhen Gao, Jennifer Gruber, Julius Natrup, Sungho Kim, and Shah Faisal. Lastly, I would like to thank my friends Christian Schneider, Alexander Segner, and Florian Stuhlmann, the members of the ``differential geometry club'' as well as my family for their support.

\newpage\thispagestyle{empty}\

\chapter{Introduction}
\label{chap:intro}
Over the last century, complexification has become a popular and versatile tool in various fields of mathematics. In group theory, for instance, complexification is quite helpful in the classification of simple Lie algebras.
Simple Lie algebras are completely described by their root systems which themselves are fully encoded in so-called Dynkin diagrams. Unfortunately, Dynkin diagrams only classify simple Lie algebras over algebraically closed fields meaning they are not suited for real Lie algebras. To rectify this, one complexifies the real Lie algebra $\mathfrak{g}$ to the complex Lie algebra $\mathfrak{g}_\C\coloneqq \mathfrak{g}\oplus i\mathfrak{g}$. The real form $\mathfrak{g}$ of $\mathfrak{g}_\C$ can then be extracted from the Dynkin diagram for $\mathfrak{g}_\C$ by the means of Satake diagrams, a special variant of Dynkin diagrams. Naturally, complexification is similarly useful for other mathematical problems, e.g. the classification of linear representations of Lie groups (cf. the notion of weights).\\
The overarching idea of this thesis is essentially to complexify mathematical objects, usually from symplectic geometry, in order to achieve one of the following goals:
\begin{enumerate}
 \item Determine which properties transfer from the real to the complex system and work out possible differences.
 \item Link two a priori different notions via a complexified structure.
 \item Gain a better understanding of already established results and reveal hidden connections using complexification.
\end{enumerate}
The main part of this thesis is comprised of two chapters, each devoted to a structure whose complexification we wish to study in detail. \autoref{chap:PHHS}\footnote{\autoref{chap:PHHS} can be found in similar form as a preprint on arXiv (cf. \cite{Wagner2023}).} focuses on Hamiltonian systems, while \autoref{chap:duality} deals with coadjoint orbits, specifically their Kähler structure.

\subsubsection*{Holomorphic Hamiltonian Systems}

A real Hamiltonian system (RHS) consists of a manifold $M$ together with a symplectic form $\omega$ and a function $H:M\to\R$ called Hamiltonian. The objects $\omega$ and $H$ allow us to define the Hamiltonian vector field $X_H$ via the formula $\iota_{X_H}\omega = -dH$. Since their inception, RHSs have attracted a vast amount of attention from the scientific community, mostly due to their physical relevance: The integral curves of $X_H$ describe the trajectories of massive point-like particles in the phase space $(M,\omega)$ subject to the dynamics dictated by the energy function $H$.\\
The complex analogue of a RHS, a \textbf{holomorphic Hamiltonian system} (HHS), is described by similar data: a complex manifold $X$, a holomorphic symplectic form $\Omega$, and a holomorphic function $\mH:X\to\C$. HHSs have been studied for the last 20 years. The research in this field is mainly focused on the interplay of HHSs with their real forms. In \cite{gerd2002}, for instance, the authors used the integrability of the Toda chain to show that the complex Toda chain as its complexification is also integrable. From this, they were able to infer that all RHSs which emerge as real forms of the complex Toda chain are integrable as well.\\
In this thesis, we take a different approach: We examine HHSs independently of their real forms and work out which properties transfer from the real to the complex case. Consider, for instance, the trajectories of Hamiltonian systems. It is common knowledge that the maximal trajectories $\gamma:I\subset\R\to M$ as integral curves of $X_H$ exist and are unique given an initial value. Locally, the same statement still holds true for HHSs: As in the real setup, $\Omega$ and $\mH$ single out a unique holomorphic vector field $X_\mH$, the holomorphic Hamiltonian vector field, defined by $\iota_{X_\mH}\Omega = -d\mH$. The holomorphicity of the vector field $X_\mH$ implies that real and imaginary parts of $X_\mH$ commute. Thus, the flows of $\text{Re} X_\mH$ and  $\text{Im} X_\mH$ commute as well. Combining the flows gives us a holomorphic map $\gamma:U\subset\C\to X$ satisfying the holomorphic integral curve equation $\gamma^\prime (z) = X_\mH (\gamma (z))$. For small enough domains $U$, $\gamma$ is well-defined and unique given an initial value. However, the maximal trajectories $\gamma:U\subset\C\to X$ are not unique in that sense. This behavior is in sharp contrast to the real case. The non-uniqueness of maximal holomorphic trajectories is caused by monodromy effects or, simply put, by the fact that the flows of $\text{Re} X_\mH$ and  $\text{Im} X_\mH$ do not commute globally. In \autoref{sec:HHS}, we use a Kepler-like central problem (cf. Example \autoref{ex:holo_cen_prob}) with Hamiltonian $\mH(Q,P)\coloneqq P^2/2 - 1/8Q^2$ to demonstrate our findings. The holomorphic trajectories of this system include square roots, e.g. the solution of the Hamilton equations
\begin{gather*}
 Q^\prime (z) = P(z),\quad P^\prime (z) = -\frac{1}{4Q^3 (z)}
\end{gather*}
with initial value $Q(0) = 1$ and $P(0) = 1/2$ is given by $Q(z) = \sqrt{z+1}$. Depending on how the domain of the square root is chosen, i.e., depending on where the branch cut of $\sqrt{\cdot}$ lies, $Q(z) = \sqrt{z+1}$ defines different maximal trajectories with initial value $(Q(0), P(0)) = (1,1/2)$.\\
It has been observed before that HHSs are affected by monodromy (cf. \cite{shanzhong2020}). Still, to the author's knowledge, it has not been shown before that the holomorphic trajectories themselves exemplify monodromy. Similarly, it is a new result that we can restore the uniqueness of maximal trajectories by promoting them to leaves of a certain foliation. Here, we again draw inspiration from RHSs: If $E\in\R$ is a regular energy value of the RHS $(M,\omega,H)$, then $X_H$ gives rise to a one-dimensional distribution on the energy hypersurface $H^{-1}(E)$. Every one-dimensional distribution is involutive, hence, $H^{-1}(E)$ admits a foliation whose leaves are tangent to $X_H$ meaning they are the maximal trajectories. In the complex case, $\text{Re} X_\mH$ and  $\text{Im} X_\mH$ span a two-dimensional distribution on $\mH^{-1}(E)$, where $E\in\C$ is now a regular value of $\mH$. The distribution is involutive, as $\text{Re} X_\mH$ and  $\text{Im} X_\mH$ commute. Again, this yields a foliation of the hypersurface $\mH^{-1}(E)$, however, the leaves are two-dimensional this time. By the holomorphic Frobenius theorem, the foliation is even holomorphic, so the leaves are, in fact, immersed Riemann surfaces. By construction, the leaf through a given point $p\in\mH^{-1}(E)$, which we can interpret as an initial value, is unique. Even though each maximal trajectory is contained in one leaf, the leaves do not agree with the maximal trajectories in general. Take, for example, the Kepler-like problem from before. The leaves in this case are the connected components of the hypersurfaces $\mH^{-1}(E)$ which $Q(z) = \sqrt{z+1}$ is clearly not.

\subsubsection*{Holomorphic Symplectic Lefschetz Fibrations}

Restoring the uniqueness of maximal trajectories is not the only advantage the holomorphic foliation has. It can also be seen as a link between Lefschetz and almost toric fibrations. Lefschetz fibrations were initially introduced by their namesake to study the topology of complex surfaces, but piqued the interest of symplectic geometers like Donaldson and Gompf in the 80s because of their relation to symplectic four-folds. Broadly speaking, they can be understood as generalized fiber bundles over a surface. Usually, the fibers of a fiber bundle are all isomorphic. Lefschetz fibrations allow for singular fibers with special local structure. If $\pi:X\to C$ is a Lefschetz fibration ($\dim_\R X = 2m$, $\dim_\R C = 2$), then $\pi$ takes the following form near singular fibers, i.e., near critical points:
\begin{gather*}
 \psi_C\circ \pi\circ\psi_X^{-1} (z_1,\ldots,z_m) = \sum^m_{j=1} z^2_j,
\end{gather*}
where $\psi_X$ and $\psi_C$ are smooth, $\C$-valued charts of $X$ and $C$, respectively.\\
Similarly, almost toric fibrations also expand the definition of a fiber bundle by singular fibers. First established by M. Symington in 2002 (cf. \cite{Symington2002}), they generalize the notion of toric fibrations, i.e., moment maps of effective Hamiltonian torus actions. Accordingly, almost toric fibrations only make sense for symplectic manifolds $(X,\omega)$. By definition, the projection $\pi:(X,\omega)\to C$ of an almost toric fibration ($\dim_\R X = 2m$, $\dim_\R C = m$) assumes the following local structure in suitable charts of $X$ and $C$ ($0\leq k\leq m$):
\begin{enumerate}
 \item $\omega = \sum^m_{j=1} dx_j\wedge dy_j$,
 \item $\pi_j = x_j$ for $1\leq j\leq k$,
 \item $\pi_j = x^2_j + y^2_j$ (toric) \textbf{or}\\
 $(\pi_j, \pi_{j+1}) = (x_jy_j + x_{j+1}y_{j+1}, x_jy_{j+1} - x_{j+1}y_j)$ (nodal) for $k<j\leq m$.
\end{enumerate}
For $k = m$, all points in the chart domain are regular and, for $k<m$, the chart domain intersects singular fibers. The key difference between toric and almost toric fibrations are the nodal points. Indeed, if all singularities are of toric type, then $\pi$ describes a toric fibration (hence the name). Toric and nodal points not only differ in their local structure, but also in their position: While toric singularities lie on the boundary of $\im\pi$, nodal points live in the interior of $\im\pi$.\\
To establish a connection between HHSs, Lefschetz fibrations, and almost toric fibrations, we introduce an object which encapsulates the essence of these three notions: a \textbf{holomorphic symplectic Lefschetz fibration}. Roughly speaking, a holomorphic Lefschetz fibration is a Lefschetz fibration $\pi:X\to C$, where $X$ and $C$ are complex manifolds and $\pi$ is a holomorphic map. We say a holomorphic Lefschetz fibration $\pi:X\to C$ is symplectic if $X$ carries a holomorphic symplectic form $\Omega$ which is compatible with $\pi$ in the sense that there are holomorphic Morse-Darboux charts $\psi_X = (z_1,\ldots,z_{2n})$ and $\psi_C$ near critical points in which $\Omega$ and $\pi$ take the following form:
\begin{gather*}
 \Omega = \sum^n_{j=1} dz_{j+n}\wedge dz_j,\quad \psi_C\circ\pi = \sum^{2n}_{j=1} z^2_j.
\end{gather*}
The existence of Morse-Darboux charts ensures that holomorphic symplectic Lefschetz fibrations are almost toric. Precisely speaking, every proper holomorphic symplectic Lefschetz fibration $\pi:(X,\Omega)\to C$ with $\dim_\R X = 4$\linebreak gives rise to two almost toric fibrations, namely $\pi:(X,\Omega_R)\to C$ and\linebreak $\pi:(X,\Omega_I)\to C$, where $\Omega\eqqcolon \Omega_R + i\Omega_I$ (cf. Proposition \autoref{prop:holo_lef_toric}). The fibration $\pi:(X,\Omega)\to C$ can also be seen as a HHS. Indeed, after choosing a holomorphic chart for $C$, $(X,\Omega,\pi)$ becomes a HHS. As explained before, the regular energy hypersurfaces of $(X,\Omega,\pi)$ admit a holomorphic foliation. In dimension $4$, the leaves of this foliation are exactly the regular fibers of the Lefschetz/almost toric fibration $\pi$. This connection between HHSs, Lefschetz fibrations, and almost toric fibrations has not been observed before.\\
During the investigation of Lefschetz and almost toric fibrations, we also tackle the question\footnote{To the author's knowledge, this question has not been examined yet by the scientific community.} whether there are obstructions for a holomorphic Lefschetz fibration $\pi:X\to C$ equipped with a holomorphic symplectic form $\Omega$ to possess Morse-Darboux charts. As it turns out, it suffices to answer this question on a real form: The given problem is local in nature, thus, we only need to consider the case $X = \C^{2n}$ and $C = \C$. If $\Omega$ and $\pi$ assume the standard form
\begin{gather*}
 \Omega = \sum^n_{j=1} dz_{j+n}\wedge dz_j,\quad \pi = \sum^{2n}_{j=1} z^2_j,
\end{gather*}
they reduce to the real-analytic tensors
\begin{gather*}
 \omega = \sum^n_{j=1} dx_{j+n}\wedge dx_j,\quad f = \sum^{2n}_{j=1} x^2_j
\end{gather*}
on the real form $\R^{2n}\subset\C^{2n}$. By unique continuation, the pair $(\omega,f)$ completely determines $(\Omega,\pi)$. Therefore, real-analytic Morse-Darboux charts on the real form $\R^{2n}$ automatically give us holomorphic Morse-Darboux charts on $\C^{2n}$ by complexification. We will show in \autoref{sec:HHS} that the problem of finding real-analytic Morse-Darboux charts can be expressed in an elegant way (cf. Theorem \autoref{thm:morse_darboux_ex}):\pagebreak

\begin{theorem*}[Existence of real-analytic Morse-Darboux charts]
 Let $(M^{2n},\omega)$ be a real-analytic symplectic manifold, let $f:M\to\R$ be a real-analytic function, and let $p\in M$. Then, the following statements are equivalent:
 \begin{enumerate}[label = (\alph*)]
  \item There is a real-analytic Morse-Darboux chart near $p$, i.e., a real-analytic chart $\psi = (x_1,\ldots, x_{2n}):U\to V\subset\R^{2n}$ of $M$ near $p$ with $\psi (p) = 0$ such that:
  \begin{gather*}
   \omega = \sum^n_{j=1} dx_{j+n}\wedge dx_j,\quad f = f(p) + \sum^{2n}_{j=1} x^2_j.
  \end{gather*}
  \item There exists a flat Kähler structure near $p$ with symplectic form $\omega$ and mixed\footnote{Confer Definition \autoref{def:mixed}.} Kähler potential $\frac{f-f(p)}{2}$, i.e., there is an open neighborhood\linebreak $U\subset M$ of $p$ and an almost complex structure $J$ on $U$ such that:
  \begin{enumerate}[label = (\roman*)]
   \item $J$ is integrable,
   \item $g\coloneqq \omega(\cdot,J\cdot)$ is a flat Riemannian metric,
   \item $\frac{f-f(p)}{2}$ is the mixed Kähler potential near $p$.
  \end{enumerate}
 \end{enumerate}
\end{theorem*}

With our previous knowledge, proving existence of holomorphic Morse-Darboux charts now reduces to the problem of finding a real form which locally exhibits the Kähler structure specified in the previous theorem (cf. Corollary \autoref{cor:morse_darboux_holo_ex}).\\
Regarding the existence of Morse-Darboux charts, we are mostly interested in the case $\dim_\C X = 2$, since only in this dimension holomorphic Lefschetz fibrations can be almost toric. Judging by the real analogue, there seems to be no obstruction in the two-dimensional case: Given a real-analytic RHS $(M,\omega,f)$ of dimension $2$ and a critical point $p\in M$ of $f$ with Morse index $\neq 1$, we can always find a real-analytic diffeomorphism $\psi:\R\to\R$ and a real-analytic chart $\psi_M = (x,y)$ of $M$ near $p$ such that $\omega = dy\wedge dx$ and $H = x^2 + y^2$, where $H\coloneqq \psi\circ f$ (cf. Lemma \autoref{lem:morse_darboux_lem_I} and  \autoref{lem:morse_darboux_lem_II}).

\subsubsection*{Action Functionals for HHSs}

The trajectories of a Hamiltonian system cannot only be interpreted as integral curves of the Hamiltonian vector field, but also as critical points of an action functional. Recall that, for an exact RHS $(M,\omega = d\lambda,H)$ and an interval $I_0 = [t_1,t_2]$, the curve $\gamma\in C^\infty(I_0,M)$ is a critical point of the action functional $\mathcal{A}^\lambda_H:C^\infty(I_0,M)\to\mathbb{R}$,
\begin{gather*}
 \mathcal{A}^\lambda_H[\gamma]\coloneqq \int\limits_{I_0} \gamma^\ast\lambda - \int\limits^{t_2}_{t_1} H\circ\gamma (t)\, dt,
\end{gather*}
if and only if $\gamma$ satisfies $\dot\gamma = X_H\circ\gamma$. Here, we take the term ``critical point'' with a grain of salt: We only consider those variations $\gamma_\varepsilon$ of $\gamma$ which keep $\gamma$ fixed at the boundary, i.e., $\gamma_\varepsilon (t_1) = \gamma (t_1)$ and $\gamma_\varepsilon (t_2) = \gamma (t_2)$ for all $\varepsilon$ (cf. Remark \autoref{rem:critical_point}). If one wishes to turn the physical trajectories into actual critical points, one can either restrict $\mathcal{A}^\lambda_H$ to periodic curves ($\gamma (t_1) = \gamma (t_2)$) or consider only those curves for which $\lambda$ vanishes at the start and end point ($\lambda_{\gamma(t_1)} = \lambda_{\gamma (t_2)} = 0$).\\
At the end of \autoref{sec:HHS}, we complexify this observation\footnote{Apparently, this has not been done before.}: We assign to each exact HHS an action functional whose critical points are the holomorphic trajectories. To construct such a functional, we first observe that an exact HHS $(X,\Omega = d\Lambda,\mH)$ gives rise to four exact RHSs. The Hamiltonian vector fields of these RHSs are -- up to minus signs and factors of $2$ -- the real and imaginary parts of the holomorphic Hamiltonian vector field $X_\mH$. Now let $\Rec\coloneqq [t_1,t_2] + i[s_1,s_2]\subset\C$ be a rectangle in the complex plane and let $\gamma\in C^\infty (\Rec,X)$ be a smooth map defined on that rectangle. The curve $\gamma$ is a holomorphic trajectory of the HHS $(X,\Omega,\mH)$ if and only if $\gamma_s$ is an integral curve\footnote{We utilize the decompositions $\Omega = \Omega_R + i\Omega_I$,  $\Lambda = \Lambda_R + i\Lambda_I$, $\mH = \mH_R + i\mH_I$, and $X_\mH = 1/2(X^R_\mH - iJ(X^R_\mH))$, where $J$ is the complex structure of $X$.} of $X^R_\mH$ for every $s\in [s_1,s_2]$ and $\gamma_t$ is an integral curve of $J(X^R_\mH)$ for every $t\in [t_1,t_2]$ ($\gamma_s (t)\coloneqq \gamma (t+is)\eqqcolon \gamma_t (s)$). Thus, assuming that $\gamma$ is a holomorphic trajectory, $\gamma_s$ and $\gamma_t$ ($s\in [s_1,s_2]$ and $t\in [t_1,t_2]$ fixed) are critical points of the action functionals $\mathcal{A}^{\Lambda_R}_{\mH_R}$ and $\mathcal{A}^{\Lambda_R}_{-\mH_I}$ associated with the RHSs $(X,\Omega_R,\mH_R)$ and $(X,\Omega_R,-\mH_I)$, respectively. In order for find functionals whose critical points fulfill this property not just for a fixed $s$ or $t$, but for all numbers in the respective intervals, we have to integrate $\mathcal{A}^{\Lambda_R}_{\mH_R}$ and $\mathcal{A}^{\Lambda_R}_{-\mH_I}$ over the remaining variable giving us the averaged action functionals:
\begin{gather*}
 \gamma\mapsto\int\limits^{s_2}_{s_1}\mathcal{A}^{\Lambda_R}_{\mH_R}[\gamma_s]ds\quad\text{and}\quad\gamma\mapsto\int\limits^{t_2}_{t_1}\mathcal{A}^{\Lambda_R}_{-\mH_I}[\gamma_t]dt.
\end{gather*}
Taking a suitable complex combination of the averaged functionals yields the desired action functional $A^{\Rec}_{\mH}:C^\infty (\Rec,X)\to\C$:
\begin{align*}
 \mathcal{A}^{\Rec}_\mH[\gamma]&\coloneqq \int\limits^{s_2}_{s_1} \mathcal{A}^{\Lambda_R}_{\mH_R}[\gamma_s] ds - i\int\limits^{t_2}_{t_1} \mathcal{A}^{\Lambda_R}_{-\mH_I}[\gamma_t] dt\\
 &= \int\limits^{t_2}_{t_1}\int\limits^{s_2}_{s_1}\left[\Lambda_{R,\gamma (t+is)}\left(2\frac{\partial\gamma}{\partial z}(t+is)\right) - \mH\circ\gamma (t+is)\right] ds\, dt\ \text{with}\\
 \frac{\partial \gamma}{\partial z}&\coloneqq \frac{1}{2}\left(\frac{\partial \gamma}{\partial t} - i\frac{\partial \gamma}{\partial s}\right)\quad\forall \gamma\in C^\infty (\Rec,X).
\end{align*}
The critical points of $\mathcal{A}^{\Rec}_\mH$ are the holomorphic trajectories of the HHS $(X,\Omega,\mH)$ with domain $\Rec$. As in the real case, the term ``critical point'' is used rather loosely here: We only look at those variations which keep $\gamma$ fixed at $\partial\Rec$. One can circumvent this by either restricting $\mathcal{A}^{\Rec}_\mH$ to those maps which are periodic in both $s$- and $t$-direction or by assuming that all curves in the domain of $\mathcal{A}^{\Rec}_\mH$ send $\partial\Rec$ to points where $\Lambda_R$ vanishes.\\
$\mathcal{A}^{\Rec}_\mH$ is not the only action functional whose critical points are holomorphic trajectories. In fact, there is a plethora of functionals exhibiting this property, some of which are even real. They differ by the shape of $\gamma$'s domain, by how the action functionals of the RHSs are averaged, and by which complex combination is taken. A large selection of action functionals is explored in \autoref{app:various_action_functionals}.\\
There is one noteworthy aspect about action functionals which describe curves defined on parallelograms $P_\alpha\coloneqq [0,t]+e^{i\alpha}[0,r]$. If such a curve is holomorphic and periodic in $t$- and $r$-direction, then it factors through a complex torus $\C/\Gamma$, where the lattice $\Gamma$ is spanned by $t$ and $re^{i\alpha}$. In general, two complex tori are not biholomorphic. Hence, an action functional for doubly-periodic curves also measures the complex structure of the curve's domain. This feature is quite remarkable, especially since it has no real analogue: Given a RHS, the domains of its periodic orbits are all isomorphic to $S^1$.

\subsubsection*{Pseudo-Holomorphic Hamiltonian Systems}

Even though doubly-periodic trajectories occur in some situations (cf. Example \autoref{ex:complex_torus}), they are exceedingly rare. For instance, any HHS on $X = \C^{2n}$ can only admit trivial doubly-periodic trajectories. This phenomenon is caused by the maximum principle: Any holomorphic map from a compact complex manifold to $\C^{2n}$ must be constant. However, the maximum principle only poses a problem if $\C^{2n}$ is equipped with the standard complex structure. Indeed, Moser showed in his beautifully written paper \cite{moser1995} that $\C^2\cong\R^4$ possesses an almost complex structure $J$ such that the standard complex tours can be pseudo-holomorphically\footnote{We say $f:(X_1,J_1)\to (X_2,J_2)$ is holomorphic if $df\circ J_1 = J_2\circ df$. If, additionally, $J_1$ or $J_2$ is not integrable, we emphasize this point by calling $f$ pseudo-holomorphic.} embedded into $(\R^4,J)$. Observe that the almost complex structure $J$ constructed by Moser is not integrable.\\
Inspired by Moser's construction, we introduce a new type of Hamiltonian system in this thesis, called \textbf{pseudo-holomorphic Hamiltonian system} (PHHS), which generalizes the notion of HHSs to almost complex manifolds $(X,J)$. A priori, it is not clear what this generalization should look like, since the complex structure $J$ only enters the definition of a HHS implicitly. Simply dropping the integrability of $J$ while keeping the other stipulations in place does not work. We will demonstrate in \autoref{sec:PHHS} that it is the closedness of $\Omega_I$ which poses a problem in this case. To explain how to obtain a reasonable notion of PHHSs, we first recall that a HHS is described by a triple $(X,\Omega,\mH)$. We can divide the triple $(X,\Omega,\mH)$ into six objects: the manifold $X$, the complex structure $J$ of $X$, and the real and imaginary parts of $\Omega$ and $\mH$, i.e., $\Omega = \Omega_R + i\Omega_I$ and $\mH = \mH_R + i\mH_I$. These objects are not independent, but rather satisfy the following relations:
\begin{gather*}
 \Omega (J\cdot,J\cdot) = -\Omega,\quad \Omega (J\cdot,\cdot) = i\Omega,\quad d\mH\circ J = i\, d\mH.
\end{gather*}
We now see that $\Omega_I$ and $\mH_I$ are redundant, as they can be recovered from $J$, $\Omega_R$, and $\mH_R$ using the relations above. This observation leads us to the idea to only use the minimal set of data to define a PHHS. In that spirit, a PHHS is a quadruple $(X,J;\Omega_R,\mH_R)$ where $X$ is a smooth manifold with almost complex structure $J$ on it, $\Omega_R$ is a symplectic form on $X$ satisfying $\Omega_R (J\cdot,J\cdot) = -\Omega_R$, and $\mH_R:X\to\R$ is a smooth function such that $d\mH_R\circ J$ is exact (cf. Definition \autoref{def:PHHS_2}).\\
In many regards, a PHHS exhibits the same properties as a HHS: By setting $\Omega_I\coloneqq -\Omega_R(J\cdot,\cdot)$ and by taking $\mH_I$ to be a primitive of $-d\mH_R\circ J$, one can define the complex tensors $\Omega\coloneqq \Omega_R + i\Omega_I$ and $\mH\coloneqq \mH_R + i\mH_I$. The form $\Omega$ is of type $(2,0)$ and non-degenerate on $T^{(1,0)}X$, while $\mH:X\to\C$ is a pseudo-holomorphic function. They allow us to define the vector field $X_\mH$ of type $(1,0)$ via $\iota_{X_\mH}\Omega = -d\mH$. One can show that the real and imaginary parts of $X_\mH$ commute allowing us to define pseudo-holomorphic trajectories $\gamma$ via the following equation:
\begin{gather*}
 \frac{\partial\gamma}{\partial z} (z)\coloneqq \frac{1}{2}\left(\frac{\partial \gamma}{\partial t} (z) - i\frac{\partial\gamma}{\partial s} (z)\right) = X_\mH (\gamma (z))\quad\text{for } z = t + is.
\end{gather*}
As for HHSs, the pseudo-holomorphic trajectories locally exist and are unique given an initial value. Again, the uniqueness does not extend to maximal trajectories, but can be restored by promoting the maximal trajectories to leaves of the foliation generated by the real and imaginary part of $X_\mH$. Furthermore, the pseudo-holomorphic trajectories obey an action principle. Indeed, it is possible to define, for instance, the action $\mathcal{A}^{\Rec}_\mH$ for PHHSs, since it only requires a primitive for $\Omega_R$, which, in contrast to $\Omega_I$, can be exact.\\
Still, there are crucial differences between HHSs and PHHSs: Since $J$ is not integrable, there is no notion of holomorphic charts or vector fields, so neither $X_\mH$ nor its induced foliation have a chance to be holomorphic. The foliation is at least pseudo-holomorphic in the sense that its leaves are pseudo-holomorphically embedded Riemann surfaces. On top of that, pseudo-holomorphic trajectories only depend smoothly on their initial value, while holomorphic trajectories do so holomorphically. This phenomenon can be traced back to the fact that neither real nor imaginary part of $X_\mH$ need to be $J$-invariant\footnote{A vector field $V$ is $J$-invariant if $L_VJ = 0$ where $L_VJ$ is the Lie derivative of $J$ w.r.t. $V$.} (cf. Proposition \autoref{prop:J-preserving_criteria}).\\
The biggest difference, however, concerns $\Omega_I$: In stark contrast to HHSs, the form $\Omega_I$ associated with a PHHS is generally not closed. One might wonder whether there are at least proper\footnote{A PHHS is proper if it is not simultaneously a HHS.} PHHSs with closed form $\Omega_I$. Surprisingly, the answer to that question is negative. At first glance, this result appears to be very strange, as one would expect the integrability of $J$, not the closedness of $\Omega_I$ to separate HHSs and PHHSs. It turns out that, when it comes to the relation between HHSs and PHHSs, both are equivalent (cf. Theorem \autoref{thm:rel_HSM_PHSM} and Corollary \autoref{cor:rel_HHS_PHHS}):

\begin{theorem*}[Relation between HHSs and PHHSs]
 Let $(X,J;\Omega_R, \mH_R)$ be a PHHS with $2$-forms $\Omega_I\coloneqq -\Omega_R (J\cdot,\cdot)$ and $\Omega\coloneqq \Omega_R + i\Omega_I$ as well as a function $\mH\coloneqq \mH_R + i\mH_I$, where $\mH_I$ is a primitive of the $1$-form $\Omega_R (J(X^{\Omega_R}_{\mH_R}),\cdot)$. Then, the following statements are equivalent:
 \begin{enumerate}
  \item $(X,\Omega, \mH)$ is a HHS with complex structure $J$.
  \item $\Omega_I$ is closed, $d\Omega_I = 0$.
  \item $J$ is integrable.
 \end{enumerate}
\end{theorem*}

\subsubsection*{Construction of PHHSs}

More interesting than the properties of a PHHS might be the study of examples. Yet, it is astonishingly difficult to find examples of proper PHHSs. Contrary to HHSs, there are no standard examples like cotangent bundles and we cannot simply complexify a RHS. To solve this problem, we will introduce a method in the first subsection of \autoref{sec:construction_and_deformation} which generates PHHSs out of HHSs (cf. Proposition \autoref{prop:constructing_PHHS_out_of_HHS}). The idea is to turn a HHS into a PHHS by twisting its complex structure $J$ with an appropriate $(1,1)$-tensor $A$. Specifically, if $(X,\Omega,\mH)$ is a HHS with complex structure $J$ and the $(1,1)$-tensor $A$ satisfies $\Omega_R (A\cdot,A\cdot) = \Omega_R$, then the twisted tensor $J_A\coloneqq AJA^{-1}$ is a generally non-integrable almost complex structure and fulfills $\Omega_R (J_A\cdot,J_A\cdot) = -\Omega_R$. Thus, $(X,J_A;\Omega_R,\mH_R)$ is a PHHS if $d\mH_R\circ J_A$ is exact.\\
With Proposition \autoref{prop:constructing_PHHS_out_of_HHS}, the problem of finding PHHSs reduces to the problem of finding suitable tensors $A$. To carry out this task, it is often convenient to assume that $A$ itself is an almost complex structure and choose $A$ by fixing the semi-Riemannian metric $g\coloneqq \Omega_R (\cdot,A\cdot)$. In \autoref{sec:construction_and_deformation}, we will demonstrate how this procedure works by applying Proposition \autoref{prop:constructing_PHHS_out_of_HHS} to the simplest non-trivial example: the HHS with $X = \C^2$, $J = i$, $\Omega = dz_2\wedge dz_1$, and $\mH = iz_1$ (cf. Example \autoref{ex:constructing_PHHS}). Here, we take $A = I_g$ to be an almost complex structure determined by the metric $g$:
\begin{gather*}
 g(\pa{x_1},\pa{x_1}) = g(\pa{x_2},\pa{x_2})^{-1} = f,\quad g(\pa{y_1},\pa{y_1}) = g(\pa{y_2},\pa{y_2})^{-1} = h,\\
 g(\pa{x_1},\pa{x_2}) = g(\pa{y_1},\pa{y_2}) = g(\pa{x_i},\pa{y_j}) = 0,
\end{gather*}
where $f,h:\mathbb{C}^2\to\mathbb{R}$ are smooth, nowhere-vanishing functions. The twisted almost complex structure $J_g\coloneqq -J_A = I_gJI_g$ then only depends on the quotient $r\coloneqq f/h$:
\begin{gather*}
 J_g (\pa{x_1}) = r\pa{y_1},\ J_g (\pa{x_2}) = r^{-1}\pa{y_2},\quad J_g (\pa{y_1}) = -r^{-1}\pa{x_1},\ J_g (\pa{y_2}) = -r\pa{x_2}.
\end{gather*}
For the given choices, $d\mH_R\circ J_g$ is exact if and only if $r$ depends solely on $x_1$. We note that the complex structure $J$ is unchanged ($J_g = J$) if we set $r\equiv 1$. In particular, the alteration of $J$ is purely local if we assume that $r$ only deviates from $1$ within a small neighborhood.\\
As it turns out, the twisting method from Proposition \autoref{prop:constructing_PHHS_out_of_HHS} is related to Hyperkähler structures. Indeed, both have a very similar setup: From the symplectic viewpoint, a Hyperkähler manifold is a symplectic manifold equipped with two anticommuting complex structures where one is compatible with the symplectic form, while the other one is anticompatible (cf. \autoref{app:kaehler}). The twisting method, on the other hand, also involves a symplectic form - the form $\Omega_R$ - and two $(1,1)$-tensors $J$ and $A$ which are usually almost complex structures. The complex structure $J$ is anticompatible with $\Omega_R$ meaning $\Omega_R (J\cdot,J\cdot) = -\Omega_R$, whereas $A = I_g$ satisfies $\Omega_R (I_g\cdot,I_g\cdot) = \Omega_R$. The setup from Proposition \autoref{prop:constructing_PHHS_out_of_HHS} and the Hyperkähler setup solely differ by the signature of the metric $g$, the integrability of $I_g$, and the commutation relation between the almost complex structures. The commutation relation is the most central difference, since if $I_g$ and $J$ commuted or anticommuted, $J$ and $J_g$ would at best differ by a sign. Still, we can interpret the twisting method as some sort of deformation of a Hyperkähler structure. Example \autoref{ex:constructing_PHHS} demonstrates this beautifully: If the quotient $r$ is just $1$, then $\Omega_R$, $I_g$, and $J$ give rise to the standard Hyperkähler structure on $\C^2\cong\Qua$. As $r$ moves away from $1$, the Hyperkähler structure breaks down, but the twist $J\mapsto J_g = I_gJI_g$ is no longer a trivial operation.

\subsubsection*{Genericity of PHHSs}

To conclude the discussion of Hamiltonian systems, we investigate how ``large'' the set of proper PHHSs is within the set of all PHHSs. We will prove in the second subsection of \autoref{sec:construction_and_deformation} that the set of proper PHHSs is open and dense implying that being a proper PHHS is a generic property (cf. Theorem \autoref{thm:generic}):

\begin{theorem*}[Proper PHHSs are generic]
 Let $X$ be a smooth manifold, then the following statements apply depending on the real dimension of $X$:
 \begin{enumerate}
  \item If $\text{\normalfont dim}_\mathbb{R}(X) = 2$: Every almost complex structure on $X$ is integrable and automatically a complex structure.
  \item If $\text{\normalfont dim}_\mathbb{R}(X) > 2$: Every complex manifold $(X,J)$ and HSM\footnote{``HSM'' and ``PHSM'' stand for ``holomorphic symplectic manifold'' and ``pseudo-holomorphic symplectic manifold''. A PHSM is a symplectic manifold $(X,\Omega_R)$ equipped with an almost complex structure $J$ such that $\Omega_R (J\cdot,J\cdot) = -\Omega_R$ holds. For the definition of a proper deformation, confer Definition \autoref{def:deformation}.} $(X,\Omega)$ admits a proper deformation. In particular, the non-integrable almost complex structures and the proper PHSMs on $X$ are generic within the set of all almost complex structures and PHSMs on $X$, respectively.
  \item If $\text{\normalfont dim}_\mathbb{R}(X) > 4$: Every HHS $(X,\Omega, \mH)$ admits a proper deformation. In particular, the proper PHHSs on $X$ are generic within the set of all PHHSs on $X$.
 \end{enumerate}
\end{theorem*}

The idea behind the proof is to show that any HHS can be turned into a proper PHHS by an arbitrarily small perturbation. To find the perturbations, we use the twisting method again. For this, we first prove that any regular HHS can locally be brought into standard form meaning $\Omega = \sum_j dz_{j+n}\wedge dz_j$ and $\mH = z_{2n}$ (cf. Lemma \autoref{lem:HHS_in_standard_form}). Afterwards, we alter the $(x_1,y_1,x_2,x_2)$-components of $J$ within this small neighborhood similar to Example \autoref{ex:constructing_PHHS}. If the dimension of $X$ is sufficiently large ($\dim_R X = 2n>4$), $\mH:X\to\C$, which only depends on the coordinates $(x_{2n},y_{2n})$, is still a pseudo-holomorphic function giving us a proper PHHS.

\subsubsection*{Holomorphic Kähler Structure of Coadjoint Orbits}

As already mentioned in the beginning, \autoref{chap:duality} is dedicated to the Kähler structure of coadjoint orbits\footnote{To be more precise, we will discuss adjoint and coadjoint orbits. Since both are isomorphic via a suitable chosen metric in the cases we are interested in, we will only talk about coadjoint orbits. We will adopt a similar convention for tangent and cotangent bundles.}. It was first noted in the 50s that coadjoint orbits of compact Lie groups are compact homogeneous Kähler manifolds. The interest in these Kähler structures was mainly fueled by two curious observations (cf. Chapter 8 of \cite{Besse2007}):
\begin{enumerate}[label = (\arabic*)]
 \item Coadjoint orbits of compact Lie groups are Kähler-Einstein manifolds with positive scalar curvature, in fact, they were the first examples to be found\footnote{Not too long ago, one of the biggest unsolved problems in the field of Kähler geometry concerned the existence of Kähler-Einstein metrics on compact Kähler manifolds with prescribed first Chern class. In the case that the first Chern class/curvature is negative or vanishes, this problem was known as the Calabi conjecture and solved by Yau in the 70s. The Fano case (positive Chern class/curvature) remained open until about ten years ago, when Chen, Donaldson, and Sun found a solution (cf. \cite{ChenI}\cite{ChenII}\cite{ChenIII}).}.
 \item All simply-connected compact homogeneous Kähler manifolds are coadjoint orbits.
\end{enumerate}
At the end of the 80s, Kronheimer rekindled interest in this topic when he posed the question whether coadjoint orbits of complex reductive groups, the complexifications of compact Lie groups, also exhibit some sort of Kähler structure. He and Kovalev were able to show that the coadjoint orbits of semisimple complex reductive groups are Hyperkähler manifolds (cf. \cite{Kronheimer1990} and \cite{Kovalev1996}). Their idea was to identify these orbits with moduli spaces of instantons, i.e., spaces of anti-self-dual connection $1$-forms modulo gauge transformations, which were known to possess Hyperkähler structures.\\
In \autoref{chap:duality}, we will show that, additionally, coadjoint orbits of complex reductive groups admit a new type of Kähler structure which we call \textbf{holomorphic Kähler structure}. To understand the notion of a holomorphic Kähler manifold, we first need to explain semi-Kähler and holomorphic semi-Kähler structures. Simply put, a \textbf{semi-Kähler manifold} $(M,\omega,J)$ is a Kähler manifold where the metric $g\coloneqq \omega (\cdot,J\cdot)$ does not need to be positive definite anymore. A \textbf{holomorphic semi-Kähler manifold} $(X,\omega,J,I)$ consists of a semi-Kähler manifold $(X,\omega,J)$ and a complex structure $I$ satisfying:
\begin{enumerate}
 \item $\omega (I\cdot,I\cdot) = -\omega$ and $IJ = JI$.
 \item $\Omega\coloneqq \omega -i\omega(I\cdot,\cdot)$ and $J$ viewed as a section of $\End(T^{(1,0)}_IX)$ are holomorphic.
\end{enumerate}
Now, a holomorphic Kähler manifold is a holomorphic semi-Kähler manifold equipped with a special real structure. We call $\sigma:X\to X$ a \textbf{real structure}\linebreak on $(X,\omega,J,I)$ if it is a smooth involution which is $J$-holomorphic,\linebreak $I$-antiholomorphic, and leaves $\omega$ invariant meaning $\sigma^\ast\omega = \omega$. The quadruple $(X,\omega,J,I)$ becomes a holomorphic Kähler manifold if the semi-Kähler structure $(X,\omega,J)$ restricts to a Kähler structure on the real form $M\coloneqq\Fix\sigma$. The reason for imposing these conditions should be obvious: They allow us to interpret $(X,\Omega,J)$ as a complexification of the Kähler manifold $(M,\omega\vert_M,J\vert_M)$ with respect to the complex structure $I$.\\
The construction of holomorphic Kähler structures on coadjoint orbits of complex reductive groups follows more or less the construction of Kähler structures on coadjoint orbits of compact groups (cf. Chapter 8 of \cite{Besse2007}). In both cases, the idea is to first define a canonical complex structure $J$ on the adjoint orbit $\mathcal{O}$ and a canonical symplectic structure $\kks$ on the coadjoint orbit $\mathcal{O}^\ast$. Afterwards, one identifies $\mathcal{O}$ with $\mathcal{O}^\ast$ via an Ad-invariant scalar product.\\
Let us begin with the compact case. Consider a compact Lie group $G$. The differential of the adjoint representation $\Ad:G\to\GL (\mathfrak{g})$ is $\ad:\mathfrak{g}\to\End (\mathfrak{g})$, $\ad_vw = [v,w]$ which implies $T_w\mathcal{O} = \im\ad_w$. Thus, it suffices to define $J_w$ on $\im\ad_w$. Every compact Lie group admits a bi-invariant Riemannian metric $g$ or, equivalently, a positive definite\footnote{For us, a scalar product is just a non-degenerate symmetric bilinear form.} $\Ad$-invariant scalar product $\skcdot$ on $\mathfrak{g}$. One infers from the $\Ad$-invariance of $\skcdot$ that $\ad_w$ is skew-symmetric meaning $\sk{\ad_wu}{v} = -\sk{u}{\ad_wv}$. The spectral theorem now gives us $\Spec\ad_w\subset i\R$ allowing us to define $J_w$:
\begin{gather*}
 J_w v\coloneqq \frac{1}{\mu} \ad_w v\quad\forall v\in E_\mu,
\end{gather*}
where we set $E_\mu\coloneqq \mathfrak{g}\cap (E_{i\mu}\oplus E_{-i\mu})$ and $E_{i\mu},\, E_{-i\mu}\subset\mathfrak{g}_\C$ are the eigenspaces of $\ad_w$ for the eigenvalues $i\mu,\, -i\mu$ ($\mu>0$), respectively. The symplectic form $\kks$ is the famous Kirillov-Kostant-Souriau form. It is induced by the canonical Poisson structure on the dual Lie algebra $\mathfrak{g}^\ast$ (cf. Proposition \autoref{prop:dual_lie}):
\begin{gather*}
 \{F,G\} (\alpha)\coloneqq \alpha \left(\left[dF_\alpha, dG_\alpha\right]\right),
\end{gather*}
where $F,G\in C^\infty (\mathfrak{g}^\ast)$ and $\alpha\in\mathfrak{g}^\ast$. Setting $X^\ast_v (\alpha)\coloneqq -\alpha\circ\ad_v$ to be the fundamental vector field of the coadjoint action with respect to the vector $v\in\mathfrak{g}$, we can express $\kks$ as:
\begin{gather*}
 \omega_{\text{KKS},\alpha} (X^\ast_v (\alpha), X^\ast_w (\alpha)) = \alpha ([v,w])\quad\forall\alpha\in\mathcal{O}^\ast\, \forall v,w\in\mathfrak{g}.
\end{gather*}
We now use the $\Ad$-invariant scalar product $\skcdot$ from before to identify $\mathcal{O}$ with $\mathcal{O}^\ast$, i.e., $\mathcal{O}\to\mathcal{O}^\ast$, $w\mapsto\sk{w}{\cdot}$. It turns out that, under this identification, $\kks (\cdot,J\cdot)$ is a Riemannian metric giving us a Kähler structure on $\mathcal{O}\cong\mathcal{O}^\ast$.\\
To generalize this construction to complex reductive groups, we have to determine where the compactness of $G$ enters the construction. Going through the construction step by step, we realize that the compactness of $G$ is only required to guarantee the existence of a bi-invariant Riemannian metric $g$ which itself is only needed in two places:
\begin{enumerate}
 \item It ensures $\Spec\ad_w\subset i\R$ which allows us to define $J_w$.
 \item It lets us identify $\mathcal{O}$ with $\mathcal{O}^\ast$.
\end{enumerate}
Hence, any coadjoint orbit admits a Kähler structure as long as the orbit fulfills (i) and (ii). If an adjoint orbit $\mathcal{O}$ satisfies (i) for one and, thus, for any point $w\in\mathcal{O}$, we call the orbit skew-symmetric. For (ii), it suffices if the Lie algebra $\mathfrak{g}$ admits a non-degenerate $\Ad$-invariant scalar product $\skcdot$. In particular, we do not require $\skcdot$ to be positive definite. However, the trade-off is that the resulting structure is only semi-Kähler (cf. Theorem \autoref{thm:semi-kaehler}):

\begin{theorem*}[Semi-Kähler structures on (co)adjoint orbits]
 Let $G$ be a Lie group with Lie algebra $\mathfrak{g}$, dual Lie algebra $\mathfrak{g}^\ast$, and $\Ad$-invariant, non-degenerate scalar product $\sk{\cdot}{\cdot}$ on $\mathfrak{g}$. Further, let $\mathcal{O}\subset\mathfrak{g}$ be a skew-symmetric adjoint orbit of $G$. Then, $\mathcal{O}\subset\mathfrak{g}$ is an immersed submanifold and carries a $G$-invariant semi-Kähler structure. Its complex structure $J$ is the canonical complex structure on skew-symmetric adjoint orbits. If we identify $\mathcal{O}$ via $\sk{\cdot}{\cdot}$ with the coadjoint orbit $\mathcal{O}^\ast\subset\mathfrak{g}^\ast$, its symplectic form becomes the Kirillov-Kostant-Souriau form $\kks$. Moreover, $\mathcal{O}$ is Kähler if $\skcdot$ is positive definite.
\end{theorem*}

If, additionally, $G$ is a complex Lie group, then the semi-Kähler structure is even holomorphic. Indeed, the complex Lie group $G$ is equipped with a complex structure $I$ which descends to the complex structure $I_e$\footnote{We denote the neutral element of $G$ by $e$. Also note that we are a bit vague here. Precisely speaking, $I_e:\mathfrak{g}\to\mathfrak{g}$ is just a complex structure of the vector space $\mathfrak{g}$. However, it restricts to the map $I_{e,w}:\im\ad_w\to\im\ad_w$ at every point $w\in\mathcal{O}$ inducing a complex structure on $\mathcal{O}$.} on the adjoint orbit $\mathcal{O}$. We will check in \autoref{sec:holo_semi-kaehler} that $(\mathcal{O},\kks,J,I_e)$ constitutes a holomorphic semi-Kähler manifold (cf. Theorem \autoref{thm:complex_orbit}):

\begin{theorem*}[Holomorphic semi-Kähler structure on (co)adjoint orbits]
 Let $G$ be a Lie group that satisfies the conditions of the previous theorem, i.e.,\linebreak $G$ admits a non-degenerate $\Ad$-invariant scalar product and a skew-symmetric adjoint orbit $\mathcal{O}\subset\mathfrak{g}$. If $G$ is a complex Lie group with complex structure $I$, then $(\mathcal{O},\kks, J, I_e)$ is a $G$-invariant holomorphic semi-Kähler manifold.
\end{theorem*}

To show that coadjoint orbits of complex reductive groups exhibit holomorphic Kähler structures, we need to verify that such groups satisfy the prerequisites of Theorem \autoref{thm:complex_orbit} and that their orbits possess Kähler manifolds as their real forms. Complex reductive groups fulfill the conditions of Theorem \autoref{thm:complex_orbit}, since their real forms, compact Lie groups, do so. Indeed, the complexification of a skew-symmetric orbit is still skew-symmetric (cf. Proposition \autoref{prop:comp_of_skew-symmetric}), while the scalar product $\skcdot$ is obtained by complexifying the bi-invariant Riemannian metric $g$ of the compact real form (cf. Proposition \autoref{prop:comp_of_metric}). To prove the existence of Kähler real forms, we will show that the real structure of a complex reductive group descends to a real structure on its coadjoint orbits. The corresponding real forms are orbits of compact Lie groups which we know are Kähler manifolds implying the desired result (cf. Theorem \autoref{thm:complex_orbit_real_form} and Corollary \autoref{cor:reductive}):

\begin{theorem*}[(Co)Adjoint orbits of complex reductive groups]
 Let $G$ be a complex reductive group with real form $G_\R$, Lie algebras $\mathfrak{g}_\R\subset\mathfrak{g}$, and dual Lie algebras\footnote{$\mathfrak{g}^\ast_\R$ is the space of all linear maps $\mathfrak{g} = \mathfrak{g}_\R\oplus I_e\mathfrak{g}_\R\to\R$ that vanish on $I_e\mathfrak{g}_\R$.} $\mathfrak{g}^\ast_\R\subset\mathfrak{g}^\ast$. Further, let $\mathcal{O}\subset\mathfrak{g}$ be an adjoint orbit of $G$ and $\mathcal{O}^\ast\subset\mathfrak{g}^\ast$ a coadjoint orbit of $G$ such that $\mathcal{O}\cap\mathfrak{g}_\R\neq\emptyset\neq\mathcal{O}^\ast\cap\mathfrak{g}^\ast_\R$. Then, $\mathcal{O}$ and $\mathcal{O}^\ast$ carry $G$-invariant holomorphic Kähler structures.
\end{theorem*}

\subsubsection*{Kähler Duality}

We name the phenomenon that a space admits Hyperkähler and holomorphic Kähler structures \textbf{Kähler duality}. Combining Kronheimer's and our results, we see that coadjoint orbits of complex reductive groups are examples of Kähler duality. Exhibiting Kähler duality is a very curious property, especially from the symplectic viewpoint: Both Hyperkähler and holomorphic Kähler manifolds can be described by quadruples $(X,\omega,J,I)$, where $I$ and $J$ are complex structures and $\omega$ is a symplectic form satisfying $\omega (J\cdot,J\cdot) = \omega$ and $\omega (I\cdot,I\cdot) = -\omega$. The sole difference between the Hyperkähler and the holomorphic Kähler case is the commutation relation of $I$ and $J$. For Hyperkähler manifolds, $I$ and $J$ anticommute, while they commute for holomorphic Kähler manifolds. On coadjoint orbits, the similarities between Hyperkähler and holomorphic Kähler manifolds are even more striking. Indeed, Hyperkähler and holomorphic Kähler manifolds can be seen as holomorphic symplectic manifolds where the complex manifold $(X,I)$ is equipped with the holomorphic symplectic form $\Omega\coloneqq \omega - i\omega(I\cdot,\cdot)$. For coadjoint orbits, the form $\Omega$ is in both cases given by the holomorphic Kirillov-Kostant-Souriau form $\KKS$.\\
Given such a strong geometrical structure as Kähler duality, one is naturally drawn to the question whether the Kähler duality of coadjoint orbits is just accidental. In this thesis, we claim that the Kähler duality of coadjoint orbits can be traced back to double cotangent bundles $T^\ast (T^\ast M)$. More precisely, we conjecture that $T^\ast (T^\ast M)$ also exhibits Kähler duality and that if $M = G_\R$ is a compact Lie group with complexification $G$, the Kähler duality of $T^\ast (T^\ast M)$ is related to the Kähler duality of coadjoint orbits of $G$ via Hyperkähler/holomorphic Kähler reduction (outlined in \autoref{app:reduction}).\\
The Kähler duality of $T^\ast (T^\ast M)$ is based on a famous result\footnote{Confer Theorem \autoref{thm:cotangent_kaehler}. Also note that the different versions of Stenzel's theorem as well as the diagram later on have to be taken with a grain of salt: In general, these structures do not exist on all of $T^\ast M$ or $T^\ast (T^\ast M)$, but only on an open neighborhood of the zero section $M\subset T^\ast M\subset T^\ast (T^\ast M)$.} due to Guillemin-Stenzel (cf. \cite{Stenzel1990} and \cite{Guillemin1991}) and Lempert-Sz{\H{o}}ke (cf. \cite{Lempert1991} and \cite{Szoeke1991}): If $(M,g)$ is a real-analytic Riemannian manifold, then $(T^\ast M,-\omega_{\can}, J_g)$ is a Kähler manifold where $\omega_{\can}$ is the canonical symplectic form on $T^\ast M$ and $J_g$ is the unique complex structure on $T^\ast M$ adapted\footnote{Confer Definition \autoref{def:adapted_comp_str}.} to $g$. Stenzel's theorem can be modified depending on which structure the base manifold $M$ carries. In the case that the metric $g$ on $M$ is Kähler, Kaledin (cf. \cite{Kaledin1997}) and Feix (cf. \cite{Feix2001}) showed that $T^\ast M$ is actually Hyperkähler. In this spirit, we claim that a Hyperkähler structure on $T^\ast M$ is given by the quadruple $(T^\ast M,-\omega_{\can}, \phi^\ast J_g, T^\ast I)$ where $(M,g,I)$ is the base Kähler manifold and $\phi:T^\ast M\to T^\ast M$ is a real-analytic, fiber-preserving diffeomorphism with $\phi\vert_M = \id_M$ (cf. Lemma \autoref{lem:stenzel_for_kaehler}):

\begin{lemma*}[Stenzel's theorem for Kähler metrics]
 Let $(M,g, I)$ be a real-analytic Kähler manifold, let $J_g$ be the complex structure adapted to $g$, and let\linebreak $\phi:T^\ast M\to T^\ast M$ be a real-analytic, fiber-preserving diffeomorphism with $\phi\vert_M = \id_M$. If $\phi^\ast J_g$ and $T^\ast I$ anticommute and $-\omega_{\can} (\cdot,\phi^\ast J_g\cdot)$ is a Riemannian metric, then $(T^\ast M,-\omega_{\can}, \phi^\ast J_g, T^\ast I)$ is a Hyperkähler manifold.
\end{lemma*}

Furthermore, we believe that there is also a version of Stenzel's theorem for pairs $(M,g)$ where $g$ is the real part of a holomorphic metric (cf. Lemma \autoref{lem:stenzel_for_holo_metrics}):

\begin{lemma*}[Stenzel's theorem for holomorphic metrics]
 Let $(M,I)$ be a complex manifold with holomorphic metric $G = g-ig(I\cdot,\cdot)$, let $J_g$ be the complex structure adapted to $g$, and let $\phi:T^\ast M\to T^\ast M$ be a real-analytic, fiber-preserving diffeomorphism with $\phi\vert_M = \id_M$. If $\phi^\ast J_g$ and $T^\ast I$ commute and $-\omega_{\can} (\cdot,\phi^\ast J_g\cdot)$ is a semi-Riemannian metric, then $(T^\ast M,-\omega_{\can}, \phi^\ast J_g, T^\ast I)$ is a holomorphic semi-Kähler manifold. If, additionally, $\sigma$ is a real structure on the complex manifold $(M,I)$ with real form $M_\R$ satisfying $\overline{\sigma^\ast G} = G$, the induced metric $g_\R$ on $M_\R$ is positive definite, and $\phi$ can be chosen such that the real structure $\sigma^\ast (\alpha)\coloneqq \alpha\circ d\sigma$ on $T^\ast M$ is $\phi^\ast J_g$-holomorphic, then $(T^\ast M,-\omega_{\can}, \phi^\ast J_g, T^\ast I)$ is even holomorphic Kähler with respect to $\sigma^\ast$.
\end{lemma*}

With the different versions of Stenzel's theorem in mind, the Kähler duality of $T^\ast (T^\ast M)$ now arises as depicted in the following diagram:
\begin{center}
 \begin{tikzcd}[column sep=3.6em,row sep=2em]
  \textbf{Hyperk.}& (T^\ast M, g_{T^\ast M}) \arrow{r}{\text{Stenzel}} & (T^\ast (T^\ast M), -\omega_{\can}, \phi^\ast_1 J_{g_{T^\ast M}}, T^\ast J_g)\\
  (M,g) \arrow{ur}[sloped,above]{\text{Stenzel}} \arrow{dr}[sloped,below]{\text{\footnotesize Complexification}} & &\\
  \textbf{Holo. K.}& (T^\ast M, g_\C) \arrow{r}[below]{\text{Stenzel}} & (T^\ast (T^\ast M), -\omega_{\can}, \phi^\ast_2 J_{g_\C}, T^\ast J_g)
 \end{tikzcd}
\end{center}
Here, the upper path illustrates how to obtain a Hyperkähler structure, while the lower one does the same for holomorphic Kähler structures.\\\\
\textbf{Hyperkähler path:}\footnote{Note that Hyperkähler structures on double cotangent bundles were already described by Bielawski (cf. \cite{Bielawski2003}).} Starting with the real-analytic Riemannian manifold $(M,g)$, we can apply Stenzel's theorem to obtain the Kähler metric\linebreak $g_{T^\ast M}\coloneqq -\omega_{\can} (\cdot,J_g\cdot) =  \omega_{\can} (J_g\cdot,\cdot)$ on $T^\ast M$. Assuming that a map $\phi_1$ exists as in Lemma \autoref{lem:stenzel_for_kaehler}, we can now apply Stenzel's theorem for Kähler\linebreak metrics to $(T^\ast M, g_{T^\ast M}, J_g)$ yielding the Hyperkähler manifold\linebreak $(T^\ast (T^\ast M), -\omega_{\can}, \phi^\ast_1 J_{g_{T^\ast M}}, T^\ast J_g)$.\\\\
\textbf{Holomorphic Kähler path:} Stenzel's theorem does not only tell us that $(T^\ast M,-\omega_{\can}, J_g)$ is a Kähler manifold, but also that the fiberwise map\linebreak $T^\ast M\to T^\ast M$, $\alpha\mapsto -\alpha$ is a real structure (cf. Theorem \autoref{thm:cotangent_kaehler}). Its real form is the zero section $M\subset T^\ast M$. Since the real form $M$ carries a real-analytic metric $g$, we can find a unique holomorphic continuation of $g$ on $T^\ast M$ (cf. \autoref{app:real_structures}). Call the real part of this holomorphic metric $g_\C$. We are now in the setup of Lemma \autoref{lem:stenzel_for_holo_metrics} (Stenzel's theorem for holomorphic metrics) giving us the holomorphic Kähler manifold $(T^\ast (T^\ast M), -\omega_{\can}, \phi^\ast_2 J_{g_\C}, T^\ast J_g)$ where we assume, of course, that a map $\phi_2$ as specified above exists.\\\\
With the Kähler duality of $T^\ast (T^\ast M)$ taken care of, we need to relate $T^\ast (T^\ast M)$ to coadjoint orbits. For this, we choose $(M,g)$ to be a compact Lie group $G_\R$ with bi-invariant Riemannian metric $g$. In this case, $T^\ast G_\R$ is isomorphic to the universal complexification $G$ of $G_\R$. Thus, the previously described process yields Hyperkähler/holomorphic Kähler structures on $T^\ast G$. After reduction, the cotangent bundle $T^\ast G$ becomes a coadjoint orbit $\mathcal{O}^\ast$ of $G$. At the same time, the symplectic form $-\omega_{\can}$ on $T^\ast G$ reduces to the Kirillov-Kostant-Souriau form $\kks$ on $\mathcal{O}^\ast$. From that perspective, it seems plausible that the Kähler structures in question are also compatible with the reduction process.\\
Unfortunately, the proofs we present in this thesis regarding the Kähler duality of $T^\ast (T^\ast M)$ and its relation to coadjoint orbits are incomplete. Precisely speaking, we leave two questions unanswered. The first one regards the existence of the map $\phi$ and the commutation relations of $\phi^\ast J_g$ and $T^\ast I$ in Lemma \autoref{lem:stenzel_for_kaehler} and \autoref{lem:stenzel_for_holo_metrics}. It is vital for the modified versions of Stenzel's theorem that $\phi^\ast J_g$ and $T^\ast I$ anticommute in the Hyperkähler case and that they commute in the holomorphic Kähler case. However, a complete proof for these commutation relations is missing (cf. Conjecture \autoref{con:commutation_relation}). Secondly, we do not carry out the reduction process in detail (cf. Conjecture \autoref{con:reduction}). To make up for that, we check the commutation relations on the zero section and for flat $g$ (cf. Lemma \autoref{lem:commutation_relation}) and sketch how reduction could possibly relate a suitable Kähler structure on $T^\ast G$ to $\mathcal{O}^\ast$ (cf. \autoref{app:reduction}).

\subsubsection*{Structural Remarks}

Before we jump into the main part, we should add a few remarks regarding the structure of this thesis. Often, papers and theses contain a chapter on preliminaries outlining the basics of the topic in question. The author of this thesis is of the opinion that these chapters have a negative impact on the text as a whole: They tend to drag on, elaborate on details with which the educated reader is already familiar, and introduce concepts at the beginning which are only needed at the end. Therefore, we take a different approach in this thesis: In an effort to tell a compelling, but cohesive story, we start our explanations \textit{in medias res}, try to give background knowledge only when needed, and develop concepts as we move along\footnote{Unfortunately, this does not cure the author of his inability to write concise passages.}. If a preliminary section cannot be avoided, we indicate its purpose (only needed to fix notations etc) and point out that the sophisticated reader may skip the section, as done, for instance, at the beginning of \autoref{sec:lie_groups}. As everything in life, this approach comes with a trade-off. Not every reader has access to the same knowledge which is why some readers may struggle with certain parts. To offset this, sections which are usually part of a preliminary chapter are moved to the appendix, for example \autoref{app:real_structures} and \autoref{app:kaehler}. These appendices are always referenced when they become relevant, so the reader may consult them if they feel they lack the needed background knowledge\footnote{Of course, the reader is always invited to read the appendix just out of curiosity.}.
Writing is not the author's claim to fame, so it is left to the reader to decide whether our approach was successfully implemented.

\chapter{Pseudo-Holomorphic Hamiltonian Systems}
\label{chap:PHHS}
Hamiltonian systems (HSs) as the mathematical model for classical mechanics have been central to the advance of modern physics and mathematics alike. In physics, HSs provide a theoretical foundation for several approaches to quantization. In mathematics, the interest in HSs has led to the study of symplectic geometry and topology. The methods developed in this study, e.g. Floer theory, have proven to be of great success for various branches of mathematics and physics, for instance celestial mechanics and string theory.\\
Simply put, a HS consists of three data: A manifold $M$, a symplectic $2$-form $\omega$ on $M$, and function $H\in C^\infty (M,\mathbb{R})$. In physical terms, the symplectic manifold $(M,\omega)$ can be understood as the phase space of the system, while the function $H$, often called Hamilton function or simply Hamiltonian, assigns to every point in phase space its energy. These data allow us to define the Hamiltonian vector field $X_H$ on $M$ via the equation $\iota_{X_H}\omega = -dH$. The dynamics of the HS $(M,\omega, H)$ is governed by the vector field $X_H$. Precisely speaking, the physical trajectories of point-like particles described by the HS $(M,\omega, H)$ are exactly the integral curves of $X_H$. The connection between the integral curve equation of $X_H$ and the Hamilton equations known from classical mechanics is given by Darboux's theorem which states that $\omega$ can locally be written as:
\begin{gather*}
 \omega = \sum^n_{i = 1} dp_i\wedge dq_i.
\end{gather*}
In such Darboux charts, the integral curve equation of $X_H$ reduces to the Hamilton equations:
\begin{gather*}
 \dot q_i(t) = \frac{\partial H}{\partial p_i},\quad \dot p_i(t) = -\frac{\partial H}{\partial q_i}\quad\forall t\in I\ \forall i\in\{1,\ldots, m\}.
\end{gather*}
Since the integral curve equation is just a first-order ODE, there exists for every initial value $(t_0,x_0)\in\R\times M$ an open interval $I\ni t_0$ and a curve $\gamma:I\to M$ solving the Hamilton equations and satisfying $\gamma (t_0) = x_0$. Furthermore, two trajectories $\gamma_1:I_1\to M$ and $\gamma_2:I_2\to M$ are identical iff they have the same domain ($I_1 = I_2\equiv I$) and attain the same value at some point $t_0\in I$. In particular, maximal integral curves are unique given an initial value and foliate energy hypersurfaces $H^{-1}(E)$ for regular values $E$ of $H$.\\
On top of that, physical trajectories obey the action principle, i.e., they can be obtained as ``critical points'' of the action functional $\mathcal{A}_H:C^\infty(I,M)\to\mathbb{R}$ assigned to an exact HS $(M,\omega = d\lambda, H)$:
\begin{gather*}
 \mathcal{A}_H[\gamma]\equiv \mathcal{A}^\lambda_H[\gamma]\coloneqq \int\limits_I \gamma^\ast\lambda - \int\limits_I H\circ\gamma (t)\, dt.
\end{gather*}
Here, ``critical point'' means that the first variation of $\mathcal{A}_H$ has to vanish at a physical trajectory $\gamma\in C^\infty(I,M)$, where we only allow for variations of $\gamma$ which keep the endpoints of $\gamma$ fixed. Sometimes, for instance in Floer theory, one wishes to view certain trajectories as actual critical points of some action functional. In this case, we have to ensure that the boundary terms vanish. There are several ways to achieve this, e.g. by putting the endpoints of a trajectory on an exact Lagrangian ($\lambda_{\gamma(t_1)} = 0$ and $\lambda_{\gamma(t_2)} = 0$ for $I = [t_1,t_2]$) or by only considering periodic trajectories (cf. Remark \autoref{rem:critical_point}).\\
Since HSs are given in terms of real manifolds $M$, forms $\omega$, and functions $H$, it is only natural to ask whether a similar construction with similar properties exists for complex manifolds $X$, forms $\Omega$, and functions $\mathcal{H}$. This question directly leads us to the notion of \textbf{holomorphic Hamiltonian systems} (HHSs). Similarly to real Hamiltonian systems\footnote{To distinguish real and complex Hamiltonian systems, we call HSs $(M,\omega,H)$ real Hamiltonian systems from now on.} (RHSs), HHSs are also described by three data (cf. \autoref{sec:HHS}): A complex manifold $X$ (implicitly defining an integrable complex structure $J$), a holomorphic symplectic $2$-form $\Omega$ on $X$, and a holomorphic function $\mH:X\to\mathbb{C}$. HHSs have been studied since the early 2000s, e.g. by Gerdjikov and Kyuldjiev \cite{gerd2001}, \cite{gerd2002}, \cite{gerd2004} or by Arathoon and Fontaine \cite{arathoon2020}. In the given references, HHSs are usually viewed as complexifications of RHSs and mostly used as a tool to study RHSs which arise as real forms\footnote{A real structure on a HHS $(X,\Omega,\mathcal{H})$ is an antiholomorphic involution $\sigma:X\to X$ such that $\overline{\sigma^\ast\Omega} = \Omega$ and $\overline{\sigma^\ast\mathcal{H}} = \mathcal{H}$. Its real form is $(M,\omega, H)$ with $M\coloneqq\Fix\sigma$, $\omega\coloneqq\iota^\ast\Omega$, and $H\coloneqq\iota^\ast\mathcal{H}$, where $\iota:M\hookrightarrow X$ denotes the natural inclusion. We say that $(X,\Omega,\mathcal{H})$ is a complexification of $(M,\omega,H)$ if $M$ is nice, i.e., $M$ intersects every connected component of $X$ non-trivially. Confer \autoref{app:real_structures} for a detailed discussion of real structures.} of HHSs. In \cite{arathoon2020}, for instance, the authors find a compact and integrable real form of the complexified spherical pendulum.\\
In this chapter, we take a different approach. We study HHSs on their own and try to recreate the results known from RHSs for HHSs. \autoref{chap:PHHS} is divided into three sections which themselves are split up into multiple subsections. We begin \autoref{sec:HHS} by introducing HHSs and examining their most basic properties. In the next subsection, we discuss the existence and uniqueness of holomorphic trajectories. Similarly to RHSs, holomorphic trajectories are defined as the holomorphic integral curves of the holomorphic Hamiltonian vector field $X_\mH$. We show that, locally, holomorphic trajectories always exist and are unique, given an initial value. In sharp contrast to RHSs, maximal holomorphic trajectories are not unique anymore, even given an initial value, due to the effects of monodromy\footnote{Recently, the monodromy of the complexified Kepler problem has been studied by Sun and You (cf. \cite{shanzhong2020}).}. Nevertheless, holomorphic trajectories still foliate the energy hypersurface $\mH^{-1}(E)$ for any regular value $E$ of $\mH$. In fact, holomorphic trajectories give us a holomorphic foliation of the entire manifold $X$, if all possible energy values are regular. This foliation shows up again in the following subsection, where we use it to establish a relation between Lefschetz fibrations and almost toric fibrations.\\
We prove in the last subsection of \autoref{sec:HHS} that the holomorphic trajectories satisfy an action principle\footnote{To the extent of the author's knowledge, an action principle for HHSs has not been formulated before in the literature.}, i.e., that they can be understood -- in some sense -- as critical points of certain action functionals. These action functionals are obtained by first decomposing a HHS $(X,\Omega,\mH)$ into four RHSs, one for each combination of real and imaginary part of $\Omega$ and $\mH$. To each RHS, we can assign the usual action functional of a RHS. Afterwards, we average each of these action functionals over the imaginary (or real) time axis and take an appropriate linear combination to obtain the action functional for the HHS $(X,\Omega,\mH)$. In fact, this method gives rise to a plethora of action functionals for the HHS $(X,\Omega,\mH)$ which differ by how one averages and takes the linear combination.\\
During the investigation of action functionals, we observe that $J$, the complex structure of $X$, poses rather strong restrictions on the existence of certain holomorphic trajectories. At the end of \autoref{sec:HHS}, we consider holomorphic trajectories whose domains are complex tori and interpret them as the complexification of periodic orbits. However, by the maximum principle, such holomorphic trajectories are always constant if the complex manifold in question is $X=\mathbb{C}^{2n}$ equipped with the standard complex structure $J = i$. The same argument does not hold anymore if we allow $J$ to be any almost complex structure. In his beautiful paper \cite{moser1995} from 1995, Moser shows that it is possible to pseudo-holomorphically embed complex tori in $\mathbb{R}^4$, where $\mathbb{R}^4$ is equipped with a suitable, not necessarily integrable almost complex structure $J$.\\
To avoid constraints imposed by the integrability of $J$, we introduce special Hamiltonian systems in \autoref{sec:PHHS} which are described by the same data as HHSs, but whose almost complex structure $J$ does not need to be integrable anymore. These Hamiltonian systems are called \textbf{pseudo-holomorphic Hamiltonian systems} (PHHSs) and exhibit, by design, the same properties as HHSs. In particular, pseudo-holomorphic trajectories of PHHSs induce foliations of regular energy hypersurfaces $\mH^{-1}(E)$ and obey an action principle (cf. the first subsection of \autoref{sec:PHHS}). At first glance, PHHSs may appear to be contrived and artificial, especially since, by definition, the imaginary part of $\Omega$ does not need to be closed anymore. However, the non-closedness of the imaginary part of $\Omega$ is an unavoidable consequence of the non-integrability of $J$, as we show in the second subsection of \autoref{sec:PHHS}. In fact, we prove that we recover a HHS from a PHHS iff $J$ is integrable or, equivalently, $\text{Im}\,\Omega$ is closed.\\
In \autoref{sec:construction_and_deformation}, we tackle the question whether PHHSs are a natural generalization of HHSs. We claim that they are and support this assertion by showing that the space of proper\footnote{A proper PHHS is a PHHS which is not simultaneously a HHS.} PHHSs is open and dense in the space of PHHSs\footnote{At least on manifolds $X$ with $\dim_\R (X) >4$.}. This implies that proper PHHSs are generic. To prove that, we first give a method to construct proper PHHSs out of HHSs (cf. the first subsection of \autoref{sec:construction_and_deformation}). The method itself is very interesting, since it relates PHHSs to Hyperkähler structures and allows us to equip the cotangent bundle of a complex manifold with the structure of a PHHS. Lastly, we use this construction to deform HHSs by proper PHHSs (cf. the second subsection of \autoref{sec:construction_and_deformation}).

\newpage
\section{Holomorphic Hamiltonian Systems}
\label{sec:HHS}
In this section, we recreate the results known from real Hamiltonian systems (laid out in the introduction of \autoref{chap:PHHS}) for holomorphic Hamiltonian systems (HHSs). We first cover the basic notion of a HHS, including the definition of a holomorphic symplectic manifold and Darboux's theorem for holomorphic symplectic manifolds. Afterwards, we discuss the properties of holomorphic trajectories which, as we will see, give rise to a holomorphic foliation. Then, we utilize this holomorphic foliation to investigate the relation between Lefschetz and almost toric fibrations. Lastly, formulate an action principle for HHSs.

\subsection*{HHS: Basic Definitions and Notions}

To begin with, we define a holomorphic symplectic manifold:

\begin{definition}[Holomorphic symplectic manifold]\label{def:holo_sym_man}
 A pair $(X,\Omega)$ is called \textbf{holomorphic symplectic manifold}\footnote{Warning: In some branches of algebraic geometry, the term ``holomorphic symplectic manifold'' is also used, but defined with additional constraints on $X$ and $\Omega$!} (HSM) if $X$ is a complex manifold and $\Omega$ is a holomorphic $2$-form on $X$ which is closed and non-degenerate on the $(1,0)$-tangent bundle $T^{(1,0)}X$ of $X$. In this setup, $\Omega$ is called the holomorphic symplectic $2$-form of $(X,\Omega)$.
\end{definition}

Let us spend some time understanding the definition of a HSM. Recall that a complex manifold $X$ is defined via an atlas $\{(\phi_\alpha, U_\alpha)\}_{\alpha\in I}$ of charts with values in $\mathbb{C}^m$ such that their transition functions are holomorphic. Pick such a holomorphic chart $\phi = (z_1,\ldots, z_m):U\to V\subset\mathbb{C}^m$. A complex-valued $2$-form $\Omega$ is holomorphic on $U$ if and only if $\Omega$ can be written as:
\begin{gather*}
 \Omega\vert_U = \sum^{m}_{i,j = 1} \Omega_{ij} dz_i\wedge dz_j,
\end{gather*}
where $\Omega_{ij}:U\to\mathbb{C}$ are holomorphic functions\footnote{A function $f:X\to\C$ defined on a complex manifold $X$ is holomorphic if $f\circ\phi^{-1}_\alpha$ is holomorphic for every chart $\phi_\alpha$ in a holomorphic atlas.}. Here, we have used the fact that one can define the exterior derivative $d$ for complex-valued functions and forms by requiring $d$ to be $\C$-linear. Since the transition functions are holomorphic, the notion of holomorphicity is independent of the choice of chart. As in the real case, closedness of $\Omega$ simply means $d\Omega = 0$. To understand the non-degeneracy in Definition \autoref{def:holo_sym_man}, recall that every complex manifold $X$ implicitly defines a complex structure\footnote{$J$ is equal to the standard complex structure $i$ in a holomorphic chart.} $J$ on $X$ and, further, that the complexified tangent and cotangent bundle of $X$, viewed as a real manifold, each decompose into a direct sum of two subbundles:
\begin{gather*}
 T_\mathbb{C}X = T^{(1,0)}X\oplus T^{(0,1)}X,\quad T^\ast_\mathbb{C}X = T^{\ast, (1,0)}X\oplus T^{\ast, (0,1)}X,
\end{gather*}
where the $(1,0)$- and $(0,1)$-bundles are fiberwise eigenspaces of $J$ (or its dual $J^\ast$) with eigenvalue $i$ and $-i$, respectively. By construction, the local forms $dz_i$ are local sections of $T^{\ast, (1,0)}X$ and, hence, map elements of $T^{(0,1)}X$ to zero. This implies that holomorphic $2$-forms can never be non-degenerate on the entire complexified tangent bundle, as they always vanish on the $(0,1)$-bundle. For a holomorphic $2$-form $\Omega$, we can at most achieve non-degeneracy on $T^{(1,0)}X$, i.e.:
\begin{gather*}
 \forall x\in X\ \forall V\in T^{(1,0)}_xX\backslash\{0\}\ \exists W\in T^{(1,0)}_xX:\quad \Omega_x(V, W) \neq 0.
\end{gather*}
In particular, non-degeneracy of $\Omega$ implies that the complex dimension $m$ of $X$ is even. Also note that, by construction, $\Omega$ and $J$ satisfy the following relations:
\begin{gather}\label{eq:J-anticompatible}
 \Omega (J\cdot, \cdot) = \Omega (\cdot, J\cdot) = i\Omega,\quad \Omega (J\cdot, J\cdot) = -\Omega.
\end{gather}
Similar to ``real'' symplectic geometry, there are two standard examples of HSMs. The first one is $X = \mathbb{C}^{2n}$ together with the standard form $\Omega = \sum^n_{j = 1} dP_j\wedge dQ_j$,\linebreak where $(Q_1,\ldots, Q_n, P_1,\ldots, P_n)\in\mathbb{C}^{2n}$. The other one is the holomorphic cotangent bundle $X = T^{\ast, (1,0)}Y$ of a complex manifold $Y$ with canonical $2$-form\linebreak $\Omega_\text{can} = d\Lambda_\text{can}$, where $\Lambda_\text{can}$ is the holomorphic Liouville $1$-form:
\begin{gather*}
 \Lambda_{\can,\alpha} (v)\coloneqq \alpha\circ d\pi (v)\quad\forall v\in T_\alpha X\ \forall \alpha\in X.
\end{gather*}
Here, $\pi:X\to Y$ is the canonical projection.\\
We know by Darboux's theorem that symplectic manifolds exhibit no local invariants, since they all are locally isomorphic to the standard symplectic manifold $(\mathbb{R}^{2n}, \sum^n_{i = 1} dp_i\wedge dq_i)$. The same statement is true for HSMs:

\begin{theorem}[Darboux's theorem for HSMs]\label{thm:holo_Darboux}
 Let $(X,\Omega)$ be a HSM of complex dimension $\text{\normalfont dim}_\mathbb{C}(X) = 2n$ $(n\in\mathbb{N})$. Then, for every point $x\in X$, there exists a holomorphic chart $\psi = (Q_1,\ldots, Q_n, P_1,\ldots, P_n):U\to V\subset\mathbb{C}^{2n}$ of $X$ near $x$ such that:
 \begin{gather*}
  \Omega\vert_{U} = \sum\limits^n_{j = 1} dP_j\wedge dQ_j.
 \end{gather*}
\end{theorem}

The correctness of Darboux's theorem for HSMs is widely accepted in the mathematical community, however, no proof has yet been formally written down, at least to the author's extent of knowledge. For completeness' sake, a proof of Darboux's theorem for HSMs is provided in \autoref{app:darboux}. We will make use of Theorem \autoref{thm:holo_Darboux} in \autoref{sec:construction_and_deformation}.\\
Now, let us turn our attention to holomorphic Hamiltonian systems:

\begin{definition}[Holomorphic Hamiltonian system]\label{def:holo_ham_sys}
 We call the triple\linebreak $(X,\Omega, \mathcal{H})$ a \textbf{holomorphic Hamiltonian system} (HHS) if $(X,\Omega)$ is a HSM and $\mH:X\to\mathbb{C}$ is a holomorphic function on $X$. In this setup, we call $\mH$ the \textbf{Hamilton function} or, simply, the Hamiltonian of the HHS $(X,\Omega, \mH)$.
\end{definition}

Examples of HHSs include the standard HSM $(\mathbb{C}^{2n}, \sum^n_{j = 1} dP_j\wedge dQ_j)$ together with any holomorphic function $\mH:\mathbb{C}^{2n}\to\mathbb{C}$ on it and holomorphic cotangent bundles $(T^{\ast, (1,0)}Y, \Omega_\text{can})$ together with natural Hamiltonians. We call a Hamiltonian on a holomorphic cotangent bundle \textbf{natural} if it can be written as a sum of kinetic and potential energy, $\mH = \mathcal{T} + \mathcal{V}$. A holomorphic function\linebreak $\mathcal{V}:T^{\ast, (1,0)}Y\to\mathbb{C}$ denotes \textbf{potential energy} if it factors through the\linebreak holomorphic projection $\pi:T^{\ast, (1,0)}Y\to Y$, i.e., can be written as\linebreak $\mathcal{V} = \mathcal{V}_0\circ\pi$ for some holomorphic function $\mathcal{V}_0:Y\to \mathbb{C}$. Furthermore, a\linebreak holomorphic function $\mathcal{T}:T^{\ast, (1,0)}Y\to\mathbb{C}$ is called \textbf{kinetic energy} if\linebreak $2\mathcal{T}(x)\equiv g^\ast (x,x)$ $\forall x\in T^{\ast, (1,0)}Y$, where $g^\ast$ is the dualization of some holomorphic metric $g$ on $Y$. A holomorphic metric $g$ on $Y$ is a holomorphic symmetric $\mathbb{C}$-bilinear form which is non-degenerate on $T^{(1,0)}Y$, i.e., for any holomorphic chart $\phi = (z_1,\ldots, z_n):U\to V\subset\mathbb{C}^n$ of $Y$, $g$ can be written as:
\begin{gather*}
 g\vert_U = \sum^n_{i,j = 1} g_{ij} dz_i\otimes dz_j,
\end{gather*}
where $g_{ij}:U\to\mathbb{C}$ are holomorphic functions satisfying $g_{ij} = g_{ji}$ and\linebreak $\text{det}(g_{ij})\neq 0$. We will study these examples in more detail in the upcoming subsections.

\subsection*{Holomorphic Trajectories}

Our next goal is to investigate the dynamics of a HHS. As for RHSs, they are determined by a Hamiltonian vector field:

\begin{definition}[Holomorphic Hamiltonian vector field]\label{def:holo_ham_field}
 Let $(X,\Omega, \mH)$ be a HHS. We call the holomorphic vector field $X_\mH$ on $X$ defined by $\iota_{X_\mH}\Omega = -d\mH$ the (holomorphic) \textbf{Hamiltonian vector field} of the HHS $(X,\Omega, \mH)$.
\end{definition}

\begin{remark}[$X_\mH$ is well-defined]\label{rem:ham_vec_field_well-def}
 Note that a complex vector field $V$ on a complex manifold $X$ is holomorphic if and only if, in any holomorphic chart $(z_1 = x_1 + iy_1,\ldots, z_{m} = x_m + iy_m):U\to V\subset\mathbb{C}^m$, it can be written as:
 \begin{gather*}
  V\vert_U = \sum^m_{j = 1} V_j\partial_{z_j} \equiv \sum^m_{j = 1} \frac{V_j}{2}\left(\pa{x_j} - i\pa{y_j}\right),
 \end{gather*}
 where $V_j:U\to\mathbb{C}$ are holomorphic functions on U. Thus, a holomorphic vector field on $X$ only attains values in the bundle $T^{(1,0)}X$. Together with the non-degeneracy of $\Omega$ on $T^{(1,0)}X$, this implies that the Hamiltonian vector field $X_\mH$ is well-defined.
\end{remark}

Before we define what a trajectory of a HHS is, it is wise to study $X_\mH$ or, better yet, holomorphic vector fields in general. The following proposition is a standard result from complex geometry:

\begin{proposition}[Holomorphic vector fields $\Leftrightarrow$ $J$-preserving vector fields]\label{prop:holo_vec_field_equiv_J_pre_vec_field}
 Let $X$ be a complex manifold with complex structure $J\in\Gamma \text{\normalfont End}(TX)$. Then, the tangent bundles\footnote{For $TX$, $X$ is viewed as a real manifold.} $TX$ and $T^{(1,0)}X$ are isomorphic as smooth complex vector bundles via:
 \begin{gather*}
  f:TX\to T^{(1,0)}X,\quad v_R\mapsto \frac{1}{2}(v_R - i\cdot J(v_R)),
 \end{gather*}
 where the complex vector space structure of the fibers of $TX$ is given by $J$.\\
 Now consider the space $\Gamma_J (TX)\coloneqq \{V_R\in\Gamma (TX)\mid L_{V_R}J = 0\}$ of smooth real $J$-preserving\footnote{$J$-preserving vector fields are called infinitesimal automorphisms in \cite{kobayashi1969}.} vector fields on $X$. Here, $\Gamma (TX)$ denotes the space of smooth real vector fields on $X$ and $L_{V_R}J$ is the Lie derivative of $J$ with respect to $V_R$. Then, $\Gamma_J (TX)$ together with the standard commutator $[\cdot, \cdot]$ of vector fields and the complex structure $J$ forms a complex Lie algebra. In fact, $(\Gamma_J (TX), [\cdot, \cdot])$ is isomorphic as complex Lie algebras to the space $(\Gamma (T^{(1,0)}X), [\cdot,\cdot])$ of holomorphic vector fields on $X$ via:
 \begin{gather*}
  F:\Gamma_J (TX)\to \Gamma (T^{(1,0)}X),\quad V_R\mapsto \frac{1}{2}(V_R - i\cdot J(V_R)).
 \end{gather*}
\end{proposition}

\begin{proof}
 Confer Proposition 2.10 and 2.11 in Chapter IX of \cite{kobayashi1969}. 
\end{proof}

One important consequence of Proposition \autoref{prop:holo_vec_field_equiv_J_pre_vec_field} which we heavily use later on is the fact that the real and imaginary part of a holomorphic vector field commute:

\begin{corollary}[$V_R$ and $J(V_R)$ commute]\label{cor:commute}
 Let $X$ be a complex manifold with complex structure $J$ and $V_R\in\Gamma_J (TX)$ be a $J$-preserving vector field. Then, $\left[V_R, J(V_R)\right] = 0$. In particular, the real and imaginary part of holomorphic vector fields on $X$ commute.
\end{corollary}

\begin{proof}
 By Proposition \autoref{prop:holo_vec_field_equiv_J_pre_vec_field}, we know that $(\Gamma_J (TX), [\cdot, \cdot])$ is a complex Lie algebra, hence $\left[V_R, J(V_R)\right] = J\left([V_R, V_R]\right) = 0$.
\end{proof}

By Proposition \autoref{prop:holo_vec_field_equiv_J_pre_vec_field}, we can associate with the Hamiltonian vector field $X_\mH$ of a HHS $(X,\Omega,\mH)$ a $J$-preserving vector field $X^R_\mH$ which is uniquely determined by:
\begin{gather*}
 X_\mH = \frac{1}{2}(X^R_\mH - i\cdot J(X^R_\mH)).
\end{gather*}
Equipped with this knowledge, there are now two ways to define holomorphic trajectories of a HHS. The first one is to simply say that holomorphic trajectories are holomorphic integral curves of the holomorphic Hamiltonian vector field $X_\mH$. The second one is to define the holomorphic trajectories as analytic continuations of the integral curves of $X^R_\mH$. Both definitions are indeed equivalent and we make use of both of them. For our purposes, we use the first one as the actual definition and the second one to construct and investigate holomorphic trajectories afterwards:

\begin{definition}[Holomorphic trajectories]\label{def:holo_traj}
 Let $(X,\Omega,\mH)$ be a HHS and\linebreak $X_\mH = 1/2(X^R_\mH - i\cdot J(X^R_\mH))$ be its Hamiltonian vector field. We call a holomorphic map $\gamma:U\to X$ a \textbf{holomorphic trajectory} of the HHS $(X,\Omega, \mH)$ if $\gamma$ satisfies the holomorphic integral curve equation:
 \begin{gather*}
  \gamma^\prime (z) = X_\mH (\gamma (z))\quad\forall z\in U,
 \end{gather*}
 where $U\subset\mathbb{C}$ is an open and connected subset and $\gamma^\prime$ is the complex derivative\footnote{The complex derivative of a curve $\gamma (t+is)$ is defined via holomorphic charts, but can also be computed by $2\gamma^\prime = \pa{t}\gamma - i\pa{s}\gamma$ (cf. proof of Proposition \autoref{prop:holo_traj}).} of $\gamma$. We call a holomorphic trajectory $\gamma:U\to X$ \textbf{maximal} if for every holomorphic trajectory $\hat\gamma:\hat U\to X$ with $U\subset \hat U$ and $\hat\gamma\vert_U = \gamma$ one has $\hat U = U$. We call the integral curves of the real vector field $X^R_\mH$ the \textbf{real trajectories} of the HHS $(X,\Omega,\mH)$.
\end{definition}

Next, let us consider the existence and uniqueness of holomorphic trajectories:

\begin{proposition}[Existence and uniqueness of holomorphic trajectories]\label{prop:holo_traj}
 Let $(X,\Omega,\mH)$ be a HHS. Then, for any $z_0\in\mathbb{C}$ and $x_0\in X$, there exists an open and connected subset $U\subset\mathbb{C}$ with $z_0\in U$ and a holomorphic trajectory\ $\gamma^{z_0, x_0}:U\to X$ of $(X,\Omega,\mH)$ with $\gamma^{z_0, x_0} (z_0) = x_0$. Two holomorphic trajectories $\gamma^{z_0, x_0}_1:U_1\to X$ and $\gamma^{z_0, x_0}_2:U_2\to X$ with $\gamma^{z_0, x_0}_1 (z_0) = x_0 = \gamma^{z_0, x_0}_2 (z_0)$ locally coincide, in particular, they are equal iff their domains $U_1$ and $U_2$ are equal. Furthermore, the holomorphic trajectory $\gamma^{z_0, x_0}$ depends holomorphically on $z_0$ and $x_0$.
\end{proposition}

\begin{proof}
 Let $(X,\Omega,\mH)$ be a HHS with $X_\mH = 1/2(X^R_\mH - iJ(X^R_\mH))$ and let $z_0\in\mathbb{C}$ and $x_0\in X$ be any points. To construct a holomorphic trajectory $\gamma^{z_0, x_0}$, we first realize that $t\mapsto\gamma^{z_0, x_0}(t + is)$ for fixed $s\in\mathbb{R}$ is a real trajectory. We can see this by taking the real part of the holomorphic integral curve equation. Thus, finding holomorphic trajectories amounts to finding analytic continuations of real trajectories. To accomplish this task, we observe that similarly\linebreak $s\mapsto\gamma^{z_0, x_0}(t + is)$ for fixed $t\in\mathbb{R}$ is an integral curve of $J(X^R_\mH)$. Naively, one hopes that $\gamma^{z_0, x_0}(t + is)$ is given by:
 \begin{gather*}
   \gamma^{z_0, x_0}(t + is) = \varphi^{J(X^R_\mH)}_{s-s_0}\circ\varphi^{X^R_\mH}_{t-t_0} (x_0),
 \end{gather*}
 where $z_0 = t_0 + is_0$ and $\varphi^{J(X^R_\mH)}_{s-s_0}$ and $\varphi^{X^R_\mH}_{t-t_0}$ are the flows of the vector fields $J(X^R_\mH)$ and $X^R_\mH$ with times $s-s_0$ and $t-t_0$, respectively. In general, however, this expression is problematic: Even though it is an integral curve of $J(X^R_\mH)$ for fixed $t$, it might not be an integral curve of $X^R_\mH$ for fixed $s$ anymore due to the composition with $\varphi^{J(X^R_\mH)}_{s-s_0}$. In order to avoid this problem, we need the composition of the flows to commute, at least for small times $t-t_0$ and $s-s_0$. This occurs if the vector fields $X^R_\mH$ and $J(X^R_\mH)$ themselves commute. In our situation, this is indeed the case (cf. Corollary \autoref{cor:commute}). Thus, we can define:
 \begin{align*}
   \gamma^{z_0, x_0}(t + is)&\coloneqq \varphi^{J(X^R_\mH)}_{s-s_0}\circ\varphi^{X^R_\mH}_{t-t_0} (x_0)\equiv \varphi^{X^R_\mH}_{t-t_0}\circ\varphi^{J(X^R_\mH)}_{s-s_0} (x_0)\\
   &\equiv\varphi^{(t-t_0)X^R_\mH + (s-s_0)J(X^R_\mH)}_1 (x_0).
 \end{align*}
 The expressions above are well-defined for $|t-t_0|,|s-s_0|<\varepsilon$ with $\varepsilon>0$ small enough and all identical due to the commutativity of $X^R_\mH$ and $J(X^R_\mH)$.\\
 Let us check that the given expressions for $\gamma^{z_0, x_0}$ indeed define a holomorphic trajectory. By construction, the map $\gamma^{z_0, x_0}$ is holomorphic, as it satisfies the Cauchy-Riemann equations. Hence, we only need to compute the complex derivative $\gamma^{z_0, x_0\ \prime}$. If $\phi = (z_1,\ldots, z_{2n}):V\to W\subset\mathbb{C}^{2n}$ is a holomorphic chart of $X$ near $x_0$, then we can define the complex derivative $\gamma^{z_0, x_0\ \prime}(z)$ for suitable $z$ using $(\gamma^{z_0, x_0}_1(z),\ldots, \gamma^{z_0, x_0}_{2n}(z))\coloneqq \phi\circ\gamma^{z_0, x_0} (z)$:
 \begin{gather*}
 \gamma^{z_0, x_0\ \prime} (z)\coloneqq\sum^{2n}_{j = 1} \gamma^{z_0, x_0\ \prime}_j (z)\cdot\left.\pa{z_j}\right\vert_{\gamma^{z_0, x_0} (z)},
 \end{gather*}
 where $\gamma^{z_0, x_0\ \prime}_j (z)$ is the usual complex derivative of a holomorphic map from $\mathbb{C}$ to $\mathbb{C}$. A straightforward calculation reveals that the complex derivative $\gamma^{z_0, x_0\ \prime}(z)$ equates to:
 \begin{gather*}
  \gamma^{z_0, x_0\ \prime}(z) = \frac{1}{2}\left(\frac{\partial \gamma^{z_0, x_0}}{\partial t}(z) - i\cdot\frac{\partial \gamma^{z_0, x_0}}{\partial s}(z)\right).
 \end{gather*}
 By definition of $\gamma^{z_0, x_0}$, we have:
 \begin{gather*}
  \frac{\partial \gamma^{z_0, x_0}}{\partial t}(z) = X^R_\mH (\gamma^{z_0, x_0}(z)),\quad \frac{\partial \gamma^{z_0, x_0}}{\partial s}(z) = J\left(X^R_\mH (\gamma^{z_0, x_0}(z))\right).
 \end{gather*}
 Putting everything together gives:
 \begin{gather*}
  \gamma^{z_0, x_0\ \prime}(z) = \frac{1}{2}\left(X^R_\mH (\gamma^{z_0, x_0}(z)) - i\cdot J\left(X^R_\mH (\gamma^{z_0, x_0}(z))\right)\right) = X_\mH (\gamma^{z_0, x_0}(z)).
 \end{gather*}
 Thus, $\gamma^{z_0, x_0}$ is indeed a holomorphic trajectory. Clearly, $\gamma^{z_0, x_0}$ satisfies\linebreak $\gamma^{z_0, x_0}(z_0) = x_0$ proving the existence in Proposition \autoref{prop:holo_traj}.\\
 To show local uniqueness given an initial value, we recall that $\gamma^{z_0, x_0}$ is just an integral curve of $X^R_\mH$ along the $t$-axis satisfying $\gamma^{z_0, x_0}(z_0) = x_0$. Hence, every other holomorphic trajectory $\hat\gamma^{z_0, x_0}$ with $\hat\gamma^{z_0, x_0}(z_0) = x_0$ agrees with $\gamma^{z_0, x_0}$ for $s = s_0$ and $t$ near $t_0$. This allows us to apply the identity theorem for holomorphic functions to the coordinates of $\gamma^{z_0, x_0}$ and $\hat\gamma^{z_0, x_0}$ in a holomorphic chart near $x_0$ giving us the local uniqueness. In order to show that $\gamma^{z_0, x_0}$ and $\hat\gamma^{z_0, x_0}$ coincide completely iff their domains are equal, we cover the images of $\gamma^{z_0, x_0}$ and $\hat\gamma^{z_0, x_0}$ with holomorphic charts and repeatedly apply the identity theorem.\\
 Lastly, we need to show that $\gamma^{z_0, x_0}$ depends analytically on $z_0$ and $x_0$. For $z_0$, this is trivial, since $\gamma^{z_1, x_0}(z)$ and $\gamma^{z_2, x_0}(z)$ for $z_1\neq z_2$ only differ by a translation in $z$. For $x_0$, this is true if and only if the flows $\varphi^{X^R_\mH}_{t-t_0}$ and $\varphi^{J(X^R_\mH)}_{s-s_0}$ of $X^R_\mH$ and $J(X^R_\mH)$ are holomorphic maps from and to $X$. As explained in Chapter IX of \cite{kobayashi1969}, the $J$-preserving vector fields on $X$ are exactly those real vector fields on $X$ whose flow is holomorphic. Remembering that, by Proposition \autoref{prop:holo_vec_field_equiv_J_pre_vec_field}, the vector fields $X^R_\mH$ and $J(X^R_\mH)$ are $J$-preserving concludes the proof.
\end{proof}

\begin{remark}\label{rem:holo_traj_alpha}
 In the last proof, we have used that a holomorphic trajectory $\gamma (t+is)$ of a HHS $(X,\Omega,\mH)$ is an integral curve of $X^R_\mH$ for fixed $s$ and an integral curve of $J(X^R_\mH)$ for fixed $t$. We can generalize this observation. If we express $t+is$ in polar coordinates, $t+is = re^{i\alpha}$, then $\gamma (re^{i\alpha})$ is an integral curve of $\cos(\alpha) X^R_\mH + \sin (\alpha)J(X^R_\mH)$ for fixed $\alpha$.
\end{remark}

The properties we have found so far seem to indicate that holomorphic trajectories of a HHS exhibit the same behavior as trajectories of a RHS. However, this is not entirely true. In sharp contrast to the real case, the maximal holomorphic trajectories, given an initial value, do \underline{not} need to be unique, as the following counterexample demonstrates.

\begin{example}[Central problem in one complex dimension]\label{ex:holo_cen_prob} Let\linebreak $X\coloneqq T^{\ast, (1,0)}\mathbb{C}^\times\cong\mathbb{C}^\times\times\mathbb{C}$ be the holomorphic cotangent bundle of $\mathbb{C}^\times\coloneqq\mathbb{C}\backslash\{0\}$ together with the standard form $\Omega = \Omega_\text{can} = dP\wedge dQ$, $(Q,P)\in X$, and the natural Hamiltonian\footnote{The given Hamiltonian is even regular, i.e., $d\mH\neq 0$ for all points of $X$.} $\mH (Q,P)\coloneqq \frac{P^2}{2} - \frac{1}{8Q^2}$. Physically speaking, the HHS $(X,\Omega, \mH)$ is the complexification of the RHS describing a single particle in one-dimensional position space subject to the almost Kepler-like central potential $V(q) = -\frac{1}{8q^2}$. The Hamiltonian vector field $X_\mH$ of the HHS $(X,\Omega,\mH)$ is:
 \begin{gather*}
  X_\mH (Q,P) = P\cdot\pa{Q} - \frac{1}{4Q^3}\cdot\pa{P}.
 \end{gather*}
 Hence, the holomorphic trajectories $\gamma (z) = (Q(z), P(z))$ of $(X,\Omega, \mH)$ satisfy
 \begin{gather*}
  Q^\prime (z) = P (z),\quad P^\prime (z) = -\frac{1}{4Q^3 (z)}
 \end{gather*}
 or, combining both equations:
 \begin{gather*}
  Q^{\prime\prime} (z) = -\frac{1}{4Q^3 (z)}.
 \end{gather*}
 We want to determine the holomorphic trajectories $\gamma$ with $\gamma (z_0) = x_0 = (Q_0,P_0)$ for $z_0, P_0\in\mathbb{C}$ and $Q_0\in\mathbb{C}^\times$. After translation in $z$, we can assume $z_0 = 0$. A straightforward computation reveals that, locally, the desired solutions are given by:
 \begin{gather*}
  Q(z) = \sqrt{Q^2_0 + 2Q_0 P_0\cdot z + 2E_0\cdot z^2},\quad P(z) = Q^\prime (z),
 \end{gather*}
 where $E_0\coloneqq \mH (Q_0, P_0)$ and $\sqrt{\cdot}$ is chosen such that $\sqrt{Q^2_0} = Q_0$. Two square roots mapping $Q^2_0$ to $Q_0$ coincide on a small neighborhood of $Q^2_0$, however, they do not need to have the same domain. Let us make this precise by choosing values for $Q_0$ and $P_0$. Pick $Q_0 = 1$ and $P_0 = \frac{1}{2}$. Then, $E_0 = 0$ and $Q(z) = \sqrt{z + 1}$. Here, all square roots are admissible that coincide with the standard square root for real positive numbers. For instance, one can choose
 \begin{alignat*}{3}
  &\sqrt{\cdot}^1:\{z\mid \text{Im}(z) \neq 0\text{ or }\text{Re}(z) > 0\}\to\mathbb{C},\ &&z = re^{i\alpha}\mapsto \sqrt{r}e^{\frac{i\alpha}{2}},\ &&\alpha \in (-\pi, \pi)\\
  &\text{or} && &&\\
  &\sqrt{\cdot}^2:\{z\mid \text{Re}(z) \neq 0\text{ or }\text{Im}(z) > 0\}\to\mathbb{C},\ &&z = re^{i\alpha}\mapsto \sqrt{r}e^{\frac{i\alpha}{2}},\ &&\alpha \in \left(-\frac{\pi}{2}, \frac{3\pi}{2}\right).
 \end{alignat*}
 Using these square roots, the holomorphic trajectories $\gamma_1: U_1\to X$ and\linebreak $\gamma_2: U_2\to X$ are given by:
 \begin{align*}
  Q_1:U_1\coloneqq\{z\in\mathbb{C}\mid \text{Im}(z+1) \neq 0\text{ or }\text{Re}(z+1) > 0\}\to\mathbb{C},\quad Q_1(z)&\coloneqq \sqrt{z + 1}^1,\\
  Q_2:U_2\coloneqq\{z\in\mathbb{C}\mid \text{Re}(z+1) \neq 0\text{ or }\text{Im}(z+1) > 0\}\to\mathbb{C},\quad Q_2(z)&\coloneqq \sqrt{z + 1}^2.
 \end{align*}
 Clearly, $\gamma_1$ and $\gamma_2$ are maximal holomorphic trajectories satisfying\linebreak $\gamma_1 (0) = (1,\frac{1}{2}) = \gamma_2 (0)$. However, their domains $U_1$ and $U_2$ differ showing that maximal trajectories are not unique, even given an initial value.\linebreak In particular, the trajectories $\gamma_1$ and $\gamma_2$ yield different values for\linebreak $z\in\{z\in\mathbb{C}\mid\text{Re}(z+1)<0\text{ and }\text{Im}(z+1)<0\}\subset U_1\cap U_2$, namely\linebreak $\gamma_1 (z) = -\gamma_2 (z)$.
\end{example}

Even though the square root, which spoils the uniqueness of maximal trajectories in the previous example, is not well-defined on all of $\mathbb{C}\backslash\{0\}$, it is well-defined on the $2:1$ covering $z\mapsto z^2$ of $\mathbb{C}\backslash\{0\}$. In the same vein, the maximal trajectories of a HHS become unique after ``passing them down to a covering''. Precisely speaking, we have to promote the maximal trajectories to leaves to make them unique. To understand this idea, we first recall the definition of a foliation:

\pagebreak

\begin{definition}[Foliation]\label{def:holo_foli}
 A $d$-dimensional \textbf{foliation} $\{L_{x_0}\}_{x_0\in I}$ ($I$: index set) of a manifold $M$ is a decomposition of $M$ into path-connected subsets $L_{x_0}\subset M$ called \textbf{leaves}, i.e. $M = \dot\bigcup_{x_0\in I} L_{x_0}$, such that for every point $p\in M$ there exists a chart $\phi = (x_1,\ldots, x_n):U\to V\subset\R^n$ of $M$ near $p$ fulfilling: For every leaf $L_{x_0}$ with $U\cap L_{x_0}\neq\emptyset$, the connected components of $U\cap L_{x_0}$ are given by $x_{d+1} = c_{d+1}$,\ldots, $x_n = c_n$ for some constants $c_j\in\R$. The foliation $\{L_{x_0}\}_{x_0\in I}$ is called \textbf{holomorphic} if $M$ is a complex manifold and the charts $\phi$ can be chosen to be holomorphic.
\end{definition}

Given a RHS $(M,\omega, H)$, we remember that if $E$ is a regular value of $H$, the energy hypersurface $H^{-1}(E)$ is foliated by maximal trajectories of $(M,\omega, H)$. Similarly, there are holomorphic foliations of regular energy hypersurfaces for HHSs. However, the leaves of the holomorphic foliation are ``more than just'' the maximal trajectories this time:

\begin{proposition}[Holomorphic foliation of a regular hypersurface]\label{prop:holo_foli}
 Let $(X,\Omega,\mH)$ be a HHS with complex structure $J$, Hamiltonian vector field\linebreak $X_\mH = 1/2 (X^R_\mH - i J(X^R_\mH))$, and regular value $E$ of $\mH$. Then, the energy hypersurface $\mH^{-1}(E)$ admits a holomorphic foliation. The leaf $L_{x_0}$ of this foliation through a point $x_0\in\mH^{-1}(E)$ is given by:
 \begin{align*}
  L_{x_0}\coloneqq \{y\in X\mid &y = \varphi^{X^R_\mH}_{t_1}\circ\varphi^{J(X^R_\mH)}_{s_1}\circ\varphi^{X^R_\mH}_{t_2}\circ\varphi^{J(X^R_\mH)}_{s_2}\circ\ldots\circ\varphi^{X^R_\mH}_{t_n}\circ\varphi^{J(X^R_\mH)}_{s_n} (x_0);\\
  &t_1,\ldots,t_n, s_1,\ldots, s_n\in\mathbb{R};\ n\in\mathbb{N}\},
 \end{align*}
 where $\varphi^{X^R_\mH}_{t_j}$ and $\varphi^{J(X^R_\mH)}_{s_j}$ are the flows of $X^R_\mH$ and $J(X^R_\mH)$ for time $t_j$ and $s_j$, respectively. Every holomorphic trajectory of $(X,\Omega,\mH)$ with energy $E$ is completely contained in one such leaf.
\end{proposition}

\begin{proof}
 Take the assumptions and notations from above. $E$ is a regular value of $\mH$, hence, $\mH^{-1}(E)$ is a complex submanifold of $X$. The tangent space of $\mH^{-1} (E)$ consists of vectors $W$ in the tangent space of $X$ satisfying $d\mH (W) = 0$. Using the holomorphicity of $\mH$, $d\mH\circ J = i\cdot d\mH$, we obtain
 \begin{align*}
  d\mH (X^R_\mH) &= d\mH (X_\mH) = -\Omega (X_\mH, X_\mH) = 0,\\
  d\mH (J(X^R_\mH)) &= i\cdot d\mH (X^R_\mH) = 0
 \end{align*}
 showing that $X^R_\mH$ and $J(X^R_\mH)$ live in the tangent space of $\mH^{-1}(E)$. This allows us to restrict $X^R_\mH$ and $J(X^R_\mH)$ to real vector fields on $\mH^{-1}(E)$.\\
 As $\mH$ is regular on $\mH^{-1} (E)$, neither $X^R_\mH$ nor $J(X^R_\mH)$ vanish on $\mH^{-1} (E)$. Furthermore, as real vector fields, they are $\mathbb{R}$-linearly independent at every point of $\mH^{-1}(E)$, since there is no real number which squares to $-1$. By Corollary \autoref{cor:commute}, $X^R_\mH$ and $J(X^R_\mH)$ also commute. This allows us to apply Frobenius' theorem to the distribution spanned by the real vector real fields $X^R_\mH$ and $J(X^R_\mH)$ giving us a foliation of $\mH^{-1}(E)$ whose leaves take the form described in Proposition \autoref{prop:holo_foli}.
 As $X^R_\mH$ and $J(X^R_\mH)$, the vector fields generating the foliation, are real and imaginary part of a holomorphic vector field, the foliation itself is holomorphic due to the holomorphic version of Frobenius' theorem (cf. Theorem 2.26 in \cite{voisin2002}). Comparing the construction of holomorphic trajectories in the proof of Proposition \autoref{prop:holo_traj} with the form of the leaves in Proposition \autoref{prop:holo_foli} reveals that a holomorphic trajectory is completely contained in one leaf concluding the proof.
\end{proof}

\begin{remark}[Flows of $X^R_\mH$ and $J(X^R_\mH)$ do not commute globally]\label{rem:flows_do_not_commute}
 One might be tempted to set $n$ in the definition of the leaves in Proposition \autoref{prop:holo_foli} to $1$, since $X^R_\mH$ and $J(X^R_\mH)$ as well as their flows commute. However, this is only locally the case. To illustrate this, consider Example \autoref{ex:holo_cen_prob} again. Choose the initial value $x_0 = (Q_0, P_0) = (1, 1/2)$ and set $n=2$,  $t_1 = 0$, $s_1 = -2$, $t_2 = -2$, and $s_2 = 1$. Using the solution $Q(z) = \sqrt{z+1}$, we find:
 \begin{gather*}
  \varphi^{X^R_\mH}_{0}\circ\varphi^{J(X^R_\mH)}_{-2}\circ\varphi^{X^R_\mH}_{-2}\circ\varphi^{J(X^R_\mH)}_{1} \left(1, \frac{1}{2}\right) = \left(\sqrt[4]{2}\cdot e^{i\frac{5\pi}{8}}, \frac{1}{2\sqrt[4]{2}}\cdot e^{-i\frac{5\pi}{8}}\right).
 \end{gather*}
 If we exchange the order of $s_1$ and $s_2$, the result differs by a sign:
 \begin{gather*}
  \varphi^{X^R_\mH}_{0}\circ\varphi^{J(X^R_\mH)}_{1}\circ\varphi^{X^R_\mH}_{-2}\circ\varphi^{J(X^R_\mH)}_{-2} \left(1, \frac{1}{2}\right) = -\left(\sqrt[4]{2}\cdot e^{i\frac{5\pi}{8}}, \frac{1}{2\sqrt[4]{2}}\cdot e^{-i\frac{5\pi}{8}}\right).
 \end{gather*}
\end{remark}

In light of Proposition \autoref{prop:holo_foli}, one might say that the leaves of a HHS $(X,\Omega, \mH)$ should be considered to be the holomorphic counterpart to the maximal trajectories of a RHS $(M,\omega, H)$. Let us investigate this statement further. To do that, we first need to generalize the notion of holomorphic trajectories in such a way that we can view any Riemann surface as the domain of a trajectory, not only subsets of $\mathbb{C}$:

\begin{definition}[Geometric trajectory]\label{def:geo_traj}
 Let $(X,\Omega,\mH)$ be a HHS with regular value $E$ of $\mH$ and foliation $L = \{L_{x_0}\}_{x_0\in I}$ ($I$: some index set) of $\mH^{-1}(E)$ as in Proposition \autoref{prop:holo_foli}. Further, let $\Sigma$ be a Riemann surface, i.e., a connected, complex one-dimensional manifold. We call a holomorphic map $\gamma:\Sigma\to X$ a \textbf{geometric trajectory} of energy $E$ if $\gamma$ is an immersion and the image of $\gamma$ is completely contained in one leaf $L_{x_0}$ of the foliation $L$.
\end{definition}

The definition of a geometric trajectory is reasonable, as every geometric trajectory is locally a holomorphic trajectory:

\begin{proposition}[Locally, geometric trajectories are holomorphic trajectories]\label{prop:geo_traj}
 Let $(X,\Omega,\mH)$ be a HHS with regular value $E$ of $\mH$, $\Sigma$ be a Riemann surface, and $\gamma:\Sigma\to X$ be a geometric trajectory of energy $E$. Then, for every $s_0\in\Sigma$, there exists an open neighborhood $V\subset\Sigma$ of $s_0$ and a holomorphic chart $\varphi:V\to U\subset\mathbb{C}$ of $\Sigma$ such that $\gamma\circ\varphi^{-1}:U\to X$ is a holomorphic trajectory of the HHS $(X,\Omega,\mH)$.
\end{proposition}

\begin{proof}
 Invoke the assumptions and notations from above. As $\gamma$ is a holomorphic immersion whose image is completely contained in one leaf and the leaves of $L$ are generated by the Hamiltonian vector field $X_\mH$, there exists for every $s\in\Sigma$ an uniquely determined vector $Y_\mH (s)\in T^{(1,0)}_s\Sigma$ such that:
 \begin{gather*}
  d\gamma_s (Y_\mH (s)) = X_\mH (\gamma (s)).
 \end{gather*}
 These vectors form a holomorphic vector field $Y_\mH$ on $\Sigma$. Now pick $s_0\in\Sigma$ and $z_0\in\mathbb{C}$. By Proposition \autoref{prop:holo_traj} and Remark \autoref{rem:holo_vec_fields}, there exists an open and connected subset $U\subset\mathbb{C}$ such that $\varphi^{-1}:U\to\Sigma$ is a holomorphic integral curve of $Y_\mH$ satisfying $\varphi^{-1}(z_0) = s_0$. After shrinking $U$ if necessary, $\varphi^{-1}$ becomes a biholomorphism onto its image $\varphi^{-1}(U)\eqqcolon V$. Thus, $\varphi:V\to U$ is a holomorphic chart of $\Sigma$ near $s_0$. Furthermore, the curve $\gamma\circ\varphi^{-1}:U\to X$ fulfills:
 \begin{gather*}
  \left(\gamma\circ\varphi^{-1}\right)^\prime (z) = d\gamma_{\varphi^{-1}(z)}\left(\varphi^{-1\,\prime} (z)\right) = d\gamma_{\varphi^{-1}(z)}\left(Y_\mH(\varphi^{-1} (z))\right) = X_\mH (\gamma\circ\varphi^{-1}(z)).
 \end{gather*}
 Hence, $\gamma\circ\varphi^{-1}$ is a holomorphic trajectory concluding the proof.
\end{proof}

Let us now assume that $(X,\Omega,\mH)$ is a HHS with regular value $E$ of $\mH$ and foliation $L = \{L_{x_0}\}_{x_0\in I}$ of $\mH^{-1}(E)$. Pick a leaf $L_{x_0}$. If $L_{x_0}$ is a complex submanifold of $X$, the inclusion $L_{x_0}\hookrightarrow X$ is clearly a geometric trajectory. If $L_{x_0}$ is not a complex submanifold of $X$, we can always equip $L_{x_0}$ with the structure of a complex manifold by choosing a suitable atlas such that the inclusion $L_{x_0}\hookrightarrow X$ becomes a geometric trajectory. The atlas in question consists of maps $\gamma^{-1}$ where $\gamma:U\to X$ is an injective holomorphic trajectory whose image is contained in $L_{x_0}$. In contrast to maximal trajectories, the geometric trajectories $L_{x_0}\hookrightarrow X$ are unique given the initial value $x_0\in \mH^{-1}(E)\subset X$. In fact, the uniqueness can be expressed as a universal property: Pick $x_0\in\mH^{-1}(E)$. Then, every geometric trajectory $\gamma:\Sigma\to X$ with initial value $x_0\in\gamma (\Sigma)$ factors uniquely through the geometric trajectory $L_{x_0}\hookrightarrow X$ and the geometric trajectory $L_{x_0}\hookrightarrow X$ is unique up to biholomorphisms with that property.\\
Often, the geometric trajectories $L_{x_0}\hookrightarrow X$ can be understood as coverings of or to be covered by maximal trajectories. Example \autoref{ex:holo_cen_prob} exemplifies this behavior. To see this, we first need to determine the leaves in Example \autoref{ex:holo_cen_prob}. We achieve this by applying the following proposition:

\begin{proposition}[Energy hypersurfaces of low-dimensional systems]\label{prop:low_dim}
 Let $(X,\Omega,\mH)$ be a HHS with $\text{\normalfont dim}_\mathbb{C}(X) = 2$ and regular value $E$ of $\mH$. Then, the leaves of $\mH^{-1} (E)$ are its connected components and, in particular, complex submanifolds of $X$.
\end{proposition}

\begin{proof}
 Take the notations and assumptions from above. $X$ is a complex two-dimensional manifold, hence, $\mH^{-1}(E)$ is a complex one-dimensional one. Likewise, the leaves of $\mH^{-1}(E)$ are one-dimensional complex manifolds, immersed in $\mH^{-1}(E)$. Thus, the leaves are open in $\mH^{-1}(E)$. By definition of a foliation, $\mH^{-1}(E)$ decomposes into a disjoint union of leaves. Therefore, the leaves are also closed in $\mH^{-1}(E)$ concluding the proof.
\end{proof}

Return to Example \autoref{ex:holo_cen_prob}. By Proposition \autoref{prop:low_dim}, the leaves in this example are just the connected components of the energy hypersurfaces $\mH^{-1}(E)$. For $E\neq 0$, the energy hypersurface is connected and only consists of one leaf, namely itself. For $E = 0$, we have:
\begin{gather*}
 \mH(Q,P) = \frac{P^2}{2} - \frac{1}{8Q^2} = \frac{1}{2}\left(P-\frac{1}{2Q}\right)\left(P+\frac{1}{2Q}\right) \stackrel{!}{=} 0.
\end{gather*}
We see that there are two connected components and, consequentially, two leaves this time: One with $Q\cdot P = \frac{1}{2}$ and another one with $Q\cdot P = -\frac{1}{2}$. We have already determined the maximal trajectories of the leaf with $Q\cdot P = \frac{1}{2}$. They are given by $Q (z) = \sqrt{z+1}$ and $P = Q^\prime (z)$ with appropriate square roots. Precomposing the maximal trajectories with the $2:1$ covering $z\mapsto z^2 -1$ gives us a well-defined map from $\mathbb{C}\backslash\{0\}$ to the leaf with $Q\cdot P = \frac{1}{2}$. In fact, this map is a biholomorphism. In that regard, the leaf can be understood as the double cover of the maximal trajectories. For an example where a leaf is covered by a maximal trajectory, confer Example \autoref{ex:complex_torus}.\\
Before we conclude the subsection on holomorphic trajectories, we quickly add two comments. The first one concerns the notion of holomorphic and geometric trajectories. Throughout the remainder of \autoref{chap:PHHS}, we will not distinguish between holomorphic and geometric trajectories anymore and use both names interchangeably. The second remark concerns holomorphic vector fields on general complex manifolds: 

\begin{remark}\label{rem:holo_vec_fields}
 The results of this subsection are not only true for holomorphic Hamiltonian vector fields and holomorphic trajectories, but for all holomorphic vector fields and holomorphic integral curves.
\end{remark}

\subsection*{Application of HHSs: Lefschetz and Almost Toric Fibrations}

One interesting aspect of HHSs is their interplay with two important structures in symplectic geometry: Lefschetz fibrations and almost toric fibrations. In this subsection, we briefly explore the connection between these structures. We begin the investigation by giving a short introduction to Lefschetz and almost toric fibrations.

\subsubsection*{Lefschetz Fibrations}

Let us first recall the definition of a Lefschetz fibration:

\begin{definition}[Lefschetz fibration]
 Let $X$ be a smooth $2m$-manifold and $C$ be a smooth $2$-manifold (both possibly with boundary). We call a smooth surjective map $\pi:X\to C$ a \textbf{Lefschetz fibration} if the following three conditions are satisfied:
 \begin{enumerate}
  \item $\partial X = \pi^{-1} (\partial C)$
  \item All points on the boundary $\partial X$ are regular points of $\pi$.
  \item For every critical point $p\in\text{Crit}(\pi)\subset\text{Int}(X)$, there exists a smooth chart $\psi_X:U_X\to V_X\subset\R^{2m}\cong \C^m$ near $p$ with $\psi_X (p) = 0$ and a smooth chart\linebreak $\psi_C:U_C\to V_C\subset\R^{2}\cong \C$ near $\pi (p)$ such that:
   \begin{gather*}
    \psi_C\circ\pi\circ\psi_X^{-1} (z_1,\ldots,z_m) = \sum^m_{j=1}z^2_j.
   \end{gather*}
 \end{enumerate}
\end{definition}

Roughly speaking, a Lefschetz fibration generalizes the notion of a fiber bundle over a surface where we now also allow singular fibers. This aspect is captured by the following proposition:

\begin{proposition}[Lefschetz fibrations as fiber bundles]\label{prop:Lef_fiber_bundles}
 Let $\pi:X\to C$ be a Lefschetz fibration and $C^\ast$ be the set of regular values of $\pi$. Further, assume that $X$ is connected. If $\pi:\pi^{-1}(C^\ast)\to C^\ast$ is proper, then $\pi^{-1}(C^\ast)\stackrel{\pi}{\to} C^\ast$ is a fiber bundle. In particular, if $X$ (and then also $C$) is compact, $\pi^{-1}(C^\ast)\stackrel{\pi}{\to} C^\ast$ is a fiber bundle.
\end{proposition}

\begin{proof}
 This follows directly from the fact\footnote{Due to Ehresmann, cf. \cite{Ehresmann1952}.} that every smooth, surjective, and proper submersion between connected manifolds is a fiber bundle.
\end{proof}

Most authors include additional conditions in the definition of a Lefschetz fibration. Often, $X$ is assumed to be a compact, connected, and oriented four-fold. On one hand, this has historic reasons: While Lefschetz introduced these fibrations to study the topology of complex surfaces, Donaldson and Gompf brought Lefschetz fibrations to the attention of the symplectic community by showing that, under mild conditions, every compact symplectic four-fold admits the structure of a Lefschetz fibration (after blowing up if necessary) and vice versa. For a comprehensive overview of the history of four-folds and Lefschetz fibrations, confer the introduction in \cite{Fuller2003}.\\
On the other hand, the case $m=2$ has a rich structure and is well understood: If $X$ is a closed four-fold, then the regular fibers of $\pi:X\to C$ are closed surfaces. Under suitable conditions\footnote{Both $X$ and $C$ are oriented and connected and, additionally, $C$ is simply connected (cf. \cite{Naylor2016}).}, the regular fibers of $\pi:X\to C$ are even oriented and connected meaning that they are surfaces of genus $g$. In this case, the singular fibers are pinched surfaces of genus $g$ (cf. \cite{Naylor2016}). In our investigation, the case $m=2$ also plays an important role.

\subsubsection*{Almost Toric Fibrations}

Almost toric fibrations generalize the notion of toric fibrations which themselves can be understood as the moment map of an effective Hamiltonian torus action (cf. \autoref{app:reduction} for the definition of a moment map). To capture this idea, let us recall the famous convexity theorem of Atiyah, Guillemin, and Sternberg (confer, for instance, \cite{Audin2004} and \cite{Symington2002}):

\begin{theorem}[Convexity theorem]
 Let $(X^{2m},\omega)$ be a connected and closed symplectic manifold with an effective Hamiltonian $T^m$-action on it. Then, the image of the associated moment map $\pi: X\to\R^m$ is a convex polytope.
\end{theorem}

The polytope has a natural stratification. The highest dimensional stratum, i.e., its interior is the set of regular values of $\pi$. The regular level sets are Lagrangian tori and one can express $\omega$ and $\pi$ in a neighborhood of any regular point as:\footnote{This is the equally famous Arnold-Liouville theorem. We use $\omega = \sum_jdx_j\wedge dy_j$ instead of $\omega = \sum_jdy_j\wedge dx_j$ here to match the conventions used in \cite{Symington2002} and \cite{Leung2003}.}
\begin{gather*}
 \omega = \sum^m_{j=1} dx_j\wedge dy_j,\quad \pi_j = x_j.
\end{gather*}
For $k<m$, the points on the $k$-dimensional stratum are critical values. Similarly to the regular case, one can show that in a neighborhood of such a critical point we can write (cf. \cite{Symington2002} and \cite{Leung2003}):
\begin{gather}
 \omega = \sum^m_{j=1} dx_j\wedge dy_j,\ \pi_j = x_j\text{ for }1\leq j\leq k,\ \pi_j = x^2_j + y^2_j\text{ for }k<j\leq m.\label{eq:toric}
\end{gather}
We see that a toric fibration is a symplectic manifold viewed as a fiber bundle whose regular fibers are Lagrangian tori and whose singular fibers take the local form described by \autoref{eq:toric}.\\
Almost toric fibrations, introduced by M. Symington in \cite{Symington2002}, exemplify a similar structure, but the local description of their singular fibers is broadened. The following definitions are taken from \cite{Symington2002} and \cite{Leung2003}:

\begin{definition}[Lagrangian and almost toric fibrations]\label{def:alm_tor}
 Let $C^m$ be a smooth and $(X^{2m},\omega)$ be a symplectic manifold (both possibly with boundary and corners). Further, let $\pi:X\to C$ be a smooth and surjective map. We call $\pi$ a \textbf{Lagrangian fibration} of $(X,\omega)$ if there exists an open and dense subset $C^\ast\subset C$ such that $\pi^{-1}(C^\ast)\stackrel{\pi}{\to} C^\ast$ is a fiber bundle with Lagrangian fibers. We call $\pi$ an \textbf{almost toric fibration} if $\pi$ is a Lagrangian fibration and every critical point $p\in X$ of $\pi$ has a neighborhood in which $\omega$ and $\pi$, after choosing charts for $X$ and $C$, take the following form ($0\leq k<m$, $x(p) = y(p) = 0$):
 \begin{enumerate}
  \item $\omega = \sum^m_{j=1} dx_j\wedge dy_j$,
  \item $\pi_j = x_j$ for $1\leq j\leq k$,
  \item $\pi_j = x^2_j + y^2_j$ (toric) \textbf{or}\\
  $(\pi_j, \pi_{j+1}) = (x_jy_j + x_{j+1}y_{j+1}, x_jy_{j+1} - x_{j+1}y_j)$ (nodal) for $k<j\leq m$.
 \end{enumerate}
\end{definition}

\begin{remark}[Fibers of almost toric fibrations]\label{rem:toric_fibers}
 If the regular fibers of an almost toric fibration are compact and connected, then they are diffeomorphic to a torus by the Arnold-Liouville theorem, explaining the name ``almost toric fibration''.
\end{remark}

\begin{remark}[Difference between toric and nodal points]\label{rem:alm_tor_crit}
 Let $\pi:X\to C$ be an almost toric fibration and let $p\in X$ be a critical point of $\pi$ contained in the interior of $X$. If the local description of $p$ has at least one toric component, then $\pi (p)$ lives in the boundary of $C$. If all components are nodal, it is easy to verify that $\pi (p)$ is an interior point of $C$.
\end{remark}

Let us now investigate the relation between Lefschetz and almost toric fibrations using HHSs. The idea is that we define a compatible holomorphic symplectic structure on a Lefschetz fibration. Such a fibration can be interpreted as a HHS. At the same time, a holomorphic symplectic Lefschetz fibration in real dimension four is, under mild conditions, an almost toric fibration, as we will show later on.\\
First, we formulate the notion of a Lefschetz fibration in a holomorphic setup:

\begin{definition}[Holomorphic Lefschetz fibration]
 Let $X$ be a complex $m$-dimensional and $C$ be a complex one-dimensional manifold (both possibly with boundary\footnote{A complex manifold with boundary is locally biholomorphic to $\{z\in\C^m\mid f(z)\leq 0\}$, where $f\in C^\infty (\C^m,\R)$ is a submersion. We have to allow for curved boundary models, as hypersurfaces in $\C^m$ are usually not biholomorphic to $\{z\in\C^m\mid \text{Im}\, z_1 = 0\}$.}). We say that a surjective holomorphic map $\pi:X\to C$ is a \textbf{holomorphic Lefschetz fibration} if the following three conditions are satisfied:
 \begin{enumerate}
  \item $\partial X = \pi^{-1} (\partial C)$
  \item All points on the boundary $\partial X$ are regular points of $\pi$.
  \item Every critical point $p\in\text{Crit}(\pi)\subset\text{Int}(X)$ is non-degenerate.
 \end{enumerate}
\end{definition}

Here, a critical point $p\in\text{Crit}(\pi)$ is called non-degenerate if its (complex) Hessian\footnote{Since $p$ is a critical point, the non-degeneracy of the Hessian is well-defined, i.e., independent of the choice of charts.} at $p$ is non-degenerate in some holomorphic charts of $X$ and $C$. As the name implies, a holomorphic Lefschetz fibration is also a (smooth) Lefschetz fibration:

\begin{proposition}[Holomorphic Lefschetz fibrations are Lefschetz fibrations]
 Let $\pi:X\to C$ be a holomorphic Lefschetz fibration. Then, $\pi:X\to C$ is a (smooth) Lefschetz fibration in the usual sense. In particular, the charts $\psi_X$ and $\psi_C$ can be chosen to be holomorphic.
\end{proposition}

\begin{proof}
 This is a direct consequence of the holomorphic Morse lemma.
\end{proof}

\begin{lemma}[Holomorphic Morse lemma]
 Let $X$ be a complex manifold, $f:X\to\C$ be a holomorphic function, and $p\in X$ be a non-degenerate critical point of $f$. Then, there exists a holomorphic chart $\psi_X:U_X\to V_X\subset\C^m$ of $X$ near $p$ with $\psi_X (p) = 0$ such that:
 \begin{gather*}
  f\circ\psi^{-1}_X (z_1,\ldots, z_m) = f(p) + \sum^m_{j=1} z^2_j\quad\forall (z_1,\ldots, z_m)\in V_X.
 \end{gather*}
\end{lemma}

\begin{remark}[Signature of Hessian]
 In contrast to the real Morse lemma, we do not need to care about the Morse index, i.e., the signature of the (complex) Hessian in the holomorphic Morse lemma, since all non-degenerate symmetric $\C$-bilinear forms on $\C^m$ are isomorphic.
\end{remark}

\begin{proof}
 The holomorphic Morse lemma can be shown in the same way as the usual Morse lemma (confer, for instance, the proof in \cite{Audin2014} on page 12 ff.), where, of course, we use the appropriate theorems from complex analysis instead of their smooth counterparts\footnote{For instance, the holomorphic implicit function theorem instead of the implicit function theorem.}.
\end{proof}

We now put a compatible symplectic structure on a holomorphic Lefschetz\linebreak fibration:

\begin{definition}[Holomorphic symplectic Lefschetz fibration]\label{def:holo_symp_Lef_fi}
 We call a holomorphic Lefschetz fibration $\pi:X\to C$ \textbf{symplectic} if $X$ admits the structure of a HSM $(X,\Omega)$ such that every critical point $p\in\text{Crit}(\pi)$ has holomorphic Morse-Darboux charts near it, i.e., there are holomorphic charts $\psi_X = (z_1,\ldots, z_{2n}):U_X\to V_X\subset\C^{2n}$ of $X$ near $p$ with $\psi_X (p) = 0$ and $\psi_C:U_C\to V_C\subset\C$ of $C$ near $\pi (p)$ with $\pi (U_X)\subset U_C$ satisfying:
 \begin{enumerate}
  \item $\Omega\vert_{U_X} = \sum^n_{j=1} dz_{j+n}\wedge dz_j$,
  \item $\psi_C\circ\pi\vert_{U_X} = \sum^{2n}_{j=1} z^2_j$.
 \end{enumerate}
\end{definition}

If $\pi:X\to C$ is a holomorphic symplectic Lefschetz fibration with underlying HSM $(X,\Omega)$, we can locally interpret $(X,\Omega,\pi)$ as a HHS after choosing a holomorphic chart $\psi_C$ of $C$. Of course, the Hamiltonian vector field of this HHS is not well-defined, since the Hamilton function depends on the choice of $\psi_C$. However, two Hamiltonian vector fields only differ by a holomorphic function:
\begin{gather*}
 X_{\hat \psi_C\circ\pi} (p) = (\hat\psi_C\circ\psi^{-1}_C)^\prime (\psi_C\circ\pi(p))\cdot X_{\psi_C\circ\pi} (p)\quad\text{for }p\in X.
\end{gather*}
Hence, the two Hamiltonian vector fields have the same trajectories, just parameterized differently. This implies that the holomorphic foliation of the regular level sets of the HHS is still well-defined, even though the Hamiltonian vector field is not. In particular, if $X$ is complex two-dimensional or, equivalently, real four-dimensional and the level sets of $\pi:X\to C$ are connected, then the leaves of this foliation are just the regular level sets themselves allowing us to interpret a holomorphic symplectic Lefschetz fibration as the holomorphic foliation of a HHS.\\
The next step is to link holomorphic symplectic Lefschetz fibrations to almost toric fibrations. We show that every holomorphic symplectic Lefschetz fibration in real dimension four is also an almost toric fibration if $\pi$ is proper.

\begin{proposition}[Holomorphic symplectic Lefschetz fibrations in real dimension four]\label{prop:holo_lef_toric}
 Let $\pi:X\to C$ be a holomorphic symplectic Lefschetz fibration with underlying HSM $(X,\Omega = \Omega_R + i\Omega_I)$. Assume that $X$ is connected and has real dimension four. If $\pi:X\to C$ is proper, then $\pi:X\to C$ is an almost toric fibration of $(X,\Omega_R)$. In particular, if $X$ is compact, then $\pi:X\to C$ is an almost toric fibration of $(X,\Omega_R)$.
\end{proposition}

\begin{proof}
 The set $C^\ast$ of regular values of $\pi$ is open and dense in $C$. By the same argument as in the proof of Proposition \autoref{prop:Lef_fiber_bundles}, $\pi^{-1}(C^\ast)\stackrel{\pi}{\to}C^\ast$ is a fiber bundle. We now consider the fibers of $\pi^{-1}(C^\ast)\stackrel{\pi}{\to}C^\ast$. The HHS $(X,\Omega,\pi)$ is complex two-dimensional, hence, the regular level sets of $\pi$ are complex Lagrangian submanifolds of $(X,\Omega)$. Thus, they are also Lagrangian submanifolds of $(X,\Omega_R)$, as $\Omega = \Omega_R + i\Omega_I$ is a holomorphic $2$-form. This implies that $\pi:X\to C$ is a Lagrangian fibration of $(X,\Omega_R)$.\\
 To conclude the proof, we need to consider a critical point $p$ of $\pi$. By definition, there are holomorphic charts $\psi_X = (z_1, z_2):U_X\to V_X\subset\C^2$ of $X$ near $p$ with $\psi_X (p) = 0$ and $\psi_C:U_C\to V_C\subset\C$ of $C$ near $\pi (p)$ satisfying:
 \begin{align*}
  \Omega\vert_{U_X} &= dz_2\wedge dz_1,\\
  \psi_C\circ\pi\vert_{U_X} &= z^2_1 + z^2_2.
 \end{align*}
 In the new coordinates $z_1\eqqcolon (\hat z_2 - i\hat z_1)/\sqrt{2}$, $z_2\eqqcolon (\hat z_1 - i\hat z_2)/\sqrt{2}$, and\linebreak $\hat \psi_C \coloneqq -\psi_C/2i$, we find:
 \begin{alignat*}{4}
  \Omega\vert_{\hat U_X} &= \, d\hat z_1\wedge d\hat z_2 &&= &&(d\hat x_1\wedge d\hat x_2 - d\hat y_1\wedge d\hat y_2) &&+ i(d\hat x_1\wedge d\hat y_2 + d\hat y_1\wedge d\hat x_2) ,\\
  \hat \psi_C\circ\pi\vert_{\hat U_X} &= \hat z_1\hat z_2 &&= &&(\hat x_1 \hat x_2 - \hat y_1\hat y_2) &&+ i(\hat x_1\hat y_2 + \hat y_1 \hat x_2),
 \end{alignat*}
 where we have used the decomposition $\hat z_j = \hat x_j + i\hat y_j$. In particular, we have\linebreak $\Omega_R\vert_{\hat U_X} = d\hat x_1\wedge d\hat x_2 - d\hat y_1\wedge d\hat y_2$ in these coordinates. Setting $x_1 \coloneqq \hat x_1$, $x_2\coloneqq -\hat y_1$, $y_1\coloneqq \hat x_2$, $y_2\coloneqq \hat y_2$, $\pi_1 \equiv \text{Re}(\hat \psi_C\circ\pi\vert_{\hat U_X})$, and $\pi_2 \equiv \text{Im}(\hat \psi_C\circ\pi\vert_{\hat U_X})$ reproduces the local structure near a nodal critical point as in the definition of an almost toric fibration.
\end{proof}

\begin{remark}[Proposition \autoref{prop:holo_lef_toric} for $(X,\Omega_I)$]
 With the same assumptions as in Proposition \autoref{prop:holo_lef_toric}, we find that $\pi:X\to C$ is also an almost toric fibration of $(X,\Omega_I)$ if $\pi:X\to C$ is proper. In that regard, a holomorphic symplectic Lefschetz fibration gives rise to two different almost toric fibrations.
\end{remark}

\begin{remark}[Critical points are nodal]
 Observe that the critical points of a holomorphic symplectic Lefschetz fibration are nodal (cf. Definition \autoref{def:alm_tor}).
\end{remark}

One might wonder how many manifolds $X$ Proposition \autoref{prop:holo_lef_toric} is applicable to. In the case that $X$ is closed, the answer is already known: Leung and Symington classified in \cite{Leung2003} all closed almost toric four-folds up to diffeomorphisms. They have shown that the only examples which are locally Lefschetz, i.e., whose critical points are non-toric are the K3 surface with base $C = S^2$ and its $\mathbb{Z}_2$-quotient, i.e., the Enriques surface with base $C = \R P^2$. Since $\R P^2$ is not orientable, it cannot admit a complex structure and, thus, the Enriques surface cannot be a holomorphic symplectic Lefschetz fibration. Hence, the only possible example of a closed holomorphic symplectic Lefschetz fibration in four dimensions is the K3 surface\footnote{By Remark \autoref{rem:toric_fibers}, the regular fibers of the K3 surface are tori.}.\\
To conclude this section, we generalize the last question and ask ourselves whether there is an obstruction for a holomorphic Lefschetz fibration to be symplectic. Of course, the space $X$ of a holomorphic Lefschetz fibration $\pi:X\to C$ needs to admit the structure of a HSM in order for $\pi:X\to C$ to be symplectic. This itself is a non-trivial condition. But even if $X$ is a HSM, it is not clear whether the Lefschetz fibration admits Morse-Darboux charts near critical points. This problem is especially interesting, since we already know that every HSM admits Darboux charts near any point and every holomorphic Lefschetz fibration admits Morse charts near critical points.\\
To tackle this question, consider the model case $X=\C^{2n}$ with:
\begin{gather*}
 \Omega = \sum^n_{j=1} dz_{j+n}\wedge dz_j,\quad \pi = \sum^{2n}_{j=1} z^2_j.
\end{gather*}
On $\R^{2n}\subset\C^{2n}$, $\Omega$ and $\pi$ become:
\begin{gather*}
 \omega = \sum^n_{j=1} dx_{j+n}\wedge dx_j,\quad f = \sum^{2n}_{j=1} x^2_j.
\end{gather*}
The holomorphic objects $\Omega$ and $\pi$ are completely determined by their values on the real form $\R^{2n}\subset\C^{2n}$, i.e., by $\omega$ and $f$, and every real-analytic pair $(\omega,f)$ on a real form gives rise, at least locally, to a unique pair $(\Omega,\pi)$ (cf. \autoref{app:real_structures}). Thus, the problem of finding holomorphic Morse-Darboux charts reduces to the problem of finding real-analytic Morse-Darboux charts on a real form. The existence of real-analytic Morse-Darboux charts is intimately linked to the local existence of flat Kähler metrics. Kähler manifolds and their various modifications are explored in \autoref{app:kaehler}. For now, it suffices to know that a Kähler manifold is a symplectic manifold $(M,\omega)$ together with a complex structure $J$ such that $g\coloneqq\omega (\cdot,J\cdot)$ is a Riemannian metric. Locally, every Kähler form $\omega$ possesses a Kähler potential\footnote{This is just a consequence of the $\partial\bar{\partial}$-lemma. To define $\partial$ and $\bar{\partial}$, recall that forms on complex manifolds admit a bigrading which allows us to split up the exterior derivative $d$, i.e., $d = \partial + \bar{\partial}$. The derivatives $\partial$ and $\bar{\partial}$ satisfy $\partial^2 = \bar{\partial}^2 = 0$ and $\partial\bar{\partial} = -\bar{\partial}\partial$. In holomorphic coordinates $(z_1,\ldots,z_n)$, $\partial$ and $\bar{\partial}$ are given by ($f$ is a $\C$-valued function):
\begin{align*}
 \partial (fdz_{r_1}\wedge\ldots\wedge dz_{r_k}\wedge d\bar{z}_{s_1}\wedge\ldots\wedge d\bar{z}_{s_l}) &=\sum^n_{j=1} \pra{f}{z_j} dz_j\wedge dz_{r_1}\wedge\ldots\wedge d\bar{z}_{s_l},\\
 \bar{\partial} (fdz_{r_1}\wedge\ldots\wedge dz_{r_k}\wedge d\bar{z}_{s_1}\wedge\ldots\wedge d\bar{z}_{s_l}) &=\sum^n_{j=1} \pra{f}{\bar{z}_j} d\bar{z}_j\wedge dz_{r_1}\wedge\ldots\wedge d\bar{z}_{s_l}.
\end{align*}
} $f$, i.e., $\omega = i\partial\bar{\partial}f$. However, $f$ is not unique. Indeed, if $g:M\to\C$ is a holomorphic function, then $\omega = i\partial\bar{\partial}\hat f$, where $\hat f\coloneqq f + g + \bar{g}$. Still, we can single out a unique Kähler potential by imposing additional constraints. For this, let $(M,J)$ be a complex manifold, $p\in M$ a point, and $f:M\to\R$ a real-analytic function with $f(p) = 0$. If we pick holomorphic coordinates $\psi = (z_1,\ldots,z_n)$ near $p$ with $\psi (p) = 0$, then we can expand $f$ in a power series of $z_j$ and $\bar{z}_j$. This allows us to write $f = h_1 + h_2 + h_3$, where $h_1$ is a power series in $z_j$, $h_2$ is a power series in $\bar{z}_j$, and $h_3$ collects all terms mixing $z_j$ and $\bar{z}_j$. The decomposition $f = h_1 + h_2 + h_3$ is independent of the choice of chart $\psi$, as long as $\psi (p) = 0$ holds. $h_1$ is holomorphic and, since $f$ is real, we have $h_2 = \bar{h}_1$. Thus, $h_3$ completely determines $\omega\coloneqq i\partial\bar{\partial}f = i\partial\bar{\partial}h_3$. Conversely, $\omega$ fixes $h_3$, as $\omega$ determines the functions $\pa{z_j}\pa{\bar{z}_k}h_3$ and every term in $h_3$ is proportional to $z_j\bar{z}_k$. We call a real-analytic function $f$ for which $h_1$ and $h_2$ vanish \textbf{mixed} near $p$. In that regard, every Kähler manifold $(M,\omega,J)$ admits a unique local Kähler potential that is mixed near $p$. We are now ready to formulate Theorem \autoref{thm:morse_darboux_ex} which shows that the existence of real-analytic Morse-Darboux charts is equivalent to the local existence of flat Kähler metrics with prescribed symplectic form $\omega$ and mixed Kähler potential $\frac{f-f(p)}{2}$:

\begin{theorem}[Existence of real-analytic Morse-Darboux charts]\label{thm:morse_darboux_ex}
 Let\linebreak $(M^{2n},\omega)$ be a real-analytic symplectic manifold, let $f:M\to\R$ be a real-analytic function, and let $p\in M$. Then, the following statements are equivalent:
 \begin{enumerate}[label = (\alph*)]
  \item There is a real-analytic Morse-Darboux chart near $p$, i.e., a real-analytic chart $\psi = (x_1,\ldots, x_{2n}):U\to V\subset\R^{2n}$ of $M$ near $p$ with $\psi (p) = 0$ such that (we drop the restrictions ``$\vert_U$''):
  \begin{gather*}
   \omega = \sum^n_{j=1} dx_{j+n}\wedge dx_j,\quad f = f(p) + \sum^{2n}_{j=1} x^2_j.
  \end{gather*}
  \item There exists a flat Kähler structure near $p$ with symplectic form $\omega$ and mixed Kähler potential $\frac{f-f(p)}{2}$, i.e., there is an open neighborhood $U\subset M$ of $p$ and an almost complex structure $J$ on $U$ such that:
  \begin{enumerate}[label = (\roman*)]
   \item $J$ is integrable,
   \item $g\coloneqq \omega(\cdot,J\cdot)$ is a flat Riemannian metric,
   \item $\frac{f-f(p)}{2}$ is the mixed Kähler potential near $p$.
  \end{enumerate}
 \end{enumerate}
\end{theorem}

\begin{remark}
 If (a) or (b) in Theorem \autoref{thm:morse_darboux_ex} is fulfilled, then $p\in M$ is a non-degenerate critical point of $f$ with Morse index\footnote{If this notion is unfamiliar, confer Definition \autoref{def:morse_index} and Remark \autoref{rem:morse_index}.} $\mu_f (p) = 0$. If we want to describe non-degenerate critical points $p$ with Morse index $\mu_f (p)\neq 0$, we have to modify Theorem \autoref{thm:morse_darboux_ex} slightly: $g$ is now a semi-Riemannian metric and minus signs need to be included in the local description of $f$, i.e.,\linebreak $f = f(p) -\sum^{\mu_f(p)}_{j=1}x^2_j + \sum^{2n}_{k=\mu_f(p)+1} x^2_k$.
\end{remark}

\begin{proof}
 ``$\text{(a)}\Rightarrow\text{(b)}$'': Define $J$ on $U$ by $J\pa{x_j}\coloneqq-\pa{x_{j+n}}$ and $J\pa{x_{j+n}}\coloneqq\pa{x_{j}}$. We verify that $J$ satisfies (i), (ii), and (iii). To see that $J$ is integrable, consider the complex coordinates $z_j\coloneqq x_{j+n} + ix_j$. $J$ turns into $i$ in these coordinates, as they fulfill $dz_j\circ J = idz_j$. Next, we compute $g$. It is easy to check that:
 \begin{gather*}
  g\coloneqq \omega(\cdot,J\cdot) = \sum^{2n}_{j=1} dx^2_j.
 \end{gather*}
 Clearly, $g$ is a flat Riemannian metric. Lastly, we calculate $\frac{i}{2}\partial\bar{\partial}f$. In the holomorphic coordinates $(z_1,\ldots,z_n)$, we find:
 \begin{gather*}
  f = f(p) + \sum^n_{j=1} z_j\bar{z}_j\quad\Rightarrow\ \frac{i}{2}\partial\bar{\partial}f = \frac{i}{2}\sum^n_{j=1}dz_j\wedge d\bar{z}_j = \sum^n_{j=1} dx_{j+n}\wedge dx_j = \omega.
 \end{gather*}
 It immediately follows from the last computation that $\frac{f-f(p)}{2}$ is the mixed Kähler potential near $p$.\\
 ``$\text{(b)}\Rightarrow\text{(a)}$'': Let $J$ and $g$ be as in Statement (b). It is a standard result from Kähler theory (cf. Lemma \autoref{lem:kaehler} or Theorem 4.17 in \cite{Ballmann2006}) that every Kähler manifold admits holomorphic normal coordinates\footnote{There are subtle differences between Riemannian and Kählerian normal coordinates. However, they coincide in the flat case (cf. \cite{Bochner1947} and \cite{Jentsch2017} for details).}. This means that there are holomorphic coordinates $\Psi = (z_1 = x_{n+1} + ix_1,\ldots, z_n = x_{2n} + ix_n)$ near $p$ with $\Psi (p) = 0$ such that $g$ takes the following form at $p$:
 \begin{gather*}
  g_p = \frac{1}{2}\sum^n_{j=1} d_pz_j\otimes d_p\bar{z}_j + d_p\bar{z}_j\otimes d_pz_j = \sum^{2n}_{j=1} d_px^2_j
 \end{gather*}
 After shrinking $U$ if necessary, we can assume that $\Psi$ is defined on all of $U$. Since $g$ is flat, the last equation not only holds for $p$, but everywhere on $U$. Thus, we find for $\omega$:
 \begin{gather*}
  \omega = g(J\cdot,\cdot) = \frac{i}{2}\sum^n_{j=1} dz_j\wedge d\bar{z}_j = \sum^n_{j=1}dx_{j+n}\wedge dx_j.
 \end{gather*}
 We can now read off the mixed Kähler potential $\hat f$ near $p$:
 \begin{gather*}
  \hat f = \frac{1}{2}\sum^n_{j=1} z_j\bar{z}_j = \frac{1}{2}\sum^{2n}_{j=1} x^2_j.
 \end{gather*}
 By assumption, $\hat f$ coincides with $\frac{f-f(p)}{2}$ implying:
 \begin{gather*}
  f = f(p) + \sum^{2n}_{j=1} x^2_j.
 \end{gather*}
 This shows that $\psi\coloneqq (x_1,\ldots,x_{2n})$ is the desired Morse-Darboux chart near $p$ concluding the proof.
\end{proof}

Let us return to the problem of finding holomorphic Morse-Darboux charts. As discussed, it suffices to find a real form admitting real-analytic Morse-Darboux charts. Thus, complexifying the notion of Kähler manifolds allows us to formulate Theorem \autoref{thm:morse_darboux_ex} in the holomorphic setting. The holomorphic version of a Kähler manifold is fittingly called \textbf{holomorphic Kähler manifold}. Holomorphic Kähler manifolds are discussed in \autoref{app:kaehler} as well as in \autoref{sec:holo_semi-kaehler} and \autoref{sec:duality}. Equipped with this knowledge, we obtain the following corollary:

\begin{corollary}[Existence of holomorphic Morse-Darboux charts]\label{cor:morse_darboux_holo_ex}
 Let $(X,\Omega)$ be a HSM, $\pi:X\to\C$ a holomorphic function, and $p\in X$ a point. Then, there is a holomorphic Morse-Darboux chart near $p$ if and only if there exists a flat holomorphic Kähler structure near $p$ with holomorphic symplectic form $\Omega$ and mixed Kähler potential $\frac{\pi - \pi(p)}{2}$ (in the sense of Theorem \autoref{thm:morse_darboux_ex}, i.e., $\frac{\pi - \pi(p)}{2}$ restricts to the mixed Kähler potential near $p$ on the real form).
\end{corollary}

\begin{proof}
 The direction ``$\Rightarrow$'' is rather obvious: In the model case $X = \C^{2n}$,\linebreak $\Omega = \sum^n_{j=1} dz_{j+n}\wedge dz_j$, $\pi = \sum^{2n}_{j=1} z^2_j$, we take $\R^{2n}\subset\C^{2n}$ to be the real form and define $J$ by $J\pa{z_j}\coloneqq -\pa{z_{j+n}}$ and $J\pa{z_{j+n}}\coloneqq \pa{z_j}$. This turns $(X,\Omega,J)$ into a holomorphic Kähler manifold with the desired properties. To prove the converse direction, we consider a flat holomorphic Kähler structure near $p$ with holomorphic symplectic form $\Omega$ and mixed Kähler potential $\frac{\pi - \pi (p)}{2}$. On its real form, it becomes a flat Kähler manifold satisfying the assertions of Theorem \autoref{thm:morse_darboux_ex}. Hence, Theorem \autoref{thm:morse_darboux_ex} gives us a real-analytic Morse-Darboux chart near $p$ on the real form. This chart now yields a holomorphic Morse-Darboux chart near $p$ by complexification (cf. Definition \autoref{def:real_form}, ``(ii)$\Rightarrow$(iii)'') concluding the proof.
\end{proof}

Corollary \autoref{cor:morse_darboux_holo_ex} does not completely answer the question in which case a HSM $(X,\Omega)$ and a holomorphic Lefschetz fibration $\pi:X\to C$ are compatible. Indeed, the projection $\pi$ is not a function with values in $\C$. To turn $\pi$ into a holomorphic function, we have to choose a holomorphic chart of $C$. This chart introduces an additional degree of freedom which makes it easier to find a Morse-Darboux chart. As we will see shortly, the only remaining obstacle to find Morse-Darboux charts, at least in the real two-dimensional case, is the unoriented Morse index:

\begin{definition}[Unoriented Morse index $\mu_f (p)$]\label{def:morse_index}
 Let $M^m$ and $L^1$ be $C^2$-manifolds, $f\in C^2 (M,L)$, and $p\in M$ be a critical point of $f$. We say that $p$ is non-degenerate if the Hessian of $\psi_L\circ f\circ \psi^{-1}_M$ at $\psi_M (p)$ is non-degenerate for some charts $\psi_M$ of $M$ near $p$ and $\psi_L$ of $L$ near $f(p)$. The \textbf{unoriented Morse index} $\mu_f (p)$ is the number $\min\{k, m-k\}$, where $k$ is the usual Morse index of $\psi_L\circ f\circ \psi^{-1}_M$ at $\psi_M (p)$, i.e, the number of negative eigenvalues of its Hessian.
\end{definition}

\begin{remark}\label{rem:morse_index}
 Even though the non-degeneracy of the Hessian is independent of the choice of charts, the number of negative eigenvalues of the Hessian is not. An orientation reversing transformation of the chart $\psi_L$, for instance $\psi_L\mapsto -\psi_L$, changes the number from $k$ to $m-k$, explaining the definition of the unoriented Morse index. Often, we will drop the adjective ``unoriented'' and simply say ``Morse index''. It will be clear from the context whether we mean the usual or the unoriented Morse index.
\end{remark}

Let us now consider Lemma \autoref{lem:morse_darboux_lem_I} and Lemma \autoref{lem:morse_darboux_lem_II}:

\pagebreak

\begin{lemma}[Morse-Darboux lemma I]\label{lem:morse_darboux_lem_I}
 Let $(M^2,\omega)$ be a real-analytic symplectic surface, $L^1$ a real-analytic $1$-manifold, $f:M\to L$ a real-analytic function, and $p\in M$ a non-degenerate critical point of $f$ with Morse index $\mu_f (p) = 0$. Further, let $T>0$ be a positive real number. Then, there exists a real-analytic chart $\psi_L:U_L\to V_L\subset\R$ of $L$ near $f(p)$ such that all non-constant trajectories near $p$ of the RHS $(U_M, \omega\vert_{U_M}, H)$ with $U_M\coloneqq f^{-1}(U_L)$ and $H\coloneqq \psi_L\circ f\vert_{U_M}$ are $T$-periodic.
\end{lemma}

\begin{proof}
 Confer \autoref{app:morse_darboux}.
\end{proof}

\begin{lemma}[Morse-Darboux lemma II]\label{lem:morse_darboux_lem_II}
 Let $(M^2,\omega)$ be a real-analytic symplectic surface and let $H:M\to\R$ be a real-analytic function on $M$ with non-degenerate critical point $p\in M$ of Morse index $\mu_H (p)\neq 1$. Further, let $T>0$ be a positive real number. Then, the following statements are equivalent:
 \begin{enumerate}
  \item There exists a real-analytic Morse-Darboux chart near $p$, i.e., a real-analytic chart $\psi_M = (x,y):U_M\to V_M\subset\R^2$ of $M$ near $p$ such that ($\psi_M (p) = 0$):
  \begin{enumerate}[label = (\alph*)]
   \item $H\vert_{U_M} = H(p) \pm \frac{\pi}{T}(x^2 + y^2)$,
   \item $\omega\vert_{U_M} = dy\wedge dx$.
  \end{enumerate}
  \item There exists an open neighborhood $U_M\subset M$ of $p$ such that all non-constant trajectories of the RHS $(U_M, \omega\vert_{U_M}, H\vert_{U_M})$ are $T$-periodic.
  \item There exists a number $E_0 > 0$ such that $\int_{U(E)} \omega = T\cdot E$ for every number $E\in [0, E_0]$, where $U(E)$ is the connected component containing $p$ of the set $\{q\in M\mid |H(q)-H(p)|\leq E\}$.
 \end{enumerate}
\end{lemma}

\begin{proof}
 Confer \autoref{app:morse_darboux}.
\end{proof}

With Lemma \autoref{lem:morse_darboux_lem_I} and \autoref{lem:morse_darboux_lem_II} in mind, it is clear how to show the existence of Morse-Darboux charts: First, we use Lemma \autoref{lem:morse_darboux_lem_I} to find a chart $\psi_L$ in which all trajectories near $p$ have period $T = \pi$ and the usual Morse index of $\psi_L\circ f$ is $0$. Afterwards, we apply Lemma \autoref{lem:morse_darboux_lem_II} to find a chart $\psi_M$ in which $\omega$ and $H = \psi_L\circ f$ assume their respective standard form.\\
Given a HSM $(X,\Omega)$ and a holomorphic Lefschetz fibration $\pi:X\to C$ with $\dim_\C X = 2$, it seems very likely, judging by the real case, that $(X,\Omega)$ is compatible with $\pi$ in the sense that holomorphic Morse-Darboux charts exist near critical points. To prove this rigorously, one has to construct a real form $M\subset X$ through a critical point $p$ such that both $\Omega$ and $\pi$ turn into real objects on $M$. Since the author is unaware how to show this statement, we stop our discussion here and move on to action functionals instead.

\subsection*{Action Functionals and Principles for HHSs}

As in the real case, we wish to link the holomorphic trajectories of a HHS $(X,\Omega,\mH)$ to critical points of some action functional. To achieve this, let us first study the holomorphic symplectic $2$-form $\Omega$ on $X$. By definition, $\Omega$ is non-degenerate on $T^{(1,0)}X$, but vanishes on $T^{(0,1)}X$. Every complex form $\Omega$ satisfies $\overline{\Omega (V,W)} = \overline{\Omega} (\overline{V},\overline{W})$ ($V$ and $W$ are complex tangent vectors), hence, the complex conjugate $\overline{\Omega}$ of the holomorphic symplectic form $\Omega$ is non-degenerate on $T^{(0,1)}X$, but vanishes on $T^{(1,0)}X$\footnote{Here, we have also used that the complex conjugation maps $T^{(1,0)}X$ to $T^{(0,1)}X$ and $T^{(0,1)}X$ to $T^{(1,0)}X$.}. In total, this implies that the real and imaginary part of $\Omega$,
\begin{gather*}
 \Omega_R = \frac{1}{2}(\Omega + \overline{\Omega})\quad\text{and}\quad\Omega_I = \frac{-i}{2}(\Omega - \overline{\Omega}),
\end{gather*}
are non-degenerate on the entire bundle $T_\mathbb{C}X = T^{(1,0)}X\oplus T^{(0,1)}X$. Clearly, $\Omega_R$ and $\Omega_I$ are real $2$-forms on $X$ and, hence, must be already non-degenerate on the real tangent bundle $TX$. As $\Omega$ is holomorphic and closed, $\Omega_R$ and $\Omega_I$ are smooth and closed. Putting everything together, we find that the real and imaginary part $\Omega_R$ and $\Omega_I$ of $\Omega$, respectively, are symplectic $2$-forms on $X$ viewed as a real manifold.\\
Let us return to the HHS $(X,\Omega, \mH)$. Obviously, the real and imaginary part of the Hamiltonian $\mH = \mH_R + i\mH_I$ are smooth real functions on $X$. Thus, any HHS $(X,\Omega,\mH)$ gives rise to four \underline{RHSs}: $(X,\Omega_R,\mH_R)$, $(X,\Omega_I,\mH_R)$, $(X,\Omega_R,\mH_I)$, and $(X,\Omega_I,\mH_I)$\footnote{As we will see later on, it might be more appropriate to say that $(X,\Omega,\mH)$ gives rise to only two RHSs, since $(X,\Omega_R,\mH_R)$ and $(X,\Omega_I,\mH_I)$ as well as $(X,\Omega_I,\mH_R)$ and $(X,\Omega_R,-\mH_I)$ are subject to the same dynamics.}. Our next task is to determine the Hamiltonian vector fields of the four RHSs. We start with $(X,\Omega_R,\mH_R)$. We write the holomorphic Hamiltonian vector field $X_\mH$ as $X_\mH = 1/2 (X^R_\mH - iJ(X^R_\mH))$ and compute $\iota_{X^R_\mH}\Omega_R$:
\begin{align*}
 \iota_{X^R_\mH}\Omega_R &= \frac{1}{2}\left(\iota_{X^R_\mH}\Omega + \iota_{X^R_\mH}\overline{\Omega}\right) = \frac{1}{2}\left(\iota_{X^R_\mH}\Omega + \iota_{\overline{X^R_\mH}}\overline{\Omega}\right) = \frac{1}{2}\left(\iota_{X^R_\mH}\Omega + \overline{\iota_{X^R_\mH}\Omega}\right)\\
 &= \frac{1}{2}\left(\iota_{X_\mH}\Omega + \overline{\iota_{X_\mH}\Omega}\right) = -\frac{1}{2}\left(d\mH + \overline{d\mH}\right) = -\frac{1}{2}d\left(\mH + \overline{\mH}\right) = - d\mH_R,
\end{align*}
where we used that $X^R_\mH$ is a real vector field on $X$, i.e., $\overline{X^R_\mH} = X^R_\mH$, and that\linebreak $\iota_{X_\mH}\Omega = \iota_{X^R_\mH}\Omega$ due to Equation \eqref{eq:J-anticompatible}. We deduce from the expression above that $X^R_\mH$ is the real Hamiltonian vector field of the RHS $(X,\Omega_R,\mH_R)$. Similarly, one can show the following proposition:

\begin{proposition}\label{prop:four_RHS}
 Let $(X,\Omega = \Omega_R + i\Omega_I,\mH = \mH_R + i\mH_I)$ be a HHS with Hamiltonian vector field $X_\mH = 1/2 (X^R_\mH - iJ(X^R_\mH))$. Then:
 \begin{enumerate}
  \item $(X,\Omega_R, \mH_R)$ is a RHS with Hamiltonian vector field $X^R_\mH$.
  \item $(X,\Omega_R, \mH_I)$ is a RHS with Hamiltonian vector field $-J(X^R_\mH)$.
  \item $(X,\Omega_I, \mH_R)$ is a RHS with Hamiltonian vector field $J(X^R_\mH)$.
  \item $(X,\Omega_I, \mH_I)$ is a RHS with Hamiltonian vector field $X^R_\mH$.
 \end{enumerate}
\end{proposition}

\begin{remark}[Cauchy-Riemann-like relations]\label{rem:cauchy-riemann}
 At first glance, one might be confused why the Hamiltonian vector fields of $(X,\Omega_R,\mH_R)$ and $(X,\Omega_I,\mH_I)$ coincide, while the Hamiltonian vector fields of $(X,\Omega_R,\mH_I)$ and $(X,\Omega_I,\mH_R)$ differ by a sign. However, this observation is just a consequence of the analyticity of the HHS $(X,\Omega,\mH)$ and one might think of it as Cauchy-Riemann-like relations:
 \begin{gather*}
  X^{\Omega_R}_{\mH_R} = X^{\Omega_I}_{\mH_I},\quad J\left(X^{\Omega_R}_{\mH_R}\right) = X^{\Omega_I}_{\mH_R} = - X^{\Omega_R}_{\mH_I},
 \end{gather*}
 where $X^{\Omega_a}_{\mH_b}$ is the Hamiltonian vector field of the RHS $(X,\Omega_a,\mH_b)$.
\end{remark}

The upshot of Proposition \autoref{prop:four_RHS} is that the real trajectories of the HHS $(X,\Omega,\mH)$ are just the trajectories of the RHSs $(X,\Omega_R,\mH_R)$ and $(X,\Omega_I,\mH_I)$. In particular, the real trajectories of $(X,\Omega,\mH)$ are critical points of an action functional, at least if $(X,\Omega = d\Lambda, \mH)$ is exact\footnote{A HHS $(X,\Omega,\mH)$ is exact iff $\Omega$ has a holomorphic primitive $\Lambda$.}. Of course, the same is true for the trajectories of $(X, \Omega_R, \mH_I)$ and $(X, \Omega_I, \mH_R)$:

\begin{proposition}[Action principle for real trajectories]\label{prop:action_fct_of_real_traj}
 Let $(X,\Omega = d\Lambda, \mH)$ be an exact HHS with Hamiltonian vector field $X_\mH = 1/2 (X^R_\mH - iJ(X^R_\mH))$ and decompositions $\Omega = \Omega_R + i\Omega_I$, $\Lambda = \Lambda_R + i\Lambda_I$, and $\mH = \mH_R + i\mH_I$. Let $I_0\subset\mathbb{R}$ be an interval. We set $\mathcal{P}_{I_0}\coloneqq C^\infty (I_0,X)$ and define the action functionals $\mathcal{A}^{\Lambda_R}_{\mH_R}:\mathcal{P}_{I_0}\to\mathbb{R}$, $\mathcal{A}^{\Lambda_R}_{-\mH_I}:\mathcal{P}_{I_0}\to\mathbb{R}$\footnote{Note the different signs in the definition of $\mathcal{A}^{\Lambda_R}_{-\mH_I}$ due to the Cauchy-Riemann-like relations.}, $\mathcal{A}^{\Lambda_I}_{\mH_R}:\mathcal{P}_{I_0}\to\mathbb{R}$, $\mathcal{A}^{\Lambda_I}_{\mH_I}:\mathcal{P}_{I_0}\to\mathbb{R}$, $\mathcal{A}^{\Lambda}_{\mH}:\mathcal{P}_{I_0}\to\mathbb{C}$, and $\mathcal{A}^{\Lambda}_{i\mH}:\mathcal{P}_{I_0}\to\mathbb{C}$ by:
 \begin{align*}
  \mathcal{A}^{\Lambda_R}_{\mH_R}[\gamma]&\coloneqq \int\limits_{I_0}\gamma^\ast\Lambda_R - \int\limits_{I_0}\mH_R\circ\gamma (t)dt,\quad \mathcal{A}^{\Lambda_R}_{-\mH_I}[\gamma]\coloneqq \int\limits_{I_0}\gamma^\ast\Lambda_R + \int\limits_{I_0}\mH_I\circ\gamma (t)dt,\\
  \mathcal{A}^{\Lambda_I}_{\mH_R}[\gamma]&\coloneqq \int\limits_{I_0}\gamma^\ast\Lambda_I - \int\limits_{I_0}\mH_R\circ\gamma (t)dt,\quad \mathcal{A}^{\Lambda_I}_{\mH_I}[\gamma]\coloneqq \int\limits_{I_0}\gamma^\ast\Lambda_I - \int\limits_{I_0}\mH_I\circ\gamma (t)dt,\\
  \mathcal{A}^{\Lambda}_{\mH}[\gamma]&\coloneqq \mathcal{A}^{\Lambda_R}_{\mH_R}[\gamma] + i\mathcal{A}^{\Lambda_I}_{\mH_I}[\gamma] = \int\limits_{I_0}\gamma^\ast\Lambda - \int\limits_{I_0}\mH\circ\gamma (t)dt,\\
  \mathcal{A}^{\Lambda}_{i\mH}[\gamma]&\coloneqq \mathcal{A}^{\Lambda_R}_{-\mH_I}[\gamma] + i\mathcal{A}^{\Lambda_I}_{\mH_R}[\gamma] = \int\limits_{I_0}\gamma^\ast\Lambda - i\int\limits_{I_0}\mH\circ\gamma (t)dt,
 \end{align*}
 where $\gamma\in\mathcal{P}_{I_0}$. Now let $\gamma\in\mathcal{P}_{I_0}$ be a smooth path in $X$. Then, $\gamma$ is a real trajectory of $(X,\Omega,\mH)$ if and only if $\gamma$ is a ``critical point'' of the action functionals $\mathcal{A}^{\Lambda_R}_{\mH_R}$, $\mathcal{A}^{\Lambda_I}_{\mH_I}$, and $\mathcal{A}^{\Lambda}_{\mH}$. Similarly, $\gamma$ is a (real) integral curve of $J(X^R_\mH)$ if and only if $\gamma$ is a ``critical point'' of the action functionals $\mathcal{A}^{\Lambda_R}_{-\mH_I}$, $\mathcal{A}^{\Lambda_I}_{\mH_R}$, and $\mathcal{A}^{\Lambda}_{i\mH}$.
\end{proposition}

\begin{proof}
 Proposition \autoref{prop:action_fct_of_real_traj} is a consequence of Proposition \autoref{prop:four_RHS} and the action principle for RHSs (cf. introduction of \autoref{chap:PHHS}).
\end{proof}

\begin{remark}[Meaning of ``critical point'']\label{rem:critical_point}
 ``Critical point'' does not denote an actual critical point. ``Critical point'' in Proposition \autoref{prop:action_fct_of_real_traj} means that the first derivative of the action functionals vanishes at $\gamma$ \underline{only} for those variations of $\gamma$ \textbf{which keep the endpoints of $\mathbf{\gamma}$ fixed!} To understand this notion, consider an exact RHS $(M,\omega = d\lambda, H)$ with associated action functional\linebreak $\mathcal{A}_H:\mathcal{P}_{I_0}\to\mathbb{R}$:
 \begin{gather*}
  \mathcal{A}_H[\gamma]\equiv \mathcal{A}^\lambda_H[\gamma]\coloneqq \int\limits_{I_0} \gamma^\ast\lambda - \int\limits_{I_0} H\circ\gamma (t)\, dt.
 \end{gather*}
 Let $I_0 = [t_1,t_2]$ and let $\gamma_\varepsilon$ be a variation of $\gamma$ which fixes the endpoints of $\gamma$, i.e., a $1$-parameter family of curves with $\gamma_0 = \gamma$ and $\gamma_\varepsilon (t_j) = \gamma (t_j)$ for $j\in\{1,2\}$. Define the vector field $\hat\gamma$ along $\gamma$ by:
 \begin{gather*}
  \hat\gamma (t)\coloneqq\left.\frac{d}{d\varepsilon}\right\vert_{\varepsilon = 0} \gamma_\varepsilon (t).
 \end{gather*}
 To calculate the derivative of $\gamma^\ast_\varepsilon\lambda$, we interpret $\gamma_\varepsilon$ as the flow of $\hat\gamma$. This allows us to express the derivative of $\gamma^\ast_\varepsilon\lambda$ as the Lie derivative $\gamma^\ast L_{\hat\gamma}\lambda$ which we compute with the help of Cartan's magic formula:
 \begin{align*}
  \left.\frac{d}{d\varepsilon}\right\vert_{\varepsilon = 0}\mathcal{A}_H[\gamma_\varepsilon] &= \left.\frac{d}{d\varepsilon}\right\vert_{\varepsilon = 0}\int\limits_{[t_1,t_2]} \gamma^\ast_\varepsilon\lambda - \left.\frac{d}{d\varepsilon}\right\vert_{\varepsilon = 0}\int\limits^{t_2}_{t_1} H\circ\gamma_\varepsilon (t)dt\\
  &= \int\limits_{[t_1,t_2]}\gamma^\ast L_{\hat\gamma}\lambda - \int\limits^{t_2}_{t_1} dH_{\gamma (t)} (\hat\gamma (t))dt\\
  &= \int\limits_{[t_1,t_2]}\gamma^\ast d\iota_{\hat\gamma}\lambda + \int\limits_{[t_1,t_2]}\gamma^\ast \iota_{\hat\gamma}d\lambda - \int\limits^{t_2}_{t_1} \omega_{\gamma (t)} (\hat\gamma (t), X_H (\gamma (t)))dt\\
  &= \int\limits_{[t_1,t_2]}d(\gamma^\ast\iota_{\hat\gamma}\lambda) + \int\limits_{[t_1,t_2]}\gamma^\ast \iota_{\hat\gamma}\omega - \int\limits^{t_2}_{t_1} \omega_{\gamma (t)} (\hat\gamma (t), X_H (\gamma (t)))dt\\
  &= \lambda_{\gamma(t_2)}(\hat\gamma (t_2)) - \lambda_{\gamma(t_1)}(\hat\gamma (t_1)) + \int\limits^{t_2}_{t_1} \omega_{\gamma (t)} (\hat\gamma (t), \dot\gamma (t)- X_H (\gamma (t)))dt\\
  &= \int\limits^{t_2}_{t_1} \omega_{\gamma (t)} (\hat\gamma (t), \dot\gamma (t)- X_H (\gamma (t)))dt\\
  &=0\quad \text{for all variations }\gamma_\varepsilon\qquad \Leftrightarrow\qquad \dot\gamma = X_H\circ\gamma.
 \end{align*}
 Here, the boundary terms vanish, since the vector field $\hat \gamma$ is zero on the boundary.\\
 Often, one wishes to view trajectories as actual critical points of some action functional, not just with fixed endpoints. There are two main ways to achieve this:
 \begin{enumerate}
  \item One can only consider paths which start and end at points where the primitive of the symplectic form vanishes, usually exact Lagrangian submanifolds of the symplectic manifold.
  \item One can only consider periodic paths such that the boundary terms in the first derivative of the action cancel each other.
 \end{enumerate}
\end{remark}

\begin{remark}[Action functional for ``tilted'' trajectories]\label{rem:tilted_traj}
 The observation that the integral curves of $X^R_\mH$ and $J(X^R_\mH)$ are linked to the ``critical points'' of the action functionals $\mathcal{A}^{\Lambda}_{\mH}$ and $\mathcal{A}^{\Lambda}_{i\mH}$, respectively, can be generalized to ``tilted'' trajectories. For any $\alpha\in\mathbb{R}$, $\gamma$ is an integral curve of $\cos (\alpha)\cdot X^R_\mH + \sin (\alpha)\cdot J(X^R_\mH)$ if and only if $\gamma$ is a ``critical point'' of the action functional $\mathcal{A}^{\Lambda}_{e^{i\alpha}\mH}:\mathcal{P}_{I_0}\to\mathbb{C}$ defined by:
 \begin{gather*}
  \mathcal{A}^{\Lambda}_{e^{i\alpha}\mH}[\gamma]\coloneqq \int\limits_{I_0}\gamma^\ast\Lambda - e^{i\alpha}\int\limits_{I_0}\mH\circ\gamma (t)dt.
 \end{gather*}
 Of course, the same action principle holds true if we only consider the real or imaginary part of $\mathcal{A}^{\Lambda}_{e^{i\alpha}\mH}$. This fact will be of great importance later on and in \autoref{app:various_action_functionals}.
\end{remark}

As the holomorphic trajectories of the HHS $(X,\Omega,\mH)$ are analytic continuations of its real trajectories, one might be content with finding action functionals for the real trajectories. However, it is also possible to construct an action functional for the holomorphic trajectories of $(X,\Omega,\mH)$ by averaging the action functionals of the four underlying RHSs:

\begin{lemma}[Action principle for holomorphic trajectories]\label{lem:holo_action_prin}
 Let\linebreak $(X,\Omega = d\Lambda,\mH)$ be an exact HHS with $\Lambda = \Lambda_R + i\Lambda_I$. Furthermore, let\linebreak $\Rec\coloneqq [t_1,t_2] + i[s_1,s_2]\subset\mathbb{C}$ be a rectangle in the complex plane with real numbers $t_1< t_2$ and $s_1< s_2$. Denote the space of smooth maps from $\Rec$ to $X$ by $\mathcal{P}_{\Rec}$ and define the action functional $\mathcal{A}^{\Rec}_\mH:\mathcal{P}_{\Rec}\to\mathbb{C}$ by:
 \begin{align*}
  \mathcal{A}^{\Rec}_\mH[\gamma]&\coloneqq \int\limits^{t_2}_{t_1}\int\limits^{s_2}_{s_1}\left[\Lambda_{R,\gamma (t+is)}\left(2\frac{\partial\gamma}{\partial z}(t+is)\right) - \mH\circ\gamma (t+is)\right] ds\, dt\ \text{with}\\
  \frac{\partial \gamma}{\partial z}&\coloneqq \frac{1}{2}\left(\frac{\partial \gamma}{\partial t} - i\frac{\partial \gamma}{\partial s}\right)\quad\forall \gamma\in\mathcal{P}_{\Rec}.
 \end{align*}
 Now let $\gamma\in\mathcal{P}_{\Rec}$ be a smooth map from ${\Rec}$ to $X$. Then, $\gamma$ is a holomorphic trajectory of $(X,\Omega,\mH)$ if and only if $\gamma$ is a ``critical point''\footnote{``Critical point'' means that only those variations are allowed which keep $\gamma$ fixed on the boundary $\partial {\Rec}$.} of $\mathcal{A}^{\Rec}_\mH$.
\end{lemma}

\begin{proof}
 Take the notations from above and decompose $\mH = \mH_R + i\mH_I$. Furthermore, let $\gamma\in\mathcal{P}_{\Rec}$ be a smooth map and let $\gamma_s:[t_1,t_2]\to X$ and $\gamma_t:[s_1,s_2]\to X$ be defined by $\gamma_s (t) = \gamma (t+is) = \gamma_t (s)$ for any $s\in [s_1,s_2]$ and $t\in [t_1,t_2]$. Recall the action functionals from Proposition \autoref{prop:action_fct_of_real_traj}. A short calculation reveals that $\mathcal{A}^{\Rec}_\mH [\gamma]$ can be expressed as:
 \begin{gather*}
  \mathcal{A}^{\Rec}_\mH [\gamma] = \int\limits^{s_2}_{s_1} \mathcal{A}^{\Lambda_R}_{\mH_R}[\gamma_s] ds - i\int\limits^{t_2}_{t_1} \mathcal{A}^{\Lambda_R}_{-\mH_I}[\gamma_t] dt.
 \end{gather*}
 $\gamma$ is a ``critical point'' of $\mathcal{A}^{\Rec}_\mH$ iff $\gamma$ is a ``critical point'' of its real and imaginary part. Now consider the real part of $\mathcal{A}^{\Rec}_\mH$. For any $s\in [s_1,s_2]$, $\gamma_s$ is a ``critical point'' of $\mathcal{A}^{\Lambda_R}_{\mH_R}$ iff $\gamma_s$ is an integral curve of $X^R_\mH$, where $X_\mH = 1/2(X^R_\mH - iJ(X^R_\mH))$ is the Hamiltonian vector field of $(X,\Omega,\mH)$. Explicitly writing down the first derivative of the functional $\gamma\mapsto\int \mathcal{A}^{\Lambda_R}_{\mH_R}[\gamma_s] ds$ shows that this property is preserved under averaging: $\gamma$ is a ``critical point'' of $\gamma\mapsto\int \mathcal{A}^{\Lambda_R}_{\mH_R}[\gamma_s] ds$ iff $\gamma_s$ is an integral curve of $X^R_\mH$ for every $s\in [s_1,s_2]$.\\
 Similarly, we find for the imaginary part of $\mathcal{A}^{\Rec}_\mH$ that $\gamma$ is a ``critical point'' of\linebreak $\gamma\mapsto\int \mathcal{A}^{\Lambda_R}_{-\mH_I}[\gamma_t] dt$ iff $\gamma_t$ is an integral curve of $J(X^R_\mH)$ for every $t\in [t_1,t_2]$. Combining our results, we find that $\gamma$ is a ``critical point'' of $\mathcal{A}^{\Rec}_\mH$ iff $\gamma_s$ is an integral curve of $X^R_\mH$ for every $s\in [s_1,s_2]$ and $\gamma_t$ is an integral curve of $J(X^R_\mH)$ for every $t\in [t_1,t_2]$. We conclude the proof by recalling from the first subsection of \autoref{sec:HHS} that holomorphic trajectories of $(X,\Omega,\mH)$ are exactly those smooth maps $\gamma$ which are integral curves of $X^R_\mH$ in $t$-direction and of $J(X^R_\mH)$ in $s$-direction.
\end{proof}

To get a better understanding of Lemma \autoref{lem:holo_action_prin}, several remarks are in order:

\begin{remark}\label{rem:several_remarks}
 {\textcolor{white}{Easter Egg}}
 \begin{enumerate}
  \item Note that the action functional $\mathcal{A}^{\Rec}_\mH$ in Lemma \autoref{lem:holo_action_prin} only uses the real part $\Lambda_R$ and \underline{not} $\Lambda_I$. Of course, a similar action functional including $\Lambda_I$ exists, but we will not use it for reasons that become apparent in \autoref{sec:PHHS}.
  \item As for RHSs, one might wish to express holomorphic trajectories as actual critical points of some functional. Again, there are two main ways to achieve this: One may only consider smooth maps $\gamma$ from the rectangle ${\Rec}$ to $X$ which\ldots
  \begin{enumerate}
   \item \dots map the boundary $\partial {\Rec}$ to points in $X$ where the $1$-form $\Lambda_R$ vanishes, usually exact Lagrangian submanifolds of $X$.
   \item \dots are doubly-periodic, i.e., periodic in both $t$- and $s$-direction.
  \end{enumerate}
  Theoretically, one can even imagine a mix of both methods, where one only considers maps $\gamma$ which are periodic in one direction and map the boundary orthogonal to the remaining direction to exact Lagrangian submanifolds of $X$.
  \item In more geometrical terms, the action $\mathcal{A}^{\Rec}_\mH$ can be expressed as:
  \begin{gather*}
   \mathcal{A}^{\Rec}_\mH[\gamma] =  \iint\limits_{\Rec}\left[\Lambda_{R,\gamma (t+is)}\left(2\frac{\partial \gamma}{\partial z}(t+is)\right) - \mH\circ\gamma (t+is)\right] dt\wedge ds,
  \end{gather*}
  where $dt\wedge ds$ is the standard area form on $\mathbb{C}\cong\mathbb{R}^2$. If $\gamma$ is a ``critical point'' of $\mathcal{A}^{\Rec}_\mH$ or simply a holomorphic curve, we find:
  \begin{gather*}
   \frac{\partial\gamma}{\partial z} (t+is) = \gamma^\prime (t+is) \in T^{(1,0)}X.
  \end{gather*}
  Using $\Lambda_R (V) = i\Lambda_I (V)$ for $V\in T^{(1,0)}X$, we obtain that the action at such $\gamma$ is given by:
  \begin{gather*}
   \mathcal{A}^{\Rec}_\mH[\gamma] =  \iint\limits_{R}\left[\Lambda\left(\gamma^\prime (t+is)\right) - \mH\circ\gamma (t+is)\right] dt\wedge ds,
  \end{gather*}
  where the expression in rectangular brackets is holomorphic in $z = t + is$.
  \item Upon closer inspection of Lemma \autoref{lem:holo_action_prin}, one might wonder whether Lemma \autoref{lem:holo_action_prin} is still true if one restricts the domain $\mathcal{P}_{\Rec}$ of $\mathcal{A}^{\Rec}_\mH$ to the holomorphic curves from ${\Rec}$ to $X$ instead of varying over all smooth curves from ${\Rec}$ to $X$. Clearly, this is not the case, as the values $\gamma$ attains at $\partial {\Rec}$ completely determine one holomorphic curve, so variation over this space is not viable. A different perspective is offered by the action functional $\mathcal{A}^{\Rec}_\mH$ itself. By writing $dt\wedge ds = i/2 \cdot dz\wedge d\bar{z}$ and recalling Point (iii), $\mathcal{A}^{\Rec}_\mH [\gamma]$ can be written as the integral of a form admitting a primitive for holomorphic $\gamma$. By Stokes' theorem, $\mathcal{A}^{\Rec}_\mH[\gamma]$ then only depends on the values of $\gamma$ on the boundary $\partial {\Rec}$. Since these values are kept fixed during variation, the action functional never changes in the variational process and gives us no information. An additional explanation for this behavior is presented in \autoref{app:various_action_functionals}.
 \end{enumerate}
\end{remark}

We can define action functionals like $\mathcal{A}^{\Rec}_\mH$ not only for rectangles, but for all kinds of domains in $\mathbb{C}$. A large selection of them is explored in \autoref{app:various_action_functionals}. Here, let us quickly introduce one generalization of $\mathcal{A}^{\Rec}_\mH$, namely the action functional $\mathcal{A}^{P_\alpha}_\mH$ for parallelograms $P_\alpha$. For the sake of simplicity, we assume that the first vector spanning the parallelogram $P_\alpha$ is parallel to the real axis such that we can write $P_\alpha = [t_1,t_2] + e^{i\alpha}\cdot[r_1,r_2]$ for some angle $\alpha\in\mathbb{R}\backslash\{n\cdot\pi\mid n\in\mathbb{Z}\}$ and some real numbers $t_1 < t_2$ and $r_1 < r_2$. Using the standard area form $dt\wedge ds$, the generalization of $\mathcal{A}^{\Rec}_\mH$ to $\mathcal{A}^{P_\alpha}_\mH$ is straightforward:
\begin{gather*}
 \mathcal{A}^{P_\alpha}_\mH[\gamma]\coloneqq \iint\limits_{P_\alpha}\left[\Lambda_{R,\gamma (t+is)}\left(2\frac{\partial\gamma}{\partial z}(t+is)\right) - \mH\circ\gamma (t+is)\right] dt\wedge ds\quad\forall\gamma\in\mathcal{P}_{P_\alpha},
\end{gather*}
where we used the coordinates $z = t + is$ and defined $\partial\gamma/\partial z$ as in Lemma \autoref{lem:holo_action_prin}. To show that $\gamma$ is a ``critical point'' of $\mathcal{A}^{P_\alpha}_\mH$ iff $\gamma$ is a holomorphic trajectory, we express $\mathcal{A}^{P_\alpha}_\mH$ in the ``tilted'' coordinates $z = t + r\cdot e^{i\alpha}$ ($\alpha$ fixed):
\begin{align*}
 \mathcal{A}^{P_\alpha}_\mH[\gamma]&= \int\limits^{t_2}_{t_1}\int\limits^{r_2}_{r_1}\left[\Lambda_{R,\gamma (t+re^{i\alpha})}\left(2\frac{\partial\gamma}{\partial z}(t+re^{i\alpha})\right) - \mH\circ\gamma (t+re^{i\alpha})\right]\, \sin(\alpha)\,dr\, dt\\
 &= \int\limits^{t_2}_{t_1}\int\limits^{r_2}_{r_1}\left[\Lambda_{R,\gamma (t+re^{i\alpha})}\left(ie^{-i\alpha}\cdot \frac{d\gamma_r}{d t}(t) - i\cdot \frac{d\gamma_t}{dr}(r)\right)\right.\\
 &\qquad\quad \left.- \left(ie^{-i\alpha}\mH_R - i\text{Re}(e^{i\alpha}\mH)\right)\circ\gamma (t+re^{i\alpha})\right]\, dr\, dt\\
 &= ie^{-i\alpha}\int\limits^{r_2}_{r_1}\mathcal{A}^{\Lambda_R}_{\mH_R} [\gamma_r] dr - i\int\limits^{t_2}_{t_1}\text{Re} (\mathcal{A}^\Lambda_{e^{i\alpha}\mH} [\gamma_t]) dt,
\end{align*}
where, for any $t\in [t_1, t_2]$ and $r\in [r_1,r_2]$, the curves $\gamma_r:[t_1,t_2]\to X$ and $\gamma_t:[r_1,r_2]\to X$ are defined by $\gamma_r (t) = \gamma (t + re^{i\alpha}) = \gamma_t (r)$, $\text{Re}(\cdot)$ denotes the real part, and $\mathcal{A}^{\Lambda_R}_{\mH_R}$ and $\mathcal{A}^{\Lambda}_{e^{i\alpha}\mH}$ are the action functionals from Proposition \autoref{prop:action_fct_of_real_traj} and Remark \autoref{rem:tilted_traj}, respectively. For $\alpha\neq n\cdot \pi$, $n\in\mathbb{Z}$, the complex numbers $ie^{-i\alpha}$ and $-i$ form a $\mathbb{R}$-linear basis of $\mathbb{C}$. Thus, $\gamma$ is a ``critical point'' of $\mathcal{A}^{P_\alpha}_\mH$ iff $\gamma$ is a ``critical point'' of the functionals $\gamma\mapsto\int\mathcal{A}^{\Lambda_R}_{\mH_R} [\gamma_r] dr$ and $\gamma\mapsto\int\text{Re} (\mathcal{A}^\Lambda_{e^{i\alpha}\mH} [\gamma_t]) dt$. The rest now follows as in proof of Lemma \autoref{lem:holo_action_prin} by exploiting Proposition \autoref{prop:action_fct_of_real_traj}, Remark \autoref{rem:tilted_traj}, and the fact that holomorphic trajectories are exactly those smooth curves which are integral curves of $X^R_\mH$ in $t$-direction and integral curves of $\cos (\alpha) X^R_\mH + \sin (\alpha) J(X^R_\mH)$ in $r$-direction ($z = t + re^{i\alpha}$) for\linebreak $\alpha\in\mathbb{R}\backslash\{n\cdot\pi\mid n\in\mathbb{Z}\}$. Summing up our results, we have just shown:

\begin{proposition}[Action principle for parallelograms]\label{prop:holo_action_prin_para}
 Let $(X,\Omega = d\Lambda,\mH)$ be an exact HHS with $\Lambda = \Lambda_R + i\Lambda_I$. For $\alpha\in\mathbb{R}\backslash\{n\cdot\pi\mid n\in\mathbb{Z}\}$, let\linebreak $P_\alpha\coloneqq [t_1,t_2] + e^{i\alpha}[r_1,r_2]\subset\mathbb{C}$ be a parallelogram in the complex plane with real numbers $t_1< t_2$ and $r_1< r_2$. Denote the space of smooth maps from $P_\alpha$ to $X$ by $\mathcal{P}_{P_\alpha}$ and define the action functional $\mathcal{A}^{P_\alpha}_\mH:\mathcal{P}_{P_\alpha}\to\mathbb{C}$ by:
\begin{align*}
 \mathcal{A}^{P_\alpha}_\mH[\gamma]&\coloneqq \iint\limits_{P_\alpha}\left[\Lambda_{R,\gamma (t+is)}\left(2\frac{\partial\gamma}{\partial z}(t+is)\right) - \mH\circ\gamma (t+is)\right] dt\wedge ds\ \text{with}\\
  \frac{\partial \gamma}{\partial z}&\coloneqq \frac{1}{2}\left(\frac{\partial \gamma}{\partial t} - i\frac{\partial \gamma}{\partial s}\right)\quad\forall \gamma\in\mathcal{P}_{P_\alpha}.
\end{align*}
 Now let $\gamma\in\mathcal{P}_{P_\alpha}$ be a smooth map from $P_\alpha$ to $X$. Then, $\gamma$ is a holomorphic trajectory of $(X,\Omega,\mH)$ iff $\gamma$ is a ``critical point''\footnote{Again, ``critical point'' means that only those variations are allowed which keep $\gamma$ fixed on the boundary $\partial P_\alpha$.} of $\mathcal{A}^{P_\alpha}_\mH$.
\end{proposition}

\begin{remark}\label{rem:action_para}
 Proposition \autoref{prop:holo_action_prin_para} is a direct generalization of Lemma \autoref{lem:holo_action_prin}, since one obtains Lemma \autoref{lem:holo_action_prin} from Proposition \autoref{prop:holo_action_prin_para} by setting $\alpha = \pi/2$.
\end{remark}

Before we conclude this subsection, let us inspect Point (ii)b of Remark \autoref{rem:several_remarks} more closely. If a holomorphic curve $\gamma:P_\alpha\to X$ whose domain is a parallelogram $P_\alpha\subset\mathbb{C}$ is doubly-periodic, i.e., periodic in $t$- and $r$-direction for $z = t + re^{i\alpha}$, then we can also view $\gamma$ as a holomorphic map from a complex torus to $X$. In this sense, we can interpret holomorphic trajectories whose domains are complex tori as the holomorphic analogue of periodic orbits. In contrast to periodic orbits of RHSs, however, the domains of two holomorphic periodic orbits do not need to be isomorphic. Indeed, the complex structure of such a torus is determined by the shape of the parallelogram $P_\alpha$. Therefore, the action functional $\mathcal{A}^{P_\alpha}_\mH$ is only sensitive to certain holomorphic periodic orbits, namely those whose domains share the complex structure induced by $P_\alpha$.\\
Non-constant holomorphic periodic orbits are rather rare and do not exist in most HHSs $(X,\Omega,\mH)$. For instance, take $X$ to be the standard example $\mathbb{C}^{2n}$. Due to the compactness of a complex torus $\mathbb{C}/\Gamma$, the maximum principle applies and any holomorphic map $\gamma:\mathbb{C}/\Gamma\to X$ has to be constant. The same result applies if $X$ is Brody hyperbolic\footnote{A complex manifold $X$ is Brody hyperbolic iff every holomorphic map $\gamma:\mathbb{C}\to X$ defined on all of $\mathbb{C}$ is constant.}. Furthermore, if $X$ is compact, then all holomorphic trajectories are constant, since all Hamiltonians are constant by the maximum principle. Still, there are examples of HHSs $(X,\Omega,\mH)$, where a plethora of holomorphic periodic orbits exists.

\begin{example}[Natural Hamiltonians on complex tori $\mathbb{C}^n/\Gamma$]\label{ex:complex_torus}
 Let $n\in\mathbb{N}$ be a natural number and $\Gamma\subset\mathbb{C}^n$ be a lattice, i.e.
 \begin{gather*}
  \Gamma\coloneqq \left\{\sum\limits^{2n}_{j=1} k_j\cdot e_j\middle| k_j\in\mathbb{Z}\right\},
 \end{gather*}
 where the vectors $e_1,\ldots, e_{2n}\in\mathbb{C}^n$ form an $\mathbb{R}$-linear basis of $\mathbb{C}^n$. Then,\linebreak $\mathbb{C}^n/\Gamma$ is a complex torus of complex dimension $n$. Now consider the\linebreak holomorphic cotangent bundle $X\coloneqq T^{\ast, (1,0)}(\mathbb{C}^n/\Gamma)\cong \mathbb{C}^n/\Gamma\times\mathbb{C}^n$ with\linebreak coordinates $([Q_1,\ldots, Q_j], P_1,\ldots, P_j)\in \mathbb{C}^n/\Gamma\times\mathbb{C}^n$ and canonical $2$-form\linebreak $\Omega = \sum^n_{j=1} dP_j\wedge dQ_j$. We want to determine all natural Hamiltonians $\mH = \mathcal{T} + \mathcal{V}$ on the HSM $(X,\Omega)$. The potential energy $\mathcal{V}$ factors through a holomorphic function on $\mathbb{C}^n/\Gamma$. As $\mathbb{C}^n/\Gamma$ is compact, all holomorphic functions on it are constant due to the maximum principle. Since changing the Hamiltonian by a constant does not change the dynamics of the system, we can set the potential energy to zero without loss of generality. To compute the kinetic energy $\mathcal{T}$, we need to classify all holomorphic metrics $g$ on $\mathbb{C}^n/\Gamma$. The projection $\mathbb{C}^n\to\mathbb{C}^n/\Gamma$ gives rise to $n$ linearly independent holomorphic $1$-forms $dQ_j$, $j = 1,\ldots, n$, on the torus $\mathbb{C}^n/\Gamma$ which we have already used to express the canonical form $\Omega$. Using these $1$-forms, we can write $g$ as
 \begin{gather*}
  g = \sum^{2n}_{i,j = 1} g_{ij}dQ_i\otimes dQ_j,
 \end{gather*}
 where $g_{ij}$ are holomorphic functions on the torus. As before, these functions have to be constant implying that the space of holomorphic metrics $g$ on $\mathbb{C}^n/\Gamma$ is isomorphic to the space of symmetric and non-degenerate $\mathbb{C}$-bilinear forms on the complex vector space $\mathbb{C}^n$. By a standard result from linear algebra, every symmetric and non-degenerate $\mathbb{C}$-bilinear form on $\mathbb{C}^n$ can be brought into the standard form $g_{ij} = \delta_{ij}$\footnote{Here, $\delta_{ij}$ is the Kronecker delta!} by a $\mathbb{C}$-linear transformation. Hence, after transforming the lattice $\Gamma$ if necessary, we can assume that the metric $g$ is given by $g = \sum^{2n}_{j = 1} dQ^2_j$. In total, it suffices to investigate the dynamics of the HHS $(\mathbb{C}^n/\Gamma\times\mathbb{C}^n, \sum^n_{j=1} dP_j\wedge dQ_j, \mH)$ with Hamiltonian
 \begin{gather*}
  \mH (Q_1,\ldots, Q_n, P_1,\ldots, P_n) = \sum^{2n}_{j = 1} \frac{P^2_j}{2}
 \end{gather*}
 for all lattices $\Gamma\subset\mathbb{C}^n$ in order to study all natural Hamiltonians on a complex torus. Clearly, the Hamilton equations related to this problem are given by:
 \begin{gather*}
  Q^\prime_j (z) = P_j (z),\quad P^\prime_j (z) = 0.
 \end{gather*}
 Given the initial value $\gamma (0) = ([Q^0_1,\ldots, Q^0_n], P^0_1,\ldots, P^0_n)$, the Hamiltonian equations are solved by the holomorphic trajectory $\gamma:\mathbb{C}\to X$:
 \begin{gather*}
  \gamma (z)\coloneqq ([Q^0_1 + z\cdot P^0_1,\ldots, Q^0_n + z\cdot P^0_n], P^0_1,\ldots, P^0_n).
 \end{gather*}
 
 Let us now define $P^0\coloneqq (P^0_1,\ldots, P^0_n)\in\mathbb{C}^n$ and consider different values for $P^0$:
 \begin{enumerate}
  \item If $P^0 = 0$, then $\gamma$ is just a constant curve.
  \item If $P^0 \neq 0$ and $z\cdot P^0\notin\Gamma$ for every $z\in\mathbb{C}\backslash\{0\}$, then $\gamma$ is a regular holomorphic trajectory with no periodicity.
  \item If $P^0\neq 0$ and $z_1\cdot P^0\in\Gamma$ for at least one $z_1\neq 0$, then $\gamma$ is a regular holomorphic trajectory which is periodic in at least one direction. In this case, we can view $\gamma$ as a holomorphic map from a complex cylinder to $X$.
  \item If $P^0\neq 0$ and $z_1\cdot P^0, z_2\cdot P^0\in\Gamma$ for two $\mathbb{R}$-linearly independent complex numbers $z_1, z_2\in\mathbb{C}$, then $\gamma$ is a regular, doubly-periodic holomorphic trajectory. In this case, $\gamma$ is holomorphic periodic orbit.
 \end{enumerate}
 We observe that the topology and the complex structure of the domain of $\gamma$ changes depending on the momentum $P^0$.
\end{example}

\begin{remark}[General Hamiltonians on a complex torus]\label{rem:ham_on_com_tor}\phantom{x}\\
 As it turns out, Example \autoref{ex:complex_torus} covers all possible Hamiltonians on the HSM\linebreak $(T^{\ast, (1,0)}(\mathbb{C}^n/\Gamma), \sum^{2n}_{j = 1} dP_j\wedge dQ_j)$. Let $\mH$ be any holomorphic function on\linebreak $T^{\ast, (1,0)}(\mathbb{C}^n/\Gamma)$. Since $T^{\ast, (1,0)}(\mathbb{C}^n/\Gamma)$ is isomorphic to $\mathbb{C}^n/\Gamma\times\mathbb{C}^n$, $\mH$ cannot depend on the $Q_j$-coordinates due to maximum principle. This allows us to repeat the discussion from Example \autoref{ex:complex_torus} by replacing $P^0_j$ with $\partial\mH/\partial P_j (P^0)$ in the solution to the Hamilton equations.
\end{remark}

\newpage
\section{Pseudo-Holomorphic Hamiltonian Systems}
\label{sec:PHHS}
In \autoref{sec:HHS}, we have seen that non-constant holomorphic periodic orbits cannot occur for a large class of HHSs $(X,\Omega, \mH)$. One possible obstruction to the existence of such orbits is the integrability of the almost complex structure $J$ of $X$. In this section, we drop the integrability of $J$ leading us to the notion of a pseudo-holomorphic Hamiltonian system (PHHS). We demonstrate in the first subsection that PHHSs exhibit, by design, almost the same properties as found for HHSs in \autoref{sec:HHS}, in particular with respect to the existence and uniqueness of pseudo-holomorphic trajectories and with respect to action functionals and principles. In the second subsection, we explore the relation between HHSs and PHHSs and show that we recover a HHS from a PHHS if we restore the integrability of its almost complex structure $J$.

\subsection*{PHHS: Basic Definitions, Notions, and Properties}

In \autoref{sec:HHS}, we have found that most HHSs $(X,\Omega,\mH)$ do not possess non-constant holomorphic periodic orbits. Often, their existence was forbidden by the maximum principle. Consider, for instance, $X = \mathbb{R}^{4}\cong \mathbb{C}^2$ with the standard complex structure $J = i$. We recall that a holomorphic periodic orbit is a holomorphic map $\gamma:\mathbb{C}/\Gamma\to X$ satisfying Hamilton's equations, where $\C/\Gamma$ is a complex torus and, thus, compact. Therefore, the maximum principle applies and every holomorphic periodic orbit in $\mathbb{C}^2$ is constant.\\
In his beautiful paper \cite{moser1995} from 1995, Moser showed that the same argument does \underline{not} apply if we equip $\mathbb{R}^4$ with a different almost complex structure $J$. Let $\gamma$ be any smooth embedding of the $2$-torus into $\mathbb{R}^4$, e.g. the inclusion\linebreak $S^1\times S^1 \subset \mathbb{R}^2\times\mathbb{R}^2\equiv \mathbb{R}^4$. Then, the image of $\gamma$ is a $2$-dimensional submanifold of $\mathbb{R}^4$ and its tangent bundle can be continued to a smooth $2$-dimensional distribution $D$ on $\mathbb{R}^4$. This can be seen as follows: The tangent bundle of $S^1\times S^1\subset\mathbb{R}^4$ is spanned by the two vector fields $V_1$ and $V_2$ on $S^1\times S^1$:
\begin{align*}
 V_1:&S^1\times S^1\to\mathbb{R}^4,\ (x_1,x_2,x_3,x_4)\mapsto (-x_2,x_1,0,0),\\ V_2:&S^1\times S^1\to\mathbb{R}^4,\ (x_1,x_2,x_3,x_4)\mapsto (0,0, -x_4, x_3).
\end{align*}
We show that there are two linearly independent vector fields $\hat V_1$ and $\hat V_2$ on $\mathbb{R}^4$ continuing $V_1$ and $V_2$, i.e., $\hat V_i\vert_{S^1\times S^1} = V_i$. To construct $\hat V_1$ and $\hat V_2$, we first define the functions $r_1\coloneqq \sqrt{x^2_1 + x^2_2}$, $r_2\coloneqq \sqrt{x^2_3 + x^2_4}$, and $R\coloneqq (1-r^2_1)^2 + (1-r^2_2)^2$. Next, define the vector fields $\hat V_1$, $W_1$, and $W_2$ on $\mathbb{R}^4$ as follows:
\begin{align*}
 \hat V_1:& \mathbb{R}^4\to\mathbb{R}^4,\ (x_1,x_2,x_3,x_4)\mapsto (-x_2,x_1,R,0),\\
 W_1:& \mathbb{R}^4\to\mathbb{R}^4,\ (x_1,x_2,x_3,x_4)\mapsto (-x_2x_4R,\, x_1x_4R,\, -r^2_1x_4,\, r^2_1x_3),\\
 W_2:& \mathbb{R}^4\to\mathbb{R}^4,\ (x_1,x_2,x_3,x_4)\mapsto (x_1,x_2,0,R).
\end{align*}
One easily checks that $\hat V_1$ is a continuation of $V_1$, vanishes nowhere, and is orthogonal to $W_1$ and $W_2$ with respect to the standard metric on $\mathbb{R}^4$. Furthermore, we notice that $W_1$ is a continuation of $V_2$. However, $W_1$ vanishes for $x_1 = x_2 = 0$ or $x_3 = x_4 = 0$. To rectify this, we take an appropriate linear combination of $W_1$ and $W_2$. For that, we first observe that $W_1$ does not vanish on $S^1\times S^1$. Hence, we can pick an open neighborhood $U\subset\mathbb{R}^4$ of $S^1\times S^1$ such that $W_1$ does not vanish on $U$. Next, we choose a partition of unity $\{f_1, f_2\}$ on $\mathbb{R}^4$ subordinate to the open covering $\{U, \mathbb{R}^4\backslash (S^1\times S^1)\}$ of $\mathbb{R}^4$, i.e, two smooth functions $f_1,f_2\in C^\infty (\mathbb{R}^4,\mathbb{R}_{\geq 0})$ satisfying:
\begin{enumerate}
 \item $f_1(x) + f_2(x) = 1\quad\forall x\in\mathbb{R}^4$,
 \item $\text{supp}(f_1)\subset U$ and $\text{supp}(f_2)\subset\mathbb{R}^4\backslash (S^1\times S^1)$.
\end{enumerate}
Now define the vector field $\hat V_2$ by $\hat V_2\coloneqq f_1\cdot W_1 + f_2\cdot W_2$. By construction, the vector field $\hat V_2$ is a continuation of $V_2$. Moreover, one can show that $\hat V_2$ vanishes nowhere by considering $\hat V_2$ separately on $S^1\times S^1$, $U\backslash (S^1\times S^1)$, and $\mathbb{R}^4\backslash U$. As $\hat V_1$ is orthogonal to $W_1$ and $W_2$, $\hat V_1$ is also orthogonal to the vector field $\hat V_2$. Two orthogonal vector fields which vanish nowhere are linearly independent, hence, the vector fields $\hat V_1$ and $\hat V_2$ are the desired continuations of $V_1$ and $V_2$. The distribution $D$ is now just the span of $\hat V_1$ and $\hat V_2$.\\
Let us return to Moser's construction. Choose a Riemannian metric $g$ on $\mathbb{R}^4$ and consider the orthogonal complement $D^\perp$ of $D$ with respect to $g$. $D^\perp$ is also a smooth $2$-dimensional distribution on $\mathbb{R}^4$. Moreover, $D$ and $D^\perp$ span the tangent bundle of $\mathbb{R}^4$, i.e., $T\mathbb{R}^4 = D\oplus D^\perp$. We can now construct the almost complex structure $J$ on $\mathbb{R}^4$ as follows: Choose orientations for $D$ and $D^\perp$ and define $J$ to be the $90^{\circ}$-rotation in $D$ and $D^\perp$ with respect to $g$ and the given orientations. After choosing a suitable complex structure $j$ on the $2$-torus, $\gamma$ becomes a pseudo-holomorphic\footnote{In most books and papers, a map $f:(X_1,J_1)\to (X_2, J_2)$ between manifolds $X_1$ and $X_2$ with almost complex structures $J_1$ and $J_2$ is called holomorphic if $df\circ J_1 = J_2\circ df$. In this thesis, we often want to emphasize the non-integrability of $J_1$ or $J_2$. To that extent, we say $f:(X_1,J_1)\to(X_2,J_2)$ is \textbf{pseudo-holomorphic} if $df\circ J_1 = J_2\circ df$ and usually save the expression ``\textbf{holomorphic}'' for the case where $J_1$ and $J_2$ are integrable.} embedding, i.e., $d\gamma\circ j = J\circ d\gamma$. The almost complex structure $J$ constructed this way is, in general, \underline{not} integrable.\\
Moser's example indicates that Hamiltonian systems with non-integrable almost complex structures $J$ might be richer than HHSs when it comes to pseudo-holomorphic periodic orbits. However, the generalization of HHSs is not straightforward, as the complex structure $J$ only enters most definitions regarding HHSs implicitly. To that end, let us recapitulate which objects and relations are essential to the definitions and discussions in \autoref{sec:HHS}. A HHS consists of six objects: A smooth manifold $X$ together with an integrable almost complex structure $J$ on it, two real $2$-forms $\Omega_R$ and $\Omega_I$ on $X$ which assemble to a holomorphic symplectic form $\Omega = \Omega_R + i\Omega_I$, and two smooth real functions $\mH_R$ and $\mH_I$ on $X$ forming a holomorphic function $\mH = \mH_R + i\mH_I$ on $X$. Closely tracing every step of \autoref{sec:HHS}, we see that these six objects need to satisfy the following relations:
\begin{enumerate}
 \item $\Omega_R$ must be closed\footnote{For the action functionals of holomorphic trajectories in \autoref{sec:HHS}, we need a primitive of $\Omega_R$, but not of $\Omega_I$!}.
 \item $J$, $\Omega_R$, and $\Omega_I$ need to satisfy the relations induced by Equation \eqref{eq:J-anticompatible}.
 \item $J$ and the Hamiltonian vector fields of the underlying RHSs have to fulfill Cauchy-Riemann-like relations formulated in Remark \autoref{rem:cauchy-riemann}.
 \item All Hamiltonian vector fields must commute reproducing Corollary \autoref{cor:commute}.
\end{enumerate}
One way to define a pseudo-holomorphic Hamiltonian system is to simply impose these relations:

\begin{definition}[Pseudo-holomorphic Hamiltonian system]\label{def:PHHS_1}
 We call a collection $(X,J;\Omega_R,\Omega_I;\mH_R,\mH_I)$ a \textbf{pseudo-holomorphic Hamiltonian system} (PHHS) if $X$ is a smooth manifold, $J$ is a (not necessarily integrable) almost complex structure on $X$, $\Omega_R\in\Omega^2(X)$ and $\Omega_I\in\Omega^2(X)$ are non-degenerate, and $\mH_R\in C^\infty(X,\R)$ and $\mH_I\in C^\infty(X,\R)$ are smooth functions satisfying:
 \begin{enumerate}
  \item $\Omega_R$ is closed, i.e., $d\Omega_R = 0$.
  \item $\Omega_R(J\cdot,\cdot) = \Omega_R(\cdot,J\cdot) = -\Omega_I,\quad \Omega_I(J\cdot,\cdot) = \Omega_I(\cdot,J\cdot) = \Omega_R,$\\
  $\Omega_R (J\cdot,J\cdot) = -\Omega_R,\qquad\qquad\quad\ \Omega_I (J\cdot,J\cdot) = -\Omega_I$.
  \item $X^{\Omega_R}_{\mH_R} = X^{\Omega_I}_{\mH_I}$ and $J\left(X^{\Omega_R}_{\mH_R}\right) = X^{\Omega_I}_{\mH_R} = - X^{\Omega_R}_{\mH_I}$, where $X^{\Omega_a}_{\mH_b}$ is defined by $\iota_{X^{\Omega_a}_{\mH_b}}\Omega_a = -d\mH_b$.
  \item $[X^{\Omega_a}_{\mH_b}, X^{\Omega_c}_{\mH_d}] = 0$ for all $a,b,c,d\in\{R,I\}$.
 \end{enumerate}
\end{definition}

\begin{remark}[Property (iii) in Definition \autoref{def:PHHS_1}]\label{rem:H_pseudo-holo}
 Note that we can replace Property (iii) in Definition \autoref{def:PHHS_1} with the condition that $\mH\coloneqq \mH_R + i\cdot \mH_I$ is pseudo-holomorphic. Indeed, Property (ii) and (iii) imply:
 \begin{align*}
  d\mH\circ J &= -\Omega_R (X^{\Omega_R}_{\mH_R}, J\cdot) - i\cdot \Omega_R (X^{\Omega_R}_{\mH_I}, J\cdot) = \Omega_I (X^{\Omega_I}_{\mH_I}, \cdot) + i\cdot\Omega_I (-X^{\Omega_I}_{\mH_R},\cdot)\\
  &= -d\mH_I + i\cdot d\mH_R = i\cdot d\mH.
 \end{align*}
 Conversely, Property (ii) and $d\mH\circ J = i\cdot d\mH$ imply the Cauchy-Riemann-like equations in Property (iii).
\end{remark}

\begin{remark}[The form $\Omega$, Part I]\label{rem:omega_I}
 As for complex manifolds, we can decompose the complexified tangent\footnote{Of course, similar remarks apply to the complexified cotangent bundle $T^{\ast}_\mathbb{C}X$ of $X$.} bundle $T_\mathbb{C}X$ of a manifold $X$ with almost complex structure $J$ into a direct sum of subbundles $T^{(1,0)}X$ and $T^{(0,1)}X$ which are fiberwise given by the eigenspaces of $J$ with eigenvalue $i$ and $-i$, respectively. The difference is that $T^{(1,0)}X$ is not involutive anymore, but merely a smooth complex vector bundle over $X$. Still, if we define the complex $2$-form $\Omega$ to be $\Omega_R + i\Omega_I$ for PHHSs, we find that $\Omega$ is of type $(2,0)$, i.e., vanishes on $T^{(0,1)}X$.
\end{remark}

One might be confused why we only require $\Omega_R$ to be closed. The reason is that if we were to include the closedness of $\Omega_I$ into the definition of a PHHS, the almost complex structure $J$ would automatically be integrable rendering our construction pointless. The proof of this statement is given in the second subsection, when we explore the relation between HHSs and PHHSs.\\
Definition \autoref{def:PHHS_1} is convoluted, redundant, and rather unwieldy. For a better approach to PHHSs, let us first define pseudo-holomorphic symplectic manifolds:

\begin{definition}[Pseudo-holomorphic symplectic manifolds]\label{def:PHSM}
 We call a triple $(X,J;\Omega_R)$ \textbf{pseudo-holomorphic symplectic manifold} (PHSM) if $(X,\Omega_R)$ is a symplectic manifold and $J$ is an almost complex structure on $X$ which is also $\Omega_R$-anticompatible, i.e. $\Omega_R (J\cdot, J\cdot) = -\Omega_R$. A PHSM $(X,J;\Omega_R)$ is called \textbf{proper} if $J$ is not integrable.
\end{definition}

\begin{remark}[The form $\Omega$, Part II]\label{rem:omega_II}
 Every PHSM $(X,J;\Omega_R)$ possesses forms $\Omega_I$ and $\Omega$ defined by:
 \begin{gather*}
  \Omega_I\coloneqq -\Omega_R (J\cdot,\cdot),\quad \Omega\coloneqq \Omega_R + i\Omega_I.
 \end{gather*}
 It is easy to see that $\Omega_I$ is a smooth, non-degenerate, alternating $2$-form on $X$ which is also anticompatible with $J$. Furthermore, $\Omega$ is also anticompatible with $J$, satisfies $\Omega (J\cdot,\cdot) = \Omega (\cdot, J\cdot) = i\Omega$, i.e., $\Omega$ is of type $(2,0)$, and is non-degenerate on $T^{(1,0)}X$. However, neither $\Omega_I$ nor $\Omega$ are necessarily closed.
\end{remark}

We can now give an alternative definition of a PHHS:

\begin{definition}[Pseudo-holomorphic Hamiltonian system]\label{def:PHHS_2}
 We call a collection $(X,J;\Omega_R, \mH_R)$ a \textbf{pseudo-holomorphic Hamiltonian system} (PHHS) if $(X,J;\Omega_R)$ is a PHSM, $\mH_R:X\to\mathbb{R}$ is a smooth function on $X$, and the $1$-form $\Omega_R (J(X^{\Omega_R}_{\mH_R}),\cdot)$ is exact, where $\Omega_R (X^{\Omega_R}_{\mH_R},\cdot)\coloneqq -d\mH_R$. We call a PHHS $(X,J;\Omega_R,\mH_R)$ \textbf{proper} if $J$ is not integrable.
\end{definition}

For both definitions to agree, it is obvious that the exactness condition is necessary, since, by Definition \autoref{def:PHHS_1}, we have $\Omega_R (J(X^{\Omega_R}_{\mH_R}),\cdot) = d\mH_I$. The following proposition ensures that the exactness condition is also sufficient:

\begin{proposition}[PHHSs well-defined]\label{prop:def_PHHS_agree}
 Definition \autoref{def:PHHS_1} and Definition \autoref{def:PHHS_2} coincide.
\end{proposition}

\begin{proof}
 Clearly, every PHHS as in Definition \autoref{def:PHHS_1} also fulfills Definition \autoref{def:PHHS_2}. Now let $(X,J;\Omega_R,\mH_R)$ be a PHHS as in Definition \autoref{def:PHHS_2}. Then, we define $\Omega_I\coloneqq -\Omega_R (J\cdot,\cdot)$ and take $\mH_I$ to be a primitive of the $1$-form $\Omega_R (J(X^{\Omega_R}_{\mH_R}),\cdot)$. We need to check that these data satisfy (i), (ii), (iii), and (iv) from Definition \autoref{def:PHHS_1}. Property (i) is trivially true by definition. Verifying Property (ii) is a short and easy computation. To check Property (iii), we recall Remark \autoref{rem:H_pseudo-holo}. It suffices to verify that the map $\mH = \mH_R + i\mH_I:X\to\mathbb{C}$ is pseudo-holomorphic which follows immediately:
 \begin{align*}
  d\mH_R\circ J &= -\Omega_R (X^{\Omega_R}_{\mH_R},J\cdot) = -\Omega_R (J(X^{\Omega_R}_{\mH_R}),\cdot) = -d\mH_I,\\
  d\mH_I\circ J &= \Omega_R (J(X^{\Omega_R}_{\mH_R}),J\cdot) = -\Omega_R (X^{\Omega_R}_{\mH_R},\cdot) = d\mH_R,\\
  \Rightarrow d\mH\circ J &= i\cdot d\mH.
 \end{align*}
 Lastly, we need to check Property (iv). Recall that any symplectic manifold $(X,\Omega_R)$ admits a Poisson bracket $\{\cdot,\cdot\}:C^{\infty}(X,\mathbb{R})\times C^{\infty}(X,\mathbb{R})\to C^{\infty}(X,\mathbb{R})$ given by:
 \begin{gather*}
  \{F,G\}\coloneqq \Omega_R (X_F, X_G),
 \end{gather*}
 where $X_F$ and $X_G$ are the Hamiltonian vector fields of the functions $F$ and $G$. Furthermore, remember that the map $X_\cdot:C^{\infty}(X,\mathbb{R})\to\Gamma (TX)$ is a Lie algebra homomorphism:
 \begin{gather*}
  X_{\{F,G\}} = [X_F, X_G]\quad\forall F,G\in C^{\infty}(X,\mathbb{R}).
 \end{gather*}
 Hence, it suffices to prove that $\{\mH_R,\mH_I\}$ vanishes in order to show that $X^{\Omega_R}_{\mH_R}$ and $X^{\Omega_R}_{\mH_I}$ commute. Let us calculate $\{\mH_R,\mH_I\}$ using Property (ii) and (iii):
 \begin{gather*}
  \{\mH_R,\mH_I\} = \Omega_R (X^{\Omega_R}_{\mH_R}, X^{\Omega_R}_{\mH_I}) = -\Omega_R (X^{\Omega_R}_{\mH_R}, J(X^{\Omega_R}_{\mH_R})) = \Omega_I (X^{\Omega_R}_{\mH_R},X^{\Omega_R}_{\mH_R}) = 0.
 \end{gather*}
 Commutativity of the remaining Hamiltonian vector fields follows from commutativity of $X^{\Omega_R}_{\mH_R}$ and $X^{\Omega_R}_{\mH_I}$ as well as Property (iii) concluding the proof.
\end{proof}

Now that we have found a compact definition of PHHSs, we should briefly mention some examples of PHHSs. Of course, every HHS is, by design, a PHHS with integrable $J$. Even though the set of proper PHHSs is much larger than the set of HHSs, finding them is a bit more involved and, thus, relegated to \autoref{sec:construction_and_deformation}. Partially, this is due to the fact that there are no ``standard'' examples of proper PHHSs like cotangent bundles\footnote{At least no canonical ones! In \autoref{sec:construction_and_deformation}, we will equip the holomorphic cotangent bundle of a complex manifold with a non-canonical PHHS structure.} as for RHSs and HHSs. On a deeper level, this is caused by the absence of a Darboux-like theorem. Clearly, there cannot be a counterpart to Darboux's theorem for PHSMs, since $J$ is usually not integrable and, hence, there are no coordinates in which $J$ assumes its standard form, let alone coordinates in which both $J$ and $\Omega = \Omega_R + i\Omega_I$ assume some standard form. Still, one can bring $J$ and $\Omega$ into a standard form using local frames:

\begin{lemma}[PHSMs in local frames]\label{lem:PHSM_in_sta_form}
 Let $(X,J;\Omega_R)$ be a PHSM with\linebreak $\Omega\coloneqq \Omega_R - i\Omega_R (J\cdot,\cdot)$ and let $x_0\in X$ be any point. Then, there exists an open neighborhood $U\subset X$ of $x_0$ and a local frame $\theta^Q_1,\ldots, \theta^Q_n,\theta^P_1,\ldots,\theta^P_n$ of the smooth complex vector bundle $T^{\ast, (1,0)}X$ on $U$ such that:
 \begin{gather*}
  \Omega\vert_U = \sum^{n}_{j=1} \theta^P_j\wedge \theta^Q_j.
 \end{gather*}
 In particular, the real dimension of $X$ is a multiple of $4$. The local frame can be chosen to be integrable (meaning the frame is induced by a chart) if and only if $J$ is integrable on $U$.
\end{lemma}

\begin{remark}[$J$ in standard form]\label{rem:J_sta}
 $J$ is also in standard form in the dual frame of $\theta^Q_1,\ldots, \theta^Q_n,\theta^P_1,\ldots,\theta^P_n$. Indeed, the real and imaginary part of the local frame $\theta^Q_1,\ldots, \theta^Q_n,\theta^P_1,\ldots,\theta^P_n$ ($\theta = \theta^x + i \theta^y$) give rise to a local frame of the real cotangent bundle $T^{\ast}X$. Its dual frame ${\hat e}^{Q, x}_{1},\ldots, {\hat e}^{P,y}_{n}$ is a local frame of the tangent bundle $TX$. By setting ${\hat e}\coloneqq 1/2({\hat e}^x - i{\hat e}^y)$, one obtains a local frame of $T^{(1,0)}X$. On $T^{(1,0)}X$, $J$ simply acts by $i$, thus, ${\hat e}^x$ and ${\hat e}^y$ satisfy $J({\hat e}^x) = {\hat e}^y$ and $J({\hat e}^y) = -{\hat e}^x$. This is the standard form of $J$.
\end{remark}

\begin{proof}
 Lemma \autoref{lem:PHSM_in_sta_form} follows from the application of the symplectic Gram-Schmidt process, which can be found in any textbook on symplectic geometry, to a local frame of $T^{(1,0)}X$. Confer Proposition 2.8 in \cite{bogomolov2020} for the complex analogue of the symplectic Gram-Schmidt process. For completeness' sake, we repeat the explicit construction here. Let $(X,J;\Omega_R)$ be a PHSM with $\Omega\coloneqq \Omega_R - i\Omega_R (J\cdot,\cdot)$ and take the real dimension of $X$ to be $\text{dim}_\mathbb{R} (X) = 2m$, $m\in\mathbb{N}$. The dimension of $X$ is even, as $X$ admits an almost complex structure $J$. Then, the complex rank of the complexified bundle $T_\mathbb{C}X$ is also given by $2m$. Now recall the decomposition $T_\mathbb{C}X = T^{(1,0)}X\oplus T^{(0,1)}X$. Since the vector bundles $T^{(1,0)}X$ and $T^{(0,1)}X$ are isomorphic via the complex conjugation $v + iw\mapsto v- iw$, their fibers have the same complex dimension, namely $m$. Now let $x_0\in X$ be any point and pick a local frame $v_1,\ldots, v_m$ of $T^{(1,0)}X$ on an open neighborhood $U\subset X$ of $x_0$. $\Omega$ is non-degenerate on $T^{(1,0)}X$ by Remark \autoref{rem:omega_II}, hence, there exists a vector ${\hat e}^Q_{1,x_0}\in T^{(1,0)}_{x_0}X$ such that $\Omega_{x_0} (v_1(x_0), {\hat e}^Q_{1,x_0})\neq 0$. $v_1,\ldots, v_m$ is a local frame of $T^{(1,0)}X$ near $x_0$, thus, we can write:
 \begin{gather*}
  {\hat e}^Q_{1,x_0} = \sum^m_{j = 1} c_j\cdot v_j(x_0),
 \end{gather*}
 where $c_j\in\mathbb{C}$. Now define the local section ${\hat e}^Q_1 \coloneqq \sum_{j} c_j\cdot v_j$ of $T^{(1,0)}X$. After shrinking $U$ if necessary, one obtains $\Omega_x (v_1(x), {\hat e}^Q_1 (x))\neq 0$ for every $x\in U$. Setting ${\hat e}^P_1\coloneqq v_1$ and changing the normalization of ${\hat e}^Q_1$ if necessary allows us to write $\Omega\vert_U ({\hat e}^P_1, {\hat e}^Q_1) = 1$.\\
 If $m = 2$, we simply define $\theta^Q_1, \theta^P_1$ to be the dual frame of ${\hat e}^Q_1, {\hat e}^P_1$. If $m> 2$, then we can pick one local section of the frame $v_1,\ldots, v_m$, say $v_2$, such that ${\hat e}^Q_1 (x)$, ${\hat e}^P_1 (x)$, and $v_2(x)$ are $\mathbb{C}$-linearly independent for every $x\in U$ (after shrinking $U$ if necessary). We set:
 \begin{gather*}
  \hat v_2\coloneqq v_2 - \Omega\vert_U (v_2, {\hat e}^Q_1)\cdot {\hat e}^P_1 + \Omega\vert_U (v_2, {\hat e}^P_1)\cdot {\hat e}^Q_1.
 \end{gather*}
 Then, ${\hat e}^Q_1$, ${\hat e}^P_1$, and $\hat v_2$ are still $\mathbb{C}$-linearly independent on $U$ and $\hat v_2$ is $\Omega$-orthogonal to ${\hat e}^Q_1$ and ${\hat e}^P_1$, i.e., $\Omega\vert_U (\hat v_2, {\hat e}^Q_1) = \Omega\vert_U (\hat v_2, {\hat e}^P_1) = 0$. Again by the non-degeneracy of $\Omega$, we can find a vector $e^Q_{2,x_0}\in T^{(1,0)}_{x_0}X$ such that $\Omega_{x_0} (\hat v_2 (x_0), e^Q_{2,x_0})\neq 0$. As before, we can write:
 \begin{gather*}
  {e}^Q_{2,x_0} = \sum^m_{j = 1} d_j\cdot v_j (x_0),
 \end{gather*}
 where $d_j\in\mathbb{C}$, and define the local section ${e}^Q_2 \coloneqq \sum_{j} d_j\cdot v_j$ of $T^{(1,0)}X$. After shrinking $U$ and changing the normalization of $e^Q_2$ if necessary, we obtain $\Omega\vert_U (\hat v_2, e^Q_2) = 1$. Now we set:
 \begin{gather*}
  {\hat e}^P_2\coloneqq \hat v_2,\quad {\hat e}^Q_2\coloneqq e^Q_2 - \Omega\vert_U (e^Q_2, {\hat e}^Q_1)\cdot {\hat e}^P_1 + \Omega\vert_U (e^Q_2, {\hat e}^P_1)\cdot {\hat e}^Q_1.
 \end{gather*}
 Proceeding inductively gives us a local frame ${\hat e}^Q_1,\ldots {\hat e}^Q_n$, ${\hat e}^P_1,\ldots, {\hat e}^P_n$ of $T^{(1,0)}X$ on some neighborhood $U$ of $x_0$ ($n\coloneqq m/2$) satisfying:
 \begin{gather*}
  \Omega\vert_U ({\hat e}^Q_i, {\hat e}^Q_j) = \Omega\vert_U ({\hat e}^P_i, {\hat e}^P_j) = 0,\quad \Omega\vert_U ({\hat e}^P_i, {\hat e}^Q_j) = \delta_{ij}.
 \end{gather*}
 Thus, the frame $\theta^Q_1,\ldots, \theta^Q_n$, $\theta^P_1,\ldots, \theta^P_n$ dual to the frame ${\hat e}^Q_1,\ldots {\hat e}^Q_n$, ${\hat e}^P_1,\ldots, {\hat e}^P_n$ is the desired local frame of $T^{\ast, (1,0)}X$ near $x_0$ in which $\Omega$ takes the form:
 \begin{gather*}
  \Omega\vert_U = \sum^{n}_{j=1} \theta^P_j\wedge \theta^Q_j.
 \end{gather*}
 In particular, the real dimension of $X$ is $4n$. Lastly, we show the equivalence.\linebreak If the frame $\theta^Q_1,\ldots, \theta^P_n$ is integrable, then there exists a chart\linebreak $\phi = (Q_1,\ldots, Q_n, P_1,\ldots, P_n):U\to V\subset\mathbb{C}^{2n}$ near $x_0$ such that $\theta^Q_j = dQ_j$ and $\theta^P_j = dP_j$. Because $\theta^Q_j$ and $\theta^P_j$ are forms of type $(1,0)$, we obtain $dQ_j\circ J = idQ_j$ and $dP_j\circ J = idP_j$. Thus, $\phi$ is a holomorphic chart and $J$ is integrable on $U$. The converse direction follows from Theorem \autoref{thm:rel_HSM_PHSM} (cf. the second subsection) and Darboux's theorem for HSMs (cf. Theorem \autoref{thm:holo_Darboux} and \autoref{app:darboux}).
\end{proof}

Next, let us investigate the dynamics of a PHHS $(X,J;\Omega_R,\mH_R)$. To do so, we need to introduce pseudo-holomorphic Hamiltonian vector fields and trajectories. We try to imitate the definitions from \autoref{sec:HHS}. For the sake of simplicity, we always associate from now on with a PHHS $(X,J;\Omega_R,\mH_R)$ the forms $\Omega_I\coloneqq -\Omega_R (J\cdot,\cdot)$ and $\Omega\coloneqq \Omega_R + i\Omega_I$ as well as the functions $\mH_I$ and $\mH\coloneqq \mH_R + i\mH_I$, where $\mH_I$ is a primitive of $\Omega_R (J(X^{\Omega_R}_{\mH_R}),\cdot)$.

\begin{definition}[Pseudo-holomorphic Hamiltonian vector fields]\label{def:pseudo-holo_ham_field_and_traj}
 Let\linebreak $(X,J;\Omega_R, \mH_R)$ be a PHHS. We call the smooth section $X_\mH$ of $T^{(1,0)}X$ defined by $\iota_{X_\mH}\Omega = -d\mH$ the (pseudo-holomorphic) \textbf{Hamiltonian vector field} of the PHHS $(X,J;\Omega_R, \mH_R)$. Furthermore, we call a smooth map $\gamma:U\to X$, $U\subset\mathbb{C}$ open and connected, a \textbf{pseudo-holomorphic trajectory} of the PHHS $(X,J;\Omega_R, \mH_R)$ if $\gamma$ satisfies the integral curve equation:
 \begin{gather*}
  \frac{\partial\gamma}{\partial z} (z)\coloneqq \frac{1}{2}\left(\frac{\partial \gamma}{\partial t} (z) - i\frac{\partial\gamma}{\partial s} (z)\right) = X_\mH (\gamma (z))\quad\forall z = t + is\in U.
 \end{gather*}
 We call a pseudo-holomorphic trajectory $\gamma:U\to X$ \textbf{maximal} if for every pseudo-holomorphic trajectory $\hat\gamma:\hat U\to X$ with $U\subset \hat U$ and $\hat\gamma\vert_U = \gamma$ one has $\hat U = U$.
\end{definition}

Alternatively, one can define pseudo-holomorphic Hamiltonian vector fields and trajectories in terms of the vector field $X^{\Omega_R}_{\mH_R}$:

\begin{proposition}[Alternative Definition of $X_\mH$ and $\gamma$]\label{prop:pseudo-holo_ham_field_and_traj}
 Let $(X,J;\Omega_R,\mH_R)$ be a PHHS with vector field $X^{\Omega_R}_{\mH_R}$ defined by $\Omega_R (X^{\Omega_R}_{\mH_R},\cdot) = -d\mH_R$. Then:
 \begin{gather*}
  X_\mH = \frac{1}{2} \left( X^{\Omega_R}_{\mH_R} - i\cdot J(X^{\Omega_R}_{\mH_R})\right).
 \end{gather*}
 Now let $\gamma:U\to X$ be a smooth map, where $U\subset\mathbb{C}$ is an open and connected subset. Then, $\gamma$ is a pseudo-holomorphic trajectory iff $\gamma_s$ defined by\linebreak $\gamma_s (t)\coloneqq \gamma (t + is)$ is an integral curve of $X^{\Omega_R}_{\mH_R}$ for every suitable $s\in\mathbb{R}$ and $\gamma:U\to X$ is pseudo-holomorphic.
\end{proposition}

\begin{proof}
 A straightforward calculation verifies that $X_\mH$ as in Definition \autoref{def:pseudo-holo_ham_field_and_traj} satisfies the equation above. If $\gamma$ is a pseudo-holomorphic trajectory, then the real part of the pseudo-holomorphic integral curve equation implies that $\gamma_s$ is an integral curve of $X^{\Omega_R}_{\mH_R}$. Moreover, $\gamma$ satisfies:
 \begin{gather*}
  \frac{\partial\gamma}{\partial s} = -2\text{Im} X_\mH = J(X^{\Omega_R}_{\mH_R}´) = J\frac{\partial\gamma}{\partial t}.
 \end{gather*}
 This equation is equivalent to $\gamma$ being pseudo-holomorphic. The converse direction works similarly.
\end{proof}

We designed PHHSs in such a way that all properties we found for HHSs in \autoref{sec:PHHS} transfer almost completely to PHHSs. For instance, we find the following PHHS-counterpart to Proposition \autoref{prop:holo_traj}:

\begin{proposition}[Existence and uniqueness of pseudo-holomorphic trajectories]\label{prop:pseudo-holo_traj}
 Let $(X,J;\Omega_R,\mH_R)$ be a PHHS. Then, for any $z_0\in\mathbb{C}$ and $x_0\in X$, there exists an open and connected subset $U\subset\mathbb{C}$ and a pseudo-holomorphic trajectory $\gamma^{z_0, x_0}:U\to X$ of $(X,J;\Omega_R,\mH_R)$ with $\gamma^{z_0, x_0} (z_0) = x_0$. Two pseudo-holomorphic trajectories $\gamma^{z_0, x_0}_1:U_1\to X$ and $\gamma^{z_0, x_0}_2:U_2\to X$ with $\gamma^{z_0, x_0}_1 (z_0) = x_0 = \gamma^{z_0, x_0}_2 (z_0)$ locally coincide, in particular, they are equal iff their domains $U_1$ and $U_2$ are equal. Furthermore, the pseudo-holomorphic trajectory $\gamma^{z_0, x_0}$ depends pseudo-holomorphically on $z_0$, but, in general, \underline{only} smoothly on $x_0$.
\end{proposition}

\begin{proof}
 The proof of Proposition \autoref{prop:pseudo-holo_traj} works very similarly to the proof of Proposition \autoref{prop:holo_traj}. By Definition \autoref{def:PHHS_1} and Proposition \autoref{prop:pseudo-holo_ham_field_and_traj}, the real and imaginary part of $X_\mH$ commute, hence, we can proceed as in the proof of Proposition \autoref{prop:holo_traj} to show that pseudo-holomorphic trajectories $\gamma^{z_0,x_0}$ given an initial value $\gamma^{z_0, x_0} (z_0) = x_0$ exist.\\
 To prove uniqueness, we observe that the formula
 \begin{align*}
  \gamma^{z_0, x_0}(z)&\coloneqq \varphi^{J(X^{\Omega_R}_{\mH_R})}_{s-s_0}\circ\varphi^{X^{\Omega_R}_{\mH_R}}_{t-t_0} (x_0)\equiv \varphi^{X^{\Omega_R}_{\mH_R}}_{t-t_0}\circ\varphi^{J(X^{\Omega_R}_{\mH_R})}_{s-s_0} (x_0)\\
  &\equiv\varphi^{(t-t_0)X^{\Omega_R}_{\mH_R} + (s-s_0)J(X^{\Omega_R}_{\mH_R})}_1 (x_0),
 \end{align*}
 where $z = t+is$, uniquely determines $\gamma^{z_0, x_0}$ on a small rectangle in $\mathbb{C}$ near $z_0 = t_0 + is_0$. The rest now follows by covering a path between $z_0$ and any point $z_1$ in $U_1\equiv U_2$ with a finite number of such rectangles.\\
 Lastly, let us consider the dependence of $\gamma^{z_0, x_0}$ on $z_0\in\mathbb{C}$ and $x_0\in X$. Again, $\gamma^{z_1, x_0} (z)$ and $\gamma^{z_2, x_0} (z)$ only differ by a translation in $z$. As pseudo-holomorphic trajectories are pseudo-holomorphic maps, the $z_0$-dependence is also pseudo-holomorphic. For the $x_0$-dependence, we need to consider the flow of $X^{\Omega_R}_{\mH_R}$ and $J(X^{\Omega_R}_{\mH_R})$. Both $X^{\Omega_R}_{\mH_R}$ and $J(X^{\Omega_R}_{\mH_R})$ are smooth vector fields, thus, their flows are smooth as well concluding the proof.
\end{proof}

\begin{remark}[$x_0$-dependence]\label{rem:x_0-dependance}
 Note that holomorphic trajectories of HHSs depend holomorphically on $x_0$, while pseudo-holomorphic trajectories of PHHSs do \underline{not} generally depend pseudo-holomorphically on $x_0$. This distinction can be traced back to the Hamiltonian vector field $X_\mH$. For HHSs, $X_\mH$ is a holomorphic vector field, in particular its real and imaginary part are $J$-preserving vector fields (cf. Proposition \autoref{prop:holo_vec_field_equiv_J_pre_vec_field}) implying that the differential of their flows commute with $J$. For PHHSs, this does not need to be the case anymore: Neither $X^{\Omega_R}_{\mH_R}$ nor $J(X^{\Omega_R}_{\mH_R})$ are required to be $J$-preserving! In fact, we study an example of a proper PHHS in \autoref{sec:construction_and_deformation}, where $X^{\Omega_R}_{\mH_R}$ is $J$-preserving, but $J(X^{\Omega_R}_{\mH_R})$ is not.
\end{remark}

As for HHSs, the maximal trajectories of a PHHS, given an initial value, do not need to be unique, however, we can still pseudo-holomorphically foliate energy hypersurfaces $\mH^{-1}(E)$ of a PHHS:

\begin{proposition}[Pseudo-holomorphic foliation of a regular hypersurface]\label{prop:pseudo-holo_foli}
 Let $(X,J;\Omega_R,\mH_R)$ be a PHHS with Hamiltonian vector field\linebreak $X_\mH = 1/2 (X^{\Omega_R}_{\mH_R} - i J(X^{\Omega_R}_{\mH_R}))$ and regular\footnote{As before, a PHHS $(X,J;\Omega_R,\mH_R)$ is regular at the energy $E\in\mathbb{C}$ if $d\mH$ or, equivalently, $d\mH_R$ does not vanish on $\mH^{-1}(E)$.} value $E$ of $\mH$. Then, the energy hypersurface $\mH^{-1}(E)$ admits a pseudo-holomorphic foliation. The leaf $L_{x_0}$ of this foliation through a point $x_0\in\mH^{-1}(E)$ is given by:
 \begin{align*}
  L_{x_0}\coloneqq \{y\in X\mid &y = \varphi^{X^{\Omega_R}_{\mH_R}}_{t_1}\circ\varphi^{J(X^{\Omega_R}_{\mH_R})}_{s_1}\circ\ldots\circ\varphi^{X^{\Omega_R}_{\mH_R}}_{t_n}\circ\varphi^{J(X^{\Omega_R}_{\mH_R})}_{s_n} (x_0);\\
  &t_1,\ldots,t_n, s_1,\ldots, s_n\in\mathbb{R};\ n\in\mathbb{N}\},
 \end{align*}
 where $\varphi^{X^{\Omega_R}_{\mH_R}}_{t_j}$ and $\varphi^{J(X^{\Omega_R}_{\mH_R})}_{s_j}$ are the flows of $X^{\Omega_R}_{\mH_R}$ and $J(X^{\Omega_R}_{\mH_R})$ for time $t_j$ and $s_j$, respectively. Every pseudo-holomorphic trajectory of $(X,J;\Omega_R,\mH_R)$ with energy $E$ is completely contained in one such leaf.
\end{proposition}

\begin{proof}
 First, we need to clarify the notion of a pseudo-holomorphic foliation. For this, recall the definition of a holomorphic foliation (cf. Definition \autoref{def:holo_foli}). Clearly, we cannot directly transfer Definition \autoref{def:holo_foli} to the non-integrable case, since generic almost complex manifolds do not admit holomorphic charts. Therefore, we call $\{L_{x_0}\}_{x_0\in I}$ a pseudo-holomorphic foliation of an almost complex manifold $(X,J)$ if $\{L_{x_0}\}_{x_0\in I}$ is a foliation of $X$ and the tangent spaces of the leaves $L_{x_0}$ are closed under the action of $J$. We can now prove the last proposition in the same way as Proposition \autoref{prop:holo_foli} by applying the Frobenius theorem to the vector fields $X^{\Omega_R}_{\mH_R}$ and $J(X^{\Omega_R}_{\mH_R})$. 
\end{proof}

Similarly to HHSs, we can also define the notion of geometric trajectories for\linebreak PHHSs. We simply copy Definition \autoref{def:geo_traj} and replace the term ``holomorphic'' with ``pseudo-holomorphic''. All results we found in \autoref{sec:HHS} for geometric trajectories of HHSs still hold in the pseudo-holomorphic case. In particular, Proposition \autoref{prop:geo_traj} is still true for PHHSs. The proof is essentially the same as in the holomorphic case. However, the vector field $Y_\mH$ on the Riemann surface $\Sigma$ is a priori only a smooth section of $T^{(1,0)}\Sigma$. $Y_\mH$ becomes a holomorphic vector field on $\Sigma$ by noting that the real and imaginary part of the Hamiltonian vector field $X_\mH$ commute by construction. Thus, the real and imaginary part of $Y_\mH$ also commute, as the push-forward of $\gamma$ is a Lie algebra homomorphism, i.e., $\gamma_\ast [V,W] = [\gamma_\ast V, \gamma_\ast W]$ for vector fields\footnote{Precisely speaking, this is not correct, since $\gamma:\Sigma\to X$ is only an immersion and not a diffeomorphism, hence, the push-forward of $\gamma$ is not well-defined. Nevertheless, the argument still holds if we consider $\gamma_\ast V$ and $\gamma_\ast W$ to be sections of the pull-back bundle $\gamma^\ast TX$ and adjust the definition of the Lie bracket accordingly.} $V$ and $W$ on $\Sigma$. Now note that, for a Riemann surface $\Sigma$, the real and imaginary part of a smooth section $V$ of $T^{(1,0)}\Sigma$ commute if and only if $V$ is a holomorphic vector field on $\Sigma$. This is easily verified in holomorphic charts of $\Sigma$. This shows the holomorphicity of $Y_\mH$. The proofs for the remaining results regarding geometric trajectories work as in the holomorphic case after adjusting the language where need be.\\
Before we conclude this subsection, we want to formulate action functionals and principles for pseudo-holomorphic trajectories. By construction of PHHSs, this can be done in the same way as in \autoref{sec:HHS}. First, we note that a PHHS $(X,J;\Omega_R,\mH_R)$ decomposes into multiple RHSs. In contrast to HHSs, we only obtain two RHSs this time, namely $(X,\Omega_R,\mH_R)$ and $(X,\Omega_R,\mH_I)$, since $\Omega_I$ is, in general, not closed. However, this suffices to find action functionals for pseudo-holomorphic trajectories, as only two of the four underlying RHSs of a HHS are subject to different dynamics. If the PHHS $(X,J;\Omega_R,\mH_R)$ is exact, i.e., $\Omega_R = d\Lambda_R$, the two RHSs $(X,\Omega_R,\mH_R)$ and $(X,\Omega_R,\mH_I)$ are also exact and possess themselves action functionals. As in the case of HHSs, we can now average these action functionals over the remaining time variable and take suitable linear combinations afterwards to find the following action functional for pseudo-holomorphic trajectories:

\pagebreak

\begin{lemma}[Action principle for pseudo-holomorphic trajectories]\label{lem:pseudo-holo_action_prin_para}
 Let $(X,J;\Omega_R = d\Lambda_R,\mH_R)$ be an exact PHHS. For $\alpha\in\mathbb{R}\backslash\{n\cdot\pi\mid n\in\mathbb{Z}\}$, let\linebreak $P_\alpha\coloneqq [t_1,t_2] + e^{i\alpha}[r_1,r_2]\subset\mathbb{C}$ be a parallelogram in the complex plane with real numbers $t_1< t_2$ and $r_1< r_2$. Denote the space of smooth maps from $P_\alpha$ to $X$ by $\mathcal{P}_{P_\alpha}$ and define the action functional $\mathcal{A}^{P_\alpha}_\mH:\mathcal{P}_{P_\alpha}\to\mathbb{C}$ by:
 \begin{align*}
  \mathcal{A}^{P_\alpha}_\mH[\gamma]&\coloneqq \iint\limits_{P_\alpha}\left[\Lambda_{R,\gamma (t+is)}\left(2\frac{\partial\gamma}{\partial z}(t+is)\right) - \mH\circ\gamma (t+is)\right] dt\wedge ds\ \text{with}\\
  \frac{\partial \gamma}{\partial z}&\coloneqq \frac{1}{2}\left(\frac{\partial \gamma}{\partial t} - i\frac{\partial \gamma}{\partial s}\right)\quad\forall \gamma\in\mathcal{P}_{P_\alpha}.
 \end{align*}
 Now let $\gamma\in\mathcal{P}_{P_\alpha}$ be a smooth map from $P_\alpha$ to $X$. Then, $\gamma$ is a pseudo-holomorphic trajectory of $(X,J;\Omega_R,\mH_R)$ iff $\gamma$ is a ``critical point''\footnote{``Critical point'' means that only those variations are allowed which keep $\gamma$ fixed on the boundary $\partial P_\alpha$.} of $\mathcal{A}^{P_\alpha}_\mH$.
\end{lemma}

\begin{proof}
 For the proof of Lemma \autoref{lem:pseudo-holo_action_prin_para}, repeat the steps from \autoref{sec:HHS} for PHHSs, in particular the proof of Proposition \autoref{prop:holo_action_prin_para}.
\end{proof}

As before, if we wish to view pseudo-holomorphic trajectories as actual critical points of some functional, we can achieve this by either mapping the boundary $\partial P_\alpha$ to an exact Lagrangian submanifold of $(X,\Omega_R)$ or by imposing periodicity on the curves $\gamma$.

\subsection*{Relation between HHSs and PHHSs}

At this point, it is not clear how PHHSs relate to HHSs and why PHHSs are a ``reasonable'' generalization of HHSs with regard to the integrability of $J$. In particular, we do not know yet why the notion of PHHSs introduced in the first subsection should coincide with the notion of HHSs when we restore the integrability of $J$. A priori, there is no reason why the $2$-form $\Omega$ associated with a PHHS $(X,J;\Omega_R,\mH_R)$ should be holomorphic or even closed for integrable $J$. The following theorem guarantees that this is indeed the case:

\begin{theorem}[Relation between HSMs and PHSMs]\label{thm:rel_HSM_PHSM}
 Let $(X,J;\Omega_R)$ be a PHSM with $2$-forms $\Omega_I\coloneqq -\Omega_R (J\cdot,\cdot)$ and $\Omega\coloneqq \Omega_R + i\Omega_I$. Further, let $x_0\in X$ be any point. Then, $J$ is integrable near $x_0$ if and only if $d\Omega_I$ vanishes near $x_0$. Moreover, the following statements are equivalent:
 \begin{enumerate}
  \item $(X,\Omega)$ is a HSM with complex structure $J$.
  \item $\Omega_I$ is closed, i.e., $d\Omega_I = 0$.
  \item $J$ is integrable.
 \end{enumerate}
\end{theorem}

\begin{proof}
 We only show the equivalence of Statement (i), (ii), and (iii). The first part of Theorem \autoref{thm:rel_HSM_PHSM} is then just a local version of the equivalence. Direction ``(i)$\Rightarrow$(ii)'' is trivially true by definition of a HSM. Implication ``(ii)$\Rightarrow$(iii)'' is due to Verbitsky (confer Theorem 3.5 in \cite{verbitsky2013} and Proposition 2.12 in \cite{bogomolov2020}). For completeness' sake, we include the proof here. Let $(X,J;\Omega_R)$ be a PHSM with $2$-forms $\Omega_I\coloneqq -\Omega_R (J\cdot,\cdot)$ and $\Omega\coloneqq \Omega_R + i\Omega_I$. Further, assume $d\Omega_I = 0$. We want to show that $J$ is integrable. By the Newlander-Nirenberg theorem, $J$ is integrable if and only if $J$ has no torsion, i.e., its Nijenhuis tensor vanishes. Now we apply Theorem 2.8 in Chapter IX of \cite{kobayashi1969}. Thus, $J$ is integrable if and only if the space of smooth sections of $T^{(0,1)}X$ is closed under the commutator $[\cdot,\cdot]$. From the first subsection, we know that $\Omega$ is non-degenerate on $T^{(1,0)}X$, but vanishes on $T^{(0,1)}X$. Hence, a complex vector field $V$ on $X$ is a smooth section of $T^{(0,1)}$ if and only if $\iota_V\Omega = 0$. Therefore, $J$ is integrable if and only if for every pair of two complex vector fields $V$ and $W$ on $X$ satisfying $\iota_V\Omega = \iota_W\Omega = 0$ one has $\iota_{[V,W]}\Omega = 0$. Now let $V$ and $W$ be two complex vector fields on $X$ with $\iota_V\Omega = \iota_W\Omega = 0$. Recall that the interior product $\iota$ applied to forms fulfills the relation (cf. Proposition 3.10 in Chapter I of \cite{kobayashi1963}):
 \begin{gather*}
  \iota_{[V,W]} = [L_V, \iota_W],
 \end{gather*}
 where $L_V$ is the Lie derivative of $V$. We can calculate $L_V\Omega$ by using Cartan's magic formula, $\iota_V\Omega = 0$, and $d\Omega = 0$:
 \begin{gather*}
  L_V\Omega = d\iota_V\Omega + \iota_Vd\Omega = 0.
 \end{gather*}
 In total, we obtain using $\iota_W\Omega = 0$:
 \begin{gather*}
  \iota_{[V,W]}\Omega = [L_V,\iota_W]\Omega = L_V(\iota_W\Omega) - \iota_W(L_V\Omega) = 0
 \end{gather*}
 proving the integrability of $J$.\\
 The remaining direction ``(iii)$\Rightarrow$(i)'' can be proven as follows: Let $(X,J;\Omega_R)$ be a PHSM with $2$-forms $\Omega_I\coloneqq -\Omega_R (J\cdot,\cdot)$ and $\Omega\coloneqq \Omega_R + i\Omega_I$. Further assume that $J$ is integrable. Then, $X$ is a complex manifold with complex structure $J$. We need to show that $\Omega$ is a closed, holomorphic $2$-form which is non-degenerate on $T^{(1,0)}X$. By Remark \autoref{rem:omega_II}, $\Omega$ is non-degenerate on $T^{(1,0)}X$ and of type $(2,0)$. Hence, in a holomorphic chart $\phi = (z_1,\ldots, z_{2m}):U\to V\subset\mathbb{C}^{2m}$ of $X$, $\Omega$ can be written as:
 \begin{gather*}
  \Omega\vert_U = \sum\limits^{2m}_{i,j = 1} \Omega_{ij} dz_i\wedge dz_j,
 \end{gather*}
 where the coefficients $\Omega_{ij} = -\Omega_{ji}:U\to\mathbb{C}\cong\mathbb{R}^2$ are smooth functions on $U$. From $d\Omega_R = 0$, we deduce:
 \begin{align*}
  0 &= 2d\Omega_R\vert_U = (\partial + \bar{\partial})(\Omega + \overline{\Omega})\vert_U\\
  &= \sum\limits^{2m}_{i,j,k = 1}\left(\frac{\partial\Omega_{ij}}{\partial z_k} dz_k\wedge dz_i\wedge dz_j + \frac{\partial\Omega_{ij}}{\partial \bar{z}_k} d\bar{z}_k\wedge dz_i\wedge dz_j\right.\\
  &\qquad\qquad + \left.\frac{\partial\overline{\Omega}_{ij}}{\partial z_k} dz_k\wedge d\bar{z}_i\wedge d\bar{z}_j + \frac{\partial\overline{\Omega}_{ij}}{\partial \bar{z}_k} d\bar{z}_k\wedge d\bar{z}_i\wedge d\bar{z}_j\right).
 \end{align*}
 This equation implies:
 \begin{gather*}
  \frac{\partial\Omega_{ij}}{\partial \bar{z}_k} = 0\quad\forall i,j,k\in\{1,\ldots, 2m\}.
 \end{gather*}
 Thus, the coefficients $\Omega_{ij}$ are holomorphic functions on $U$. As the last argument can be repeated for any holomorphic chart of $X$, the form $\Omega$ itself is holomorphic. Therefore, its exterior derivative $d\Omega$ is also a holomorphic form. In particular, $d\Omega$ satisfies:
 \begin{gather*}
  d\Omega (J\cdot,\cdot,\cdot) = i\cdot d\Omega.
 \end{gather*}
 As the exterior derivative is $\mathbb{C}$-linear, the decomposition of $d\Omega$ into real and imaginary part amounts to $d\Omega = d\Omega_R + id\Omega_I$. Combining this decomposition with the previous equation gives us:
 \begin{gather*}
  d\Omega_I = -d\Omega_R (J\cdot,\cdot,\cdot) = 0,
 \end{gather*}
 where we have used the closedness of $\Omega_R$ again. This shows that $\Omega$ has the desired properties concluding the proof.
\end{proof}

\begin{remark}[Closedness of $\Omega_R$]\label{rem:closedness_of_Omega_R}
 Note that the closedness of $\Omega_R$ is crucial for Theorem \autoref{thm:rel_HSM_PHSM}: If $(X,\Omega = \Omega_R + i\Omega_I)$ is a HSM and $f:X\to\mathbb{R}_+$ is a positive smooth function on $X$, then $f\cdot\Omega_R$ is still a non-degenerate $2$-form on $X$ which is anticompatible with the integrable complex structure $J$, however, neither $f\cdot \Omega_R$ nor $f\cdot \Omega_I$ are necessarily closed. In fact, $f\cdot \Omega$ is, in general, not even holomorphic.
\end{remark}

Of course, we can also formulate Theorem \autoref{thm:rel_HSM_PHSM} for Hamiltonian systems:

\begin{corollary}[Relation between HHSs and PHHSs]\label{cor:rel_HHS_PHHS}
 Let $(X,J;\Omega_R, \mH_R)$ be a PHHS with $2$-forms $\Omega_I\coloneqq -\Omega_R (J\cdot,\cdot)$ and $\Omega\coloneqq \Omega_R + i\Omega_I$ as well as a function $\mH\coloneqq \mH_R + i\mH_I$, where $\mH_I$ is any primitive of the $1$-form $\Omega_R (J(X^{\Omega_R}_{\mH_R}),\cdot)$. Then, the following statements are equivalent:
 \begin{enumerate}
  \item $(X,\Omega, \mH)$ is a HHS with complex structure $J$.
  \item $\Omega_I$ is closed, $d\Omega_I = 0$.
  \item $J$ is integrable.
 \end{enumerate}
\end{corollary}

\begin{proof}
 Corollary \autoref{cor:rel_HHS_PHHS} is a direct consequence of Theorem \autoref{thm:rel_HSM_PHSM} and Remark \autoref{rem:H_pseudo-holo}.
\end{proof}

We can interpret Theorem \autoref{thm:rel_HSM_PHSM} and Corollary \autoref{cor:rel_HHS_PHHS} as follows: The integrability of the almost complex structure $J$ of a PHSM $(X,J;\Omega_R)$ or a PHHS $(X,J;\Omega_R,\mH_R)$ is completely measured by the closedness of the imaginary part $\Omega_I\coloneqq -\Omega_R (J\cdot,\cdot)$ and vice versa. Moreover, these quantities are the only local invariants of a PHSM or a PHHS: We know by Darboux's theorem for HSMs (cf. Theorem \autoref{thm:holo_Darboux}) that any HSM can locally be brought into standard form. Similarly, we will see in \autoref{sec:construction_and_deformation} that also (regular) HHSs can locally be brought into standard form. The existence of coordinates in which some geometrical object assumes a standard form implies that said geometrical object exhibits no local invariant. In this sense, Theorem \autoref{thm:rel_HSM_PHSM} and Corollary \autoref{cor:rel_HHS_PHHS} state that the Nijenhuis tensor $N_J$ of $J$ or, equivalently, the exterior derivative $d\Omega_I$ are the only local invariants of PHSMs and (regular) PHHSs. For general PHHSs, the Nijenhuis tensor and the behavior of the Hamiltonian near singular points are the only local invariants.

\newpage
\section[Construction of Proper PHHSs and Deformation of HHSs]{Construction of Proper PHHSs and Deformation of HHSs\sectionmark{Construction and Deformation}}
\sectionmark{Construction and Deformation}
\label{sec:construction_and_deformation}
\autoref{sec:construction_and_deformation} is divided into two subsections. The aim of the first subsection is to provide examples of proper PHHSs. In fact, we present a general method for constructing PHHSs out of HHSs (cf. Proposition \autoref{prop:constructing_PHHS_out_of_HHS}), which allows us, for instance, to equip the holomorphic cotangent bundle of a complex manifold with a (non-canonical) PHHS structure. In the second subsection, we study the ``size'' of the set of proper PHHSs within the set of all PHHSs. The main result of this investigation is that proper PHHSs are generic. To prove this result, we deform HHSs by proper PHHSs (cf. Theorem \autoref{thm:generic}).

\subsection*{Constructing proper PHHSs out of HHSs}

The goal of this subsection is to find examples of proper PHHSs. The basic idea is to start with a HHS $(X,\Omega,\mH)$ and turn it into a PHHS by twisting its complex structure $J$ with a $\Omega_R$-compatible $(1,1)$-tensor $A$, usually an almost complex structure. To elaborate on this idea, we first explain the notion of compatibility:

\begin{definition}[$\omega$-compatible $(1,1)$-tensor]
 Let $(M,\omega)$ be a symplectic manifold. We call a $(1,1)$-tensor $A\in\Gamma\End (TM)$ \textbf{compatible} with $\omega$ if $\omega_p (A_pv, A_pw) = \omega_p(v,w)$ holds for all vectors $v,w\in T_pM$ and every point $p\in M$.
\end{definition}

\begin{remark}
 We add two facts regarding $\omega$-compatible $(1,1)$-tensors:
 \begin{enumerate}
  \item Since $\omega$ is non-degenerate, every $\omega$-compatible $(1,1)$-tensor $A$ is invertible, i.e., $A^{-1}\in\Gamma\End (TM)$ exists with $AA^{-1} = A^{-1}A = \mathds{1}$. It is easy to see that $A^{-1}$ is also compatible with $\omega$.
  \item In the literature on symplectic geometry, the term ``compatibility'' is reserved for almost complex structures and used in a slightly different manner: Given an almost complex structure $I$ on $M$, one usually says that $I$ is $\omega$-compatible if $g\coloneqq \omega (\cdot,I\cdot)$ is a Riemannian metric which implies $\omega (I\cdot, I\cdot) = \omega$. The connection between this and our notion is that a $\omega$-compatible $(1,1)$-tensor $A$ satisfies $A^2 = -\mathds{1}$ if and only if $\omega (\cdot, A\cdot)$ is a semi-Riemannian metric. Still, the commonly used notion of ``compatibility'' is stronger than our notion, since $\omega (\cdot, A\cdot)$ does not have to be positive definite in our case.
 \end{enumerate}
\end{remark}

Now pick a HHS $(X,\Omega = \Omega_R + i\Omega_I,\mH = \mH_R + i\mH_I)$ with complex structure $J$ and choose a $\Omega_R$-compatible $(1,1)$-tensor $A$. We consider the $(1,1)$-tensor $J_A\coloneqq AJA^{-1}$. Clearly, $J_A$ is an almost complex structure:
\begin{gather*}
 J^2_A = (AJA^{-1})^2 = AJA^{-1}AJA^{-1} = AJ^2A^{-1} = -AA^{-1} = -\mathds{1}.
\end{gather*}
Moreover, $J_A$ is still $\Omega_R$-anticompatible:
\begin{align*}
 \Omega_R (J_A\cdot, J_A\cdot) &= \Omega_R (AJA^{-1}\cdot, AJA^{-1}\cdot) = \Omega_R (JA^{-1}\cdot, JA^{-1}\cdot) = -\Omega_R(A^{-1}\cdot,A^{-1}\cdot)\\
 &= -\Omega_R.
\end{align*}
Thus, $(X,J_A;\Omega_R)$ is a PHSM. If $A$ was chosen such that $\Omega_R (J_A (X^{\Omega_R}_{\mH_R}),\cdot)$ is exact, then $(X,J_A;\Omega_R,\mH_R)$ is even a PHHS. The PHHS constructed this way is generally proper. To check that, it suffices by Corollary \autoref{cor:rel_HHS_PHHS} to compute the exterior derivative of $\Omega^A_I\coloneqq -\Omega_R (J_A\cdot,\cdot)$. The following proposition collects our findings:

\begin{proposition}[Constructing PHHSs]\label{prop:constructing_PHHS_out_of_HHS}
 Let $(X,\Omega,\mH)$ be a HHS with complex structure $J$ and decompositions $\Omega = \Omega_R + i\Omega_I$ and $\mH = \mH_R + i\mH_I$. Further, let $A$ be a $\Omega_R$-compatible $(1,1)$-tensor. Then, $(X,J_A;\Omega_R)$ is a PHSM, where $J_A\coloneqq AJA^{-1}$. If, additionally, the $1$-form $\Omega_R (J_A (X^{\Omega_R}_{\mH_R}),\cdot)$ is exact, then $(X,J_A;\Omega_R,\mH_R)$ is even a PHHS. $(X,J_A;\Omega_R)$ and $(X,J_A;\Omega_R,\mH_R)$ are proper if and only if $d\Omega^A_I\neq 0$, where $\Omega^A_I\coloneqq -\Omega_R (J_A\cdot,\cdot)$.
\end{proposition}

\begin{remark}\label{rem:constructing_PHHS_out_of_HHS}
 If the $\Omega_R$-compatible tensor $A$ is an almost complex structure or, equivalently, $g\coloneqq\Omega_R (\cdot, A\cdot)$ is a semi-Riemannian metric, we emphasize this circumstance by writing $I_g\coloneqq A$, $J_g\coloneqq -J_A = I_gJI_g$, and $\Omega^g_I\coloneqq -\Omega_R (J_g\cdot,\cdot)$. As before, $(X,J_g;\Omega_R)$ is a PHSM, $(X,J_g;\Omega_R,\mH_R)$ is a PHHS if $\Omega_R (J_g (X^{\Omega_R}_{\mH_R}),\cdot)$ is exact, and both are proper iff $d\Omega^g_I\neq 0$.
\end{remark}

In practice, one can find $\Omega_R$-compatible tensors $A$ by requiring $A = I_g$ to be an almost complex structure and picking a suitable semi-Riemannian metric $g$. Given any HHS $(X,\Omega,\mH)$, however, an arbitrary $\Omega_R$-compatible almost complex structure $I_g$ is usually not compatible with $\mH_R$ in the sense explained above. Most often, it is simpler to first pick an almost complex structure $I_g$ and afterwards pick a suitable real function $\mH_R$. Since $X$ is contractible in most examples we wish to study, e.g. local considerations and $X = \mathbb{C}^{2m}$, one can often find suitable $\mH_R$ by solving the differential equation $d[\Omega_R (J_g (X^{\Omega_R}_{\mH_R}),\cdot)] = 0$. To illustrate the construction, let us consider the simplest non-trivial example:

\begin{example}[PHHS on $X = \mathbb{C}^2$]\label{ex:constructing_PHHS}\normalfont
 Let $(X,\Omega,\mH)$ be the HHS with $X = \mathbb{C}^2$, $\Omega = dz_2\wedge dz_1$, where $(z_1,z_2)\in\mathbb{C}^2$, and $\mH (z_1,z_2) = i\cdot z_1$. With the decomposition $z_j = x_j + i\cdot y_j$, we obtain:
 \begin{gather*}
  \Omega_R = dx_2\wedge dx_1 - dy_2\wedge dy_1,\quad \mH_R = -y_1.
 \end{gather*}
 The complex structure $J$ of the HHS $(X,\Omega,\mH)$ is given by $i$:
 \begin{gather*}
  J (\partial_{x_j})\coloneqq \partial_{y_j},\quad J (\partial_{y_j})\coloneqq -\partial_{x_j}.
 \end{gather*}
 Now pick the following semi-Riemannian metric $g$ on $\mathbb{C}^2$:
 \begin{gather*}
  g(\pa{x_1},\pa{x_1}) = g(\pa{x_2},\pa{x_2})^{-1} = f,\quad g(\pa{y_1},\pa{y_1}) = g(\pa{y_2},\pa{y_2})^{-1} = h,\\
  g(\pa{x_1},\pa{x_2}) = g(\pa{y_1},\pa{y_2}) = g(\pa{x_i},\pa{y_j}) = 0,
 \end{gather*}
 where $f,h:\mathbb{C}^2\to\mathbb{R}$ are smooth, nowhere-vanishing functions. The corresponding almost complex structure $I_g$ is given by:
 \begin{gather*}
  I_g (\pa{x_1}) = -f\pa{x_2},\ I_g (\pa{x_2}) = f^{-1}\pa{x_1},\quad I_g (\pa{y_1}) = h\pa{y_2},\ I_g (\pa{y_2}) = -h^{-1}\pa{y_1}.
 \end{gather*}
 We now employ the notation $r\coloneqq f/h$. Computing $J_g = I_gJI_g$ in the standard basis yields:
 \begin{gather*}
  J_g (\pa{x_1}) = r\pa{y_1},\ J_g (\pa{x_2}) = r^{-1}\pa{y_2},\quad J_g (\pa{y_1}) = -r^{-1}\pa{x_1},\ J_g (\pa{y_2}) = -r\pa{x_2}.
 \end{gather*}
 We find for the induced $2$-form $\Omega^g_I= -\Omega_R (J_g\cdot,\cdot)$:
 \begin{gather*}
  \Omega^g_I = r^{-1}dx_2\wedge dy_1 + rdy_2\wedge dx_1.
 \end{gather*}
 Hence, the exterior derivative of $\Omega^g_I$ amounts to:
 \begin{gather*}
  d\Omega^g_I 
  = dr\wedge\left( dy_2\wedge dx_1 - r^{-2}dx_2\wedge dy_1\right).
 \end{gather*}
 In general, the exterior derivative $d\Omega^g_I$ does not vanish. For instance, we can set\linebreak $f (z_1, z_2) = 1$ and $h(z_1, z_2) = e^{x_1}$ resulting in $r (z_1, z_2) = e^{-x_1}$. This choice yields:
 \begin{gather*}
  d\Omega^g_I = e^{x_1} dx_1\wedge dx_2\wedge dy_1.
 \end{gather*}
 This form does not vanish at any point of $X$. Thus, $(X,J_g;\Omega_R)$ is a proper PHSM for the presented choice of $f$ and $h$.\\
 Next, we check whether the $1$-form $\Omega_R (J_g (X^{\Omega_R}_{\mH_R}),\cdot)$ is exact. Since $X = \mathbb{C}^2$ is contractible, it suffices to check whether $\Omega_R (J_g (X^{\Omega_R}_{\mH_R}),\cdot)$ is closed. For the sake of generality, we first perform the computation for any smooth real function $H:\mathbb{C}^2\to\mathbb{R}$ and afterwards insert $H = \mH_R = -y_1$. We start by determining the Hamiltonian vector field $X_H\equiv X^{\Omega_R}_H$. It can be written as:
 \begin{gather*}
  X_H = (\pa{x_2}H)\cdot\pa{x_1} - (\pa{x_1}H)\cdot\pa{x_2} - (\pa{y_2}H)\cdot\pa{y_1} + (\pa{y_1}H)\cdot\pa{y_2}.
 \end{gather*}
 Applying $J_g$ to $X_H$ yields:
 \begin{gather*}
  J_g (X_H) = (r^{-1}\pa{y_2}H)\cdot \pa{x_1} - (r\pa{y_1}H)\cdot \pa{x_2} + (r\pa{x_2}H)\cdot \pa{y_1} - (r^{-1}\pa{x_1}H)\cdot\pa{y_2}.
 \end{gather*}
 Contracting $J_g (X_H)$ with $\Omega_R$ gives:
 \begin{align*}
  \Omega_R (J_g (X_H),\cdot) &= -(r\pa{y_1}H)\cdot dx_1 - (r^{-1}\pa{y_2}H)\cdot dx_2\\
  &\phantom{=}\ + (r^{-1}\pa{x_1}H)\cdot dy_1 + (r\pa{x_2}H)\cdot dy_2\\
  &\equiv w_{x_1} dx_1 + w_{x_2} dx_2 + w_{y_1} dy_1 + w_{y_2} dy_2.
 \end{align*}
 For $\Omega_R (J_g (X_H),\cdot)$ to be closed, the coefficients $w$ need to satisfy the following conditions:
 \begin{gather*}
  \pa{x_i} w_{x_j} = \pa{x_j} w_{x_i},\quad \pa{y_i} w_{y_j} = \pa{y_j} w_{y_i},\quad \pa{x_i} w_{y_j} = \pa{y_j} w_{x_i}.
 \end{gather*}
 Let us now set $H = \mH_R = -y_1$. Then, all coefficients except for $w_{x_1} = r$ vanish. In particular, $\Omega_R (J_g (X_H),\cdot)$ is exact if and only if $r$ only depends on $x_1$. In this case, the primitive $\mH^g_I$ of $\Omega_R (J_g (X_H),\cdot)$ is given by $R$, where $R$ only depends on $x_1$ and satisfies $\pa{x_1}R = r$. For instance, we can again set $f (z_1, z_2) = 1$ and $h(z_1, z_2) = e^{x_1}$ leading to $r (z_1, z_2) = e^{-x_1}$. Then, $(X,J_g;\Omega_R,\mH_R)$ is a proper PHHS for this choice, where $\mH^g_I = -e^{-x_1} + c$, $c\in\R$, is the imaginary part of the Hamiltonian.
\end{example}

Of course, Proposition \autoref{prop:constructing_PHHS_out_of_HHS} is not the only way to obtain a PHSM out of a HSM. For instance, the simple structure of the example $(X = \mathbb{C}^2, \Omega = dz_2\wedge dz_1)$ allows us to transform the standard complex structure $J$ into a new $\Omega_R$-anticompatible almost complex structure $J_\varphi$ by rotation of the axes:
\begin{alignat*}{2}
 J_\varphi (\pa{x_1})&\coloneqq \cos (\varphi)\pa{y_1} - \sin (\varphi)\pa{y_2},\quad &&J_\varphi (\pa{x_2})\coloneqq \sin (\varphi)\pa{y_1} + \cos (\varphi)\pa{y_2},\\
 J_\varphi (\pa{y_1})&\coloneqq -\cos (\varphi)\pa{x_1} - \sin (\varphi)\pa{x_2},\quad &&J_\varphi (\pa{y_2})\coloneqq \sin (\varphi)\pa{x_1} - \cos (\varphi)\pa{x_2},
\end{alignat*}
where $\varphi:\mathbb{C}^2\to\mathbb{R}$ is any smooth function. In general, the resulting PHSM $(X, J_\varphi; \Omega_R)$ is proper.\\
The reason why we focus our attention on the construction given in Proposition \autoref{prop:constructing_PHHS_out_of_HHS} is that it possesses an interesting connection to holomorphic Kähler and Hyperkähler manifolds. The various flavors of Kähler manifolds are discussed in \autoref{app:kaehler}. The quick rundown is that a Kähler manifold $(X,g,I)$ is a Riemannian manifold $(X,g)$ with complex structure $I$ such that $\omega\coloneqq g(I\cdot,\cdot)$ is a symplectic form, that a holomorphic Kähler manifold is a complexification of a Kähler manifold, and that a Hyperkähler manifold $(X,g,I,J,K)$ is a quaternionic (i.e. $I^2 = J^2 = K^2 = IJK = -\mathds{1}$) manifold such that $(X,g,I)$, $(X,g,J)$, and $(X,g,K)$ are Kähler. In \autoref{app:kaehler}, we discover that holomorphic Kähler and Hyperkähler manifolds are quite similar from a symplectic viewpoint (cf. Lemma \autoref{lem:hyper_vs_holo}): Both are described by the data $(X,\omega,I,J)$, where $\omega$ is a symplectic form on the manifold $X$ and $I,J$ are complex structures on $X$ satisfying $\omega (I\cdot,I\cdot) = \omega$ and $\omega(J\cdot,J\cdot) = -\omega$. In particular, $\Omega^J\coloneqq \omega -i\omega (J\cdot,\cdot)$ is a holomorphic\footnote{The superscript $J$ of $\Omega^J$ indicates that $\Omega^J$ is holomorphic with respect to $J$.} symplectic form in both situations. The only difference between holomorphic Kähler and Hyperkähler manifolds is that $I$ and $J$ commute ($IJ = JI$) in the holomorphic Kähler case, while they anticommute ($IJ = -JI$) for Hyperkähler manifolds.\\
We now realize that the symplectic picture of holomorphic Kähler/Hyperkähler manifolds is quite similar to the setup of Proposition \autoref{prop:constructing_PHHS_out_of_HHS}, at least in the case $A = I_g$ (cf. Remark \autoref{rem:constructing_PHHS_out_of_HHS}): We also have a symplectic form $\Omega_R$, a complex structure $J$ anticompatible with $\Omega_R$, and an almost complex structure $I_g$ compatible with $\Omega_R$. Solely the non-integrability of $I_g$ and the missing commutation relation between $I_g$ and $J$ differ from the holomorphic Kähler/Hyperkähler case. Clearly, it does not make sense to impose a commutation relation for the purposes of Proposition \autoref{prop:constructing_PHHS_out_of_HHS}. Indeed, if we imposed a commutation relation, the resulting structure $J_g = I_gJI_g = \mp J$ would only differ from $J$ by sign. In particular, we would not obtain a proper PHHS this way.\\
Still, we might think of Proposition \autoref{prop:constructing_PHHS_out_of_HHS} as a generalization or deformation of holomorphic Kähler/Hyperkähler structures in some sense. To explore this statement further, consider Example \autoref{ex:constructing_PHHS} again. If we choose $g$ to be the standard Euclidean metric $\delta$ on $\mathbb{C}^2\cong\mathbb{R}^4$, i.e., $f \equiv h\equiv 1$, then $I_\delta$ anticommutes with the standard complex structure $J$, thus, $J_\delta = J$. In fact, $(\C^2,\delta, I_\delta, J, K\coloneqq I_\delta J)$ constitutes a Hyperkähler manifold and $dz_2\wedge dz_1$ comes from this Hyperkähler structure:
\begin{gather*}
 dz_2\wedge dz_1 = \delta (I_\delta\cdot,\cdot) - i\delta (K\cdot,\cdot).
\end{gather*}
Now recall that the metrics $g$ in Example \autoref{ex:constructing_PHHS} are parameterized by two smooth nowhere-vanishing functions $f$ and $h$ on $\C^2$. The space of pairs $(f,h)$ has four connected components which are isomorphic via $(f,h)\mapsto (-f,h)$ and\linebreak $(f,h)\mapsto (f,-h)$. Furthermore, each connected component is contractible. Thus, we can reach any pair $(f,h)$ in this space by a path starting at the point $(1,1)$, whose corresponding metric is $\delta$, and a discrete transformation. Next, consider such a path in this space, i.e., a smooth $1$-parameter family of nowhere-vanishing functions $f^\varepsilon$ and $h^\varepsilon$ on $\mathbb{C}^2$ such that $f^0\equiv h^0 \equiv 1$. This family induces a $1$-parameter family $J^\varepsilon$ of almost complex structures on $\mathbb{C}^2$ satisfying $J^0 = J$. By construction, $J^\varepsilon$ itself is defined via a $1$-parameter family of almost complex structures $I^\varepsilon$ satisfying $I^0 = I_\delta$. With our previous knowledge, we can interpret $I^\varepsilon$ as a deformation of the Hyperkähler manifold $(X = \mathbb{C}^2,\delta, I_\delta, J, K)$. Hence, we can say that every almost complex structure $J_g$ in Example \autoref{ex:constructing_PHHS} giving rise to a PHSM comes from such a deformation up to a discrete transformation. In general, this seems to be the best application of our construction: Pick a HSM $(X,\Omega)$ coming from a Hyperkähler manifold $(X,g,I,J,K)$ and then deform $I$ as described to find PHSMs.\\
Next, we show that Proposition \autoref{prop:constructing_PHHS_out_of_HHS} is applicable to a rather large class of complex manifolds $X$, namely the class of holomorphic cotangent bundles\linebreak $X = T^{\ast, (1,0)}Y$ of complex manifolds $Y$. Recall that the complex structure of $Y$ induces a canonical complex structure $J$ on $X$ and that $X$ as a holomorphic cotangent bundle possesses a canonical HSM structure with $2$-form $\Omega_\text{can} = \Omega_R + i\Omega_I$ (cf. \autoref{sec:HHS}). We now observe that the real part $\Omega_R$ of $\Omega_\text{can}$ can be identified with the canonical symplectic $2$-form $\omega_\text{can}$ of the real cotangent bundle $T^\ast Y$:

\begin{proposition}[$F^\ast\omega_\text{can} = \Omega_R$]\label{prop:F}
 Let $Y$ be a complex manifold. Then, the map $F:T^{\ast, (1,0)}Y\to T^\ast Y$, $F(\alpha)\coloneqq \text{\normalfont Re} (\alpha)$, where $\text{\normalfont Re} (\alpha)$ denotes the real part of $\alpha$, is a bundle isomorphism between smooth real vector bundles. In particular, $F$ is a diffeomorphism between the smooth manifolds $T^{\ast, (1,0)}Y$ and $T^\ast Y$ satisfying $F^\ast \omega_\text{can} = \Omega_R$.
\end{proposition}

\begin{proof}
 Let $\psi = (Q_1 = Q_{x,1} + iQ_{y,1},\ldots, Q_n = Q_{x,n} + iQ_{y,n})$ be a holomorphic chart of $Y$. Then, $T^{\ast, (1,0)}\psi$ is a holomorphic chart of $T^{\ast, (1,0)}Y$ and $T^{\ast}\psi$ is a smooth chart of $T^\ast Y$. Identifying $\C^{2n}$ with $\R^{4n}$ allows us to view $T^{\ast, (1,0)}\psi$ as a real chart of $T^{\ast, (1,0)}Y$ giving us the expressions:\footnote{We suppress the base point $\psi^{-1}(\tilde{Q}_{1},\ldots, \tilde{Q}_{n})$ in the second and third line.}
 \begin{align*}
  &\left(T^{\ast, (1,0)}\psi\right)^{-1} (\tilde{Q}_{x,1},\ldots, P_{y,n}) \coloneqq \sum^n_{j=1} P_j dQ_{j,\psi^{-1}(\tilde{Q}_{1},\ldots, \tilde{Q}_{n})}\\
  &\phantom{\left(T^{\ast, (1,0)}\psi\right)^{-1} (\tilde{Q}_{x,1},\ldots, P_{y,n})\ }= \sum^n_{j = 1} P_{x,j}dQ_{x,j} - P_{y,j}dQ_{y,j}\\
  &\phantom{\left(T^{\ast, (1,0)}\psi\right)^{-1} (\tilde{Q}_{x,1},\ldots, P_{y,n})=\,}+ i\sum^n_{j=1}P_{x,j}dQ_{y,j} + P_{y,j}dQ_{x,j},\\
  &\left(T^{\ast}\psi\right)^{-1} (\tilde{q}_{x,1},\ldots, p_{y,n}) \coloneqq \sum^n_{j = 1} p_{x,j}dQ_{x,j,\psi^{-1}(\tilde{q}_{1},\ldots, \tilde{q}_{n})} + p_{y,j}dQ_{y,j,\psi^{-1}(\tilde{q}_{1},\ldots, \tilde{q}_{n})}.
 \end{align*}
 Denote the components of $T^{\ast, (1,0)}\psi$ by $T^{\ast, (1,0)}\psi = (Q_{x,1},\ldots, P_{y,n})$ and the components of $T^{\ast}\psi$ by $T^{\ast}\psi = (q_{x,1},\ldots, p_{y,n})$. In these coordinates, $\Omega_\text{can}$ and $\omega_\text{can}$ become:
 \begin{align*}
  \Omega_\text{can} &= \sum^n_{j = 1} dP_j\wedge dQ_j = \sum^n_{j=1} dP_{x,j}\wedge dQ_{x,j} - dP_{y,j}\wedge dQ_{y,j}\\
  &\phantom{= \sum^n_{j = 1} dP_j\wedge dQ_j =\,}+ i\sum^n_{j=1}dP_{x,j}\wedge dQ_{y,j} + dP_{y,j}\wedge dQ_{x,j},\\
  \omega_\text{can} &= \sum^n_{j=1} dp_{x,j}\wedge dq_{x,j} + dp_{y,j}\wedge dq_{y,j}.
 \end{align*}
 Expressing $F$ in these coordinates gives:
 \begin{gather*}
  T^\ast\psi\circ F\circ \left(T^{\ast, (1,0)}\psi\right)^{-1} (\tilde{Q}_{x,j}, \tilde{Q}_{y,j}, P_{x,j}, P_{y,j}) = (\tilde{Q}_{x,j}, \tilde{Q}_{y,j}, P_{x,j}, -P_{y,j}).
 \end{gather*}
 In total, we obtain:
 \begin{gather*}
  F^\ast\omega_\text{can} = \sum^n_{j=1} dP_{x,j}\wedge dQ_{x,j} - dP_{y,j}\wedge dQ_{y,j} = \text{Re} (\Omega_\text{can}) = \Omega_R.
 \end{gather*}
\end{proof}

Let us now choose a semi-Riemannian metric $g$ on $Y$. In \autoref{app:almost_complex_structures}, we show that any metric $g$ induces an almost complex structure $J^\ast_{\nabla^{g}}$ on $T^\ast Y$ which is $\omega_\text{can}$-compatible. To summarize the construction of $J^\ast_{\nabla^{g}}$, one first splits\linebreak $TY = H\oplus V$ into horizontal and vertical bundle with the help of the Levi-Civita connection $\nabla^g$, afterwards defines a complex structure on $TY$ by $H\oplus V\to H\oplus V$, $(w_1,w_2)\mapsto (w_2,-w_1)$ and lastly identifies $TY$ with $T^\ast Y$ via $g$. Using $F$, we can transfer this almost complex structure from $T^\ast Y$ to $X$. We denote the result by $I_g\coloneqq dF^{-1}\circ J^\ast_{\nabla^{g}}\circ dF$. By Proposition \autoref{prop:F}, $I_g$ is also $\Omega_R$-compatible. Thus, $J_g\coloneqq I_gJI_g$ is $\Omega_R$-anticompatible and $(X,J_g;\Omega_R)$ is a PHSM.\\
In general, this PHSM is proper, since we have not imposed any relation between $g$ and the complex structure on $Y$. However, there are two special cases in which $J_g$ is integrable, namely if $g = h_R$ is the real part of a holomorphic metric $h = h_R + ih_I$ or if $g$ is a Kähler metric. Indeed, $I_g$ and $J$ commute for $g = h_R$ and anticommute for Kähler metrics $g$. A very similar statement is proven in \autoref{chap:duality} in a slightly different context (cf. Lemma \autoref{lem:commutation_relation}). Here, we are content with outlining the proof for the case $g = h_R$ (cf. Lemma \autoref{lem:just_one_commutation_relation}). The Kähler case then works analogously to Lemma \autoref{lem:just_one_commutation_relation} and Lemma \autoref{lem:commutation_relation}:

\begin{lemma}[$g = h_R\ \Rightarrow\ J$ and $I_g$ commute]\label{lem:just_one_commutation_relation}
 Let $Y$ be a complex manifold with holomorphic metric $h = h_R + ih_I$. Then, $I_g$ obtained from the construction above for $g = h_R$ commutes with the complex structure $J$ on $X\coloneqq T^{\ast, (1,0)}Y$. In particular, $J_g \coloneqq I_gJI_g = -J$.
\end{lemma}

\begin{proof}
 The idea of the proof is to choose coordinates in which $J$ and $I_g$ take a simple form. Let $p\in Y$ be any point. We start by choosing holomorphic coordinates $\psi = (Q_1 = Q_{x,1} + iQ_{y,1},\ldots, Q_n = Q_{x,n} + iQ_{y,n})$ of $Y$ near $p$ which, at the same time, are normal coordinates of $(Y,h_R)$ near $p$. In \autoref{app:holo_connection}, we show that such coordinates exist by considering normal coordinates of the holomorphic Levi-Civita connection $\nabla^h$. The holomorphic normal coordinates $\psi$ then give rise to holomorphic coordinates $T^{\ast, (1,0)}\psi = (Q_1,\ldots, P_n = P_{x,n} + iP_{y,n})$ of $X$. As $J$ is the complex structure of $X$, $J$ takes the following form in $T^{\ast, (1,0)}\psi$ ($\alpha\in T^{\ast, (1,0)}_p Y$):
 \begin{alignat*}{2}
  J\left(\partial_{Q_{x,j},\alpha}\right) &= \partial_{Q_{y,j},\alpha},\quad J\left(\partial_{Q_{y,j},\alpha}\right) &&= -\partial_{Q_{x,j},\alpha},\\
  J\left(\partial_{P_{x,j},\alpha}\right) &= \partial_{P_{y,j},\alpha},\quad J\left(\partial_{P_{y,j},\alpha}\right) &&= -\partial_{P_{x,j},\alpha}.
 \end{alignat*}
 In \autoref{app:almost_complex_structures}, we show that $I_g$ takes the following form\footnote{To be precise, we have determined the form of $J^\ast_{\nabla^{h_R}}$ in coordinates $T^\ast \psi$ in \autoref{app:almost_complex_structures}. We obtain the form of $I_g\coloneqq dF^{-1}\circ J^\ast_{\nabla^{h_R}}\circ dF$ by applying $F$. Mind the change of signs for $P_{y,j}$ due to $F$.} in $T^{\ast, (1,0)}\psi$:
 \begin{alignat*}{2}
  I_g\left(\partial_{Q_{x,j},\alpha}\right) &= -\partial_{P_{x,j},\alpha},\quad I_g\left(\partial_{Q_{y,j},\alpha}\right) &&= -\partial_{P_{y,j},\alpha},\\
  I_g\left(\partial_{P_{x,j},\alpha}\right) &= \partial_{Q_{x,j},\alpha},\quad I_g\left(\partial_{P_{y,j},\alpha}\right) &&= \partial_{Q_{y,j},\alpha}.
 \end{alignat*}
 With these formulas, one easily verifies by direct computation that $J$ and $I_g$ commute at the point $\alpha$. Since the same argument can be repeated for every $\alpha\in T^{\ast, (1,0)}_pY$ and $p\in Y$, the result is true for all points of $X$ concluding the proof.
\end{proof}

Before we conclude this subsection and turn our attention to the deformation of HHSs, we should briefly illustrate one feature of proper PHHSs which is absent in the integrable case: The fact that the Hamiltonian vector fields $X^{\Omega_R}_{\mH_R}$ and $J(X^{\Omega_R}_{\mH_R})$ of a PHHS $(X,J;\Omega_R,\mH_R)$ do not need to be $J$-preserving (cf. Remark \autoref{rem:x_0-dependance}). Example \autoref{ex:constructing_PHHS} beautifully demonstrates that feature. To see this, let us first develop a criterion that tells us in which cases the vector fields $X^{\Omega_R}_{\mH_R}$ and $J(X^{\Omega_R}_{\mH_R})$ are $J$-preserving:

\begin{proposition}\label{prop:J-preserving_criteria}
 Let $(X,J;\Omega_R,\mH_R)$ be a PHHS with Hamiltonian vector field $X_\mH = 1/2 (X^{\Omega_R}_{\mH_R} - iJ(X^{\Omega_R}_{\mH_R}))$ and $2$-form $\Omega_I = -\Omega_R (J\cdot,\cdot)$. Then, $X^{\Omega_R}_{\mH_R}$ is $J$-preserving if and only if $d\Omega_I (X^{\Omega_R}_{\mH_R},\cdot,\cdot) \equiv 0$. Similarly, $J(X^{\Omega_R}_{\mH_R})$ is $J$-preserving if and only if $d\Omega_I (J(X^{\Omega_R}_{\mH_R}),\cdot,\cdot) \equiv 0$.
\end{proposition}

\begin{proof}
 Take the notations from above and let $V\in\{X^{\Omega_R}_{\mH_R}, J(X^{\Omega_R}_{\mH_R})\}$. We want to determine in which cases $V$ is $J$-preserving. Thereto, we use Proposition 2.10 in Chapter IX of \cite{kobayashi1969}: A vector field $V$ on a smooth manifold $X$ with almost complex structure $J$ is $J$-preserving if and only if the following equation is fulfilled for every vector field $W$:
 \begin{gather*}
  [V, J(W)] = J([V,W]).
 \end{gather*}
 Since $\Omega_R$ is non-degenerate, we find that $V$ is $J$-preserving if and only if
 \begin{gather*}
  \iota_{[V, J(W)]}\Omega_R - \iota_{J([V,W])}\Omega_R = 0
 \end{gather*}
 holds for every vector field $W$ on $X$. Let us now compute the left-hand side of this equation. First, we obtain by definition of $\Omega_I$:
 \begin{gather*}
  \iota_{J([V,W])}\Omega_R = -\iota_{[V,W]}\Omega_I.
 \end{gather*}
 Now remember the relation (cf. Proposition 3.10 in Chapter I of \cite{kobayashi1963}):
 \begin{gather*}
  \iota_{[V,W]} = [L_V,\iota_W].
 \end{gather*}
 Furthermore, recall that by Definition \autoref{def:PHHS_1} the $1$-forms $\Omega_R (V,\cdot)$ and\linebreak $\Omega_I (V,\cdot) = -\Omega_R (J(V),\cdot)$ are exact. Together with $d\Omega_R = 0$ and Cartan's magic formula, this implies:
 \begin{gather*}
  L_V\Omega_R = 0,\quad L_V\Omega_I = \iota_V d\Omega_I,
 \end{gather*}
 where $L_V$ denotes the Lie derivative of $V$. These relations allow us to compute:
 \begin{align*}
  \iota_{[V,J(W)]}\Omega_R &= [L_V,\iota_{J(W)}]\Omega_R = L_V\left(\iota_{J(W)}\Omega_R\right) - \iota_{J(W)}\left(L_V\Omega_R\right)\\
  &= -L_V\left(\iota_W\Omega_I\right) = -[L_V,\iota_W]\Omega_I - \iota_W\left(L_V\Omega_I\right)\\
  &= -\iota_{[V,W]}\Omega_I - \iota_W\iota_V d\Omega_I = \iota_{J([V,W])}\Omega_R - \iota_W\iota_V d\Omega_I.
 \end{align*}
 In total, we find that $V$ is $J$-preserving if and only if
 \begin{gather*}
  \iota_{[V, J(W)]}\Omega_R - \iota_{J([V,W])}\Omega_R = - \iota_W\iota_V d\Omega_I = 0
 \end{gather*}
 holds for every vector field $W$ on $X$ concluding the proof.
\end{proof}

Now let us consider Example \autoref{ex:constructing_PHHS} again with $f = 1$, $h = e^{x_1}$, and\linebreak $H = \mH_R = -y_1$. Then, the Hamiltonian vector fields and the exterior derivative of $\Omega^g_I$ are given by:
\begin{gather*}
 X^{\Omega_R}_{\mH_R} = -\pa{y_2},\quad J_g(X^{\Omega_R}_{\mH_R}) = e^{-x_1}\pa{x_2},\quad d\Omega^g_I = e^{x_1} dx_1\wedge dx_2\wedge dy_1.
\end{gather*}
Using Proposition \autoref{prop:J-preserving_criteria}, we immediately see that $X^{\Omega_R}_{\mH_R}$ is $J_g$-preserving, while $J_g(X^{\Omega_R}_{\mH_R})$ is not. By choosing $H = \mH_R = ay_1 + be^{-x_1}$ ($a,b\in\R$), Example \autoref{ex:constructing_PHHS} even shows that neither $X^{\Omega_R}_{\mH_R}$ nor $J(X^{\Omega_R}_{\mH_R})$ of a PHHS $(X,J;\Omega_R,\mH_R)$ need to be $J$-preserving.

\subsection*{Deforming HHSs by proper PHHSs}

In this subsection, we examine the question: ``How `large' is the set of proper PHHSs within the set of all PHHSs?'' The answer and the main result of this subsection is that proper PHHSs are generic if $\text{dim}_\mathbb{R}(X)>4$. To prove this, we first reduce the genericity of proper PHHSs to the claim that every HHS can be deformed by proper PHHSs. Afterwards, we show this claim in two steps. The first step is to locally bring every HHS into standard form. Secondly, we deform the HHS standard form by proper PHHSs within a small neighborhood.\\
We start by recalling the definition of a generic property. For a topological space $B$ and a subset $A\subset B$, we call the property that an element is contained in $A$ \textbf{generic} if $B\backslash A$ is a meager subset of $B$ in the sense of Baire. In particular, the property $a\in A$ is generic if $A$ is an open and dense subset of $B$. Now let $X$ be a smooth manifold. We introduce the following notations for the set of almost complex structures, of PHSMs, and of PHHSs on $X$:
\begin{align*}
 \mathcal{J}_\text{a.c.} (X)&\coloneqq \{J\in\Gamma \text{End}(TX)\mid J^2 = -1\},\\
 \text{PHSM} (X)&\coloneqq \{(J,\Omega_R)\in\mathcal{J}_\text{a.c.}(X)\times\Omega^2 (X)\mid (X,J;\Omega_R)\text{ is a PHSM}\},\\
 \text{PHHS} (X)&\coloneqq \{(J,\Omega_R,\mH_R)\in\mathcal{J}_\text{a.c.}(X)\times\Omega^2 (X)\times C^\infty (X,\mathbb{R})\mid\\
 &\phantom{\coloneqq \{\ }(X,J;\Omega_R,\mH_R)\text{ is a PHHS}\}.
\end{align*}
Similarly, we denote the set of complex structures, of HSMs, and of HHSs on $X$ by:
\begin{align*}
 \mathcal{J}_\text{c} (X)&\coloneqq \{J\in\mathcal{J}_\text{a.c.}(X)\mid J\text{ is integrable}\},\\
 \text{HSM} (X)&\coloneqq \{(J,\Omega_R)\in \text{PHSM} (X)\mid J\text{ is integrable}\},\\
 \text{HHS} (X)&\coloneqq \{(J,\Omega_R,\mH_R)\in \text{PHHS} (X)\mid J\text{ is integrable}\}.
\end{align*}
Lastly, we write for the set of proper, i.e., non-integrable almost complex structures, of proper PHSMs, and of proper PHHSs on $X$:
\begin{align*}
 \mathcal{J}_\text{p} (X)&\coloneqq \mathcal{J}_\text{a.c.} (X)\backslash \mathcal{J}_\text{c} (X),\\
 \text{PHSM}_\text{p} (X)&\coloneqq \text{PHSM} (X)\backslash \text{HSM} (X),\\
 \text{PHHS}_\text{p} (X)&\coloneqq \text{PHHS} (X)\backslash \text{HHS} (X).
\end{align*}
We equip every set with the topology induced by the compact-open topology. We wish to show that $\mathcal{J}_\text{p} (X)\subset \mathcal{J}_\text{a.c.} (X)$, $\text{PHSM}_\text{p} (X)\subset \text{PHSM} (X)$, and $\text{PHHS}_\text{p} (X)\subset \text{PHHS} (X)$ are open and dense subsets if $\text{dim}_\mathbb{R} (X)>4$. Clearly, the sets in question are open subsets, since all former mentioned subsets contain exactly those elements from their respective supersets for which the Nijenhuis tensor $N_J$ of $J$ does not vanish. The argument is completed by noting that $N_J\neq 0$ is an open condition.\\
It is left to show that the subsets are dense in their respective supersets. This can be done by showing that the ``integrable'' sets $\mathcal{J}_\text{c} (X)$, $\text{HSM} (X)$, and $\text{HHS}(X)$ are contained in the boundary of the ``proper'' sets $\mathcal{J}_\text{p}$, $\text{PHSM}_\text{p} (X)$, and $\text{PHHS}_\text{p} (X)$, respectively. The last statement is true if the ``integrable'' elements can be deformed by ``proper'' elements:

\begin{definition}[Proper deformation]\label{def:deformation}
 Let $X$ be a smooth manifold with almost complex structure $J$ on it. $(X,J^\varepsilon)$ is called a \textbf{deformation} of $(X,J)$ if $J^\varepsilon$ describes a smooth\footnote{Here, smooth path means that the map $J^{\cdot}:X\times\mathbb{R}\to \text{\normalfont End}(TX)$, $(x,\varepsilon)\mapsto J^{\varepsilon}_x$ is smooth. Similar remarks apply to $\Omega^\varepsilon_R$ and $\mH^\varepsilon_R$.} path in $\mathcal{J}_\text{\normalfont a.c.} (X)$ with start point $J^0 = J$. Now let $\Omega_R$ be a $2$-form on $X$ such that $(X,J;\Omega_R)$ is a PHSM. Then, $(X,J^\varepsilon; \Omega^\varepsilon_R)$ is a deformation of $(X,J;\Omega_R)$ if $(J^\varepsilon, \Omega^\varepsilon_R)$ describes a smooth path in $\text{\normalfont PHSM}(X)$ with start point $(J^0, \Omega^0_R) = (J,\Omega_R)$. If, additionally, $\mH_R$ is a function on $X$ such that $(X,J;\Omega_R,\mH_R)$ is a PHHS, we say $(X,J^\varepsilon; \Omega^\varepsilon,\mH^\varepsilon_R)$ is a deformation of $(X,J;\Omega_R,\mH_R)$ if $(J^\varepsilon, \Omega^\varepsilon_R, \mH^\varepsilon_R)$ describes a smooth path in $\text{\normalfont PHHS}(X)$ with start point $(J^0, \Omega^0_R,\mH^0_R) = (J,\Omega_R,\mH_R)$. We call a deformation \textbf{proper} if the corresponding $J^\varepsilon$ is not integrable for $\varepsilon\neq 0$.
\end{definition}

Before we continue with the proof of genericity, we shall quickly address one aspect concerning our definitions: In the definition of a deformation, we have neglected the imaginary parts $\Omega_I$ and $\mH_I$, even though they are crucial for the definition of HSMs and HHSs. One might wonder whether this is justified, i.e., whether every deformation of a PHSM or a PHHS automatically gives us a suitable $1$-parameter family $\Omega^\varepsilon_I$ and $\mH^\varepsilon_I$ of imaginary parts. Regarding the form $\Omega_I$, this is certainly true by simply setting $\Omega^\varepsilon_I\coloneqq -\Omega^\varepsilon_R (J^\varepsilon\cdot,\cdot)$. However, $\mH_I$ defined as a primitive of $\Omega_R (J(X^{\Omega_R}_{\mH_R}),\cdot)$ might be more problematic. It is not obvious why a smooth $1$-parameter family of exact $1$-forms $\alpha^\varepsilon$ should be the differential of a smooth $1$-parameter family of functions $f^\varepsilon$. Nevertheless, this is still true in the case of $X$ being contractible:

\begin{proposition}[Primitive of $\alpha^\varepsilon$]\label{prop:primitive}
 Let $X$ be a smooth contractible manifold, $\alpha^\varepsilon$ a smooth $1$-parameter family of exact $1$-forms on $X$, and\linebreak $f\in C^\infty (X,\mathbb{R})$ such that $\alpha^0 = df$. Then, there exists a smooth $1$-parameter family of functions $f^\varepsilon$ on $X$ such that $\alpha^\varepsilon = df^\varepsilon$ and $f^0 = f$.
\end{proposition}

\begin{proof}
 Let $X$, $\alpha^\varepsilon$, and $f$ be as above. $X$ is contractible, so there is a smooth map $F:[0,1]\times X\to X$ such that $F(0,x) = x_0$ and $F(1,x) = x$ for every $x\in X$ and some $x_0\in X$. For fixed $\varepsilon$, $F^\ast\alpha^\varepsilon$ is a closed $1$-form on $[0,1]\times X$ and can be expressed as (using $F_t:X\to X$, $F_t (x)\coloneqq F (t,x)$ for every $t\in [0,1]$):
 \begin{gather*}
  F^\ast\alpha^\varepsilon = F^\ast_t\alpha^\varepsilon + \beta^\varepsilon_t dt,
 \end{gather*}
 where $\beta^\varepsilon_t$ is a function on $X$ smoothly depending on $\varepsilon$ and $t$. For fixed $t$ and $\varepsilon$, $F^\ast_t\alpha^\varepsilon$ can be understood as a closed $1$-form on $X$. Applying the exterior derivative to $F^\ast\alpha^\varepsilon = F^\ast_t\alpha^\varepsilon + \beta^\varepsilon_t dt$ and using the closedness of $F^\ast\alpha^\varepsilon$ and $F^\ast_t\alpha^\varepsilon$ yields:
 \begin{gather*}
  \frac{d}{dt}F^\ast_t\alpha^\varepsilon = d\beta^\varepsilon_t.
 \end{gather*}
 In total, we obtain:
 \begin{align*}
  \alpha^\varepsilon = \text{id}^\ast_X\alpha^\varepsilon - \text{const}^\ast_{x_0}\alpha^\varepsilon = F^\ast_1\alpha^\varepsilon - F^\ast_0\alpha^\varepsilon = \int\limits^1_0 \frac{d}{dt}(F^\ast_t\alpha^\varepsilon)\, dt = d\int\limits^1_0 \beta^\varepsilon_t dt\eqqcolon df^\varepsilon.
 \end{align*}
 Shifting $f^\varepsilon$ by a constant to match $f$ at $\varepsilon = 0$ concludes the proof.
\end{proof}

To show that properness is generic, we will only use local deformations. In particular, we can always shrink the neighborhood in which we deform an integrable structure such that it becomes contractible. Thus, we can rightfully neglect deformations of the imaginary parts $\Omega_I$ and $\mH_I$ here.\\
Let us return to the proof of genericity. As explained, it suffices to find a proper deformation of every ``integrable'' element to show genericity. We construct the desired proper deformations in two steps:
\begin{enumerate}
 \item We bring the ``integrable'' elements locally into standard form.
 \item We deform the standard form ``properly'' within a small neighborhood.
\end{enumerate}
Regarding the first step, it is clear what the standard forms of complex manifolds and HSMs are and how to achieve them. For complex manifolds, the standard form is just $X = \mathbb{C}^m$ with complex structure $J = i$ and every complex manifold can be brought into standard form by holomorphic charts. For HSMs, the standard form is $X = \mathbb{C}^{2n}$ together with $\Omega = \sum_j dP_j\wedge dQ_j$ and every HSM can be brought into standard form by holomorphic Darboux charts (cf. Theorem \autoref{thm:holo_Darboux} and \autoref{app:darboux}).\\
For a HHS $(X,\Omega,\mH)$, we observe that $\mH$ is either locally constant or regular at some point, i.e., there exists $x_0\in X$ such that $d\mH_{x_0}\neq 0$. If the given HHS is locally constant near a point, then holomorphic Darboux coordinates describe the standard form of the HHS in question near that point. If the given HHS is regular near a point, then the following lemma brings it into standard form:

\begin{lemma}[Regular HHSs in standard form]\label{lem:HHS_in_standard_form}
 Let $(X,\Omega,\mH)$ be a HHS and $x_0\in X$ be a point such that $d\mH_{x_0}\neq 0$. Then, there exists a holomorphic chart $\phi = (z_1,\ldots, z_{2n}):U\to V\subset\mathbb{C}^{2n}$ of $X$ near $x_0\in U$ such that\linebreak $\Omega\vert_U = \sum^n_{j = 1} dz_{j+n}\wedge dz_j$ and $\mH\vert_U = z_{2n}$.
\end{lemma}

\begin{proof}
 Let $(X,\Omega,\mH)$ be a HHS and $x_0\in X$ be a point such that $d\mH_{x_0}\neq 0$. Without loss of generality, we can assume $\mH (x_0) = 0$. The proof consists of three steps:
 \begin{enumerate}
  \item Construct a holomorphic function $G$ defined locally near $x_0$ satisfying $\{\mH, G\}\coloneqq \Omega (X_\mH, X_G) = 1$.
  \item Find a holomorphic chart $\phi_3 = (z^\#_1,\ldots, z^\#_{2n}):U_3\to V_3$ of $X$ near $x_0$ such that $G\vert_{U_3} = z^\#_n$, $\mH\vert_{U_3} = z^\#_{2n}$, $X_\mH\vert_{U_3} = \pa{z^\#_n}$, and $X_G\vert_{U_3} = -\pa{z^\#_{2n}}$.
  \item We have $\Omega\vert_{U_3} = dz^\#_{2n}\wedge dz^\#_{n} + \Sigma$, where $\Sigma$ is a closed $2$-form only depending on $z^\#_1,\ldots, z^\#_{n-1}, z^\#_{n+1},\ldots, z^\#_{2n-1}$. Bring $\Sigma$ into standard form using Darboux's theorem for HSMs (cf. Theorem \autoref{thm:holo_Darboux} and \autoref{app:darboux}).
 \end{enumerate}
 \textbf{Step 1:} We pick a small open neighborhood $U_1$ of $x_0$ such that $U_1$ is the domain of a holomorphic chart $\phi_1 = (\hat z_1,\ldots, \hat z_{2n}):U_1\to V_1\subset\mathbb{C}^{2n}$ with $\phi_1 (x_0) = 0$ and $d\mH$ does not vanish on $U_1$. $\pa{\hat z_1},\ldots, \pa{\hat z_{2n}}$ form a holomorphic frame of the tangent bundle $T^{(1,0)}U_1$. Since the holomorphic Hamiltonian vector field $X_\mH$ of the HHS $(X,\Omega,\mH)$ does not vanish on $U_1$, we can replace one coordinate vector field, say $\pa{\hat z_n}$, in the collection $\pa{\hat z_1},\ldots, \pa{\hat z_{2n}}$ with $X_\mH$ and still obtain a frame of $T^{(1,0)}U_1$ after shrinking $U_1$ if necessary.\\
 Now let $\varphi^{X_\mH}_z:U_1\to X$ be the complex flow of the holomorphic Hamiltonian vector field $X_\mH$ for suitable $z\in\mathbb{C}$. $\varphi^{X_\mH}_z$ is constructed as follows: For $x\in X$, we denote by $\gamma_x$ the holomorphic curve from $\mathbb{C}$ to $X$ solving the following initial value problem:
 \begin{gather*}
  \gamma^\prime_x (z) = X_\mH (\gamma_x (z)),\quad \gamma_x (0) = x.
 \end{gather*}
 By Proposition \autoref{prop:holo_traj}, the curves $\gamma_x$ exist on a small neighborhood of $0\in\C$ and are locally unique. With this, we can set the flow of $X_\mH$ to be $\varphi^{X_\mH}_z (x)\coloneqq \gamma_x (z)$.\\
 Let us now consider the map $\phi^{-1}_2$ from $\mathbb{C}^{2n}$ to $X$ defined by:
 \begin{gather*}
  \phi^{-1}_2(\tilde z_1,\ldots, \tilde z_{2n})\coloneqq \varphi^{X_\mH}_{\tilde z_n}\circ \phi^{-1}_1 (\tilde z_1,\ldots, \tilde z_{n-1}, 0, \tilde z_{n+1},\ldots, \tilde z_{2n}).
 \end{gather*}
 The differential of $\phi^{-1}_2$ at $0\in\mathbb{C}^{2n}$ maps the standard complex basis of $\mathbb{C}^{2n}$ to the vector fields $\pa{\hat z_1},\ldots, \pa{\hat z_{n-1}}, X_{\mH}, \pa{\hat z_{n+1}}, \ldots, \pa{\hat z_{2n}}$ at $x_0\in X$. These vector fields form a complex basis at $x_0$, hence, $\phi^{-1}_2$ is a local biholomorphism near $0$ by the holomorphic version of the inverse function theorem. This gives us the holomorphic chart $\phi_2 = (\tilde z_1,\ldots, \tilde z_{2n}):U_2\to V_2$, where $U_2$ is a small open neighborhood of $x_0$.\\
 Next, we set $G:U_2\to\mathbb{C}$ to be the coordinate $\tilde z_{n}$. We find:
 \begin{gather*}
  dG_x (X_\mH (x)) = \left.\frac{d}{dz}\right\vert_{z = 0} (G\circ \varphi^{X_\mH}_z (x)) = 1\quad\forall x\in U_2.
 \end{gather*}
 This implies for the Poisson bracket of $\mH$ and $G$:
 \begin{gather*}
  \{\mH, G\} \coloneqq \Omega (X_\mH, X_G) = dG (X_\mH) = 1,
 \end{gather*}
 where $X_G$ is the holomorphic Hamiltonian vector field of the HHS $(U_2,\Omega\vert_{U_2}, G)$.\\\\
 \textbf{Step 2:} Taking the Hamiltonian vector field of a holomorphic function is a Lie algebra homomorphism, hence, we get:
 \begin{gather*}
  [X_\mH, X_G] = X_{\{\mH, G\}} = 0.
 \end{gather*}
 As in the real case, the commutativity of the vector fields implies the commutativity of their flows. This allows us to find a holomorphic chart\linebreak $\phi_3 = (z^\#_1,\ldots, z^\#_{2n}):U_3\to V_3$ of $X$ near $x_0$ such that $G\vert_{U_3} = z^\#_n$, $\mH\vert_{U_3} = z^\#_{2n}$, $X_\mH\vert_{U_3} = \pa{z^\#_n}$, and $X_G\vert_{U_3} = -\pa{z^\#_{2n}}$. The construction of $\phi_3$ makes use of the regular value theorem for complex manifolds. Consider the holomorphic map $f\coloneqq(H\vert_{U_2},G):U_2\to\mathbb{C}^2$. The map $f$ is a submersion due to $\Omega (X_\mH, X_G) = 1$. By the regular value theorem, the level sets of $f$ are complex submanifolds of $U_2$. The tangent space of a level set is given by the kernel of $df$. Now pick the level set $W\coloneqq f^{-1}(\mH (x_0), G(x_0)) = f^{-1}(0,0)$ containing $x_0$ and a holomorphic chart $\psi$ of $W$ near $x_0$ with $\psi (x_0) = 0$. Next, we define the map $\phi^{-1}_3$ via:
 \begin{gather*}
  \phi^{-1}_3 (z^\#_1,\ldots, z^\#_{2n})\coloneqq \varphi^{X_\mH}_{z^\#_n}\circ \varphi^{X_G}_{-z^\#_{2n}}\circ \psi^{-1} (z^\#_1,\ldots, z^\#_{n-1}, z^\#_{n+1}, \ldots, z^\#_{2n-1}).
 \end{gather*}
 The differential of $\phi^{-1}_3$ at $0\in\mathbb{C}^{2n}$ maps the standard complex basis of $\mathbb{C}^{2n}$ to the vector fields $X_\mH$, $-X_G$, and the coordinate vector fields $v_j$ of $\psi$ at $x_0$. $v_j$ satisfy $d\mH (v_j) = dG (v_j) = 0$ and, hence, are orthogonal to $X_\mH$ and $X_G$ with respect to the symplectic form $\Omega$. Because $\Span\{X_\mH, X_G\}\subset T^{(1,0)}_{x_0}X$ is a symplectic subspace, we have:
 \begin{gather*}
  T^{(1,0)}_{x_0}X = \Span\{X_\mH, X_G\}\oplus \Span\{X_\mH, X_G\}^{\perp\Omega} = \Span\{X_\mH, X_G\}\oplus \Span\{v_j\}.
 \end{gather*}
 Thus, $X_\mH$, $-X_G$, and $v_j$ form a complex basis of $T^{(1,0)}_{x_0}X$. We can again apply the inverse function theorem to show that $\phi^{-1}_3$ is locally biholomorphic near $0$. From $\phi^{-1}_3$, we obtain the holomorphic chart $\phi_3 = (z^\#_1,\ldots, z^\#_{2n}):U_3\to V_3$ of $X$ near $x_0$. Using the commutativity of $\varphi^{X_\mH}_{z^\#_n}$ and $\varphi^{X_G}_{-z^\#_{2n}}$, we easily compute $\pa{z^\#_n} = X_\mH\vert_{U_3}$ and $\pa{z^\#_{2n}} = -X_G\vert_{U_3}$. By construction, $\mH$ only depends on $z^\#_{2n}$ and $G$ only depends on $z^\#_{n}$. Integrating $dG(X_\mH) = d\mH (-X_G) = \Omega (X_\mH, X_G) = 1$ gives us $G = z^\#_{n}$ and $\mH = z^\#_{2n}$.\\\\
 \textbf{Step 3:} Let us inspect $\Omega$ in the chart $\phi_3$ more closely. The coordinates of $\phi_3$ satisfy $\iota_{\pa{z^\#_{n}}}\Omega\vert_{U_3} = -dz^\#_{2n}$ and $\iota_{\pa{z^\#_{2n}}}\Omega\vert_{U_3} = dz^\#_{n}$. Therefore, $\Omega\vert_{U_3}$ takes the following form:
 \begin{gather*}
  \Omega\vert_{U_3} = dz^\#_{2n}\wedge dz^\#_n + \sum\limits_{i,j\neq n,2n} \Omega_{ij} dz^\#_{i}\wedge dz^\#_{j}\eqqcolon dz^\#_{2n}\wedge dz^\#_n + \Sigma.
 \end{gather*}
 $\Omega\vert_{U_3}$ and $dz^\#_{2n}\wedge dz^\#_n$ are closed $2$-forms on $U_3$, thus, $\Sigma$ is also a closed $2$-form on $U_3$. $d\Sigma = 0$ implies that the partial derivatives of $\Omega_{ij}$ with respect to $z^\#_{n}$ and $z^\#_{2n}$ have to vanish for every $i,j\neq n,2n$. Thus, the restriction of $\Sigma$ to a hyperplane $z^\#_n = c_1$ and $z^\#_{2n} = c_2$ does not depend on the values $c_1$ and $c_2$. Furthermore, the restriction of $\Sigma$ to such a hyperplane is a holomorphic symplectic $2$-form. This is a direct consequence of $\Omega\vert_{U_3}$ being a holomorphic symplectic $2$-form on $U_3$ and the hyperplane being a complex submanifold of $U_3$. Therefore, we can apply Theorem \autoref{thm:holo_Darboux} to $\Sigma$ giving us a holomorphic chart $\Psi$ of the hyperplane $z^\#_n = 0$ and $z^\#_{2n} = 0$ in which $\Sigma$ assumes the standard form. Replacing the coordinates $z^\#_1,\ldots, z^\#_{n-1}, z^\#_{n+1}, \ldots, z^\#_{2n-1}$ of $\phi_3$ with the coordinates of $\Psi$ yields the holomorphic chart $\phi = (z_1,\ldots, z_{2n}):U\to V$ of $X$ near $x_0$ (choose $\Psi (x_0) = 0$):
 \begin{gather*}
  \phi^{-1} (z_1,\ldots, z_{2n})\coloneqq \varphi^{X_\mH}_{z_n}\circ \varphi^{X_G}_{-z_{2n}}\circ \Psi^{-1} (z_1,\ldots, z_{n-1}, z_{n+1}, \ldots, z_{2n-1}).
 \end{gather*}
 In this chart, $\Omega$ and $\mH$ take the form $\Omega\vert_U = \sum^n_{j = 1} dz_{j+n}\wedge dz_j$ and $\mH\vert_U = z_{2n}$ concluding the proof.
\end{proof}

\begin{corollary}[Every regular HHS is locally integrable as a Hamiltonian system]\label{cor:integrable}
 Let $(X,\Omega,\mH)$ be a HHS and $x_0\in X$ be a point such that $d\mH_{x_0}\neq 0$. Then, there exists an open neighborhood $U$ of $x_0$ and holomorphic functions $F_i,G_i:U\to\mathbb{C}$ for $i\in\{1,\ldots, n\}$, $2n\coloneqq \text{\normalfont dim}_\mathbb{C}(X)$, such that:
 \begin{gather*}
  F_n = \mH\vert_U,\quad \{F_i, G_j\} = \delta_{ij},\quad \{F_i, F_j\} = \{G_i,G_j\} = 0\quad\forall i,j\in\{1,\ldots, n\}.
 \end{gather*}
 In particular, every regular HHS is locally integrable as a Hamiltonian system.
\end{corollary}

\begin{proof}
 Take the chart $\phi$ from Lemma \autoref{lem:HHS_in_standard_form} and set $F_j\coloneqq z_{j + n}$ as well as $G_j\coloneqq z_j$.
\end{proof}

\begin{remark}[Every regular RHS is locally integrable]\label{rem:integrable}
 Lemma \autoref{lem:HHS_in_standard_form} and Corollary \autoref{cor:integrable} are also true for regular RHSs with almost exactly the same proofs. In particular, every regular RHS is locally integrable.
\end{remark}

With the ``integrable'' elements in standard form, let us now turn our attention to the second step. We need to construct local proper deformations of the standard form. Here, we only show how to deform HHSs explicitly. The deformations of HSMs and complex manifolds can be obtained in a similar way.

\begin{proposition}[Deformation of standard HHSs]\label{prop:deformation_of_standard_HHS}
 Let $(X,\Omega,\mH)$ be a HHS with complex structure $J$ and decompositions $\Omega = \Omega_R + i\Omega_I$ and $\mH = \mH_R + i\mH_I$, where $X = \mathbb{C}^{2n}\ (n>1)$, $J = i$, $\Omega = \sum^n_{j = 1} z_{j+n}\wedge z_j$, and $\mH \equiv c$ for some constant $c\in\mathbb{C}$ or $\mH = z_{2n}$. Furthermore, let $U\subset \mathbb{C}^{2n}$ be any non-empty open subset. Then, there exists a proper deformation $(X,J^\varepsilon;\Omega^\varepsilon_R,\mH^\varepsilon_R)$ of the HHS $(X,\Omega,\mH)$ such that $J^\varepsilon\vert_{X\backslash U} = J\vert_{X\backslash U}$, $\Omega^\varepsilon_R = \Omega_R$, and $\mH^\varepsilon_R = \mH_R$ for every $\varepsilon\in\mathbb{R}$.
\end{proposition}

\begin{proof}
 The idea for the deformation is based on Example \autoref{ex:constructing_PHHS}. Let $(X,\Omega,\mH)$ be a HHS as above with non-empty open subset $U\subset X$. Now pick a non-constant smooth function $f:X\cong\mathbb{R}^{4n}\to\mathbb{R}$ satisfying:
 \begin{gather*}
  f(x)\geq 0\ \forall x\in X,\quad f(x) = 0\ \forall x\in X\backslash U.
 \end{gather*}
 We define the $1$-parameter family of smooth functions $r^\varepsilon$ on $X$ as follows:
 \begin{gather*}
  r^\varepsilon (x)\coloneqq 1 + \varepsilon^2 f(x)\quad \forall x\in X\ \forall \varepsilon\in\mathbb{R}.
 \end{gather*}
 Using $r^\varepsilon$, we define the $1$-parameter family of almost complex structures $J^\varepsilon$:
 \begin{align*}
  J^\varepsilon (\pa{x_1}) &\coloneqq r^\varepsilon\pa{y_1},\ J^\varepsilon (\pa{x_{n+1}}) \coloneqq \frac{1}{r^\varepsilon}\pa{y_{n+1}},\\
  J^\varepsilon (\pa{y_1}) &\coloneqq -\frac{1}{r^\varepsilon}\pa{x_1},\ J^\varepsilon (\pa{y_{n+1}}) \coloneqq -r^\varepsilon\pa{x_{n+1}},\\
  J^\varepsilon (\pa{x_j}) &\coloneqq \pa{y_j},\quad J^\varepsilon (\pa{y_j}) \coloneqq -\pa{x_j}\qquad\qquad\quad \forall j\in\{2,\ldots, n, n+2,\ldots, 2n\},
 \end{align*}
 where $\pa{x_j}$ and $\pa{y_j}$ are the vector fields coming from the coordinates\linebreak $(z_1 = x_1 + iy_1,\ldots, z_{2n} = x_{2n} + iy_{2n})\in\mathbb{C}^{2n}$. Clearly, $J^\varepsilon$ coincides with $J = i$ for $\varepsilon = 0$. Moreover, one easily checks that $J^\varepsilon$ satisfies $J^\varepsilon\vert_{X\backslash U} = J\vert_{X\backslash U}$ and is $\Omega_R$-anticompatible for every $\varepsilon\in\mathbb{R}$.\\
 Next, we set $\Omega^\varepsilon_R\equiv \Omega_R$ and $\Omega^\varepsilon_I\coloneqq -\Omega^\varepsilon_R (J^\varepsilon\cdot,\cdot)$ for every $\varepsilon\in\mathbb{R}$. As in Example \autoref{ex:constructing_PHHS}, the exterior derivative of $\Omega^\varepsilon_I$ can be expressed as:
 \begin{gather*}
  d\Omega^\varepsilon_I = \varepsilon^2 df\wedge\left( dy_{n+1}\wedge dx_1 - \left(\frac{1}{r^\varepsilon}\right)^2dx_{n+1}\wedge dy_1\right).
 \end{gather*}
 We see that for every $\varepsilon\neq 0$ there exists a point $x_0\in U$ satisfying $d\Omega^\varepsilon_{I,x_0}\neq 0$, as $f$ is non-constant. Thus, $(X,J^\varepsilon;\Omega^\varepsilon_R)$ becomes a proper PHSM for every $\varepsilon\neq 0$ due to Theorem \autoref{thm:rel_HSM_PHSM}.\\
 The only thing left to check is that $\mH^\varepsilon_R\equiv \mH_R$ is compatible with the PHSM $(X,J^\varepsilon;\Omega^\varepsilon_R)$, i.e., $d[\Omega^\varepsilon_R (J^\varepsilon (X^{\Omega^\varepsilon_R}_{\mH^\varepsilon_R}),\cdot)] = 0$ for every $\varepsilon$. In the case of $\mH\equiv c$, this is trivially true, because then the Hamiltonian vector field $X^{\Omega^\varepsilon_R}_{\mH^\varepsilon_R}$ vanishes. For $\mH = z_{2n}$, we first realize that the equation above is equivalent to the condition that $\mH^\varepsilon_R$ is the real part of some pseudo-holomorphic function $\mH^\varepsilon$ (cf. Remark \autoref{rem:H_pseudo-holo}). We already know that $\mH^\varepsilon\equiv \mH$ is pseudo-holomorphic with respect to $J = i$. Since $\mH^\varepsilon$ only depends on the last component $z_{2n}$ and both $J^\varepsilon$ and $J = i$ act in the same way on the last components of $X$ for every $\varepsilon$ (as $n>1$), $\mH^\varepsilon\equiv\mH$ is also pseudo-holomorphic with respect to $J^\varepsilon$ for every $\varepsilon$. This turns $(X,J^\varepsilon;\Omega^\varepsilon_R,\mH^\varepsilon_R)$ into the desired deformation of $(X,\Omega,\mH)$ concluding the proof.
\end{proof}

\begin{remark}[The case $n=1$]\label{rem:constant_Ham}
 Upon closer inspection, we note that the proof also works for $n = 1$ if $\mH$ is constant or if we disregard $\mH$ entirely, i.e., if we are only interested in proper deformations of HSMs or complex manifolds. Indeed, we only need the condition $n>1$ to ensure that $\mH^\varepsilon_R$ is the real part of some pseudo-holomorphic function $\mH^\varepsilon$. However, the constant function is pseudo-holomorphic with respect to any almost complex structure.
\end{remark}

Combining Lemma \autoref{lem:HHS_in_standard_form} with Proposition \autoref{prop:deformation_of_standard_HHS} proves the following theorem:

\begin{theorem}[Proper PHHSs are generic]\label{thm:generic}
 Let $X$ be a smooth manifold, then the following statements apply depending on the real dimension of $X$:
 \begin{enumerate}
  \item If $\text{\normalfont dim}_\mathbb{R}(X) = 2$: Every almost complex structure on $X$ is integrable and automatically a complex structure.
  \item If $\text{\normalfont dim}_\mathbb{R}(X) > 2$: Every complex manifold $(X,J)$ and HSM $(X,\Omega)$ admits a proper deformation. In particular, the non-integrable almost complex structures and the proper PHSMs on $X$ are generic within the set of all almost complex structures and PHSMs on $X$, respectively.
  \item If $\text{\normalfont dim}_\mathbb{R}(X) > 4$: Every HHS $(X,\Omega, \mH)$ admits a proper deformation. In particular, the proper PHHSs on $X$ are generic within the set of all PHHSs on $X$.
 \end{enumerate}
\end{theorem}

\begin{proof}
 For (i), it suffices to check that the Nijenhuis tensor always vanishes in two dimensions which is a simple calculation. In the case of (ii), we first find a holomorphic chart for $(X,J)$ or a holomorphic Darboux chart for $(X,\Omega)$ to locally bring these manifolds into standard form. Afterwards, we apply Proposition \autoref{prop:deformation_of_standard_HHS} (disregarding $\mH$) to locally deform the manifolds in the chosen charts. For (iii), we use Lemma \autoref{lem:HHS_in_standard_form} to locally bring $(X,\Omega,\mH)$ into standard form and then employ Proposition \autoref{prop:deformation_of_standard_HHS} to find a local deformation finishing the proof.
\end{proof}

Before we conclude this section, let us quickly comment on PHHSs $(X,J;\Omega_R,\mH_R)$ in dimension $\text{dim}_\mathbb{R}(X) = 4$. These systems are the only geometrical object studied in this subsection to which Theorem \autoref{thm:generic} does not apply. The reason is that Proposition \autoref{prop:deformation_of_standard_HHS} fails for $\text{dim}_\mathbb{R}(X) = 4$, since $\mH_R = x_2$ cannot be the real part of a complex function which is pseudo-holomorphic with respect to deformations $J^\varepsilon$ as chosen in the proof of Proposition \autoref{prop:deformation_of_standard_HHS}. Nevertheless, Statement 3 of Theorem \autoref{thm:generic} might still be true for $\text{dim}_\mathbb{R}(X) = 4$. One possible way to prove this could be to modify Proposition \autoref{prop:deformation_of_standard_HHS}. Instead of choosing $\mH^\varepsilon_R$ to be independent of $\varepsilon$, we could allow for general deformations $\mH^\varepsilon_R$ of $\mH_R$. Finding such deformations $\mH^\varepsilon_R$ involves solving a second order PDE with boundary conditions. However, the existence of non-trivial solutions to the given problem might be forbidden by the Nijenhuis tensor. We elaborate on this thought: Let $X$ be a smooth manifold with almost complex structure $J$, $N_J$ be the Nijenhuis tensor of $J$, and $f:X\to\mathbb{C}$ be a pseudo-holomorphic map, i.e., $df\circ J = i\cdot df$. A straightforward calculation reveals that $df\left(N_J (V,W)\right) = 0$ for all vector fields $V$ and $W$ on $X$. Thus, the image of $N_{J,x}$ is contained within the kernel of $df_x$ for any pseudo-holomorphic function $f$ and any point $x\in X$. This implies that there are at most $1/2(\text{dim}_\mathbb{R}(X)-r_x)$ pseudo-holomorphic functions on $X$ whose differentials at a given point $x\in X$ are $\mathbb{C}$-linearly independent, where $r_x$ is the rank\footnote{The Nijenhuis tensor satisfies the relation $N_J (J(V),W) = -JN_J (V,W)$, hence, its rank is even.} of $N_{J,x}$. In four dimensions, the rank of $N_J$ alone does not exclude the existence of non-trivial pseudo-holomorphic functions: Let $V_1$, $V_2$, $J(V_1)$, and $J(V_2)$ be a local frame of $X$. One easily sees that, because of the symmetries of $N_J$, i.e., $N_J (V,W) = -N_J(W,V)$ and $N_J(J(V), W) = -JN_J (V,W)$, $N_J (V_1, V_2)$ and $N_J (V_1, J(V_2)) = -JN_J (V_1,V_2)$ are the only two components of $N_J$ in the local frame $V_1$, $V_2$, $J(V_1)$, and $J(V_2)$ which are not redundant. Thus, the rank of the Nijenhuis tensor is at most $2$ in four dimensions. For instance, the deformation $J^\varepsilon$ from Proposition \autoref{prop:deformation_of_standard_HHS} for $n=1$,
 \begin{gather*}
  J^\varepsilon (\pa{x_1}) \coloneqq r^\varepsilon\pa{y_1},\ J^\varepsilon (\pa{x_{2}}) \coloneqq \frac{1}{r^\varepsilon}\pa{y_{2}},\quad J^\varepsilon (\pa{y_1}) \coloneqq -\frac{1}{r^\varepsilon}\pa{x_1},\ J^\varepsilon (\pa{y_{2}}) \coloneqq -r^\varepsilon\pa{x_{2}},
 \end{gather*}
 yields for $r^\varepsilon\in C^\infty(\mathbb{R}^4,\mathbb{R}_{+})$:
 \begin{gather*}
  N_{J^\varepsilon} (\pa{x_1},\pa{x_2}) = \frac{1}{r^\varepsilon}J^\varepsilon\left(N_{J^\varepsilon} (\pa{x_1},\pa{y_2})\right) = \sum\limits_{a\in\{x,y\}}\sum^2_{i,j = 1,\ i\neq j}\left(\pa{a_i}\ln (r^\varepsilon)\right)\cdot \pa{a_j}.
 \end{gather*}
 Hence, the (real) rank of $N_{J^\varepsilon,x}$ is $2$ for $dr^\varepsilon_x \neq 0$ and $0$ for $dr^\varepsilon_x = 0$.\\
 Nevertheless, the rank of $N_J$ is a rather weak bound for the number of independent pseudo-holomorphic functions. The exact number is given by\linebreak $1/2(\text{dim}_\mathbb{R}(X) - k)$, where $k$ is the rank of the IJ-bundle\footnote{Confer \cite{muskarov1986} for the definition of the IJ-bundle and a detailed investigation of the relation between the Nijenhuis tensor and pseudo-holomorphic functions.} over $X$ containing the image of $N_J$. Even though there are non-integrable almost complex structures $J$ in four dimensions whose IJ-bundle does not have full rank, e.g. $J^\varepsilon$ with $r^\varepsilon = e^{-x_1}$ as in Example \autoref{ex:constructing_PHHS}, it is not clear whether and why this should also apply to the almost complex structures $J^\varepsilon$ as above for general functions $r^\varepsilon$.

\chapter[Kähler Duality of Complex Coadjoint Orbits]{Kähler Duality of Complex Coadjoint Orbits\chaptermark{Kähler Duality}}
\chaptermark{Kähler Duality}
\label{chap:duality}
Kähler structures on coadjoint\footnote{Most statements in this chapter apply to both adjoint and coadjoint orbits, as they are isomorphic via a suitably chosen metric $g$ in the cases we consider. For the sake of brevity, we often only talk about coadjoint orbits. We adopt a similar convention for tangent and cotangent bundles.} orbits of compact Lie groups have been studied since the 50s. The interest in this topic stems from the fact that these orbits are compact homogeneous Kähler manifolds which were the first manifolds known to admit a Kähler-Einstein metric with positive scalar curvature\footnote{One major question in the field of Kähler-Einstein manifolds regards the existence of Kähler-Einstein metrics on compact spaces. For negative or vanishing scalar curvature/first Chern class, this question was answered by Yau in the 70s (cf. the Calabi conjecture). The Fano case (positive scalar curvature) remained an open problem for the longest time and was only solved about ten years ago by Chen-Donaldson-Sun (cf. \cite{ChenI}\cite{ChenII}\cite{ChenIII}).} (cf. \cite{Besse2007}). In fact, the coadjoint orbits of compact Lie groups classify all simply-connected compact homogeneous Kähler manifolds (cf. \cite{Besse2007}). At the end of the 80s and the beginning of the 90s, the question arose whether these Kähler structures persist in the complex category. This problem was solved by Kronheimer (cf. \cite{Kronheimer1990}) and Kovalev (cf. \cite{Kovalev1996}): They proved that the coadjoint orbits of semisimple complex reductive groups, i.e, complex Lie groups with semisimple compact real forms, admit the structure of a Hyperkähler manifold.\\
In this chapter, we show that these orbits are not only Hyperkähler, but also possess a holomorphic Kähler structure. From a symplectic viewpoint, Hyperkähler and holomorphic Kähler manifolds are quite similar\footnote{Confer \autoref{app:kaehler} for details, in particular Lemma \autoref{lem:hyper_vs_holo}.}. Both consist of a symplectic form $\omega$ and two complex structures $I$ and $J$ where $I$ is $\omega$-anticompatible and $J$ is $\omega$-compatible in the sense that $\omega (I\cdot, I\cdot) = -\omega$ and $\omega (J\cdot,J\cdot) = \omega$. The only difference is the commutation relation of $I$ and $J$: $I$ and $J$ anticommute for Hyperkähler manifolds ($IJ = -JI$), while they commute for holomorphic Kähler manifolds ($IJ = JI$). On coadjoint orbits, the similarities between Hyperkähler and holomorphic Kähler structures go even further. To explain this, we first note that $\omega$ and $I$ always induce a holomorphic symplectic structure given by the closed holomorphic two-form $\Omega\coloneqq \omega - i\omega (I\cdot,\cdot)$. In general, there is no reason why the form $\Omega$ of a Hyperkähler structure and the form $\Omega$ of a holomorphic Kähler structure should be related in any way. On a coadjoint orbit, however, the Hyperkähler and holomorphic Kähler structure in question possess the same form $\Omega$, namely the holomorphic Kirillov-Kostant-Souriau form $\KKS$. Therefore, these two structures only differ by $J$.\\
We say a space $X$ exhibits \textbf{Kähler duality} if $X$ admits Hyperkähler and holomorphic Kähler structures in a somewhat natural manner. As discussed, coadjoint orbits exhibit Kähler duality. Of course, we instinctively ask whether the Kähler duality of coadjoint orbits is just mere coincidence or caused by some deeper reason. We believe that the Kähler duality of coadjoint orbits originates from a similar structure on double cotangent bundles. Precisely speaking, we claim in this chapter that double cotangent bundles naturally exhibit Kähler duality and that this Kähler duality is connected to coadjoint orbits via a suitable reduction process (Hyperkähler/holomorphic Kähler reduction, cf. \autoref{app:reduction}).\\
In order to sketch how Kähler duality on double cotangent bundles occurs, let us recall a famous result\footnote{Confer Theorem \autoref{thm:cotangent_kaehler}.} due to Guillemin-Stenzel (cf. \cite{Stenzel1990} and \cite{Guillemin1991}) and Lempert-Sz{\H{o}}ke (cf. \cite{Lempert1991} and \cite{Szoeke1991}): If $(M,g)$ is a real-analytic Riemannian manifold, then $(T^\ast M,-\omega_{\can}, J_g)$ is a Kähler manifold where $\omega_{\can}$ is the canonical symplectic form on $T^\ast M$ and $J_g$ is the unique complex structure on $T^\ast M$ adapted to $g$. The Kähler duality on $T^\ast (T^\ast M)$ now arises as depicted in the following diagram\footnote{This diagram has to be taken with a grain of salt: In general, these structures do not exist on all of $T^\ast M$ or $T^\ast (T^\ast M)$, but only on an open neighborhood of the zero section $M\subset T^\ast M\subset T^\ast (T^\ast M)$.}:
\begin{center}
 \begin{tikzcd}[column sep=3.6em,row sep=2em]
  \textbf{Hyperk.}& (T^\ast M, g_{T^\ast M}) \arrow{r}{\text{Stenzel}} & (T^\ast (T^\ast M), -\omega_{\can}, \phi^\ast_1 J_{g_{T^\ast M}}, T^\ast J_g)\\
  (M,g) \arrow{ur}[sloped,above]{\text{Stenzel}} \arrow{dr}[sloped,below]{\text{\footnotesize Complexification}} & &\\
  \textbf{Holo. K.}& (T^\ast M, g_\C) \arrow{r}[below]{\text{Stenzel}} & (T^\ast (T^\ast M), -\omega_{\can}, \phi^\ast_2 J_{g_\C}, T^\ast J_g)
 \end{tikzcd}
\end{center}
Here, the upper path illustrates how to obtain a Hyperkähler structure, while the lower one does the same for holomorphic Kähler structures.\\
\textbf{Hyperkähler path:} Starting with $(M,g)$, we can apply Stenzel's theorem to obtain the Kähler metric $g_{T^\ast M}\coloneqq -\omega_{\can} (\cdot,J_g\cdot) =  \omega_{\can} (J_g\cdot,\cdot)$ on $T^\ast M$. This gives us another real-analytic Riemannian manifold $(T^\ast M,g_{T^\ast M})$ allowing us to apply Stenzel's theorem again which results in the Kähler manifold $(T^\ast (T^\ast M),-\omega_{\can}, J_{g_{T^\ast M}})$. To construct the second\footnote{This complex structure plays the role of $I$.} complex structure on $T^\ast (T^\ast M)$, we first observe that the cotangent bundle of a complex manifold naturally inherits a complex structure from its base manifold. It is a standard result from complex geometry that the cotangent bundle together with its inherited complex structure and the (negative) canonical symplectic structure forms a holomorphic symplectic manifold. If we choose $(T^\ast M, J_g)$ to be our complex base manifold, we denote the induced complex structure on $T^\ast(T^\ast M)$ by $T^\ast J_g$. We know at this point that $-\omega_{\can}$ is compatible with $J_{g_{T^\ast M}}$ and anticompatible with $T^\ast J_g$. Thus, $(T^\ast (T^\ast M), -\omega_{\can}, J_{g_{T^\ast M}}, T^\ast J_g)$ is a Hyperkähler\footnote{The existence of a Hyperkähler structure on the cotangent bundle of a Kähler manifold was already observed by Kaledin (cf. \cite{Kaledin1997}) and Feix (cf. \cite{Feix2001}). Also note that Hyperkähler structures on double cotangent bundles were already described by Bielawski (cf. \cite{Bielawski2003}).} manifold if $T^\ast J_g$ and $J_{g_{T^\ast M}}$ anticommute. However, it may occur, depending on the choice of $g$, that $T^\ast J_g$ and $J_{g_{T^\ast M}}$ do not anticommute. In this case, we need to modify $J_{g_{T^\ast M}}$ by a diffeomorphism $\phi_1:T^\ast (T^\ast M)\to T^\ast (T^\ast M)$ to obtain $\phi^\ast_1 J_{g_{T^\ast M}}$ (cf. Conjecture \autoref{con:commutation_relation} and Lemma \autoref{lem:stenzel_for_kaehler}).\\
\textbf{Holomorphic Kähler path:} Stenzel's theorem does not only tell us that $(T^\ast M,-\omega_{\can}, J_g)$ is a Kähler manifold, but also that the fiberwise map\linebreak $T^\ast M\to T^\ast M$, $\alpha\mapsto -\alpha$ is a real structure (cf. Theorem \autoref{thm:cotangent_kaehler}). Its real form is the zero section $M\subset T^\ast M$. Since the real form $M$ carries a real-analytic metric $g$, we can find a unique holomorphic continuation of $g$ on $T^\ast M$ (cf. \autoref{app:real_structures}). Call the real part of this holomorphic metric $g_\C$. $g_\C$ is a real-analytic semi-Riemannian metric on $T^\ast M$, thus, we can again apply Stenzel's theorem, but this time to $(T^\ast M, g_\C)$\footnote{Stenzel's theorem also works for semi-Riemannian manifolds (cf. \cite{Szoeke2004}).}. As before, the resulting complex structure $J_{g_\C}$ on $T^\ast (T^\ast M)$ might not satisfy the appropriate commutation relation with $T^\ast J_g$, so we have to modify $J_{g_\C}$ by a diffeomorphism $\phi_2$ to obtain $\phi^\ast_2 J_{g_\C}$ (cf. Conjecture \autoref{con:commutation_relation} and Lemma \autoref{lem:stenzel_for_holo_metrics}). Hence, $(T^\ast (T^\ast M), -\omega_{\can}, \phi^\ast_2 J_{g_\C}, T^\ast J_g)$ is a holomorphic Kähler manifold if $T^\ast J_g$ and $\phi^\ast_2 J_{g_\C}$ commute.\\
To relate $T^\ast (T^\ast M)$ to coadjoint orbits, we choose $(M,g)$ to be a compact Lie group $G_\R$ with bi-invariant Riemannian metric $g$. In this case, $T^\ast G_\R$ is isomorphic to the universal complexification\footnote{Confer Definition \autoref{def:group_real_str_main} and the explanation afterwards.} $G$ of $G_\R$. Thus, the previously described process yields Hyperkähler/holomorphic Kähler structures on $T^\ast G$. As explained in \autoref{app:reduction}, the cotangent bundle $T^\ast G$ becomes a coadjoint orbit $\mathcal{O}^\ast$ after reduction. At the same time, the symplectic form $-\omega_{\can}$ on $T^\ast G$ reduces to the Kirillov-Kostant-Souriau form $\kks$ on $\mathcal{O}^\ast$. From that perspective, it seems plausible that the Kähler structures in question are also compatible with the reduction process.\\
The proofs we present in this thesis regarding the Kähler duality of $T^\ast (T^\ast M)$ and its relation to coadjoint orbits are incomplete: A complete proof for the commutation relations of $T^\ast J_g$, $\phi^\ast_1 J_{g_{T^\ast M}}$, and $\phi^\ast_2 J_{g_\C}$ is missing. Moreover, we do not carry out the reduction process in detail. To make up for that, we check the commutation relations for flat $g$ and sketch how reduction could possibly relate a suitable Kähler structure on $T^\ast G$ to $\mathcal{O}^\ast$ (cf. \autoref{app:reduction}).\\
\autoref{chap:duality} is divided into four parts: \autoref{sec:lie_groups} offers a short introduction to Lie groups and discusses the most important constructions we use throughout this chapter. In \autoref{sec:semi-kaehler}, we recall the construction of Kähler structures on coadjoint orbits of compact Lie groups following Chapter 8 in \cite{Besse2007}. As we want to apply this construction to complex Lie groups afterwards, we also generalize it to skew-symmetric orbits of groups carrying a bi-invariant semi-Riemannian metric. In \autoref{sec:holo_semi-kaehler}, we adapt the results of \autoref{sec:semi-kaehler} to the complex category. In particular, we prove that coadjoint orbits of complex reductive groups admit holomorphic Kähler structures. The last section (\autoref{sec:duality}) is concerned with the Kähler duality of coadjoint orbits and its origin. Barring the commutation relation, we show that double cotangent bundles exhibit Kähler duality. Furthermore, we investigate the case $M = G_\R$ and relate the resulting Kähler structures on $T^\ast G$ to Bremigan's construction (cf. \cite{Bremigan2000}).

\newpage
\section[Preliminaries on Lie Groups]{Preliminaries on Lie Groups\sectionmark{Lie Groups}}
\sectionmark{Lie Groups}
\label{sec:lie_groups}
In this section, we briefly recapitulate the basics of Lie groups. Its purpose is two-fold: On one hand, we fix the notations and recall the theorems we use throughout \autoref{chap:duality}. On the other hand, this section shall serve as a short introduction to Lie groups for those unfamiliar with the topic. The reader acquainted with Lie groups may choose to skip \autoref{sec:lie_groups}.\\
There are several great introductory books on Lie groups which offer an in-depth analysis of the topic, for instance \cite{Bump2004}, \cite{Varadarajan1984}, or the first few chapters in \cite{Hamilton2017}. Here, we are only interested in a small selection of constructions and statements starting with the definition of a Lie group:

\begin{definition}[Lie group]\label{def:lie_group}
 We call $G$ \textbf{Lie group} if $G$ is a group and a smooth manifold such that the multiplication $G\times G\to G,\, (g,h)\mapsto gh$ and the inversion $G\to G,\, g\mapsto g^{-1}$ are smooth maps. We denote the left multiplication with a group element $g\in G$ by $L_g:G\to G$ and the right multiplication by $R_g:G\to G$, i.e.:
 \begin{gather*}
  L_g (h)\coloneqq gh,\quad R_g (h)\coloneqq hg\quad\forall g,h\in G.
 \end{gather*}
\end{definition}

Closely related to the notion of a Lie group is the concept of a Lie algebra:

\begin{definition}[Lie algebra]\label{def:lie_algebra}
 Let $G$ be a Lie group. The \textbf{Lie algebra} of $G$ is the tangent space $\mathfrak{g}\coloneqq T_e G$ of $G$ at the neutral element $e$. It naturally comes with a Lie bracket. To construct this Lie bracket, we consider left-invariant\footnote{One can also define the Lie bracket by identifying $\mathfrak{g}$ with the space of right-invariant vector fields. Note that the two conventions do not lead to the same Lie bracket: The Lie brackets differ by a sign.} vector fields $X$, i.e., smooth vector fields $X\in\Gamma (TG)$ satisfying:
 \begin{gather*}
  (L_g)_\ast X = X\quad\forall g\in G.
 \end{gather*}
 Denote the set of left-invariant vector fields by $\Gamma_G (TG)$. Observe that the commutator of two left-invariant vector fields is again left-invariant:
 \begin{gather*}
  (L_g)_\ast [X,Y] = [(L_g)_\ast X, (L_g)_\ast Y] = [X,Y]\quad\forall g\in G\, \forall X,Y\in\Gamma_G (TG).
 \end{gather*}
 Identifying $\mathfrak{g}$ with $\Gamma_G (TG)$ via the map $\cdot^L:\mathfrak{g}\to\Gamma_G (TG)$ given by
 \begin{gather*}
  X^L (g)\coloneqq dL_{g,e} X\quad\forall X\in\mathfrak{g}
 \end{gather*}
 now allows us to define the Lie bracket $[\cdot,\cdot]:\mathfrak{g}\times\mathfrak{g}\to\mathfrak{g}$:
 \begin{gather*}
  [X,Y]^L\coloneqq [X^L,Y^L]\quad\forall X,Y\in\mathfrak{g}.
 \end{gather*}
\end{definition}

Lie groups (and algebras) form categories. Their morphisms are called Lie group (or algebra) homomorphisms:

\begin{definition}[Homomorphism]\label{def:homo}
 Let $G$ and $H$ be Lie groups with Lie algebras $\mathfrak{g}$ and $\mathfrak{h}$, respectively. A \textbf{Lie group homomorphism} from $G$ to $H$ is a smooth map $R:G\to H$ satisfying $R(g_1g_2) = R(g_1)R(g_2)$ for all $g_1,g_2\in G$. Similarly, a \textbf{Lie algebra homomorphism} from $\mathfrak{g}$ to $\mathfrak{h}$ is a linear map $\rho:\mathfrak{g}\to\mathfrak{h}$ satisfying $\rho ([X,Y]) = [\rho (X), \rho (Y)]$ for all $X,Y\in\mathfrak{g}$.
\end{definition}

It is easy to show that, if $R:G\to H$ is a Lie group homomorphism, the differential $dR_e:\mathfrak{g}\to\mathfrak{h}$ is a Lie algebra homomorphism (cf. \cite{Hamilton2017} for the proof). This intimate relation is reflected by the following important example:

\begin{example}[Conjugation and adjoint action]\label{ex:conj_adj}
 Let $G$ be a Lie group. For every $g\in G$, the \textbf{conjugation} $c_g:G\to G$, $c_g (h)\coloneqq ghg^{-1}$ is a Lie group homomorphism. Its differential $\Ad (g):\mathfrak{g}\to\mathfrak{g},\ \Ad (g)\coloneqq dc_{g,e}$ is a Lie algebra homomorphism called the \textbf{adjoint action} or the adjoint representation. The differential of the map $\Ad:G\to\GL (\mathfrak{g})$ is often denoted by $\ad:\mathfrak{g}\to\End (\mathfrak{g})$. It is straightforward to verify that $\ad$ is just the Lie bracket (cf. \cite{Hamilton2017}):
 \begin{gather*}
  \ad_X Y\coloneqq \ad (X)(Y) = [X,Y]\quad\forall X,Y\in\mathfrak{g}.
 \end{gather*}
 Since $\Ad (g)$ is a Lie algebra homomorphism, it satisfies:
 \begin{gather}
  \Ad (g)[X,Y] = [\Ad (g)X, \Ad (g)Y]\quad\forall g\in G\, \forall X,Y\in\mathfrak{g}.\label{eq:adj_rep}
 \end{gather}
 Taking the derivative with respect to $g$ now yields for the map $\ad_Z$:
 \begin{gather}
  \ad_Z [X,Y] = [\ad_ZX, Y] + [X,\ad_ZY]\quad\forall X,Y,Z\in\mathfrak{g}.\label{eq:jacobi}
 \end{gather}
 In light of $\ad_X Y = [X,Y]$, we can interpret \autoref{eq:jacobi} as the Jacobi identity. We will make heavy use of \autoref{eq:adj_rep} and \eqref{eq:jacobi} in \autoref{sec:semi-kaehler} and \autoref{sec:holo_semi-kaehler}.
\end{example}

The next concept we introduce is the exponential map for Lie groups:

\begin{definition}[Exponential Map]\label{def:exp_lie}
 Let $G$ be a Lie group. The \textbf{exponential map} $\exp:\mathfrak{g}\to G$ is defined by
 \begin{gather*}
  \exp (X)\coloneqq \gamma_X (1),
 \end{gather*}
 where $\gamma_X$ is the unique integral curve of the left-invariant vector field $X^L$ satisfying $\gamma_X (0) = e$.
\end{definition}

As the definition indicates, the exponential map satisfies the usual flow properties:
\begin{enumerate}
 \item $\exp (0) = e$,
 \item $\exp ((s+t)X) = \exp (sX)\exp (tX)$,
 \item $\exp (tX)^{-1} = \exp (-tX)$.
\end{enumerate}
One easily verifies that $\exp$ is smooth and well-defined on all of $\mathfrak{g}$. Its derivative at the neutral element $e$ is $d\exp_e = \id_{\mathfrak{g}}$. Hence, by the inverse function theorem, $\exp$ is a local diffeomorphism near $e$ and gives rise to charts of $G$ after choosing a basis of $\mathfrak{g}$. The usefulness of the exponential map lies in the fact that it allows us to easily describe tangent vectors as curves. To understand this statement, we note that any tangent space $T_gG$ can be identified with $\mathfrak{g}$ via:
\begin{gather*}
 dL_{g,e}:\mathfrak{g}\to T_g G,\ X\mapsto X^L (g).
\end{gather*}
Hence, $g\exp (tX)$ is a curve through $g$ with tangent vector $X^L (g)$.

\pagebreak

Let us now turn our attention to metrics on Lie groups:

\begin{definition}[Bi-invariant metric]\label{def:bi-invariant}
 Let $G$ be a Lie group and $h$ a semi-Riemannian metric on $G$. We call $h$ left-invariant if it satisfies $L^\ast_g h = h$ for all $g\in G$. Similarly, we say $h$ is right-invariant if it satisfies $R^\ast_g h = h$ for all $g\in G$. $h$ is called \textbf{bi-invariant} if $h$ is both left- and right-invariant.
\end{definition}

Left-invariant metrics as well as right-invariant metrics are in one-to-one correspondence with non-degenerate scalar products on $\mathfrak{g}$: If $\skcdot$ is a non-degenerate scalar product on $\mathfrak{g}$, then $\skcdot^L$ is left-invariant and $\skcdot^R$ is right-invariant, where $\skcdot^L$ and $\skcdot^R$ are defined by:
\begin{gather*}
 \skcdot^L_g\coloneqq \sk{dL_{g^{-1},g}\cdot}{dL_{g^{-1},g}\cdot},\quad \skcdot^R_g\coloneqq \sk{dR_{g^{-1},g}\cdot}{dR_{g^{-1},g}\cdot} \quad\forall g\in G.
\end{gather*}
Conversely, any left- or right-invariant metric $h$ can be written as $\skcdot^L$ or $\skcdot^R$ where $\skcdot$ is given by $h_e$.\\
As it turns out, the left- or right-invariant metric induced by $\skcdot$ is bi-invariant if and only if $\skcdot$ is \textbf{Ad-invariant}, i.e.:
\begin{gather}
 \sk{\Ad (g)\cdot}{\Ad (g)\cdot} = \skcdot\quad\forall g\in G.\label{eq:Ad-invariance}
\end{gather}
The proof of this equivalence is fairly simple: First, recall that the conjugation with $g\in G$ can be written as:
\begin{gather*}
 c_g = L_g\circ R_{g^{-1}} = R_{g^{-1}}\circ L_g.
\end{gather*}
Thus, the adjoint action as its differential satisfies:
\begin{gather*}
 \Ad (g) = dc_{g,e} = dL_{g,g^{-1}}\circ dR_{g^{-1},e} = dR_{g^{-1},g}\circ dL_{g,e}.
\end{gather*}
This shows that a left- or right-invariant metric $h$ is bi-invariant if and only if $h_e$ is $\Ad$-invariant.\\
One important example of an $\Ad$-invariant two-form is given by the \textbf{Killing form} $K:\mathfrak{g}\times\mathfrak{g}\to\R$, $K(u,v)\coloneqq-\tr\left(\ad_u\circ\ad_v\right)$. One easily verifies that $K$ is bilinear, symmetric, and $\Ad$-invariant. However, $K$ is usually not positive definite, in fact, it is often not even non-degenerate. It is a famous fact, known as Cartan's criterion, that $K$ is non-degenerate if and only if $\mathfrak{g}$ is semi-simple. Moreover, $K$ is positive semi-definite if $\mathfrak{g}$ is the Lie algebra of a compact Lie group.\\
At this point, it should be noted that $\ad_Z$ is skew-symmetric with respect to any $\Ad$-invariant scalar product $\skcdot$:
\begin{gather}
 \sk{\ad_ZX}{Y} = -\sk{X}{\ad_ZY}\quad\forall X,Y,Z\in\mathfrak{g}.\label{eq:ad-invariant}
\end{gather}
\autoref{eq:ad-invariant} directly follows from \autoref{eq:Ad-invariance} by taking the derivative with respect to $g$. In fact, \autoref{eq:Ad-invariance} and \eqref{eq:ad-invariant} are equivalent if $G$ is connected (cf. \autoref{sec:holo_semi-kaehler}).\\
\autoref{eq:ad-invariant} poses a rather heavy condition on the existence of bi-invariant metrics: Milnor has shown that a connected Lie group $G$ admits a positive definite, bi-invariant metric if and only if $G$ is isomorphic to the Cartesian product of a compact Lie group and an Abelian Lie group (cf. \cite{Milnor1976}). If we drop the assumption that the metric needs to be positive definite, we find a larger class of Lie groups that admit bi-invariant metrics. For instance, complex reductive Lie groups always admit bi-invariant metrics of signature $(n,n)$ (cf. \autoref{sec:holo_semi-kaehler} and \autoref{app:complex_lie_groups}).\\
Before we conclude \autoref{sec:lie_groups}, we want to present two famous theorems from the theory of Lie groups -- Cartan's subgroup theorem and a special case of Godement's theorem -- and afterwards apply them to the adjoint action. Let us begin with Cartan's subgroup theorem. To formulate it, we first need to review the two different notions of subgroups of Lie groups:

\begin{definition}[Lie subgroup]
 Let $G$ and $H$ be two Lie groups. We say that $H$ is an \textbf{immersed subgroup} of $G$ if there exists a Lie group homomorphism $R:H\to G$ which is also an injective immersion. $H$ is a \textbf{Lie subgroup} or an embedded subgroup of $G$ if there exists a Lie group homomorphism $R:H\to G$ which is also an embedding. If $H\subset G$, we call $H$ an immersed (or embedded) subgroup of $G$ if the inclusion $H\hookrightarrow G$ is a Lie group homomorphism and an immersion (or embedding).
\end{definition}

Obviously, every embedded subgroup is also an immersed one. The converse, however, is false. The most prominent counterexample is the one-parameter subgroup
\begin{gather*}
 T_{p,q}\coloneqq\{(e^{ip\alpha}, e^{iq\alpha})\mid \alpha\in\R\}
\end{gather*}
of the torus $T^2 = S^1\times S^1$, where the ratio of $p\in\R$ and $q\in\R$ is irrational. The problem in this case is that $T_{p,q}$ is not closed in the topology of $T^2$. By Cartan's subgroup theorem, this is the only obstacle:

\begin{theorem}[Cartan's subgroup theorem]\label{thm:cartan}
 Let $G$ be a Lie group and $H\subset G$ be a subset. $H$ is a Lie subgroup of $G$ if and only if $H$ is a subgroup of $G$ and closed in the topology of $G$.
\end{theorem}

\begin{proof}
 The direction ``$\Rightarrow$'' is a rather straightforward computation. The converse direction is much more involved. We refer to \cite{Hamilton2017} for the proof.
\end{proof}

Next, we wish to formulate Godement's theorem for Lie groups. To do so, we have to cover group actions and quotients:

\begin{definition}[Group action]
 Let $G$ be a Lie group and $M$ a smooth manifold. A left $G$-\textbf{action} on $M$ is a smooth map $G\times M\to M$, $(g,p)\mapsto gp$ satisfying:
 \begin{gather*}
  g(hp) = (gh)p\quad\text{and}\quad ep = p\quad\forall g,h\in G\ \forall p\in M.
 \end{gather*}
 The \textbf{orbit} of $G$ through $p\in M$ is defined by:
 \begin{gather*}
  \mathcal{O}\coloneqq\{gp\mid g\in G\}\subset M.
 \end{gather*}
 $p\sim gp$ defines an equivalence relation whose equivalence classes are the orbits of $G$. The \textbf{quotient} $M/G$ is defined to be the set of all equivalence classes, i.e., the space of all orbits:
 \begin{gather*}
  M/G\coloneqq \{\mathcal{O}\subset M\mid \mathcal{O}\text{ orbit of }G\}.
 \end{gather*}
 The \textbf{stabilizer} of $p\in M$ is defined as follows:
 \begin{gather*}
  G_p\coloneqq \{g\in G\mid gp = p\}\subset G.
 \end{gather*}
 A $G$-action on $M$ is called \textbf{free} if $G_p = \{e\}$ for all points $p\in M$. It is called \textbf{proper} if the map
 \begin{gather*}
  G\times M\to M\times M,\ (g,p)\mapsto (gp,p)
 \end{gather*}
 is proper.
\end{definition}

\begin{remark}[Right action]
 Naturally, there also exists the notion of a right action which is just a smooth map $M\times G\to M$, $(p,g)\mapsto pg$ satisfying:
 \begin{gather*}
  (pg)h = p(gh)\quad\text{and}\quad pe = p\quad\forall g,h\in G\ \forall p\in M.
 \end{gather*}
 Orbits, stabilizers, and so on are defined analogously for right actions. The terms ``left action'' and ``right action'' are interchangeable, since every left action gives rise to a right action via $pg\coloneqq g^{-1}p$ and vice versa. In light of this observation, we often just write ``action'' and only specify ``left'' or ``right'' if the need arises.
\end{remark}

It is natural to ask at this point whether the quotient $M/G$ admits the structure of a smooth manifold. Godement's theorem tells us that this is the case if the $G$-action is free and proper:

\begin{theorem}[Godement's theorem\footnote{The standard version of Godement's theorem deals with the question whether, given a manifold $M$ with equivalence relation $R$ on it, the set $M/R$ of equivalence classes is a manifold. In our case, the relation $R$ is given by the $G$-action, namely $p\sim gp$.} for Lie groups]\label{thm:godement}
 Let $G$ be a Lie group and $M$ be a smooth manifold with a $G$-action on it. If the action is free and proper, then $M/G$ is a smooth manifold and the canonical projection\linebreak $\pi:M\to M/G$ is a surjective submersion.
\end{theorem}

\begin{proof}
 Again, we refer to \cite{Hamilton2017} for the proof.
\end{proof}

\begin{remark}\label{rem:godement}\ \vspace{-0.2cm}
 \begin{enumerate}
  \item If $G$ is a Lie group with Lie subgroup $H\subset G$, then $G\times H\to G$, $(g,h)\mapsto gh$ is a free and proper right $H$-action (cf. \cite{Hamilton2017}). Hence, by\linebreak Theorem \autoref{thm:godement}, $G/H$ is a smooth manifold carrying the left $G$-action $g_1 [g_2]\coloneqq [g_1g_2]$. Furthermore, the canonical projection $\pi:G\to G/H$ is a surjective submersion which is also equivariant with respect to the left $G$-actions on $G$ and $G/H$, i.e.:
  \begin{gather*}
   \pi (g_1g_2) = [g_1g_2] = g_1[g_2] = g_1\pi (g_2)\quad\forall g_1,g_2\in G.
  \end{gather*}
  \item Usually, it can be quite challenging to determine whether a $G$-action is proper. However, if $G$ is a compact Lie group, then it is rather simple to show that any $G$-action is proper.
  \item Let $G$ be a Lie group acting on a manifold $M$. If the action is free and proper, then an additional consequence of Theorem \autoref{thm:godement} is that\linebreak $\pi:M\to M/G$ is a $G$-principal bundle. In fact, a $G$-action on $M$ is free and proper if and only if $\pi:M\to M/G$ is a $G$-principal bundle with respect to that action (cf. \cite{Hamilton2017}).
 \end{enumerate}
\end{remark}

Two important examples of group actions on manifolds are the adjoint action $\Ad: G\times\mathfrak{g}\to\mathfrak{g}$ (cf. Example \autoref{ex:conj_adj}) and its dual, the coadjoint action\linebreak $\Ad^\ast:G\times\mathfrak{g}^\ast\to\mathfrak{g}^\ast$ ($\mathfrak{g}^\ast\coloneqq\{\alpha\mid\alpha:\mathfrak{g}\to\R\text{ linear}\}$):
\begin{gather*}
 \Ad^\ast (g)\alpha\coloneqq \alpha\circ \Ad (g^{-1})\quad\forall g\in G\, \forall \alpha\in\mathfrak{g}^\ast.
\end{gather*}
Both are $G$-actions, since the conjugation $c_{\cdot}:G\to\Aut (G)$ is a group homomorphism. In \autoref{sec:semi-kaehler} and \autoref{sec:holo_semi-kaehler}, we explore the geometrical structure of the (co)adjoint orbits, i.e., the orbits of $\Ad$ ($\Ad^\ast$). Before we can analyze the orbits in detail, we need to study the question whether the (co)adjoint orbits admit a manifold structure. The answer to this question is positive, as one can show with the help of Theorem \autoref{thm:cartan} and \autoref{thm:godement}:

\begin{lemma}[$G/G_p\cong\mathcal{O}\subset M$ immersed submanifold]\label{lem:immersed_orbit}
 Let $G$ be a Lie group, $M$ be a smooth manifold with a left $G$-action on it, and $p\in M$ be a point. Furthermore, let $G_p\subset G$ be the stabilizer of $p$ and $\mathcal{O}\subset M$ be the orbit of $G$ through $p$. Then, the map
 \begin{gather*}
  f_p: G/G_p\to M,\ [g]\mapsto gp
 \end{gather*}
 is an injective immersion whose image is $\mathcal{O}$. Moreover, $f_p$ is equivariant with respect to the left $G$-actions on $G/G_p$ and $M$, i.e.:
 \begin{gather*}
  f_p (g_1[g_2]) = g_1 f_p([g_2])\quad\forall g_1, g_2\in G.
 \end{gather*}
\end{lemma}

\begin{proof}
 Consider the evaluation map $\ev_p:G\to M$, $\ev_p (g)\coloneqq gp$. Clearly, $\ev_p$ is continuous. Thus, the stabilizer $G_p\subset G$ is closed, since $G_p = \ev_p^{-1}(\{p\})$. As $G_p$ is obviously a subgroup of $G$, we can use Cartan's subgroup theorem (cf. Theorem \autoref{thm:cartan}) to show that $G_p$ is a Lie subgroup of $G$. This allows us to apply Theorem \autoref{thm:godement} to $G_p\subset G$ (cf. Remark \autoref{rem:godement}) proving that $G/G_p$ is a smooth manifold and $\pi:G\to G/G_p$ is a surjective submersion.\\
 Now consider the following commuting diagram:
 \begin{center}
 \begin{tikzcd}
  G \arrow[r, "\ev_p"] \arrow[d,"\pi"]
  & M\\
  G/G_p \arrow[ur, "f_p", dashed]
 \end{tikzcd}
 \end{center}
 The smoothness of $\ev_p$ together with the universal property of quotient spaces (cf. Lemma 3.7.5 in \cite{Hamilton2017}) implies the smoothness of $f_p$. Clearly, $f_p$ is injective and its image is $\mathcal{O}$. Moreover, one directly verifies that $f_p$ is $G$-equivariant:
 \begin{gather*}
  f_p (g_1[g_2]) = f_p ([g_1g_2]) = (g_1g_2)p = g_1(g_2p) = g_1f_p ([g_2])\quad\forall g_1,g_2\in G.
 \end{gather*}
 It remains to be shown that $f_p$ is an immersion. Since $f_p$ is $G$-equivariant, it suffices to show that $df_{p,[e]}$ is injective. However, this directly follows from the previous commuting diagram and the observations (cf. \cite{Hamilton2017} for details):
 \begin{gather*}
  \ker d\ev_{p,e} = \mathfrak{g}_p\quad\text{and}\quad T_{[e]} (G/G_p) = \mathfrak{g}/\mathfrak{g}_p,
 \end{gather*}
 where $\mathfrak{g}_p$ is the Lie algebra of $G_p$.
\end{proof}

\newpage
\section[Semi-Kähler Structure of Coadjoint Orbits]{Semi-Kähler Structure of Coadjoint Orbits\sectionmark{Semi-Kähler Structure}}
\sectionmark{Semi-Kähler Structure}
\label{sec:semi-kaehler}
The existence of Kähler structures on (co)adjoint orbits is well-known and can be traced back to the 50s and 60s (cf. the introduction of Chapter 8 in \cite{Besse2007}). The construction of these Kähler structures (cf. \cite{Besse2007}) relies heavily on the fact that the groups in question are compact. This poses a problem: Eventually, our goal is to construct special Kähler structures on specific complex Lie groups (cf. \autoref{sec:holo_semi-kaehler}). However, complex Lie groups are often not compact\footnote{The adjoint action of any connected, compact, and complex Lie group is trivial due to the maximum principle. Consequently, such a group is Abelian implying that the identity component of every compact complex Lie group is just a torus.}. To circumvent this obstacle, we generalize the construction in \cite{Besse2007} to skew-symmetric\footnote{Confer Definition \autoref{def:skew} and Remark \autoref{rem:skew}.} orbits of Lie groups admitting a bi-invariant semi-Riemannian metric. The trade-off is that the constructed Kähler metric is usually not positive definite. Precisely speaking, we show that, under the mild conditions listed above, a (co)adjoint orbit of a Lie group carries a semi-Kähler\footnote{Confer Definition \autoref{def:kaehler_in_main} or \autoref{app:kaehler}.} structure\linebreak (Theorem \autoref{thm:semi-kaehler} and \autoref{thm:semi-kaehler_unique}).\\
The semi-Kähler structures we develop are canonical to some extent: The complex structure of an adjoint orbit does not depend on any choices, while the symplectic structure is canonical for coadjoint orbits. However, we need to identify adjoint and coadjoint orbits to construct the semi-Kähler structures in question. This identification is not unique and depends on the choice of an $\Ad$-invariant scalar product.\\
Throughout this section, Chapter 8 of \cite{Besse2007} shall serve as a guideline for developing the main result. \autoref{sec:semi-kaehler} is divided into four parts: First, we construct the canonical complex structure of adjoint orbits. Afterwards, we recall the canonical symplectic structure of coadjoint orbits, given by the Kirillov-Kostant-Souriau form. In the third part, we combine the complex and the symplectic structure via an $\Ad$-invariant, non-degenerate scalar product to obtain a semi-Kähler structure on (co)adjoint orbits. Lastly, we cross-check our construction against \cite{Besse2007} by applying it to compact Lie groups.

\subsection*{Complex Structure $\mathbf{J}$ of Adjoint Orbits}

The construction of $J$ is carried out in three steps: First, we define a complex structure $J_w$ on the vector space $\im\ad_w\subset\mathfrak{g}$ for each skew-symmetric element $w\in\mathfrak{g}$. Afterwards, we use the identification $T_w\mathcal{O}\cong\im\ad_w$ to equip a skew-symmetric adjoint orbit $\mathcal{O}$ with the complex structure $J$. Lastly, we compute the Nijenhuis tensor $N_J$ in order to show that $J$ is integrable.

\subsubsection*{Pointwise Construction of $\mathbf{J}$}

An adjoint orbit only admits the canonical complex structure $J$ if the orbit contains a skew-symmetric element, as we will soon see. Since the notion of skew-symmetric elements is purely algebraic, we begin by reviewing some facts about Lie algebras:

\pagebreak

\begin{definition}[Skew-symmetric elements]\label{def:skew}
 Let $(\mathfrak{g},[\cdot,\cdot])$ be a real Lie algebra. An element $w\in\mathfrak{g}$ is called \textbf{skew-symmetric} if the map $\ad_w\coloneqq [w,\cdot]\in\End (\mathfrak{g})$ satisfies the following properties:
 \begin{enumerate}
  \item Its complexification\footnote{Set $\mathfrak{g}_\C\coloneqq\mathfrak{g}\otimes_\R \C$ and extend $\ad_w$ to $\mathfrak{g}_\C$ by $\C$-linearity. Also note that we denote both the map on $\mathfrak{g}$ and its complexification by $\ad_w$.} $\ad_w\in\End (\mathfrak{g}_\C)$ is diagonalizable.
  \item The non-vanishing eigenvalues of $\ad_w\in\End (\mathfrak{g}_\C)$ are purely imaginary.
 \end{enumerate}
\end{definition}

We call this property ``skew-symmetric'', since $\ad_w$ is skew-symmetric with respect to some positive definite scalar product $\sk{\cdot}{\cdot}$ on $\mathfrak{g}$ if and only if $w\in\mathfrak{g}$ is a skew-symmetric element. One can easily verify this equivalence with the help of the spectral theorem.\\
Before we associate complex structures to skew-symmetric elements, note the following remarks:

\begin{remark}\label{rem:skew}\ \vspace{-0.3cm}
 \begin{enumerate}
  \item Condition (i) in Definition \autoref{def:skew} implies $\mathfrak{g} = \ker\ad_w\oplus\im\ad_w$ for every skew-symmetric element $w\in\mathfrak{g}$.
  \item If $w\in\mathfrak{g}$ is skew-symmetric, then the non-vanishing eigenvalues of $\ad_w$ come in pairs $(i\mu, -i\mu)$ with $\mu>0$, as $\ad_w$ is real, i.e., $\ad_w (\bar{v}) = \overline{\ad_w (v)}$, where $\bar{\cdot}$ denotes the complex conjugation.
  \item If $G$ is a Lie group, $\mathfrak{g}$ its Lie algebra, and $w\in\mathfrak{g}$ a skew-symmetric element, then every element in the adjoint orbit through $w$ is also skew-symmetric. This is a direct consequence of \autoref{eq:adj_rep}:
  \begin{gather*}
   \ad_{\Ad (g) w} = \Ad (g)\circ \ad_{w}\circ \Ad (g)^{-1}\quad\forall g\in G.
  \end{gather*}
  Hence, if an adjoint orbit contains a skew-symmetric element, we call the entire orbit skew-symmetric.
 \end{enumerate}
\end{remark}

The reason why skew-symmetric elements $w$ are important for our discussion is that the image of $\ad_w$ naturally comes with a complex structure $J_w$: 

\begin{proposition}[Pointwise construction of $J$]\label{prop:can_comp_str}
 Let $\mathfrak{g}$ be a Lie algebra and $w\in\mathfrak{g}$ be skew-symmetric. Then, $V\coloneqq\im\ad_w\subset\mathfrak{g}$ carries a canonical complex structure $J_w$ given by:
 \begin{gather*}
  J_w v\coloneqq \frac{1}{\mu} \ad_w v\quad\forall v\in E_\mu,
 \end{gather*}
 where we set $E_\mu\coloneqq V\cap (E_{i\mu}\oplus E_{-i\mu})$ and $E_{i\mu},\, E_{-i\mu}\subset\mathfrak{g}_\C$ are the eigenspaces of $\ad_w$ for the eigenvalues $i\mu,\, -i\mu$ ($\mu>0$), respectively.
\end{proposition}

\begin{proof}
 Clearly, the linear map $J_w:V\to V$ is well-defined. Indeed, one has $J_w (E_\mu)\subset E_\mu$ and the decomposition
 \begin{gather*}
  V = \bigoplus_{\mu>0} E_\mu,
 \end{gather*}
 where the direct sum is taken over the norm $\mu$ of non-vanishing\linebreak eigenvalues of $\ad_w$.\\
 We now show $J^2_w = -\id_V$. By definition of $J_w$, if suffices to prove $J^2_w (v) = -v$ for all $v\in E_\mu$ and $\mu>0$. Any $v\in E_\mu$ can be written as $v = u + \bar{u}$, where $u\in E_{i\mu}$ and, consequently, $\bar{u}\in E_{-i\mu}$. Applying $J^2_w$ to $v$ yields:
 \begin{gather*}
  J^2_w (v) = \frac{1}{\mu^2}\left(\ad^2_w (u) + \ad^2_w(\bar{u})\right) = i^2 u + (-i)^2 \bar{u} = -(u+\bar{u}) = -v.
 \end{gather*}
\end{proof}

\begin{remark}\label{rem:J_w}
 If $V$ is a real vector space with map $J\in\End (V)$ satisfying $J^2 = -\id_V$, we can decompose $V_\C$ into the spaces
 \begin{align*}
  V^{(1,0)}&\coloneqq E_i \equiv \left\{ v\in V_\C \mid J(v) = iv\right\}\quad\text{and}\\
  V^{(0,1)}&\coloneqq E_{-i} \equiv \left\{ v\in V_\C \mid J(v) = -iv\right\}.
 \end{align*}
 In the case of $J_w$, the spaces $V^{(1,0)}$ and $V^{(0,1)}$ are given by:
 \begin{gather*}
  V^{(1,0)} = \bigoplus_{\mu > 0} E_{i\mu},\quad V^{(0,1)} = \bigoplus_{\mu > 0} E_{-i\mu}.
 \end{gather*}
 In particular, restricting $J_w$ to $E_\mu$ gives:
  \begin{gather*}
   E_\mu^{(1,0)} = E_{i\mu},\quad E_\mu^{(0,1)} = E_{-i\mu}\qquad\forall \mu >0.
  \end{gather*}
\end{remark}

\subsubsection*{Global Construction of $\mathbf{J}$}

So far, we have constructed the complex structure $J_w$ on the vector space $\im\ad_w$. In order to equip a skew-symmetric adjoint orbit $\mathcal{O}$ with a complex structure, we need to identify $T_w\mathcal{O}$ with $\im\ad_w$ for every point $w\in\mathcal{O}$. For this, we first recall that $\mathcal{O}\subset\mathfrak{g}$ is an immersed submanifold (cf. Lemma \autoref{lem:immersed_orbit}):

\begin{corollary}\label{cor:immersed}
 Let $G$ be a Lie group, $\mathfrak{g}$ its Lie algebra, $w\in\mathfrak{g}$ a point, and $G_{w}$ the stabilizer of $w$. Then, the map $f_{w}:G/G_{w}\to\mathfrak{g}$, $[g]\mapsto \Ad (g)w$ is a well-defined, $G$-equivariant, and injective immersion whose image is the adjoint orbit $\mathcal{O}$ of $G$ through $w$.
\end{corollary}

\begin{remark}
 Even though there are some examples in which $\mathcal{O}\subset\mathfrak{g}$ is even an embedded submanifold, for instance if $G$ is compact, this is not always the case. Note, however, that every immersion is locally an embedding which is why this subtlety will not pose any problems when we discuss the integrability of $J$.
\end{remark}

Corollary \autoref{cor:immersed} allows us to identify the tangent spaces of $\mathcal{O}$ with subspaces\linebreak of $\mathfrak{g}$. Indeed, the map ($v\in\mathfrak{g}$)
\begin{gather*}
 \gamma_{v,w} (t)\coloneqq f_w ([\exp (tv)]) = \Ad (\exp (tv))w
\end{gather*}
is a smooth curve in $\mathcal{O}$ through $w$. Its derivative at $t=0$ is the fundamental vector field $X_v$:
\begin{gather*}
X_v (w)\coloneqq\dot\gamma_{v,w} (0) = \ad_v w = -\ad_w v.
\end{gather*}
Since all vectors tangent to $\mathcal{O}$ can be represented by such a curve, we find:
\begin{gather*}
 T_w\mathcal{O}\cong \im\ad_w\subset\mathfrak{g}\quad\forall w\in\mathcal{O}.
\end{gather*}
This gives us the following corollary:

\begin{corollary}[Global construction of $J$]\label{cor:can_sec}
 Let $G$ be a Lie group, $\mathfrak{g}$ its Lie algebra, and $\mathcal{O}\subset\mathfrak{g}$ a skew-symmetric adjoint orbit of $G$. Then, the section $J\in\Gamma\End (T\mathcal{O})$ defined pointwise by $J_w\in\End (T_w\mathcal{O})\cong\End (\im\ad_w)$ as in Proposition \autoref{prop:can_comp_str} satisfies $J^2_w = -\id_{T_w\mathcal{O}}$.
\end{corollary}

For $J$ to become an almost complex structure on $\mathcal{O}$, we need to show that the section $J$ is smooth. This is a consequence of the $G$-invariance of $J$:

\begin{proposition}[$J$ is an almost complex structure]\label{prop:J_smooth}
 Let $G$ be a Lie group, $\mathfrak{g}$ its Lie algebra, and $\mathcal{O}\subset\mathfrak{g}$ a skew-symmetric adjoint orbit of $G$. Then, the section $J\in\Gamma\End (T\mathcal{O})$ is $G$-invariant, i.e.:
 \begin{gather}
  (gJ)_w \coloneqq \Ad (g)\circ J_{\Ad (g^{-1})w} \circ \Ad (g^{-1}) = J_w\quad\forall g\in G\ \forall w\in\mathcal{O}.\label{eq:G-invariance}
 \end{gather}
 In particular, $J$ is a smooth almost complex structure on $\mathcal{O}$.
\end{proposition}

\begin{proof}
 First, we show that $J$ is $G$-invariant: The complex structure\linebreak $J_w\in\End (\im\ad_w)$ can be written as the composition $r_w\circ\ad_w$. Here,\linebreak $r_w\in\End (\im\ad_w)$ is defined by:
 \begin{gather*}
  r_w (v)\coloneqq\frac{v}{\mu}\quad\forall v\in E_\mu\ \forall\,\text{norms }\mu>0\text{ of eigenvalues},
 \end{gather*}
 where we employ the notations from Proposition \autoref{prop:can_comp_str}. The identification $T_w\mathcal{O}\cong\im\ad_w$ allows us to interpret $w\mapsto r_w$ and $w\mapsto \ad_w$ as sections of $\End (T\mathcal{O})$. If these sections are $G$-invariant, then $J$ is $G$-invariant as well. The $G$-invariance of $\ad_w$ is an immediate consequence of \autoref{eq:adj_rep}:
 \begin{gather*}
  \ad_{\Ad (g)w} = \Ad (g)\circ\ad_w\circ\Ad (g^{-1}).
 \end{gather*}
 The last equation implies that the decomposition
 \begin{gather*}
  T_w\mathcal{O}\cong\im\ad_w = \bigoplus_{\mu >0} E_\mu
 \end{gather*}
 is $G$-invariant. Therefore, $r_w$ is also $G$-invariant.\\
 Next, we demonstrate how the $G$-invariance of $J$ infers its smoothness: First, observe that, by Theorem \autoref{thm:godement} and Lemma \autoref{lem:immersed_orbit}, the map
 \begin{gather*}
  \ev_{w_0} = f_{w_0}\circ \pi: G\to\mathcal{O},\ \ev_{w_0} (g) = \Ad (g)w_0
 \end{gather*}
 is a surjective submersion, where $w_0\in\mathcal{O}$ is fixed. Hence, for every $w\in\mathcal{O}$, there exists an embedded submanifold $U_{w}\subset G$ such that $W_{w}\coloneqq\ev_{w_0} (U_{w})\subset\mathcal{O}$ is an open subset containing $w$ and $\ev_{w_0}\vert_{U_{w}}$ is a diffeomorphism onto its image $W_{w}$. We now rewrite \autoref{eq:G-invariance} as follows:
 \begin{gather*}
  J_{\Ad (g)w_0} = \Ad (g)\circ J_{w_0}\circ\Ad (g^{-1}).
 \end{gather*}
 Hence, we find with the help of $\ev_{w_0}$:
 \begin{gather*}
  J_{w^\prime} = \Ad (g(w^\prime))\circ J_{w_0}\circ\Ad (g(w^\prime)^{-1})\quad\forall w^\prime\in W_w\forall w\in\mathcal{O},
 \end{gather*}
 where $g(w^\prime)\coloneqq\ev_{w_0}\vert^{-1}_{U_w} (w^\prime)$. This equation proves that $J$ is smooth, because the adjoint action on $\mathcal{O}$ is smooth.
\end{proof}

\subsubsection*{Integrability of $\mathbf{J}$}

This part is devoted to the integrability of $J$. By the Newlander-Nirenberg theorem (cf. \cite{kobayashi1969}), it suffices to compute the Nijenhuis tensor $N_J$ to check the integrability of $J$. To perform this calculation, we extend $J$ to a local section of $\End (T\mathfrak{g})$ and use Proposition \autoref{prop:nijenhuis}.\\
To formulate Proposition \autoref{prop:nijenhuis}, we introduce the following notation: Let $U\subset\R^n$ be an open subset, $X = (x_1,\ldots, x_n)\in C^\infty (U,\R^n)$, and $Y\in\R^n$. We set:
\begin{gather*}
 dX (Y)\coloneqq (dx_1 (Y),\ldots, dx_n (Y)).
\end{gather*}
Similarly, we set for $A = (a_{ij})_{i,j = 1,\ldots, n}\in C^\infty (U,\R^{n\times n})$:
\begin{gather*}
 dA (Y)\coloneqq (da_{ij}(Y))_{i,j = 1,\ldots, n}.
\end{gather*}

\begin{proposition}\label{prop:nijenhuis}
 Let $U\subset\R^n$ be an open subset and $A$ a $(1,1)$-tensor\linebreak on $U$, i.e., $A\in C^\infty (U,\R^{n\times n})$. Then, the Nijenhuis tensor $N_A$ of $A$ is given by:
 \begin{gather*}
  N_A (X,Y) = A [X,Y]_{dA} - [X,Y]_{dA,A}\quad\forall X,Y\in\R^n,
 \end{gather*}
 where the brackets $[\cdot, \cdot]_{dA}$ and $[\cdot, \cdot]_{dA,A}$ are defined by ($X,Y\in\R^n$):
 \begin{gather*}
  [X,Y]_{dA} \coloneqq dA(X) Y - dA (Y) X,\quad [X,Y]_{dA,A} \coloneqq dA(AX) Y - dA (AY) X.
 \end{gather*}
\end{proposition}

\begin{proof}
 Let $p\in U$ be a point and $V, W$ vector fields on $U$ such that\linebreak $V(p)\equiv V_p = X$ and $W(p)\equiv W_p = Y$. Then, the Nijenhuis tensor $N_A$\linebreak at $p$ is defined as follows:
 \begin{gather}\label{eq:nijenhuis}
  N_{A,p} (X,Y) = -A^2_p [V,W]_p + A_p\left( [AV, W]_p + [V, AW]_p\right) - [AV, AW]_p.
 \end{gather}
 The Lie bracket of two vector fields on $U$ is given by:
 \begin{gather*}
  [V,W] = dW (V) - dV (W).
 \end{gather*}
 Hence, we obtain by applying the product rule:
 \begin{align*}
  [V,W]_p &= dW_p (X) - dV_p (Y),\\
  [AV,W]_p &= dW_p (A_p X) - A_p dV_p (Y) - dA_p (Y) X,\\
  [V,AW]_p &= A_p dW_p (X) - dV_p (A_p Y) + dA_p (X) Y,\\
  [AV,AW]_p &= dA_p (A_p X) Y + A_p dW_p \left(A_p X\right) - dA_p (A_p Y) X - A_p dV_p (A_p Y).
 \end{align*}
 Inserting these equations into Equation \eqref{eq:nijenhuis} finishes the proof.
\end{proof}

Proposition \autoref{prop:nijenhuis} allows us to compute $N_A$ of a $(1,1)$-tensor $A$ whose underlying manifold is a vector space. However, the almost complex structure $J$ is not defined on a vector space, but only on an immersed submanifold. We can circumvent this problem by locally extending $J$:

\begin{definition}[Local continuation]\label{def:local_continuation}
 Let $M,N$ be two smooth manifolds and $f:M\to N$ an immersion. Further, let $g$ be a function, $X$ a vector field, and $A$ a $(1,1)$-tensor on $M$. Then, $\hat g$ ($\hat X$, $\hat A$) is a \textbf{local continuation} of $g$ ($X$, $A$) with respect to the immersion $f:M\to N$ if there exist non-empty open subsets $U\subset M$ and $V\subset N$ such that $f(U)\subset V$ and\ldots
 \begin{itemize}
  \item \dots $\hat g$ is a function on $V$ satisfying $\hat g\circ f (p) = g(p)$ for all $p\in U$,
  \item \dots $\hat X$ is a vector field on $V$ satisfying $\hat X (f(p)) = df_p X(p)$ for all $p\in U$,
  \item \dots $\hat A$ is a $(1,1)$-tensor on $V$ satisfying $\hat A_{f(p)}\circ df_p = df_p\circ A_p$ for all $p\in U$.
 \end{itemize}
 Local continuations of general tensors are defined analogously.
\end{definition}

\begin{remark}[Local continuations always exist]
 Given any immersion\linebreak $f:M\to N$, tensor $T$ on $M$, and point $p\in M$, we can always find a local continuation $\hat T$ of $T$ where $U$ from Definition \autoref{def:local_continuation} is a neighborhood of $p$. Indeed, any immersion is locally an embedding, i.e., there exists an open neighborhood $U$ of $p$ such that $f\vert_U$ is an embedding. After shrinking $U$ if necessary, we can assume that $f(U)\subset N$ completely lies in one submanifold chart. It is now trivial to extend $T$ in this submanifold chart.
\end{remark}

A priori, we do not know how the Nijenhuis tensor of $A$ and its local continuation $\hat A$ are related. Intuitively, we expect $N_{\hat A}$ to be a local continuation of $N_A$. The next proposition confirms our intuition:

\begin{proposition}[$N_{\hat A}$ is local continuation of $N_A$]\label{prop:local_continuation_nijenhuis}
 Let $f:M\to N$ be an immersion between two manifolds. Further, let $A$ be a $(1,1)$-tensor on $M$ and $\hat A$ a local continuation of $A$ with respect to $f$. Then, $N_{\hat A}$ is a local continuation of $N_A$ with respect to $f$. In particular, $N_{A,p}$ vanishes if and only if $N_{\hat A, f(p)}$ vanishes on $df_p(T_pM)$.
\end{proposition}

\begin{proof}
 Proposition \autoref{prop:local_continuation_nijenhuis} immediately follows from the definition of $N_A$ (cf. \autoref{eq:nijenhuis}) and Proposition \autoref{prop:local_continuation_commutator}.
\end{proof}

\begin{proposition}\label{prop:local_continuation_commutator}
 Let $f:M\to N$ be an immersion between two manifolds. Further, let $V,W$ be vector fields on $M$ and $A$ a $(1,1)$-tensor on $M$ with local continuations $\hat V$, $\hat W$, and $\hat A$, respectively. Then, $\hat A \hat V$ is a local continuation of $AV$ and $[\hat V,\hat W]$ is a local continuation of $[V,W]$.
\end{proposition}

\begin{proof}
 It clearly follows from Definition \autoref{def:local_continuation} that $\hat A\hat V$ is a local continuation of $AV$. Now consider $[\hat V,\hat W]$. Since the commutator can be computed locally and every immersion is locally an embedding, we can assume without loss of generality that $M = \R^m$, $N = \R^n$, and that $f$ is given by the inclusion\linebreak $\R^m\hookrightarrow\R^n = \R^m\times\R^{n-m}$. If we denote the coordinates of $\R^n$ by $x_1,\ldots, x_n$, we can express the vector fields at hand as follows:
 \begin{alignat*}{2}
  V &= \sum^m_{i=1} v_i\pa{x_i},\quad W &&= \sum^m_{i=1} w_i\pa{x_i},\\
  \hat V &= \sum^n_{i=1} \hat v_i\pa{x_i},\quad \hat W &&= \sum^n_{i=1} \hat w_i\pa{x_i},
 \end{alignat*}
 where the functions $v_i,w_i\in C^\infty (\R^m)$ and $\hat v_i,\hat w_i\in C^\infty (\R^n)$ satisfy:
 \begin{align*}
  \hat v_i (x_1,\ldots, x_m, 0,\ldots, 0) &= \begin{cases} v_i (x_1,\ldots, x_m)&\text{ for }i\in\{1,\ldots, m\},\\ 0&\text{ else,}\end{cases}\\
  \hat w_i (x_1,\ldots, x_m, 0,\ldots, 0) &= \begin{cases} w_i (x_1,\ldots, x_m)&\text{ for }i\in\{1,\ldots, m\},\\ 0&\text{ else.}\end{cases}
 \end{align*}
 Thus, $[\hat V, \hat W]$ becomes:
 \begin{alignat*}{2}
  [\hat V,\hat W] &= \phantom{+} \sum^m_{i,j = 1} [\hat v_i\pa{x_i}, \hat w_j\pa{x_j}] &&+ \sum^m_{i=1}\sum^n_{j = m+1} [\hat v_i\pa{x_i}, \hat w_j\pa{x_j}]\\
  &\phantom{=}+ \sum^n_{i=m+1}\sum^m_{j = 1} [\hat v_i\pa{x_i}, \hat w_j\pa{x_j}] &&+ \sum^n_{i,j = m+1} [\hat v_i\pa{x_i}, \hat w_j\pa{x_j}].
 \end{alignat*}
 Evaluating the first term on $\R^m\times\{0\}$ gives:
 \begin{align*}
  \sum^m_{i,j = 1} [\hat v_i\pa{x_i}, \hat w_j\pa{x_j}] (x_1,\ldots, x_m,0,\ldots,0) &= \sum^m_{i,j = 1} [v_i\pa{x_i}, w_j\pa{x_j}] (x_1,\ldots, x_m)\\
  &= [V,W](x_1,\ldots, x_m).
 \end{align*}
 The remaining terms vanish, since $v_j$, $w_j$, $\pa{x_i}v_j$, and $\pa{x_i}w_j$ evaluated on\linebreak $\R^m\times\{0\}\subset\R^n$ are identically zero for $i\in\{1,\ldots, m\}$ and $j\in\{m+1,\ldots, n\}$. Therefore, we have:
 \begin{gather*}
  [\hat V,\hat W](x_1,\ldots, x_m,0,\ldots,0) = [V,W](x_1,\ldots, x_m),
 \end{gather*}
 concluding the proof.
\end{proof}

Let us return to the almost complex structure $J$. Our strategy is to compute $N_J$ by applying Proposition \autoref{prop:nijenhuis} to a local continuation $\hat J$. To do that, we need to determine the differential $d\hat J$. In general, this might be tricky, however, it becomes much simpler if we chose $\hat J$ to be $G$-invariant (cf. \autoref{eq:G-invariance}):

\begin{proposition}\label{prop:dJ}
 Let $G$ be a Lie group with Lie algebra $\mathfrak{g}$ and skew-symmetric adjoint orbit $\mathcal{O}$. Further, let $\hat J$ be a local, $G$-invariant continuation of $J$ with respect to the immersion $\mathcal{O}\subset\mathfrak{g}$. For $v\in\mathfrak{g}$, define the associated fundamental vector field as follows:
 \begin{gather*}
  X_v (w)\equiv X_v\coloneqq -\ad_w v\in \im\ad_w\cong T_w\mathcal{O}\quad\forall w\in\mathcal{O}.
 \end{gather*}
 Then, the differential $d\hat J$ is given by:
 \begin{gather*}
  d\hat J_w (X_u) X_v = [u, \hat J_w X_v] - \hat J_w [u, X_v]\quad\forall u,v\in\mathfrak{g}\ \forall w\in\mathcal{O}.
 \end{gather*}
\end{proposition}

\begin{proof}
 For $w\in\mathcal{O}$ and $u\in\mathfrak{g}$, define the curve:
 \begin{gather*}
  \gamma_{u,w} (t)\coloneqq \Ad (\exp (tu))w.
 \end{gather*}
 Its derivative at $t=0$ is $\dot\gamma_{u,w} (0) = X_u$.\\
 Now consider $d\hat J$. The $G$-invariance of $\hat J$ implies ($v\in\mathfrak{g}$):
 \begin{align*}
  d\hat J_w (X_u) X_v &= \left.\frac{d}{dt}\right\vert_{t=0} \hat J_{\gamma_{u,w}(t)} X_v\\
  &= \left.\frac{d}{dt}\right\vert_{t=0}\left\{\Ad (\exp (tu))\circ \hat J_w\circ\Ad (\exp (-tu)) (X_v)\right\}\\
  &= [u, \hat J_w X_v] - \hat J_w [u, X_v].
 \end{align*}
\end{proof}

\begin{remark}\vspace{-0.3cm}\ 
 \begin{enumerate}
  \item $J$ possesses a local, $G$-invariant continuation $\hat J$ near any point $w\in\mathcal{O}$, because $J$ itself is $G$-invariant (cf. Proposition \autoref{prop:J_smooth}). $\hat J$ can be constructed using the orbit slice theorem (cf. Theorem I.2.1 in \cite{Audin2004}).
  \item Of course, the differential $d\hat J$ depends on the chosen continuation $\hat J$. However, the Nijenhuis tensor $N_{\hat J, w}\vert_{T_w\mathcal{O}\times T_w\mathcal{O}}$ obtained from $d\hat J$ via Proposition \autoref{prop:nijenhuis} is independent of the choice of $\hat J$ (cf. Proposition \autoref{prop:local_continuation_nijenhuis}). In light of this observation, we drop the hat symbol from now on and denote both $J$ and $\hat J$ simply by $J$. It is clear from the context whether we mean the almost complex structure or the local continuation.
 \end{enumerate}
\end{remark}

Let us now use Proposition \autoref{prop:nijenhuis} and \autoref{prop:dJ} to compute the Nijenhuis tensor. Since we have the decomposition $T_w\mathcal{O}\cong\im\ad_w = \bigoplus_{\mu>0} E_\mu$, it suffices to evaluate $N_J$ on the spaces $E_\mu$. We find:

\begin{proposition}\label{prop:N_J}
 Let $G$ be a Lie group with Lie algebra $\mathfrak{g}$ and skew-symmetric adjoint orbit $\mathcal{O}\subset\mathfrak{g}$. Then, the Nijenhuis tensor $N_J$ of the almost complex structure $J$ of $\mathcal{O}$ is given by:
 \begin{gather*}
  N_{J,w} (X_u, X_v) = (\lambda + \mu)\left(J_w [u,v] - [J_wu, v] - [u, J_wv] - J_w [J_w u, J_w v]\right),
 \end{gather*}
 where $w\in\mathcal{O}$, $u\in E_\lambda$, and $v\in E_\mu$.
\end{proposition}

\begin{proof}
 By Proposition \autoref{prop:nijenhuis}, $N_{J,w}$ amounts to:
 \begin{align*}
  N_{J,w} (X_u, X_v) &= J_w\left\{dJ_w (X_u) X_v - dJ_w (X_v) X_u\right\}\\
  &\phantom{=} - dJ_w (J_w X_u) X_v + dJ_w (J_w X_v)X_u.
 \end{align*}
 With the help of Proposition \autoref{prop:dJ}, we evaluate the differential $dJ_w (X_u) X_v$:
 \begin{gather*}
  dJ_w (X_u) X_v = [u, J_w X_v] - J_w [u, X_v]
 \end{gather*}
 A similar expression holds for $dJ_w (X_v) X_u$. To compute the remaining differentials in $N_{J,w}$, we need to rewrite $J_w X_u$ and $J_w X_v$ in the form $X_Y$ for some vectors $Y\in\mathfrak{g}$. By definition of $J_w$, we have:
 \begin{gather*}
  J_w u = \frac{1}{\lambda}\ad_w u,\quad J_w v = \frac{1}{\mu}\ad_w v.
 \end{gather*}
 This implies:
 \begin{alignat*}{6}
  X_u &= -\ad_w u &&= -\lambda J_w u\ \Rightarrow\ && X_{J_wu} &= -\lambda J^2_wu &&= \lambda u &&= J_w X_u,\\
  X_v &= -\ad_w v &&= -\mu J_w v\ \Rightarrow\ && X_{J_wv} &= -\mu J^2_wv &&= \mu v &&= J_w X_v.
 \end{alignat*}
 Inserting these identities into the differentials yields:
 \begin{align*}
  dJ_w (X_u) X_v &= \mu\left([u,v] + J_w[u,J_wv]\right),\\
  -dJ_w (X_v) X_u &= \lambda\left([u,v] + J_w[J_w u,v]\right),\\
  -dJ_w (J_w X_u) X_v &= -\mu\left([J_w u, v] + J_w [J_w u, J_w v]\right),\\
  dJ_w (J_w X_v) X_u &= -\lambda\left([u, J_w v] + J_w [J_w u, J_w v]\right).
 \end{align*}
 Lastly, we put everything together to find:
 \begin{gather*}
  N_{J,w} (X_u, X_v) = (\lambda + \mu)\left(J_w [u,v] - [J_wu, v] - [u, J_wv] - J_w [J_w u, J_w v]\right),
 \end{gather*}
 where we used $J^2_w s = -s$ for $s\coloneqq \mu [u, J_w v] + \lambda [J_w u, v]$. This equation is valid, since $s$ is an element of $\im\ad_w$, as one application of \autoref{eq:jacobi} shows:
 \begin{gather*}
  \mu [u, J_w v] + \lambda [J_w u, v] = [u,\ad_w v] + [\ad_w u, v] = \ad_w [u,v].
 \end{gather*}
\end{proof}

Proposition \autoref{prop:N_J} implies that $J$ is integrable. This can be seen by complexifying the Nijenhuis tensor:

\begin{proposition}\label{prop:J_integrable}
 Let $G$ be a Lie group with Lie algebra $\mathfrak{g}$ and skew-symmetric adjoint orbit $\mathcal{O}\subset\mathfrak{g}$. Then, the Nijenhuis tensor $N_J$ of the almost complex structure $J$ of $\mathcal{O}$ vanishes implying that $J$ is integrable.
\end{proposition}

\begin{proof}
 Complexifying all maps that occur during the computation of $N_J$ ($\ad_w$, $J_w$, $[\cdot, \cdot]$, and so on) allows us to apply Proposition \autoref{prop:N_J} not only to vectors $u\in E_\lambda$ and $v\in E_\mu$, but also to vectors $u\in E_{\lambda,\C}$ and $v\in E_{\mu,\C}$. By Remark \autoref{rem:J_w}, these spaces admit the following decomposition:
 \begin{gather*}
  E_{\lambda,\C} = E^{(1,0)}_\lambda\oplus E^{(0,1)}_\lambda = E_{i\lambda}\oplus E_{-i\lambda},\quad E_{\mu,\C} = E^{(1,0)}_\mu\oplus E^{(0,1)}_\mu = E_{i\mu}\oplus E_{-i\mu}.
 \end{gather*}
 Now let $u^{(1,0)}\in E_{i\lambda}$, $u^{(0,1)}\in E_{-i\lambda}$, $v^{(1,0)}\in E_{i\mu}$, and $v^{(0,1)}\in E_{-i\mu}$ be any vectors. Exploiting $J_w u^{(1,0)} = iu^{(1,0)}$, $J_w u^{(0,1)} = -iu^{(0,1)}$, and similar equations for $v$ allows us to evaluate the formula in Proposition \autoref{prop:N_J}:
 \begin{alignat*}{3}
  N_{J, w} (X_{u^{(1,0)}}, X_{v^{(1,0)}}) &= 2(\lambda + \mu) (J_w - i)[u^{(1,0)}, v^{(1,0)}],\ && N_{J, w} (X_{u^{(0,1)}}, X_{v^{(1,0)}}) &&= 0,\\
  N_{J, w} (X_{u^{(0,1)}}, X_{v^{(0,1)}}) &= 2(\lambda + \mu) (J_w + i)[u^{(0,1)}, v^{(0,1)}],\ && N_{J, w} (X_{u^{(1,0)}}, X_{v^{(0,1)}}) &&= 0.
 \end{alignat*}
 It remains to be shown:
 \begin{gather*}
  J_w[u^{(1,0)}, v^{(1,0)}] = i[u^{(1,0)}, v^{(1,0)}],\quad J_w[u^{(0,1)}, v^{(0,1)}] = -i[u^{(0,1)}, v^{(0,1)}].
 \end{gather*}
 However, this follows immediately from \autoref{eq:jacobi}:
 \begin{align*}
  \ad_w [u^{(1,0)}, v^{(1,0)}] &= [\ad_w u^{(1,0)}, v^{(1,0)}] + [u^{(1,0)}, \ad_w v^{(1,0)}]\\
  &= i(\lambda + \mu)[u^{(1,0)}, v^{(1,0)}],\\
  \ad_w [u^{(0,1)}, v^{(0,1)}] &= [\ad_w u^{(0,1)}, v^{(0,1)}] + [u^{(0,1)}, \ad_w v^{(0,1)}]\\
  &= -i(\lambda + \mu)[u^{(0,1)}, v^{(0,1)}].
 \end{align*}
\end{proof}

To conclude the discussion about $J$, we collect our results in the following lemma (cf. Section B in Chapter 8 of \cite{Besse2007}, in particular Proposition 8.39):

\begin{lemma}[Canonical complex structure]\label{lem:can_comp_str}
 Let $G$ be a Lie group with Lie algebra $\mathfrak{g}$ and let $\mathcal{O}\subset\mathfrak{g}$ be a skew-symmetric adjoint orbit of $G$. Then, $\mathcal{O}\subset\mathfrak{g}$ is an immersed submanifold and carries a canonical $G$-invariant complex\linebreak structure $J$.
\end{lemma}

\subsection*{Symplectic Structure of Coadjoint Orbits}

Let us now recall the construction of the Kirillov-Kostant-Souriau form $\kks$ on coadjoint orbits induced by the natural Poisson structure on $\mathfrak{g}^\ast$. Since these concepts are common knowledge among symplectic geometers, we keep the review brief. For a detailed discussion of Poisson manifolds and the KKS form, confer any textbook on symplectic geometry, for instance \cite{Marsden1999}.\\
We begin by reviewing basic facts about Poisson manifolds:

\begin{definition}\label{def:Poisson}
 A pair $(M,\{\cdot,\cdot\})$ is called \textbf{Poisson manifold} if $M$ is a smooth manifold and $\{\cdot,\cdot\}:C^\infty (M)\times C^\infty (M)\to C^\infty (M)$ is a bilinear map satisfying:
 \begin{enumerate}
  \item $\{\cdot,\cdot\}$ is skew-symmetric, i.e., $\{F,G\} = -\{G,F\}$ for $F,G\in C^\infty (M)$.
  \item $\{\cdot,\cdot\}$ is a derivation:
  \begin{gather*}
   \{FG, H\} = F\{G,H\} + G\{F,H\}\quad\forall F,G,H\in C^\infty (M).
  \end{gather*}
  \item $\{\cdot, \cdot\}$ fulfills the Jacobi identity:
  \begin{gather*}
   \{\{F, G\}, H\} + \{\{H, F\}, G\} + \{\{G, H\}, F\} = 0\quad\forall F,G,H\in C^\infty (M).
  \end{gather*}
 \end{enumerate}
\end{definition}

Condition (i) and (ii) allow us to rewrite any Poisson bracket $\{\cdot, \cdot\}$ in terms of a skew-symmetric bivector field $\pi\in\Gamma \left(\Lambda^2 TM\right)$, namely:
\begin{gather}\label{eq:bivector}
 \{F, G\} \eqqcolon \pi (dF, dG)\quad\forall F,G\in C^\infty (M).
\end{gather}
As every cotangent vector is the differential of a function at some point, the bivector field $\pi$ is unique. Observe that, conversely, every bracket $\{\cdot, \cdot\}$ defined via a skew-symmetric bivector field $\pi$ as in Equation \eqref{eq:bivector} automatically satisfies Condition (i) and (ii) from Definition \autoref{def:Poisson}. Condition (iii) is satisfied if and only if the Schouten bracket of $\pi$ with itself vanishes, i.e., $[\pi, \pi] = 0$.\\
Similar to symplectic manifolds, Poisson manifolds $(M,\{\cdot,\cdot\})$ come with a musical map $\# : T^\ast M\to TM$. It is defined by $\# (\alpha)\coloneqq \iota_\alpha\pi$. However, $\#$ is usually not an isomorphism. The image of $\#$ is a distribution on $M$ whose rank is, in general, non-constant. By Condition (iii), the distribution $\im\#$ is involutive. Thus, the leaves of $\im\#$ are immersed submanifolds of $M$ due to the generalized Frobenius theorem. They also carry a symplectic structure. To see this, let $L\subset M$ be a leaf of $\im\#$ and $p\in L$ be a point. The tangent space $T_p L$ is equal to $\im \#_p$, hence, we can define $\omega\in\Omega^2 (L)$ by:
\begin{gather*}
 \omega_p (\#_p (\alpha), \#_p (\beta))\coloneqq \pi_p (\alpha,\beta)\quad\forall \alpha, \beta \in T^\ast_p M\ \forall p\in M.
\end{gather*}
$\omega$ is well-defined, since we have $\alpha-\alpha^\prime\in\ker\#_p$ if $\#_p (\alpha) = \#_p (\alpha^\prime)$ which implies $\pi_p (\alpha,\cdot) = \pi_p (\alpha^\prime, \cdot)$. By construction, $\omega$ is non-degenerate. Also, one easily verifies that $\omega$ is smooth. Furthermore, a straightforward, but tedious computation reveals that $\omega$ is closed proving that $\omega$ is a symplectic form on $L$. This gives us the following proposition:

\begin{proposition}\label{prop:symp_leaves}
 Let $(M,\{\cdot,\cdot\})$ be a Poisson manifold with induced distribution $\im\#$. Then, $\im\#$ is involutive and the leaves of $\im\#$ are immersed symplectic submanifolds of $M$.
\end{proposition}

Note that the distribution $\im\#$ coincides with the entire tangent bundle $TM$ if $\#$ is an isomorphism. In this case, the leaves of $\im\#$ are the connected components of $M$ and the Poisson bracket $\{\cdot,\cdot\}$ agrees with the Poisson bracket induced by the symplectic structure on $M$. In this sense, we can interpret Poisson manifolds as a generalization of symplectic manifolds.\\
We now turn our attention to the most classical example of Poisson manifolds:

\begin{proposition}\label{prop:dual_lie}
 Let $(\mathfrak{g}, [\cdot,\cdot])$ be a real Lie algebra and $\mathfrak{g}^\ast$ its dual. Then, $\mathfrak{g}^\ast$ carries a canonical Poisson bracket defined by:
 \begin{gather*}
  \{F,G\} (\alpha)\coloneqq \alpha \left(\left[dF_\alpha, dG_\alpha\right]\right),
 \end{gather*}
 where $F,G\in C^\infty (\mathfrak{g}^\ast)$, $\alpha\in\mathfrak{g}^\ast$, and we use the canonical identifications\linebreak $T^\ast_\alpha\mathfrak{g}^\ast\cong\mathfrak{g}^{\ast\ast}\cong \mathfrak{g}$.
\end{proposition}

\begin{proof}
 Condition (i) and (ii) of Definition \autoref{def:Poisson} are easily checked. To verify Condition (iii), we first note that the Jacobi identity of the Lie bracket $[\cdot, \cdot]$ implies the Jacobi identity of the Poisson bracket $\{\cdot, \cdot\}$ for linear functions on $\mathfrak{g}^\ast$. Together with Condition (ii), this shows that the Jacobi identity also holds for polynomials on $\mathfrak{g}^\ast$. Now observe that the Jacobi identity only involves nested Poisson brackets of the form $\{\{F,G\}, H\}$ and that such terms only contain derivatives up to second order. Thus, we can replace the functions on $\mathfrak{g}^\ast$ by their second-order Taylor expansion to verify the Jacobi identity. Since these are just polynomials, the Jacobi identity is fulfilled for all functions on $\mathfrak{g}^\ast$ concluding the proof.
\end{proof}

Our next task is to determine the symplectic leaves of the distribution induced by the Poisson structure on $\mathfrak{g}^\ast$. To do so, we first need to compute the musical\linebreak map $\#$. As $\{F,G\}(\alpha)$ only depends on the differentials $dF_\alpha$ and $dG_\alpha$, it suffices to evaluate the Poisson bracket only for linear functions on $\mathfrak{g}^\ast$ in order to calculate $\#$. For $v\in\mathfrak{g}$, we define the linear function $F_v\in C^\infty (\mathfrak{g}^\ast)$ by $F_v (\alpha)\coloneqq \alpha (v)$. Then:
\begin{gather*}
 \{F_v, F_w\} (\alpha) = \alpha ([v,w])\quad\forall v,w\in\mathfrak{g}\ \forall\alpha\in\mathfrak{g}^\ast.
\end{gather*}
Hence, the musical map $\#_\alpha :\mathfrak{g}\cong T^\ast_\alpha\mathfrak{g}^\ast\to\mathfrak{g}^\ast\cong T_\alpha\mathfrak{g}^\ast$ amounts to:
\begin{gather*}
 \#_\alpha (v) = \alpha\circ\ad_v \eqqcolon -\ad^\ast_v (\alpha) = \ad^\ast_{-v} (\alpha)\quad\forall v\in\mathfrak{g}.
\end{gather*}
Now assume that $\mathfrak{g}$ is the Lie algebra of a Lie group $G$. The Lie group $G$ acts on its dual Lie algebra $\mathfrak{g}^\ast$ by the coadjoint representation, i.e.:
\begin{gather*}
 \Ad^\ast (g)\alpha \coloneqq \alpha\circ\Ad (g^{-1})\quad\forall g\in G\,\forall\alpha\in\mathfrak{g}^\ast.
\end{gather*}
We now associate a fundamental vector field $X^\ast_v$ on $\mathfrak{g}^\ast$ to every vector $v\in\mathfrak{g}$:
\begin{gather*}
 X^\ast_v (\alpha)\coloneqq\left.\frac{d}{dt}\right\vert_{t=0}\Ad^\ast (\exp (tv))\alpha = -\alpha\circ\ad_v = \ad^\ast_v (\alpha)\quad\forall\alpha\in\mathfrak{g}^\ast.
\end{gather*}
Similar to the adjoint case, the fundamental vector fields $X^\ast_v$ span the tangent spaces of the coadjoint orbits of $G$. Moreover, the last equations reveal that $\# (v)$ is the fundamental vector field $X^\ast_{-v}$. Hence, the distribution $\im\#$ consists of the tangent spaces of the coadjoint orbits and its symplectic leaves are the connected components of the coadjoint orbits. The symplectic form on the leaves is called the \textbf{Kirillov-Kostant-Souriau} form $\kks$. On a coadjoint orbit $\mathcal{O}^\ast\subset\mathfrak{g}^\ast$, it is given by:
\begin{gather}
 \omega_{\text{KKS},\alpha} (X^\ast_v (\alpha), X^\ast_w (\alpha)) = \alpha ([v,w])\quad\forall\alpha\in\mathcal{O}^\ast\, \forall v,w\in\mathfrak{g},\label{eq:kks_form}
\end{gather}
where we identify $T_\alpha\mathcal{O}^\ast$ with the span of the fundamental vector fields.\\
Before we conclude this subsection, we quickly observe that $\kks$ is $G$-invariant:
\begin{align*}
 &\phantom{=}\,\ \omega_{\text{KKS}, \Ad^\ast (g)\alpha} \left(\Ad^\ast (g) (\ad^\ast_v \alpha), \Ad^\ast (g) (\ad^\ast_w \alpha)\right)\\
 &= \omega_{\text{KKS}, \Ad^\ast (g)\alpha} \left(\ad^\ast_{\Ad (g)v} (\Ad^\ast (g) \alpha), \ad^\ast_{\Ad (g)w} (\Ad^\ast (g) \alpha)\right)\\
 &= (\Ad^\ast (g) \alpha) \left([\Ad (g)v, \Ad (g)w]\right) = \alpha\left(\Ad (g^{-1}) [\Ad (g)v, \Ad (g)w]\right)\\
 &= \alpha ([v,w]) = \omega_{\text{KKS},\alpha} (\ad^\ast_v \alpha, \ad^\ast_w \alpha),
\end{align*}
where we used \autoref{eq:adj_rep} to get from line $1$ to $2$ and from line $3$ to $4$.\\
The following lemma summarizes our findings:

\begin{lemma}[Canonical symplectic structure]\label{lem:kks}
 Let $G$ be a Lie group with dual Lie algebra $\mathfrak{g}^\ast$ and let $\mathcal{O}^\ast\subset\mathfrak{g}^\ast$ be a coadjoint orbit of $G$. Then,\linebreak $\mathcal{O}^\ast\subset\mathfrak{g}^\ast$ is an immersed submanifold and carries a canonical $G$-invariant symplectic structure given by the Kirillov-Kostant-Souriau form $\kks$.
\end{lemma}

\subsection*{Semi-Kähler Structure of (Co)Adjoint Orbits}

The main goal of \autoref{sec:semi-kaehler} and especially this subsection is to construct semi-Kähler structures on (co)adjoint orbits. Simply put, a semi-Kähler manifold is a Kähler manifold whose Kähler metric does not need to be positive definite (cf. \autoref{app:kaehler} for details):

\begin{definition}[Semi-Kähler manifolds]\label{def:kaehler_in_main}
 A \textbf{pre-semi-Kähler} manifold is a triple $(M,\omega, J)$ where $M^{2n}$ is a smooth manifold and the tensors $\omega\in\Omega^2 (M)$ and $J\in\Gamma\End (TM)$ satisfy:
 \begin{enumerate}
  \item $\omega$ is non-degenerate, i.e., $\omega^n_p\neq 0$ for all $p\in M$,
  \item $J$ is an almost complex structure, i.e., $J^2_p = -\id_{T_pM}$ for all $p\in M$,
  \item $\omega$ and $J$ are compatible in the sense that $\omega (J\cdot , J\cdot) = \omega$.
 \end{enumerate}
 We drop the prefix ``pre'' if $\omega$ is closed and $J$ is integrable, i.e., if $(M,\omega, J)$ satisfies the integrability conditions $d\omega = 0$ and $N_J = 0$.\\
 We drop the prefix ``semi'' if the semi-Riemannian metric $g\coloneqq \omega (\cdot, J\cdot)$ is positive definite.
\end{definition}

We are now ready to formulate the central theorem of \autoref{sec:semi-kaehler}:

\begin{theorem}[Semi-Kähler structures on (co)adjoint orbits]\label{thm:semi-kaehler}
 Let $G$ be a Lie group with Lie algebra $\mathfrak{g}$, dual Lie algebra $\mathfrak{g}^\ast$, and $\Ad$-invariant, non-degenerate scalar product $\sk{\cdot}{\cdot}$ on $\mathfrak{g}$. Further, let $\mathcal{O}\subset\mathfrak{g}$ be a skew-symmetric adjoint orbit of $G$. Then, $\mathcal{O}\subset\mathfrak{g}$ is an immersed submanifold and carries a $G$-invariant semi-Kähler structure. Its complex structure $J$ is the canonical complex structure on skew-symmetric adjoint orbits. If we identify $\mathcal{O}$ via $\sk{\cdot}{\cdot}$ with the coadjoint orbit $\mathcal{O}^\ast\subset\mathfrak{g}^\ast$, its symplectic form becomes the Kirillov-Kostant-Souriau form $\kks$. Moreover, $\mathcal{O}$ is Kähler if $\skcdot$ is positive definite.
\end{theorem}

We have already investigated the complex structure on skew-symmetric adjoint orbits and the symplectic structure on coadjoint orbits. The only missing ingredient for Theorem \autoref{thm:semi-kaehler} is the relation between adjoint and coadjoint orbits. In our setup, we link adjoint and coadjoint orbits by an $\Ad$-invariant, non-degenerate scalar product $\skcdot$ (cf. \autoref{sec:lie_groups}, in particular Definition \autoref{def:bi-invariant} and \autoref{eq:Ad-invariance}):

\begin{proposition}[$\mathcal{O}\cong\mathcal{O}^\ast$ via $\skcdot$]\label{prop:ad_coad}
 Let $G$ be a Lie group with Lie\linebreak algebra $\mathfrak{g}$ and dual Lie algebra $\mathfrak{g}^\ast$. Further, let $\sk{\cdot}{\cdot}:\mathfrak{g}\times\mathfrak{g}\to\R$ be an $\Ad$-invariant, non-degenerate scalar product, i.e.:
 \begin{gather*}
  \sk{\Ad (g)v}{\Ad (g)w} = \sk{v}{w}\quad\forall g\in G\, \forall v,w\in\mathfrak{g}.
 \end{gather*}
 Then, the map $b:\mathfrak{g}\to\mathfrak{g}^\ast$, $w\mapsto \sk{w}{\cdot}$ restricts to an $\Ad$-$\Ad^\ast$-equivariant diffeomorphism from adjoint to coadjoint orbits.
\end{proposition}

\begin{proof}
 Clearly, the map $b:\mathfrak{g}\to\mathfrak{g}^\ast$ is a diffeomorphism. $\Ad$-invariance of $\sk{\cdot}{\cdot}$ implies $\Ad$-$\Ad^\ast$-equivariance of $b$:
 \begin{gather*}
  b\circ\Ad (g) = \Ad^\ast (g)\circ b\quad\forall g\in G.
 \end{gather*}
 Now let $\mathcal{O}\subset\mathfrak{g}$ be an adjoint orbit of $G$ and $w\in\mathcal{O}$ a point. Denote the coadjoint orbit through $b(w)$ by $\mathcal{O}^\ast\subset\mathfrak{g}^\ast$. Since $b$ is $\Ad$-$\Ad^\ast$-equivariant, we have $b(\mathcal{O})\subset \mathcal{O}^\ast$ and $b^{-1}(\mathcal{O}^\ast)\subset \mathcal{O}$. This implies $b (\mathcal{O}) = \mathcal{O}^\ast$. As $\mathcal{O}\subset\mathfrak{g}$ and $\mathcal{O}^\ast\subset\mathfrak{g}^\ast$ are immersed submanifolds, the maps $b\vert_{\mathcal{O}}:\mathcal{O}\to\mathcal{O}^\ast$ and $b^{-1}\vert_{\mathcal{O}^\ast}:\mathcal{O}^\ast\to\mathcal{O}$ are smooth and, hence, diffeomorphisms.
\end{proof}

\begin{remark}\label{rem:Ad-invariant_metric}
 The $\Ad$-invariance of $\sk{\cdot}{\cdot}$ directly infers the skew-symmetry of $\ad_w$ with respect to $\sk{\cdot}{\cdot}$ (cf. \autoref{eq:ad-invariant}):
 \begin{gather*}
  \sk{\ad_w u}{v} = -\sk{u}{\ad_w v}\quad\forall u,v,w\in\mathfrak{g}.
 \end{gather*}
 A consequence of this observation is that every element $w\in\mathfrak{g}$ is skew-symmetric if $\skcdot$ is positive definite (cf. the beginning of \autoref{sec:semi-kaehler}, in particular Definition \autoref{def:skew}). However, this relation breaks down for indefinite scalar products. As a counterexample, consider the Minkowskian metric $\sk{\cdot}{\cdot}$ on $\R^2$ given by $\sk{e_j}{e_k} = (\delta_{1j}-\delta_{2j})\delta_{jk}$ in standard basis and the matrix $A = \begin{pmatrix}0 & 1\\ 1 & 0\end{pmatrix}$. $A$ is skew-symmetric with respect to $\sk{\cdot}{\cdot}$, nevertheless, the eigenvalues of $A$ are $\pm 1$ and, thus, real.
\end{remark}

To obtain a semi-Kähler structure, we can use $\sk{\cdot}{\cdot}$ to either transfer $J$ from $\mathcal{O}$ to $\mathcal{O}^\ast$ or pull back $\kks$ on $\mathcal{O}^\ast$ to $\mathcal{O}$. The two semi-Kähler manifolds found this way are isomorphic which is why we focus our attention on the adjoint orbit $\mathcal{O}$ for the remainder of \autoref{sec:semi-kaehler}. Pulling back $\kks$ to $\mathcal{O}$ yields the form $\omega$:

\begin{proposition}\label{prop:adjoint_omega}
 Let $G$ be a Lie group with Lie algebra $\mathfrak{g}$ and $\Ad$-invariant, non-degenerate scalar product $\sk{\cdot}{\cdot}$ on $\mathfrak{g}$. Further, let $\mathcal{O}\subset\mathfrak{g}$ be an adjoint orbit of $G$. Then, the symplectic form $\omega \coloneqq b^\ast\kks$ on $\mathcal{O}$ is given by
 \begin{gather*}
  \omega_w (X_u(w) , X_v(w) ) = \sk{w}{[u,v]}\quad\forall w\in\mathcal{O}\, \forall u,v\in\mathfrak{g},
 \end{gather*}
 where $b$ is the musical map induced by $\sk{\cdot}{\cdot}$ and $X_u(w) \coloneqq \ad_u w$.
\end{proposition}

\begin{proof}
 The skew-symmetry of $\ad_u$ implies:
 \begin{gather*}
  b(X_u (w)) = \sk{\ad_u w}{\cdot} = -\sk{w}{\ad_u\cdot} = \ad^\ast_u b(w) = X^\ast_u (b(w))\quad\forall u,w\in\mathfrak{g}.
 \end{gather*}
 Thus, we obtain:
 \begin{align*}
  \omega_w (X_u(w) , X_v(w) ) &= \omega_{\text{KKS}, b(w)} (b(X_u (w)), b(X_v (w)))\\
  &= \omega_{\text{KKS}, b(w)} (X^\ast_u (b(w)), X^\ast_v (b(w)))\\
  &= b(w)\left([u,v]\right) = \sk{w}{[u,v]}.
 \end{align*}
\end{proof}

We now have all tools at hand to prove Theorem \autoref{thm:semi-kaehler}:

\begin{proof}[Proof of Theorem \autoref{thm:semi-kaehler}]
 We only need to show the $\omega$-compatibility of $J$ and, for positive definite $\skcdot$, the positive definiteness of $g\coloneqq \omega (\cdot, J\cdot)$.\\
 $\omega = \omega (J\cdot,J\cdot)$ is equivalent to $\omega (J\cdot,\cdot) = -\omega (\cdot, J\cdot)$. To verify the last equation, we evaluate $\omega_w (J_wX_u, X_v)$. Since $w\in\mathcal{O}$ is skew-symmetric, $\mathfrak{g}$ decomposes into $\ker\ad_w$ and the spaces $E_\mu$ for $\mu>0$. If $u$ or $v$ is an element of $\ker\ad_w$, we have:
 \begin{gather*}
  \omega_w (J_wX_u, X_v) = 0 = -\omega_w (X_u, J_wX_v).
 \end{gather*}
 Let us now assume $u\in E_\lambda$ and $v\in E_\mu$. As in the proof of Proposition \autoref{prop:N_J}, we have $J_w X_u = X_{J_wu}$, thus:
 \begin{gather}\label{eq:compatible}
  \omega_w (J_w X_u, X_v) = \omega_w (X_{J_wu}, X_v) = \sk{w}{[J_wu, v]}.
 \end{gather}
 Next, consider the expression:
 \begin{gather*}
  s_{u,v}\coloneqq [J_wu, v] + [u,J_wv].
 \end{gather*}
 We want to show that $s_{u,v}$ lies in $\im\ad_w$. To achieve that, we distinguish two cases: If $\lambda = \mu$, we find:
 \begin{gather*}
  s_{u,v} = \frac{1}{\mu}\left([\ad_wu, v] + [u,\ad_wv]\right) = \frac{1}{\mu}\ad_w[u,v],
 \end{gather*}
 where we exploited \autoref{eq:jacobi}. If $\lambda\neq\mu$, we first note that it suffices to prove $[E_\lambda, E_\mu]\subset\im\ad_w$ in order to show $s_{u,v}\in\im\ad_w$, since $J_w (E_\lambda)\subset E_\lambda$ and $J_w (E_\mu)\subset E_\mu$. To show $[E_\lambda, E_\mu]\subset\im\ad_w$, we recall the decomposition:
 \begin{gather*}
  E_{\lambda,\C} = E_{i\lambda}\oplus E_{-i\lambda},\quad E_{\mu,\C} = E_{i\mu}\oplus E_{-i\mu}.
 \end{gather*}
 Exploiting \autoref{eq:jacobi} now yields:
 \begin{gather*}
  \ad_w [u,v] = [\ad_wu,v] + [u,\ad_w v] = i (\lambda+\mu) [u,v]\quad\forall u\in E_{i\lambda}\, \forall v\in E_{i\mu}
 \end{gather*}
 implying that $[E_{i\lambda}, E_{i\mu}]$ is zero or contained in the eigenspace $E_{i(\lambda+\mu)}$. In both cases, $[E_{i\lambda}, E_{i\mu}]$ is a subset of $\im\ad_w$. Similarly, one can show that $[E_{-i\lambda}, E_{i\mu}]$, $[E_{i\lambda}, E_{-i\mu}]$, and $[E_{-i\lambda}, E_{-i\mu}]$ are subspaces of $\im\ad_w$. This proves\linebreak $[E_\lambda, E_\mu]\subset\im\ad_w$ and, therefore, $s_{u,v}\in\im\ad_w$.\\
 We can infer from $s_{u,v}\in\im\ad_w$ that $w$ and $s_{u,v}$ are orthogonal:
 \begin{gather*}
  \sk{w}{\ad_w \xi} = -\sk{\ad_w w}{\xi} = 0\quad\forall \xi\in\mathfrak{g},
 \end{gather*}
 where we employed \autoref{eq:ad-invariant}. $\sk{w}{s_{u,v}} = 0$ implies:
 \begin{gather*}
  \sk{w}{[J_wu, v]} = -\sk{w}{[u, J_wv]}
 \end{gather*}
 Combining the last equation with \autoref{eq:compatible} yields:
 \begin{align*}
  \omega_w (J_w X_u, X_v) &= \sk{w}{[J_wu,v]} = -\sk{w}{[u,J_wv]}\\
  &= -\omega_w (X_u, J_w X_v)
 \end{align*}
 showing the $\omega$-compatibility of $J$.\\
 If $\sk{\cdot}{\cdot}$ is positive definite, $(\mathcal{O},\omega, J)$ is even a Kähler manifold. To see that, we first rewrite $\omega$:
 \begin{gather*}
  \omega_w (X_u, X_v) = \sk{w}{[u,v]} = \sk{w}{\ad_u v} = \sk{\ad_u w}{-v} = \sk{X_u}{-v}\quad\forall u,v\in\mathfrak{g}.
 \end{gather*}
 We can now express $g\coloneqq \omega (\cdot, J\cdot)$ as:
 \begin{align*}
  g_w(X_u, X_v) &= \omega_w (X_u, J_wX_v) = \omega_w (X_u, X_{J_w v})\\
  &= \sk{X_v}{-J_wv} = \frac{\sk{X_u}{X_v}}{\mu}\quad\forall u\in\mathfrak{g}\,\forall v\in E_\mu.
 \end{align*}
 Since $g_w$ is symmetric, the last equation shows that the spaces $E_\mu$ are orthogonal with respect to $g_w$ and that, restricted to $E_\mu$, $\mu g_w$ coincides with $\sk{\cdot}{\cdot}$. In particular, $g_w$ is positive definite if $\skcdot$ is positive definite.
\end{proof}

\begin{remark}
 The signature of $g_w$ is determined by the signature of $\sk{\cdot}{\cdot}$ restricted to $\im\ad_w$, as the previous proof shows. In general, $\sk{\cdot}{\cdot}$ does not need to be positive definite on $\im\ad_w$. This occurs, for example, in the case of complex reductive groups (cf. \autoref{sec:holo_semi-kaehler}).
\end{remark}

Before we end the discussion of semi-Kähler structures on (co)adjoint orbits, we add some comments regarding the uniqueness of those structures. In \cite{Besse2007}, it is shown that every $G$-invariant, closed two-form on an adjoint orbit $\mathcal{O}$ of a compact Lie group $G$ is the image of an $\Ad$-equivariant map\footnote{Note that, in \cite{Besse2007}, the maps $s$ are denoted by $\sigma$ and assumed to be $G$-invariant sections of the vector bundle $\mathfrak{s}\to\mathcal{O}$, where the fiber $\mathfrak{s}_w$ is the center of $\ker\ad_w$ for every $w\in\mathcal{O}$. As we will see, every map $s:\mathcal{O}\to\mathfrak{g}$ satisfying $s(\Ad (g)w) = \Ad (g)s(w)$ for all $g\in G$ and $w\in\mathcal{O}$ automatically fulfills these criteria.} $s:\mathcal{O}\to\mathfrak{g}$ by transgression and vice versa. This means that the set of $G$-invariant, closed two-forms is equal to the set of forms $\omega_s$ defined as follows:
\begin{gather*}
 \omega_{s,w} (X_u, X_v)\coloneqq \sk{s(w)}{[u,v]}\quad\forall w\in\mathcal{O}\ \forall u,v\in\mathfrak{g},
\end{gather*}
where $\skcdot$ is a fixed, $\Ad$-invariant, positive definite scalar product and $s$ is a map from $\mathcal{O}$ to $\mathfrak{g}$ satisfying $s(\Ad (g)w) = \Ad (g)s(w)$ for all $g\in G$ and $w\in\mathcal{O}$. Furthermore, it is proven in \cite{Besse2007} that every form $\omega_s$ is compatible with the canonical complex structure $J$ in the sense that the equation $\omega_s (J\cdot, J\cdot) = \omega_s$ holds. If $s$ is chosen such that $\omega_s$ is non-degenerate, $(\mathcal{O},\omega_s, J)$ forms a semi-Kähler manifold. It is Kähler if $\sk{ds_w\cdot}{\cdot}$ restricted to $\im\ad_w$ is positive definite. This observation follows from a calculation similar to the one at the end of the proof of Theorem \autoref{thm:semi-kaehler}. Simply put, the forms\footnote{The form $\omega$ from Proposition \autoref{prop:adjoint_omega} is the image of the map $s(w)\coloneqq w$ by transgression. If $\omega^\prime$ is the pullback of $\kks$ with respect to another $\Ad$-invariant scalar product $\skcdot^\prime$, then $\omega^\prime$ is the image of $s = \#\circ b^\prime$ by transgression.} $\omega_s$ exhaust all $G$-invariant semi-Kähler structures on $\mathcal{O}$ with complex structure $J$.\\
If we allow $G$ to be non-compact and assume that $\skcdot$ is only a non-degenerate, $\Ad$-invariant scalar product, then the image $\omega_s$ of an $\Ad$-equivariant map\linebreak $s:\mathcal{O}\to\mathfrak{g}$ is still $G$-invariant, as one easily checks. One can perform the same computation as in the compact case (cf. Lemma 8.67 and 8.68 in \cite{Besse2007}) to show that $\omega_s$ is also closed. Replacing $\omega$ with $\omega_s$ in the proof of Theorem \autoref{thm:semi-kaehler} reveals that $\omega_s$ further satisfies $\omega_s (J\cdot,J\cdot) = \omega_s$. Indeed, the only vital ingredient for that proof is $s(w)\in\ker\ad_w$ for all $w\in\mathcal{O}$. This property, however, is an immediate consequence of the $\Ad$-equivariance of $s$: Taking the derivative of $\Ad (g)s(w) = s(\Ad (g)w)$ with respect to $g$ gives us:
\begin{gather*}
 \ad_u s(w) = ds_w (\ad_uw)\quad\forall u\in\mathfrak{g}\ \forall w\in\mathcal{O}.
\end{gather*}
For $u=w$, this yields $\ad_w s(w) = 0$ and, therefore, $s(w)\in\ker\ad_w$. In fact, we find $\ad_u s(w) = 0$ for all $u\in\ker\ad_w$ which implies that $s(w)$ lies in the center $\mathfrak{s}_w$ of the Lie subalgebra $\ker\ad_w\subset\mathfrak{g}$.\\
Combining these results, we find that, as in the compact case, $(\mathcal{O},\omega_s, J)$ is semi-Kähler if $\omega_s$ is non-degenerate and Kähler if $\sk{ds_w\cdot}{\cdot}$ restricted to $\im\ad_w$ is positive definite. The only difference between the compact and non-compact case is that we do not know whether the forms $\omega_s$ exhaust all $G$-invariant, closed two-forms and, therefore, all $G$-invariant semi-Kähler structures on $\mathcal{O}$ with complex structure $J$. In the compact case, the proof of this fact requires an in-depth analysis of the root system which may not apply to non-compact groups.\\
To conclude this part, we collect our findings in the following theorem:

\begin{theorem}[Uniqueness of $G$-invariant semi-Kähler structures]\label{thm:semi-kaehler_unique}
 Let $G$ be a Lie group with fixed $\Ad$-invariant, non-degenerate scalar product $\skcdot$ on its Lie algebra $\mathfrak{g}$ and skew-symmetric adjoint orbit $\mathcal{O}\subset\mathfrak{g}$. Further, let $J$ be the canonical complex structure on $\mathcal{O}$. For any $\Ad$-equivariant map $s:\mathcal{O}\to\mathfrak{g}$, the image $\omega_s$ of $s$ by transgression, i.e.
 \begin{gather*}
  \omega_{s,w} (X_u, X_v)\coloneqq \sk{s(w)}{[u,v]}\quad\forall w\in\mathcal{O}\ \forall u,v\in\mathfrak{g},
 \end{gather*}
 is a closed, $G$-invariant two-form on $\mathcal{O}$. $(\mathcal{O},\omega_s, J)$ is a semi-Kähler manifold if $\omega_s$ is non-degenerate and a Kähler manifold if $\sk{ds_w\cdot}{\cdot}$ restricted to $\im\ad_w$ is positive definite. If $G$ is compact, the forms $\omega_s$ exhaust all closed, $G$-invariant two-forms and all $G$-invariant semi-Kähler structures on $\mathcal{O}$ with complex structure $J$.
\end{theorem}

\subsection*{Application: Compact Lie Groups}

In the last subsection of \autoref{sec:semi-kaehler}, we cross-check our construction against \cite{Besse2007} by applying Theorem \autoref{thm:semi-kaehler} to compact Lie groups. Our goal is to prove the following theorem (cf. \cite{Besse2007}):

\begin{theorem}[Kähler structures on (co)adjoint orbits]\label{thm:kaehler}
 Let $G$ be a compact Lie group with Lie algebra $\mathfrak{g}$ and dual Lie algebra $\mathfrak{g}^\ast$. Further, let $\mathcal{O}\subset\mathfrak{g}$ be an adjoint orbit of $G$ and $\mathcal{O}^\ast\subset\mathfrak{g}^\ast$ a coadjoint orbit of $G$. Then, $\mathcal{O}\subset\mathfrak{g}$ and $\mathcal{O}^\ast\subset\mathfrak{g}^\ast$ are embedded submanifolds and carry $G$-invariant Kähler structures. The complex structure of $\mathcal{O}$ and the symplectic structure of $\mathcal{O}^\ast$ are the canonical ones of adjoint and coadjoint orbits, respectively.
\end{theorem}

\begin{proof}
 To obtain Theorem \autoref{thm:kaehler} from Theorem \autoref{thm:semi-kaehler}, we need to answer the following questions for compact Lie groups:
 \begin{enumerate}
  \item[1.] Why are all (co)adjoint orbits embedded submanifolds?
  \item[2.] Why are all adjoint orbits skew-symmetric?
  \item[3.] Why does the Lie algebra always admit an $\Ad$-invariant metric?
  \item[4.] Why are the orbits not only semi-Kähler, but also Kähler manifolds?
 \end{enumerate}
 The first question is pretty easy to answer: By Lemma \autoref{lem:immersed_orbit}, the maps\linebreak $G/G_{w_0}\to\mathfrak{g}$, $[g]\mapsto \Ad (g)w_0$ and $G/G_{\alpha_0}\to\mathfrak{g}^\ast$, $[g]\mapsto \Ad^\ast (g)\alpha_0$ are injective immersions with image $\mathcal{O}$ and $\mathcal{O}^\ast$, respectively. Since $G$ is compact, the quotients $G/G_{w_0}$ and $G/G_{\alpha_0}$ are compact as well. Thus, the maps $G/G_{w_0}\to\mathfrak{g}$ and $G/G_{\alpha_0}\to\mathfrak{g}^\ast$ are proper. It is a standard result from differential geometry that proper injective immersions are embeddings answering the first question.\\
 The remaining questions can be answered all at once: By Remark \autoref{rem:Ad-invariant_metric} and Theorem \autoref{thm:semi-kaehler}, the existence of an $\Ad$-invariant, positive definite scalar product $\skcdot$ guarantees that every orbit is skew-symmetric and that the associated semi-Kähler structure is Kähler. The existence of $\skcdot$ is ensured by Proposition \autoref{prop:averaging}.
\end{proof}

To construct the desired scalar product, we use a standard trick from representation theory. We turn any positive definite scalar product on $\mathfrak{g}$ into an $\Ad$-invariant scalar product by averaging over $G$:

\begin{proposition}\label{prop:averaging}
 Let $G$ be a compact Lie group with Lie algebra $\mathfrak{g}$ and positive definite scalar product $\sk{\cdot}{\cdot}$ on $\mathfrak{g}$. Then, the averaged scalar product $\sk{\cdot}{\cdot}_G$ defined by
 \begin{gather*}
  \sk{u}{v}_G\coloneqq\frac{1}{\vol G}\int\limits_G \sk{\Ad (g)u}{\Ad (g)v} d\mu_G (g)\quad\forall u,v\in\mathfrak{g}
 \end{gather*}
 is an $\Ad$-invariant, positive definite scalar product on $\mathfrak{g}$, where $d\mu_G$ is the Haar measure on $G$ and $\vol G$ is the volume of $G$ with respect to the Haar measure.
\end{proposition}

\begin{proof}
 $\sk{\cdot}{\cdot}_G$ is clearly a positive definite scalar product. We only need to check that $\sk{\cdot}{\cdot}_G$ is $\Ad$-invariant. We compute:
 \begin{align*}
  \sk{\Ad (g)u}{\Ad (g)v}_G &= \frac{1}{\vol G}\int\limits_G \sk{\Ad (g^\prime g)u}{\Ad (g^\prime g)v} d\mu_G (g^\prime)\\
  &= \frac{1}{\vol G}\int\limits_G \sk{\Ad (h)u}{\Ad (h)v} d\mu_G (hg^{-1})\\
  &= \frac{1}{\vol G}\int\limits_G \sk{\Ad (h)u}{\Ad (h)v} d\mu_G (h)\\
  &= \sk{u}{v}_G\quad\forall g\in G\,\forall u,v\in\mathfrak{g},
 \end{align*}
 where we set $h = g^\prime g$ and used the fact that the Haar measure on a compact Lie group is right invariant.
\end{proof}

\newpage
\section[Holomorphic Semi-Kähler Structure of Complex Coadjoint Orbits]{Holomorphic Semi-Kähler Structure of Complex Coadjoint Orbits\sectionmark{Holomorphic Semi-Kähler Structure}}
\sectionmark{Holomorphic Semi-Kähler Structure}
\label{sec:holo_semi-kaehler}
As demonstrated in \autoref{sec:semi-kaehler}, the skew-symmetric orbits of a Lie group $G$ admitting a bi-invariant metric carry $G$-invariant semi-Kähler structures. The goal of this section is to show that the semi-Kähler structures are even holomorphic\footnote{Confer Definition \autoref{def:holo_kaehler_main} or \autoref{app:kaehler}.} if $G$ is a complex Lie group. We especially investigate the situation where $G$ is the complexification of a real Lie group $G_\R$ with bi-invariant metric and skew-symmetric orbit. In this case, the holomorphic semi-Kähler structure naturally comes with a real structure\footnote{Confer Definition \autoref{def:comp_kaehler_real_str_main} or \autoref{app:kaehler}.} whose real form consists of skew-symmetric orbits of $G_\R$. Examples of such $G$ include complex reductive groups, i.e., complex Lie groups with compact real forms.\\
\autoref{sec:holo_semi-kaehler} is split up into three parts: The first part offers a short explanation of the most important objects we use throughout the section, particularly holomorphic semi-Kähler structures and complex Lie groups. In the second part, we construct holomorphic semi-Kähler structures on skew-symmetric orbits of complex Lie groups $G$ admitting a bi-invariant metric (Theorem \autoref{thm:complex_orbit} and Corollary \autoref{cor:complex_orbit_transgression}). Lastly, we consider the special case of $G$ being the complexification of a real group $G_\R$ (Theorem \autoref{thm:complex_orbit_real_form} and Corollary \autoref{cor:transgression_complexification}) and apply the results to complex reductive groups (Corollary \autoref{cor:reductive}).

\subsection*{Preliminaries}

This subsection introduces the concepts relevant to \autoref{sec:holo_semi-kaehler}, especially holomorphic semi-Kähler manifolds and complex Lie groups. A detailed account of these notions can be found in \autoref{app:kaehler} and \autoref{app:complex_lie_groups}.\\
Let us begin with holomorphic semi-Kähler manifolds:

\begin{definition}[Holomorphic semi-Kähler manifolds]\label{def:holo_kaehler_main}
 A \textbf{complexified pre-semi-Kähler} manifold is a collection $(X,\omega, J, I)$ where $X^{4n}$ is a smooth manifold and the tensors $\omega\in\Omega^2 (X)$ and $I,J\in\Gamma\End (TX)$ satisfy:
 \begin{enumerate}
  \item $(X,\omega, J)$ is a pre-semi-Kähler manifold,
  \item $I$ is an almost complex structure, i.e., $I^2_p = -\id_{T_pX}$ for all $p\in X$,
  \item $I$ is anticompatible with $\omega$ and commutes with $J$:
  \begin{gather*}
   \omega (I\cdot, I\cdot) = -\omega\quad\text{and}\quad IJ = JI.
  \end{gather*}
 \end{enumerate}
 We say that $(X,\omega,J,I)$ is \textbf{holomorphic} instead of complexified if $I$ is integrable and $\Omega\coloneqq \omega -i\omega (I\cdot,\cdot)$ as well as $J$ viewed as a section\footnote{Here, the subscript indicates that the decomposition $T_\C X = T^{(1,0)}X\oplus T^{(0,1)}X$ is understood with respect to $I$.} of $\End (T^{(1,0)}_IX)$ is holomorphic. We drop the prefix ``pre'' if $(X,\omega,J)$ is semi-Kähler.
\end{definition}

Here, the prefix ``semi'' has a slightly different meaning than for semi-Kähler manifolds. To give a precise explanation, we first need to define real structures\footnote{The general notion of real structures on complex manifolds is discussed in \autoref{app:real_structures}.} on holomorphic semi-Kähler manifolds:

\begin{definition}[Real structure of Kähler manifolds]\label{def:comp_kaehler_real_str_main}
 Let $(X,\omega,J,I)$ be a complexified pre-semi-Kähler manifold. A \textbf{real structure} $\sigma$ on $(X,\omega,J,I)$ is a smooth involution on $X$ satisfying:
 \begin{enumerate}
  \item $\sigma$ preserves $\omega$, i.e., $\sigma^\ast\omega = \omega$,
  \item $\sigma$ is $J$-holomorphic, i.e., $J\circ d\sigma = d\sigma\circ J$,
  \item $\sigma$ is $I$-antiholomorphic, i.e., $I\circ d\sigma = -d\sigma\circ I$.
 \end{enumerate}
 The fixed point set $M\coloneqq\Fix\sigma$ is called \textbf{real form}. $M$ is \textbf{nice} if it meets every connected component of $X$. We drop the prefix ``semi'' if $g\coloneqq\omega (\cdot,J\cdot)$ restricts to a Riemannian metric on a non-empty real form $M$.
\end{definition}

As the name suggests, the real form $M$ of a complexified pre-semi-Kähler\linebreak manifold $X$ is a pre-semi-Kähler manifold (cf. \autoref{app:kaehler}). Per definition, $M$ is pre-Kähler if and only if $X$ is a complexified pre-Kähler manifold.\\
Next, we introduce the concept of complex Lie groups: A complex Lie group is a Lie group equipped with a complex structure such that the group operations (left and right multiplication as well as inversion) are holomorphic. It turns out that all constructions and results from \autoref{sec:lie_groups} ``adapt properly to the complex setup'': The associated geometrical objects (Lie algebra, (co)tangent bundle, (co)adjoint orbit, and so on) naturally inherit a complex structure from the underlying complex Lie group. The maps constructed in \autoref{sec:lie_groups} are compatible with these complex structures in a suitable sense, for instance the Lie bracket $[\cdot,\cdot]$ is $\C$-bilinear, the exponential map $\exp$ is holomorphic, the adjoint action $\Ad (g)$ commutes with the complex structure on $\mathfrak{g}$ and is holomorphic in $g$, and so on.\\
Similar to Kähler manifolds, Lie groups may also possess real structures:

\begin{definition}[Real structure of Lie groups]\label{def:group_real_str_main}
 Let $G$ be a complex Lie group. A \textbf{real structure} on $G$ is an antiholomorphic involution $\sigma:G\to G$ which is also a group homomorphism. Its fixed point set $G_\R\coloneqq\Fix\sigma$ is called \textbf{real form}. We say $G_\R$ is \textbf{nice} if $G_\R$ meets every connected component of $G$. In this case, we call $G$ a \textbf{complexification} of $G_\R$. A complex Lie group is called \textbf{reductive} if it admits a nice compact real form.
\end{definition}

The real form $G_\R$ of a complex Lie group $G$ is itself a Lie group with dimension $\dim_\R G_\R = \dim_\C G$ (cf. \autoref{app:complex_lie_groups}). Most real Lie groups $G_\R$ admit a special kind of complexification, the \textbf{universal complexification}. The universal complexification $G$ of a Lie group $G_\R$ is defined by the universal property that every Lie group homomorphism $f:G_\R\to H$ into a complex Lie group $H$ can be uniquely extended to a complex Lie group homomorphism $F:G\to H$. As a universal object, $G$ is unique up to isomorphisms. It is a well-established fact that every compact Lie group $G_\R$ possesses a universal complexification $G$ constructed via the polar decomposition (cf. \cite{Bump2004} and Lemma \autoref{lem:kaehler_reductive_groups}).

\subsection*{Holomorphic Semi-Kähler structure of (Co)Adjoint Orbits}

The main goal of this part is to show Theorem \autoref{thm:complex_orbit} which states that the semi-Kähler structures on adjoint orbits of complex Lie groups are holomorphic:

\begin{theorem}[Holomorphic semi-Kähler structure on adjoint orbits]\label{thm:complex_orbit}
 Let $G$ be a Lie group that satisfies the conditions of Theorem \autoref{thm:semi-kaehler}, i.e., $G$ admits a bi-invariant semi-Riemannian metric and a skew-symmetric adjoint orbit $\mathcal{O}\subset\mathfrak{g}$. Denote the $G$-invariant semi-Kähler structure on $\mathcal{O}$ from Theorem \autoref{thm:semi-kaehler}\linebreak by $(\mathcal{O}, \omega, J)$. If $G$ is a complex Lie group with complex structure $I$, then $(\mathcal{O},\omega, J, I_e)$\footnote{For the sake of simplicity, we denote the complex structure on $\mathcal{O}$ by $I_e$, even though this is somewhat inaccurate. The precise definition of the complex structure on $\mathcal{O}$ is given in \autoref{app:complex_lie_groups} and in the proof of Theorem \autoref{thm:complex_orbit}.} is a $G$-invariant holomorphic semi-Kähler manifold.
\end{theorem}

\begin{remark}
 Of course, a similar result holds for suitable coadjoint orbits.
\end{remark}

\begin{proof}
 Evaluating the complex structure $I$ of $G$ at the neutral element $e\in G$ gives us a complex structure on the Lie algebra $(\mathfrak{g}, [\cdot,\cdot])$. Restricting $I_e$ to $\im\ad_w\cong T_w\mathcal{O}$ allows us to view $I_e$ as a section of $\End (T\mathcal{O})$. Clearly, this section is an almost complex structure on $\mathcal{O}$. As one would expect, $I_e$ is even a complex structure, i.e., integrable (cf. \autoref{app:complex_lie_groups}). Moreover, $I_e$ is $G$-invariant meaning that it commutes with $\Ad (g)$ for every $g\in G$. Thus, we only need to check three properties: $I_e$ and $J$ commute, $J$ viewed as a section of $\End (T^{(1,0)}_{I_e}\mathcal{O})$ is holomorphic, and $I_e$ is $\omega$-anticompatible\footnote{Precisely speaking, we also need to check that $\Omega\coloneqq \omega -i\omega (I_e\cdot,\cdot)$ is $I_e$-holomorphic. However, this is not necessary, since the form $\Omega\coloneqq \omega -i\omega (I\cdot,\cdot)$ is automatically $I$-holomorphic if $\omega$ is closed, $I$ is integrable, and $\omega$ is anticompatible with $I$ (cf. Theorem \autoref{thm:rel_HSM_PHSM} and Remark \autoref{rem:comp_kaehler}).}.\\
 To verify $I_eJ = JI_e$, we first recall the definition of $J$:
 \begin{gather*}
  J_wv\coloneqq \frac{1}{\mu}\ad_wv\quad\forall v\in E_\mu,
 \end{gather*}
 where $\mu>0$ is the norm of a complex eigenvalue of $\ad_w$, $E_\mu$ is the space $\im\ad_w\cap (E_{i\mu}\oplus E_{-i\mu})$, and $E_{i\mu}, E_{-i\mu}$ are the complex eigenspaces of $\ad_w$ with respect to the eigenvalues $i\mu, -i\mu$, respectively. Since $I_e$ turns $(\mathfrak{g}, [\cdot,\cdot])$ into a complex Lie algebra, we have the following relation:
 \begin{gather*}
  I_e [u,v] = [I_eu, v] = [u,I_ev]\quad\forall u,v\in\mathfrak{g}.
 \end{gather*}
 In particular, this means that $I_e$ and $\ad_w$ commute. Hence, $I_e$ maps eigenspaces of $\ad_w$ to eigenspaces of $\ad_w$ implying $I_e (E_\mu)\subset E_\mu$. It is now obvious from the definition of $J_w$ that $I_e$ and $J_w$ also commute.\\
 Next, we need to check that the section $J\in\End (T^{(1,0)}_{I_e}\mathcal{O})$ is $I_e$-holomorphic. This is a rather easy task: We simply note that the adjoint action of a complex Lie group is holomorphic and repeat the proof of Proposition \autoref{prop:J_smooth} in the complex category.\\
 Lastly, we have to consider how $\omega$ and $I_e$ interact in order to show\linebreak $\omega (I_e\cdot,I_e\cdot) = -\omega$. Recall that $\omega$ is defined as follows:
 \begin{gather*}
  \omega_w (X_u (w), X_v (w)) = \sk{w}{[u,v]}\quad\forall w\in\mathcal{O}\ \forall u,v\in\mathfrak{g},
 \end{gather*}
 where $X_u (w) = -\ad_w u$ and $X_v (w) = -\ad_w v$ are the fundamental vector fields associated with $u$ and $v$, respectively, and $\skcdot$ is the $\Ad$-invariant scalar product on $\mathfrak{g}$ determined by the bi-invariant metric on $G$. We now use the fact that $I_e$ is the complex structure of the complex Lie algebra $(\mathfrak{g}, [\cdot,\cdot])$:
 \begin{align*}
  \omega_w (I_eX_u (w), I_eX_v (w)) &= \omega_w (X_{I_eu} (w), X_{I_ev} (w))\\
  &= \sk{w}{[I_eu, I_ev]}\\
  &= \sk{w}{I^2_e[u,v]} = -\sk{w}{[u,v]}\\
  &= -\omega_w (X_u (w), X_v (w)).
 \end{align*}
 This shows $\omega (I_e\cdot,I_e\cdot) = -\omega$ concluding the proof.
\end{proof}

Let us quickly comment on one surprising fact, namely that the holomorphic semi-Kähler manifolds $(\mathcal{O}, \omega, J, I_e)$ and $(\mathcal{O}^\ast, \kks, b_\ast J, I^\ast_e)$ are isomorphic, even though $b:\mathfrak{g}\to\mathfrak{g}^\ast$ is not necessarily holomorphic:

\begin{remark}
 In \autoref{sec:semi-kaehler}, we have stated that the semi-Kähler structures on $\mathcal{O}$ and $\mathcal{O}^\ast$ are isomorphic. The reasoning behind this statement is that the musical map $b:\mathfrak{g}\to\mathfrak{g}^\ast$ derived from the scalar product $\skcdot$ gives rise to diffeomorphisms between adjoint and coadjoint orbits. If $G$ is a complex Lie group, $b:(\mathcal{O},I_e)\to(\mathcal{O}^\ast,I^\ast_e)$ is even a biholomorphism. To see this, observe that the (co)adjoint action is a holomorphic and transitive action on $\mathcal{O}$ ($\mathcal{O}^\ast$). As $b:\mathcal{O}\to\mathcal{O}^\ast$ intertwines adjoint and coadjoint action, it must be holomorphic as well.\\
 However, this fact does not imply that $b:\mathfrak{g}\to\mathfrak{g}^\ast$ is holomorphic. Indeed, $b$ is holomorphic if and only if $b$ satisfies:
 \begin{gather*}
  b(I_e v) = I^\ast_e(b(v)) \equiv b(v)\circ I_e\quad\forall v\in\mathfrak{g}.
 \end{gather*}
 By unfolding the definition of $b$, we can rephrase the preceding equation as $\sk{I_e\cdot}{\cdot} = \sk{\cdot}{I_e\cdot}$ or, equivalently, $\sk{I_e\cdot}{I_e\cdot} = -\skcdot$. One has $\sk{I_e\cdot}{I_e\cdot} = -\skcdot$ if and only if $\sk{\cdot}{\cdot}$ is the real part of a $\C$-bilinear two-form on $\mathfrak{g}$.\\
 Nevertheless, not every $\Ad$-invariant scalar product $\skcdot$ is the real part of a $\C$-bilinear two-form. To construct a counterexample, consider the case $G = G^\prime\times H$ where $G^\prime$ is some complex Lie group admitting $\Ad$-invariant scalar products and $H$ is a complex Abelian Lie group. We obtain a scalar product $\skcdot$ on $\mathfrak{g} = \mathfrak{g}^\prime\oplus\mathfrak{h}$ by choosing scalar products $\skcdot_{G^\prime}$ on $\mathfrak{g}^\prime$ and $\skcdot_H$ on $\mathfrak{h}$. If $\skcdot_{G^\prime}$ is $\Ad$-invariant, then $\skcdot$ is $\Ad$-invariant as well, since $\skcdot_H$ is automatically $\Ad$-invariant (the adjoint action of an Abelian group is trivial). However, $\skcdot$ cannot be the real part of a $\C$-bilinear two-form if $\skcdot_H$ is not one which is possible. We will see that $\skcdot$ is naturally the real part of a $\C$-bilinear two-form if $G$ is a complexification of a real Lie group $G_\R$.
\end{remark}

For a real Lie group $G$, $(\mathcal{O},\omega,J)$ is not the only semi-Kähler structure we can construct on an adjoint orbit $\mathcal{O}$ of $G$. As exemplified in \autoref{sec:semi-kaehler}, adjoint orbits admit a large class of semi-Kähler structures where we simply replace $\omega$ in $(\mathcal{O},\omega, J)$ with the image $\omega_s$ of an $\Ad$-equivariant map $s:\mathcal{O}\to\mathfrak{g}$ by transgression, i.e.:
\begin{gather*}
 \omega_{s,w} (X_u (w), X_v (w))\coloneqq\sk{s(w)}{[u,v]}\quad\forall w\in\mathcal{O}\ \forall u,v\in\mathfrak{g}.
\end{gather*}
Of course, we immediately ask at this point whether $(\mathcal{O}, \omega_s, J, I_e)$ is also a holomorphic semi-Kähler manifold for complex Lie groups $G$. To answer that question, we only need to check whether $I_e$ is $\omega_s$-anticompatible. This is obviously the case:
 \begin{align*}
  \omega_{s,w} (I_eX_u (w), I_eX_v (w)) &= \omega_{s,w} (X_{I_eu} (w), X_{I_ev} (w))\\
  &= \sk{s(w)}{[I_eu, I_ev]}\\
  &= \sk{s(w)}{I^2_e[u,v]} = -\sk{s(w)}{[u,v]}\\
  &= -\omega_{s,w} (X_u (w), X_v (w)).
 \end{align*}
Thus, we have shown the following corollary:

\begin{corollary}\label{cor:complex_orbit_transgression}
 Let $G$ be a Lie group admitting a bi-invariant semi-Riemannian metric. Furthermore, let $\mathcal{O}\subset\mathfrak{g}$ be a skew-symmetric adjoint orbit of $G$ and $s:\mathcal{O}\to\mathfrak{g}$ be an $\Ad$-equivariant map. Moreover, let $\omega_s$ be the closed two-form induced by $s$ (cf. Theorem \autoref{thm:semi-kaehler_unique}). If $\omega_s$ is non-degenerate and $G$ is a complex Lie group with complex structure $I$, then $(\mathcal{O},\omega_s, J, I_e)$ is a $G$-invariant holomorphic semi-Kähler manifold.
\end{corollary}

\subsection*{Application: Complex Reductive Groups}

Let us now turn our attention to complexifications. We first want to show Theorem \autoref{thm:complex_orbit_real_form} which, among others, states that $G$ automatically satisfies the conditions of Theorem \autoref{thm:complex_orbit} if $G$ is a complexification of a real Lie group $G_\R$ fulfilling the conditions of Theorem \autoref{thm:semi-kaehler}. Afterwards, we consider complex reductive groups as an example of such $G$ (Corollary \autoref{cor:reductive}).\\
To prove Theorem \autoref{thm:complex_orbit_real_form}, we need to answer two questions:
\begin{enumerate}[label = {(\arabic*)}]
 \item How does the metric on $G_\R$ become a metric on $G$?
 \item Why are complexifications of skew-symmetric orbits again skew-symmetric?
\end{enumerate}
Proposition \autoref{prop:comp_of_metric} covers the first question, while Proposition \autoref{prop:comp_of_skew-symmetric} answers the second one:

\begin{proposition}\label{prop:comp_of_metric}
 Let $G$ be a complex Lie group with complex structure $I$, nice\footnote{Confer Definition \autoref{def:group_real_str_main}.} real form $G_\R\subset G$ and Lie algebras $\mathfrak{g}_\R\subset\mathfrak{g}$. Further, let $\skcdot_\R$ be a $G_\R$-invariant\footnote{We sometimes say $G$-invariant instead of $\Ad$-invariant to specify the group with respect to which the scalar product at hand is preserved.}, non-degenerate scalar product on $\mathfrak{g}_\R$. Then, there exists a unique $G$-invariant, non-degenerate scalar product $\skcdot$ on $\mathfrak{g}$ which restricts to $\skcdot_\R$ on $\mathfrak{g}_\R$, satisfies $\sk{I_e\cdot}{I_e\cdot} = -\skcdot$ and with respect to which the decomposition $\mathfrak{g} = \mathfrak{g}_\R\oplus I_e\mathfrak{g}_\R$ is orthogonal.
\end{proposition}

\begin{proof}
 \textbf{Existence:} We define the scalar product $\skcdot$ on $\mathfrak{g} = \mathfrak{g}_\R\oplus I_e\mathfrak{g}_\R$ by:
 \begin{gather}
  \sk{u_1 + I_eu_2}{v_1 + I_ev_2}\coloneqq \sk{u_1}{v_1}_\R - \sk{u_2}{v_2}_\R\quad\forall u_1,u_2,v_1,v_2\in\mathfrak{g}_\R.\label{eq:comp_metric}
 \end{gather}
 Clearly, $\skcdot$ restricts to $\skcdot_\R$ on $\mathfrak{g}_\R$ and satisfies $\sk{I_e\cdot}{I_e\cdot} = -\skcdot$. It is obvious that the decomposition $\mathfrak{g} = \mathfrak{g}_\R\oplus I_e\mathfrak{g}_\R$ is orthogonal with respect to $\skcdot$. Since $\skcdot_\R$ is non-degenerate, $\skcdot$ is non-degenerate as well.\\
 It remains to be shown that $\skcdot$ is $G$-invariant. The $G_\R$-invariance of $\skcdot_\R$ implies:
 \begin{gather*}
  \sk{\ad_wu}{v}_\R = -\sk{u}{\ad_wv}_\R\quad\forall u,v,w\in\mathfrak{g}_\R.
 \end{gather*}
 Recall that $I_e$ commutes with $\ad_w$ for every $w\in\mathfrak{g}$. Thus, we obtain:
 \begin{gather}
  \sk{\ad_wu}{v} = -\sk{u}{\ad_wv}\quad\forall w\in\mathfrak{g}_\R\ \forall u,v\in\mathfrak{g}.\label{eq:skew_for_comp}
 \end{gather}
 Next, we note:
 \begin{gather*}
  \ad_{I_ew}u = -\ad_u I_ew = -I_e\ad_uw = I_e\ad_wu\quad\forall u,w\in\mathfrak{g}
 \end{gather*}
 which implies that \autoref{eq:skew_for_comp} also holds for $w\in\mathfrak{g} = \mathfrak{g}_\R\oplus I_e\mathfrak{g}_\R$. Integrating \autoref{eq:skew_for_comp} now gives us:
 \begin{gather}
  \sk{\Ad (\exp (tw))u}{\Ad (\exp (tw))v} = \sk{u}{v}\quad\forall t\in\R\ \forall u,v,w\in\mathfrak{g}.\label{eq:inv_local}
 \end{gather}
 Hence, $\skcdot$ is $\Ad (g)$-invariant for group elements $g$ in an open neighborhood of $e$. To extend the invariance to all of $G$, we distinguish two cases:\\\\
 \textbf{Case 1}: $G$ is connected.\\\\
 It is a basic fact from Lie group theory that every element $g$ in a connected Lie group $G$ can written as a product of exponentials:
 \begin{gather*}
  g = \exp(w_1)\exp (w_2)\ldots\exp (w_n).
 \end{gather*}
 This decomposition now allows us to repeatedly apply \autoref{eq:inv_local} to prove that $\skcdot$ is $G$-invariant.\\\\
 \textbf{Case 2}: $G$ is not connected.\\\\
 Denote the connected component of $G$ containing $e$ by $G_0$ and the remaining components by $G_i$. For every component $G_i$, fix one element $g_i\in G_i$. Then, we can write every element $g\in G$ as $g_ig_0$, where $g_i$ lies in the connected component containing $g$ and $g_0$ is some element in $G_0$. As we have seen in Case 1, the scalar product $\skcdot$ is $G_0$-invariant. Therefore, we only need to show that $\skcdot$ is invariant with respect to the fixed elements $g_i$. For this, we use the fact that $G_\R$ is a nice real form of $G$. It implies that $G_\R$ intersects every component $G_i$ non-trivially. Hence, we can chose the fixed elements $g_i$ to be contained in $G_\R\cap G_i$. We can directly infer from the definition of $\skcdot$ (together with the fact that $\skcdot_\R$ is $G_\R$-invariant and that $\Ad (g)$ commutes with $I_e$) that $\skcdot$ is $G_\R$-invariant. In particular, it is invariant with respect to the fixed elements $g_i$ proving that $\skcdot$ is $G$-invariant.
 \\\\
 \textbf{Uniqueness:} Let $\skcdot$ be a scalar product satisfying the properties specified in Proposition \autoref{prop:comp_of_metric}. The decomposition $\mathfrak{g} = \mathfrak{g}_\R\oplus I_e\mathfrak{g}_\R$ is orthogonal with respect to $\skcdot$, thus, we can write:
 \begin{gather*}
  \sk{u_1 + I_eu_2}{v_1 + I_ev_2}\coloneqq \sk{u_1}{v_1}_{\mathfrak{g}_\R} + \sk{u_2}{v_2}_{I_e\mathfrak{g}_\R}\quad\forall u_1,u_2,v_1,v_2\in\mathfrak{g}_\R,
 \end{gather*}
 where $\skcdot_{\mathfrak{g}_\R}$ and $\skcdot_{I_e\mathfrak{g}_\R}$ are scalar products on $\mathfrak{g}_\R$. The relation\linebreak $\sk{I_e\cdot}{I_e\cdot} = -\skcdot$ enforces $\skcdot_{\mathfrak{g}_\R} = -\skcdot_{I_e\mathfrak{g}_\R}$. Because $\skcdot$ restricts to $\skcdot_\R$ on $\mathfrak{g}_\R$, the scalar product $\skcdot_{\mathfrak{g}_\R}$ must coincide with $\skcdot_\R$. This shows that $\skcdot$ agrees with scalar product defined by \autoref{eq:comp_metric} finishing the proof.
\end{proof}

\begin{remark}
 The scalar product $\skcdot$ from Proposition \autoref{prop:comp_of_metric} is the real part of the non-degenerate $\C$-bilinear two-form:
 \begin{gather*}
  \skcdot_\C\coloneqq \skcdot - i\sk{I_e\cdot}{\cdot}.
 \end{gather*}
 Consequently, the signature of $\skcdot$ is $(n,n)$ where $n\coloneqq\dim_\R\mathfrak{g}_\R$.
\end{remark}

Regarding the second question, it suffices to show that complexifications of diagonalizable linear maps with purely imaginary eigenvalues are also diagonalizable with purely imaginary eigenvalues:

\begin{proposition}\label{prop:comp_of_skew-symmetric}
 Let $V$ be a real vector space equipped with $I\in\End (V)$ satisfying $I^2 = -\id_V$ and let $V_\R\subset V$ be a subspace such that $V = V_\R\oplus I(V_\R)$.\linebreak Further, let $A\in\End (V)$ be a $\R$-linear map such that $A(V_\R)\subset V_\R$,\linebreak $AI = IA$, and $A\in\End (V_{\R,\C})$\footnote{We set $V_{\R,\C}\coloneqq V_\R\otimes_\R\C$. For the sake of simplicity, we also denote the restrictions and complexifications of $A$ by $A$.} is diagonalizable with pairwise distinct eigenvalues $\pm i\mu_1,\ldots,\pm i\mu_n$ and eigenspaces $E^{V_{\R,\C}}_{\pm i\mu_1},\ldots, E^{V_{\R,\C}}_{\pm i\mu_n}\subset V_{\R,\C}$ ($\mu_1,\ldots,\mu_n>0$). Then, the complexification $A\in\End (V_\C)$ is also diagonalizable with eigenvalues $\pm i\mu_1,\ldots,\pm i\mu_n$. The eigenspaces $E^{V_\C}_{\pm i\mu_1},\ldots, E^{V_\C}_{\pm i\mu_n}\subset V_\C$ of $A\in\End (V_\C)$ satisfy:
 \begin{gather*}
  E^{V_\C}_{\pm i\mu_j} = E^{V_{\R,\C}}_{\pm i\mu_j}\oplus I(E^{V_{\R,\C}}_{\pm i\mu_j})\quad\forall j\in\{1,\ldots,n\}.
 \end{gather*}
\end{proposition}
 
\begin{proof}
 If $v\in V_{\R,\C}$ is an eigenvector of $A\in\End (V_{\R,\C})$ w.r.t. the eigenvalue $\pm i\mu_j$, then $v,Iv\in V_\C$ are eigenvectors of $A\in\End (V_\C)$ w.r.t. the same eigenvalue:
 \begin{gather*}
  Av = \pm i\mu_jv\Rightarrow AIv = IAv = \pm i\mu_j Iv.
 \end{gather*}
\end{proof}

Combining Proposition \autoref{prop:comp_of_metric} and \autoref{prop:comp_of_skew-symmetric} now gives us Theorem \autoref{thm:complex_orbit_real_form}:

\begin{theorem}[Adjoint orbits of complexifications]\label{thm:complex_orbit_real_form}
 Let $G$ be a complex Lie group with real structure $\sigma$ and nice real form $G_\R$. Further, assume that $G_\R$ admits a bi-invariant semi-Riemannian metric and a skew-symmetric adjoint orbit $\mathcal{O}_\R\subset\mathfrak{g}_\R$. Denote the adjoint orbit of $G$ containing $\mathcal{O}_\R$ by $\mathcal{O}$. Then, $\mathcal{O}$ admits a $G$-invariant holomorphic semi-Kähler structure given by $(\mathcal{O},\omega,J,I_e)$ from Theorem \autoref{thm:complex_orbit} together with a compatible real structure $\Sigma\coloneqq d\sigma_e\vert_{\mathcal{O}}$. If the metric on $G_\R$ is positive definite, $(\mathcal{O},\omega,J,I_e)$ forms with $\Sigma$ a holomorphic Kähler manifold.
\end{theorem}

\begin{remark}
 Again, a similar result holds for suitable coadjoint orbits.
\end{remark}

\begin{proof}
 First, we check that we can apply Theorem \autoref{thm:complex_orbit} to $G$. For this, we need to find an $\Ad$-invariant scalar product on $\mathfrak{g}$ and show that the orbit $\mathcal{O}$ is skew-symmetric. By Proposition \autoref{prop:comp_of_metric}, the bi-invariant metric on $G_\R$ or, equivalently, the $\Ad$-invariant scalar product on $\mathfrak{g}_\R$ gives rise to a unique $\Ad$-invariant scalar product on $\mathfrak{g}$. To prove that $\mathcal{O}$ is skew-symmetric, we first note that $\mathcal{O}$ contains the skew-symmetric orbit $\mathcal{O}_\R$ of $G_\R$. Thus, if we pick an element $w\in\mathcal{O}_\R\subset\mathcal{O}$, the complexification of $\ad_w\in\End (\mathfrak{g}_\R)$ is diagonalizable with purely imaginary eigenvalues. We can now apply Proposition \autoref{prop:comp_of_skew-symmetric} to the map $\ad_w$ which shows that the complexification of $\ad_w\in\End (\mathfrak{g})$ is also diagonalizable with purely imaginary eigenvalues. Hence, $\mathcal{O}$ contains a skew-symmetric element making $\mathcal{O}$ itself skew-symmetric. With this, all conditions of Theorem \autoref{thm:complex_orbit} are satisfied giving us the holomorphic semi-Kähler structure $(\mathcal{O},\omega,J,I_e)$.\\
 Now consider $\Sigma$. In \autoref{app:complex_lie_groups}, we show that $\Sigma$ is a real structure of the complex manifold $(\mathcal{O},I_e)$ with real form $\mathcal{O}\cap\mathfrak{g}_\R$. Thus, we only need to check that $\Sigma$ preserves $\omega$ and is $J$-holomorphic in order to prove that $\Sigma$ is a real structure of $(\mathcal{O},\omega,J,I_e)$. Let us first compute the pullback $\Sigma^\ast\omega$:
 \begin{align*}
  (\Sigma^\ast\omega)_w \left(X_u (w), X_v (w)\right) &= \omega_{\Sigma (w)}\left(d\Sigma_w X_u (w), d\Sigma_w X_v (w)\right)\\
  &= \omega_{d\sigma_e w}\left(d\sigma_e X_u (w), d\sigma_e X_v (w)\right)\\
  &= \omega_{d\sigma_e w}\left(X_{d\sigma_e u} (d\sigma_e w), X_{d\sigma_e v} (d\sigma_e w)\right)\\
  &= \sk{d\sigma_e w}{\left[d\sigma_e u, d\sigma_e v\right]}\\
  &= \sk{d\sigma_e w}{d\sigma_e [u,v]}\\
  &= \sk{w}{[u,v]}\\
  &= \omega_w \left(X_u (w), X_v (w)\right)\quad\forall w\in\mathcal{O}\,\forall u,v\in\mathfrak{g},
 \end{align*}
 where $\skcdot$ is the $\Ad$-invariant scalar product on $\mathfrak{g}$. To get from Line 2 to 3 and from Line 4 to 5, we used the fact that $d\sigma_e$ as the differential of a Lie group homomorphism is a Lie algebra homomorphism:
 \begin{gather*}
  d\sigma_e [u,v] = [d\sigma_e u, d\sigma_e v]\quad\forall u,v\in\mathfrak{g}.
 \end{gather*}
 To get from Line 5 to 6, we exploited the fact that $d\sigma_e$ preserves the scalar product $\skcdot$. To verify this statement, recall that $\skcdot$ is defined in Proposition \autoref{prop:comp_of_metric} as follows:
 \begin{gather*}
  \sk{u_1 + I_eu_2}{v_1 + I_ev_2}\coloneqq \sk{u_1}{v_1}_\R - \sk{u_2}{v_2}_\R\quad\forall u_1,u_2,v_1,v_2\in\mathfrak{g}_\R,
 \end{gather*}
 where $\mathfrak{g}$ decomposes into $\mathfrak{g}_\R\oplus I_e\mathfrak{g}_\R$ and $\skcdot_\R$ is the $\Ad$-invariant scalar product on $\mathfrak{g}_\R$. Since $d\sigma_e$ fixes $\mathfrak{g}_\R$ and anticommutes with $I_e$, $d\sigma_e$ acts on $\mathfrak{g}_\R\oplus I_e\mathfrak{g}_\R$ as $(u_1,u_2)\mapsto (u_1,-u_2)$. It is now clear from the definition of $\skcdot$ that $d\sigma_e$ preserves $\skcdot$.\\
 Next, we show that $\Sigma$ is $J$-holomorphic, i.e, that $d\Sigma$ commutes with $J$. We again use the fact that $d\sigma_e$ is a Lie algebra homomorphism:
 \begin{gather*}
  d\Sigma_w\circ \ad_w = d\sigma_e\circ\ad_w = \ad_{d\sigma_e w}\circ d\sigma_e = \ad_{\Sigma (w)}\circ d\Sigma_w\quad\forall w\in\mathcal{O}.
 \end{gather*}
 The last equation states that the differential $d\Sigma$ commutes with the $(1,1)$-tensor field\footnote{$w\mapsto\ad_w$ is a section of $\End (T\mathcal{O})$ and, therefore, a $(1,1)$-tensor field.} $w\mapsto \ad_w$. In particular, $d\Sigma$ preserves the ``eigenbundles'' $E_\mu$ of $w\mapsto \ad_w$. It is now evident from the definition of $J$,
 \begin{gather*}
  J_w v\coloneqq\frac{1}{\mu}\ad_w v\quad\forall v\in E_\mu\,\forall\text{ suitable } \mu>0,
 \end{gather*}
 that $d\Sigma$ also commutes with $J$.\\
 Lastly, we discuss the special case of $\skcdot_\R$ being positive definite. Our goal is to prove that $(\mathcal{O},\omega,J,I_e)$ forms with $\Sigma$ a holomorphic Kähler manifold in this case. In \autoref{app:kaehler}, we show that the real form of a holomorphic semi-Kähler manifold is itself a semi-Kähler manifold. In our situation, the semi-Kähler structure on the real form $\Fix\Sigma = \mathcal{O}\cap\mathfrak{g}_\R$ coincides with the one obtained from applying Theorem \autoref{thm:semi-kaehler} to $G_\R$. If $\skcdot_\R$ is positive definite, this semi-Kähler structure is even Kähler concluding the proof.
\end{proof}

Before we discuss complex reductive groups, we want to answer the question what happens if we replace $\omega$ in Theorem \autoref{thm:complex_orbit_real_form} with the image $\omega_s$ of an $\Ad$-equivariant map $s:\mathcal{O}\to\mathfrak{g}$ by transgression. We expect that the resulting holomorphic semi-Kähler structure is compatible with the real structure $\Sigma$. Indeed, the only condition we need to check for this is that $\Sigma$ also preserves $\omega_s$. By modifying the proof of Theorem \autoref{thm:complex_orbit_real_form}, we see that this is the case if $s$ commutes with $d\sigma_e$. It is straightforward to verify that the commutation relation $d\sigma_e\circ s = s\circ d\sigma_e$ is equivalent to $s(\mathcal{O}\cap\mathfrak{g}_\R)\subset\mathfrak{g}_\R$. The latter statement allows us to view $s$ as a complexification of the $G_\R$-equivariant map $s_\R\coloneqq s\vert_{\mathcal{O}\cap\mathfrak{g}_\R}$ into $\mathfrak{g}_\R$. Together with Theorem \autoref{thm:semi-kaehler_unique}, this gives us the following corollary:

\begin{corollary}[Complexifications of transgressions]\label{cor:transgression_complexification}
 Let $G$ be a complex Lie group with real structure $\sigma$ and nice real form $G_\R$. Further, assume that $G_\R$ admits a bi-invariant semi-Riemannian metric and a skew-symmetric adjoint orbit $\mathcal{O}_\R\subset\mathfrak{g}_\R$. Denote the adjoint orbit of $G$ containing $\mathcal{O}_\R$ by $\mathcal{O}$. Now let\linebreak $s:\mathcal{O}\to\mathfrak{g}$ be an $\Ad$-equivariant map satisfying $s(\mathcal{O}\cap\mathfrak{g}_\R)\subset \mathfrak{g}_\R$. If $\omega_s$ from Theorem \autoref{thm:semi-kaehler_unique} is non-degenerate, then $(\mathcal{O},\omega_s,J,I_e)$ is a $G$-invariant holomorphic semi-Kähler manifold compatible with the real structure $\Sigma\coloneqq d\sigma_e\vert_{\mathcal{O}}$. If $\sk{ds_{\R,w}\cdot}{\cdot}_\R$ restricted to $\im\ad_w$ is positive definite, $(\mathcal{O},\omega_s,J,I_e)$ forms with $\Sigma$ a holomorphic Kähler manifold.
\end{corollary}

Let us now apply Theorem \autoref{thm:complex_orbit_real_form} to complex reductive groups, i.e. complex Lie groups with nice compact real form:

\begin{corollary}[(Co)Adjoint orbits of complex reductive groups]\label{cor:reductive}
 Let $G$ be a complex reductive group with real form $G_\R$, Lie algebras $\mathfrak{g}_\R\subset\mathfrak{g}$, and dual Lie algebras\footnote{$\mathfrak{g}^\ast_\R$ is the space of all linear maps $\mathfrak{g} = \mathfrak{g}_\R\oplus I_e\mathfrak{g}_\R\to\R$ that vanish on $I_e\mathfrak{g}_\R$.} $\mathfrak{g}^\ast_\R\subset\mathfrak{g}^\ast$. Further, let $\mathcal{O}\subset\mathfrak{g}$ be an adjoint orbit of $G$ and $\mathcal{O}^\ast\subset\mathfrak{g}^\ast$ a coadjoint orbit of $G$ such that $\mathcal{O}\cap\mathfrak{g}_\R\neq\emptyset\neq\mathcal{O}^\ast\cap\mathfrak{g}^\ast_\R$. Then, $\mathcal{O}$ and $\mathcal{O}^\ast$ carry $G$-invariant holomorphic Kähler structures. The first complex structure ($J$) of $\mathcal{O}$ and the symplectic structure of $\mathcal{O}^\ast$ are the canonical ones of adjoint and coadjoint orbits, respectively. The second complex structure ($I$) of $\mathcal{O}$ and $\mathcal{O}^\ast$ is the canonical one induced by the complex structure of $G$.
\end{corollary}

\begin{proof}
 It suffices to show that the conditions of Theorem \autoref{thm:complex_orbit_real_form} are satisfied, i.e., that\dots
 \begin{enumerate}
  \item $\dots G_\R$ admits a bi-invariant Riemannian metric,
  \item $\dots \mathcal{O}$ ($\mathcal{O}^\ast$) contains a skew-symmetric\footnote{We say a coadjoint orbit is skew-symmetric if it is isomorphic to a skew-symmetric adjoint orbit via an $\Ad$-invariant scalar product.} orbit of $G_\R$.
 \end{enumerate}
 \textbf{Condition (i):} In \autoref{sec:semi-kaehler}, we have shown that every compact Lie group admits a bi-invariant Riemannian metric (cf. Proposition \autoref{prop:averaging}), so this also holds for the compact group $G_\R$.\\
 \textbf{Condition (ii):} As $\mathcal{O}$ ($\mathcal{O}^\ast$) intersects $\mathfrak{g}_\R$ ($\mathfrak{g}^\ast_\R$) non-trivially, $\mathcal{O}$ ($\mathcal{O}^\ast$) contains an orbit of the compact group $G_\R$. We have shown in \autoref{sec:semi-kaehler} that every (co)adjoint orbit of a compact Lie group is skew-symmetric concluding the proof.
\end{proof}

\begin{remark}
 Obviously, Corollary \autoref{cor:reductive} is still true if we replace the symplectic form $\omega$ by the image $\omega_s$ of an $\Ad$-equivariant map $s:\mathcal{O}\to\mathfrak{g}$ by transgression where $s$ satisfies $s(\mathcal{O}\cap\mathfrak{g}_\R)\subset\mathfrak{g}_\R$ (cf. Corollary \autoref{cor:transgression_complexification}).
\end{remark}

\newpage
\section{Kähler Duality}
\label{sec:duality}
It is known since the 90s that coadjoint orbits of semisimple complex reductive groups $G$ are not only Kählerian, but also admit the structure of a Hyperkähler manifold\footnote{Here, a complex reductive group $G$ denotes the universal complexification of its compact real form $G_\R$ (cf. Definition \autoref{def:group_real_str_main} and the remarks afterwards). The Hyperkähler structures on coadjoint orbits were first discovered by Kronheimer (cf. \cite{Kronheimer1990}) and later generalized by Kovalev (cf. \cite{Kovalev1996}). A sketch of their construction can be found in \autoref{app:hyp_orb}.}. As in the Kähler case, the Hyperkähler structure is not unique and depends on certain choices\footnote{These dependencies include the choice of a triple $(\tau_1,\tau_2,\tau_3)$ of Lie algebra elements and the choice of a homomorphism $\rho$ (cf. \cite{Kronheimer1990} and \cite{Kovalev1996} for details on $(\tau_1,\tau_2,\tau_3)$ and $\rho$). To the best of the author's knowledge, the Hyperkähler structures in question should also depend on the choice of an $\Ad$-invariant scalar product $\skcdot_\R$, as explained in \autoref{app:hyp_orb}. However, this dependence is nowhere explicitly mentioned.}. The existence of these Hyperkähler structures is remarkable, especially because their underlying orbits also possess a family of holomorphic Kähler structures, as we have seen in \autoref{sec:holo_semi-kaehler}. If a space admits a family of both Hyperkähler and holomorphic Kähler structures in a somewhat natural way, we say this space exhibits \textbf{Kähler duality}. In that regard, coadjoint orbits are an example of Kähler duality.\\
If we encounter such a rich and powerful structure like Kähler duality, we are naturally drawn to the question whether this structure is just coincidental or originates from some deeper and hidden theory. We conjecture that the Kähler duality of coadjoint orbits can be traced back to double cotangent bundles. Precisely speaking, we claim that double cotangent bundles admit a family of Hyperkähler and holomorphic Kähler structures and that, for Lie groups, these Kähler structures become, via a suitable reduction process (cf. \autoref{app:reduction}), the previously mentioned Kähler structures on coadjoint orbits.\\
The main goal of \autoref{sec:duality} is to construct the two different kinds of Kähler structures on double cotangent bundles. The following diagram illustrates the idea behind the procedure:
\begin{center}
 \begin{tikzcd}[column sep=3.6em,row sep=2em]
  \textbf{Hyperk.}& (T^\ast M, g_{T^\ast M}) \arrow{r}{\text{Stenzel}} & (T^\ast (T^\ast M), -\omega_{\can}, \phi^\ast_1 J_{g_{T^\ast M}}, T^\ast J_g)\\
  (M,g) \arrow{ur}[sloped,above]{\text{Stenzel}} \arrow{dr}[sloped,below]{\text{\footnotesize Complexification}} & &\\
  \textbf{Holo. K.}& (T^\ast M, g_\C) \arrow{r}[below]{\text{Stenzel}} & (T^\ast (T^\ast M), -\omega_{\can}, \phi^\ast_2 J_{g_\C}, T^\ast J_g)
 \end{tikzcd}
\end{center}
Our construction is based on Stenzel's theorem (cf. Theorem \autoref{thm:cotangent_kaehler}) which states that the cotangent bundle $T^\ast M$ of a Riemannian manifold $(M,g)$ admits a canonical Kähler structure $(T^\ast M,-\omega_{\can}, J_g)$. Here, $J_g$ is the complex structure adapted to $g$ (cf. Definition \autoref{def:adapted_comp_str}). As shown in the diagram above, the plan is to apply Stenzel's theorem to two natural metrics on $T^\ast M$: One is the Kähler metric $g_{T^\ast M}\coloneqq-\omega_{\can}(\cdot,J_g\cdot)$ coming from Stenzel's theorem, the other metric $g_\C$ is the complexification of $g$. For this idea to work, the complex structures $J_{g_{T^\ast M}}$ and $J_{g_\C}$ need to be twisted with diffeomorphisms $\phi_1,\phi_2:T^\ast(T^\ast M)\to T^\ast (T^\ast M)$ to ensure the right commutation relations (cf. Conjecture \autoref{con:commutation_relation}).\\
\autoref{sec:duality} is divided into four subsections: First, we summarize the history of Kähler structures on cotangent bundles. The purpose of this part is to provide historical context and explain which results are already known. Afterwards, we prove Stenzel's theorem. During this process, we introduce the notion of adapted complex structures. In the third subsection, we discuss the main result of \autoref{sec:duality}: We show -- barring a complete proof for the commutation relations (cf. Conjecture \autoref{con:commutation_relation}) -- that double cotangent bundles exhibit Kähler duality (cf. Theorem \autoref{thm:kaehler_duality}). Lastly, we apply Theorem \autoref{thm:kaehler_duality} to compact Lie groups $M = G_\R$ with bi-invariant metrics $g$ resulting in a family of Hyperkähler and holomorphic Kähler structures on $T^\ast G$, where $G$ is the universal complexification of $G_\R$ (cf. Theorem \autoref{thm:kaehler_duality_lie}). We conjecture that one can reduce these structures to similar structures on coadjoint orbits of $G$ (cf. Conjecture \autoref{con:reduction} and \autoref{app:reduction}).

\subsection*{Historical Background}

Our story begins in the late 80s when Matthew B. Stenzel started working on his PhD thesis under the supervision of his advisor Victor W. Guillemin. They tackled the question whether the cotangent bundle $T^\ast M$ of a real-analytic Riemannian manifold $(M,g)$ admits a unique complex structure $J_g$ such that $(T^\ast M, -\omega_{\can} = -d\lambda_{\can})$\footnote{Throughout \autoref{sec:duality}, we always have this inconvenient minus sign in front of $\omega_{\can}$. There are two ways to get rid of it: One can either change the definition of the Kähler metric ($g\coloneqq \omega (J\cdot,\cdot)$ instead of $g\coloneqq \omega (\cdot,J\cdot)$) or pick a different sign convention for $\omega_{\can}$. In light of \autoref{app:reduction}, where we see that $(T^\ast G,\omega_{\can})$ reduces to $(\mathcal{O}^\ast,-\kks)$, changing the sign convention for $\omega_{\can}$ might be the more natural choice. Still, we do not change the conventions to keep the notation consistent with \autoref{chap:PHHS}.} becomes a Kähler manifold with Kähler potential\linebreak $U_g (\alpha)\coloneqq g^\ast (\alpha,\alpha)$, where $g^\ast$ is the metric dual to $g$. It was known at that time that every real-analytic manifold $M$ can be embedded\footnote{This is the famous Bruhat-Whitney embedding theorem (cf. \cite{Whitney1959}).} as a totally real submanifold (cf. \autoref{app:real_structures}) into a complex manifold $X$ and that one can choose $X$ to be the cotangent bundle of $M$, where the real structure on $T^\ast M$ is just the fiberwise inversion $\sigma (\alpha)\coloneqq -\alpha$. However, the choice $X = T^\ast M$ does not fix the complex structure on $T^\ast M$. In fact, there are many complex structures on $T^\ast M$, some of which even turn $(T^\ast M, -\omega_{\can})$ into a Kähler manifold (cf. \cite{Stenzel1990}). The idea of Stenzel and Guillemin was to single out one complex structure $J_g$ by fixing the Kähler potential to be $U_g$. When Stenzel handed in his thesis in 1990, he was only able to give a formal solution to this problem in terms of a power series. Even though this series implies uniqueness of $J_g$ if it exists and Stenzel could show existence of $J_g$ for a plethora of examples, he was unable to prove convergence of the power series for all cases in his thesis (cf. \cite{Stenzel1990}). About one year later, Stenzel and Guillemin published a paper (cf. \cite{Guillemin1991}) in which they found a complete answer to the problem at hand: Abandoning the power series\footnote{In \cite{Guillemin1991}, Guillemin and Stenzel claim that they found a ``long and complicated proof'' for the convergence of the power series. However, they never published the proof, because, as they state, their ``efforts were completely misguided'': The other approach is much simpler and implicitly contained in \cite{Burns1982}, but they lacked ``the perspecality to notice so earlier''.}, they reduced the construction of $J_g$ to a Monge-Amp\`{e}re problem and solved the Monge-Amp\`{e}re equation by using the distance function on $M$ induced by $g$ and ideas of Dan Burns (cf. \cite{Burns1982}).\\
Independently of Stenzel and Guillemin, L{\'a}szl{\'o} Lempert and R{\'{o}}bert Sz{\H{o}}ke found a different way to define and construct $J_g$ (cf. \cite{Lempert1991} and \cite{Szoeke1991}). Initially, they set out to study the question to what extent global conditions determine the solutions of the Monge-Amp\`{e}re equation. In their setup, Lempert and Sz{\H{o}}ke were able to reduce this question to the problem of finding complex structures adapted to a metric $g$. As it turns out, the complex structure $J_g$ adapted to $g$ is unique and coincides with the complex structure constructed by Stenzel and Guillemin (cf. \cite{Szoeke1991}). Later on, Roger Bielawski and, independently, Sz{\H{o}}ke realized that adapted complex structures can also be defined for semi-Riemannian\footnote{Even though Stenzel's formal solution also works for semi-Riemannian manifolds $(M,g)$ and coincides with the adapted complex structure $J_g$, the results of the paper \cite{Guillemin1991} do not transfer to the semi-Riemannian case, as semi-Riemannian metrics $g$ do not give rise to a distance function.} manifolds $(M,g)$ and even Koszul manifolds $(M,\nabla)$, i.e., manifolds with connections (cf. \cite{Bielawski2003} and \cite{Szoeke2004}).\\
Shortly after the works of Guillemin-Stenzel and Lempert-Sz{\H{o}}ke, people started to wonder what kind of geometry $T^\ast M$ has if the semi-Riemannian manifold $(M,g)$ admits additional structures. Dmitry Kaledin (1997) and, independently, Birte Feix (in her PhD thesis, handed in in 1999) investigated the case where $M$ is a Kähler manifold. Each with their own method, they were able to show that $T^\ast M$ carries a Hyperkähler structure in this scenario: Kaledin used Hodge manifold theory (cf. \cite{Kaledin1997}), while Feix employed a twistor construction (cf. \cite{Feix2001}). A couple of years later, Bielawski discovered a similar result, namely a Hyperkähler structure on double cotangent bundles. As Bielawski points out (cf. \cite{Bielawski2003}), his construction is a special case of Kaledin's and Feix's theorem where the Kähler manifold $M$ is chosen to be a cotangent bundle equipped with Stenzel's Kähler structure.\\
Also around the turn of the millennium, Ralph Bremigan examined the geometry of cotangent bundles of semisimple Lie groups $M = G_\R$ equipped with the Killing form $K$ as their metric (cf. \cite{Bremigan2000}). The Kähler structure he found on $T^\ast G_\R$ agrees with Stenzel's one (cf. the end of \autoref{sec:duality}) and is invariant under both left and right translations.

\subsection*{Adapted Complex Structures}

In this subsection, we introduce the notion of an adapted complex structure and work out its most important properties following \cite{Bielawski2003} and \cite{Szoeke2004}. We utilize these structures to formulate and prove Stenzel's theorem (cf. Theorem \autoref{thm:cotangent_kaehler}). The main result of this subsection are the commutation relations which we can only show on the zero section and for flat $g$ (cf. Lemma \autoref{lem:commutation_relation} and Conjecture \autoref{con:commutation_relation}). We start with the definition of an adapted complex structure:

\begin{definition}[Adapted complex structure]\label{def:adapted_comp_str}
 Let $M$ be a manifold with connection\footnote{By a connection $\nabla$ on $M$, we mean a linear connection $\nabla$ on the vector bundle $TM\to M$.} $\nabla$ on it and let $U\subset TM$ be an open neighborhood of the zero section\footnote{We identify the zero section of a vector bundle $V\to M$ with $M$ via $M\to V$, $p\mapsto 0_p\in V_p$.} $M\subset TM$. We say a complex structure $J_\nabla$ on $U$ is \textbf{adapted} to $\nabla$ if for every $\nabla$-geodesic $\gamma:I\to M$ the differential
 \begin{gather*}
  d\gamma: (TI,i)\to (U,J_\nabla)
 \end{gather*}
 is a holomorphic curve. Here, we view $TI$ as a subset of $\C$ by taking the base point to be the real part and the fiber component to be the imaginary part. If $\nabla = \nabla^g$ is the Levi-Civita connection of a semi-Riemannian metric $g$ on $M$, we write $J_g$ instead of $J_{\nabla^g}$ and call $J_g$ adapted to $g$.
\end{definition}

To gain some insight into this rather abstract definition, it is helpful to consider two canonical vector fields on $TM$: The vertical vector field $X^{\ver}$ and the geodesic vector field $X^{\geo}$. The vertical vector field $X^{\ver}:W\to TW$ is defined for any vector bundle $W\stackrel{\pi}{\to}M$. Recall that the vertical bundle $V\subset TW$ is the kernel of $d\pi$. For any $w\in W$, there exists a canonical isomorphism between $W_{\pi (w)}$ and $V_w$:
\begin{gather*}
 \ver_w: W_{\pi (w)}\to V_w,\quad v\mapsto \left.\frac{d}{dh}\right\vert_{h=0} w + hv.
\end{gather*}
The vertical vector field is now given by $X^{\ver}(w)\coloneqq \ver_w w$. The geodesic vector field $X^{\geo}$, on the other hand, is only defined for Koszul manifolds $(M,\nabla)$. Let $v\in TM$, $q\coloneqq \pi (v)\in M$, and $\gamma_v:I\to M$ be the unique maximal $\nabla$-geodesic with $\gamma_v (0) = q$ and $\dot\gamma_v (0) = v$. Consider the curve $\dot\gamma_v:I\to TM$. We define the geodesic vector field $X^{\geo}:TM\to T(TM)$ as follows:
\begin{gather*}
 X^{\geo} (v)\coloneqq \left.\frac{d}{dh}\right\vert_{h=0} \dot\gamma_v (h).
\end{gather*}
Coordinates $(q_1,\ldots, q_n):Q\subset M\to\R^n$ of $M$ induce coordinates\linebreak $(q_1,\ldots,q_n,v_1,\ldots,v_n):TQ\subset TM\to\R^{2n}$ of $TM$. In these coordinates, $X^{\ver}$ and $X^{\geo}$ are given by:
\begin{gather*}
 X^{\ver} (q,v) = \sum^n_{j=1} v_j \pa{v_j},\quad X^{\geo} (q,v) = \sum^n_{j=1} v_j\pa{q_j} - \sum^n_{i,j,k=1}\Gamma^k_{ij} (q)v_iv_j\pa{v_k},
\end{gather*}
where $\Gamma^k_{ij}$ are the Christoffel symbols associated with $\nabla$. It is now obvious that $X^{\ver}$ and $X^{\geo}$ vanish on the zero section $M\subset TM$, but are linearly independent for points in $TM\backslash M$. Furthermore, one easily checks with the help of these formulas that $[X^{\ver},X^{\geo}] = X^{\geo}$. Thus, $X^{\ver}$ and $X^{\geo}$ span a two-dimensional, involutive distribution on $TM\backslash M$. By the Frobenius theorem, this distribution gives rise to a foliation, called the Koszul foliation (cf. \cite{Szoeke2004}).\\
Next, we want to show that the leaves of the Koszul foliation are complex submanifolds of $(M,J_\nabla)$ if $J_\nabla$ is a complex structure adapted to $\nabla$. To do so, we pick a geodesic $\gamma:I\to M$ and introduce the map $u(t+is)\coloneqq d\gamma_t (s) = s\dot\gamma (t)$. We need to compute $\pa{t}u (t_0 + is_0)$ and $\pa{s}u (t_0 + is_0)$ for $s_0\in\R\backslash\{0\}$. It is straightforward to calculate $\pa{s}u$:
\begin{align*}
 \pa{s}u(t_0+is_0) &= \left.\frac{d}{dh}\right\vert_{h=0} (s_0+ h)\dot\gamma (t_0) = \left.\frac{d}{dh}\right\vert_{h=0} s_0\dot\gamma (t_0) + \frac{h}{s_0}s_0\dot\gamma (t_0)\\
 &= \frac{1}{s_0}\ver_{s_0\dot\gamma (t_0)} s_0\dot\gamma (t_0) = \frac{1}{s_0}X^{\ver}(u (t_0 + is_0)).
\end{align*}
To determine $\pa{t}u(t_0 + is_0)$, we first observe:
\begin{align*}
 \pa{t}u(t_0+is_0) &= \left.\frac{d}{dh}\right\vert_{h=0} s_0\dot\gamma (t_0 + h) = \left.\frac{d}{dh}\right\vert_{h=0} s_0\dot\gamma^\prime (h)\\
 &= dm_{s_0}X^{\geo}(v),
\end{align*}
where $m_{s_0}:TM\to TM, m_{s_0}(v)\coloneqq s_0 v$ is the fiberwise multiplication by $s_0$ and $\gamma^\prime (h)\coloneqq \gamma (t_0 + h)$ is the unique geodesic with $\gamma^\prime (0) = \gamma (t_0) \eqqcolon p$ and $\dot\gamma^\prime (0) = \dot\gamma (t_0) \eqqcolon v$. Now note that $\gamma^\prime_{s_0} (h)\coloneqq \gamma^\prime (s_0h)$ is the unique geodesic with $\gamma^\prime_{s_0} (0) = p$ and $\dot\gamma^\prime_{s_0} (0) = s_0 v$. With this, we find:
\begin{gather*}
 X^{\geo} (s_0v) = \left.\frac{d}{dh}\right\vert_{h=0}\dot\gamma^\prime_{s_0} (h) = \left.\frac{d}{dh}\right\vert_{h=0}s_0\dot\gamma^\prime (s_0h) = s_0dm_{s_0}X^{\geo} (v).
\end{gather*}
Thus, we obtain:
\begin{gather*}
 \pa{t}u(t_0+is_0) = \frac{1}{s_0}X^{\geo} (s_0v) = \frac{1}{s_0}X^{\geo} (u(t_0+is_0)).
\end{gather*}
If $J_\nabla$ is an adapted complex structure, then $s_0\pa{t}u = X^{\geo}$ and $s_0\pa{s}u = X^{\ver}$ imply:
\begin{gather}
 J_\nabla X^{\geo} = X^{\ver},\quad J_\nabla X^{\ver} = -X^{\geo}.\label{eq:koszul_leaves}
\end{gather}
\autoref{eq:koszul_leaves} shows that the distribution spanned by $X^{\ver}$ and $X^{\geo}$ is\linebreak $J_\nabla$-invariant which entails that the Koszul leaves are complex submanifolds.\\
The Koszul foliation helps us to answer the question of uniqueness for adapted complex structures. Imagine for a moment that $J_\nabla$ is not defined on an open neighborhood $U$ of $M\subset TM$, but on $TM\backslash M$. By our previous consideration, $J_\nabla$ is an adapted complex structure away from the zero section if and only if $J_\nabla$ preserves the Koszul distribution. In simpler terms, this means that the Koszul foliation fixes one holomorphic vector field on $(TM\backslash M,J_\nabla)$, namely $X^{\ver} + iX^{\geo}$. This vector field does not determine $J_\nabla$, however, as we still have the freedom to choose the remaining holomorphic vector fields (cf. \cite{Szoeke2004}).\\
Nevertheless, $J_\nabla$ is completely determined on the zero section $M\subset TM$. This seemingly contradicts our previous statement. Yet, being an adapted complex structure on a neighborhood of the zero section means more than just fulfilling \autoref{eq:koszul_leaves}. Indeed, $X^{\geo}$ and $X^{\ver}$ vanish on the zero section, while $\pa{t}u(t_0 + is_0)$ and $\pa{s}u(t_0+is_0)$ are not necessarily zero for $s_0 = 0$. To make this idea more precise, recall that a linear connection $\nabla$ on a vector bundle $W\to M$ gives rise to a horizontal bundle $H\subset TW$ which is complementary to the vertical bundle $V$: $TW = H\oplus V$. While the vertical bundle $V$ is canonical, the horizontal bundle $H$ is not and depends on the choice of $\nabla$. Despite that, $H_p$ is independent of $\nabla$ for points $p\in M\subset W$ of the zero section. Indeed, if $z:M\to W$ is the zero section, then $H_p$ is just $\im dz_p$ and we can identify $H_p$ with $T_pM$ via $dz_p$. Returning to $W = TM$, we see that the tangent space $T_p(TM)$ admits a canonical decomposition into $T_pM\oplus T_pM$ for points $p\in M\subset TM$:
\begin{gather*}
 T_p(TM) = H_p\oplus V_p\cong T_pM\oplus T_pM,
\end{gather*}
where we identified $V_p$ with $T_pM$ via $\ver_p:T_p M\to V_p$ as before. Now consider $\pa{t}u(t_0+is_0)$ and $\pa{s}u(t_0+is_0)$ for $s_0 = 0$. To highlight that $u$ depends on the geodesic $\gamma$, we write $u_\gamma$ instead of $u$. The following relations are easy to prove ($p\in M\subset TM$):
\begin{align*}
 H_p &= \{\pa{t}u_\gamma (t_0)\mid \gamma\text{ is a geodesic with }\gamma (t_0) = p\},\\
 V_p &= \{\pa{s}u_\gamma (t_0)\mid \gamma\text{ is a geodesic with }\gamma (t_0) = p\}.
\end{align*}
Equipped with this knowledge, we see that the adapted complex\linebreak structure $J_\nabla$ is completely determined at $p$ by the Cauchy-Riemann equations\linebreak $J_\nabla\pa{t}u_\gamma (t_0) = \pa{s}u_\gamma (t_0)$ and $J_\nabla\pa{s}u_\gamma (t_0) = -\pa{t}u_\gamma (t_0)$. In fact, we can even describe how $J_\nabla$ acts on $T_p(TM)$. Unfolding the definitions and keeping the identifications in mind, the Cauchy-Riemann equations become\footnote{This equation is reminiscent of the equation defining the almost complex structure constructed in \autoref{app:almost_complex_structures} (note that the minus sign is in a different position). We will see that the adapted complex structure $J_\nabla$ and the almost complex structure from \autoref{app:almost_complex_structures} have a lot in common, even though they are usually different objects.}:
\begin{gather}
 J_\nabla (v_1,v_2) = (-v_2,v_1)\quad\forall (v_1,v_2)\in T_p(TM)\ \forall p\in M\subset TM.\label{eq:adapted_at_zero}
\end{gather}
As Sz{\H{o}}ke shows in \cite{Szoeke2004}, the uniqueness of $J_\nabla$ at the zero section transfers to other points $p\in U$ if $p$ is connected to the zero section by a geodesic. This leads to the notion of an admissible domain (cf. \cite{Szoeke2004}):

\begin{definition}[Admissible domain]\label{def:adm_dom}
 Let $(M,\nabla)$ be a Koszul manifold. We say an open neighborhood $U\subset TM$ of $M\subset TM$ is an \textbf{admissible domain} if for every maximal $\nabla$-geodesic $\gamma:I\to M$ the set $(d\gamma)^{-1}(U)\subset\C$ is connected.
\end{definition}

\begin{remark}
 Using the notion of a felicitous domain (cf. \cite{Szoeke2004}), one easily sees that every neighborhood of $M\subset TM$ contains an admissible domain.
\end{remark}

\begin{proposition}[Uniqueness of adapted complex structures]\label{prop:unique_adapted}
 Let $(M,\nabla)$ be a Koszul manifold. If $J_{\nabla,1}$ and $J_{\nabla,2}$ are two adapted complex structures on the admissible domain $U\subset TM$, then $J_{\nabla,1} = J_{\nabla,2}$. In particular, if $J_{\nabla,1}$ and $J_{\nabla,2}$ are adapted complex structures on the neighborhoods $U_1$ and $U_2$ of $M\subset TM$, respectively, then there exists an open neighborhood $U\subset U_1\cap U_2$ of $M\subset TM$ such that $J_{\nabla,1}\vert_U = J_{\nabla,2}\vert_U$.
\end{proposition}

\begin{proof}
 This is Theorem 0.2 in \cite{Szoeke2004}.
\end{proof}

Let us now discuss existence of adapted complex structures. The following proposition ensures that $J_\nabla$ exists if $(M,\nabla)$ is real-analytic:

\begin{proposition}[Existence of adapted complex structures]\label{prop:existence_adapted}
 Let $(M,\nabla)$ be a real-analytic Koszul manifold. Then, there exists an open neighborhood $U$ of $M\subset TM$ and an adapted complex structure $J_\nabla$ on $U$. Furthermore, the fiberwise inversion $\sigma: U\to U$, $v\mapsto -v$ is a real structure on $(U,J_\nabla)$ with real form $M\subset U$.
\end{proposition}

\begin{proof}
 This statement is proven in \cite{Bielawski2003} (Proposition 1.1) and in \cite{Szoeke2004} (Theorem 0.3). Both proofs are very similar and essentially use the exponential map to construct $J_\nabla$. Still, we follow the proof of Roger Bielawski here, as it is more concise. Let $(M,\nabla)$ be a real-analytic Koszul manifold. By Bruhat-Whitney's embedding theorem, there exists a complex manifold $(X,I)$ with nice real structure $\tau$ on it (cf. \autoref{app:real_structures}) such that $\Fix\tau = M$. Since $\nabla$ is a real-analytic connection on the real form $M$, we can extend\footnote{Strictly speaking, it might only be possible to extend $\nabla$ to an open neighborhood $U$ of $M\subset X$. In this case, we replace $X$ by $U$.} $\nabla$ to a $\tau$-invariant holomorphic connection $\nabla^X$ on $X$ (cf. Lemma \autoref{lem:holo_continuation_of_connections}). Now consider a point $p\in M\subset X$ and the holomorphic tangent space $T^{(1,0)}_pX$. It carries a real structure given by $v\mapsto d\tau_p (\bar{v})$. This real structure gives rise to the decomposition
 \begin{gather*}
  T^{(1,0)}_pX = E^+_p\oplus E^-_p,
 \end{gather*}
 where $E^{\pm}_p$ is the eigenspace of the real structure w.r.t. the eigenvalue $\pm 1$. $E^+_p$ is isomorphic to $E^-_p$ via $i:E^+_p\to E^-_p$, $v\mapsto iv$ and isomorphic to $T_pM$ via $A_p:T_pM\to E^+_p$, $v\mapsto 1/2(v -iI_pv)$. Hence, $TM\cong E^+\cong E^-$.\\
 Next, we define the map $f:TM\to X$ as follows:
 \begin{gather*}
  f(v)\coloneqq \exp^X_p (iA_p v)\quad\forall v\in T_pM\ \forall p\in M,
 \end{gather*}
 where $\exp^X_p:T^{(1,0)}_pX\to X$ is the exponential map of the holomorphic connection $\nabla^X$ at the point $p\in M\subset X$ (cf. Definition \autoref{def:holo_exp_map}). $f$ satisfies $f(0_p) = p$ for every zero $0_p\in T_pM$, so it is injective on the zero section $M\subset TM$. Furthermore, we have $d\exp^X_{p,0_p} = \id_{T^{1,0}_pX}$ which entails that $f$ is a local diffeomorphism near any point of the zero section. Together, this shows that there is an open neighborhood $U\subset TM$ of $M\subset TM$ such that $f\vert_U$ is a diffeomorphism onto its image. We now define the complex structure $J_\nabla$ to be the pullback of $I$ under $f\vert_U$. It remains to be shown that $J_\nabla$ is adapted to $\nabla$. First, we note that the following relation holds:
 \begin{gather}
  \exp^X_p((z_1+z_2)v) = \exp^X_{\exp^X_p (z_2v)}\left(z_1\left.\frac{d}{dz}\right\vert_{z=z_2}\exp^X_p(zv)\right),\label{eq:bielawski_rel}
 \end{gather}
 where $p\in X$, $z_1,z_2\in\C$, and $v\in T^{(1,0)}_p X$. To see this, consider the curves:
 \begin{gather*}
  \alpha_{z_2} (z_1)\coloneqq \exp^X_p((z_1+z_2)v),\ \beta_{z_2} (z_1)\coloneqq\exp^X_{\exp^X_p (z_2v)}\left(z_1\left.\frac{d}{dz}\right\vert_{z=z_2}\exp^X_p(zv)\right).
 \end{gather*}
 One easily verifies that the curves $\alpha_{z_2}$ and $\beta_{z_2}$ as well as their first derivatives agree for $z_1 = 0$. As both $\alpha_{z_2}$ and $\beta_{z_2}$ are geodesics of $\nabla^X$, they must coincide for all $z_1$. Now let $\exp^M_p$ be the exponential map of $\nabla$ at $p\in M$. Then for any $v\in T_pM$, $\gamma (t)\coloneqq \exp^M_p (tv)$ is a geodesic of $\nabla$. We compute $f\circ u(t_0 + is_0)$, where $u(t+is)\coloneqq s\dot\gamma (t)$:\footnote{Caution is advised: The derivative in the second line is a complex derivative! That is why $A_p$ is absent in the second line.}
 \begin{align*}
  f\circ u (t_0+is_0) &= \exp^X_{\exp^M_p (t_0v)}\left(is_0 A_p\left.\frac{d}{dt}\right\vert_{t=t_0}\exp^M_p (tv)\right)\\
  &= \exp^X_{\exp^X_p (t_0A_pv)}\left(is_0\left.\frac{d}{dz}\right\vert_{z=t_0}\exp^X_p (zA_pv)\right)\\
  &= \exp^X_p\left((t_0+is_0)A_p v\right).
 \end{align*}
 Here, we have used \autoref{eq:bielawski_rel} and that $\exp^X_p$ is the holomorphic continuation of $\exp^M_p$, i.e., $\exp^M_p (tv) = \exp^X_p (tA_pv)$. Due to the holomorphicity of $\exp^X_p$, $f\circ u$ is a holomorphic curve proving that $J_\nabla$ is adapted to $\nabla$.\\
 Lastly, we show that $\sigma$ is a real structure on $U$. For this, it suffices to prove $\tau\circ f = f\circ\sigma$. Since the holomorphic connection $\nabla^X$ is $\tau$-invariant, $\exp^X_p$ is compatible with $\tau$ in the following way ($v\in T_pM$, $p\in M$):
 \begin{gather*}
  \tau\circ f (v) = \tau\circ\exp^X_p (iA_pv) = \exp^X_p (d\tau_p(\overline{iA_pv})) = \exp^X_p (-iA_p v) = f\circ\sigma (v).
 \end{gather*}
\end{proof}

Next, we consider the case where $M$ carries a semi-Riemannian metric $g$. In this scenario, $g$ furnishes a natural identification of $TM$ with $T^\ast M$. This diffeomorphism allows us to pullback the adapted complex structure $J_g$ from $T M$ to $T^\ast M$. For the sake of convenience, we also denote the resulting complex structure by $J_g$. The following theorem states that $T^\ast M$ equipped with $-\omega_{\can}$ and $J_g$ is, in fact, a semi-Kähler manifold:

\begin{theorem}[Stenzel's theorem]\label{thm:cotangent_kaehler}
 Let $(M,g)$ be a real-analytic\linebreak semi-Riemannian manifold. Furthermore, let $\sigma:T^\ast M\to T^\ast M$, $\sigma (\alpha)\coloneqq -\alpha$\linebreak be the fiberwise inversion, $\lambda_{\can}\in\Omega^1 (T^\ast M)$ the canonical one-form, and\linebreak $U_g:T^\ast M\to \R$ the function defined by $U_g (\alpha)\coloneqq g^\ast (\alpha,\alpha)$, where $g^\ast$ is the metric dual to $g$. Then, there exists an open neighborhood $U\subset T^\ast M$ of $M\subset T^\ast M$ and an almost complex structure $J_g$ on $U$ such that:
 \begin{enumerate}
  \item $J_g$ is integrable, i.e. $N_{J_g} = 0$,
  \item $\sigma$ is $J_g$-antiholomorphic, i.e. $d\sigma\circ J_g = -J_g\circ d\sigma$,
  \item $dU_g\circ J_g = 2\lambda_{\can}$.
 \end{enumerate}
 $J_g$ is unique in the following sense: If there is another almost complex structure $J^\prime_g$ on a neighborhood $U^\prime\subset T^\ast M$ of $M$ satisfying (i), (ii), and (iii), then there exists a neighborhood $U^{\prime\prime}\subset U\cap U^\prime$ of $M$ such that $\left. J^\prime_g\right\vert_{U^{\prime\prime}} = \left. J_g\right\vert_{U^{\prime\prime}}$. Under the identification $TM\cong T^\ast M$, $J_g$ becomes the complex structure on $TM$ adapted to $g$. Moreover, $(U,-\omega_{\can} = -d\lambda_{\can},J_g)$ is a semi-Kähler manifold with potential $U_g$ and signature $(2r,2s)$ of its semi-Kähler metric $g_{T^\ast M}$, where $(r,s)$ is the signature of $g$. In particular, $(U,-\omega_{\can}, J_g)$ is Kähler if $g$ is positive definite.
\end{theorem}

\begin{proof}
 This theorem is mostly due to Stenzel (cf. \cite{Stenzel1990}), but the proof given in his PhD thesis is incomplete. Consider the complex structure on $TM$ adapted to $g$ (exists due to Proposition \autoref{prop:existence_adapted}) and pull it back to $T^\ast M$ via $v\mapsto g(v,\cdot)$. Call the result $J_g$. By Proposition \autoref{prop:existence_adapted}, $J_g$ is defined on an open neighborhood $U\subset T^\ast M$ of the zero section $M\subset T^\ast M$ and satisfies (i) and (ii). (iii) is shown in \cite{Lempert1991} (cf. Corollary 5.5) and \cite{Szoeke2004} (cf. Theorem 0.6\footnote{The proof of Theorem 0.6 is correct. However, Equation (0.4) is false: It should be ``$2\text{Im}\bar{\partial}E = \Theta$'' instead of ``$\bar{\partial}E = i\Theta$''.}). Uniqueness is proven by Stenzel using his formal power series\footnote{Stenzel only formulated his power series for Riemannian manifolds $(M,g)$. Nevertheless, it still make sense for semi-Riemannian manifolds $(M,g)$.} (cf. \cite{Stenzel1990}). To show that $(U,-\omega_{\can},J_g)$ is a semi-Kähler manifold, we first observe:
 \begin{align*}
  i\partial\bar{\partial} U_g &= i(\partial + \bar{\partial})\bar{\partial} U_g = id\bar{\partial} U_g = \frac{i}{2}d(dU_g + idU_g\circ J_g) = \frac{i^2}{2}d\left(dU_g\circ J_g\right)\\
  &\stackrel{\text{(iii)}}{=} -d\lambda_{\can} = -\omega_{\can}.
 \end{align*}
 It is a standard result from complex geometry that every two-form $\omega$ on a complex manifold $(X,J)$ satisfying $\omega = i\partial\bar{\partial}f$ for some function $f\in C^\infty (X,\R)$ fulfills $\omega (J\cdot,J\cdot) = \omega$ (in this case, we say $f$ is the potential for $\omega$). Thus, $(U,-\omega_{\can}, J_g)$ is a semi-Kähler manifold with potential $U_g$.\\
 Lastly, we have to determine the signature of the semi-Kähler metric. By shrinking $U$ if necessary, we can assume that every point in $U$ is connected to the zero section $M\subset T^\ast M$. Since the signature of a semi-Riemannian metric does not change within a connected component, it suffices to compute the signature on $M\subset T^\ast M$. As before, we have the canonical decomposition $T_p(T^\ast M) = H_p\oplus V_p$ for points $p\in M$. While the horizontal space $H_p$ is still isomorphic to $T_p M$ via the zero section, the vertical space $V_p$ is now isomorphic to $T^\ast_p M$. Hence, we have $T_p (T^\ast M)\cong T_p M\oplus T^\ast_p M$. Under this isomorphism, $-\omega_{\can}$ and $J_g$ become (cf. \autoref{eq:adapted_at_zero}):
 \begin{gather*}
  -\omega_{\can} \left((v_1,\alpha_1), (v_2,\alpha_2)\right) = \alpha_2 (v_1) - \alpha_1 (v_2),\qquad J_g (v,\alpha) = (-\#(\alpha), b(v)),
 \end{gather*}
 where $b$ and $\#$ are the musical maps associated with $g$. Thus, the semi-Kähler metric $g_{T^\ast M}\coloneqq -\omega_{\can} (\cdot,J_g\cdot)$ takes the following form on $M\subset T^\ast M$:
 \begin{gather*}
  g_{T^\ast M} \left((v_1,\alpha_1), (v_2,\alpha_2)\right) = b(v_2)(v_1) + \alpha_1 (\#(\alpha_2)) = g(v_1,v_2) + g^\ast (\alpha_1,\alpha_2).
 \end{gather*}
 Both $g$ and $g^\ast$ have signature $(r,s)$ concluding the proof.
\end{proof}

Before we turn our attention to the main result of this subsection, it is instructive to consider an important class of adapted complex structures:

\begin{example}[Adapted complex structures for flat $g$]\label{ex:flat_adapted}
 In \autoref{app:almost_complex_structures}, we construct the almost complex structure $J^\ast_{\nabla^g}$ on $T^\ast M$ which bears great resemblance to the adapted complex structure $J_g$: With the formulas given in Theorem \autoref{thm:alm_cpx_str_on_tan}, one easily verifies that $-J^\ast_{\nabla^g}$ satisfies Property (ii) and (iii) of Theorem \autoref{thm:cotangent_kaehler}. Furthermore, $-J^\ast_{\nabla^g}$ is always compatible with $-\omega_{\can}$ (cf. \autoref{app:almost_complex_structures}). $-J^\ast_{\nabla^g}$ even has an advantage over $J_g$ when it comes to existence: While existence of $J_g$ is only ensured on a neighborhood $U$ of $M\subset T^\ast M$ and only for real-analytic metrics $g$, $-J^\ast_{\nabla^g}$ can be defined on all of $T^\ast M$ and for all metrics $g$. However, this comes with a trade-off: $-J^\ast_{\nabla^g}$ is usually not integrable. Indeed, as stated in Theorem \autoref{thm:alm_cpx_str_on_tan}, $-J^\ast_{\nabla^g}$ is integrable if and only if $g$ is flat. In this case, $-J^\ast_{\nabla^g}$ fulfills all properties from Theorem \autoref{thm:cotangent_kaehler} and, therefore, agrees with $J_g$. Since we have explicit formulas for $J^\ast_{\nabla^g}$ (cf. Theorem \autoref{thm:alm_cpx_str_on_tan}), this observation allows us to compute $J_g$ directly if $g$ is flat.
\end{example}

Let us now assume that the semi-Riemannian manifold $(M,g)$ carries a complex structure $I$. In this case, $I$ induces an additional complex structure $T^\ast I$ on $T^\ast M$. Naturally, we are drawn to the question how $J_g$ and $T^\ast I$ interact with each other. The next proposition suggests that the commutation relation of $J_g$ and $T^\ast I$ depends on whether $I$ is compatible or anticompatible with $g$:

\begin{lemma}[Commutation relations]\label{lem:commutation_relation}
 Let $(M,g)$ be a real-analytic semi-Riemannian manifold and $I$ a complex structure on $M$. Furthermore, let $T^\ast I$ be the complex structure on $T^\ast M$ induced by $I$ and $J_g$ the complex structure adapted to $g$ defined on a neighborhood $U$ of $M\subset T^\ast M$. Then, $T^\ast I$ and $J_g$ satisfy the following commutation relations on the zero section $M\subset T^\ast M$:
 \begin{enumerate}
  \item $T^\ast I$ and $J_g$ commute, i.e. $T^\ast I\circ J_g = J_g\circ T^\ast I$, if $G\coloneqq g-ig(I\cdot,\cdot)$ is a holomorphic $\C$-bilinear two-form on $M$.
  \item $T^\ast I$ and $J_g$ anticommute, i.e. $T^\ast I\circ J_g = -J_g\circ T^\ast I$, if $(M,g,I)$ is a semi-Kähler manifold.
 \end{enumerate}
 If $g$ is flat, then the commutation relations hold for all points in $U$ (after shrinking $U$ if necessary).
\end{lemma}

\begin{proof}
 The idea is to find coordinates of $T^\ast M$ in which $T^\ast I$ and $J_g$ take the following form:
 \begin{gather}
  T^\ast I = \begin{pmatrix} I_0 & 0\\ 0 & -I_0\end{pmatrix} ,\quad J_g = \begin{pmatrix} 0 & -G^{-1}_0\\ G_0 & 0\end{pmatrix}.\label{eq:matrix_notation}
 \end{gather}
 Here, $I_0$ and $G_0$ are $2n\times 2n$-matrices ($n\coloneqq\dim_\C M$). $I_0$ is always the standard complex structure on $\R^{2n}$, i.e.:
 \begin{gather*}
  I_0 = \begin{pmatrix} 0 & -\mathds{1}_{n\times n}\\ \mathds{1}_{n\times n} & 0\end{pmatrix},
 \end{gather*}
 whereas $G_0$ differs between Statement (i) and (ii). $I_0$ and $G_0$ anticommute for Statement (i), while they commute for Statement (ii). Given \autoref{eq:matrix_notation}, it is trivial to verify the commutation relations:
 \begin{align*}
  T^\ast I\circ J_g &= \begin{pmatrix} I_0 & 0\\ 0 & -I_0\end{pmatrix}\begin{pmatrix} 0 & -G^{-1}_0\\ G_0 & 0\end{pmatrix} = \begin{pmatrix} 0 & -I_0G^{-1}_0\\ -I_0G_0 & 0\end{pmatrix}\\
  &= \begin{pmatrix} 0 & \pm G^{-1}_0I_0\\ \pm G_0I_0 & 0\end{pmatrix} = \pm\begin{pmatrix} 0 & -G^{-1}_0\\ G_0 & 0\end{pmatrix}\begin{pmatrix} I_0 & 0\\ 0 & -I_0\end{pmatrix}\\
  &= \pm J_g\circ T^\ast I,
 \end{align*}
 Thus, it suffices to show \autoref{eq:matrix_notation} in order to prove Lemma \autoref{lem:commutation_relation}. All coordinates of $T^\ast M$ we construct are induced by holomorphic coordinates of $M$. In this case, $T^\ast I$ automatically fulfills \autoref{eq:matrix_notation}. To see this, consider holomorphic coordinates $\psi = (q_1,\ldots, q_n):V\subset M\to \C^n$ of $M$. They give rise to holomorphic coordinates $\Psi:T^{\ast, (1,0)}V\to \C^{2n}$ of $T^{\ast, (1,0)}M$, namely:
 \begin{gather*}
  \Psi^{-1} (\tilde q_1,\ldots, \tilde q_n, \tilde p_1,\ldots, \tilde p_n) \coloneqq \sum^n_{j=1}\tilde p_j dq_{j,\psi^{-1}(\tilde q_1,\ldots, \tilde q_n)}.
 \end{gather*}
 Now consider the diffeomorphism $f:T^{\ast, (1,0)}M\to T^\ast M$, $f(\alpha)\coloneqq\text{Re} (\alpha)$. We define $T^\ast I$ by declaring $f$ to be biholomorphic. This means that $f$ sends holomorphic charts of $T^{\ast, (1,0)}M$ to holomorphic charts of $T^\ast M$. Decomposing $q_j = q_{x,j} + iq_{y,j}$ and $\tilde p_j = \tilde p_{x,j} + i\tilde p_{y,j}$ into real and imaginary part, we obtain:
 \begin{gather*}
  \beta = f\circ\Psi^{-1}(\tilde q_1,\ldots, \tilde q_n, \tilde p_1,\ldots, \tilde p_n) = \sum^n_{j=1} \tilde p_{x,j}dq_{x,j} - \tilde p_{y,j}dq_{y,j},
 \end{gather*}
 where we dropped the base point $\psi^{-1}(\tilde q_1,\ldots, \tilde q_n)$ for the sake of clarity. Comparing this to the coordinates $(q_{x,j}, q_{y,j}, p_{x,j}, p_{y,j}):T^\ast V\to \R^{4n}$ of $T^\ast M$ induced by $\psi$, i.e.
 \begin{gather*}
  \beta = \sum^n_{j=1} p_{x,j}dq_{x,j} + p_{y,j}dq_{y,j},
 \end{gather*}
 we realize that the $p_{y,j}$-components pick up an additional minus sign to become holomorphic. Thus, the complex structure $T^\ast I$ can be written as:
 \begin{alignat*}{3}
  T^\ast I (\pa{q_{x,j}}) &= \pa{q_{y,j}},\quad && T^\ast I (\pa{q_{y,j}}) &&= -\pa{q_{x,j}},\\
  T^\ast I (\pa{p_{x,j}}) &= -\pa{p_{y,j}},\quad && T^\ast I (\pa{p_{y,j}}) &&= \pa{p_{x,j}}.
 \end{alignat*}
 In matrix notation, this is just \autoref{eq:matrix_notation}.\\
 The next step is to find ``good'' coordinates for $J_g$. We have to distinguish two cases:\\
 
 \textbf{Case 1:} $J_g$ on the zero section $M\subset T^\ast M$\\
 
 Fix a point $p\in M$ and pick holomorphic coordinates $(q_1,\ldots, q_n):V\to\C^n$ of $M$ near $p$. We want to show that, in the coordinates $(q_{x,j}, q_{y,j}, p_{x,j}, p_{y,j})$ induced by $(q_1,\ldots, q_n)$, $J_g$ is described by \autoref{eq:matrix_notation} at the point $p\in M\subset T^\ast M$ of the zero section. We recall that, similar to \autoref{eq:adapted_at_zero}, $J_g$ takes the following simple form on the zero section:
 \begin{gather*}
  J_g (v,\alpha) = (-\#(\alpha), b(v)) = (-b^{-1}(\alpha), b(v)),
 \end{gather*}
 where $\#$ and $b$ are the musical maps associated with $g$ and we used the canonical decomposition $T_p(T^\ast M) = H_p\oplus V_p\cong T_pM\oplus T^\ast _pM$ into horizontal and vertical space.
 Under the given identifications, the musical map $b$ turns into the matrix $G_0$. We check that $I_0$ and $G_0$ satisfy the appropriate commutation relations. $g(I\cdot,I\cdot) = \pm g$ implies:
 \begin{gather*}
  -I_0G_0I_0 = I^t_0G_0I_0 = \pm G_0\quad\Rightarrow G_0 I_0 = \pm I_0 G_0.
 \end{gather*}
 For the sake of clarity, we shall give explicit formulas for $G_0$. Given the setup of (i), we can always achieve that $G = g - ig(I\cdot,\cdot)$ assumes its standard form at $p\in M$, i.e.
 \begin{gather*}
  G_p = \sum^n_{j=1} dq^2_{j,p}\quad\Rightarrow g_p = \sum^n_{j=1} d_p q^2_{x,j} - d_p q^2_{y,j},
 \end{gather*}
 by applying a $\C$-linear transformation to $(q_1,\ldots,q_n)$. We can now directly read off $G_0$:
 \begin{gather}
  G_0 = G^{-1}_0 = \begin{pmatrix}\mathds{1}_{n\times n} & 0\\ 0 & -\mathds{1}_{n\times n}\end{pmatrix}.\label{eq:standard_G_0_i}
 \end{gather}
 \autoref{eq:standard_G_0_i} allows us to explicitly verify $I_0G_0 = -G_0I_0$.\\
 For Statement (ii), we choose holomorphic coordinates $(q_1,\ldots, q_n)$ such that $2g_p(\pa{q_j, p},\overline{\pa{q_k, p}}) = \pm \delta_{jk}$. We can always achieve this by applying the Gram-Schmidt procedure to the basis $\{\pa{q_{x,j}}\}$ of the complex vector space $(T_pM,I_p)$ equipped with the semi-Hermitian metric $g_p + i g_p (I_p\cdot,\cdot)$. $\{\pa{q_{x,j}}$, $\pa{q_{y,j}}\}$ now forms a Sylvester basis of $(T_pM,g_p)$, where $g_p(\pa{q_{x,j}},\pa{q_{x,j}}) = g_p(\pa{q_{y,j}},\pa{q_{y,j}}) = \pm 1$, i.e.:
 \begin{gather*}
   g_p = \sum^r_{j=1} d_p q^2_{x,j} + d_p q^2_{y,j} - \sum^{n = r+s}_{j=r+1} d_p q^2_{x,j} + d_p q^2_{y,j}.
 \end{gather*}
 Here, $(2r,2s)$ is the signature of $g$. Hence, $G_0$ takes the following form:
 \begin{gather}
  G_0 = G^{-1}_0 = \begin{pmatrix}G^\prime_0 & 0\\ 0 & G^\prime_0\end{pmatrix}\quad\text{with}\quad G^\prime_0 = \begin{pmatrix}\mathds{1}_{r\times r} & 0\\ 0 & -\mathds{1}_{s\times s}\end{pmatrix}.\label{eq:standard_G_0_ii}
 \end{gather}
 Given these formulas, it is trivial to check that $I_0$ and $G_0$ commute.\\
 
 \textbf{Case 2:} $g$ is flat\\
 
 Recall that for flat $g$ the adapted complex structure $J_g$ is just the complex structure $-J^\ast_{\nabla^g}$ from \autoref{app:almost_complex_structures} (after shrinking $U$ if necessary). Theorem \autoref{thm:alm_cpx_str_on_tan} now tells us that $J^\ast_{\nabla^g}$ takes a simple form in normal coordinates of $(M,g)$. We would like to choose the normal coordinates in such a way that they are also holomorphic with respect to $I$ in order to preserve the form of $T^\ast I$.\\
 For Statement (i), the theory of holomorphic connections, developed in \autoref{app:holo_connection}, guarantees that holomorphic normal coordinates exist: By Lemma \autoref{lem:normal_of_h_and_h_R}, there are holomorphic coordinates $(q_1,\ldots, q_n):V\subset M\to \C^n$ near any point $p\in M$ such that the corresponding real coordinates $(q_{x,j},q_{y,j}):V\to\R^{2n}$ are normal coordinates of $(M,g)$ near $p$. In these coordinates, $g$ takes the form $g_p = \sum^n_{j=1} d_pq^2_{x,j} - d_pq^2_{y,j}$ at $p$. We can now apply Theorem \autoref{thm:alm_cpx_str_on_tan} to express $J_g = -J^\ast_{\nabla^g}$ in the coordinates $(q_{x,j}, q_{y,j}, p_{x,j}, p_{y,j}):T^\ast V\to \R^{4n}$ induced by $(q_{x,j},q_{y,j})$:
 \begin{alignat*}{3}
  J_g (\pa{q_{x,j}}) &= \pa{p_{x,j}},\quad && J_g (\pa{p_{x,j}}) &&= -\pa{q_{x,j}},\\
  J_g (\pa{q_{y,j}}) &= -\pa{p_{y,j}},\quad && J_g (\pa{p_{y,j}}) &&= \pa{q_{y,j}}.
 \end{alignat*}
 In matrix notation, this is just \autoref{eq:matrix_notation} with $G_0$ from \autoref{eq:standard_G_0_i}.\\
 For semi-Kähler manifolds (Statement (ii)), the existence of holomorphic normal coordinates is a well-known fact (cf. Lemma \autoref{lem:kaehler} or Theorem 4.17\footnote{This theorem only shows existence of holomorphic normal coordinates for Kähler manifolds. However, the same proof still applies in the semi-Kähler case. Also note that there are subtle differences between Riemannian and Kählerian normal coordinates. Nevertheless, they coincide in the flat case (cf. \cite{Bochner1947} and \cite{Jentsch2017} for details).} in \cite{Ballmann2006}). We can again use Theorem \autoref{thm:alm_cpx_str_on_tan} to express $J_g$ in these coordinates. If done correctly, one sees that $J_g$ assumes the form given in \autoref{eq:matrix_notation}, where $G_0$ is given by \autoref{eq:standard_G_0_ii} concluding the proof.
\end{proof}

In light of Lemma \autoref{lem:commutation_relation}, one might expect that the commutation relations hold everywhere on $U$ for all suitable metrics $g$. However, results by Dancer and Sz{\H{o}}ke (cf. \cite{Dancer1997}) suggest that the commutation relations are only satisfied if we modify $J_g$ by an appropriate diffeomorphism $\phi:T^\ast M\to T^\ast M$. The paper \cite{Dancer1997} together with Lemma \autoref{lem:commutation_relation} motivates the following conjecture:

\begin{conjecture}[Commutation relations]\label{con:commutation_relation}
 Let $(M,g)$ be a real-analytic semi-Riemannian manifold with signature $(r,s)$, $J_g$ the complex structure adapted to $g$ defined on a neighborhood $U$ of $M\subset T^\ast M$, and $I$ a complex structure on $M$ inducing the complex structure $T^\ast I$ on $T^\ast M$. Then, there exists a real-analytic, fiber-preserving diffeomorphism $\phi:U\to U$ (shrink $U$ if needed) with $\phi\vert_M = \id_M$ such that $-\omega_{\can}(\cdot,\phi^\ast J_g\cdot)$ is a metric with signature $(2r,2s)$ and:
 \begin{enumerate}
  \item $T^\ast I$ and $\phi^\ast J_g$ commute, i.e. $T^\ast I\circ \phi^\ast J_g = \phi^\ast J_g\circ T^\ast I$, if $G\coloneqq g-ig(I\cdot,\cdot)$ is a holomorphic $\C$-bilinear two-form on $M$. If, additionally, $(M,I)$ admits a real structure $\sigma$ with $\sigma^\ast g = g$, then $\phi$ can be chosen in such a way that the induced real structure $\sigma^\ast (\alpha)\coloneqq \alpha\circ d\sigma$ on $T^\ast M$ is $\phi^\ast J_g$-holomorphic and the signature of $-\omega_{\can}(\cdot,\phi^\ast J_g\cdot)$ restricted to the real form $\Fix\sigma^\ast$ is $(2r^\prime,2s^\prime)$, where $(r^\prime,s^\prime)$ is the signature of $g$ restricted to $\Fix\sigma$.
  \item $T^\ast I$ and $\phi^\ast J_g$ anticommute, i.e. $T^\ast I\circ \phi^\ast J_g = -\phi^\ast J_g\circ T^\ast I$, if $(M,g,I)$ is a semi-Kähler manifold.
 \end{enumerate}
\end{conjecture}

\begin{remark}
 The second statement in Case (i) concerning real structures may seem odd at first, especially since there is no analogous statement in Lemma \autoref{lem:commutation_relation}. Still, we include it here, because we will need this statement later on for the proof of Lemma \autoref{lem:stenzel_for_holo_metrics}. As a sanity check, one can consider the case $\phi = \id_U$. In this case, it immediately follows from the definition of $J_g$ as well as the fact that $\sigma$ is an isometry that $\sigma^\ast$ is $J_g$-holomorphic and that $-\omega_{\can}(\cdot,J_g\cdot)$ restricted to the real form $\Fix\sigma^\ast$ has the appropriate signature.
\end{remark}

As we will see in the next subsection, Conjecture \autoref{con:commutation_relation} allows us to show that double cotangent bundles exhibit Kähler duality in a natural way. Unfortunately, the author does not have a complete proof for Conjecture \autoref{con:commutation_relation} as of the time of writing.

\subsection*{Kähler Duality on $T^\ast(T^\ast M)$}

In the previous subsection, we have seen that the cotangent bundle $T^\ast M$ of a Riemannian manifold $(M,g)$ admits a unique Kähler structure. We now want to tackle the question how this Kähler structure changes if the base manifold $M$ carries additional structures. Here, two cases are particularly interesting for us: In the first case, $M$ itself possesses a Kähler structure, while the metric on $M$ is holomorphic in the second case. The following two lemmata show that $T^\ast M$ is Hyperkähler in the first case and holomorphic Kähler in the second:

\begin{lemma}[Stenzel's theorem for Kähler manifolds]\label{lem:stenzel_for_kaehler}
 Let $(X,\omega, J)$ be a real-analytic semi-Kähler manifold with metric $g\coloneqq\omega (\cdot,J\cdot)$ and let $J_g$ be the complex structure adapted to $g$ defined on a neighborhood $U\subset T^\ast X$\linebreak of $X\subset T^\ast X$. If Conjecture \autoref{con:commutation_relation} holds, then $(U,-\omega_{\can}, \phi^\ast J_g, T^\ast J)$ is\linebreak a semi-Hyperkähler manifold, where $\phi:U\to U$ is the diffeomorphism from Conjecture \autoref{con:commutation_relation}. If, additionally, $g$ is positive definite, then $(U,-\omega_{\can}, \phi^\ast J_g, T^\ast J)$ is Hyperkähler.
\end{lemma}

\begin{proof}
 By Lemma \autoref{lem:hyper_vs_holo}, $(X,\omega,J,I)$ is a semi-Hyperkähler manifold if $I$ and $J$ are complex structures on $X$ such that they anticommute, $I$ is anticompatible with the symplectic form $\omega$, and $J$ is compatible with $\omega$. Clearly, $-\omega_{\can}$ is a symplectic form and both $\phi^\ast J_g$ and $T^\ast J$ are complex structures. If Conjecture \autoref{con:commutation_relation} is true, then $\phi^\ast J_g$ and $T^\ast J$ anticommute. The compatibility of $-\omega_{\can}$ with $\phi^\ast J_g$ also follows from Conjecture \autoref{con:commutation_relation}. The anticompatibility of $-\omega_{\can}$ with $T^\ast J$ is a direct consequence of the following observation: Given a complex manifold $(X,J)$, the space $T^{\ast,(1,0)}X$ possesses a canonical holomorphic symplectic form $\Omega_{\can}$. Under the biholomorphism $T^{\ast, (1,0)}X\to T^\ast X$, $\alpha\mapsto\text{Re}(\alpha)$, the real part of $\Omega_{\can}$ becomes $\omega_{\can}$ which implies that $-\omega_{\can}$ is anticompatible with $T^\ast J$. We conclude the proof by noting that $g_{T^\ast X}\coloneqq -\omega_{\can}(\cdot, \phi^\ast J_g\cdot)$ is positive definite if $g$ is positive definite (cf. Conjecture \autoref{con:commutation_relation}).
\end{proof}

It has already been established by Kaledin (cf. \cite{Kaledin1997}) and Feix (cf. \cite{Feix2001}) that cotangent bundles of Kähler manifolds admit Hyperkähler structures. Their approach differs greatly from ours, so one might wonder whether our\linebreak Hyperkähler structure is isomorphic to the one constructed by Feix and Kaledin.\linebreak A recently released preprint by Su-Jen Kan (cf. \cite{Kan2024}) suggests that they are not isomorphic. However, we will not pursue this idea any further.
Let us now turn our attention to the second case:

\pagebreak

\begin{lemma}[Stenzel's theorem for holomorphic metrics]\label{lem:stenzel_for_holo_metrics}
 Let $(X,J)$ be a complex manifold with holomorphic metric $G = g-ig(J\cdot,\cdot)$ and let $J_g$ be the complex structure adapted to $g$ defined on a neighborhood $U\subset T^\ast X$ of\linebreak $X\subset T^\ast X$. If Conjecture \autoref{con:commutation_relation} is true, then $(U,-\omega_{\can}, \phi^\ast J_g, T^\ast J)$ is a holomorphic semi-Kähler manifold, where $\phi:U\to U$ is the map from Conjecture \autoref{con:commutation_relation}. If, additionally, $\sigma$ is a real structure on $X$ with $\sigma^\ast g = g$ and the induced metric $g_M$ on the real form $M\coloneqq\Fix\sigma$ is positive definite, then $(U,-\omega_{\can}, \phi^\ast J_g, T^\ast J)$ is even holomorphic Kähler.
\end{lemma}

\begin{proof}
 The proof of the first statement in Lemma \autoref{lem:stenzel_for_holo_metrics} is quite similar to the proof of Lemma \autoref{lem:stenzel_for_kaehler}: By Lemma \autoref{lem:hyper_vs_holo}, $(X,\omega,J,I)$ is a holomorphic semi-Kähler manifold if $I$ and $J$ are complex structures on $X$ such that they commute, $I$ is anticompatible with the symplectic form $\omega$, and $J$ is compatible with $\omega$. Most properties can be shown as in Lemma \autoref{lem:stenzel_for_kaehler}. Only the commutation relation for $\phi^\ast J_g$ and $T^\ast J$ differs which is accounted for by Conjecture \autoref{con:commutation_relation}.\\
 To verify the second assertion, we first have to find a real structure of the holomorphic semi-Kähler manifold $(U,-\omega_{\can}, \phi^\ast J_g, T^\ast J)$. It is easy to prove, for instance by going into holomorphic charts, that $\sigma^\ast:T^\ast X\to T^\ast X$, $\alpha\mapsto \alpha\circ d\sigma$\linebreak is a real structure of $(T^\ast X, T^\ast J)$ with real form $T^\ast M$. If $U$ and $-\omega_{\can}$ are invariant under $\sigma^\ast$ and $\sigma^\ast$ is $\phi^\ast J_g$-holomorphic, then $\sigma^\ast$ is also a real structure\linebreak of $(U,-\omega_{\can}, \phi^\ast J_g,T^\ast J)$. By Conjecture \autoref{con:commutation_relation}, we know that $\sigma^\ast$ is\linebreak $\phi^\ast J_g$-holomorphic. After shrinking $U$, we can assume $\sigma^\ast (U) = U$. To show that $\sigma^\ast$ preserves $-\omega_{\can}$, recall the following fact: If $f:X\to X$ is a diffeomorphism of a manifold $X$, then $f$ gives rise to the following symplectomorphism of $(T^\ast X,\omega_{\text{can}})$:
 \begin{gather*}
  f_{T^\ast X}:T^\ast X\to T^\ast X,\quad \alpha\in T^\ast_p X\mapsto \alpha\circ df^{-1}\in T^\ast_{f(p)}X.
 \end{gather*}
 $\sigma^\ast$ is a map of that form and, therefore, leaves $-\omega_{\can}$ invariant.\\
 As shown in \autoref{app:kaehler}, the real form of a holomorphic semi-Kähler manifold is itself a semi-Kähler manifold meaning that the real form of $(U,-\omega_{\can}, \phi^\ast J_g, T^\ast J)$ carries a semi-Kähler structure. Since we assume the metric $g_M$ to be positive definite, the Kähler metric of the real form $U\cap T^\ast M$ is positive definite as well (cf. Conjecture \autoref{con:commutation_relation}). This shows that $(U,-\omega_{\can}, \phi^\ast J_g, T^\ast J)$ is a holomorphic Kähler manifold finishing the proof.
\end{proof}

%

Lemma \autoref{lem:stenzel_for_kaehler} and \autoref{lem:stenzel_for_holo_metrics} now allow us to construct Hyperkähler and holomorphic Kähler structures on double cotangent bundles of real-analytic Riemannian manifolds $(M,g)$. The idea is depicted in the following diagram:\footnote{Both $\phi_1$ and $\phi_2$ denote $\phi$-maps as in Conjecture \autoref{con:commutation_relation}. The different indices merely emphasize the point that we do not need to choose the same $\phi$-map in the Hyperkähler and holomorphic Kähler case.}

\begin{center}
 \begin{tikzcd}[column sep=3.6em,row sep=2em]
  \textbf{Hyperk.}& (T^\ast M, g_{T^\ast M}) \arrow{r}{\text{Stenzel}} & (T^\ast (T^\ast M), -\omega_{\can}, \phi_1^\ast J_{g_{T^\ast M}}, T^\ast J_g)\\
  (M,g) \arrow{ur}[sloped,above]{\text{Stenzel}} \arrow{dr}[sloped,below]{\text{\footnotesize Complexification}} & &\\
  \textbf{Holo. K.}& (T^\ast M, g_\C) \arrow{r}[below]{\text{Stenzel}} & (T^\ast (T^\ast M), -\omega_{\can}, \phi_2^\ast J_{g_\C}, T^\ast J_g)
 \end{tikzcd}
\end{center}
The construction is based on the observation that there are two natural metrics on $T^\ast M$: One is the (real-analytic) Kähler metric $g_{T^\ast M}$ obtained from Stenzel's theorem (Theorem \autoref{thm:cotangent_kaehler}). To define the other metric, recall that $M$ is the real form of the complex manifold $(T^\ast M,J_g)$ (cf. Proposition \autoref{prop:existence_adapted}). $g$ is a real-analytic metric on the real form $M$, so, by Lemma \autoref{lem:holo_continuation_of_tensors}, there exists a holomorphic continuation of $g$ on $T^\ast M$. Let $g_\C$ be the real part of this holomorphic continuation. Clearly, $g_\C$ is a real-analytic semi-Riemannian metric on $T^\ast M$. Thus, we can apply Lemma \autoref{lem:stenzel_for_kaehler} to $(T^\ast M, g_{T^\ast M})$ and Lemma \autoref{lem:stenzel_for_holo_metrics} to $(T^\ast M, g_\C)$ showing:

\begin{theorem}[Kähler duality on $T^\ast (T^\ast M)$]\label{thm:kaehler_duality}
 Let $(M,g)$ be a real-analytic\linebreak semi-Riemannian manifold, $J_g$ the complex structure adapted to $g$, and $g_{T^\ast M}$ the semi-Kähler metric from Theorem \autoref{thm:cotangent_kaehler}. Furthermore, let $g_\C$ be the real part of the holomorphic continuation of $g$, $T^\ast J_g$ the complex structure induced by $J_g$, and let $J_{g_{T^\ast M}}$ and $J_{g_\C}$ be the complex structures adapted to $g_{T^\ast M}$ and $g_\C$, respectively, defined on a neighborhood $U\subset T^\ast (T^\ast M)$ of $M$. If Conjecture \autoref{con:commutation_relation} holds, then $(U,-\omega_{\can}, \phi^\ast_1 J_{g_{T^\ast M}}, T^\ast J_g)$ is a semi-Hyperkähler manifold and $(U,-\omega_{\can}, \phi^\ast_2 J_{g_\C}, T^\ast J_g)$ is a holomorphic semi-Kähler manifold. We can drop the prefix ``semi'' if $g$ is positive definite.
\end{theorem}

Hyperkähler structures on double (co)tangent bundles have already been observed by Bielawski: In \cite{Bielawski2003}, he uses adapted complex structures and a twistor construction to equip $TTM$ with a Hyperkähler structure. It is possible that the first part of Theorem \autoref{thm:kaehler_duality} reproduces Bielawski's construction after a suitable identification of $TTM$ with $T^\ast (T^\ast M)$. We shall not investigate this idea any further, but rather apply Theorem \autoref{thm:kaehler_duality} to Lie groups $M = G_\R$ with bi-invariant metrics $g$.

\subsection*{Application: $M = G_\R$}

The purpose of this subsection is to describe the results of Theorem \autoref{thm:cotangent_kaehler} and \autoref{thm:kaehler_duality} when applied to Lie groups. In particular, we want to work out that Stenzel's theorem applied to groups is a generalization of Bremigan's construction (cf. \cite{Bremigan2000}), that the universal complexification of a compact group carries a family of Kähler structures (cf. Lemma \autoref{lem:kaehler_reductive_groups}), and that cotangent bundles of such groups exhibit Kähler duality (cf. Theorem \autoref{thm:kaehler_duality_lie}). This leads us to the conjecture that the Kähler duality on cotangent bundles reduces to the Kähler duality on coadjoint orbits (cf. Conjecture \autoref{con:reduction}). We start by discussing connections on Lie groups:

\begin{definition}[Left-invariant connection]\label{def:left_connection}
 Let $G$ be a Lie group and $\nabla$ a connection on $G$. We call $\nabla$ \textbf{left-invariant} if $\nabla_X Y$ is left-invariant for all left-invariant vector fields $X,Y$ on $G$.
\end{definition}

The following properties of left-invariant connections are fairly obvious:

\pagebreak

\begin{proposition}\label{prop:properties_of_left_con}
 Let $G$ be a Lie group.
 \begin{enumerate}
  \item The set of left-invariant connections $\nabla$ on $G$ is isomorphic to the set of bilinear maps $\alpha_\nabla:\mathfrak{g}\times\mathfrak{g}\to\mathfrak{g}$.
  \item The map $\alpha_\nabla$ associated with a left-invariant connection $\nabla$ on $G$ decomposes into a symmetric part $\alpha^S_\nabla$ and an antisymmetric part $\alpha^A_\nabla$:
  \begin{align*}
   \alpha^S_\nabla (X,Y) &= \frac{1}{2} \left(\alpha_\nabla (X,Y) + \alpha_\nabla (Y,X)\right) = \frac{1}{2}\left(\nabla_X Y + \nabla_Y X\right),\\
   \alpha^A_\nabla (X,Y) &= \frac{1}{2} \left(\alpha_\nabla (X,Y) - \alpha_\nabla (Y,X)\right) = \frac{1}{2}\left(\nabla_X Y - \nabla_Y X\right).
  \end{align*}
  \item The symmetric part $\alpha^S_\nabla$ of a left-invariant connection $\nabla$ vanishes if and only if all immersed $1$-dimensional Lie subgroups of $G$ are geodesics of $\nabla$.
  \item A left-invariant connection $\nabla$ is torsion-free if and only if the following equation holds:
  \begin{gather*}
   \alpha^A_\nabla (X,Y) = \frac{1}{2}[X,Y]\quad\forall X,Y\in\mathfrak{g}.
  \end{gather*}
 \end{enumerate}
\end{proposition}

\begin{proof}
 For (i), let $X_i$ be a basis of the space of left-invariant vector fields. Observe that every connection $\nabla$ defined by
 \begin{gather*}
  \nabla_{X_i} X_j \coloneqq \sum^n_{k=1} \Gamma^k_{ij} X_k
 \end{gather*}
 with constant Christoffel symbols $\Gamma^k_{ij}\in\R$ is left-invariant and that every left-invariant connection can be written in such a way. Identifying $\mathfrak{g}$ with the space of left-invariant vector fields now shows (i). (ii) is trivial. To prove (iii), recall that every immersed $1$-dimensional Lie subgroup is the integral curve of some left-invariant vector field through the neutral element $e\in G$. Let $\gamma:\R\to G$ be the integral curve of the left-invariant vector field $X$ with $\gamma (0) = e$ and suppose that $\gamma$ is also a geodesic. Then, the geodesic equation enforces:
 \begin{gather*}
  \nabla_X X = \nabla_{\dot \gamma}\dot \gamma = 0.
 \end{gather*}
 Thus, all immersed $1$-dimensional Lie subgroups of $G$ are geodesics of $\nabla$ if and only if $\alpha_\nabla (X,X) = \alpha^S_\nabla (X,X) = 0$ for every $X\in\mathfrak{g}$. By the polarization identity, this is equivalent to $\alpha^S_\nabla = 0$ showing (iii). Lastly, we verify (iv). The torsion tensor $T_\nabla$ of a left-invariant connection $\nabla$ is given by:
 \begin{gather*}
  T_\nabla (X,Y) = \nabla_X Y - \nabla_YX - [X,Y] = 2\alpha^A_\nabla (X,Y) - [X,Y].
 \end{gather*}
 $\nabla$ is torsion-free if and only if $T_\nabla$ vanishes finishing the proof.
\end{proof}

Proposition \autoref{prop:properties_of_left_con} infers that every Lie group $G$ admits a canonical connection:

\begin{corollary}[Cartan connection]\label{cor:cartan_connection}
 Let $G$ be a Lie group. Then, there exists a unique left-invariant connection $\nabla^G$ on $G$ which is torsion-free and whose symmetric form $\alpha^S_\nabla$ vanishes.
 We call $\nabla^G$ the \textbf{Cartan connection} of $G$ and its adapted complex structure $J_{\nabla^{G}}$ \textbf{Cartan structure}. If $X$ and $Y$ are left-invariant vector fields on $G$, $\nabla^G$ is given by:
 \begin{gather*}
  \nabla^G_X Y = \frac{1}{2}[X,Y].
 \end{gather*}
\end{corollary}

We are interested in the Cartan connection, because it agrees with the Levi-Civita connection $\nabla^g$ of a bi-invariant metric $g$:

\begin{proposition}\label{prop:cartan_levi-civita}
 Let $G$ be a Lie group with left-invariant metric $g$. Denote the Cartan connection on $G$ by $\nabla^G$ and the Levi-Civita connection of $g$ by $\nabla^g$. If $g$ is bi-invariant, then $\nabla^G = \nabla^g$. If $G$ is connected, this is an equivalence.
\end{proposition}

\begin{proof}
 Recall that one can compute the Levi-Civita connection $\nabla^g$ with the help of the Koszul formula:
 \begin{align*}
  2g (\nabla^g_XY,Z) = &\phantom{+}\ X(g(Y,Z)) + Y(g(X,Z)) - Z(g(X,Y))\\
  &+ g([X,Y],Z) - g([X,Z],Y) - g([Y,Z],X),
 \end{align*}
 where $X,Y,Z$ are vector fields on $G$. If $X,Y,Z$ are left-invariant, the first three terms in the Koszul formula vanish. From this, we see that
 \begin{gather*}
  \nabla^g_XY = \frac{1}{2}[X,Y] = \nabla^G_X Y
 \end{gather*}
 holds if and only if the last two terms in the Koszul formula cancel each other, i.e., if and only if:
 \begin{gather*}
  g(\ad_Z X, Y) = -g([X,Z],Y) = g([Y,Z],X) = -g(X,\ad_Z Y)\quad\forall X,Y,Z\in\mathfrak{g}.
 \end{gather*}
 If $G$ is connected, the last equation is equivalent to $g$ being bi-invariant, as we have seen in \autoref{sec:lie_groups}. This concludes the proof.
\end{proof}

Our next task is to compute the Cartan structure. To do that, we first consider the polar decomposition for Lie groups: Let $G$ be a complex Lie group with real form $G_\R$. Remember that the tangent bundle $TG_\R$ can be trivialized by left translation, i.e., via the map $G_\R\times\mathfrak{g}_\R\to TG_\R$, $(g,v)\mapsto dL_g (v)$, where $L_g (h)\coloneqq gh$. We now define the \textbf{polar decomposition} $P:TG_\R\to G$ of $G$ using this trivialization:
\begin{gather*}
 P(g,v) \coloneqq g\exp (iv),
\end{gather*}
where $\exp$ is the exponential map of $G$. $P$ allows us to formulate the following statement which is due to Sz\H{o}ke (cf. Proposition 7.3 in \cite{Szoeke2004}):

\begin{proposition}\label{prop:adapted_to_cartan}
 Let $G$ be a complex Lie group with complex structure $I$ and real form $G_\R$. Further, let $P:TG_\R\to G$ be the polar decomposition of $G$ and let $\gamma:\R\to G_\R$ be a geodesic of the Cartan connection $\nabla^{G_\R}$ with $g\coloneqq \gamma (0)$ and $v\coloneqq dL_{g^{-1}}\circ\dot\gamma (0)$. Define the map $f_\gamma:\C\to TG_\R$ by $f_\gamma (t+is)\coloneqq s\dot\gamma (t)$. Then, the following equation holds:
 \begin{gather}
  P\circ f_\gamma (z) = g\exp (zv)\quad\forall z\in\C,\label{eq:polar_geo}
 \end{gather}
 where $\exp$ is the exponential map of $G$. In particular, $P\circ f_\gamma$ is holomorphic and the complex structure $P^\ast I$, defined on an open and dense subset $U\subset TG_\R$, is adapted to $\nabla^{G_\R}$.
\end{proposition}

\begin{proof}
 By Proposition \autoref{prop:properties_of_left_con} (iii), the geodesics of the Cartan connection through $e$ are integral curves of left-invariant vector fields. Since $L_g:G_\R\to G_\R$ leaves $\nabla^{G_\R}$ invariant, the last statement extends to all geodesics, particularly $\gamma$:
 \begin{gather*}
  \gamma (t) = g\exp (tv)\quad\forall t\in\R.
 \end{gather*}
 We now verify \autoref{eq:polar_geo}:
 \begin{gather*}
  P\circ f_\gamma (t+is) = P(\gamma (t), sv) = g\exp (tv)\exp (isv) = g\exp ((t+is)v)\quad\forall t,s\in\R.
 \end{gather*}
 The holomorphicity of $P\circ f_\gamma$ directly follows from the holomorphicity of $\exp$. To prove the last statement, note the following, well-known fact (cf. Theorem 2.1 in \cite{Bremigan2000} and Proposition 7.4 in \cite{Szoeke2004}):
 \begin{gather*}
  U\coloneqq\{(g,v)\in TG_\R\mid dP_{(g,v)}\text{ invertible}\} = \{(g,v)\mid \Spec (\ad_v)\cap\pi\Z\subset\{0\}\}.
 \end{gather*}
 It is easy to see that $U$ is an open and dense subset of $TG_\R$ containing the zero section. On $U$, we can define the pullback complex structure $P^\ast I$. The holomorphicity of $P\circ f_\gamma$ implies that $P^\ast I$ is adapted to $\nabla^{G_\R}$ finishing the proof.
\end{proof}

The existence of the Cartan structure is not very surprising: It is a deep, but famous result from Lie group theory that every topological Lie group is automatically real-analytic\footnote{This is Hilbert's fifth problem which was solved by Gleason (cf. \cite{Gleason1952}), Montgomery, and Zippin (cf. \cite{Montgomery1952}).} which infers that the Cartan connection is also real-analytic. Proposition \autoref{prop:existence_adapted} now dictates that a complex structure adapted to the Cartan connection must exist. The more interesting aspect of Proposition \autoref{prop:adapted_to_cartan} is its description of the Cartan structure which allows us to express $J_{\nabla^{G_\R}}$ as $P^\ast I$.\\
The upshot of Proposition \autoref{prop:cartan_levi-civita} and \autoref{prop:adapted_to_cartan} is that Stenzel's theorem can be understood as a generalization of Ralph Bremigan's construction. In \cite{Bremigan2000}, Bremigan equipped the open and dense subset $U\subset TG_\R$ from Proposition \autoref{prop:adapted_to_cartan} with a semi-Kähler structure. For his construction, he used the complex structure $P^\ast I$ as well as the symplectic form $-\omega_{\can}$ after identifying $T G_\R$ with $T^\ast G_\R$ via the Killing form $K$. As we have seen in \autoref{sec:lie_groups}, the Killing form $K$ is a bi-invariant metric for semisimple Lie groups $G_\R$ which allows us to apply Theorem \autoref{thm:cotangent_kaehler} to $(G_\R, K)$. In light of Proposition \autoref{prop:cartan_levi-civita} and \autoref{prop:adapted_to_cartan}, we see that this is just Bremigan's construction. While Bremigan only constructed one semi-Kähler structure on $U\subset TG_\R$, Stenzel's approach generates a family of semi-Kähler structures by replacing\footnote{Strictly speaking, this also works for Bremigan's construction: Most statements in \cite{Bremigan2000} are still correct if one exchanges $K$ for a bi-invariant metric $g$.} $K$ with any bi-invariant metric $g$. One noteworthy point regarding this construction is that both $TG_\R$ and $T^\ast G_\R$ admit canonical structures: $P^\ast I$ on $TG_\R$ and $\omega_{\can}$ on $T^\ast G_\R$. Yet, the resulting semi-Kähler is not unique, because the identification of $TG_\R$ with $T^\ast G_\R$ depends on the choice of a bi-invariant metric $g$.\\
If the polar decomposition $P:TG_\R\to G$ is a diffeomorphism, then we can say that the complex group $G$ carries this family of semi-Kähler structures. This occurs, for example, if $G$ is the universal complexification of the compact Lie\linebreak group $G_\R$:

\begin{lemma}[Kähler structures on complex reductive groups $G$]\label{lem:kaehler_reductive_groups}
 Let $G$ be the universal complexification of a compact Lie group $G_\R$ with complex structure $I$. Then, the polar decomposition $P:TG_\R\to G$ is a diffeomorphism. In particular, $G$ admits a family of Kähler structures $(G,\omega_g,I)$. The symplectic form $\omega_g$ depends on the choice of a bi-invariant Riemannian metric $g$ on $G_\R$. Specifically, $\omega_g$ is the pullback of $-\omega_{\can}$ on $T^\ast G_\R$ where the identification\linebreak $G\cong TG_\R\cong T^\ast G_\R$ is induced by $P$ and $g$.
\end{lemma}

\begin{remark}
 The existence of Kähler structures on complex reductive groups is a well-established fact and was already known to Kronheimer (cf. the introduction of \cite{Kronheimer2004}) and Stenzel (cf. Section 4.1 in \cite{Stenzel1990}).
\end{remark}

\begin{proof}
 We begin by observing that the domain $U\subset TG_\R$ from Proposition \autoref{prop:adapted_to_cartan} is equal to $TG_\R$, because all eigenvalues of $\ad_v$, $v\in\mathfrak{g}_\R$, are purely imaginary (cf. \autoref{sec:semi-kaehler}). Thus, it suffices to show that $P:TG_\R\to G$ is bijective. There are several ways to prove this. The easiest method is to first consider the case $G_\R = \U (n)$. The universal complexification of $\U (n)$ is $\GL (n,\C)$. It is a standard result from linear algebra that the polar decomposition $P:TG_\R\to G$ is bijective for $G_\R = \U (n)$ and $G = \GL (n,\C)$. To infer the general case, we use a powerful fact from representation theory, namely that every compact Lie group admits a faithful, finite-dimensional, and unitary representation (cf. Theorem 6.1.2 in \cite{Kowalski2014}). Simply put, we can take $G_\R$ to be a closed subgroup of $\U (n)$. The bijectivity of $P$ for matrices now implies the bijectivity of $P$ for $G_\R$. The remaining statements of Lemma \autoref{lem:kaehler_reductive_groups} follow from Theorem \autoref{thm:cotangent_kaehler} as well as Proposition \autoref{prop:cartan_levi-civita} and \autoref{prop:adapted_to_cartan}.
\end{proof}

Lastly, we want to discuss what happens if we apply Theorem \autoref{thm:kaehler_duality} to compact groups $(G_\R,g)$. Since $T^\ast G_\R$ is isomorphic to the universal complexification $G$, Theorem \autoref{thm:kaehler_duality} will give us a family of Hyperkähler and holomorphic Kähler structures on the cotangent bundle $T^\ast G$. Now recall that the idea behind Theorem \autoref{thm:kaehler_duality} is to apply Stenzel's theorem to two natural metrics on $T^\ast G_\R\cong G$: One is obtained by complexifying $g$. We have already explained in \autoref{sec:holo_semi-kaehler} that the complexification $g_\C$ of $g$ is a bi-invariant metric on $G$ with signature $(n,n)$, where $n\coloneqq\dim_\R G_\R$ (cf. Proposition \autoref{prop:comp_of_metric}). We know that applying Stenzel's theorem to Lie groups with bi-invariant metrics agrees with Bremigan's construction (when $K$ is replaced by $g$), so we can take the holomorphic Kähler structures on $T^\ast G$ to be a result of Bremigan's construction\footnote{This statement has to be taken with a grain of salt: In the scenario at hand, Stenzel's theorem for holomorphic metrics (cf. Lemma \autoref{lem:stenzel_for_holo_metrics}) only coincides with Bremigan's construction if the map $\phi$ is chosen to be the identity which might not always be possible.}. The other natural metric on $G$ is the Kähler metric of $(G,\omega_g,I)$. Usually, this Kähler metric is not bi-invariant, because every group which admits a bi-invariant Riemannian metric must be a product of a compact and an Abelian group which, in general, $G$ is not. In any case, the resulting Hyperkähler structures on $T^\ast G$ should coincide with the ones mentioned in \autoref{app:hyp_orb}, even though we lack a proof for this statement. Let us summarize our findings:

\begin{theorem}[Kähler duality for Lie groups]\label{thm:kaehler_duality_lie}
 Let $G$ be the universal complexification of a compact Lie group $G_\R$. If Conjecture \autoref{con:commutation_relation} holds, then there exists a neighborhood $U\subset T^\ast G$ of $G\subset T^\ast G$ admitting a family of Hyperkähler and holomorphic Kähler structures depending on the choice of a bi-invariant Riemannian metric $g$ on $G_\R$. These structures are constructed by applying Theorem \autoref{thm:kaehler_duality} to $(G_\R,g)$. Specifically, the Hyperkähler structures on $U$ are obtained by applying Stenzel's theorem (cf. Lemma \autoref{lem:stenzel_for_kaehler}) to the Kähler structures from Lemma \autoref{lem:kaehler_reductive_groups}, while the holomorphic Kähler structures on $U$ are the result of Bremigan's construction (cf. \cite{Bremigan2000}).
\end{theorem}

To conclude \autoref{sec:duality}, we tackle the initial question concerning the origin of Kähler duality on coadjoint orbits. The following conjecture outlines how the Kähler duality of $T^\ast G$ and, thereby, of double cotangent bundles might possibly be related to the Kähler duality of coadjoint orbits:

\begin{conjecture}[Reduction of Kähler duality]\label{con:reduction}
 Let $G$ be the universal complexification of a compact Lie group $G_\R$. Then, the Kähler duality on $U\subset T^\ast G$ reduces (in the sense of \autoref{app:reduction}) to the Kähler duality on coadjoint orbits $\mathcal{O}^\ast$ of $G$. Precisely speaking, the Hyperkähler structures on $T^\ast G$ are compatible with the Hyperkähler structures of $\mathcal{O}^\ast$ described in \autoref{app:hyp_orb}, while the holomorphic Kähler structures on $T^\ast G$ are compatible with the holomorphic Kähler structures on $\mathcal{O}^\ast$ described in \autoref{sec:holo_semi-kaehler}.
\end{conjecture}

\chapter{Outlook}
\label{chap:outlook}
Future research on the topics presented in this thesis can head in several directions. One possibility is trying to solve the problems which occurred in the course of our investigations, but remained without a solution. By our count, there are five of those problems: The first one concerns holomorphic Morse-Darboux charts. During the study of Lefschetz and almost toric fibrations in the middle of \autoref{sec:HHS}, we have proven that two-dimensional real-analytic RHSs $(M^2,\omega,H)$ always admit real-analytic Morse-Darboux charts near non-degenerate critical points of $H$ with Morse index $\neq 1$ (cf. Lemma \autoref{lem:morse_darboux_lem_I} and \autoref{lem:morse_darboux_lem_II}). To extend this statement to the holomorphic setting, it suffices by Theorem \autoref{thm:morse_darboux_ex} and Corollary \autoref{cor:morse_darboux_holo_ex} to show that every complex two-dimensional HHS $(X,\Omega,\mH)$ locally admits a real form near non-degenerate critical points of $\mH$ such that $(X,\Omega,\mH)$ restricts to a RHS $(M,\omega,H)$ as above on that real form. Unfortunately, we were unable to find such a real form. Solving this problem would have great consequences, as it would mean that every proper holomorphic Lefschetz fibration $\pi:X\to C$ equipped with a holomorphic symplectic form $\Omega$ is automatically an almost toric fibration in dimension two (cf. Definition \autoref{def:holo_symp_Lef_fi} and Proposition \autoref{prop:holo_lef_toric}).\\
The second problem regards Theorem \autoref{thm:generic} which states that proper PHHSs are generic, as long as their real dimension is larger than four. We do not know whether this result still applies in dimension four, so it would be a good starting point for further work. The third and fourth problem are both listed as conjectures in \autoref{sec:duality}. Conjecture \autoref{con:commutation_relation} is about the commutation relations of the complex structures $\phi^\ast J_g$ and $T^\ast I$ needed for the modified versions of Stenzel's theorem, while Conjecture \autoref{con:reduction} asserts that the Kähler dualities of double cotangent bundles and coadjoint orbits are linked via reduction. A full proof of the commutation relations as well as a complete reduction analysis are still missing. The final open problem can also be found in \autoref{sec:duality} and is more vague than the others. It poses the question how exactly our Hyperkähler structures on double cotangent bundles are related to those already available in the literature (Bielawski \cite{Bielawski2003}, Feix \cite{Feix2001}, Kaledin \cite{Kaledin1997}, and Kronheimer \cite{Kronheimer2004}).\\
Another direction for future work concerns the generalization of HHSs with respect to the Hamiltonian $\mH$ instead of the complex structure\footnote{This path led us to the notion of a PHHS.} $J$. We have always assumed in this thesis that the Hamiltonian of a HHS $(X,\Omega,\mH)$ is a holomorphic map from $X$ to the complex plane $\C$. However, we can also view holomorphic maps from $X$ to a Riemann surface $C$ as Hamiltonians. In this case, the triple $(X,\Omega,\mH)$ becomes locally a HHS in the sense of Definition \autoref{def:holo_ham_sys} if we choose a holomorphic chart for the Riemann surface $C$. Given two charts of $C$, the associated holomorphic Hamiltonian vector fields $X_{\mH}$ only differ by multiplication with a holomorphic function. Thus, the distributions spanned by these vector fields agree (cf. \autoref{sec:HHS}, in particular the proof of Proposition \autoref{prop:holo_foli}), even though the Hamiltonian vector fields might be different, giving us a global distribution on $X$. Since this distribution is integrable, it induces a foliation, called holomorphic Hamiltonian foliation, whose leaves can be interpreted as maximal holomorphic trajectories, as we have seen in \autoref{sec:HHS} (cf. Definition \autoref{def:holo_foli} and the explanation afterwards). In this sense, the holomorphic Hamiltonian foliation can be seen as a generalization\footnote{Of course, pseudo-holomorphic Hamiltonian foliations also exist for PHHSs, giving us a generalization of HHSs with respect to $J$ and $\mH$ at the same time.} of HHSs. The advantage of holomorphic Hamiltonian foliations is that they give us easier access to periodic orbits. Consider, for instance, the complex torus $X\coloneqq\C^{2n}/\Z^{2n}+i\Z^{2n}$ parameterized by $[z_1,\ldots,z_{2n}]\in X$ and equipped with the holomorphic symplectic form $\Omega \coloneqq \sum^{n}_{j=1} dz_{j}\wedge dz_{j+n}$. Since $X$ is compact, all holomorphic maps from $X$ to $\C$ are constant, so all ``classical'' HHSs are trivial. However, the holomorphic Hamiltonian foliation associated with $\mH_j:X\to\C/\Z+i\Z$, $[z_1,\ldots,z_{2n}]\mapsto [z_j]$ ($j\in \{1,\ldots, n\}$) is non-trivial and its leaves are tori, since $X_\mH = \pa{z_{j+n}}$ is the corresponding Hamiltonian vector field (in the usual torus coordinates). There are several interesting and open questions regarding holomorphic Hamiltonian foliations, e.g. under which conditions they admit periodic orbits, i.e. under which conditions these foliations possess compact leaves.\\
To answer this question or similar questions for HHSs and PHHSs, one may follow a Floer-theoretical approach. In \autoref{sec:HHS} and \autoref{app:various_action_functionals}, we have constructed various actions functionals\footnote{In certain cases, one can also write down action functionals for holomorphic Hamiltonian foliations. Take, for instance, the torus example from above. Here, we can use more or less the same formulas from \autoref{sec:HHS} and \autoref{app:various_action_functionals} to write down action functionals which are now, of course, torus-valued.} for the (pseudo-)holomorphic trajectories of (P)HHSs. One might wonder whether it is possible to use these action functionals as a basis for a Floer-like theory. As of now, it is hard to imagine how such a theory should look like and, if it exists, whether it gives meaningful results. There are mainly two complications: First of all, most action functionals we have constructed are complex-valued which is a problem, since we do not know how to cook up a Morse theory for $\C$-valued Morse functions. The second issue concerns the remaining, real action functionals developed in \autoref{app:various_action_functionals}. They all have counterintuitive boundary conditions and, most likely, will just yield trivial results (cf. \autoref{app:various_action_functionals}).\\
Last but not least, one can study examples of HHSs and PHHSs to gain more insight into their dynamics. This path seems to be quite promising, especially because HHSs and PHHSs exhibit properties like monodromy and non-trivial structures of the domains of their trajectories which are absent in the case of RHSs.

\begin{appendix}

\counterwithin{definition}{chapter}


 \chapter{Real Structures on Complex Manifolds}
 \label{app:real_structures}
 This part is concerned with real structures on complex manifolds. It serves as a reminder to those who are already familiar with multidimensional complex analysis and as a brief introduction to those who are new to the topic. Aside from explaining the notion of a real structure and deriving some of its basic properties, we also discuss how real structures interact with tensors.\\
In a nutshell, a real structure $\sigma$ on a complex manifold $X$ is an antiholomorphic involution. As we will see soon, real structures are intimately linked to the notion of real forms which are the fixed point sets of real structures. In fact, there are three equivalent definitions of real forms:

\begin{definition}[Real form]\label{def:real_form}
 Let $X$ be a complex manifold with complex structure $J$. Further, let $M\subset X$ be a non-empty subset. $M$ is called \textbf{real form} of $X$ if one of the following, equivalent conditions is satisfied:
 \begin{enumerate}
  \item $M$ is the fixed point set of a \textbf{real structure}, i.e., there exists an open neighborhood $U\subset X$ of $M$ and an antiholomorphic involution $\sigma$ on $U$ such that $M = \Fix\sigma$.
  \item $M\subset X$ is a \textbf{totally real}, real-analytic submanifold, i.e., $M$ satisfies\linebreak $\dim_\R M = \dim_\C X$ and $J_p(T_pM)\cap T_pM = \{0\}$ for every point $p\in M$.
  \item For every point $p\in M$ there exists an open neighborhood $V\subset X$ of $p$ and a holomorphic chart $(z_1,\ldots, z_n):V\to\C^n$ such that ($z_j = x_j + iy_j$):
  \begin{gather*}
   V\cap M = \{q\in V\mid y_1 (q) = \ldots = y_n (q) = 0\}.
  \end{gather*}
 \end{enumerate}
\end{definition}

\begin{proof}
 (i)$\Rightarrow$(ii): Let $\sigma$ be as above and take a point $p\in M = \Fix\sigma$. Without loss of generality, we can assume $U = X$. Since $\sigma$ is an involution, its differential $d\sigma_p:T_pX\to T_pX$ is an involution as well. Thus, the tangent space $T_pX$ splits into the direct sum $E_1\oplus E_{-1}$ where $E_{\pm 1}$ is the eigenspace of $d\sigma_p$ with respect to the eigenvalue $\pm 1$. Now consider a real-analytic Riemannian metric $g$ on a $\sigma$-invariant neighborhood $V\subset X$ of $p$, i.e., $\sigma (V) = V$. $\sigma$ is an isometry with respect to the real-analytic metric $g_\sigma\coloneqq g + \sigma^\ast g$ on $V$. Hence, the exponential map $\exp_p:W\subset T_pX\to V$ of $g_\sigma$ at $p$ restricts to a map $\exp_p:W\cap E_1\to \Fix\sigma$ and gives rise to a real-analytic submanifold chart for $\Fix\sigma = M$. In particular, $M\subset X$ is a real-analytic submanifold and we find $T_pM = E_1$.\\
 Now recall that $\sigma$ is antiholomorphic, i.e., $d\sigma_p\circ J_p = -J_p\circ d\sigma_p$. Thus, we have $E_{-1} = J_p(E_1)$. This yields:
 \begin{gather*}
  J_p (T_pM)\cap T_pM = J_p (E_1)\cap E_1 = E_{-1}\cap E_1 = \{0\}\quad\text{and}\\
  \dim_\R M = \dim_\R T_pM = \dim_\R E_1 = \dim_\R E_{-1} = \frac{1}{2}\dim_\R T_pX = \dim_\C X
 \end{gather*}
 proving that $M\subset X$ is a totally real, real-analytic submanifold.\\\\
 (ii)$\Rightarrow$(iii): Let $M\subset X$ be a totally real, real-analytic submanifold. Pick a point $p\in M$, an open neighborhood $V\subset X$ of $p$, a holomorphic chart $\psi:V\to\C^n$ of $X$, and a real-analytic chart $\phi:V\cap M\to\R^n$ of $M$. Then, $\psi\circ\phi^{-1}$ is a real-analytic map from $\R^n$ to $\C^n\cong \R^{2n}$. Thus, there exists a unique holomorphic map $\Psi$ from $\C^n$ to $\C^n$ coinciding with $\psi\circ\phi^{-1}$ on $\R^n$. To illustrate this, we consider the case $n=1$:\\ Let $f = a+ib:\R\to\C\cong\R^2$ be a real-analytic map, i.e., the maps $a,b:\R\to\R$ are real-analytic. This means that the Taylor series of $a$ (and $b$) locally converges and coincides with the function $a$ (and $b$, respectively). Expressing $a$ and $b$ as
 \begin{gather*}
  a(x) = \sum^\infty_{k=0}a_k(x-x_0)^k\quad\text{and}\quad b(x) = \sum^\infty_{k=0}b_k(x-x_0)^k
 \end{gather*}
 allows us to locally define the unique holomorphic continuation $F$ of $f$ by setting:
 \begin{gather*}
  F(z)\coloneqq \sum^\infty_{k=0}(a_k + ib_k)(z-x_0)^k.
 \end{gather*}
 Since $F$ can locally be written as a power series, it is holomorphic. Furthermore, it is unique, as any other holomorphic continuation of $f$ must have the same power series near points $x_0\in\R$.\\
 Let us return to the map $\Psi$. We want to show that $\Psi$ is a local biholomorphism near $\phi (p)$. To prove this, it suffices to show that the differential\linebreak $d\Psi_{\phi (p)}:T_{\phi (p)}\C^n\cong \C^n\to T_{\psi (p)}\C^n\cong \C^n$ is bijective. Let $\pa{x_1,p},\ldots, \pa{x_n,p}$ be the basis vectors of $T_p M\subset T_pX$ associated to the chart $\phi$. We find:
 \begin{align*}
  d\Psi_{\phi (p)} (e_j) &= d\psi_p (\pa{x_j,p})\quad\text{and}\\
  d\Psi_{\phi (p)} (ie_j) &= id\Psi_{\phi (p)} (e_j) = id\psi_p (\pa{x_j,p}) = d\psi_p (J_p\pa{x_j,p}),
 \end{align*}
 where $e_1,\ldots, e_n$ denotes the standard basis of $\C^n$ and we exploited the holomorphicity of $\psi$ and $\Psi$. Since $M\subset X$ is a totally real submanifold, the vectors $\pa{x_1,p},\ldots, \pa{x_n,p}, J_p\pa{x_1,p},\ldots, J_p\pa{x_n,p}$ form a basis of $T_pX$. Moreover, $d\psi_p$ maps bases to bases, as $\psi$ is a chart. Together with the previous equations, this tells us that $d\Psi_{\phi (p)}$ maps a real basis of $\C^n$ to another basis, i.e., that $d\Psi_{\phi (p)}$ is bijective.\\
 We can now define the holomorphic chart $(z_1,\ldots,z_n)\coloneqq\Psi^{-1}\circ\psi:V\to\C^n$ after shrinking $V$ if necessary. By construction, $\Psi^{-1}\circ\psi$ coincides with $\phi$ on $V\cap M$ and only assumes real values on $V\cap M$ proving ($z_j = x_j + iy_j$):
 \begin{gather*}
  V\cap M = \{q\in V\mid y_1 (q) = \ldots = y_n (q) = 0\}.
 \end{gather*}
 (iii)$\Rightarrow$(i): Pick a collection $\{(V_\alpha,\psi_\alpha)\}_\alpha$ of holomorphic charts as in (iii) covering $M$. Without loss of generality, we can assume that $c(\psi_\alpha (V_\alpha)) = \psi_\alpha (V_\alpha)$, where $c:\C^n\to \C^n$ is the complex conjugation. Set $U\coloneqq \cup_\alpha V_\alpha$ and define the map $\sigma:U\to U$ as follows:
 \begin{gather*}
  \sigma (p)\coloneqq \psi^{-1}_\alpha\circ c\circ \psi_\alpha (p)\quad\forall p\in V_\alpha\forall \alpha.
 \end{gather*}
 It is clear from the definition that $\sigma$ is an antiholomorphic involution if it is well-defined. To show well-definedness, it suffices to prove that the maps $\psi_\alpha\circ\psi^{-1}_\beta$ commute with the complex conjugation $c$ for every $\alpha$ and $\beta$. For this, we first note that, by (iii), the points in $V_\alpha\cap M$ are exactly those points which are mapped to $\R^n$ under $\psi_\alpha$. Hence, the map $\psi_\alpha\circ\psi^{-1}_\beta$ sends $\R^n$ to $\R^n$. Now consider the maps $\psi_\alpha\circ\psi^{-1}_\beta$ and $c\circ\psi_\alpha\circ\psi^{-1}_\beta\circ c$. Both are holomorphic and coincide on $\R^n$. By the identity theorem\footnote{To apply the identity theorem, it is necessary that every connected component of $V_\alpha\cap V_\beta$ is met by some point in $M$. We can always achieve this by choosing a suitable collection $\{(V_\alpha,\psi_\alpha)\}_\alpha$.} for holomorphic maps, they must be equal giving us:
 \begin{gather*}
  c\circ\psi_\alpha\circ\psi^{-1}_\beta\circ c = \psi_\alpha\circ\psi^{-1}_\beta\quad\forall \alpha,\beta.
 \end{gather*}
 Since $c$ is an involution, the last equation gives us the desired commutation relation between $\psi_\alpha\circ\psi^{-1}_\beta$ and $c$.\\
 Lastly, we need to check that $M$ is the fixed point set of $\sigma$. By definition of $\sigma$, the fixed points of $\sigma$ are exactly those points in $U$ which are mapped to $\R^n$ by some $\psi_\alpha$. By (iii), these are exactly the points in $M$.
\end{proof}

Given a real form $M$, one can ask whether the corresponding real structure is unique, i.e., whether there are two different real structures whose real forms are $M$. The following proposition answers this question:

\begin{proposition}[Uniqueness of real structures]\label{prop:real_structure_unique}
 Let $X$ be a complex manifold and $\sigma_1,\sigma_2:X\to X$ be two real structures with the same non-empty fixed point set $M$. Then, there exists an open neighborhood $U\subset X$ of $M$ such that $\sigma_1\vert_U = \sigma_2\vert_U$. If every connected component of $X$ meets $M$ in some point, one has $\sigma_1 = \sigma_2$.
\end{proposition}

\begin{proof}
 Let $\psi$ be a holomorphic chart of $X$ near a point $p\in M$ as in (iii) of Definition \autoref{def:real_form} and consider the map $\sigma_1\circ\sigma_2$. In the chart $\psi$, $\sigma_1\circ\sigma_2$ becomes a holomorphic map from $\C^n$ to $\C^n$ which restricts to the identity on $\R^n$. By the identity theorem, this map must be the identity on all of $\C^n$. Hence, $\sigma_1\circ\sigma_2$ coincides with the identity on a neighborhood of the point $p$. Repeating this argument for every point $p\in M$ shows that there exists an open neighborhood $U\subset X$ of $M$ such that $\sigma_1\circ\sigma_2\vert_U = \id_U$ or, equivalently, $\sigma_1\vert_U = \sigma_2\vert_U$.\\
 If every connected component of $X$ meets $M$ in some point, then every connected component of $X$ contains some open subset on which $\sigma_1$ and $\sigma_2$ coincide. Applying the identity theorem again shows $\sigma_1 = \sigma_2$.
\end{proof}

\begin{remark}[Nice real structure]\label{rem:nice_real_structures}
 In light of Proposition \autoref{prop:real_structure_unique}, we call a real structure $\sigma$ on $X$ \textbf{nice} if its fixed point set meets every connected component of $X$ in some point. The notion of a nice real form is defined analogously.
\end{remark}

At this point, it should be noted that real-analyticity is necessary in (ii) of Definition \autoref{def:real_form}. As a counterexample, take $X$ to be the complex plane $\C$ with complex structure $i$ and $M$ to be the graph of a smooth function $f:\R\to\R$ with compact support satisfying $f(0)\neq 0$. Clearly, $M\subset\C$ is a totally real, smooth submanifold of $\C$. Now assume that $M$ is the fixed point set of some real structure $\sigma$ on an open neighborhood of $M$. Outside a compact set, $M$ coincides with the real line $\R$. Now note that the complex conjugation on $\C$ is a real structure whose real form is $\R$. This allows us to apply the identity theorem as in the proof of Proposition \autoref{prop:real_structure_unique} to show that $\sigma$ coincides with the complex conjugation. Hence, $M$ must be the real line contradicting $f(0)\neq 0$.

\begin{remark}[Real structures on almost complex manifolds]\label{rem:almost_real_structure}
 If $X$ is just an almost complex manifold, we can still use (i) or (ii)\footnote{In this case, it is sensible to drop the condition that $M$ has to be real-analytic, since $X$ itself might not be real-analytic.} from Definition \autoref{def:real_form} to define real forms $M\subset X$. In this case, it is still true that every real form in the sense of (i) is also a real form in the sense of (ii). However, the converse might fail, as the previous example suggests.
\end{remark}

The remainder of \autoref{app:real_structures} is devoted to the relation between tensors and real structures. Before we dive into this discussion, it is wise to introduce the concepts of $\sigma$-charts:

\begin{definition}[$\sigma$-Chart]
 Let $X$ be a complex manifold with real structure $\sigma$ on it. Further, let $\psi:V\to\C^n$ be a holomorphic chart of $X$. We say that $\psi$ is a $\mathbf{\sigma}$\textbf{-chart} if the following two conditions are satisfied:
 \begin{enumerate}
  \item $\sigma (V) = V$.
  \item $\psi\circ\sigma = c\circ\psi$, where $c:\C^n\to\C^n$, $c(z)\coloneqq \bar{z}$ denotes complex conjugation.
 \end{enumerate}
\end{definition}

We can interpret $\sigma$-charts as charts in which $\sigma$ assumes its standard form.\linebreak $\sigma$-Charts are enormously helpful, as they simplify a vast amount of upcoming computations. The following proposition ensures that we can always find an atlas of $\sigma$-charts (cf. \cite{Kulkarni1978}):

\begin{proposition}[$\sigma$-Atlas]\label{prop:sigma-atlas}
 Let $X$ be a complex manifold with real structure $\sigma$ on it. Then, $X$ admits an atlas of $\sigma$-charts.
\end{proposition}

\begin{proof}
 Set $M\coloneqq\Fix\sigma$ and pick a point $p\in X\backslash M$. Let $\psi:V\to\C^n$ be a holomorphic chart of $X$ near $p$. By shrinking $V$ if necessary, we can achieve that $\sigma (V)\cap V = \emptyset$. Now define the holomorphic chart $\psi_\sigma:V\cup\sigma (V)\to\C^n$ as follows:
 \begin{gather*}
  \psi_\sigma (q)\coloneqq\begin{cases}\psi (q) & \text{ if }q\in V,\\c\circ\psi\circ\sigma (q) & \text{ if }q\in\sigma (V).\end{cases}
 \end{gather*}
 By construction, $\psi_\sigma$ is a $\sigma$-chart near $p$.\\
 If $M$ is empty, the proof is already finished. If $M$ is not empty, it remains to be shown that we can find $\sigma$-charts near any point $p\in M$. However, in the proof of Definition \autoref{def:real_form}, we have shown that the charts from (iii) in Definition \autoref{def:real_form} are, by construction of $\sigma$, $\sigma$-charts. This concludes the proof.
\end{proof}

Let us now begin the discussion of the relation between tensors and real structures. We start by recalling that any manifold $X$ with an almost complex structure $J$ on it possesses a canonical splitting of its complexified (co)tangent spaces into $(1,0)$- and $(0,1)$-spaces:
\begin{gather*}
 T_{p,\C}X = T^{(1,0)}_pX \oplus T^{(0,1)}_pX\quad \left(T^\ast_{p,\C}X = T^{\ast, (1,0)}_pX \oplus T^{\ast, (0,1)}_pX\right).
\end{gather*}
The elements $v\in T^{(1,0)}_pX$ ($\alpha\in T^{\ast, (1,0)}_pX$) are defined by the equation:
\begin{gather*}
 J_p v = iv\quad \left(\alpha\circ J_p = i\alpha\right),
\end{gather*}
while the elements $w\in T^{(0,1)}_pX$ ($\beta\in T^{\ast, (0,1)}_pX$) are defined by:
\begin{gather*}
 J_p w = -iw\quad \left(\beta\circ J_p = -i\beta\right).
\end{gather*}
If $J$ is integrable, i.e., $X$ is a complex manifold, we can pick a holomorphic chart $(z_1 = x_1 + iy_1,\ldots,z_n = x_n + iy_n):V\to\C^n$ of $X$. The vectors (forms)
\begin{gather*}
 \pa{z_j}\coloneqq\frac{1}{2}\left(\pa{x_j} - i\pa{y_j}\right)\quad \left(dz_j\coloneqq dx_j + idy_j\right)
\end{gather*}
form a basis of the $(1,0)$-spaces, while the vectors (forms)
\begin{gather*}
 \pa{\bar{z}_j}\coloneqq\frac{1}{2}\left(\pa{x_j} + i\pa{y_j}\right)\quad \left(d\bar{z}_j\coloneqq dx_j - idy_j\right)
\end{gather*}
form a basis of the $(0,1)$-spaces. If, additionally, $X$ carries a real structure $\sigma$ and the chosen chart is a $\sigma$-chart, one easily checks the following relation for the corresponding coordinate vectors (forms):
\begin{gather*}
 \sigma_\ast\pa{z_j} = \pa{\bar{z}_j}\quad\left(\sigma^\ast dz_j = d\bar{z}_j\right).
\end{gather*}
Equipped with this knowledge, we are now able to define holomorphic tensors:

\begin{definition}[Holomorphic tensors]\label{def:holo_tensor}
 Let $X$ be a complex manifold and $T$ be a smooth complex tensor field on $X$. We say that $T$ is of \textbf{holomorphic type} if $T$ takes the following form in every holomorphic chart $(z_1,\ldots, z_n):V\to\C^n$ of $X$:
 \begin{gather*}
  T = T^{j_1\ldots j_r}_{k_1\ldots k_s} \pa{z_{j_1}}\otimes\ldots\otimes\pa{z_{j_r}}\otimes dz_{k_1}\otimes\ldots\otimes dz_{k_s},
 \end{gather*}
 where we employ the Einstein sum convention and the coefficients\linebreak $T^{j_1\ldots j_r}_{k_1\ldots k_s}:V\to\C$ are smooth functions on $V$. The tensor field $T$ is \textbf{holomorphic} if it is of holomorphic type and the coefficients $T^{j_1\ldots j_r}_{k_1\ldots k_s}$ are holomorphic functions. Now assume that $X$ additionally carries a real structure $\sigma$. We call the tensor field $T$ $\mathbf{\sigma}$\textbf{-invariant} if it satisfies\footnote{For a diffeomorphism $\phi:L\to N$ and a tensor field $T$ on $N$, $\phi^\ast T$ is defined by pulling back the form part of $T$ with $\phi$ and pushing forward the vector part of $T$ with $\phi^{-1}$. If $T$ is a tensor field on $L$, $\phi_\ast T$ is defined similarly, but the roles of $\phi$ and $\phi^{-1}$ are reversed. Since the diffeomorphism $\phi = \sigma$ in Definition \autoref{def:holo_tensor} is an involution, these subtle differences are not important.} $\overline{\sigma^\ast T} = T$.
\end{definition}

\begin{remark}\label{rem:after_holo_tensor}
 In the case that the almost complex structure of $X$ is not integrable, we can still define tensors of holomorphic type\footnote{Naturally, we need to employ a different, but equivalent definition for tensors of holomorphic type in this case, namely, that $T$ vanishes if we plug any $(0,1)$-vectors or -forms into it.} and $\sigma$-invariant tensors. The integrability is only needed for the definition of holomorphic tensors.
\end{remark}

The key observation about $\sigma$-invariant tensors of holomorphic type is that they restrict to real tensors on the real form $M$:

\begin{proposition}\label{prop:invariant_holo_tensors}
 Let $X$ be a complex manifold with real structure $\sigma$ and non-empty real form $M\coloneqq\Fix\sigma$. Further, let $T$ be a $\sigma$-invariant tensor field on $X$ of holomorphic type. Then, $T$ induces a real tensor field $T_M$ on $M$. If $T$ is holomorphic, $T_M$ is real-analytic.
\end{proposition}

\begin{proof}
 Pick a point $p\in M$ and a $\sigma$-chart $(z_1,\ldots, z_n):V\to\C^n$ near $p$. In this $\sigma$-chart, we can write $T$ as:
 \begin{gather*}
  T = T^{j_1\ldots j_r}_{k_1\ldots k_s} \pa{z_{j_1}}\otimes\ldots\otimes\pa{z_{j_r}}\otimes dz_{k_1}\otimes\ldots\otimes dz_{k_s},
 \end{gather*}
 where we employ the same conventions as in Definition \autoref{def:holo_tensor}. Using this expression, the equation $\overline{\sigma^\ast T} = T$ becomes:
 \begin{gather*}
  \overline{T^{j_1\ldots j_r}_{k_1\ldots k_s}\circ\sigma} = T^{j_1\ldots j_r}_{k_1\ldots k_s}.
 \end{gather*}
 Thus, the coefficients $T^{j_1\ldots j_r}_{k_1\ldots k_s}(q)$ are real for points $q\in V\cap M$. We are now able to define the real tensor $T_M$ on $V\cap M$:
 \begin{gather}\label{eq:T_M}
  T_M(q)\coloneqq T^{j_1\ldots j_r}_{k_1\ldots k_s} (q) \pa{x_{j_1},q}\otimes\ldots\otimes\pa{x_{j_r},q}\otimes dx_{k_1,q}\otimes\ldots\otimes dx_{k_s,q},
 \end{gather}
 where $q\in V\cap M$ and $z_j = x_j + iy_j$. Recall that $\sigma$-charts near points of $M$ are the holomorphic charts from (iii) in Definition \autoref{def:real_form}, i.e., we have:
 \begin{gather*}
  V\cap M = \{q\in V\mid y_1 (q) = \ldots = y_n (q) = 0\}.
 \end{gather*}
 In particular, this means that the coordinate functions $x_1,\ldots, x_n$ constitute a real-analytic submanifold chart for $M$. Hence, the object defined by \autoref{eq:T_M} is a well-defined real tensor field on $V\cap M$.\\
 To obtain a real tensor field on all of $M$, we repeat the previous step for every point $p\in M$. This is possible, since the definition of $T_M$ given in \autoref{eq:T_M} is independent of the choice of $\sigma$-chart, as one can easily verify. To conclude the proof, we observe that the coefficients $T^{j_1\ldots j_r}_{k_1\ldots k_s}$ are real-analytic if $T$ is holomorphic.
\end{proof}

Given a real tensor $T_M$, one can ask whether the corresponding $\sigma$-invariant tensor $T$ of holomorphic type is unique, i.e., whether there are two different $\sigma$-invariant tensors of holomorphic type that induce the same tensor $T_M$. In general, the tensor $T$ is not unique, however, it becomes unique if we additionally require $T$ to be holomorphic:

\begin{lemma}[Uniqueness of holomorphic tensors]\label{lem:uniqueness_of_holomorphic_tensors}
 Let $X$ be a complex manifold with real structure $\sigma$ and non-empty real form $M\coloneqq\Fix\sigma$. Further, let $S$ and $T$ be two $\sigma$-invariant holomorphic tensor fields on $X$ which induce the same real tensor field $T_M$ on $M$. Then, there exists an open neighborhood $U\subset X$ of $M$ such that $S\vert_U = T\vert_U$. If $M$ is nice, one even has $S = T$.
\end{lemma}

\begin{proof}
 Pick a point $p\in M$ and a $\sigma$-chart $(z_1,\ldots, z_n) = \psi:V\to\C^n$ near $p$ such that $V$ is connected. In this $\sigma$-chart, we can write $S$ and $T$ as:
 \begin{align*}
  S &= S^{j_1\ldots j_r}_{k_1\ldots k_s} \pa{z_{j_1}}\otimes\ldots\otimes\pa{z_{j_r}}\otimes dz_{k_1}\otimes\ldots\otimes dz_{k_s},\\
  T &= T^{j_1\ldots j_r}_{k_1\ldots k_s} \pa{z_{j_1}}\otimes\ldots\otimes\pa{z_{j_r}}\otimes dz_{k_1}\otimes\ldots\otimes dz_{k_s},
 \end{align*}
 where we employ the same conventions as in Definition \autoref{def:holo_tensor}. Since $S$ and $T$ induce the same real tensor $T_M$, their coefficients $S^{j_1\ldots j_r}_{k_1\ldots k_s}$ and $T^{j_1\ldots j_r}_{k_1\ldots k_s}$ coincide on $V\cap M$. By precomposing them with $\psi^{-1}$, we can view the coefficients as holomorphic functions from $\C^n$ to $\C$ which coincide on $\R^n$. Hence, by the identity theorem, they must coincide everywhere, i.e., on all of $V$. This shows that $S$ and $T$ agree on $V$. We can now repeat this argument for every point $p\in M$ to find a neighborhood $U\subset X$ of $M$ such that $S\vert_U = T\vert_U$. If $M$ is nice, we can apply the identity theorem to every connected component of $X$ to prove $S = T$.
\end{proof}

Lastly, we want to discuss the question of existence, i.e., given a real tensor $T_M$, is there always a $\sigma$-invariant holomorphic tensor $T$ inducing $T_M$? The answer is a resounding ``Yes!'' if $T_M$ is real-analytic:

\begin{lemma}[Holomorphic continuation of tensors]\label{lem:holo_continuation_of_tensors}
 Let $X$ be a complex manifold with real structure $\sigma$ and non-empty real form $M\coloneqq\Fix\sigma$. Further, let $T_M$ be a real-analytic tensor field on $M$. Then, there exists an open neighborhood $U\subset X$ of $M$ and a $\sigma$-invariant holomorphic tensor field $T$ on $U$ which induces $T_M$. $T$ is unique in the sense of Lemma \autoref{lem:uniqueness_of_holomorphic_tensors}.
\end{lemma}

\begin{proof}
 Pick a point $p\in M$ and a $\sigma$-chart $(z_1,\ldots, z_n) = \psi:V\to\C^n$ near $p$. The real part of the $\sigma$-chart furnishes a real-analytic submanifold chart for $M$. In this chart, $T_M$ becomes:
 \begin{gather*}
  T_M = T^{\phantom{M,}\, j_1\ldots j_r}_{M,\, k_1\ldots k_s} \pa{x_{j_1}}\otimes\ldots\otimes\pa{x_{j_r}}\otimes dx_{k_1}\otimes\ldots\otimes dx_{k_s},
 \end{gather*}
 where we employ the same conventions as before. The coefficients $T^{\phantom{M,}\, j_1\ldots j_r}_{M,\, k_1\ldots k_s}$ are real-analytic maps from $V\cap M$ to $\R\subset\C$. By precomposing them with $\psi^{-1}\vert_{\R^n}$, we can view the coefficients as real-analytic maps from $\R^n$ to $\C$. As in the proof of Definition \autoref{def:real_form}, there exists a unique holomorphic continuation of these maps. This allows us to interpret the coefficients $T^{\phantom{M,}\, j_1\ldots j_r}_{M,\, k_1\ldots k_s}$ as holomorphic functions on $V$ (after shrinking $V$ if necessary). Denote these holomorphic coefficients by $T^{j_1\ldots j_r}_{k_1\ldots k_s}$. We can now define the tensor field $T$ on $V$ by:
 \begin{gather}\label{eq:def_of_T}
  T \coloneqq T^{j_1\ldots j_r}_{k_1\ldots k_s} \pa{z_{j_1}}\otimes\ldots\otimes\pa{z_{j_r}}\otimes dz_{k_1}\otimes\ldots\otimes dz_{k_s}.
 \end{gather}
 $T$ is holomorphic, because the coefficients $T^{j_1\ldots j_r}_{k_1\ldots k_s}$ are holomorphic. Furthermore, $T$ is also $\sigma$-invariant. To see this, note that the coefficients $T^{j_1\ldots j_r}_{k_1\ldots k_s}$ assume real values on $V\cap M$. With this, we can apply the identity theorem to $\overline{T^{j_1\ldots j_r}_{k_1\ldots k_s}\circ\sigma}$ and $T^{j_1\ldots j_r}_{k_1\ldots k_s}$ to obtain:
 \begin{gather*}
  \overline{T^{j_1\ldots j_r}_{k_1\ldots k_s}\circ\sigma} = T^{j_1\ldots j_r}_{k_1\ldots k_s}.
 \end{gather*}
 Inserting the last equation into \autoref{eq:def_of_T} and exploiting the properties of $\sigma$-charts now yields $\overline{\sigma^\ast T} = T$.\\
 To define $T$ on a neighborhood $U$ of $M$, we want to repeat the last step for every point $p\in M$ and afterwards combine the various coordinate patches. This is possible if the definition of $T$ given in \autoref{eq:def_of_T} is independent of the choice of $\sigma$-chart. However, the independence of the choice of chart directly follows from Lemma \autoref{lem:uniqueness_of_holomorphic_tensors} concluding the proof.
\end{proof}

To end \autoref{app:real_structures}, we point out that real structures are not exclusive to the category of complex manifolds, but exist in all kinds of complex subcategories, for instance in the category of complex Lie groups, in the category of holomorphic symplectic manifolds, in the category of holomorphic Kähler manifolds, and so forth. In these cases, the antiholomorphic involution $\sigma$ has to preserve the structure of the subcategory, e.g., has to be a group homomorphism,  has to preserve the holomorphic symplectic form $\Omega$ (i.e. $\overline{\sigma^\ast\Omega} = \Omega$), has to preserve the Kähler structure, and so forth. This enforces that the corresponding real form $M$ carries the remaining non-complex structure of the category, e.g., is a Lie group, is a symplectic manifold, is a Kähler manifold, and so forth.

 \chapter{Proof of Darboux's Theorem for HSMs}
 \label{app:darboux}
 In this part of the appendix, we want to show that there is a holomorphic counterpart to Darboux's theorem:

\begin{theorem}[Darboux's theorem for HSMs]
 Let $(X,\Omega)$ be a holomorphic symplectic manifold (HSM) of complex dimension $\text{\normalfont dim}_\mathbb{C}(X) = 2n$\linebreak $(n\in\mathbb{N})$. Then, for every point $x\in X$, there is a holomorphic chart\linebreak $\psi = (Q_1,\ldots, Q_n, P_1,\ldots, P_n):U\to V\subset\mathbb{C}^{2n}$ of $X$ near $x$ such that
 \begin{gather*}
  \Omega\vert_{U} = \sum\limits^n_{j = 1} dP_j\wedge dQ_j.
 \end{gather*}
\end{theorem}

\begin{proof}
 As in the real setup, we employ Moser's trick to prove Darboux's theorem for HSMs. A detailed transcription of Moser's trick to complex manifolds is given by Soldatenkov and Verbitsky (cf. \cite{soldatenkov2021}).\\
 Let $(X,\Omega)$ be a HSM of complex dimension $\text{\normalfont dim}_\mathbb{C}(X) = 2n$ ($n\in\mathbb{N}$), $x\in X$ be any point of $X$, and $(\hat U, \hat\psi = (\hat z_1,\ldots, \hat z_{2n}))$ be a holomorphic chart of $X$ near $x$. First, we observe that every complex symplectic form $\omega$ on a complex vector space $V$ of dimension $\text{\normalfont dim}_\mathbb{C}(V) = 2n$ can be brought into standard form, i.e., can be written as
 \begin{gather*}
  \omega = \sum\limits^n_{j=1} \theta_{j+n}\wedge\theta_{j}
 \end{gather*}
 for a basis $(\theta_1,\ldots,\theta_{2n})$ dual to some complex basis $(e_1,\ldots, e_{2n})$ of $V$. Thus, we can assume that $\Omega$ at $x$ in the chart $(\hat U,\hat\psi)$ takes the form
 \begin{gather*}
  \Omega_x = \sum\limits^n_{j = 1} d\hat z_{j+n,x}\wedge d\hat z_{j,x}
 \end{gather*}
 by applying a $\mathbb{C}$-linear transformation to $(\hat U,\hat\psi)$ if necessary. Next, we define the following $2$-form on $\hat U$:
 \begin{gather*}
  \Omega_1\coloneqq\sum\limits^n_{j=1} d\hat z_{j+n}\wedge d\hat z_{j} \in\Omega^{(2,0)}(\hat U).
 \end{gather*}
 Clearly, both $\Omega_0\coloneqq \Omega\vert_{\hat U}$ and $\Omega_1$ are holomorphic symplectic $2$-forms on $\hat U$. Now, we define the $2$-form $\Omega_t$ on $\hat U$ as the interpolation of $\Omega_0$ and $\Omega_1$:
 \begin{gather*}
  \Omega_t\coloneqq \Omega_0 + t(\Omega_1 - \Omega_0)\quad\forall t\in\mathbb{R}.
 \end{gather*}
 For every $t\in\mathbb{R}$, the form $\Omega_t$ is holomorphic and closed, as $\Omega_0$ and $\Omega_1$ are holomorphic and closed. Further, we observe that $\Omega_{t,x} = \Omega_{0,x}$ for every $t\in\mathbb{R}$, as $\Omega_0$ and $\Omega_1$ coincide at $x$ by construction. This means that $\Omega_{t,x}$ is non-degenerate for every $t$. Because non-degeneracy is an open property, we can find an open neighborhood $U^\prime\subset\hat U$ of $x$ such that $\Omega_t\vert_{U^\prime}$ is a non-degenerate $2$-form for every $t\in [0,1]$, where we have also used the fact that $\Omega_t$ depends continuously on $t$ and $[0,1]$ is a compact interval. This turns ($U^\prime$, $\Omega_t\vert_{U^\prime}$) into a HSM for every $t\in [0,1]$. Moreover, we can assume that $U^\prime$ is contractible by shrinking $U^\prime$ if necessary. For the sake of simplicity and ease of notation, we assume from now on that $U^\prime = \hat U$. As $\hat U$ is contractible, its cohomology is trivial. By the Poincar\'{e} lemma, we can find a $1$-form $\alpha$ on $\hat U$ such that
 \begin{gather*}
  \frac{d}{dt}\Omega_t = \Omega_1 - \Omega_0 = d\alpha.
 \end{gather*}
 We can write $\alpha = \alpha^{(1,0)} + \alpha^{(0,1)}$, where $\alpha^{(1,0)}$ and $\alpha^{(0,1)}$ are of type $(1,0)$ and $(0,1)$, respectively. Since the $2$-form $\Omega_1 - \Omega_0$ is holomorphic, it is of type $(2,0)$ implying $d\alpha = d\alpha^{(1,0)} = \partial\alpha^{(1,0)}$. Thus, we can choose $\alpha$ such that it is of type $(1,0)$ and fulfills $\bar{\partial}\alpha = 0$ turning $\alpha$ into a holomorphic $1$-form. Without loss of generality, we can assume that $\alpha$ satisfies $\alpha_x = 0$ by replacing $\alpha$ with $\alpha^\prime = \alpha - \beta$ if necessary, where $\beta$ is a holomorphic $1$-form with $\beta_x = \alpha_x$ and $d\beta = 0$. For every $t\in [0,1]$, define the vector field $V_t$ on $\hat U$ with values in $T^{(1,0)}\hat U$ by $\iota_{V_t}\Omega_t = -\alpha$. Note that $V_t$ is well-defined, as $\alpha$ is a $1$-form of type $(1,0)$ and $\Omega_t$ is non-degenerate on $T^{(1,0)}\hat U$, thus, $\Omega_t$ gives rise to an isomorphism from $T^{(1,0)}\hat U$ to $T^{\ast, (1,0)}\hat U$. Since $\alpha$ and $\Omega_t$ are holomorphic, $V_t$ is a holomorphic vector field on $\hat U$ for every $t\in [0,1]$. Using the closedness of $\Omega_t$ and Cartan's magic formula, we calculate the Lie derivative $L_{V_t}\Omega_t$:
 \begin{gather*}
  L_{V_t}\Omega_t = d\iota_{V_t}\Omega_t + \iota_{V_t}d\Omega_t = d\iota_{V_t}\Omega_t = -d\alpha = -\frac{d}{dt}\Omega_t\quad\forall t\in [0,1].
 \end{gather*}
 Recall\footnote{Confer Proposition \autoref{prop:holo_vec_field_equiv_J_pre_vec_field} and \cite{kobayashi1969} for details.} that every holomorphic vector field $V$ can be uniquely written as\linebreak $1/2(V^R - i\cdot J(V^R))$, where $J$ is the complex structure of the underlying complex manifold and $V^R$ is a real $J$-preserving vector field. Now let $V^R_t$ be the real $J$-preserving vector field corresponding to $V_t$ for every $t\in [0,1]$. Next, we want to show the following equation:
 \begin{gather*}
  L_{V^R_t}\Omega_t = L_{V_t}\Omega_t = -\frac{d}{dt}\Omega_t.
 \end{gather*}
 We do this by proving a more general auxiliary lemma:
 
 \begin{auxiliarylemma}\label{lem:aux_lem_darboux_1}
  Let $X$ be a complex manifold with complex structure $J\in\Gamma (\text{\normalfont End}(TX))$. Further, let $V^R\in\Gamma_J (TX)$ be a $J$-preserving vector field on $X$ with corresponding holomorphic vector field $V\coloneqq 1/2(V^R -iJ(V^R))$ and let $T$ be a holomorphic tensor field on $X$, then the Lie derivatives of $T$ with respect to $V^R$ and $V$ coincide, i.e.:
  \begin{gather*}
   L_{V^R} T = L_{V}T.
  \end{gather*}
 \end{auxiliarylemma}

 \begin{proof}[Proof of Lemma \autoref{lem:aux_lem_darboux_1}]
  Let $X$, $J$, $V^R$, and $V$ be as above. Further, let $T$ be a holomorphic $(k,l)$-tensor field on $X$. The Lie derivative is complex linear, thus, we have by definition of $V$:
  \begin{gather*}
   L_V T = \frac{1}{2}\left(L_{V^R}T - i\cdot L_{J(V^R)}T\right).
  \end{gather*}
  Hence, it suffices to show:
  \begin{gather}\label{eq:lie_i}
   L_{J(V^R)}T = i\cdot L_{V^R}T.
  \end{gather}
  In a holomorphic chart $\phi = (z_1,\ldots, z_n):U\to V\subset\mathbb{C}^n$ of $X$, $T$ can be expressed as
  \begin{gather*}
   T\vert_U = \sum^n_{i_1\ldots i_k, j_1\ldots j_l = 1} T^{i_1\ldots i_k}_{j_1\ldots j_l}\cdot dz_{j_1}\otimes\ldots\otimes dz_{j_l}\otimes \partial_{z_{i_1}}\otimes\ldots\otimes \partial_{z_{i_k}},
  \end{gather*}
  where $T^{i_1\ldots i_k}_{j_1\ldots j_l}:U\to\mathbb{C}$ are holomorphic functions on $U$. Since the Lie derivative can be computed locally and satisfies the Leibniz rule, it suffices to show Equation \eqref{eq:lie_i} for $T$ being $T^{i_1\ldots i_k}_{j_1\ldots j_l}$, $dz_{i}$, and $\partial_{z_j}$. As the Lie derivative also commutes with the exterior differential $d$ and $\partial_{z_j}$ is a (local) holomorphic vector field, it is sufficient to prove Equation \eqref{eq:lie_i} for $T$ being a holomorphic function $h$ and holomorphic vector field $W$. For $T = h$, we find:
  \begin{gather*}
   L_{J(V^R)} h = dh\left( J (V^R)\right) = i\cdot dh (V^R) = i\cdot L_{V^R}h,
  \end{gather*}
  where we used that $h$ is holomorphic, i.e. $dh\circ J = i\cdot dh$. For $T = W$, we can use Proposition \autoref{prop:holo_vec_field_equiv_J_pre_vec_field} to obtain:
  \begin{align*}
   L_{J(V^R)}W &= \left[J(V^R), W\right] = \frac{1}{2}\left([J(V^R), W^R] - i[J(V^R), J(W^R)]\right)\\
   &= \frac{1}{2}\left([V^R, J(W^R)] - i[V^R, J^2(W^R)]\right) = \frac{i}{2}\left([V^R, W^R] - i[V^R, J(W^R)]\right)\\
   &= i\left[V^R, W\right] = i\cdot L_{V^R}W,
  \end{align*}
  concluding the proof.
 \end{proof}
 
 Let us return to the proof of Darboux's theorem for HSMs. Let $\varphi_t$ be the (possibly local) flow of the real time-dependent vector field $V^R_t$. The pull-back $\varphi^\ast_t\Omega_t$ is a solution of the initial value problem:
 \begin{gather*}
  \frac{d}{dt}(\varphi^\ast_t \Omega_t) = \varphi^\ast_t (L_{V^R_t}\Omega_t + \frac{d}{dt}\Omega_t) = 0,\quad \varphi^\ast_0\Omega_0 = \Omega_0,
 \end{gather*}
 where we suppressed the fact that $\varphi_t$ might not be defined on all of $\hat U$ for every $t\in [0,1]$ in our notation. Clearly, $\Omega_0$ is also a solution to the same initial value problem. Because the solution to the given initial value problem is unique, we obtain:
 \begin{gather*}
  \varphi^\ast_t\Omega_t = \Omega_0.
 \end{gather*}
 We have to show that the last equation holds true for every $t\in [0,1]$ in some open neighborhood $U\subset\hat U$ of $x$. For this, we recall that $\alpha_x = 0$ by construction. Thus, we have $V^R_t (x) = V_t (x) = 0$ for every $t\in [0,1]$. This implies that the flow $\varphi_t$ is stationary at $x$, i.e., $\varphi_t (x) = x$ for every $t\in [0,1]$. We can deduce from this that there exists an open neighborhood $U\subset\hat U$ of $x$ such that the flow $\varphi_t: U\to \varphi_t(U)\subset\hat U$ is a well-defined diffeomorphism for every $t\in [0,1]$. In particular, the time-$1$-map $\varphi_1:U\to\varphi_1 (U)$ satisfies:
 \begin{gather*}
  \varphi^\ast_1\Omega_1 = \Omega_0.
 \end{gather*}
 Hence, $(U,\psi\coloneqq\hat \psi\circ\varphi_1)$ is a smooth chart of $X$ near $x$ which satisfies:
 \begin{gather*}
  \psi^{-1\,\ast}\Omega\vert_U = {\hat\psi}^{-1\,\ast}\left(\varphi^{-1\,\ast}_1 \Omega_0\right) = {\hat\psi}^{-1\,\ast}\Omega_1 = \sum^n_{j=1} \theta_{j+n}\wedge \theta_{j},
 \end{gather*}
 where $\sum \theta_{j+n}\wedge\theta_{j}$ is the standard symplectic form on $\mathbb{C}^{2n}$. Hence, the holomorphic symplectic form $\Omega$ takes the following form on $(U,\psi \equiv (z_1,\ldots, z_{2n}))$:
 \begin{gather*}
  \Omega\vert_{U} = \sum\limits^n_{j = 1} dz_{j+n}\wedge dz_{j}.
 \end{gather*}
 $(U,\psi)$ is a good candidate for the desired Darboux chart. To conclude the proof, we need to show that $(U,\psi)$ is also a holomorphic chart of $X$. For this, it suffices to prove that the map $\varphi_1: U\to \varphi_1 (U)$ is locally biholomorphic. The idea behind this proof is simple: In Chapter IX of \cite{kobayashi1969}, it is shown that the time-independent $J$-preserving\footnote{$J$-preserving vector fields are called infinitesimal automorphisms in \cite{kobayashi1969}.} vector fields $V^R$ on a complex manifold $X$ are exactly those real vector fields whose flow $\varphi^{V^R}_t$ is locally biholomorphic. However, we cannot directly apply this statement to $V^R_t$, as, in general, $V^R_t$ carries a non-trivial time-dependence. To account for this, we relate $V_t$ to a time-independent holomorphic vector field $V$ on $\hat U \times O \ni (x,t)$, where $O\subset\mathbb{C}$.\\
 First, we generalize the definition of $\Omega_t$ and allow for complex times $t\equiv\tau\in\mathbb{C}$. By the same arguments as before and by shrinking $\hat U$ if necessary, we find an open neighborhood $O\subset\mathbb{C}$ of $[0,1]$ such that $\Omega_\tau$ is non-degenerate for every $\tau\in O$. This allows us to generalize the definition of $V_\tau$ to all $\tau\in O$. Instead of viewing $\Omega_\tau$ as a time-dependent $(2,0)$-form on $\hat U$, we can also take it to be a time-independent $(2,0)$-form on $\hat U\times O$. As the time-dependence of $\Omega_\tau$ is clearly holomorphic, $\Omega_\tau$ as a form on $\hat U\times O$ is also holomorphic. Thus, $V_\tau$ understood as a vector field on $\hat U\times O$ is also holomorphic. Now consider the time-independent vector field $V$ on $\hat U\times O$:
 \begin{gather*}
  V(y,\tau^\prime)\coloneqq V_{\tau^\prime}(y) + \partial_\tau(y,\tau^\prime)\quad\forall (y,\tau^\prime)\in\hat U\times O,
 \end{gather*}
 where we interpret $\tau$ as a coordinate of $\hat U\times O$. As $V_\tau$ and $\partial_\tau$ are holomorphic, $V$ is also holomorphic. Now let $V^R$ be the real $J$-preserving vector field corresponding to $V$ and let $\Gamma:[0,1]\to\hat U\times O$, $\Gamma (r)\equiv (\gamma (r), \rho (r))$ be a smooth curve in $\hat U\times O$ with $\rho (0) = 0$. Then, we have the following auxiliary lemma:
 
 \begin{auxiliarylemma}\label{lem:aux_lem_darboux_2}
  $\Gamma$ is an integral curve of $V^R$ if and only if $\gamma:[0,1]\to\hat U$ is an integral curve of $V^R_t$ and $\rho (r) = r$ for every $r\in [0,1]$.
 \end{auxiliarylemma}
 
 \begin{proof}[Proof of Lemma \autoref{lem:aux_lem_darboux_2}]
  This follows from a quick computation: Let $\tau = t + is$ be the decomposition of $\tau$ into real and imaginary part, then the vector field $V^R$ is given by:
  \begin{gather*}
   V^R(y,\tau^\prime) = V^R_{\tau^\prime}(y) + \partial_t(y,\tau^\prime)\quad\forall (y,\tau^\prime)\in\hat U\times O.
  \end{gather*}
  Further, let $\rho = \rho_R + i\rho_I$ be the decomposition of $\rho$ into real and imaginary part, then the integral curve equation of $V^R$ for $\Gamma$ can be written as ($r\in [0,1]$):
  \begin{gather*}
   \dot\gamma (r) = V^R_{\rho (r)} (\gamma (r)),\quad \dot \rho_R (r)\cdot\partial_t(\Gamma (r)) + \dot \rho_I (r)\cdot\partial_s(\Gamma (r)) = \partial_t(\Gamma (r)).
  \end{gather*}
  As $\rho$ needs to satisfy the initial condition $\rho (0) = 0$, $\rho$ is given by $\text{id}_{[0,1]}$ if $\Gamma$ is an integral curve of $V^R$. In this case, $\gamma$ has to satisfy the following differential equation:
  \begin{gather*}
   \dot\gamma (r) = V^R_{r} (\gamma (r))\quad r\in [0,1].
  \end{gather*}
  Therefore, $\gamma$ is an integral curve of the time-dependent vector field $V^R_t$ on $\hat U$ if $\Gamma$ is an integral curve of $V^R$. The converse direction follows similarly.
 \end{proof}
 
 From Lemma \autoref{lem:aux_lem_darboux_2}, it follows that the flow $\varphi^{V^R}_t$ of $V^R$ and the flow $\varphi_t$ of the time-dependent vector field $V^R_t$ are related in the following way:
 \begin{gather*}
  \varphi^{V^R}_t (y,0) = \left(\varphi_t (y), t\right)
 \end{gather*}
 for every suitable $y\in\hat U$. As discussed earlier, $\varphi^{V^R}_t$ is the flow of a holomorphic vector field and, hence, locally biholomorphic. This implies that $\varphi_t$ is also locally biholomorphic concluding the proof.
\end{proof}

 \chapter{Three Flavors of Kähler Structures}
 \label{app:kaehler}
 We briefly explore the various notions of Kähler structures in this part. In particular, we introduce three different flavors of Kähler structures: Kähler manifolds (first section), Hyperkähler manifolds (second section), and holomorphic Kähler manifolds (third section). The main result of \autoref{app:kaehler} is that, in the symplectic picture, Hyperkähler and holomorphic Kähler manifolds only differ by a sign in the commutation relation of the complex structures $I$ and $J$.

\section*{Kähler Manifolds}

\begin{definition}[Kähler manifolds]\label{def:kaehler}
 A \textbf{pre-semi-Kähler} manifold is a triple $(M,\omega, J)$ where $M^{2n}$ is a smooth manifold and the tensors $\omega\in\Omega^2 (M)$ and $J\in\Gamma\End (TM)$ satisfy:
 \begin{enumerate}
  \item $\omega$ is non-degenerate, i.e., $\omega^n_p\neq 0$ for all $p\in M$,
  \item $J$ is an almost complex structure, i.e., $J^2_p = -\id_{T_pM}$ for all $p\in M$,
  \item $\omega$ and $J$ are compatible in the sense that $\omega (J\cdot , J\cdot) = \omega$.
 \end{enumerate}
 We drop the prefix ``pre'' if $(M,\omega, J)$ satisfies the following integrability conditions:
 \begin{enumerate}[label = (\arabic*)]
  \item $\omega$ is closed, i.e., $d\omega = 0$,
  \item $J$ is integrable, i.e., the Nijenhuis tensor $N_J$ vanishes.
 \end{enumerate}
 We drop the prefix ``semi'' if the semi-Riemannian metric $g\coloneqq \omega (\cdot, J\cdot)$ is positive definite.
\end{definition}

\begin{remark}\label{rem:kaehler}\ \vspace{-0.2cm}
 \begin{enumerate}
  \item It is easy to verify that, given a pre-semi-Kähler manifold $(M,\omega, J)$, the tensor $g\coloneqq \omega (\cdot, J\cdot)$ is a semi-Riemannian metric satisfying $g(J\cdot,J\cdot) = g$. Conversely, if $M$ is a smooth manifold with semi-Riemannian metric $g$ and almost complex structure $J$ satisfying $g(J\cdot, J\cdot) = g$, then the triple $(M,\omega, J)$ is a pre-semi Kähler manifold where $\omega\coloneqq g(J\cdot, \cdot)$. In light of this observation, we call both triples $(M,\omega, J)$ and $(M,g,J)$ pre-semi-Kähler manifolds. Whether the symplectic or Riemannian picture is used, is clear from the symbols and the context.
  \item A manifold $M$ admits a pre-Kähler structure iff $M$ admits an almost complex structure iff $M$ admits a non-degenerate two-form. In particular, every pre-semi-Kähler manifold admits a pre-Kähler structure. To see this, we first note that $M$ admits an almost complex (or almost symplectic or pre-Kähler) structure iff the frame bundle of $M$ admits a $\GL (n,\C)$- (or $\Sp (2n,\R)$- or $\U (n)$-)reduction. The rest now follows from\linebreak $\Sp (2n,\R)\cap \GL (n,\C) = \U (n)$ and the fact that the inclusions\linebreak $\U (n)\hookrightarrow \GL (n,\C)$ and $\U (n)\hookrightarrow \Sp (2n)$ are strong deformation retracts as well as the observation that $G$-principal bundles admit $H$-reductions if $H\hookrightarrow G$ is a strong deformation retract of Lie groups.
 \end{enumerate}
\end{remark}

Kähler manifolds are one of the most studied geometries in mathematics. As such, a plethora of equivalent descriptions of Kähler manifolds has been found over the years (cf. \cite{Ballmann2006} and \cite{Moroianu2007} for a small selection). Here, we are only interested in the following equivalences:

\begin{lemma}[Equivalent description of integrability conditions]\label{lem:kaehler}
 Let $(M,\omega,J)$ be a pre-semi-Kähler manifold and $\nabla$ be the Levi-Civita connection of\linebreak $g\coloneqq\omega (\cdot, J\cdot)$. Then, the following statements are equivalent:
 \begin{enumerate}
  \item $(M,\omega, J)$ is a semi-Kähler manifold.
  \item $\nabla J = 0$.
  \item $\nabla\omega = 0$.
  \item $M$ admits an atlas of holomorphic normal coordinates, i.e., for every point $p\in M$ there exists a holomorphic chart $(z_1 = x_1 + iy_1,\ldots, z_n = x_n + iy_n)$ near $p$ such that
  \begin{gather*}
   g_p = \sum^k_{j=1} dx^2_{j,p} + dy^2_{j,p} - \sum^n_{j = k+1} dx^2_{j,p} + dy^2_{j,p}
  \end{gather*}
  and the first derivatives of the coefficients of $g$ vanish at $p$ in this chart.
  \item $J$ is integrable and $\omega$ possesses local Kähler potentials $f$ near any point, i.e., for every $p\in M$ there exists an open neighborhood $U$ of $p$ and a function $f\in C^\infty (U,\R)$ such that $\omega\vert_U = i\partial\bar{\partial}f$.
 \end{enumerate}
\end{lemma}

\begin{proof}
 ``(i)$\Leftrightarrow$(ii)'': This equivalence is Lemma 4.2.5 in \cite{mcduff2017}.\\
 ``(ii)$\Leftrightarrow$(iii)'': We compute for vector fields $X,Y,Z\in\Gamma (TM)$:
 \begin{align*}
  X(\omega (Y,Z)) &= (\nabla_X\omega) (Y,Z) + \omega (\nabla_X Y, Z) + \omega (Y,\nabla_XZ)\\
  = X(g(JY, Z)) &= (\nabla_Xg) (JY,Z) + g(\nabla_X(JY), Z) + g (JY,\nabla_XZ)\\
  &= g((\nabla_X J)Y, Z) + g (J\nabla_XY,Z) + g (JY,\nabla_XZ)\\
  &= g((\nabla_XJ)Y,Z) + \omega (\nabla_X Y, Z) + \omega (Y,\nabla_XZ),
 \end{align*}
 where we exploited the metric compatibility of $\nabla$. Subtracting the first from the last line yields $(\nabla_X\omega) (Y,Z) = g((\nabla_XJ)Y,Z)$ proving the equivalence.\\
 The remaining directions are proven in \cite{Ballmann2006} (cf. Theorem 4.17). Note that the equivalences are only shown for positive definite $g$ in \cite{mcduff2017} and \cite{Ballmann2006}. However, it is straightforward to check that the proofs still work in the semi-Riemannian case.
\end{proof}

\begin{remark}[Signature of $g$]\label{rem:kaehler_sig}
 A consequence of Lemma \autoref{lem:kaehler} is that the signature of a semi-Kähler metric $g$ is even, i.e., $g$ has signature $(2k, 2(n-k))$. One easily verifies that this is also true for pre-semi-Kähler manifolds $(M,g,J)$.
\end{remark}

The Kähler potential $f$ in Lemma \autoref{lem:kaehler} is not unique. For instance, adding the real part of a holomorphic function $g$ to $f$ generates another Kähler potential, i.e., if $g:U\to\C$ is a holomorphic function, then $\omega = i\partial\bar{\partial}\hat f$ where:

\begin{gather*}
\hat f\coloneqq f + \text{Re}\, g = f + \frac{1}{2}(g + \bar{g}).
\end{gather*}

To single out a unique Kähler potential, we need to impose additional conditions on $f$. One such condition is to require $f$ to be mixed near a given point $p$:

\begin{definition}[Mixed near $p$]\label{def:mixed}
 Let $M$ be a complex manifold and $f:M\to\R$ a real-analytic function on $M$. Fix a point $p\in M$ and a holomorphic chart $\psi = (z_1,\ldots,z_n)$ near $p$ with $\psi (p) = 0$. Expanding $f$ in a power series of $z_j$ and $\bar{z}_k$ gives us the decomposition:
 \begin{gather*}
  f = h_0 + h_1 + h_2 + h_3,
 \end{gather*}
 where $h_0 \equiv f(p)$ is a constant function, $h_1$ only contains powers of $z_j$, $h_2$ only contains powers of $\bar{z}_k$, and $h_3$ contains the terms mixing $z_j$ and $\bar{z}_k$. We call $f$ \textbf{mixed} near $p$ if $h_0 = h_1 = h_2 = 0$.
\end{definition}

\begin{remark}\label{rem:mixed}\ \vspace{-0.2cm}
 \begin{enumerate}
  \item The decomposition $f = h_0 + h_1 + h_2 + h_3$ is independent of the choice of chart $\psi$, as long as $\psi (p) = 0$ holds.
  \item If $f$ is mixed near $p$, then $p$ is a critical point of $f$. Indeed, all terms contained in $f$ are of quadratic or higher order in the coordinates $z_j$ and $\bar{z}_k$.
 \end{enumerate}
\end{remark}

Definition \autoref{def:mixed} now allows us to find a unique Kähler potential near a given point $p$:

\begin{proposition}[Mixed Kähler potential]\label{prop:mixed_kaehler_pot}
 Let $M$ be a complex manifold and $\omega\in\Omega^2 (M)$ a closed real-analytic $\R$-valued form of type $(1,1)$. Fix a point $p\in M$. Then, there exists a contractible neighborhood $U$ of $p$ and a unique real-analytic function $f:U\to\R$ such that $f$ is mixed near $p$ and $\omega = i\partial\bar{\partial}f$ holds. In particular, every real-analytic semi-Kähler manifold $(M,\omega,J)$ possesses a unique local Kähler potential $f$ which is mixed near $p$.
\end{proposition}

\begin{proof}
 By the $\partial\bar{\partial}$-lemma, there is a contractible neighborhood $U$ of $p$ and a function $f:U\to\R$ such that $\omega = i\partial\bar{\partial}f$. Since $\omega$ is real-analytic, we can also choose $f$ to be real-analytic. After shrinking $U$, $U$ becomes the domain of a chart $(z_1,\ldots,z_n)$ and we can write $f = h_0 + h_1 + h_2 + h_3$ as in Definition \autoref{def:mixed}. $h_1$ is holomorphic and, since $f$ is real, we have $h_2 = \bar{h}_1$. In particular, $dh_0 = \bar{\partial}h_1 = \partial h_2 = 0$ holds and $h_3$ satisfies $\omega = i\partial\bar{\partial}h_3$. Thus, $h_3$ is the mixed Kähler potential for $\omega$. Conversely, $\omega$ fixes $h_3$, as $\omega$ determines the functions $\pa{z_j}\pa{\bar{z}_k}h_3$ and every term in $h_3$ is proportional to $z_j\bar{z}_k$ concluding the proof.
\end{proof}

\section*{Hyperkähler Manifolds}

\begin{definition}[Hyperkähler manifolds]\label{def:hyperkaehler}
 A \textbf{pre-semi-Hyperkähler} manifold is a collection $(M,g,I,J,K)$ where $M^{4n}$ is a smooth manifold and the tensors $g\in\Sym^2 (M)$ and $I,J,K\in\Gamma\End (TM)$ satisfy:
 \begin{enumerate}
  \item $g$ is a semi-Riemannian metric,
  \item $I$, $J$, and $K$ form an almost quaternionic structure:
  \begin{gather*}
   I^2_p = J^2_p = K^2_p = I_pJ_pK_p = -\id_{T_pM}\quad\forall p\in M,
  \end{gather*}
  \item $g$ is compatible with $I$, $J$, and $K$ in the following sense:
  \begin{gather*}
   g(I\cdot,I\cdot) = g(J\cdot,J\cdot) = g(K\cdot,K\cdot) = g.
  \end{gather*}
 \end{enumerate}
 We drop the prefix ``pre'' if $(M,g,I)$, $(M,g,J)$, and $(M,g,K)$ are semi-Kähler.\\
 We drop the prefix ``semi'' if $g$ is positive definite.
\end{definition}

\begin{remark}\label{rem:hyperkaehler}\ \vspace{-0.2cm}
 \begin{enumerate}
  \item Given a pre-semi-Hyperkähler manifold $(M,g,I,J,K)$, every almost complex structure $L\in\{I,J,K\}$ corresponds to a non-degenerate $2$-form $\omega_L\coloneqq g(L\cdot,\cdot)$. The forms $\omega_I$, $\omega_J$, and $\omega_K$ allow us to construct the complex forms:
  \begin{gather*}
   \Omega_I\coloneqq \omega_J + i\omega_K,\ \Omega_J\coloneqq \omega_K + i\omega_I,\ \Omega_K\coloneqq \omega_I + i\omega_J.
  \end{gather*}
  For any $L\in\{I,J,K\}$, $\Omega_L$ is of type $(2,0)$ with respect to the almost complex structure $L$, i.e.:
  \begin{gather*}
   \Omega_L (L\cdot,\cdot) = \Omega_L (\cdot, L\cdot) = i\Omega_L.
  \end{gather*}
  If $(M,g,I,J,K)$ is semi-Hyperkähler, then $\Omega_L$ is a holomorphic symplectic form with respect to $L\in\{I,J,K\}$ (cf. Theorem \autoref{thm:rel_HSM_PHSM}).
  \item A manifold $M$ admits a pre-Hyperkähler structure iff $M$ admits an almost quaternionic structure. In particular, every pre-semi-Hyperkähler manifold admits a pre-Hyperkähler structure. As in the Kähler case, the proof relies on the fact that the inclusion $\GL (n,\Qua)\cap\Or (4n)\hookrightarrow \GL (n,\Qua)$ is a strong deformation retract.
 \end{enumerate}
\end{remark}

Definition \autoref{def:hyperkaehler} is somewhat redundant: To define Hyperkähler structures, we introduced three almost complex structures, however, two anticommuting ones completely suffice. Indeed, if $I$, $J$, and $K$ form an almost quaternionic structure on $M$, one easily concludes that $I$ and $J$ anticommute. The converse is also true:

\begin{proposition}\label{prop:anticommuting_acs}
 Let $M$ be a smooth manifold with two almost complex structures $I$ and $J$ and set $K\coloneqq IJ$. If $IJ = -JI$, then $I$, $J$, and $K$ form an almost quaternionic structure on $M$. If, additionally, $g$ is a semi-Riemannian metric on $M$ such that $(M,g,I)$ and $(M,g,J)$ are pre-semi-Kähler manifolds, then $(M,g,I,J,K)$ is a pre-semi-Hyperkähler manifold.
\end{proposition}

\begin{proof}
 We verify:
 \begin{gather*}
  K^2 = IJK = IJIJ = -I^2J^2 = -\mathds{1}.
 \end{gather*}
 The $K$-compatibility of $g$ is a direct consequence of the $I$- and $J$-compatibility:
 \begin{gather*}
  g(K\cdot,K\cdot) = g(IJ\cdot,IJ\cdot) = g(J\cdot,J\cdot) = g.
 \end{gather*}
\end{proof}

Not only the third almost complex structure $K$ is redundant, but also the third integrability condition:

\begin{proposition}\label{prop:int_con_hyper}
 Let $(M,g,I,J,K)$ be a pre-semi-Hyperkähler manifold. Then, the following statements are equivalent:
 \begin{enumerate}
  \item $(M,g,I,J,K)$ is semi-Hyperkähler.
  \item $(M,g,I)$ and $(M,g,J)$ are semi-Kähler.
 \end{enumerate}
\end{proposition}

\begin{proof}
 The direction ``(i)$\Rightarrow$(ii)'' is trivial, thus, we only consider ``(ii)$\Rightarrow$(i)''. If $(M,g,I)$ and $(M,g,J)$ are semi-Kähler, then, by Lemma \autoref{lem:kaehler}, we have\linebreak $\nabla I = \nabla J = 0$ for the Levi-Civita connection $\nabla$ of $g$. Therefore, the almost complex structure $K = IJ$ satisfies:
 \begin{gather*}
  \nabla K = (\nabla I)J + I\nabla J = 0
 \end{gather*}
 Applying Lemma \autoref{lem:kaehler} again now shows that $(M,g,K)$ is semi-Kähler concluding the proof.
\end{proof}

Proposition \autoref{prop:int_con_hyper} is a fundamental and well-known result in the field of Hyperkähler geometry. However, Theorem \autoref{thm:rel_HSM_PHSM} allows us to improve Proposition \autoref{prop:int_con_hyper}:

\begin{lemma}[Equivalent description of integrability conditions]\label{lem:hyperkaehler}
 Let \linebreak $(M,g,I,J,K)$ be a pre-semi-Hyperkähler manifold with corresponding two-forms $\omega_I$, $\omega_J$, and $\omega_K$. If two almost complex structures in $\{I,J,K\}$ are integrable and two forms in $\{\omega_I,\omega_J,\omega_K\}$ are closed, then $(M,g,I,J,K)$ is semi-Hyperkähler.
\end{lemma}

\begin{proof}
 If the integrable almost complex structures correspond to the closed forms, we can apply Proposition \autoref{prop:int_con_hyper}. Let us now assume that the integrable almost complex structures do not correspond to the closed forms. Without loss of generality, we can take $I$ and $J$ to be integrable and $\omega_I$ and $\omega_K$ to be closed. First, consider $i\Omega_I = -\omega_K + i\omega_J$. It is a form of type $(2,0)$ with respect to the complex structure $I$ whose real part is closed. Therefore, by Theorem \autoref{thm:rel_HSM_PHSM}, the imaginary part $\omega_J$ must be closed as well. Now consider $\Omega_K = \omega_I + i\omega_J$. It is a closed form of type $(2,0)$ with respect to $K$. Again by Theorem \autoref{thm:rel_HSM_PHSM}, $K$ must be integrable concluding the proof.
\end{proof}

As in the Kähler case, we cannot only describe a Hyperkähler manifold in a Riemannian setup, but also in a symplectic setup:

\begin{lemma}[Symplectic picture of Hyperkähler manifolds]\label{lem:symplectic_picture}
 Let $(M,\omega,J)$ be a pre-semi-Kähler manifold. Further, let $I$ be an almost complex structure on $M$ which is anticompatible with $\omega$ and anticommutes with $J$, i.e.:
 \begin{gather*}
  \omega (I\cdot,I\cdot) = -\omega\quad\text{and}\quad IJ = -JI.
 \end{gather*}
 Then, $(M,g,I,J,K)$ is a pre-semi-Hyperkähler manifold where $g\coloneqq\omega (\cdot,J\cdot)$ and $K\coloneqq IJ$. If $(M,\omega,J)$ is semi-Kähler and $I$ is integrable, then $(M,g,I,J,K)$ is semi-Hyperkähler.
\end{lemma}

\begin{proof}
 We already know from Proposition \autoref{prop:anticommuting_acs} that $I$, $J$, and $K$ form an almost quaternionic structure. On top of that, it directly follows from the definition of a pre-semi-Kähler manifold that $g$ is compatible with $J$, i.e., $g$ satisfies\linebreak $g(J\cdot,J\cdot) = g$. Since $I$ is anticompatible with $\omega$ and anticommutes with $J$, we also have $g(I\cdot,I\cdot) = g$ which in turn implies $g(K\cdot,K\cdot) = g$. This shows that $(M,g,I,J,K)$ is a pre-semi-Hyperkähler manifold.\\
 Now assume that$(M,\omega,J)$ is semi-Kähler and $I$ is integrable. By Theorem \autoref{thm:rel_HSM_PHSM}, the form
 \begin{gather*}
  \omega_K\coloneqq g(K\cdot,\cdot) = \omega (IJ\cdot, J\cdot) = -\omega (I\cdot,\cdot)
 \end{gather*}
 is closed. This allows us to apply Lemma \autoref{lem:hyperkaehler} (where $\omega_J\equiv\omega$) showing that $(M,g,I,J,K)$ is semi-Hyperkähler.
\end{proof}

To conclude \autoref{app:kaehler}, we complexify the notion of Kähler manifolds:

\section*{Holomorphic Kähler Manifolds}

\begin{definition}[Holomorphic semi-Kähler manifolds]\label{def:holo_kaehler}
 A \textbf{complexified pre-semi-Kähler} manifold is a collection $(X,\omega, J, I)$ where $X^{4n}$ is a smooth manifold and the tensors $\omega\in\Omega^2 (X)$ and $I,J\in\Gamma\End (TX)$ satisfy:
 \begin{enumerate}
  \item $(X,\omega, J)$ is a pre-semi-Kähler manifold,
  \item $I$ is an almost complex structure, i.e., $I^2_p = -\id_{T_pX}$ for all $p\in X$,
  \item $I$ is anticompatible with $\omega$ and commutes with $J$:
  \begin{gather*}
   \omega (I\cdot, I\cdot) = -\omega\quad\text{and}\quad IJ = JI.
  \end{gather*}
 \end{enumerate}
 We say that $(X,\omega,J,I)$ is \textbf{holomorphic} instead of complexified if $I$ is integrable and $\Omega\coloneqq \omega -i\omega (I\cdot,\cdot)$ as well as $J$ viewed as a section\footnote{Here, the subscript indicates that the decomposition $T_\C X = T^{(1,0)}X\oplus T^{(0,1)}X$ is understood with respect to $I$.} of $\End (T^{(1,0)}_IX)$ is holomorphic. We drop the prefix ``pre'' if $(X,\omega,J)$ is semi-Kähler.
\end{definition}

\begin{remark}\label{rem:comp_kaehler}\ \vspace{-0.2cm}
 \begin{enumerate}
  \item For every complexified pre-semi-Kähler manifold $(X,\omega, J, I)$, the forms $\Omega\coloneqq \omega -i\omega (I\cdot,\cdot)$ and $G\coloneqq \Omega (\cdot,J\cdot)\equiv g -ig (I\cdot,\cdot)$ where $g\coloneqq\omega (\cdot, J\cdot)$ are of type $(2,0)$ with respect to $I$:
  \begin{gather*}
   \Omega (I\cdot,\cdot) = \Omega (\cdot,I\cdot) = i\Omega\quad\text{and}\quad G (I\cdot,\cdot) = G (\cdot,I\cdot) = iG.
  \end{gather*}
  Furthermore, $\Omega$ and $G$ are non-degenerate on $T^{(1,0)}_IX$ implying that the real dimension of $X$ is a multiple of $4$.
  \item If $\omega$ is closed and $I$ is integrable, $\Omega$ is closed and holomorphic due to Theorem \autoref{thm:rel_HSM_PHSM}. Hence, every holomorphic semi-Kähler manifold is also a holomorphic symplectic manifold (cf. \autoref{sec:HHS}).
  \item The symmetric two-form $G$ is $\C$-bilinear with respect to $I$ and\linebreak non-degenerate on $T^{(1,0)}_IX$. Thus, its real part $g$ is indefinite with signature $(2n,2n)$ where $4n\coloneqq \dim_\R X$. The indefiniteness of $g$ is the reason why the prefix ``semi'' carries a different meaning for holomorphic Kähler manifolds than for Kähler and Hyperkähler manifolds (cf. Definition \autoref{def:semi}).
  \item In the lowest possible dimension, i.e. $\dim_\R X = 4$, the almost complex structure $J$ of a holomorphic pre-semi-Kähler manifold is automatically integrable. To see this, recall that the Nijenhuis tensor of an almost complex structure on a two-dimensional manifold naturally vanishes due to its symmetries. Similarly, the Nijenhuis tensor of a holomorphic $(1,1)$-tensor $J$ with $J^2 = -\mathds{1}$ on a complex two-dimensional manifold vanishes due to its symmetries.
  \item Also note that for every holomorphic pre-semi-Kähler manifold $(X,\omega,J,I)$ in dimension $\dim_\R X = 4$ the form $\Omega$ is a holomorphic top degree form and, thus, closed. Together with the previous remark, this implies that holomorphic pre-semi-Kähler manifolds in dimension 4 are automatically semi-Kähler.
 \end{enumerate}
\end{remark}

The rich structure of a holomorphic semi-Kähler manifold allows us to find a simple description of $I$ and $J$ in suitable coordinates:

\begin{proposition}[Local structure of $I$ and $J$]\label{prop:I-J-coordinates}
 Let $X^{2m}$ be a smooth manifold with commuting almost complex structures $I$ and $J$. If $I$ and $J$ are integrable and $J$ viewed as a section of $\End (T^{(1,0)}_IX)$ is holomorphic, then there are holomorphic coordinates $(z_1,\ldots,z_m)$ of $(X,I)$ near any point and a number $p$ with $0\leq p \leq m$ such that:
 \begin{gather*}
  J\pa{z_j} = i\pa{z_j}\quad\text{and}\quad J\pa{z_k} = -i\pa{z_k},
 \end{gather*}
 where $j\in\{1,\ldots, p\}$ and $k\in\{p+1,\ldots,m\}$. In particular, every holomorphic semi-Kähler manifold admits such coordinates.
\end{proposition}

\begin{remark}\label{rem:I-J-coordinates}\ \vspace{-0.2cm}
 \begin{enumerate}
  \item The number $p$ is constant on connected components of $X$.
  \item We call coordinates as in Proposition \autoref{prop:I-J-coordinates} $\mathbf{I}$\textbf{-}$\mathbf{J}$\textbf{-coordinates}.
  \item A consequence of Proposition \autoref{prop:I-J-coordinates} is that $I$ viewed as a section of\linebreak $\End (T^{(1,0)}_JX)$ is holomorphic. Indeed, $(z_1,\ldots, z_p, \bar{z}_{p+1},\ldots, \bar{z}_m)$ are holomorphic coordinates of $(X,J)$ in which the coefficients of $I$ are constant and, thus, holomorphic.
  \item For every holomorphic semi-Kähler manifold $(X^{4n},\omega,J,I)$, $\Omega$ takes the following form in $I$-$J$-coordinates:
  \begin{gather*}
   \Omega = \sum^p_{j=1}\sum^{2n}_{k=p+1}\Omega_{jk}dz_j\wedge dz_k.
  \end{gather*}
  This is due to the fact that all other combinations of $dz_r\wedge dz_s$ are not compatible with $J$. The given form of $\Omega$ enforces $p=n$, since, otherwise, $\Omega$ would be degenerate.
 \end{enumerate}
\end{remark}

\begin{proof}
 The idea of the proof is to apply the holomorphic Frobenius theorem to suitable holomorphic subbundles of $T^{(1,0)}_IX$ (the following version of this theorem together with its explanation is directly taken from \cite{voisin2002}):
 
 \begin{theorem}[Holomorphic Frobenius theorem, Theorem 2.26 in \cite{voisin2002}]
  Let $X$ be a complex manifold of dimension $n$ and let $E$ be a holomorphic distribution of rank $k$ over $X$, i.e., a holomorphic vector subbundle of rank $k$ of the holomorphic tangent bundle $T^{(1,0)}X$. Then, $E$ is integrable in the holomorphic sense if and only if we have the integrability condition
  \begin{gather*}
   [E,E]\subset E.
  \end{gather*}
 \end{theorem}
 Here, the integrability in the holomorphic sense means that $X$ is covered by open sets $U$ such that there exists a holomorphic submersive map
 \begin{gather*}
  \phi_U:U\to\C^{n-k}
 \end{gather*}
 satisfying
 \begin{gather*}
  E_u = \ker d\phi_{U,u}
 \end{gather*}
 for every $u\in U$.\\\\
 We first note that $J$ restricts to a section of $\End (T^{(1,0)}_IX)$, as $I$ and $J$ commute. Due to $J^2=-\mathds{1}$, the bundle $T^{(1,0)}_IX$ splits into the subbundles $E_i$ and $E_{-i}$ where $E_{\pm i}\subset T^{(1,0)}_IX$ is the eigenbundle of $J$ with respect to the eigenvalue $\pm i$. $J$ is holomorphic, hence, $E_i$ and $E_{-i}$ are holomorphic subbundles.\\
 To apply the holomorphic Frobenius theorem, we need to show that the subbundles $E_i$ and $E_{-i}$ are involutive. For this, we observe that $E_i$ and $E_{-i}$ can be expressed as:
 \begin{gather*}
  E_i = T^{(1,0)}_IX\cap T^{(1,0)}_JX,\quad E_{-i} = T^{(1,0)}_IX\cap T^{(0,1)}_JX.
 \end{gather*}
 Now recall that integrability of an almost complex structure $K$ on $X$ is equivalent to
 \begin{gather*}
  \left[T^{(1,0)}_KX,T^{(1,0)}_KX \right]\subset T^{(1,0)}_KX
 \end{gather*}
 which itself is equivalent to
 \begin{gather*}
  \left[T^{(0,1)}_KX,T^{(0,1)}_KX \right]\subset T^{(0,1)}_KX.
 \end{gather*}
 Thus, the integrability of $I$ and $J$ implies $[E_i, E_i]\subset E_i$ and $[E_{-i},E_{-i}]\subset E_{-i}$. This allows us to apply the holomorphic Frobenius theorem to $E_i$ and $E_{-i}$ giving us the desired charts, where $p$ is the rank of $E_i$.
\end{proof}

As already alluded to in Statement (iii) of Remark \autoref{rem:comp_kaehler}, the prefix ``semi'' carries a different meaning than in the Kähler and Hyperkähler case. To give a precise explanation, we need to introduce the notion of real structures on complexified pre-semi-Kähler manifolds:

\begin{definition}[Real structure]\label{def:comp_kaehler_real_str}
 Let $(X,\omega,J,I)$ be a complexified pre-semi-Kähler manifold. A \textbf{real structure} $\sigma$ on $(X,\omega,J,I)$ is a smooth involution on $X$ satisfying:
 \begin{enumerate}
  \item $\sigma$ preserves $\omega$, i.e., $\sigma^\ast\omega = \omega$,
  \item $\sigma$ is $J$-holomorphic, i.e., $J\circ d\sigma = d\sigma\circ J$,
  \item $\sigma$ is $I$-antiholomorphic, i.e., $I\circ d\sigma = -d\sigma\circ I$.
 \end{enumerate}
 The fixed point set $M\coloneqq\Fix\sigma$ is called \textbf{real form}. $M$ is \textbf{nice} if it meets every connected component of $X$.
\end{definition}

The name suggests that a real form of a complexified (pre-)semi-Kähler manifold is itself a (pre-)semi-Kähler manifold. The next proposition confirms this idea:

\begin{proposition}[Induced Kähler structure on real forms]\label{prop:induced_kaehler_on_real_forms}
 Let $(X,\omega,J,I)$ be a complexified (pre-)semi-Kähler manifold with real structure $\sigma$ and non-empty real form $M=\Fix\sigma$. Further, let $\iota:M\hookrightarrow X$ be the inclusion and define $\hat\omega\coloneqq\iota^\ast\omega$ as well as $\hat J\coloneqq J\vert_{TM}$. Then, $(M,\hat\omega,\hat J)$ is a (pre-)semi-Kähler manifold. If $(X,\omega,J,I)$ is holomorphic, $(M,\hat\omega,\hat J)$ is real-analytic.
\end{proposition}

\begin{proof}
 We begin by proving that $(M,\hat\omega,\hat J)$ is a pre-semi-Kähler manifold. As in \autoref{app:real_structures}, one can show that $M\subset X$ is a smooth submanifold with dimension $\dim_\R X = 2\dim_\R M$ and decomposition $T_pX = E_1\oplus E_{-1} = T_pM\oplus E_{-1}$ for every $p\in M$, where $E_{\pm 1}$ is the eigenspace of $d\sigma_p$ with respect to the eigenvalue $\pm 1$. We now need to prove that $\hat\omega$ is non-degenerate. For this, take a point $p\in M$ and a vector $v\in T_pM = E_1$. Since $\omega$ is non-degenerate, there exists a vector $w\in T_pX$ such that $\omega_p (v,w)\neq 0$. Because of $T_pX = E_1\oplus E_{-1}$, we can write $w = w_1 + w_{-1}$, where $w_\lambda\in E_\lambda$. Now observe that $\omega_p (v,w_{-1})$ vanishes, since:
 \begin{gather*}
  \omega_p (v,w_{-1}) = (\sigma^\ast\omega)_p (v,w_{-1}) = \omega_p (d\sigma_p v, d\sigma_p w_{-1}) = -\omega_p (v,w_{-1}).
 \end{gather*}
 This implies:
 \begin{gather*}
  \hat\omega_p (v,w_1) = \omega_p (v,w_1) = \omega_p (v,w)\neq 0
 \end{gather*}
 proving the non-degeneracy of $\hat\omega$.\\
 The next step is to show that $\hat J$ is a well-defined almost complex structure on $M$. It suffices to show that for every $p\in M$ and $v\in T_pM$ one has $J_pv\in T_pM$. This is an immediate consequence of the commutativity of $d\sigma_p$ and $J_p$ and the fact that $T_pM$ and $E_1$ coincide:
 \begin{gather*}
  d\sigma_p J_p v = J_p d\sigma_p v = J_p v\quad\Rightarrow\quad J_pv\in E_1 = T_pM.
 \end{gather*}
 The compatibility of $\hat\omega$ and $\hat J$ directly follows from the compatibility of $\omega$ and $J$, completing the proof of $(M,\hat\omega,\hat J)$ being pre-semi-Kähler.\\
 If $(X,\omega,J)$ is semi-Kähler, $\omega$ is closed and $J$ is integrable. We deduce from this that $\hat\omega$ is closed and $\hat J$ is integrable. Hence, $(M,\hat\omega,\hat J)$ is semi-Kähler in this case.\\
 Lastly, we note that the tensors $\Omega\coloneqq\omega - i\omega (I\cdot,\cdot)$ and $J$ satisfy $\overline{\sigma^\ast\Omega} = \Omega$ and $\overline{\sigma^\ast J} = J$. Thus, if $(X,\omega,J,I)$ is holomorphic, we can apply Proposition \autoref{prop:invariant_holo_tensors} to $\Omega$ and $J$ viewed as a holomorphic section of $\End (T^{(1,0)}_IX)$. The tensors induced by $\Omega$ and $J$ according to Proposition \autoref{prop:invariant_holo_tensors} are $\hat\omega$ and $\hat J$, respectively. In particular, Proposition \autoref{prop:invariant_holo_tensors} tells us that the induced tensors $\hat\omega$ and $\hat J$ are real-analytic concluding the proof.
\end{proof}

Now, we have all tools at hand to specify the meaning of ``semi'':

\begin{definition}\label{def:semi}
 Let $(X,\omega,J,I)$ be a complexified pre-semi-Kähler manifold. We drop the prefix ``semi'' if there exists a real structure $\sigma$ on $(X,\omega,J,I)$ with non-empty real form $M$ such that $(M,\hat\omega,\hat J)$ is pre-Kähler, i.e., the metric\linebreak $\hat g\coloneqq\hat\omega (\cdot,\hat J\cdot)$ is positive definite.
\end{definition}

Definition \autoref{def:semi} sheds some light on the chosen naming conventions: A complexified (pre-)Kähler manifold is, by definition, just the complexification of a (pre-)Kähler manifold. If the complexification is holomorphic, it is unique in the following sense:

\begin{lemma}[Uniqueness of holomorphic Kähler manifolds]\label{lem:uniqueness_of_holo_kaehler}
 Let $(X,\omega_1,J_1,I)$ and $(X,\omega_2,J_2,I)$ be two holomorphic pre-semi-Kähler manifolds with real structures $\sigma_1$ and $\sigma_2$ and non-empty real forms $M_1 = \Fix\sigma_1$ and $M_2=\Fix\sigma_2$, respectively. If $\sigma_1$ and $\sigma_2$ induce the same pre-semi-Kähler manifold, i.e., $(M_1,\hat\omega_1,\hat J_1)$ and $(M_2,\hat\omega_2,\hat J_2)$ coincide, then there exists an open neighborhood $U\subset X$ of $M_1 = M_2$ on which $(\omega_1,J_1,\sigma_1)$ and $(\omega_2,J_2,\sigma_2)$ agree. If, additionally, $M_1 = M_2$ is nice, one even has $\omega_1 = \omega_2$, $J_1 = J_2$, and $\sigma_1 = \sigma_2$.
\end{lemma}

\begin{proof}
 Lemma \autoref{lem:uniqueness_of_holo_kaehler} immediately follows from Proposition \autoref{prop:real_structure_unique} (applied to $\sigma_1$ and $\sigma_2$) and Lemma \autoref{lem:uniqueness_of_holomorphic_tensors} (applied to $\Omega_1$ and $\Omega_2$ as well as $J_1$ and $J_2$ viewed as holomorphic sections of $\End (T^{(1,0)}_I X)$).
\end{proof}

At this point, it is natural to ask the converse question: Does every (pre-\linebreak semi-)Kähler manifold admit a holomorphic complexification? If the Kähler manifold in question is real-analytic, the answer is positive:

\begin{lemma}[Complexification of Kähler manifolds]\label{lem:comp_of_kaehler}
 Let $(X,I)$ be a complex manifold with real structure $\sigma$ and non-empty real form $M\coloneqq\Fix\sigma$. Further, let $(M,\hat\omega,\hat J)$ be a real-analytic (pre-semi-)Kähler manifold. Then, there exists an open neighborhood $U\subset X$ of $M$ and a holomorphic (pre-semi-)Kähler manifold $(U,\omega,J,I)$ such that $\sigma$ is a real structure on $(U,\omega,J,I)$ and its induced Kähler structure is $(M,\hat\omega,\hat J)$. $(U,\omega,J,I)$ is unique in the sense of Lemma \autoref{lem:uniqueness_of_holo_kaehler}.
\end{lemma}

\begin{proof}
 Essentially, Lemma \autoref{lem:comp_of_kaehler} is a consequence of Lemma \autoref{lem:uniqueness_of_holomorphic_tensors} and \autoref{lem:holo_continuation_of_tensors}. By Lemma \autoref{lem:holo_continuation_of_tensors}, the real-analytic tensors $\hat\omega$ and $\hat J$ on $M$ possess holomorphic continuations $\Omega$ and $J^\prime$ on an open neighborhood $U\subset X$ of $M$. We take $\omega$ to be the real part of $\Omega$ and set:
 \begin{gather*}
  Jv\coloneqq\begin{cases} J^\prime v&\text{ for }v\in T^{(1,0)}_I U,\\ \overline{J^\prime \bar{v}}&\text{ for }v\in T^{(0,1)}_I U.\end{cases}
 \end{gather*}
 In order to show that $(U,\omega,J,I)$ is a complexified pre-semi-Kähler manifold, we need to check that $\omega$ is $I$-anticompatible and non-degenerate, that $J$ is an almost complex structure commuting with $I$, and that $\omega$ and $J$ are compatible. We begin with the $I$-anticompatibility and non-degeneracy of $\omega$:\\
 $\Omega$ is a holomorphic form and, thus, $I$-anticompatible. Accordingly, its real part $\omega$ is also $I$-anticompatible. Now recall that $\hat\omega$ is non-degenerate. Hence, $\Omega$ as its holomorphic continuation is non-degenerate on $T^{(1,0)}_IU$ after shrinking $U$ if necessary. Therefore, its real part $\omega$ is non-degenerate on $TU$.\\
 Let us now consider $J$. As a holomorphic tensor, $J^\prime$ commutes with $I$ and, thus, $J$ also commutes with $I$. Furthermore, one easily checks that $J$ is real, i.e., $Jv\in TU$ for every $v\in TU$. To check that $J$ is an almost complex structure, it suffices to show that $J^\prime$ is one, i.e., $(J^\prime)^2 = -\mathds{1}$. However, we already know that $(J^\prime)^2 = -\mathds{1}$ holds on $M$, since $J^\prime$ is the holomorphic continuation of an almost complex structure. Thus, by Lemma \autoref{lem:uniqueness_of_holomorphic_tensors}, $(J^\prime)^2 = -\mathds{1}$ must hold everywhere on $U$ after shrinking $U$ if necessary.\\
 Let us now show that $\omega$ and $J$ are compatible. As before, we apply Lemma \autoref{lem:uniqueness_of_holomorphic_tensors} to achieve that: The forms $\hat\omega$ and $\hat\omega (\hat J\cdot,\hat J\cdot)$ coincide, therefore, their holomorphic continuations $\Omega$ and $\Omega (J^\prime\cdot, J^\prime\cdot)$ agree on $M$. By Lemma \autoref{lem:uniqueness_of_holomorphic_tensors}, they must agree on all of $U$ after shrinking $U$ if necessary. The compatibility of $\omega$ and $J$ now follows from $\Omega = \Omega (J^\prime\cdot,J^\prime\cdot)$.\\
 So far, we have shown that $(U,\omega,J,I)$ is a holomorphic\footnote{By construction, the tensors $\Omega$ and $J^\prime = J\vert_{T^{(1,0)}_IU}$ are holomorphic. Therefore, $(U,\omega,J,I)$ is holomorphic.} pre-semi-Kähler manifold if $(M,\hat\omega,\hat J)$ is a pre-semi-Kähler manifold. The next step is to check that the same statement is true if we drop the prefix ``pre'' or ``semi''. For ``semi'', this is true by definition. For ``pre'', we need to check the integrability conditions of $(U,\omega,J,I)$. In particular, we have to verify that $\omega$ is closed and $J$ is integrable. This is the case if the exterior derivative of $\Omega$ and the Nijenhuis tensor of $J^\prime$ vanish. However, $d\Omega$ and $N_{J^\prime}$ are just the holomorphic continuations of $d\hat\omega$ and $N_{\hat J}$ which are zero. Therefore, by Lemma \autoref{lem:uniqueness_of_holomorphic_tensors},  $d\Omega$ and $N_{J^\prime}$ are zero finishing the proof.
\end{proof}

One useful application of Lemma \autoref{lem:uniqueness_of_holo_kaehler} and \autoref{lem:comp_of_kaehler} is the observation that it suffices to check the integrability conditions of a holomorphic pre-semi-Kähler manifold on its real form:

\begin{corollary}[Integrability conditions on real forms]
 Let $(X,\omega,J,I)$ be a holomorphic pre-semi-Kähler manifold with real structure $\sigma$ and nice real form $M = \Fix\sigma$. Then, the following statements are equivalent:
 \begin{enumerate}
  \item $(X,\omega,J)$ is semi-Kähler. 
  \item $(M,\hat\omega,\hat J)$ is semi-Kähler.
 \end{enumerate}
\end{corollary}

\begin{proof}
 The direction ``(i)$\Rightarrow$(ii)'' is trivial, so we only consider the converse direction. If $(M,\hat\omega,\hat J)$ is semi-Kähler, then, by Lemma \autoref{lem:comp_of_kaehler}, its complexification is semi-Kähler as well. By Lemma \autoref{lem:uniqueness_of_holo_kaehler}, this complexification must coincide with $(X,\omega,J,I)$ finishing the proof.
\end{proof}

Before we finish \autoref{app:kaehler}, we quickly want to examine how Hyperkähler and holomorphic Kähler manifolds are connected. The deep relation between Hyperkähler and holomorphic Kähler manifolds is most apparent in the symplectic picture (cf. Lemma \autoref{lem:symplectic_picture}). Indeed, semi-Hyperkähler and holomorphic semi-Kähler manifolds only differ by a sign in this picture:

\begin{lemma}[Hyperkähler vs. holomorphic Kähler]\label{lem:hyper_vs_holo}
 Let $(X,\omega,J)$ be a semi-Kähler manifold and $I$ be an integrable complex structure on $X$ satisfying $\omega (I\cdot,I\cdot) = -\omega$.
 \begin{enumerate}
  \item If $IJ = -JI$, then $(X,g,I,J,K)$ is a semi-Hyperkähler manifold where $g\coloneqq\omega (\cdot,J\cdot)$ and $K\coloneqq IJ$.
  \item If $IJ = JI$ and $J$ viewed as a section of $\End (T^{(1,0)}_IX)$ is holomorphic, then $(X,\omega,J,I)$ is a holomorphic semi-Kähler manifold.
 \end{enumerate}
\end{lemma}

\begin{proof}
 Lemma \autoref{lem:hyper_vs_holo} follows directly from Lemma \autoref{lem:symplectic_picture} and Definition \autoref{def:holo_kaehler}.
\end{proof}

 \chapter{Proof of Morse-Darboux Lemmata}
 \label{app:morse_darboux}
 Our goal in this part is to prove Lemma \autoref{lem:morse_darboux_lem_I} and \autoref{lem:morse_darboux_lem_II}. We start by showing Lemma \autoref{lem:morse_darboux_lem_I} for smooth manifolds:

\begin{lemma}[Morse-Darboux lemma I]\label{lem:morse_darboux_lem_I_smooth}
 Let $(M^2,\omega)$ be a smooth symplectic manifold, $L^1$ a smooth $1$-manifold, $f\in C^\infty (M,L)$, and $p\in M$ a non-degenerate critical point of $f$ with Morse index $\mu_f (p) = 0$. Further, let $T>0$ be a positive real number. Then, there exists a $C^1$-chart $\psi_L:U_L\to V_L\subset\R$ of $L$ near $f(p)$ which is smooth on $U_L\backslash\{f(p)\}$ such that all non-constant trajectories near $p$ of the RHS $(U_M, \omega\vert_{U_M}, H)$ with $U_M\coloneqq f^{-1}(U_L)$ and $H\coloneqq \psi_L\circ f\vert_{U_M}$ are $T$-periodic.
\end{lemma}

\begin{proof}
 The proof consists of three steps:
 \begin{enumerate}
  \item First, we convince ourselves that the non-constant trajectories $\gamma$ near $p$ are indeed periodic.
  \item Afterwards, we compute the period $\hat T(r)$ of a trajectory $\gamma$ near $p$ with $f\circ \gamma = r^2$ ($L=\R$, $f(p) = 0$) to show that $\hat T(r)$ is defined for $r\in (-\varepsilon,\varepsilon)$ ($\varepsilon>0$), depends smoothly on $r$, and is bounded from below by a positive constant.
  \item Lastly, we use these properties of $\hat T(r)$ to define a $C^1$-diffeomorphism\linebreak $\psi_L:U_L\to V_L\subset\R$ such that the trajectories of the rescaled RHS $(U_M,\omega\vert_{U_M}, H)$ with $U_M$ and $H$ as above have fixed period $T>0$.
 \end{enumerate}
 \vspace{0.3cm}
 \textbf{Step 1}\\\\
 Without loss of generality, we can assume, after choosing appropriate charts, that $L = \R$, $f(p) = 0\in\R$, and that the (usual) Morse index $\mu_f (p)$ of $f:M\to\R$ is $0$. Now, we apply the Morse lemma to find a chart $\hat\psi_M = (\hat x, \hat y):\hat U_M\to \hat V_M$ of $M$ near $p$ with $\hat \psi_M (p) = 0$ such that $f\vert_{\hat U_M} = {\hat x}^2 + {\hat y}^2$. In this chart, we have $\omega\vert_{\hat U_M} = \hat v\, d\hat x\wedge d\hat y$, where $\hat v\in C^\infty (\hat U_M,\R)$. Since $\omega$ is non-degenerate, we can assume $\hat v >0$ (after permuting $\hat x$ and $\hat y$ if necessary). Now consider the RHS $(M,\omega, f)$ and its trajectories $\gamma$ near $p$. $f$ is constant along $\gamma$, so for $r>0$ small enough the trajectory $\gamma$ near $p$ with $f\circ \gamma = r^2$ moves along the circle
 \begin{gather*}
  f^{-1}(r^2) = \hat \psi^{-1}_M (\{(\hat x,\hat y)\in\R^2\mid \hat x^2 + \hat y^2 = r^2\}) \cong S^1
 \end{gather*}
 with velocity $\dot\gamma\neq 0$. Hence, the trajectories near $p$ are periodic.\\\\
 \textbf{Step 2}\\\\
 Denote the period of $\gamma$ with $f\circ \gamma = r^2$ ($r>0$) by $\hat T(r)$. We calculate $\hat T(r)$ by going into polar coordinates:
 \begin{gather*}
  (\hat x,\hat y) = (r\cos (\varphi), r\sin (\varphi))\in\R^2.
 \end{gather*}
 Define $v \coloneqq \hat v\circ \hat\psi^{-1}_M$, then the Hamiltonian vector field $X_f$ is given by:
 \begin{gather*}
  (\hat\psi_M)_\ast X_f = \frac{2}{v(\hat x, \hat y)}
  \begin{pmatrix}
   -\hat y\\ \hat x
  \end{pmatrix} = \frac{2}{v(r\cos (\varphi), r\sin (\varphi))}
  \begin{pmatrix}
   -r\sin (\varphi)\\ r\cos (\varphi)
  \end{pmatrix}.
 \end{gather*}
 Now parameterize an integral curve $\gamma:\R\to M$ of $X_f$ by $r_\gamma, \varphi_\gamma:\R\to\R$ in the following way:
 \begin{gather*}
  \hat \psi_M\circ \gamma (t) =
  \begin{pmatrix}
   r_\gamma (t)\cos (\varphi_\gamma (t))\\
   r_\gamma (t)\sin (\varphi_\gamma (t))
  \end{pmatrix}.
 \end{gather*}
 The integral curve equation $\dot\gamma = X_f (\gamma)$ now yields:
 \begin{gather*}
  \dot r_\gamma = 0,\quad \dot \varphi_\gamma = \frac{2}{v(r_\gamma\cos (\varphi_\gamma), r_\gamma\sin (\varphi_\gamma))}.
 \end{gather*}
 Thus, $r_\gamma$ is constant and, since $f\circ\gamma = r^2$, given by $r$. This allows us to define $\Phi:\R\to \R$ by:
 \begin{align*}
  \Phi (t)&\coloneqq \varphi_\gamma (t) - \varphi_\gamma (0) = \int\limits^t_0 \dot\varphi_\gamma (t^\prime) dt^\prime\\
  &= \int\limits^t_0 \frac{2}{v(r\cos (\varphi_\gamma (t^\prime)), r\sin (\varphi_\gamma(t^\prime)))} dt^\prime.
 \end{align*}
 $\Phi$ is an orientation preserving diffeomorphism, since $\dot\Phi = \frac{2}{v}>0$. Hence, $\Phi^{-1}$ exists and is given by
 \begin{gather*}
  \Phi^{-1} (\alpha) = \frac{1}{2}\int\limits^{\varphi_0+\alpha}_{\varphi_0} v(r\cos (\varphi), r\sin (\varphi))\, d\varphi
 \end{gather*}
 with $\varphi_0\coloneqq \varphi_\gamma (0)$, as:
 \begin{align*}
  \frac{d\Phi^{-1}}{d\alpha} (\alpha) &\stackrel{\phantom{\varphi_\gamma = \varphi_0 + \Phi}}{=} \frac{1}{\dot\Phi (\Phi^{-1}(\alpha))} = \frac{1}{2}v (r\cos (\varphi_\gamma (\Phi^{-1}(\alpha))), r\sin (\varphi_\gamma (\Phi^{-1}(\alpha))))\\
  &\stackrel{\varphi_\gamma = \varphi_0 + \Phi}{=} \frac{1}{2} v(r\cos (\varphi_0 + \alpha), r\sin (\varphi_0 + \alpha)).
 \end{align*}
 We can now use $\Phi$ and $\Phi^{-1}$ to compute $\hat T(r)$:
 \begin{align*}
  &\varphi_\gamma (t + \hat T(r)) = \varphi_\gamma (t) + 2\pi\quad\forall t\in\R\quad \Rightarrow \Phi (\hat T(r)) = 2\pi\\
  \Rightarrow &\hat T(r) = \Phi^{-1} (2\pi) = \frac{1}{2}\int\limits^{\varphi_0+2\pi}_{\varphi_0} v(r\cos (\varphi), r\sin (\varphi))\, d\varphi\\
  &\phantom{\hat T(r)}\, = \frac{1}{2}\int\limits^{2\pi}_{0} v(r\cos (\varphi), r\sin (\varphi))\, d\varphi.
 \end{align*}
 As we can see, $\hat T(r)$ depends smoothly on $r>0$. In fact, this formula allows us to define $\hat T(r)$ smoothly for $r\leq 0$ as well. It turns out that the function $\hat T(r)$ is even:
 \begin{align*}
  \hat T(-r) &= \frac{1}{2}\int\limits^{2\pi}_{0} v(-r\cos (\varphi), -r\sin (\varphi))\, d\varphi = \frac{1}{2}\int\limits^{2\pi}_{0} v(r\cos (\varphi+\pi), r\sin (\varphi+\pi))\, d\varphi\\
  &= \frac{1}{2}\int\limits^{2\pi}_{0} v(r\cos (\varphi), r\sin (\varphi))\, d\varphi = \hat T(r).
 \end{align*}
 After shrinking $\hat U_M$ if necessary, we can assume that $\hat U_M$ has a compact neighborhood. Thus, $\hat v\in C^\infty (\hat U_M,\R)$ is bounded from below by a positive constant $v_{\min}>0$. Therefore, $v\coloneqq \hat v\circ\hat\psi^{-1}_M$ is also bounded from below by $v_{\min}$. This implies:
 \begin{gather*}
  \hat T(r) = \frac{1}{2}\int\limits^{2\pi}_{0} v(r\cos (\varphi), r\sin (\varphi)) \geq \pi v_{\min} > 0\quad\forall r.
 \end{gather*}
 \textbf{Step 3}\\\\
 Lastly, we use these properties of $\hat T(r)$ to define a $C^1$-diffeomorphism\linebreak $\psi_L:U_L\to V_L\subset\R$ on a neighborhood $U_L\subset L$ of $0\in\R = L$ which rescales the periods of the trajectories $\gamma$ to a fixed period $T > 0$. Again, consider a trajectory $\gamma$ near $p$ with $f\circ \gamma = r^2$. We want to define $\psi_L (s)$ with $s = r^2$ in such a way that the rescaled trajectory $\Gamma (t)\coloneqq \gamma (\hat T(r)\cdot t/T)$ is a trajectory of the rescaled RHS $(U_M,\omega\vert_{U_M},H)$ ($U_M$ and $H$ as above). Hence, we want $\Gamma$ to be an integral curve of $X_H$:
 \begin{align*}
  \dot\Gamma (t) &= \frac{\hat T(r)}{T}\dot\gamma \left(\frac{\hat T(r)}{T}t\right) = \frac{\hat T(r)}{T} X_f (\Gamma (t))\\
  &\stackrel{!}{=} X_H (\Gamma (t)) = \frac{d\psi_L}{ds} (f\circ\Gamma (t)) X_f (\Gamma (t)) = \frac{d\psi_L}{ds} (r^2) X_f (\Gamma (t))
 \end{align*}
 Thus, we obtain the following condition for $\psi_L$:
 \begin{gather}
  \frac{d\psi_L}{ds} (r^2) = \frac{\hat T(r)}{T}.\label{eq:psi_L}
 \end{gather}
 To solve Equation \eqref{eq:psi_L}, we define the function $g:\R\to\R$ by $g(s)\coloneqq \sqrt{|s|}$. $g$ is continuous on $\R$ and smooth on $\R\backslash\{0\}$. Thus, $\hat T\circ g$ is also continuous on a neighborhood $U_L$ of $0$ and smooth on $U_L\backslash\{0\}$. Now we define $\psi_L$ by:
 \begin{gather*}
  \psi_L (s)\coloneqq \frac{1}{T}\int\limits^s_0 \hat T\circ g (s^\prime) ds^\prime = \frac{1}{T} \int\limits^s_0 \hat T (\sqrt{|s^\prime|}) ds^\prime.
 \end{gather*}
 Since $\hat T\circ g$ is continuous, $\psi_L$ exists, is a $C^1$-function on $U_L$, and smooth on $U_L\backslash\{0\}$. Furthermore:
 \begin{gather*}
  \frac{d\psi_L}{ds} (s) = \frac{\hat T (\sqrt{|s|})}{T}\geq \frac{\pi v_{\min}}{T} >0.
 \end{gather*}
 Therefore, $\psi_L$ is also a $C^1$-diffeomorphism. The last equation together with the fact that $\hat T (r)$ is even shows Equation \eqref{eq:psi_L} concluding the proof:
 \begin{gather*}
  \frac{d\psi_L}{ds} (r^2) = \frac{\hat T (|r|)}{T} = \frac{\hat T(r)}{T}.
 \end{gather*}
\end{proof}

\begin{remark}[No regularity issues in a real-analytic setup]\label{rem:no_reg_issue}
 If all objects in Lemma \autoref{lem:morse_darboux_lem_I_smooth} are real-analytic instead of smooth (cf. Lemma \autoref{lem:morse_darboux_lem_I}), then the regularity issues do not occur, i.e., the chart $\psi_L$ can be chosen to be real-analytic. We can see this as follows: By similar arguments as before, the map $\hat T:\R\to\R$ assigning the period $\hat T(r)$ to each radius $r$ is given by
 \begin{gather*}
  \hat T(r) = \frac{1}{2}\int\limits^{2\pi}_{0} v(r\cos (\varphi), r\sin (\varphi))\, d\varphi
 \end{gather*}
 and, hence, real-analytic, as $v$ is real-analytic. Thus, $\hat T$ can be written as a power series near $r = 0$:
 \begin{gather*}
  \hat T (r) = \sum^\infty_{k=0} a_k r^k.
 \end{gather*}
 Now recall that $\hat T$ is even, therefore, only even powers occur in the power series of $\hat T$:
 \begin{gather*}
  \hat T (r) = \sum^\infty_{k=0} a_{2k} r^{2k}.
 \end{gather*}
 This allows us to define the real-analytic function $\hat t$ by:
 \begin{gather*}
  \hat t (s) \coloneqq \sum^\infty_{k=0} a_{2k} s^k.
 \end{gather*}
 Obviously, $\hat T$ and $\hat t$ satisfy the relation: $\hat t (r^2) = \hat T (r)$. Now we define the real-analytic chart $\psi_L$ by:
 \begin{gather*}
  \psi_L (s)\coloneqq \frac{1}{T}\int\limits^s_0 \hat t (s^\prime) ds^\prime.
 \end{gather*}
 As in the proof of Lemma \autoref{lem:morse_darboux_lem_I_smooth}, all non-constant trajectories of the Hamiltonian $\psi_L\circ f$ are $T$-periodic, since $\psi_L$ fulfills the equation
 \begin{gather*}
  \frac{d\psi_L}{ds} (r^2) = \frac{\hat t (r^2)}{T} = \frac{\hat T (r)}{T}.
 \end{gather*}
\end{remark}

Next, we prove Lemma \autoref{lem:morse_darboux_lem_II} in the smooth case:

\begin{lemma}[Morse-Darboux lemma II]\label{lem:morse_darboux_lem_II_smooth}
 Let $(M^2,\omega)$ be a smooth symplectic manifold and let $H\in C^\infty (M,\R)$ be a smooth function on $M$ with non-degenerate critical point $p\in M$ of Morse index $\mu_H (p)\neq 1$. Further, let $T>0$ be a positive real number. Then, the following statements are equivalent:
 \begin{enumerate}
  \item There exists a $C^0$-chart $\psi_M = (x,y):U_M\to V_M\subset\R^2$ of $M$ near $p$ which is smooth on $U_M\backslash\{p\}$ such that ($\psi_M (p) = 0$):
  \begin{enumerate}[label = (\alph*)]
   \item $H\vert_{U_M} = H(p) \pm \frac{\pi}{T}(x^2 + y^2)$,
   \item $\omega\vert_{U_M} = dy\wedge dx$.
  \end{enumerate}
  \item There exists an open neighborhood $U_M\subset M$ of $p$ such that all non-constant trajectories of the RHS $(U_M, \omega\vert_{U_M}, H\vert_{U_M})$ are $T$-periodic.
  \item There exists a number $E_0 > 0$ such that $\int_{U(E)} \omega = T\cdot E$ for every number $E\in [0, E_0]$, where $U(E)$ is the connected component containing $p$ of the set $\{q\in M\mid |H(q)-H(p)|\leq E\}$.
 \end{enumerate}
\end{lemma}

\begin{proof}
 The idea of the proof is simple: The implications ``(i)$\Rightarrow$(ii)'' and ``(ii)$\Rightarrow$(iii)'' follow from straightforward computations. To show the remaining implication ``(iii)$\Rightarrow$(i)'', we first choose a Morse chart $\hat \psi_M = (\hat x, \hat y):\hat U_M\to \hat V_M$ such that $H\vert_{\hat U_M} = H(p) + \varepsilon\frac{\pi}{T} (\hat x^2 + \hat y^2)$, where $\varepsilon\in\{-1,+1\}$. In general, $\hat \psi_M$ is not a Darboux chart for $\omega$. Still, the trajectories of the RHS $(\hat U_M, \omega\vert_{\hat U_M}, H\vert_{\hat U_M})$ are circles, in particular periodic orbits. Solely the angular velocity of these circles might not be constant. To rectify this, we go into polar coordinates $(r,\varphi)$ and apply an appropriately chosen diffeomorphism to $\varphi$. This operation results in a new chart $\psi_M$. Since we have not changed the radius $r$, $\psi_M$ is still a Morse chart. However, the change of the angle coordinate turns $\psi_M$ into a Darboux chart.\\\\
 $\boxed{\text{(i)}\Rightarrow\text{(ii)}}$\\\\
 In the chart $\psi_M$, we find for the Hamiltonian vector field $X_H$:
 \begin{gather*}
  (\psi_M)_\ast X_H = \pm\frac{2\pi}{T}\begin{pmatrix}y\\ -x\end{pmatrix}.
 \end{gather*}
 Hence, the integral curves $\gamma$ of $X_H$ are given by:
 \begin{gather*}
  \psi_M\circ\gamma (t) = \begin{pmatrix}r_0 \cos (\varphi_0 \mp \frac{2\pi}{T}t)\\ r_0 \sin (\varphi_0 \mp \frac{2\pi}{T}t)\end{pmatrix}.
 \end{gather*}
 Thus, the trajectories $\gamma$ near $p$ are $T$-periodic.
 
 \pagebreak
 
 $\boxed{\text{(ii)}\Rightarrow\text{(iii)}}$\\\\
 $p$ is a non-degenerate critical point of $H$ with Morse index $\mu_H (p) \neq 1$, hence, we can find a Morse chart $\hat \psi_M = (\hat x,\hat y):\hat U_M\to \hat V_M$ of $M$ near $p$ such that ($\hat\psi_M (p) = 0$):
 \begin{gather*}
  H\vert_{\hat U_M} = H(p) + \varepsilon\frac{\pi}{T}(\hat x^2 + \hat y^2)\quad\text{with}\quad\varepsilon\in\{-1,+1\}.
 \end{gather*}
 In this chart, we have\footnote{We write $d\hat x\wedge d\hat y$ instead of $d\hat y\wedge d\hat x$ here to match the convention from the proof of Lemma \autoref{lem:morse_darboux_lem_I_smooth} and to ensure that $\omega$ agrees with the standard orientation of $\R^2$.} $\omega\vert_{\hat U_M} = \hat v\cdot d\hat x\wedge d\hat y$ for $\hat v\in C^\infty (\hat U_M,\R)$. After permuting $\hat x$ and $\hat y$ if necessary, we can assume that $\hat v >0$. Let $\gamma$ be an integral curve of $X_H$. In polar coordinates $(\hat x, \hat y) = (r\cos (\varphi), r\sin (\varphi))$, we can parameterize $\gamma$ via $r_\gamma, \varphi_\gamma:\R\to\R$ as follows:
 \begin{gather*}
  \hat\psi_M\circ\gamma (t) = \begin{pmatrix}
                r_\gamma (t) \cos (\varphi_\gamma (t))\\
                r_\gamma (t) \sin (\varphi_\gamma (t))
               \end{pmatrix}.
 \end{gather*}
 With $v\coloneqq \hat v\circ \hat\psi^{-1}_M$, the integral curve equation becomes:
 \begin{gather*}
  \dot r_\gamma = 0,\quad \dot \varphi_\gamma = \frac{2\pi\varepsilon}{Tv(r_\gamma\cos (\varphi_\gamma), r_\gamma\sin (\varphi_\gamma))}.
 \end{gather*}
 Thus, $r_\gamma\equiv r$ is constant. This allows us to define $\Phi:\R\to \R$ by:
 \begin{align*}
  \Phi (t)&\coloneqq \varphi_\gamma (t) - \varphi_\gamma (0) = \int\limits^t_0 \dot\varphi_\gamma (t^\prime) dt^\prime\\
  &= \int\limits^t_0 \frac{2\pi\varepsilon}{T v(r\cos (\varphi_\gamma (t^\prime)), r\sin (\varphi_\gamma(t^\prime)))} dt^\prime.
 \end{align*}
 As in the proof of Lemma \autoref{lem:morse_darboux_lem_I_smooth}, $\Phi^{-1}$ exists and is given by:
 \begin{gather*}
  \Phi^{-1} (\alpha) = \frac{T}{2\pi\varepsilon}\int\limits^{\varphi_0+\alpha}_{\varphi_0} v(r\cos (\varphi), r\sin (\varphi))\, d\varphi
 \end{gather*}
 with $\varphi_0\coloneqq \varphi_\gamma (0)$. We now use the fact that, by assumption, $\gamma$ is $T$-periodic, so $\Phi^{-1}$ satisfies:
 \begin{gather*}
  \Phi^{-1} (2\pi) = \varepsilon T\quad\Rightarrow \int\limits^{2\pi}_0 v(r\cos (\varphi), r\sin (\varphi)) d\varphi = 2\pi.
 \end{gather*}
 Observe that the last equation holds for all $r>0$ small enough. 
 It allows us to compute the symplectic area $\int_{U(E)}\omega$:
 \begin{align*}
  \int\limits_{U(E)} \omega &= \int\limits^{2\pi}_0\int\limits^{\sqrt{\frac{TE}{\pi}}}_0  v(r\cos(\varphi), r\sin(\varphi))\, rdr\, d\varphi\\
  &= \int\limits^{\sqrt{\frac{TE}{\pi}}}_0 \left(\int\limits^{2\pi}_0 v(r\cos(\varphi), r\sin(\varphi))\, d\varphi\right) rdr\\
  &= \int\limits^{\sqrt{\frac{TE}{\pi}}}_0 2\pi r\, dr = T\cdot E.
 \end{align*}
 $\boxed{\text{(iii)}\Rightarrow\text{(i)}}$\\\\
 As in ``$\text{(ii)}\Rightarrow\text{(iii)}$'', we can find a Morse chart $\hat \psi_M = (\hat x,\hat y):\hat U_M\to \hat V_M$ of $M$ near $p$ such that ($\hat\psi_M (p) = 0$):
 \begin{gather*}
  H\vert_{\hat U_M} = H(p) + \varepsilon\frac{\pi}{T}(\hat x^2 + \hat y^2),
 \end{gather*}
 where $\varepsilon\in\{-1,+1\}$ and $\omega\vert_{\hat U_M} = \hat v\cdot d\hat x\wedge d\hat y$ for $\hat v\in C^\infty (\hat U_M,\R)$ with $\hat v >0$. By taking the derivative of $\int_{U(E)}\omega = T\cdot E$ with respect to $E$, we deduce that
 \begin{gather*}
  \int\limits^{2\pi}_0 v(r\cos (\varphi), r\sin (\varphi))\, d\varphi = 2\pi
 \end{gather*}
 holds for $r>0$ small enough and $v\coloneqq \hat v\circ \hat\psi^{-1}_M$. The last equation implies that the map $P:(0,\varepsilon_0)\times S^1\to (0,\varepsilon_0)\times S^1$ given by
 \begin{gather*}
  P(r, [\varphi])\coloneqq \left(r, \left[\int\limits^{\varphi}_0 v(r\cos (\varphi^\prime), r\sin (\varphi^\prime))\, d\varphi^\prime\right]\right)
 \end{gather*}
 is well-defined for $\varepsilon_0 >0$ small enough. $P$ is a smooth diffeomorphism, since\linebreak the determinant of $dP$ is $v>0$. Denote the polar coordinates by\linebreak $S:\R_+ \times S^1\to \R^2\backslash\{0\}$, i.e., $S(r, [\varphi])\coloneqq (r\cos (\varphi), r\sin (\varphi))$, and consider the map $S\circ P\circ S^{-1}:\mathring D_{\varepsilon_0}\to \mathring D_{\varepsilon_0}$, where $\mathring D_{\varepsilon_0}\coloneqq \{x\in\R^2\backslash\{0\}\mid ||x|| < \varepsilon_0\}$.\linebreak $S\circ P\circ S^{-1}$ is a smooth diffeomorphism, since both $S$ and $P$ are smooth diffeomorphisms. Furthermore, $S\circ P\circ S^{-1}$ maps circles of radius $r$ to circles of radius $r$, hence, we can extend $S\circ P\circ S^{-1}$ to a homeomorphism on $D_{\varepsilon_0}\coloneqq\{x\in\R^2\mid ||x|| < \varepsilon_0\}$ by setting $S\circ P\circ S^{-1} (0)\coloneqq 0$.\\
 Now consider the map $(x,y)\equiv \psi_M\coloneqq S\circ P\circ S^{-1}\circ\hat \psi_M:U_M\to V_M$. $\psi_M$ is a $C^0$-chart of $M$ near $p$ and smooth on $U_M\backslash\{p\}$. Recall that $\hat \psi_M$ is a Morse chart for $H$, thus, the level sets of $H$ are circles in the chart $\hat \psi_M$. Since the charts $\psi_M$ and $\hat \psi_M$ only differ by postcomposition with $S\circ P\circ S^{-1}$ which preserves circles, $\psi_M$ is also a Morse chart for $H$:
 \begin{gather*}
  H\vert_{U_M} = H(p) + \varepsilon\frac{\pi}{T}(x^2 + y^2).
 \end{gather*}
 Furthermore, the fact that the Jacobian of $P$ is $v$ implies:
 \begin{gather*}
  \omega\vert_{U_M} = dx\wedge dy.
 \end{gather*}
 We conclude the proof by permuting $x$ and $y$ to bring $\omega$ into the form $dy\wedge dx$.
\end{proof}

\begin{remark}[No regularity issues in a real-analytic setup]\label{rem:no_reg_issue_II}
 As for Lemma \autoref{lem:morse_darboux_lem_I_smooth}, the chart $\psi_M:U_M\to V_M$ is real-analytic on all of $U_M$ if all objects in Lemma \autoref{lem:morse_darboux_lem_II_smooth} are chosen to be real-analytic (cf. Lemma \autoref{lem:morse_darboux_lem_II}). To prove this, it suffices to show that the map $S\circ P\circ S^{-1}:D_{\varepsilon_0}\to D_{\varepsilon_0}$ is a real-analytic diffeomorphism. To do so, we employ the notations from above and define:
 \begin{gather*}
  \hat\varphi\coloneqq \int\limits^{\varphi}_0 \left(v(r\cos (\varphi^\prime), r\sin (\varphi^\prime))-1\right)\, d\varphi^\prime = \int\limits^{\varphi}_0 v(r\cos (\varphi^\prime), r\sin (\varphi^\prime))\, d\varphi^\prime - \varphi.
 \end{gather*}
 Using polar coordinates $(\hat x, \hat y) = (r\cos (\varphi),r\sin (\varphi))\in\mathring D_{\varepsilon_0}$, we obtain:
 \begin{align*}
  S\circ P\circ S^{-1}(\hat x,\hat y) &= \begin{pmatrix}r\cos (\hat\varphi + \varphi)\\ r\sin (\hat\varphi + \varphi)\end{pmatrix} = \begin{pmatrix}\cos (\hat\varphi) & -\sin (\hat\varphi)\\ \sin (\hat\varphi) & \cos (\hat\varphi)\end{pmatrix} \begin{pmatrix}r\cos (\varphi)\\ r\sin (\varphi)\end{pmatrix}\\
  &= \begin{pmatrix}\cos (\hat\varphi) & -\sin (\hat\varphi)\\ \sin (\hat\varphi) & \cos (\hat\varphi)\end{pmatrix}\begin{pmatrix}\hat x\\ \hat y\end{pmatrix},
 \end{align*}
 where we used the angle addition theorems for sine and cosine. We see that $S\circ P\circ S^{-1}$ is real-analytic on $D_{\varepsilon_0}$ if $\hat\varphi$ can be expressed as a real-analytic function on $D_{\varepsilon_0}$. We show the analyticity of $\hat\varphi$ by first recalling that $v$ satisfies:
 \begin{gather}
  \int\limits^{2\pi}_0 v(r\cos(\varphi), r\sin(\varphi))\, d\varphi = 2\pi.\label{eq:v_int_eq}
 \end{gather}
 In the case of $r=0$, the last equation becomes $v(0,0) = 1$. Now note that $v(\hat x,\hat y)$ is real-analytic in $\hat x$ and $\hat y$, hence, we can write ($v_{k_1k_2}\in\R$):
 \begin{gather*}
  v(\hat x,\hat y) - 1 = \sum^\infty_{k_1 + k_2 >0} v_{k_1k_2} \hat x^{k_1}\hat y^{k_2} = \sum^\infty_{n = 1} r^n\sum^n_{k=0} v_{k(n-k)}\cos^k (\varphi)\sin^{n-k}(\varphi).
 \end{gather*}
 Inserting the last equation into the definition of $\hat\varphi$ gives:
 \begin{gather*}
  \hat\varphi = \sum^\infty_{n = 1} r^n\sum^n_{k=0} v_{k(n-k)}\int\limits^{\varphi}_0\cos^k (\varphi^\prime)\sin^{n-k}(\varphi^\prime)\, d\varphi^\prime.
 \end{gather*}
 Next, we want to show that there are coefficients\footnote{One can even show that the coefficients $v^\prime_{nk}$ are unique!} $v^\prime_{nk}$ such that:
 \begin{gather}
  \sum^n_{k=0} v_{k(n-k)}\int\limits^{\varphi}_0\cos^k (\varphi^\prime)\sin^{n-k}(\varphi^\prime)\, d\varphi^\prime = \sum^n_{k=0}v^\prime_{nk}\cos^k(\varphi)\sin^{n-k}(\varphi).\label{eq:v_coeff_eq}
 \end{gather}
 We observe that $\hat\varphi (\varphi)$ is a $2\pi$-periodic $C^1$-function due to \autoref{eq:v_int_eq}. Therefore, $\frac{1}{n!}\pa{r}^n\hat\varphi\vert_{r=0} (\varphi)$ is also a $2\pi$-periodic $C^1$-function. In particular, its Fourier series exists and converges uniformly to $\frac{1}{n!}\pa{r}^n\hat\varphi\vert_{r=0} (\varphi)$. This allows us to write:
 \begin{gather*}
  \frac{1}{n!}\pa{r}^n\hat\varphi\vert_{r=0} (\varphi) = \sum^n_{k=0} v_{k(n-k)}\int\limits^{\varphi}_0\cos^k (\varphi^\prime)\sin^{n-k}(\varphi^\prime)\, d\varphi^\prime = \sum_{m\in\Z}a_m e^{im\varphi},
 \end{gather*}
 where $a_m\in\C$ are complex numbers. It is easy to check using\linebreak $2\cos (\varphi) = e^{i\varphi} + e^{-i\varphi}$ and $2i\sin (\varphi) = e^{i\varphi} - e^{-i\varphi}$ that $a_m$ vanishes for $|m|>n$. Let us now take the derivative of the last equation with respect to $\varphi$:
 \begin{gather*}
  \sum^n_{k=0} v_{k(n-k)}\cos^k(\varphi)\sin^{n-k}(\varphi) = \sum_{|m|\leq n} ima_m e^{im\varphi}.
 \end{gather*}
 We notice that the left-hand side changes by the sign $(-1)^n$ under the transformation $\varphi\mapsto\varphi + \pi$. On the other hand, only terms of the right-hand side with odd $m$ pick up a minus sign under this transformation, while the terms with $m$ even are invariant. Thus, we find\footnote{We always have $a_0 = 0$, since the integral in $\frac{1}{n!}\pa{r}^n\hat\varphi\vert_{r=0} (\varphi)$ does not generate constant terms.} $a_m = 0$ for $m \equiv n+1\mod 2$. Plugging in\linebreak $e^{im\varphi} = (\cos(\varphi) + i\sgn(m)\sin(\varphi))^{|m|}$ into the formula for $\frac{1}{n!}\pa{r}^n\hat\varphi\vert_{r=0} (\varphi)$ now yields ($a^\prime_{mk}\in\C$):
 \begin{gather*}
  \frac{1}{n!}\pa{r}^n\hat\varphi\vert_{r=0} (\varphi) = \sum^n_{\substack{m=0\\ n-m\,\text{even}}}\sum^m_{k=0}a^\prime_{mk} \cos^k(\varphi)\sin^{m-k}(\varphi).
 \end{gather*}
 We now multiply each term in the last equation by\nolinebreak $(\cos^2 (\varphi) + \sin^2 (\varphi))^{(n-m)/2} = 1$ to show \autoref{eq:v_coeff_eq}.\\
 Let us return to $\hat\varphi$. \autoref{eq:v_coeff_eq} allows us to write:
 \begin{align*}
  \hat\varphi (\hat x,\hat y) &= \sum^\infty_{n = 1} r^n\sum^n_{k=0}v^\prime_{nk}\cos^k(\varphi)\sin^{n-k}(\varphi)\\
  &= \sum^\infty_{n = 1} \sum^n_{k=0}v^\prime_{nk}(r\cos(\varphi))^k(r\sin(\varphi))^{n-k}\\
  &= \sum^\infty_{n = 1} \sum^n_{k=0}v^\prime_{nk}\hat x^k\hat y^{n-k}.
 \end{align*}
 The last equation shows that $\hat\varphi$ is a real-analytic function on a neighborhood of $(0,0)$ concluding the proof.
\end{remark}

 \chapter[Various Action Functionals for HHSs and PHHSs]{Various Action Functionals for HHSs and PHHSs\chaptermark{Various Action Functionals}}
 \chaptermark{Various Action Functionals}
 \label{app:various_action_functionals}
 In \autoref{sec:HHS} and \autoref{sec:PHHS}, we have defined and studied action functionals for HHSs and PHHSs. The ``critical points'' of these action functionals gave us (pseudo-)holomorphic trajectories of the system under consideration whose domains are parallelograms in the complex plane. However, these action functionals are not the only functionals whose ``critical points'' can be linked to (pseudo-)holomorphic trajectories. There is, in fact, an abundance of action functionals that differ in the domain of their trajectories and the way the ``one-dimensional'' action functionals of their underlying RHSs are integrated. In this part of the appendix, we present and examine a large selection of such action functionals. First, we only formulate and explore action functionals for HHSs. Afterwards, we explain how these action functionals need to be modified in order to give action functionals for PHHSs. Hereby, we realize that the presented action functionals for PHHSs are all real-valued. From this point of view, a Floer-like theory for PHHSs revolving around these real-valued functionals might be possible.\\
We begin by defining an action functional for holomorphic trajectories whose domains are disks $D^{z_0}_R\subset\mathbb{C}$ of radius $R>0$ centered at $z_0\in\mathbb{C}$. To do that, we first need to partition the disk $D^{z_0}_R$ into lines. We choose the partition consisting of lines starting at the center $z_0$ and ending at any boundary point $z\in\partial D^{z_0}_R$. For every such radial line, we consider the action functional $\mathcal{A}^{\Lambda}_{e^{i\alpha}\mH}$ from Remark \autoref{rem:tilted_traj}. We now obtain an action functional for HHSs by integrating the action $\mathcal{A}^{\Lambda}_{e^{i\alpha}\mH}$ over all radial lines, i.e., $\alpha\in [0,2\pi]$:

\begin{proposition}[Action functional $\actdiski$]\label{prop:actdiski}
 Let $(X,\Omega = d\Lambda, \mH)$ be an exact\linebreak HHS, $D^{z_0}_R\coloneqq\{z\in\mathbb{C}\mid |z-z_0|\leq R\}$ be a disk of radius $R>0$ centered at $z_0\in\mathbb{C}$,\linebreak $\mathcal{P}_{D^{z_0}_R}\coloneqq C^\infty (D^{z_0}_R, X)$ be the set of smooth maps from $D^{z_0}_R$ to $X$, and\linebreak $\actdiski:\mathcal{P}_{D^{z_0}_R}\to\mathbb{C}$ be the action functional defined by
 \begin{gather*}
  \actdiski [\gamma]\coloneqq \frac{1}{2\pi}\int\limits^{2\pi}_0\int\limits^R_0\left[\Lambda_{\gamma_\alpha (r)}\left(\frac{d\gamma_\alpha}{dr}(r)\right) - e^{i\alpha}\cdot \mH\circ\gamma_\alpha (r)\right]dr\, d\alpha\quad\forall \gamma \in\mathcal{P}_{D^{z_0}_R},
 \end{gather*}
 where $\gamma_\alpha:[0,R]\to X$ is defined by $\gamma_\alpha (r)\coloneqq \gamma (z_0 + re^{i\alpha})$. Now let $\gamma\in\mathcal{P}_{D^{z_0}_R}$. Then, $\gamma$ is a holomorphic trajectory of the HHS $(X,\Omega,\mH)$ iff $\gamma$ is a ``critical point'' of $\actdiski$. Here, ``critical points'' means that we only allow for those variations of $\gamma$ which keep $\gamma$ fixed at the boundary $\partial D^{z_0}_R$ and the \underline{center} $z_0$.
\end{proposition}

\begin{proof}
 Take the notations from above. Using Remark \autoref{rem:tilted_traj} and writing $\actdiski$ as
 \begin{gather*}
  \actdiski [\gamma] = \frac{1}{2\pi}\int\limits^{2\pi}_0\mathcal{A}^{\Lambda}_{e^{i\alpha}\mH} [\gamma_\alpha] d\alpha,
 \end{gather*}
 we can show as in the proof of Lemma \autoref{lem:holo_action_prin} that $\gamma$ is a ``critical point'' of $\actdiski$ iff\linebreak $\gamma_\alpha:[0,R]\to X$ is a (real) integral curve of $\cos (\alpha)\cdot X^R_\mH + \sin (\alpha)\cdot J(X^R_\mH)$ for every $\alpha\in [0,2\pi]$, where $X_\mH = 1/2 (X^R_\mH - iJ(X^R_\mH))$ is the Hamiltonian vector field of $(X,\Omega,\mH)$. Thus, a ``critical point'' $\gamma$ is uniquely determined, given an initial value $x_0\coloneqq\gamma (z_0)$, by:
 \begin{gather*}
   \gamma(z_0 + re^{i\alpha}) = \varphi^{\cos (\alpha)\cdot X^R_\mH + \sin (\alpha)\cdot J(X^R_\mH)}_r (x_0) = \varphi^{r\cos (\alpha)\cdot X^R_\mH + r\sin (\alpha)\cdot J(X^R_\mH)}_1 (x_0),
 \end{gather*}
 where $\varphi^V_t$ is the time-$t$-flow of a real vector field $V$ on $X$. Comparing the last equation with the formula for the holomorphic trajectory $\gamma^{z_0,x_0}$ satisfying $\gamma^{z_0,x_0} (z_0) = x_0$ given in the proof of Proposition \autoref{prop:holo_traj} shows that $\gamma$ is a holomorphic trajectory of the HHS $(X,\Omega,\mH)$ iff $\gamma$ is a ``critical point'' of $\actdiski$.\\
 Lastly, we have to explain why a ``critical point'' $\gamma$ of $\actdiski$ needs to fix the variations of $\gamma$ at the boundary $\partial D^{z_0}_R$ and the center $z_0$. Recall the variation of $\mathcal{A}^{\Lambda}_{e^{i\alpha}\mH}$ at $\gamma_\alpha$. In general, the variation of this functional also includes terms associated with the boundary of the image of $\gamma_\alpha$. This boundary consists of two points, namely the center $z_0$ and one boundary point $z\in\partial D^{z_0}_R$. To get rid of these boundary terms in the variation of $\actdiski$, we have to keep $\gamma$ fixed at $z_0$ and $\partial D^{z_0}_R$.
\end{proof}

In \autoref{sec:HHS}, we have given two reasons why we need to vary over all smooth curves $\gamma$ and cannot simple restrict the variational problem to the space of holomorphic curves $\gamma$. The new-found action functional offers an additional perspective on that matter. It maps every holomorphic curve $\gamma$ to zero, hence, only varying it over the space of holomorphic curves is meaningless:

\begin{proposition}\label{prop:actdiski_value}
 Take the assumptions and notations from Proposition \autoref{prop:actdiski}. Further, let $\gamma:D^{z_0}_R\to X$ be any holomorphic map from $D^{z_0}_R$ to $X$. Then:
 \begin{gather*}
  \actdiski [\gamma] = 0.
 \end{gather*}
\end{proposition}

\begin{proof}
 Take the assumptions and notations from Proposition \autoref{prop:actdiski} and let\linebreak $\gamma:D^{z_0}_R\to X$ be holomorphic. Using the relation $\Lambda\circ J = i\cdot \Lambda$ for holomorphic $1$-forms, we find:
 \begin{gather*}
  \Lambda_{\gamma_\alpha(r)}\left(\frac{d\gamma_\alpha}{dr}(r)\right) = e^{i\alpha}\cdot\Lambda_{\gamma_\alpha (r)}\left( \gamma^{\prime}(z_0 + re^{i\alpha})\right),
 \end{gather*}
 where $\gamma^{\prime}$ is the complex derivative of $\gamma$. With this, we obtain:
 \begin{align*}
  \actdiski [\gamma] &= \frac{1}{2\pi}\int\limits^{2\pi}_0\int\limits^R_0\left[\Lambda_{\gamma_\alpha (r)}\left(\frac{d\gamma_\alpha}{dr}(r)\right) - e^{i\alpha}\cdot \mH\circ\gamma_\alpha (r)\right]dr\, d\alpha\\
  &= \frac{1}{2\pi}\int\limits^{2\pi}_0\int\limits^R_0\left[\Lambda_{\gamma (z_0 + re^{i\alpha})}\left(\gamma^\prime (z_0 + re^{i\alpha})\right) - \mH\circ\gamma (z_0 + re^{i\alpha})\right]dr\, e^{i\alpha}\, d\alpha\\
  &= \int\limits^R_0 \frac{1}{2\pi i}\oint\limits_{|z-z_0| = r}\left[\Lambda_{\gamma (z)}\left(\gamma^\prime (z)\right) - \mH\circ\gamma (z)\right]dz\, \frac{dr}{r}\\
  &= \int\limits^R_0 \sum_{p\in D^{z_0}_r}\text{\normalfont Res}_{z = p}\left[\Lambda_{\gamma (z)}\left(\gamma^\prime (z)\right) - \mH\circ\gamma (z)\right] \frac{dr}{r},
 \end{align*}
 where $\text{\normalfont Res}_{z = p} [f(z)]$ denotes the residue of a meromorphic function $f(z)$ at\linebreak $z = p$. In the last line of the computation, we have used Cauchy's theorem. Clearly, the function $f$ has no residues inside $D^{z_0}_r$ ($r\in [0,R]$) in our case. Thus, the action vanishes for holomorphic curves concluding the proof.
\end{proof}

Even though the ``critical values'' of $\actdiski$ are nice and easy to understand, the action functional itself does not appear to be particularly useful. Often, we want to modify action functionals such that trajectories become actual critical points. The standard ways to achieve this are to either put the boundary of the trajectory on exact Lagrangian submanifolds or to impose periodicity. Both methods do not appear to be meaningful here. For the presented action functional, periodicity means periodicity of the radial lines. Thus, a ``periodic'' curve\linebreak $\gamma:D^{z_0}_R\to X$ needs to attain the same value on its boundary as on its center. However, the only holomorphic maps $\gamma:D^{z_0}_R\to X$ exhibiting such a behavior are constant curves by the identity theorem.\\
The other method, mapping the ``boundary'' to exact Lagrangian submanifolds, takes an unnatural and downright ugly form here, namely mapping $z_0$ and $\partial D^{z_0}_R$ to exact Lagrangian submanifolds. The action functional $\actdiskii$ improves on $\actdiski$ in that regard. To avoid boundary terms associated with $z_0$, which are at the center\footnote{Cum grano salis.} of our problem, we now partition the disk $D^{z_0}_R$ into lines starting at $z_0 - z$ and ending at $z_0 + z$ ($|z| = R$). In order to account for the doubled length of the radial lines, we only integrate over the angles $\alpha\in [0,\pi]$ this time:

\begin{proposition}[Action functional $\actdiskii$]\label{prop:actdiskii}
 Let $(X,\Omega = d\Lambda, \mH)$ be an exact\linebreak HHS, $D^{z_0}_R\coloneqq\{z\in\mathbb{C}\mid |z-z_0|\leq R\}$ be a disk of radius $R>0$ centered at $z_0\in\mathbb{C}$,\linebreak $\mathcal{P}_{D^{z_0}_R}\coloneqq C^\infty (D^{z_0}_R, X)$ be the set of smooth maps from $D^{z_0}_R$ to $X$, and\linebreak $\actdiskii:\mathcal{P}_{D^{z_0}_R}\to\mathbb{C}$ be the action functional defined by
 \begin{gather*}
  \actdiskii [\gamma]\coloneqq \frac{i}{4R}\int\limits^{\pi}_0\int\limits^R_{-R}\left[\Lambda_{\gamma_\alpha (r)}\left(\frac{d\gamma_\alpha}{dr}(r)\right) - e^{i\alpha}\cdot \mH\circ\gamma_\alpha (r)\right]dr\, d\alpha\quad\forall \gamma \in\mathcal{P}_{D^{z_0}_R},
 \end{gather*}
 where $\gamma_\alpha:[-R,R]\to X$ is defined by $\gamma_\alpha (r)\coloneqq \gamma (z_0 + re^{i\alpha})$. Now let $\gamma\in\mathcal{P}_{D^{z_0}_R}$. Then, $\gamma$ is a holomorphic trajectory of the HHS $(X,\Omega,\mH)$ iff $\gamma$ is a ``critical point'' of $\actdiskii$. Here, ``critical points'' means that we only allow for those variations of $\gamma$ which keep $\gamma$ fixed at the boundary $\partial D^{z_0}_R$.
\end{proposition}

\begin{proof}
 The proof works as the proof of Proposition \autoref{prop:actdiski} by writing $\actdiskii$ as
 \begin{gather*}
  \actdiskii [\gamma] = \frac{i}{4R}\int\limits^{\pi}_0\mathcal{A}^{\Lambda}_{e^{i\alpha}\mH} [\gamma_\alpha] d\alpha.
 \end{gather*}
 Here, the variations of $\gamma$ only need to keep $\gamma$ fixed at the boundary $\partial D^{z_0}_R$, since the radial lines start and end at $\partial D^{z_0}_R$.
\end{proof}

\begin{remark}[No Proposition \autoref{prop:actdiski_value} for $\actdiskii$]\label{rem:actdiskii_value}
 Proposition \autoref{prop:actdiski_value} does not apply to $\actdiskii$. In fact, the normalization in Proposition \autoref{prop:actdiskii} is chosen such that the action of constant curves is given by the Hamilton function:
 \begin{gather*}
  \actdiskii [\gamma_{x_0}] = \mH (x_0),
 \end{gather*}
 where $\gamma_{x_0} (z)\coloneqq x_0\in X$ for every $z\in D^{z_0}_{R}$. Thus, any singular point $x_0$ of $\mH$ with $\mH (x_0)\neq 0$ provides a counterexample to Proposition \autoref{prop:actdiski_value} for $\actdiskii$.
\end{remark}

If we modify $\actdiskii$ such that the holomorphic trajectories become actual critical points, we see that this action is a bit more reasonable. In the Lagrangian case, we now restrict the space of smooth curves $\gamma:D^{z_0}_R\to X$ to the space of those curves which only map the boundary $\partial D^{z_0}_R$ to exact Lagrangian submanifolds, as one would expect. However, the modification via periodicity still only gives trivial results. One can see this as follows: Now, periodicity means periodicity of radial lines starting and ending at $\partial D^{z_0}_R$. In this sense, we say $\gamma:D^{z_0}_R\to X$ is ``periodic'' if it assigns the same value to opposite points on the boundary $\partial D^{z_0}_R$. For the sake of simplicity, let us now assume $z_0 = 0$. For such a ``periodic'' $\gamma$, define $\gamma_-$ by $\gamma_- (z)\coloneqq \gamma (-z)$. If $\gamma$ is holomorphic, then $\gamma_-$ is also holomorphic and, by assumption, attains the same values on $\partial D^{0}_R$ as $\gamma$. Hence, by the identity theorem, $\gamma$ and $\gamma_-$ denote the same map. However, if $\gamma$ is a trajectory, then $\gamma$ is an integral curve of the Hamiltonian vector field $X_\mH$ and we have:
\begin{gather*}
 X_\mH (\gamma (z)) = \frac{d}{dz}\gamma (z) = \frac{d}{dz} \gamma (-z) = - X_\mH (\gamma (-z)) = -X_\mH (\gamma (z)).
\end{gather*}
Thus, the Hamiltonian vector field vanishes in this case and $\gamma$ is a constant curve.\\
We cannot only formulate $\actdiski$ and $\actdiskii$ for disks $D^{z_0}_R$, but for any bounded star-shaped domain\footnote{Here, a domain $D\subset\mathbb{C}$ is a path-connected set with non-empty interior $D^\circ$ dense in $D$.} $D\subset\mathbb{C}$ with smooth boundary\footnote{The boundary $b$ is parameterized by the polar angle $\alpha$, i.e,  $b(\alpha) = z_0 + R(\alpha)e^{i\alpha}\in \partial D$.} $b:\mathbb{R}/2\pi\mathbb{Z}\to \partial D$:

\begin{proposition}[Action functionals $\acti$ and $\actii$ for bounded star-shaped domains] \label{prop:action_starshaped}
 Let $(X,\Omega = d\Lambda, \mH)$ be an exact HHS, let $D\subset\mathbb{C}$ be a bounded domain in $\mathbb{C}$ which is star-shaped with respect to $z_0$ and has smooth boundary $b:\mathbb{R}/2\pi\mathbb{Z}\to\partial D$, and let $\mathcal{P}_{D}\coloneqq C^\infty (D, X)$ be the set of smooth maps from $D$ to $X$. Then, we can define the action functionals $\acti:\mathcal{P}_{D}\to\mathbb{C}$ and $\actii:\mathcal{P}_D\to \mathbb{C}$ by
 \begin{align*}
  \acti [\gamma]&\coloneqq \frac{1}{2\pi}\int\limits^{2\pi}_0\int\limits^{R(\alpha)}_0\left[\Lambda_{\gamma_\alpha (r)}\left(\frac{d\gamma_\alpha}{dr}(r)\right) - e^{i\alpha}\cdot \mH\circ\gamma_\alpha (r)\right]dr\, d\alpha,\\
  \actii [\gamma]&\coloneqq \frac{i}{4\hat{R}}\int\limits^{\pi}_0\int\limits^{R(\alpha)}_{-R(\alpha - \pi)}\left[\Lambda_{\gamma_\alpha (r)}\left(\frac{d\gamma_\alpha}{dr}(r)\right) - e^{i\alpha}\cdot \mH\circ\gamma_\alpha (r)\right]dr\, d\alpha,
 \end{align*}
 where $\gamma\in\mathcal{P}_D$, $R:\mathbb{R}/2\pi\mathbb{Z}\to\mathbb{R}$ is defined by $R(\alpha)\coloneqq |b(\alpha)-z_0|$,\linebreak $\gamma_\alpha:[-R(\alpha - \pi), R(\alpha)]\to X$ is given by $\gamma_\alpha (r)\coloneqq \gamma(z_0 + re^{i\alpha})$, and $\hat{R}$ is defined by
 \begin{gather*}
  \hat{R}\coloneqq \frac{i}{4}\left[\int\limits^{2\pi}_{\pi} R(\alpha)e^{i\alpha} d\alpha - \int\limits^{\pi}_{0} R(\alpha)e^{i\alpha} d\alpha\right].
 \end{gather*}
 Now let $\gamma\in\mathcal{P}_{D}$. Then, $\gamma$ is a holomorphic trajectory of the HHS $(X,\Omega,\mH)$ iff $\gamma$ is a ``critical point''\footnote{In the sense of Proposition \autoref{prop:actdiski}.} of $\acti$ iff $\gamma$ is a ``critical point''\footnote{In the sense of Proposition \autoref{prop:actdiskii}.} of $\actii$.
\end{proposition}

\begin{proof}
 Confer the proofs of Proposition \autoref{prop:actdiski} and \autoref{prop:actdiskii}.
\end{proof}

\begin{remark}[Normalization of $\acti$ and $\actii$]\label{rem:normalization}
 The normalization of $\acti$ and $\actii$ are chosen such that they coincide with our previous definitions for $D = D^{z_0}_R$ being a disk. In particular, $\actii$ agrees with the Hamilton function $\mH$ for constant curves $\gamma$.
\end{remark}

One might wonder how the action functionals $\acti$ and $\actii$ are related to the action functional $\mathcal{A}^{P_\alpha}_\mH$ for parallelograms $P_\alpha$ from \autoref{sec:HHS} and \autoref{sec:PHHS}, especially because a parallelogram $P_\alpha$ is also a bounded star-shaped domain. When we modify the functionals $\acti$ and $\actii$ to describe general PHHSs, we will see that $\acti$ and $\actii$ differ a lot from $\mathcal{A}^{P_\alpha}_\mH$. To compare $\mathcal{A}^{P_\alpha}_\mH$ directly with $\acti$ and $\actii$, let us express $\acti$ and $\actii$ in the same coordinates as $\mathcal{A}^{P_\alpha}_\mH$, namely Cartesian coordinates $z = t + is$:

\begin{proposition}[$\acti$ and $\actii$ in Cartesian coordinates]\label{prop:cartesian_coordinates}
 Employ the assumptions and notations from Proposition \autoref{prop:action_starshaped}. For $\gamma\in\mathcal{P}_D$, we define the following derivatives in Cartesian coordinates $z = t + is \in D$:
 \begin{gather*}
  \frac{\partial\gamma}{\partial z} (z)\coloneqq \frac{1}{2}\left(\frac{\partial\gamma}{\partial t} (z) - i\frac{\partial\gamma}{\partial s} (z)\right),\quad \frac{\partial\gamma}{\partial \bar{z}} (z)\coloneqq \frac{1}{2}\left(\frac{\partial\gamma}{\partial t} (z) + i\frac{\partial\gamma}{\partial s} (z)\right).
 \end{gather*}
 Furthermore, define the complex functions $f,g:D\to\mathbb{C}$ by:
 \begin{gather*}
  f(z)\coloneqq \Lambda_{\gamma (z)}\left(\frac{\partial\gamma}{\partial z}(z)\right) - \mH\circ\gamma (z),\quad g(z)\coloneqq \Lambda_{\gamma (z)}\left(\frac{\partial\gamma}{\partial \bar{z}}(z)\right).
 \end{gather*}
 Then, the action functionals $\acti$ and $\actii$ in Cartesian coordinates are given by:
 \begin{align*}
  \acti [\gamma] &= \frac{1}{2\pi}\iint\limits_{D}\left[\frac{f(z)}{\bar{z} - \bar{z_0}} + \frac{g(z)}{z - z_0}\right]dt\wedge ds,\\
  \actii [\gamma] &= \frac{i}{4\hat{R}}\iint\limits_{D^+}\left[\frac{f(z)}{\bar{z} - \bar{z_0}} + \frac{g(z)}{z - z_0}\right]dt\wedge ds - \frac{i}{4\hat{R}}\iint\limits_{D^-}\left[\frac{f(z)}{\bar{z} - \bar{z_0}} + \frac{g(z)}{z - z_0}\right]dt\wedge ds,
 \end{align*}
 where $z = t+is\in D$, $\bar{\cdot}$ denotes the complex conjugation,\linebreak $D^+\coloneqq \{z\in D\mid \text{\normalfont Im}(z-z_0)\geq 0\}$, and $D^-\coloneqq \{z\in D\mid \text{\normalfont Im}(z-z_0)\leq 0\}$.
\end{proposition}

\begin{proof}
 Take the assumptions and notations from above. We only show Proposition \autoref{prop:cartesian_coordinates} for $D = D^0_R\equiv D_R$ being a disk of radius $R>0$ centered at the origin. The general case can be shown similarly. Using the derivatives defined above, we can write:
 \begin{gather*}
  \frac{\partial\gamma}{\partial t} (z) = \frac{\partial\gamma}{\partial z} (z) + \frac{\partial\gamma}{\partial \bar{z}} (z),\quad \frac{\partial\gamma}{\partial s} (z) = i\left(\frac{\partial\gamma}{\partial z} (z) - \frac{\partial\gamma}{\partial \bar{z}} (z)\right).
 \end{gather*}
 Now consider the map $\gamma_\alpha:[-R,R]\to X$ defined in polar coordinates $z = re^{i\alpha}$ by $\gamma_\alpha (r) = \gamma (re^{i\alpha})$ for every $\alpha\in [0,2\pi]$. $\Lambda$ applied to the derivative of $\gamma_\alpha$ gives:
 \begin{alignat*}{2}
  \Lambda_{\gamma_\alpha (r)}\left(\frac{d\gamma_\alpha}{dr} (r)\right) &= \cos (\alpha)\cdot \Lambda_{\gamma (re^{i\alpha})}\left(\frac{\partial\gamma}{\partial t} (re^{i\alpha})\right) &&+ \sin (\alpha)\cdot \Lambda_{\gamma (re^{i\alpha})}\left(\frac{\partial\gamma}{\partial s} (re^{i\alpha})\right)\\
  &= e^{i\alpha}\cdot \Lambda_{\gamma (re^{i\alpha})}\left(\frac{\partial\gamma}{\partial z} (re^{i\alpha})\right) &&+ e^{-i\alpha}\cdot \Lambda_{\gamma (re^{i\alpha})}\left(\frac{\partial\gamma}{\partial \bar{z}} (re^{i\alpha})\right).
 \end{alignat*}
 Recalling the definition of $f$ and $g$, this allows us to write:
 \begin{align*}
  \Lambda_{\gamma_\alpha (r)}\left(\frac{d\gamma_\alpha}{dr}(r)\right) - e^{i\alpha}\cdot \mH\circ\gamma_\alpha (r) &= e^{i\alpha}\cdot f(re^{i\alpha}) + e^{-i\alpha}\cdot g(re^{i\alpha})\\
  &= r\cdot\left[\frac{f(re^{i\alpha})}{re^{-i\alpha}} + \frac{g(re^{i\alpha})}{re^{i\alpha}}\right] = r\cdot\left[\frac{f(z)}{\bar{z}} + \frac{g(z)}{z}\right].
 \end{align*}
 The expression for $\acti$ is now obtained by inserting the last equation into the defining formula for $\acti$ and using $r\cdot dr\wedge d\alpha = dt\wedge ds$ for $re^{i\alpha} = z = t + is$. Observing that the integrands of $\acti$ and $\actii$ agree on $D^+_R$ and differ by a sign on $D^-_R$ concludes the proof.
\end{proof}


The form of $\acti$ and $\actii$ in Cartesian coordinates is rather remarkable. In fact, we can express the action functional $\mathcal{A}^{P_\alpha}_\mH$ from \autoref{sec:HHS} in the same form as $\acti$ for $D = P_\alpha$, just with different complex functions $f$ and $g$, namely by setting:
\begin{gather*}
 f(z) = 2\pi(\bar{z} - \bar{z_0})\cdot\left[\Lambda_{R,\gamma (z)}\left(2\frac{\partial\gamma}{\partial z}(z)\right) - \mH\circ\gamma (z)\right],\quad g(z) = 0.
\end{gather*}
The similarities between $\acti$ and $\mathcal{A}^{P_\alpha}_\mH$ become even more apparent if we evaluate $\mathcal{A}^{P_\alpha}_\mH$ at holomorphic curves $\gamma:P_\alpha\to X$. As in Remark \autoref{rem:several_remarks}, Point 3, we find in this case:
\begin{gather*}
 \mathcal{A}^{P_\alpha}_\mH[\gamma] =  \iint\limits_{P_\alpha}\left[\Lambda\left(\frac{\partial\gamma}{\partial z}(z)\right) - \mH\circ\gamma (z)\right] dt\wedge ds.
\end{gather*}
Because the function $g$ as defined in Proposition \autoref{prop:cartesian_coordinates} vanishes for holomorphic $\gamma$, the only difference between $\acti$ and $\mathcal{A}^{P_\alpha}_\mH$ is now the factor $2\pi(\bar{z} - \bar{z_0})$ in the integrand. However, this seemingly small difference is rather impactful, as we will shortly see.\\
Lastly, we want to modify the actions $\acti$ and $\actii$ in such a way that they also apply to PHHSs. Recall that for an exact PHHS $(X,J;\Omega_R = d\Lambda_R,\mH_R)$ the induced $2$-form $\Omega_I$ is, in general, not closed. Hence, only the parts of $\acti$ and $\actii$ that do not include $\Lambda_I$ are well-defined for PHHSs. Precisely speaking, these are the real part of $\acti$ and the real part of $-4i\hat R\cdot \actii$. Still, these \underline{real-valued} functionals satisfy an action principle with respect to the pseudo-holomorphic trajectories of a PHHS. In fact, this can be shown in the same way as Proposition \autoref{prop:actdiski} by simply observing that Remark \autoref{rem:tilted_traj}, which is crucial for the proof of Proposition \autoref{prop:actdiski}, is also valid for the real part of the functional $\mathcal{A}^\Lambda_{e^{i\alpha}\mH}$.\\
The generalization of $\acti$ and $\actii$ to PHHSs has now revealed the most striking difference between $\acti$ and $\mathcal{A}^{P_\alpha}_\mH$: While the action principle related to $\acti$ still applies if we only consider the real part of $\acti$, both the real \underline{and} imaginary part of $\mathcal{A}^{P_\alpha}_\mH$ are crucial for the validity of the action principle related to $\mathcal{A}^{P_\alpha}_\mH$. In particular, there might exist a Floer-like theory related to the real part of $\acti$ or $-4i\hat R\cdot \actii$. For $\mathcal{A}^{P_\alpha}_\mH$, we have no intuition on how such a theory should look like, since we do not know how to interpret a complex function as a Morse function in the sense of Morse homology. Since we can turn $\acti$ and $\actii$ into real-valued action functionals, the same objections do not apply to them. Nevertheless, the question remains whether the real part of $\acti$ or $-4i\hat R\cdot \actii$ are indeed suitable Morse functions and whether the resulting Floer theories, if they exist, give any non-trivial result. At least for the (conjectured) Hamiltonian\footnote{In Hamiltonian Floer theory, one only considers periodic orbits.} Floer theory related to $\acti$ and $\actii$, it is most likely that it only gives trivial results, since all trajectories, which are ``periodic'' in a sense suitable for $\acti$ and $\actii$ as explained above, are automatically constant.

 \chapter[Almost Complex Structures on $\mathbf{TM}$ and $\mathbf{T^\ast M}$]{Almost Complex Structures on $\mathbf{TM}$ and $\mathbf{T^\ast M}$\chaptermark{Almost Complex Structures}}
 \chaptermark{Almost Complex Structures}
 \label{app:almost_complex_structures}
 In this part, we explain how a connection $\nabla$ on a manifold $M$ induces an almost complex structure\footnote{In general, $J_\nabla$ is \underline{not} the complex structure adapted to $\nabla$, even though both share a deep relation (cf. \autoref{sec:duality}, in particular Example \autoref{ex:flat_adapted}).} $J_\nabla$ on its tangent bundle $TM$ (cf. \cite{Dombrowski1962} and \cite{Tachibana1961}). If $\nabla = \nabla^g$ is the Levi-Civita connection of a semi-Riemannian metric $g$ on $M$, then $\nabla^g$ also defines an almost complex structure $J^\ast_{\nabla^g}$ on the cotangent bundle $T^\ast M$ via the bundle isomorphism $G:TM\to T^\ast M, v\mapsto \iota_v g$. In this case, $J^\ast_{\nabla^g}$ is compatible with the canonical symplectic form $\omega_{\can}$ on $T^\ast M$ in the sense that $\omega_{\can} (\cdot, J^\ast_{\nabla^g}\cdot)$ is a semi-Riemannian metric on $T^\ast M$ of signature $(2s, 2t)$, where $(s,t)$ is the signature of $g$. Furthermore, we will see that $J_{\nabla^g}$ or, equivalently, $J^\ast_{\nabla^g}$ is integrable if and only if $g$ is flat.\\
Let $M$ a smooth manifold of dimension $n$, $\pi_E:E\to M$ be a (smooth) vector bundle, and $\nabla$ be a linear connection\footnote{Sometimes, the term ``affine connection'' is used instead of ``linear connection''.} on the vector bundle $E\to M$. The (fiberwise) kernel of the differential $d\pi_E: TE\to TM$ yields the vertical subbundle $VE$ of $TE$, while $\nabla$ defines the horizontal subbundle $HE$ of $TE$. In fact, the notion of a horizontal subbundle $HE$ is equivalent to the notion of a linear connection $\nabla$ on $E\to M$. To see this, let $K:TE = VE\oplus HE\to E$ be the vertical projection\footnote{
 Technically speaking, the map $\hat K:TE = VE\oplus HE\to VE$ is the vertical projection. To obtain $K$ from $\hat K$, we have already exploited the fact that $E\to M$ is a vector bundle allowing us to identify the fibers of $VE$ with the fibers of $E$ via the linear isomorphism
 \begin{gather*}
  E_p\to V_wE, v\mapsto \left.\frac{d}{dt}\right\vert_{t = 0} (w + vt)
 \end{gather*}
 for $p\in M$ and $w\in E_p = \pi_E^{-1}(p)$.
} (cf. \cite{Dombrowski1962} or \cite{Eliasson1967} for the construction of $K$). The data $HE$ and $K$ are equivalent, since, given the horizontal subbundle $HE$, we can always define the vertical projection $K$ and, given the map $K$, we can always define the horizontal bundle $HE$ to be the (fiberwise) kernel of $K$. Likewise, the data $\nabla$ and $K$ are equivalent. Their relation is encoded in the following formula:
\begin{gather*}
 \nabla_X Y = K\circ dY (X)\quad\forall X\in TM,
\end{gather*}
where the section $Y\in\Gamma (E)$ is viewed as a smooth map $Y:M\to E$.

\pagebreak

Now observe that the map $f^\nabla \coloneqq (d\pi_E, K):TE\to TM\oplus E$ is a bundle map over\linebreak $\pi_E:E\to M$, i.e., the diagram
\begin{center}
 \begin{tikzcd}
  TE \arrow[r, "f^\nabla"] \arrow[d]
  & TM\oplus E \arrow[d] \\
  E \arrow[r, "\pi_E"]
  & M
 \end{tikzcd}
\end{center}
commutes, where the vertical arrows are the base point projections of the vector bundles $TE$ and $TM\oplus E$. Fiberwise, the bundle map $f^\nabla$ is a linear isomorphism. Thus, $TE$ is isomorphic to the pullback bundle $\pi_E^\ast (TM\oplus E)$.\\
Now take $E\to M$ to be the tangent bundle $TM\to M$ with base point projection $\pi\equiv \pi_{TM}:TM\to M$. Then, $f^\nabla$ allows us to define the almost complex structure $J_\nabla$ on the manifold $TM$ by ``pulling back'' the almost complex structure\linebreak $J_{TM\oplus TM}:TM\oplus TM \to TM\oplus TM, (w_1, w_2)\mapsto (w_2, -w_1)$ of the vector bundle\linebreak $TM\oplus TM\to M$ to the vector bundle $T(TM)\to TM$. Explicitly speaking, the almost complex structure $J_\nabla$ is completely determined by the following equations:
\begin{gather*}
 d\pi\circ J_\nabla = K,\quad K\circ J_\nabla = -d\pi.
\end{gather*}
Next, we wish to express $J_\nabla$ in local coordinates. For this, we need to parameterize the vertical and horizontal subspaces first. Choose a point $p\in M$ and normal coordinates $\psi = (x_1,\ldots, x_n):U\to V\subset\R^n$ of $(M,\nabla)$ near $p$, i.e., a chart $\psi$ in which all lines through the origin $\psi (p) = 0$ are geodesics\footnote{Specifically, $\psi$ is given by $l^{-1}\circ\exp_p^{-1}$, where $l:\R^n\to T_pM$ is a linear isomorphism and $\exp_p:T_pM\to M$ is the exponential map near $p$ associated with $\nabla$. The exponential map $\exp_p$ is defined by geodesics through $p$ meaning $\exp_p (v)\coloneqq \gamma_{p,v}(1)$ with $\gamma_{p,v}:[0,1]\to M$ being the unique map that fulfills $\nabla_{\dot\gamma_{p,v}}\dot\gamma_{p,v} = 0$, $\gamma_{p,v} (0) = p$, and $\dot\gamma_{p,v} (0) = v$. Restricted to small neighborhoods of $0\in T_pM$ and $\exp_p (0) = p\in M$, $\exp_p$ becomes a well-defined diffeomorphism.}. $\psi$ induces coordinates of $TM$ near any point $w\in T_pM$ which we denote by $T\psi = (\hat x_1,\ldots, \hat x_n, v_1,\ldots, v_n)$:
\begin{gather*}
 (T\psi)^{-1} (\hat x_1,\ldots, \hat x_n, v_1,\ldots, v_n)\coloneqq \sum^n_{k = 1}v_k\partial_{x_k, \psi^{-1}(\hat x_1,\ldots, \hat x_n)}.
\end{gather*}
We find for the coordinate vector fields $\pa{v_k}$:
\begin{gather*}
 d\pi (\partial_{v_k}) = \left.\frac{d}{dt}\right\vert_{t = 0} \pi\left(\sum_{l\neq k}v_l\partial_{x_l} + (v_k + t)\partial_{x_k}\right) = \left.\frac{d}{dt}\right\vert_{t = 0} \pi\left(\sum^n_{l = 1}v_l\partial_{x_l}\right) = 0.
\end{gather*}
Thus, the vector fields $\partial_{v_1},\ldots, \partial_{v_n}$ span the vertical subspaces.\\
To parameterize the horizontal subspaces, we consider the local vector field\linebreak $X_c\coloneqq \sum_k c_k\partial_{x_k}$ with constants $c_k\in\R$. We find:
\begin{align*}
 \nabla_{\partial_{x_k}} X_c (p) &= K\circ dX_{c,p} (\partial_{x_k, p}) = K\left(\left.\frac{d}{dt}\right\vert_{t = 0}\left(X_c\circ \psi^{-1}(\psi (p) + t\hat e_k)\right)\right)\\
 &= K(\partial_{\hat x_k, X_c (p)}).
\end{align*}
Hence, the vectors $\partial_{\hat x_1, w},\ldots, \partial_{\hat x_n, w}$ span the horizontal subspaces\linebreak $H(TM) = \ker (K)$ for any $w = X_c (p)\in T_pM$ if and only if the equation
\begin{gather*}
 \nabla_{\partial_{x_i}}\partial_{x_j} (p) = 0
\end{gather*}
holds for all $i,j\in\{1,\ldots, n\}$. The last equation is satisfied for all normal coordinates $(x_1,\ldots, x_n)$ near any point $p\in M$ if and only if $\nabla$ is symmetric\footnote{This is due to the fact that the geodesic equation $\nabla_{\dot\gamma}\dot\gamma = 0$ only ``sees'' the symmetric part of $\nabla$.}, i.e., satisfies:
\begin{gather*}
 \nabla_X Y - \nabla_Y X = [X,Y]\quad \forall X,Y\in\Gamma (TM).
\end{gather*}
Thus, we shall henceforth assume that the connection $\nabla$ is symmetric.\\
So far, we have found that the vertical subspaces at $w\in T_pM$ are spanned by the vectors $\partial_{v_1,w},\ldots, \partial_{v_n,w}$, while the horizontal subspaces at $w\in T_pM$ are spanned by the vectors $\partial_{\hat x_1,w},\ldots, \partial_{\hat x_n,w}$ for normal coordinates $\psi = (x_1,\ldots, x_n)$ of $(M,\nabla)$ near $p\in M$ with $T\psi = (\hat x_1,\ldots, \hat x_n, v_1,\ldots, v_n)$. If we express an arbitrary vector $u\in T_w (TM)$ and its image $J_\nabla(u)\in T_w (TM)$ as
\begin{gather*}
 u = \sum^n_{k = 1}a_k\partial_{v_k,w} + b_k\partial_{\hat x_k,w}\quad\text{and}\quad J_\nabla(u) = \sum^n_{k = 1}c_k\partial_{v_k,w} + d_k\partial_{\hat x_k,w},
\end{gather*}
we can compute the coefficients $c_k$ and $d_k$ in terms of $a_k$ and $b_k$:
\begin{alignat*}{2}
 -\sum^n_{k = 1}b_k\partial_{x_k,p} &= -d\pi (u) = K\circ J_\nabla (u) = \sum^n_{k = 1}c_k\partial_{x_k,p}\quad &\Rightarrow c_k = -b_k\\
 \sum^n_{k = 1} a_k\partial_{x_k,p} &= K(u) = d\pi\circ J_\nabla (u) = \sum^n_{k =1}d_k\partial_{x_k,p}\quad &\Rightarrow d_k = a_k,
\end{alignat*}
where we used $d\pi (\partial_{\hat x_k,w}) = \partial_{x_k,p} = K (\partial_{v_k,w})$, $V(TM) = \ker (d\pi)$, and\linebreak $H(TM) = \ker (K)$. This gives us:
\begin{gather*}
 J_\nabla (\partial_{v_k,w}) = \partial_{\hat x_k,w},\quad J_\nabla (\partial_{\hat x_k,w}) = -\partial_{v_k,w}.
\end{gather*}
We now see that the almost complex structure $J_\nabla$ assumes the standard form in normal coordinates near $p\in M$ for points $w\in T_pM$. This does not mean, however, that $J_\nabla$ is integrable, since the last equation is not necessarily true for all points $w$ within the chart domain $TU$. Clearly, this is the case if $\nabla_{\partial_{x_i}}\partial_{x_j} \equiv 0$, i.e., if $\nabla$ is flat. As it turns out, $J_\nabla$ is integrable if and only if $\nabla$ is symmetric and flat (this result was first written down by Dombrowski, cf. \cite{Dombrowski1962}, and independently by Tachibana and Okumura, cf. \cite{Tachibana1961}).\\
Next, we want to transfer the almost complex structure $J_\nabla$ from $TM$ to $T^\ast M$. In general, we can pick any bundle isomorphism $TM\to T^\ast M$, which is then also a diffeomorphism between the manifolds $TM$ and $T^\ast M$, and translate $J_\nabla$ using this diffeomorphism. However, there is no canonical choice of bundle isomorphism for generic manifolds $M$ with connection $\nabla$. The situation is different if $M$ is equipped with a semi-Riemannian metric $g$ and $\nabla$ is the Levi-Civita connection $\nabla^g$. In this case, we can choose $G:TM\to T^\ast M$, $v\mapsto \iota_v g$ as our bundle isomorphism and define the almost complex structure $J^\ast_{\nabla^g}$ on $T^\ast M$ via:
\begin{gather*}
 J^\ast_{\nabla^g}\coloneqq dG\circ J_{\nabla^g}\circ dG^{-1}.
\end{gather*}
Again, our goal is to express $J^\ast_{\nabla^g}$ in normal coordinates $\psi = (x_1,\ldots, x_n)$ of $(M, g)$ near $p\in M$. To achieve this, we first note that normal coordinates of a semi-Riemannian manifold $(M,g)$ near $p\in M$ satisfy:
\begin{gather*}
 g(p) = \sum^s_{k = 1} dx^2_{k,p} - \sum^n_{k = s+1} dx^2_{k,p}\quad\text{and}\quad \partial_{x_k} g_{lm} (p) = 0,
\end{gather*}
where $(s,t)$ is the signature of $g$ and $n = s+t$ is the dimension of $M$. As before, we introduce the notation $T\psi = (\hat x_1,\ldots, \hat x_n, v_1,\ldots, v_n)$ for the induced coordinates on $TM$. Similarly, we employ the notation $T^\ast \psi = (q_1,\ldots, q_n, p_1,\ldots, p_n)$ for the induced coordinates on $T^\ast M$:
\begin{gather*}
 (T^\ast\psi)^{-1} (q_1,\ldots, q_n, p_1,\ldots, p_n) \coloneqq \sum^n_{k = 1} p_k dx_{k, \psi^{-1}(q_1,\ldots, q_n)}.
\end{gather*}
In these coordinates, $G$ is given by:
\begin{align*}
 T^\ast\psi\circ G\circ (T\psi)^{-1} (0,\ldots, 0, v_1,\ldots, v_n) = (0,\ldots, 0, v_1,\ldots, v_s, -v_{s+1},\ldots, -v_n).
\end{align*}
Together with $\partial_{x_k} g_{lm} (p) = 0$, this implies:
\begin{gather*}
 dG_w (\partial_{\hat x_k, w}) = \partial_{q_k, G(w)},\ dG_w (\partial_{v_k, w}) = \begin{cases}\partial_{p_k, G(w)}\text{ for }1\leq k\leq s\\ -\partial_{p_k, G(w)}\text{ for } k>s\end{cases}\ \forall w\in T_pM.
\end{gather*}
This allows us to compute $J^\ast_{\nabla^g}$ in coordinates for points $\alpha\in T^\ast_p M$:
\begin{alignat*}{3}
 &J^\ast_{\nabla^g} (\partial_{q_k,\alpha}) = -\partial_{p_k, \alpha},\quad &&J^\ast_{\nabla^g} (\partial_{p_k, \alpha}) = \phantom{-} \partial_{q_k, \alpha}\quad &&(1\leq k\leq s),\\
 &J^\ast_{\nabla^g} (\partial_{q_k, \alpha}) = \phantom{-} \partial_{p_k, \alpha},\quad &&J^\ast_{\nabla^g} (\partial_{p_k, \alpha}) = -\partial_{q_k, \alpha}\quad &&(k>s).
\end{alignat*}
Again, this does not imply that $J^\ast_{\nabla^g}$ is integrable, as the equations above are only true for points $\alpha\in T^\ast_p M$ and not necessarily the entire chart domain $T^\ast U$. Indeed, $J^\ast_{\nabla^g}$ is integrable if and only if $g$ is flat:
\begin{alignat*}{2}
 &\phantom{\Leftrightarrow} && J^\ast_{\nabla^g}\text{ is integrable.}\\
 &\Leftrightarrow\ && J_{\nabla^g}\text{ is integrable.}\\
 &\Leftrightarrow\ && \nabla^g\text{ is symmetric and flat.}\\
 &\Leftrightarrow\ && \nabla^g\text{ is flat.}\\
 &\Leftrightarrow\ && g\text{ is flat.}
\end{alignat*}
The reason why we are interested in the almost complex structure $J^\ast_{\nabla^g}$ is the curious fact that $J^\ast_{\nabla^g}$ is naturally compatible with the canonical symplectic form $\omega_{\can}$ on $T^\ast M$, as one easily checks: In the coordinates $(q_1,\ldots, q_n, p_1,\ldots, p_n)$ of $T^\ast M$ from above, $\omega_{\can}$ is given by:
\begin{gather*}
 \omega_{\can}\vert_{T^\ast U} = \sum^n_{k = 1} dp_k\wedge dq_k.
\end{gather*}
Thus, $\omega_{\can} (\cdot, J^\ast_{\nabla^g}\cdot)$ is a semi-Riemannian metric on $T^\ast M$ of signature $(2s, 2t)$:
\begin{gather*}
 \omega_{\can, \alpha} (\cdot, J^\ast_{\nabla^g}\cdot) = \sum^s_{k = 1} dq^2_{k,\alpha} + dp^2_{k, \alpha} - \sum^{n}_{k = s+1} dq^2_{k,\alpha} + dp^2_{k,\alpha}\quad\forall \alpha\in T_pM.
\end{gather*}

\pagebreak

As we can find normal coordinates near any point $p\in M$, we have proven the following theorem:
\begin{theorem}[Almost complex structures on $TM$ and $T^\ast M$]\label{thm:alm_cpx_str_on_tan}
 Let $M$ be a smooth, $n$-dimensional manifold together with a connection $\nabla$ on it. Then, there exists a unique almost complex structure $J_\nabla$ on the tangent bundle $TM$ such that:
 \begin{gather*}
  d\pi\circ J_\nabla = K,\quad K\circ J_\nabla = -d\pi,
 \end{gather*}
 where $\pi:TM\to M$ is the base point projection of $TM$ and $K:T(TM)\to TM$ is the vertical projection corresponding to $\nabla$. $J_\nabla$ is integrable if and only if $\nabla$ is symmetric and flat.\\
 If $\nabla$ is symmetric, then $J_\nabla$ can be expressed as:
 \begin{gather*}
  J_\nabla (\partial_{v_k, w}) = \partial_{\hat x_k, w},\quad J_\nabla (\partial_{\hat x_k, w}) = -\partial_{v_k, w},
 \end{gather*}
 where $w\in TM$ is a point, $\psi = (x_1,\ldots, x_n)$ are normal coordinates of $(M,\nabla)$ near $p = \pi (w)$, and $T\psi = (\hat x_1,\ldots, \hat x_n, v_1,\ldots, v_n)$ are the induced coordinates on $TM$.\\
 If $\nabla = \nabla^g$ is the Levi-Civita connection of a semi-Riemannian metric $g$ on $M$ of signature $(s,t)$, then there exists a (canonical) almost complex structure $J^\ast_{\nabla^g}$ on $T^\ast M$ such that $\omega_{\can} (\cdot, J^\ast_{\nabla^g}\cdot)$ is a semi-Riemannian metric on $T^\ast M$ of signature $(2s, 2t)$, where $\omega_{\can}$ is the canonical symplectic form on $T^\ast M$. In coordinates, $J^\ast_{\nabla^g}$ is given by:
 \begin{alignat*}{3}
  &J^\ast_{\nabla^g} (\partial_{q_k, \alpha}) = -\partial_{p_k, \alpha},\quad &&J^\ast_{\nabla^g} (\partial_{p_k, \alpha}) = \phantom{-} \partial_{q_k, \alpha}\quad &&(1\leq k\leq s),\\
  &J^\ast_{\nabla^g} (\partial_{q_k, \alpha}) = \phantom{-} \partial_{p_k, \alpha},\quad &&J^\ast_{\nabla^g} (\partial_{p_k, \alpha}) = -\partial_{q_k, \alpha}\quad &&(k>s),
 \end{alignat*}
 where $\alpha\in T^\ast_p M$ is a point, $\psi = (x_1,\ldots, x_n)$ are normal coordinates of $(M,g)$ near $p\in M$, and $T^\ast \psi = (q_1,\ldots, q_n, p_1,\ldots, p_n)$ are the induced coordinates of $T^\ast M$. Furthermore, $J^\ast_{\nabla^g}$ is integrable if and only if $g$ is flat.
\end{theorem}

 \chapter{Holomorphic Connections}
 \label{app:holo_connection}
 In this part, we introduce the notion of a (linear/affine) holomorphic connection on a complex manifold $X$ (cf. the end of Section 4.2 in \cite{huybrechts2005}). In particular, we relate holomorphic connections to real connections (first section), study their interaction with real structures (second section), and define the holomorphic Levi-Civita connection $\nabla^h$ induced by a holomorphic metric $h = h_R + ih_I$ on $X$ (third section). Moreover, we show in the third section that the standard Levi-Civita connections $\nabla^{h_R}$ and $\nabla^{h_I}$ of the real and imaginary part $h_R$ and $h_I$ agree with each other, that $\nabla^{h_R} = \nabla^{h_I}$ is the connection associated (cf. first section) with the holomorphic Levi-Civita connection $\nabla^h$, and that holomorphic normal coordinates of $h$ are also normal coordinates of $h_R$.

\section*{Holomorphic and Associated Connections}

We begin by recalling the standard definition of a (linear/affine) connection $\nabla$ on a smooth manifold $X$:

\begin{definition}[Real linear connection]\label{def:rlc}
 Let $X$ be a smooth manifold. A \textbf{real linear connection}\footnote{We call them real connections to differentiate them from complex connections.} (\rlc for short) on $X$ is a $\R$-bilinear map\footnote{Here, $\Gamma (E)$ denotes the space of sections of the bundle $E\to B$. If $E\to B$ is a holomorphic bundle, $\Gamma (E)$ denotes the space of holomorphic sections, while $\Gamma_{C^\infty}(E)$ denotes the space of smooth sections.}
 \begin{gather*}
  \nabla: \Gamma (TX)\times\Gamma (TX)\to \Gamma (TX)
 \end{gather*}
 satisfying for all $f\in C^\infty (X)$ and $Y,Z\in\Gamma (TX)$\dots
 \begin{enumerate}
  \item \dots tensoriality in the first component: $\nabla_{fY}Z = f\nabla_Y Z$,
  \item \dots the Leibniz rule in the second component: $\nabla_Y fZ = Y(f)Z + f\nabla_YZ$.
 \end{enumerate}
\end{definition}

Our goal is to define holomorphic connections, thus, we must also allow for complex-valued connections:

\begin{definition}[Complex linear connection]\label{def:clc}
 Let $X$ be a smooth manifold. A \textbf{complex linear connection} (\clc for short) on $X$ is a $\C$-bilinear map\footnote{$T_\C X$ denotes the complexified tangent bundle.}
 \begin{gather*}
  \nabla: \Gamma (T_\C X)\times\Gamma (T_\C X)\to \Gamma (T_\C X)
 \end{gather*}
 satisfying for all $f\in C^\infty (X, \C)$ and $Y,Z\in\Gamma (T_\C X)$\dots
 \begin{enumerate}
  \item \dots tensoriality in the first component: $\nabla_{fY}Z = f\nabla_Y Z$,
  \item \dots the Leibniz rule in the second component: $\nabla_Y fZ = Y(f)Z + f\nabla_YZ$.
 \end{enumerate}
 A \clc $\nabla$ is called \textbf{real} if $\overline{\nabla_Y Z} = \nabla_{\overline{Y}}\overline{Z}$ holds for all complex vector fields $Y,Z\in\Gamma (T_\C X)$.
\end{definition}

Obviously, a \rlc $\nabla$ induces a \clc by complexification:
\begin{gather*}
 \nabla_{U+iV} (Y+iZ)\coloneqq \nabla_UY - \nabla_VZ + i\left(\nabla_VY + \nabla_UZ\right)\quad\forall U,V,Y,Z\in\Gamma (TX).
\end{gather*}
In light of this observation, one naturally asks which \clc are induced by \rlc As the name suggests, these are exactly the real ones:

\begin{proposition}[R.l.c. $\Leftrightarrow$ real c.l.c]\label{prop:real_clc}
 Let $\nabla$ be a \clc on a smooth manifold $X$. Then, we have the following equivalence:
 \begin{gather*}
  \nabla\text{ is induced by a \rlc} \quad\Leftrightarrow\quad \nabla\text{ is real}.
 \end{gather*}
\end{proposition}

\begin{proof}
 The direction ``$\Rightarrow$'' is trivial, since $\overline{\nabla_Y Z} = \nabla_{\overline{Y}}\overline{Z}$ is true by construction. Now consider the converse direction. Let $\nabla$ be a real \clc Then, we have $\overline{\nabla_Y Z} = \nabla_Y Z$ for all vector fields $Y,Z\in\Gamma (TX)\subset\Gamma (T_\C X)$, since $\overline{Y} = Y$ and $\overline{Z} = Z$. In particular, we find $\nabla_Y Z\in\Gamma (TX)$ for $Y,Z\in\Gamma (TX)$. This allows us to define a \rlc $\nabla^r$ by setting:
 \begin{gather*}
  \nabla^r_Y Z \coloneqq \nabla_Y Z\quad\forall Y,Z\in\Gamma (TX).
 \end{gather*}
 It is easy to check that $\nabla^r$ induces $\nabla$ concluding the proof.
\end{proof}

\begin{remark}\label{rem:real_clc}
 Considering Proposition \autoref{prop:real_clc}, it is not necessary to distinguish between \rlc and real \clc Thus, we use both notions interchangeably from now on.
\end{remark}

Before we move on to holomorphic connections, let us quickly recall that \rlc and \clc are local in the sense that, given vector fields $Y$ and $Z$ as well as a point $p\in X$, the expression $\nabla_Y Z (p)$ only depends on $Y(p)$ and $Z$ in a neighborhood of $p$. Indeed, this fact directly follows from the definition of a connection and the existence of smooth partitions of unity. The locality of a connection allows us to introduce the concept of Christoffel symbols which (locally) completely determine the connection in question:

\begin{definition}[Christoffel symbols]\label{def:christoffel_symbols}
 Let $X$ be a smooth manifold, $\nabla$ be a \clc on $X$, and $U\subset X$ be an open subset. Further, let $\{V_j\}$ be a complex frame of $T_\C U$. Then, the \textbf{Christoffel symbols} $\Gamma^l_{jk}\in C^\infty (U,\C)$ of $\nabla$ with respect to the frame $\{V_j\}$ are defined by:
 \begin{gather*}
  \nabla_{V_j}V_k\eqqcolon\sum_l\Gamma^l_{jk}V_l.
 \end{gather*}
\end{definition}

\begin{remark}\label{rem:real_christoffel_symbols}
 If the connection $\nabla$ and the frame $\{V_j\}$ are real in Definition \autoref{def:christoffel_symbols}, then the Christoffel symbols $\Gamma^l_{jk}$ are real as well and coincide with the standard notion of Christoffel symbols.
\end{remark}

Now let $X$ be a complex manifold. We wish to define a holomorphic connection $\nabla$ on $X$. However, we cannot just copy Definition \autoref{def:rlc}. The problem is that there is no holomorphic partition of unity which in turn means that locality does not immediately follow from the definition anymore, but needs to be enforced manually. This leads us to the following definition (cf. \cite{huybrechts2005}):

\begin{definition}[Holomorphic connection]\label{def:holo_connection}
 Let $X$ be a complex manifold. We call $\nabla$ a \textbf{holomorphic connection} on $X$ if $\nabla$ is a collection
 \begin{gather*}
  \{\nab{U}:\hvec{U}\times\hvec{U}\to\hvec{U}\mid U\subset X\text{\normalfont\ open}\}
 \end{gather*}
 of $\C$-bilinear maps satisfying for all open subsets $U^\prime\subset U\subset X$, all holomorphic functions $f:U\to\C$, and all holomorphic vector fields $Y,Z\in\hvec{U}$\dots
 \begin{enumerate}
  \item \dots tensoriality in the first component: $\nab{U}_{fY} Z = f\nab{U}_Y Z$,
  \item \dots the Leibniz rule in the second component: $\nab{U}_Y fZ = Y(f) Z + f\nab{U}_Y Z$,
  \item \dots the presheaf property: $(\nab{U}_Y Z)\vert_{U^\prime} = \nab{U^\prime}_{Y\vert_{U^\prime}} Z\vert_{U^\prime}$.
 \end{enumerate}
 We say that $\nabla$ is \textbf{induced} by a \clc $\nabla^\prime$ on $X$ if $\nab{U}_YZ = \nabla^\prime_YZ$ for all holomorphic vector fields $Y,Z\in\hvec{U}\subset\Gamma (T_\C U)$ and open $U\subset X$.
\end{definition}

\begin{remark}[Presheaf property]\label{rem:presheaf}
 To understand the presheaf property, consider a compact complex manifold $X$, e.g. a complex torus $\C^n/\Lambda$ ($\Lambda\subset\C^n$ is a lattice). Since holomorphic functions on compact complex manifolds are locally constant (due to the maximum principle), the Leibniz rule reduces to tensoriality for compact $X$. Hence, without the presheaf property, $\nabla = 0$ would be a holomorphic connection on compact complex manifolds contradicting our intuition and the local picture.
\end{remark}

Similar to \rlc and c.l.c., a holomorphic connection can be computed locally:

\begin{proposition}[Holomorphic connections are local]\label{prop:holo_con_local}
 Let $X$ be a complex manifold with holomorphic connection $\nabla$. For all points $p\in X$, all open neighborhoods $U^\prime\subset U\subset X$ of $p$, and all holomorphic vector fields $Y,Z\in \hvec{U}$, the expression $\nab{U}_Y Z (p)$ is completely determined by $Y(p)$ and $Z\vert_{U^\prime}$.
\end{proposition}

\begin{proof}
 Take the notations from above, then:
 \begin{gather*}
  \nab{U}_Y Z (p) = (\nab{U}_Y Z)\vert_{U^\prime} (p) \stackrel{\text{(iii)}}{=} \nab{U^\prime}_{Y\vert_{U^\prime}} Z\vert_{U^\prime} (p)
 \end{gather*}
 Thus, $\nab{U}_Y Z(p)$ only depends on $Y\vert_{U^\prime}$ and $Z\vert_{U^\prime}$. By going into holomorphic charts of $X$ near $p$ and using Property (i), we see that $\nab{U}_Y Z (p)$ only depends on $Y(p)$ and $Z\vert_{U^\prime}$.
\end{proof}

Again, locality allows us to introduce Christoffel symbols which (locally) completely determine the connection:

\begin{definition}[Holomorphic Christoffel symbols]\label{def:holo_christoffel_symbols}
 Let $X$ be a complex manifold, $\nabla$ be a holomorphic connection on $X$, and $U\subset X$ be an open subset. Further, let $\{V_j\}$ be a holomorphic frame of $T^{(1,0)} U$. Then, the \textbf{holomorphic Christoffel symbols} $\Gamma^l_{jk}:U\to\C$ of $\nabla$ with respect to the frame $\{V_j\}$ are holomorphic functions on $U$ defined by:
 \begin{gather*}
  \nab{U}_{V_j}V_k\eqqcolon\sum_l\Gamma^l_{jk}V_l.
 \end{gather*}
\end{definition}

Let us now study the relation between \clc and holomorphic connections. Our goal is to associate a unique \clc with each holomorphic connection. The difficulty of this task lies in the fact that multiple \clc may induce the same holomorphic connection, as the following example demonstrates:

\begin{example}[Holomorphic connection on $\C^n$]\label{ex:holo_connection}
 Consider the complex manifold $X = \C^n$ with global chart $\id = (z_1,\ldots, z_n):X\to\C^n$. Then, $\{\pa{z_j}\}_j$ is a holomorphic frame of $T^{(1,0)}\C^n$. Now pick holomorphic functions $\Gamma^l_{jk}:\C^n\to\C$. Using these functions as Christoffel symbols, we can define a holomorphic connection $\nabla$ on $\C^n$:
 \begin{gather*}
  \nab{\C^n}_{\pa{z_j}}\pa{z_k} = \sum_l \Gamma^l_{jk} \pa{z_l}.
 \end{gather*}
 Obviously, the remaining maps $\nab{U}$ ($U\subset\C^n$ open) of the holomorphic connection $\nabla$ are defined by the presheaf property.\\
 Now observe that $\{\pa{z_j},\pa{\bar{z}_j}\}_j$ is a smooth frame of $T^{(1,0)}\C^n\oplus T^{(0,1)}\C^n = T_\C\C^n$ and pick functions $f^l_{jk},g^l_{jk}\in C^\infty (\C^n,\C)$. With these, we can define the \clc $\nabla^1$, $\nabla^2$, $\nabla^3$, and $\nabla^4$ on $\C^n$:
 \begin{enumerate}[label = (\arabic*)]
  \item $\nabla^1_{\pa{z_j}}\pa{z_k} = \sum_l \Gamma^l_{jk} \pa{z_l}$ and $\nabla^1_{\pa{\bar{z}_j}}\pa{z_k} = \nabla^1_{\pa{z_j}}\pa{\bar{z}_k} = \nabla^1_{\pa{\bar{z}_j}}\pa{\bar{z}_k} = 0$,
  \item $\nabla^2_{\pa{z_j}}\pa{z_k} = \overline{ \nabla^2_{\pa{\bar{z}_j}}\pa{\bar{z}_k}} = \sum_l \Gamma^l_{jk} \pa{z_l}$ and $\nabla^2_{\pa{\bar{z}_j}}\pa{z_k} = \nabla^2_{\pa{z_j}}\pa{\bar{z}_k} = 0$,
  \item $\nabla^3_{\pa{z_j}}\pa{z_k} = \sum_l \Gamma^l_{jk} \pa{z_l}$, $\nabla^3_{\pa{\bar{z}_j}}\pa{z_k} = \sum_l f^l_{jk}\pa{z_l} + g^l_{jk}\pa{\bar{z}_l}$, and\\ $\nabla^3_{\pa{z_j}}\pa{\bar{z}_k} = \nabla^3_{\pa{\bar{z}_j}}\pa{\bar{z}_k} = 0$,
  \item $\nabla^4_{\pa{z_j}}\pa{z_k} = \overline{ \nabla^4_{\pa{\bar{z}_j}}\pa{\bar{z}_k}} = \sum_l \Gamma^l_{jk} \pa{z_l}$ and\\ $\nabla^4_{\pa{\bar{z}_j}}\pa{z_k} = \overline{ \nabla^4_{\pa{z_j}}\pa{\bar{z}_k}} = \sum_l f^l_{jk}\pa{z_l} + g^l_{jk}\pa{\bar{z}_l}$.
 \end{enumerate}
 It is clear from the construction that $\nabla^1$, $\nabla^2$, $\nabla^3$, and $\nabla^4$ induce the same holomorphic connection $\nabla$. In general, however, none of the connections coincide.
\end{example}

To single out a unique \clc inducing a holomorphic connection, we need to impose additional constraints on the \clc One of those constraints is that the \clc is real. Another one is compatibility with the complex structure:

\begin{definition}[$\nabla$ compatible with $J$]\label{def:nabla_compatible}
 Let $\nabla$ be a \clc on a smooth manifold $X$ with almost complex structure $J$. $\nabla$ is \textbf{compatible} with $J$ if $\nabla J = 0$.
\end{definition}

By the next proposition, the notion of compatibility is equivalent to $\nabla$ preserving the $(1,0)$- and $(0,1)$-type of the second component:

\begin{proposition}\label{prop:nabla_compatible}
 Let $X$ be a smooth manifold with almost complex structure $J$. A \clc $\nabla$ on $X$ is compatible with $J$ if and only if $\nabla$ satisfies:
 \begin{enumerate}
  \item $\nabla_Y Z\in \Gamma_{C^\infty} (T^{(1,0)}X)\quad \forall Y\in \Gamma (T_\C X)\ \forall Z\in \Gamma_{C^\infty} (T^{(1,0)}X)$,
  \item $\nabla_Y Z\in \Gamma_{C^\infty} (T^{(0,1)}X)\quad \forall Y\in \Gamma (T_\C X)\ \forall Z\in \Gamma_{C^\infty} (T^{(0,1)}X)$.
 \end{enumerate}
\end{proposition}

\begin{proof}
 ``$\Rightarrow$'': Pick $(k,l)\in\{(1,0),\, (0,1)\}$. Further, take vector fields $Y\in \Gamma (T_\C X)$ and $Z\in \Gamma_{C^\infty} (T^{(k,l)}X)$. Then, we can write $Z = 1/2 (\hat Z + (-1)^kiJ\hat Z)$ for a unique real vector field $\hat Z$ on $X$. Thus, we obtain:
 \begin{align*}
  \nabla_Y Z &\stackrel{\phantom{\nabla J = 0}}{=} \frac{1}{2}\left(\nabla_Y\hat Z + (-1)^ki\nabla_Y (J\hat Z)\right)\\
  &\stackrel{\nabla J = 0}{=} \frac{1}{2}\left(\nabla_Y\hat Z + (-1)^kiJ (\nabla_Y \hat Z)\right)\in\Gamma_{C^\infty} (T^{(k,l)}X).
 \end{align*}
 ``$\Leftarrow$'': Take smooth complex vector fields $Y, Z\in \Gamma (T_\C X)$. Then, there are unique vector fields $Z^{(1,0)}\in\Gamma_{C^\infty} (T^{(1,0)}X)$ and $Z^{(0,1)}\in\Gamma_{C^\infty} (T^{(0,1)}X)$ such that\linebreak $Z = Z^{(1,0)} + Z^{(0,1)}$. By assumption, we have:
 \begin{gather*}
  \nabla_Y Z^{(1,0)}\in\Gamma_{C^\infty} (T^{(1,0)}X) \quad\text{and}\quad \nabla_Y Z^{(0,1)}\in\Gamma_{C^\infty} (T^{(0,1)}X).
 \end{gather*}
 We now compute $(\nabla_Y J) Z$:
 \begin{align*}
  (\nabla_Y J)Z &= \nabla_Y (JZ) - J\nabla_Y Z\\
  &= \nabla_Y (iZ^{(1,0)} - iZ^{(0,1)}) - i\nabla_YZ^{(1,0)} + i\nabla_Y Z^{(0,1)} = 0.
 \end{align*}
\end{proof}

\begin{remark}
 If $\nabla$ is real, we can simplify Proposition \autoref{prop:nabla_compatible} as follows: A \rlc $\nabla$ on $X$ is compatible with $J$ if and only if $\nabla$ satisfies:
 \begin{gather*}
  \nabla_Y Z\in \Gamma_{C^\infty} (T^{(1,0)}X)\quad \forall Y\in \Gamma (T_\C X)\ \forall Z\in \Gamma_{C^\infty} (T^{(1,0)}X).
 \end{gather*}
 Of course, the same statement is true if we replace ``$(1,0)$'' with ``$(0,1)$''.
\end{remark}

Unfortunately, those two constraints are still insufficient to single out a unique \clc given a holomorphic connection. To illustrate this, consider Example \autoref{ex:holo_connection} again. In general, the connections from Example \autoref{ex:holo_connection} satisfy the following properties:
\begin{enumerate}[label = (\arabic*)]
 \item $\nabla^1$ is compatible with the complex structure, but not real.
 \item $\nabla^2$ is both compatible with the complex structure and real.
 \item $\nabla^3$ is neither compatible with the complex structure nor real.
 \item $\nabla^4$ is real, but not compatible with the complex structure.
\end{enumerate}
However, if all functions $g^l_{jk}$ are identically zero, then both $\nabla^2$ and $\nabla^4$ are real and compatible with the complex structure. In this case, they still might not be equal, since the functions $f^l_{jk}$ do not need to vanish.\\
The problem is that compatibility as a constraint is too weak. Even though compatible \rlc preserve the $(1,0)$- and $(0,1)$-type of the second component, they still allow for a mixing of types between the first and second component. Excluding this phenomenon gives us a unique \clc for every holomorphic connection:

\pagebreak

\begin{lemma}[Associated connection $\nabla^\prime$]\label{lem:associated_connection}
 Let $X$ be a complex manifold with holomorphic connection $\nabla$ on it. Then, there exists a unique \clc $\nabla^\prime$ on $X$ such that:
 \begin{enumerate}
  \item $\nabla^\prime$ induces $\nabla$.
  \item $\nabla^\prime$ is real.
  \item $\nabla^\prime_{\overline{Y}} Z = 0$ for all $Y,Z\in\Gamma (T^{(1,0)}U)\subset\Gamma (T_\C U)$ and open $U\subset X$.
 \end{enumerate}
 We call $\nabla^\prime$ the connection \textbf{associated} with $\nabla$.
\end{lemma}

\begin{proof}
 The idea behind the proof is to first show uniqueness of $\nabla^\prime$ and afterwards use connection $\nabla^2$ from Example \autoref{ex:holo_connection} as a local model to construct $\nabla^\prime$.\\\\
 \textbf{Uniqueness:} Let $\nabla^\prime$ and $\nabla^{\prime\prime}$ be two \clc associated with $\nabla$. Further, pick a point $p\in X$ and a holomorphic chart $(z_1,\ldots, z_n):V\to\C^n$ of $X$ near $p$. Now observe that $\{\pa{z_j}, \pa{\bar{z}_j}\}_j$ is a frame of $T_\C V$. If $\nabla^\prime$ and $\nabla^{\prime\prime}$ map the frame $\{\pa{z_j}, \pa{\bar{z}_j}\}_j$ to the same vector fields, they coincide on the neighborhood $V$ of $p$. Since both $\nabla^\prime$ and $\nabla^{\prime\prime}$ induce $\nabla$, we obtain:
 \begin{gather*}
  \nabla^\prime_{\pa{z_j}}\pa{z_k} = \nab{V}_{\pa{z_j}}\pa{z_k} = \nabla^{\prime\prime}_{\pa{z_j}}\pa{z_k}.
 \end{gather*}
 Using the fact that $\nabla^\prime$ and $\nabla^{\prime\prime}$ are also real yields:
 \begin{gather*}
  \nabla^\prime_{\pa{\bar{z}_j}}\pa{\bar{z}_k} = \overline{\nabla^\prime_{\pa{z_j}}\pa{z_k}} = \overline{\nabla^{\prime\prime}_{\pa{z_j}}\pa{z_k}} = \nabla^{\prime\prime}_{\pa{\bar{z}_j}}\pa{\bar{z}_k}.
 \end{gather*}
 Lastly, combining Property (iii) with realness gives:
 \begin{gather*}
  \nabla^\prime_{\pa{\bar{z}_j}}\pa{z_k} = \overline{\nabla^\prime_{\pa{z_j}}\pa{\bar{z}_k}} = 0 =  \overline{\nabla^{\prime\prime}_{\pa{z_j}}\pa{\bar{z}_k}} = \nabla^{\prime\prime}_{\pa{\bar{z}_j}}\pa{z_k}.
 \end{gather*}
 Hence, $\nabla^\prime$ and $\nabla^{\prime\prime}$ coincide on a neighborhood of $p$. Repeating the previous argument for every point $p\in X$ proves uniqueness.\\\\
 \textbf{Existence:} Pick a holomorphic chart $(z_1,\ldots, z_n):V\to\C^n$ of $X$. We observe that $\{\pa{z_j}\}_j$ is a holomorphic frame of $T^{(1,0)}V$. Let $\Gamma^l_{jk}$ be the holomorphic Christoffel symbols of $\nabla$ with respect to that frame. We now define the \clc $\nabla^\prime$ on $V$:
 \begin{gather*}
  \nabla^\prime_{\pa{z_j}}\pa{z_k} = \overline{\nabla^\prime_{\pa{\bar{z}_j}}\pa{\bar{z}_k}} = \sum_l \Gamma^l_{jk} \pa{z_l},\quad \nabla^\prime_{\pa{\bar{z}_j}}\pa{z_k} = \nabla^\prime_{\pa{z_j}}\pa{\bar{z}_k} = 0
 \end{gather*}
 Clearly, $\nabla^\prime$ satisfies Property (i), (ii), and (iii) on $V$. We now extend the definition of $\nabla^\prime$ to all of $X$ by repeating this construction for any holomorphic chart of $X$. Since the connection associated with $\nabla$ is unique, the definition of $\nabla^\prime$ agrees on different coordinate patches concluding the proof.
\end{proof}

\begin{remark}[Associated connection]\label{rem:associated_connection}\ 
 \begin{itemize}
  \item Note that the associated connection $\nabla^\prime$ in Lemma \autoref{lem:associated_connection} is compatible with the complex structure of $X$. Indeed, it is easy to check in holomorphic charts that $\nabla^\prime$ preserves the $(1,0)$-type of the second component (cf. the proof of Lemma \autoref{lem:associated_connection}).
  \item Usually, we denote a holomorphic connection and its associated connection by the same symbol.
 \end{itemize}
\end{remark}

\section*{Holomorphic Connections on Manifolds with Real Structures}

Next, we investigate the interaction of holomorphic connections $\nabla$ with real structures $\sigma$ (cf. \autoref{app:real_structures}). In particular, we show that $\sigma$-invariant connections $\nabla$ induce connections on real forms, that two holomorphic connections inducing the same connection on a real form coincide, and that every real-analytic connection on a real form gives rise to a unique holomorphic connection on a neighborhood of the real form.\\
To state and prove these results, we first need to study the relation between connections and diffeomorphisms:

\begin{definition}[$\phi^\ast\nabla$, $\phi_\ast\nabla$, and $\phi$-invariance]\label{def:phi-invariance}
 Let $M,N$ be smooth manifolds, $\phi:M\to N$ be a diffeomorphism, and $\nabla$ be a \clc on $N$. Then, the \textbf{pullback connection} $\phi^\ast\nabla$ on $M$ is defined by:
 \begin{gather*}
  (\phi^\ast\nabla)_Y Z\coloneqq \phi^{-1}_\ast\left(\nabla_{\phi_\ast Y} \phi_\ast Z\right)\quad\forall Y,Z\in\Gamma (T_\C M).
 \end{gather*}
 If $\nabla$ is a \clc on $M$, then the \textbf{pushforward connection} $\phi_\ast\nabla$ on $N$ is defined by:
 \begin{gather*}
  (\phi_\ast\nabla)_Y Z\coloneqq ((\phi^{-1})^\ast\nabla)_Y Z = \phi_\ast\left(\nabla_{\phi^{-1}_\ast Y} \phi^{-1}_\ast Z\right)\quad\forall Y,Z\in\Gamma (T_\C N).
 \end{gather*}
 In the case of $M = N$, we say that $\nabla$ is $\mathbf{\phi}$-\textbf{invariant} if $\phi^\ast\nabla = \nabla$ or, equivalently, $\phi_\ast\nabla = \nabla$. For complex manifolds $M = N$ and holomorphic connections $\nabla$ on $M$, $\nabla$ is called $\phi$-invariant if the connection associated with $\nabla$ is $\phi$-invariant.
\end{definition}

We can now explain how $\sigma$-invariant connections induce connections on real forms:

\begin{lemma}[Induced connection on real form]\label{lem:induced_connection_on_real_form}
 Let $X$ be a complex manifold with real structure $\sigma$ and non-empty real form $M \coloneqq \Fix\sigma$. Further, let $\nabla$ be a $\sigma$-invariant \rlc on $X$. Then, $\nabla$ induces a \rlc $\nabla^M$ on $M$. In particular, if $\nabla$ is the connection associated with a $\sigma$-invariant holomorphic connection on $X$, $\nabla^M$ is real-analytic.
\end{lemma}

\begin{proof}
 The idea of the proof is to use the Christoffel symbols of $\nabla$ in $\sigma$-charts (cf. \autoref{app:real_structures}) to construct the induced connection on $M$. Pick a point $p\in M$ and a $\sigma$-chart $(z_1,\ldots, z_n):V\to\C^n$ near $p$. In the frame $\{\pa{z_j},\pa{\bar{z}_j}\}_j$,\linebreak $\nabla$ is given by:
 \begin{gather*}
  \nabla_{\pa{z_j}}\pa{z_k} = \sum_l a^l_{jk}\pa{z_l} + b^l_{jk}\pa{\bar{z}_l},\quad \nabla_{\pa{\bar{z}_j}}\pa{z_k} = \sum_l c^l_{jk}\pa{z_l} + d^l_{jk}\pa{\bar{z}_l},
 \end{gather*}
 where $a^l_{jk}, b^l_{jk}, c^l_{jk}, d^l_{jk}\in C^\infty (V,\C)$ are Christoffel symbols. The $\sigma$-invariance of $\nabla$ now implies $a^l_{jk} = \overline{a^l_{jk}\circ\sigma}$ and similar equations for $b^l_{jk}, c^l_{jk}, d^l_{jk}$, as the following computation illustrates:
 \begin{align*}
  \sum_l a^l_{jk}\pa{z_l} + b^l_{jk}\pa{\bar{z}_l} &= \nabla_{\pa{z_j}}\pa{z_k} = (\sigma^\ast\nabla)_{\pa{z_j}}\pa{z_k} = \sigma_\ast\left(\nabla_{\sigma_\ast\pa{z_j}}\sigma_\ast\pa{z_k}\right)\\
  &= \sigma_\ast\left(\nabla_{\pa{\bar{z}_j}}\pa{\bar{z}_k}\right) = \sigma_\ast\overline{\nabla_{\pa{z_j}}\pa{z_k}} = \sigma_\ast\sum_l \overline{a^l_{jk}}\pa{\bar{z}_l} + \overline{b^l_{jk}}\pa{z_l}\\
  &= \sum_l \overline{a^l_{jk}\circ\sigma}\pa{z_l} + \overline{b^l_{jk}\circ\sigma}\pa{\bar{z}_l}.
 \end{align*}
 From this, we deduce that the Christoffel symbols $a^l_{jk},\ldots,d^l_{jk}$ are real when evaluated on points $q\in V\cap M$. We can now define the \rlc $\nabla^M$ on $V\cap M$ via\footnote{For the sake of readability, we drop the restriction ``$\vert_{V\cap M}$'' in the definition of $\nabla^M$.}:
 \begin{gather*}
  \nabla^M_{\pa{x_j}}\pa{x_k} \coloneqq \sum_l \Gamma^l_{jk}\pa{x_l},
 \end{gather*}
 where $\Gamma^l_{jk}\coloneqq a^l_{jk} + b^l_{jk} + c^l_{jk} + d^l_{jk}$, we used the decomposition $z_j = x_j + iy_j$ and exploited the fact that the real part of a $\sigma$-chart gives rise to a real-analytic chart of $M$. We extend the definition of $\nabla^M$ to all of $M$ by repeating the construction for any point $p\in M$. This is possible, because the definition of $\nabla^M$ is independent of the choice of $\sigma$-chart. To see this, we first observe that, given holomorphic functions $f_k:V\to\C$, the following identities hold:
 \begin{gather*}
  \nabla_{\pa{z_j}}\sum_kf_k\pa{z_k} = \sum_k \frac{\partial f_k}{\partial z_j}\pa{z_k} + f_k\nabla_{\pa{z_j}}\pa{z_k},\quad \nabla_{\pa{\bar{z}_j}}\sum_kf_k\pa{z_k} = \sum_k f_k\nabla_{\pa{\bar{z}_j}}\pa{z_k}.
 \end{gather*}
 Here, we have exploited tensoriality in the second component as well as the fact that holomorphicity of $f_k$ implies $\pa{\bar{z}_j}f_k = 0$. These identities tell us that the coefficients $a^l_{jk}$ transform like holomorphic Christoffel symbols under a change of $\sigma$-chart, i.e.:
 \begin{gather*}
  \left(a^\prime\right)^{l^\prime}_{j^\prime k^\prime} = \sum_l \frac{\partial^2 z_l}{\partial z^\prime_{j^\prime}\partial z^\prime_{k^\prime}}\frac{\partial z^\prime_{l^\prime}}{\partial z_l} + \sum_{j,k,l}a^l_{jk}\frac{\partial z_j}{\partial z^\prime_{j^\prime}}\frac{\partial z_k}{\partial z^\prime_{k^\prime}}\frac{\partial z^\prime_{l^\prime}}{\partial z_l},
 \end{gather*}
 where $\left(a^\prime\right)^{l^\prime}_{j^\prime k^\prime}$ are Christoffel symbols of $\nabla$ in a different $\sigma$-chart\linebreak $(z^\prime_1,\ldots, z^\prime_n):V^\prime\to\C^n$ near $p$, while the remaining coefficients $b^l_{jk}, c^l_{jk}, d^l_{jk}$ transform like tensors, i.e.:
 \begin{align*}
  \left(b^\prime\right)^{l^\prime}_{j^\prime k^\prime} &= \sum_{j,k,l}b^l_{jk}\frac{\partial z_j}{\partial z^\prime_{j^\prime}}\frac{\partial z_k}{\partial z^\prime_{k^\prime}}\frac{\partial \bar{z}^\prime_{l^\prime}}{\partial \bar{z}_l},\\
  \left(c^\prime\right)^{l^\prime}_{j^\prime k^\prime} &= \sum_{j,k,l}c^l_{jk}\frac{\partial \bar{z}_j}{\partial \bar{z}^\prime_{j^\prime}}\frac{\partial z_k}{\partial z^\prime_{k^\prime}}\frac{\partial z^\prime_{l^\prime}}{\partial z_l},\\
  \left(d^\prime\right)^{l^\prime}_{j^\prime k^\prime} &= \sum_{j,k,l}d^l_{jk}\frac{\partial \bar{z}_j}{\partial \bar{z}^\prime_{j^\prime}}\frac{\partial z_k}{\partial z^\prime_{k^\prime}}\frac{\partial \bar{z}^\prime_{l^\prime}}{\partial \bar{z}_l},
 \end{align*}
 where $\left(b^\prime\right)^{l^\prime}_{j^\prime k^\prime}$, $\left(c^\prime\right)^{l^\prime}_{j^\prime k^\prime}$, and $\left(d^\prime\right)^{l^\prime}_{j^\prime k^\prime}$ are defined similarly to $\left(a^\prime\right)^{l^\prime}_{j^\prime k^\prime}$. For the definition of $\nabla^M$ to be independent of the choice of $\sigma$-chart, we need to show that the sum $\Gamma^l_{jk}$ evaluated on points of $M$ transforms as follows\footnote{Again, we drop the restriction to $M$.}:
 \begin{gather*}
  \left(\Gamma^\prime\right)^{l^\prime}_{j^\prime k^\prime} = \sum_l \frac{\partial^2 x_l}{\partial x^\prime_{j^\prime}\partial x^\prime_{k^\prime}}\frac{\partial x^\prime_{l^\prime}}{\partial x_l} + \sum_{j,k,l}\Gamma^l_{jk}\frac{\partial x_j}{\partial x^\prime_{j^\prime}}\frac{\partial x_k}{\partial x^\prime_{k^\prime}}\frac{\partial x^\prime_{l^\prime}}{\partial x_l}.
 \end{gather*}
 However, this transformation behavior is a direct consequence of the previous transformation rules and the fact that two $\sigma$-charts $(z_1,\ldots, z_n):V\to\C^n$ and $(z^\prime_1,\ldots, z^\prime_n):V^\prime\to\C^n$ near $p$ satisfy:
 \begin{gather*}
  \frac{\partial^r z_l}{\partial z^\prime_{s_1}\ldots\partial z^\prime_{s_r}} (q) = \frac{\partial^r \bar{z}_l}{\partial \bar{z}^\prime_{s_1}\ldots\partial \bar{z}^\prime_{s_r}} (q) = \frac{\partial^r x_l}{\partial x^\prime_{s_1}\ldots\partial x^\prime_{s_r}} (q)\quad\forall q\in V\cap V^\prime\cap M.
 \end{gather*}
 If $\nabla$ is a $\sigma$-invariant holomorphic connection on $X$, then its associated connection $\nabla$ is a $\sigma$-invariant \rlc on $X$. Hence, we can construct $\nabla^M$ as before. Since the Christoffel symbols $\Gamma^l_{jk} = a^l_{jk}$ are holomorphic and, thus, real-analytic in this case, $\nabla^M$ is real-analytic.
\end{proof}

\begin{remark}
 Note that a $\sigma$-invariant \rlc $\nabla$ on $X$ induces not just $\nabla^M$, but a plethora of \rlc on $M$. Indeed, the construction in the proof of Lemma \autoref{lem:induced_connection_on_real_form} still works if we substitute $\Gamma^l_{jk}$ for the affine combination
 \begin{gather*}
  a^l_{jk} + r_bb^l_{jk} + r_cc^l_{jk} + r_dd^l_{jk},
 \end{gather*}
 where $r_b,r_c,r_d\in\R$ are any real numbers. Among those connections, $\nabla^M$ stands out, as it satisfies:
 \begin{gather*}
  \nabla_{\pa{x_j}}\pa{x_k} (q) = \nabla^M_{\pa{x_j}}\pa{x_k} (q)\quad\forall q\in V\cap M,
 \end{gather*}
 where $(z_1 = x_1 + iy_1,\ldots, z_n = x_n + iy_n):V\to M$ is a $\sigma$-chart. Also note that if $\nabla$ is associated with a holomorphic connection, the coefficients $b^l_{jk},c^l_{jk},d^l_{jk}$ vanish and $\nabla$ only induces $\nabla^M$.
\end{remark}

To conclude this section, we want to show that two holomorphic connections inducing the same \rlc are equal on a neighborhood of the real form and that any real-analytic connection on a real form gives rise to a holomorphic connection. These results are the ``connection counterparts'' of Lemma \autoref{lem:uniqueness_of_holomorphic_tensors} and \autoref{lem:holo_continuation_of_tensors}. Indeed, the proofs are very similar and employ the same techniques which is why we keep them brief:

\begin{lemma}[Uniqueness of holomorphic connections]\label{lem:uniqueness_of_holomorphic_connections}
 Let $X$ be a complex manifold with real structure $\sigma$ and non-empty real form $M\coloneqq\Fix\sigma$. Further, let $\nabla^1$ and $\nabla^2$ be two $\sigma$-invariant holomorphic connections on $X$ which induce the same \rlc $\nabla^M$ on $M$. Then, there exists an open neighborhood $U\subset X$ of $M$ such that $\nab{U}^1 = \nab{U}^2$. If $M$ is nice, one even has $\nabla^1 = \nabla^2$.
\end{lemma}

\begin{proof}
 Pick a point $p\in M$ and a $\sigma$-chart $(z_1,\ldots, z_n):V\to\C^n$ near $p$ such that $V$ is connected. Now let $\left(\Gamma^1\right)^l_{jk}$ and $\left(\Gamma^2\right)^l_{jk}$ be the holomorphic Christoffel symbols of $\nabla^1$ and $\nabla^2$, respectively, in this $\sigma$-chart. Since $\nabla^1$ and $\nabla^2$ induce the same \rlc $\nabla^M$, their Christoffel symbols $\left(\Gamma^1\right)^l_{jk}$ and $\left(\Gamma^2\right)^l_{jk}$ coincide on $V\cap M$. Hence, by the identity theorem, they must coincide on all of $V$. This shows that $\nab{V}^1$ and $\nab{V}^2$ are equal. We can now repeat this argument for every point $p\in M$ to find a neighborhood $U\subset X$ of $M$ such that $\nab{U}^1 = \nab{U}^2$. If $M$ is nice, we can apply the identity theorem to every connected component of $X$ to prove $\nabla^1 = \nabla^2$.
\end{proof}

\begin{lemma}[Holomorphic continuation of connections]\label{lem:holo_continuation_of_connections}
 Let $X$ be a complex manifold with real structure $\sigma$ and non-empty real form $M\coloneqq\Fix\sigma$. Further, let $\nabla^M$ be a real-analytic \rlc on $M$. Then, there exists an open neighborhood $U\subset X$ of $M$ and a $\sigma$-invariant holomorphic connection $\nabla$ on $U$ which induces $\nabla^M$. $\nabla$ is unique in the sense of Lemma \autoref{lem:uniqueness_of_holomorphic_connections}.
\end{lemma}

\begin{proof}
 Pick a point $p\in M$ and a $\sigma$-chart $(z_1,\ldots, z_n):V\to\C^n$ near $p$. The real part of the $\sigma$-chart furnishes a real-analytic submanifold chart for $M$. Let $\left(\Gamma^M\right)^l_{jk}$ be the real-analytic Christoffel symbols of $\nabla^M$ in this submanifold chart. As in the proof of Lemma \autoref{lem:holo_continuation_of_tensors}, we can find holomorphic continuations $\Gamma^l_{jk}:V\to\C$ of $\left(\Gamma^M\right)^l_{jk}$ (after shrinking $V$ if necessary) and interpret them as holomorphic Christoffel symbols. We can now define the holomorphic connection $\nabla$ on $V$ by:
 \begin{gather*}
  \nab{V}_{\pa{z_j}}\pa{z_k} \coloneqq \sum_l \Gamma^l_{jk}\pa{z_l}.
 \end{gather*}
 It is easy to check that $\nabla$ is a $\sigma$-invariant holomorphic connection on $V$ inducing $\nabla^M$. To define $\nabla$ on a neighborhood $U$ of $M$, we repeat the last step for every point $p\in M$ and afterwards glue together the various coordinate patches. This is possible, since the construction is independent of the choice of chart which directly follows from Lemma \autoref{lem:uniqueness_of_holomorphic_connections}.
\end{proof}

\section*{Holomorphic Levi-Civita Connection}

In this section, we want to define the holomorphic Levi-Civita connection $\nabla^h$ of a holomorphic metric $h$. As in the real case, $\nabla^h$ should be the unique holomorphic connection on $X$ which is torsion-free and compatible with $h$. Thus, we first need to define the torsion $T$ of a holomorphic connection $\nabla$:

\begin{definition}[Torsion $T$]
 Let $X$ be a complex manifold with holomorphic connection $\nabla$. The \textbf{torsion}\footnote{Sometimes called torsion tensor or torsion (tensor) field.} $T$ of $\nabla$ is the collection of maps
 \begin{gather*}
  \{\tor{U}:\hvec{U}\times\hvec{U}\to\hvec{U}\mid U\subset X\text{\normalfont\ open}\}
 \end{gather*}
 defined by
 \begin{gather*}
  \tor{U}(Y,Z)\coloneqq \nab{U}_Y Z - \nab{U}_Z Y - [Y,Z]\quad\forall Y,Z\in\hvec{U}\ \forall U\subset X\text{\normalfont\ open}.
 \end{gather*}
 $\nabla$ is said to be \textbf{torsion-free} if $\tor{U}\equiv 0$ for all open subsets $U\subset X$.
\end{definition}

\begin{remark}[$T$ is a tensor]
 $T$ is bitensorial and satisfies the presheaf property.
\end{remark}

We also need to say what it means for a holomorphic connection $\nabla$ to be compatible with a holomorphic metric $h$:

\begin{definition}[Metric compatibility]
 Let $X$ be a complex manifold with holomorphic connection $\nabla$. Furthermore, let $h$ be a holomorphic metric on $X$. The \textbf{metric compatibility tensor} $\nabla h$ of $\nabla$ and $h$ is the collection of maps
 \begin{gather*}
  \{\nab{U}h:\hvec{U}^3\to\mathcal{O}_U\coloneqq\{f:U\to\C\mid f\text{\normalfont\ holomorphic}\}\mid U\subset X\text{\normalfont\ open}\}
 \end{gather*}
 defined by
 \begin{gather*}
  \nab{U}h(Z_1,Z_2,Z_3)\coloneqq Z_1(h\vert_U (Z_2,Z_3)) - h\vert_U (\nab{U}_{Z_1} Z_2, Z_3) - h\vert_U (Z_2, \nab{U}_{Z_1} Z_3)
 \end{gather*}
 for all holomorphic vector fields $Z_1,Z_2,Z_3\in\hvec{U}$ and all open subsets $U\subset X$. We say that $\nabla$ is \textbf{compatible} with $h$ if $\nab{U}h\equiv 0$ for every open subset $U\subset X$.
\end{definition}

\begin{remark}[$\nabla h$ is a tensor.]
 $\nabla h$ is tritensorial and satisfies the presheaf property.
\end{remark}

We are now ready to define the holomorphic Levi-Civita connection $\nabla^h$:

\begin{lemma}[Levi-Civita connection $\nabla^h$]\label{lem:holo_levi-civita}
 Let $X$ be a complex manifold together with a holomorphic metric $h$. Then, there exists exactly one holomorphic connection $\nabla^h$ on $X$, the holomorphic \textbf{Levi-Civita connection}, which is torsion-free and compatible with $h$.
\end{lemma}

\begin{proof}
 Take the notations from above. The proof works very similarly to the real case by exploiting the Koszul formula.\\\\
 \noindent \textbf{Uniqueness:} Let $\nabla^{h,1}$ and $\nabla^{h,2}$ be two holomorphic connections which are torsion-free and compatible with $h$. As in the real case, one can show that $\nabla^{h,1}$ and $\nabla^{h,2}$ satisfy the Koszul formula ($U\subset X$ open; $Z_1,Z_2,Z_3\in\hvec{U}$):
 \begin{align*}
  2h\vert_U (\nab{U}^{h,1}_{Z_1} Z_2, Z_3) &= Z_1(h\vert_U (Z_2,Z_3)) + Z_2(h\vert_U (Z_3,Z_1)) - Z_3 (h\vert_U (Z_1,Z_2))\\
  &\phantom{=}+ h\vert_U ([Z_1,Z_2], Z_3) - h\vert_U ([Z_2,Z_3], Z_1) + h\vert_U([Z_3,Z_1],Z_2)\\
  &= 2h\vert_U (\nab{U}^{h,2}_{Z_1} Z_2, Z_3).
 \end{align*}
 Pick any point $p\in U$ and a holomorphic chart $(z_1,\ldots, z_n):U^\prime\to\C^n$ of $X$ near $p$ such that $U^\prime\subset U$. Due to the presheaf property of $\nabla^{h,i}$, we have:
 \begin{gather*}
  \nab{U}^{h,i}_{Z_1} Z_2 (p) = \nab{U^\prime}^{h,i}_{Z_1\vert_{U^\prime}} Z_2\vert_{U^\prime} (p).
 \end{gather*}
 We now set $Z_3\vert_{U^\prime} = \partial_{z_j}$ in the Koszul formula (evaluated at $p$) and obtain using the previous equation:
 \begin{gather*}
  h\vert_{U^\prime} (\nab{U^\prime}^{h,1}_{Z_1\vert_{U^\prime}} Z_2\vert_{U^\prime}, \partial_{z_j}) (p) = h\vert_{U^\prime} (\nab{U^\prime}^{h,2}_{Z_1\vert_{U^\prime}} Z_2\vert_{U^\prime}, \partial_{z_j}) (p).
 \end{gather*}
 Since the last equation holds for every $j$, we find:
 \begin{gather*}
  \nab{U}^{h,1}_{Z_1} Z_2 (p) = \nab{U^\prime}^{h,1}_{Z_1\vert_{U^\prime}} Z_2\vert_{U^\prime} (p) = \nab{U^\prime}^{h,2}_{Z_1\vert_{U^\prime}} Z_2\vert_{U^\prime} (p) = \nab{U}^{h,2}_{Z_1} Z_2 (p).
 \end{gather*}
 As this argument can be repeated for every $p\in U$, we get:
 \begin{gather*}
  \nab{U}^{h,1}_{Z_1} Z_2 = \nab{U}^{h,2}_{Z_1} Z_2.
 \end{gather*}
 Again, we can repeat the argument for all $Z_1,Z_2\in\hvec{U}$ and any open subset $U\subset X$ implying $\nabla^{h,1} = \nabla^{h,2}$.\\\\
 \textbf{Existence:} The Koszul formula gives us an expression for
 \begin{gather*}
  h\vert_U (\nab{U}^{h}_{Z_1} Z_2, Z_3).
 \end{gather*}
 A priori, it is not cleared whether this expression completely determines $\nab{U}^h_{Z_1} Z_2$, as $U$ might be ``too large'' to admit enough linearly independent, holomorphic vector fields $Z$. However, we can always shrink $U$ by the presheaf property. Especially, if we shrink $U$ to be a chart domain, then $U$ admits enough holomorphic vector fields. Thus, we can define $\nab{U}^h_{Z_1} Z_2 (p)$ for any point $p\in U$ via the Koszul formula in a small chart near $p$. Since the Koszul formula is independent of the choice of charts, the resulting holomorphic connection $\nabla^h$ is well-defined. One easily checks that $\nabla^h$ defined this way is torsion-free and compatible with $h$.
\end{proof}

Next, we want to consider normal coordinates of a holomorphic connection. To do so, we need to define geodesics first:

\begin{definition}[Complex derivative $\nabla/dz$ along $\gamma$ and geodesics]
 Let $X$ be a complex manifold with holomorphic connection $\nabla$. Further, let $D\subset\C$ be an open and connected subset, $\gamma:D\to X$ be a holomorphic curve, and $Y$ be a holomorphic vector field along $\gamma$, i.e., a holomorphic section of the pullback bundle $\gamma^\ast T^{(1,0)}X$. Analogously to the real case, we can define the \textbf{complex derivative} $\nabla/dz$ along $\gamma$. In holomorphic charts $(z_1,\ldots, z_n):U^\prime\to\C^n$ of $X$, $\nabla/dz$ is defined by:
 \begin{gather*}
  \left(\frac{\nabla}{dz} Y\right)_l (z)\coloneqq Y_l^\prime (z) + \sum^n_{j,k=1}\gamma_j^\prime (z)\cdot Y_k (z)\cdot \Gamma^l_{jk} (\gamma (z)),
 \end{gather*}
 where $f^\prime$ denotes the usual complex derivative of a holomorphic function\linebreak $f:\C\to\C$ and $\Gamma^l_{jk}$ are the holomorphic Christoffel symbols of $\nabla$. We call $\gamma$ a \textbf{geodesic} of $\nabla$ if $\nabla/dz\ \gamma^\prime \equiv 0$.
\end{definition}

The following existence and uniqueness results regarding geodesics of $\nabla$ are reminiscent of the real case:

\begin{proposition}[Existence and uniqueness of geodesics]
 Let $X$ be a complex manifold with holomorphic connection $\nabla$. Further, let $p\in X$ and\linebreak $v\in T^{(1,0)}_p X$. Then, for every $z_0\in\C$, there exists an open and connected neighborhood $D\subset\C$ of $z_0$ and a geodesic $\gamma:D\to X$ of $\nabla$ such that $\gamma (z_0) = p$ and $\gamma^\prime (z_0) = v$. Moreover, if $\gamma_1, \gamma_2:D\to X$ are two geodesics of $\nabla$ with $\gamma_1 (z_0) = \gamma_2 (z_0) = p$ and $\gamma^\prime_1 (z_0) = \gamma^\prime_2 (z_0) = v$ for any open and connected neighborhood $D\subset\C$ of $z_0$, then $\gamma_1 \equiv \gamma_2$.
\end{proposition}

\begin{proof}
 We observe that the geodesic equation is locally a second order complex differential equation. By introducing new variables $v_i\coloneqq \gamma^\prime_i$, thus, doubling the number of equations, we can rewrite the second order differential equation into a first order one. Since the geodesic equation has no explicit time-dependence ($z$-dependence), we can interpret the first order differential equation obtained this way as the integral curve equation of a holomorphic vector field, namely the geodesic vector field. The existence and uniqueness results for geodesics now follow from the existence and uniqueness results for integral curves of holomorphic vector fields (cf. Proposition \autoref{prop:holo_traj}).\\
 There is also an alternative way to prove uniqueness: Write the geodesic equation in holomorphic charts near $p$ and expand the coordinates of $\gamma_1$, $\gamma_2$ in a power series around $z_0$. This way, the geodesic equation becomes a recursive formula for the coefficients of the power series. Furthermore, we realize that this recursive formula completely determines all coefficients if we fix the first two terms in each power series. Hence, fixing $\gamma_1 (z_0) = \gamma_2 (z_0) = p$ and $\gamma^\prime_1 (z_0) = \gamma^\prime_2 (z_0) = v$ uniquely determines the coefficients of the power series. The rest now follows from the identity theorem.
\end{proof}

\begin{remark}[Geodesics depend holomorphically on initial values]\label{rem:geodesic_initial_value}
 Proposition \autoref{prop:holo_traj} also shows that a geodesic $\gamma$ of a holomorphic connection $\nabla$ with $\gamma (z_0) = p$ and $\gamma^\prime (z_0) = v$ depends holomorphically on $z_0$, $p$, and $v$.
\end{remark}

We can now use the geodesics of $\nabla$ to define the exponential map of $\nabla$:

\begin{definition}[Exponential map]\label{def:holo_exp_map}
 Let $X$ be a complex manifold with holomorphic connection $\nabla$  and let $p\in X$ be any point. The \textbf{exponential map} of $\nabla$, $\exp_p:V\to X$ with a suitable open neighborhood $V\subset T^{(1,0)}_p X$ of $0$, is defined by:
 \begin{gather*}
  \exp_p (v) \coloneqq \gamma^v_p (1),
 \end{gather*}
 where $\gamma^v_p$ is a geodesic of $\nabla$ satisfying $\gamma^v_p (0) = p$ and $\gamma^{v\prime}_p (0) = v$.
\end{definition}

\begin{remark}[Exponential map] We add two remarks concerning the exponential map:
 \begin{enumerate}
  \item $V\subset T^{(1,0)}_p X$ is chosen small enough such that $\gamma^v_p (1)$ is uniquely defined and unaffected by monodromy effects (cf. \autoref{sec:HHS}, in particular Proposition \autoref{prop:holo_traj} and Example \autoref{ex:holo_cen_prob}). For instance, choose $V\subset T^{(1,0)}_pX$ such that the domain of the geodesic $\gamma^v_p:D\to X$ can be chosen to be $D = \{z\in\C\mid |z|<2\}$ for every $v\in V$.
  \item By Remark \autoref{rem:geodesic_initial_value}, the exponential map $\exp_p$ of a holomorphic connection $\nabla$ is holomorphic.
 \end{enumerate}
\end{remark}

We can now shrink the domain of the exponential map to obtain a biholomorphism:

\begin{proposition}
 Let $X$ be a complex manifold with holomorphic connection $\nabla$  and let $p\in X$ be any point. Then, there exists an open neighborhood $V\subset T^{(1,0)}_p X$ of $0$ such that $\exp_p: V\to \exp_p (V)$ is a biholomorphism.
\end{proposition}

\begin{proof}
 Apply the holomorphic inverse function theorem to $d\exp_{p,0} = \text{id}_{T^{1,0}_p X}$.
\end{proof}

Lastly, we use this biholomorphism to define normal coordinates:

\begin{definition}[Normal coordinates]
 Let $X$ be a complex manifold with holomorphic connection $\nabla$  and let $p\in X$ be any point. Further, choose a\linebreak $\C$-linear isomorphism $l:\C^n\to T^{(1,0)}_p X$. We call the coordinates $(z_1,\ldots, z_n)$ of the holomorphic chart $l^{-1}\circ \exp_p^{-1}$ \textbf{holomorphic normal coordinates} of $\nabla$ near $p$. If $\nabla = \nabla^h$ is the Levi-Civita connection of a holomorphic metric $h$ on $X$, we additionally require that the vectors $v_1\coloneqq l (\hat e_1),\ldots, v_n\coloneqq l (\hat e_n)$\linebreak ($\hat e_1,\ldots, \hat e_n$ denotes the standard basis of $\C^n$) satisfy:
 \begin{gather*}
  h (v_i, v_j) = \delta_{ij}.
 \end{gather*}
\end{definition}

\begin{remark}
 Note that in the complex case all non-degenerate, symmetric bilinear forms are isomorphic, while in the real case two non-degenerate, symmetric bilinear forms are isomorphic if and only if they have the same signature.
\end{remark}

Holomorphic normal coordinates satisfy properties similar to their real counterparts:

\pagebreak

\begin{proposition}[Properties of holomorphic normal coordinates]
 Let $X$ be a complex manifold with holomorphic metric $h$ and corresponding Levi-Civita connection $\nabla^h$. Further, let $p\in X$ be any point and $\phi = (z_1,\ldots, z_n):U\to\C^n$ be holomorphic normal coordinates of $\nabla^h$ near $p$. In the coordinates $(z_1,\ldots, z_n)$, we have:
 \begin{enumerate}[label = (\arabic*)]
  \item $h_{ij} (p) = \delta_{ij}$
  \item $\partial_{z_k}h_{ij} (p) = 0$
  \item $\Gamma^k_{ij} (p) = 0$
 \end{enumerate}
\end{proposition}

\begin{proof}
 We only show (3), as (1) is obvious and (2) follows from (3) by\linebreak exploiting the vanishing of the metric compatibility tensor. Define for\linebreak $x = (x_1,\ldots, x_n)\in\C^n$ the curve $\gamma (z)\coloneqq \phi^{-1}(z\cdot x)$. By definition of holomorphic normal coordinates, $\gamma$ is a geodesic of $\nabla^h$ with $\gamma^\prime (0) = d\phi^{-1}_0 (x)$. In coordinates $(z_1,\ldots, z_n)$, the geodesic equation reads:
 \begin{gather*}
  \gamma^{\prime\prime}_{j}(z) + \sum^n_{k,l = 1} \Gamma^j_{kl} (\gamma (z))\gamma^\prime_k (z)\gamma^\prime_l (z) = 0
 \end{gather*}
 For $z = 0$, this gives:
 \begin{gather*}
  \sum^n_{k,l = 1} \Gamma^j_{kl} (p)x_kx_l = 0
 \end{gather*}
 The last equation is true for any $x\in\C^n$, thus, it enforces $\Gamma^j_{kl} (p) + \Gamma^j_{lk}(p) = 0$. Since $\nabla^h$ is symmetric, $\Gamma^j_{kl} (p)$ is symmetric in $k$ and $l$. This implies $2\Gamma^j_{kl} (p) = 0$ concluding the proof.
\end{proof}

The final goal of \autoref{app:holo_connection} is to relate the holomorphic Levi-Civita connection $\nabla^{h}$ to the standard Levi-Civita connections $\nabla^{h_R}$ and $\nabla^{h_I}$. In particular, we aim to relate the holomorphic normal coordinates of $\nabla^h$ to the standard normal coordinates of $\nabla^{h_R}$ and $\nabla^{h_I}$ for holomorphic metrics $h = h_R + ih_I$. To achieve that, we first prove the following proposition:

\begin{proposition}\label{prop:h_R_in_coordinates}
 Let $X$ be a complex manifold with holomorphic metric $h = h_R + ih_I$ and corresponding Levi-Civita connections $\nabla^h$, $\nabla^{h_R}$, and $\nabla^{h_I}$. Further, let $p\in X$ be any point and $(z_1 = x_1 + iy_1,\ldots, z_n = x_n + iy_n):U\to\C^n$ be holomorphic normal coordinates of $\nabla^h$ near $p$. Then:
 \begin{enumerate}[label = (\arabic*)]
  \item $h_{R,p} (\partial_{x_j,p}, \partial_{x_k,p}) = -h_{R,p} (\partial_{y_j,p}, \partial_{y_k,p}) = \delta_{jk},\quad h_{I,p} (\partial_{x_j,p}, \partial_{y_k,p}) = \delta_{jk}$\\
  $h_{R,p} (\partial_{x_j,p}, \partial_{y_k,p}) \, = h_{I,p} (\partial_{x_j,p}, \partial_{x_k,p})\ \ \, = h_{I,p} (\partial_{y_j,p}, \partial_{y_k,p})\qquad\ \ \, = 0$
  \item Derivatives of $h_R$ and $h_I$ vanish at $p$ in coordinates $(x_1,\ldots, x_n, y_1,\ldots, y_n)$.
  \item All Christoffel symbols of $\nabla^{h_R}$ and $\nabla^{h_I}$ vanish at $p$ in coordinates\linebreak $(x_1,\ldots, x_n, y_1,\ldots, y_n)$.
 \end{enumerate}
\end{proposition}

\begin{proof}
 (1) follows from (1) of the previous proposition:
 \begin{align*}
  h_{R,p} + i h_{I,p} &= h_p = \sum^n_{j = 1} dz^2_{j,p}\\
  &= \sum^n_{j = 1}\left(dx^2_{j,p} - dy^2_{j,p}\right) + i\sum^n_{j = 1}\left(dx_{j,p}\otimes dy_{j,p} + dy_{j,p}\otimes dx_{j,p}\right).
 \end{align*}
 (2) follows from (2) of the previous proposition and the fact that the components $h_{ij}$ are holomorphic, i.e., $\partial_{\bar{z}_k} h_{ij} (p) = 0$.\\
 (3) follows from (2) of the proposition at hand and the Koszul formula for Christoffel symbols of standard Levi-Civita connections.
\end{proof}

The last proposition implies that the Levi-Civita connections $\nabla^{h_R}$ and $\nabla^{h_I}$ describe the same connection.

\begin{corollary}[$\nabla^{h_R} = \nabla^{h_I}$]
 Let $X$ be a complex manifold with holomorphic metric $h = h_R + ih_I$. Then, the (standard) Levi-Civita connections $\nabla^{h_R}$ and $\nabla^{h_I}$ coincide, i.e., $\nabla^{h_R} = \nabla^{h_I}$.
\end{corollary}

\begin{proof}
 This is an immediate consequence of the previous proposition: For every point $p\in X$, there are coordinates in which the Christoffel symbols of $\nabla^{h_R}$ and $\nabla^{h_I}$ agree at $p$, hence, $\nabla^{h_R} = \nabla^{h_I}$ by definition of the Christoffel symbols.
\end{proof}

We want to show now that $\nabla^{h_R} = \nabla^{h_I}$ is the connection associated with the holomorphic Levi-Civita connection $\nabla^h$. To do that, we need to investigate how the complex structure $J$ interacts with the connection $\nabla^{h_R}$. It turns out that $\nabla^{h_R}$ is compatible with $J$:

\begin{proposition}[$\nabla^{h_R}$ compatible with $J$]\label{prop:levi_compatible_with_J}
 Let $X$ be a complex manifold with complex structure $J$ and holomorphic metric $h = h_R + ih_I$. Then, the connection $\nabla^{h_R}$ is compatible with $J$.
\end{proposition}

\begin{proof}
 We need to show:
 \begin{gather*}
  (\nabla^{h_R}_Y J)Z = \nabla^{h_R}_Y (JZ) - J\nabla^{h_R}_Y Z = 0\quad\forall Y,Z\in\Gamma (TX).
 \end{gather*}
 The object $(\nabla^{h_R}_Y J) (Z)$ is tensorial in $Y$ and $Z$, hence, it suffices to evaluate $\nabla^{h_R} J$ for the coordinate vector fields of the holomorphic normal coordinates $(z_1 = x_1 + iy_1,\ldots, z_n = x_n + iy_n)$ near any point $p\in X$. By Proposition \autoref{prop:h_R_in_coordinates}, we find:
 \begin{align*}
  \nabla^{h_R}_{\partial_{x_i}}\left( J \partial_{x_j}\right) (p) &= \nabla^{h_R}_{\partial_{x_i}} \partial_{y_j} (p) = 0,\quad \nabla^{h_R}_{\partial_{x_i}}\left( J \partial_{y_j}\right) = \ldots\\
  J\nabla^{h_R}_{\partial_{x_i}} \partial_{x_j} (p) &= 0,\quad J\nabla^{h_R}_{\partial_{x_i}} \partial_{y_j} (p) =\ldots
 \end{align*}
 Hence, $\nabla^{h_R} J = 0$.
\end{proof}

$\nabla^{h_R}$ is not just compatible with $J$, it also satisfies the following useful formula:

\begin{proposition}\label{prop:helpful_formula}
 Let $X$ be a complex manifold with holomorphic metric $h = h_R + i h_I$ and Levi-Civita connection $\nabla^{h_R}$. Then:
 \begin{gather*}
  \nabla^{h_R}_{JY} Z = i\nabla^{h_R}_Y Z\quad\forall Y\in\Gamma (TU)\ \forall Z\in\Gamma (T^{(1,0)}U)\ \forall U\subset X\text{ open}.
 \end{gather*}
\end{proposition}

\begin{proof}
 Take the notations from above and pick any point $p\in U$. Let\linebreak $(z_1 = x_1 + iy_1,\ldots, z_n = x_n + iy_n)$ be holomorphic normal coordinates of $\nabla^h$ near $p$ and write
 \begin{gather*}
  Z = \sum^n_{j = 1} c_j\pa{z_j}
 \end{gather*}
 for some locally defined holomorphic functions $c_j$. Then:
 \begin{align*}
  \nabla^{h_R}_{JY} Z (p) &= \sum^n_{j = 1} c_j (p)\nabla^{h_R}_{JY}\pa{z_j} (p) + dc_j (JY) (p)\pa{z_j,p}\\
  &= \sum^n_{j = 1} c_j (p)\nabla^{h_R}_{JY}\pa{z_j} (p) + i dc_j (Y) (p)\pa{z_j,p}\\
  &= \sum^n_{j = 1} idc_j (Y) (p)\pa{z_j,p}\\
  &= i\sum^n_{j = 1} c_j (p) \nabla^{h_R}_{Y}\pa{z_j} (p) + dc_j (Y) (p)\pa{z_j,p} = i\nabla^{h_R}_Y Z (p),
 \end{align*}
 where we used that $\nabla^{h_R}$ is $\C$-linear in both components, tensorial in the first component, and $\nabla^{h_R}_{\pa{x_i}} \pa{x_j} (p) = \nabla^{h_R}_{\pa{x_i}} \pa{y_j} (p) = \ldots = 0$ (cf. Proposition \autoref{prop:h_R_in_coordinates}).
\end{proof}

\begin{remark}
 The formula in Proposition \autoref{prop:helpful_formula} not only applies to $\nabla^{h_R}$, but to any connection $\nabla$ associated with a holomorphic connection. Indeed, if\linebreak $p\in U\subset X$ is any point and $Y\in\Gamma (TU)$ is a smooth real vector field, then there exists a holomorphic vector field $Y^\prime$ (possibly only defined on a neighborhood of $p$) such that $Y^\prime (p) = Y(p) - iJ_pY (p)$. Exploiting Property (iii) from Lemma \autoref{lem:associated_connection} then yields:
 \begin{gather*}
  0 = \nabla_{\overline{Y^\prime}}Z (p) = \nabla_Y Z (p) +i \nabla_{JY}Z (p)
 \end{gather*}
 for any holomorphic vector field $Z\in\hvec{U}$.
\end{remark}

We now have all tools at hand to show that $\nabla^{h_R}$ is the connection associated with $\nabla^h$:

\begin{lemma}[$\nabla^{h_R}$ associated with $\nabla^h$]\label{lem:associated_levi}
 Let $X$ be a complex manifold with holomorphic metric $h = h_R + ih_I$. Then, the standard Levi-Civita connection $\nabla^{h_R} = \nabla^{h_I}$ is the connection associated with the holomorphic Levi-Civita connection $\nabla^h$.
\end{lemma}

\begin{proof}
 For real vector fields $Z_1,Z_2,Z_3\in\Gamma (TX)$, we can express the terms
 \begin{gather*}
  2h_R \left(\nabla^{h_R}_{Z_1} Z_2, Z_3\right)\quad\text{and}\quad 2h_I \left(\nabla^{h_I}_{Z_1} Z_2, Z_3\right)
 \end{gather*}
 with the help of the Koszul formula as in the proof of Lemma \autoref{lem:holo_levi-civita}. We realize that the Koszul formula also holds for complex vector fields $Z_1,Z_2,Z_3\in\Gamma (T_\C X)$ if we complexify the Koszul formula by $\C$-linearity in $Z_1,Z_2,Z_3$. In particular, for holomorphic vector fields $Z_1,Z_2,Z_3\in\Gamma (T^{(1,0)}X)$, this gives:
 \begin{align*}
  2h (\nabla^{h_R}_{Z_1} Z_2, Z_3) &= 2h_R (\nabla^{h_R}_{Z_1} Z_2, Z_3) + i 2h_I (\nabla^{h_R}_{Z_1} Z_2, Z_3)\\
  &= 2h_R (\nabla^{h_R}_{Z_1} Z_2, Z_3) + i 2h_I (\nabla^{h_I}_{Z_1} Z_2, Z_3)\\
  &= (\text{Koszul formula for $h_R$}) + i (\text{Koszul formula for $h_I$})\\
  &= \text{Koszul formula for $h$} = 2h (\nabla^h_{Z_1} Z_2, Z_3).
 \end{align*}
 Thus, we have:
 \begin{gather}\label{eq:tired}
  h(\nabla^{h}_{Z_1} Z_2 - \nabla^{h_R}_{Z_1} Z_2, Z_3) = 0.
 \end{gather}
 Recall that, by Proposition \autoref{prop:levi_compatible_with_J}, $\nabla^{h_R}$ is compatible with the complex structure $J$ of $X$. In particular, $\nabla^{h_R}$ preserves the $(1,0)$-type of the second component (cf. Proposition \autoref{prop:nabla_compatible}). Thus, the vector field $\nabla^{h}_{Z_1} Z_2 - \nabla^{h_R}_{Z_1} Z_2$ assumes values in the bundle $T^{(1,0)}X$. Since $h$ is non-degenerate on $T^{(1,0)}X$, \autoref{eq:tired} proves that the connection $\nabla^{h_R}$ induces $\nabla^h$. Clearly, $\nabla^{h_R}$ is real. Thus, it only remains to be shown that $\nabla^{h_R}$ satisfies:
 \begin{gather*}
  \nabla^{h_R}_{\overline{Y}}Z = 0\quad\forall Y,Z\in\hvec{U}\ \forall U\subset X\text{ open}.
 \end{gather*}
 Let $U\subset X$ be an open subset and $Y,Z$ be two holomorphic vector field on $U$. Then, there exists a unique real vector field $\hat Y$ on $U$ such that $Y = \hat Y - iJ\hat Y$. Proposition \autoref{prop:helpful_formula} now yields the desired property concluding the proof:
 \begin{gather*}
  \nabla^{h_R}_{\overline{Y}}Z = \nabla^{h_R}_{\hat Y + iJ\hat Y}Z = \nabla^{h_R}_{\hat Y}Z + i\nabla^{h_R}_{J\hat Y}Z = \nabla^{h_R}_{\hat Y}Z + i^2\nabla^{h_R}_{\hat Y}Z = 0.
 \end{gather*}
\end{proof}

To conclude \autoref{app:holo_connection}, we examine the relation between holomorphic normal coordinates of $\nabla^h$ and normal coordinates of $\nabla^{h_R}$. For this, we need to link holomorphic geodesics of $\nabla^h$ to geodesics of $\nabla^{h_R}$:

\begin{proposition}\label{prop:link_geodesics}
 Let $X$ be a complex manifold with holomorphic metric\linebreak $h = h_R + ih_I$  and Levi-Civita connections $\nabla^h$ and $\nabla^{h_R}$. Further, let\linebreak $U\coloneqq [t_0, t_1] + i[s_0, s_1]$ be a domain in $\C$ and $\gamma:U\to X$ be a holomorphic curve. Define the curves $\gamma_s:[t_0,t_1]\to X$ and $\gamma_t: [s_0, s_1]\to X$ by $\gamma_s (t) \coloneqq \gamma (t+is) \eqqcolon \gamma_t (s)$. Then, $\gamma$ is a (holomorphic) geodesic of $\nabla^h$ iff $\gamma_s$ is a geodesic of $\nabla^{h_R} = \nabla^{h_I}$ for every $s\in[s_0, s_1]$ iff $\gamma_t$ is a geodesic of $\nabla^{h_R} = \nabla^{h_I}$ for every $t\in[t_0, t_1]$.
\end{proposition}

\begin{proof}
 This statement is mostly a consequence of Proposition \autoref{prop:levi_compatible_with_J} and \autoref{prop:helpful_formula} as well as Lemma \autoref{lem:associated_levi}: For $z = t+is$, we can write:
 \begin{gather*}
  \gamma^\prime (z) = \frac{1}{2}\left(\pa{t}\gamma (z) - iJ\pa{t}\gamma (z)\right) = \frac{1}{2}\left(-J\pa{s}\gamma (z) - i\pa{s}\gamma (z)\right).
 \end{gather*}
 Hence:
 \begin{align*}
  \frac{\nabla^h}{dz}\gamma^\prime &= \nabla^h_{\gamma^\prime} \gamma^\prime = \nabla^{h_R}_{\gamma^\prime} \gamma^\prime = \nabla^{h_R}_{1/2 (\pa{t}\gamma - iJ\pa{t}\gamma)} \gamma^\prime = \nabla^{h_R}_{\pa{t}\gamma}\gamma^\prime\\
  &= \frac{1}{2}\left(\nabla^{h_R}_{\pa{t}\gamma} \pa{t}\gamma - iJ\left(\nabla^{h_R}_{\pa{t}\gamma} \pa{t}\gamma\right)\right)\\
  &= \frac{1}{2}\left(\frac{\nabla^{h_R}}{dt}\frac{d\gamma_s}{dt} - iJ\left(\frac{\nabla^{h_R}}{dt}\frac{d\gamma_s}{dt}\right)\right).
 \end{align*}
 A similar expression can be found for $\gamma_t$ concluding the proof.
\end{proof}

An immediate consequence of Proposition \autoref{prop:link_geodesics} is that the exponential maps of $\nabla^h$ and $\nabla^{h_R}$ coincide:

\begin{corollary}\label{cor:exp}
 Let $X$ be a complex manifold with holomorphic metric\linebreak $h = h_R + ih_I$  and Levi-Civita connections $\nabla^h$ and $\nabla^{h_R}$. Further, let $p\in X$ be a point and $\exp^h_p$ and $\exp^{h_R}_p$ be the exponential maps of $\nabla^h$ and $\nabla^{h_R}$, respectively. Take a vector $v\in T^{(1,0)}_p X$ such that $\exp^h_p (v)$ exists and write it as $v = 1/2 (\hat v - iJ_p\hat v)$ for a unique vector $\hat v\in T_pX$. Then, $\exp^{h_R}_p (\hat v)$ exists and we have:
 \begin{gather*}
  \exp^{h_R}_p (\hat v) = \exp^h_p (v).
 \end{gather*}
\end{corollary}

\begin{proof}
 Take the notations from above. $\exp^h_p (v)$ is defined as $\gamma (1)$, where $\gamma$ is a holomorphic geodesic of $\nabla^h$ satisfying $\gamma (0) = p$ and $\gamma^\prime (0) = v$. By Proposition \autoref{prop:link_geodesics}, $\gamma (t)$, $t\in\R$, is a geodesic of $\nabla^{h_R}$ satisfying $\gamma (0) = p$ and $\dot\gamma (0) = \hat v$. Since $\gamma (1)$ exists, $\exp^{h_R}_p (\hat v)$ exists as well and we have:
 \begin{gather*}
  \exp^{h_R}_p (\hat v) = \gamma (1) = \exp^h_p (v).
 \end{gather*}
\end{proof}

\begin{remark}
 Observe that Proposition \autoref{prop:link_geodesics} and Corollary \autoref{cor:exp} are still true if we replace $\nabla^h$ by any holomorphic connection and $\nabla^{h_R}$ by its associated connection. Indeed, Corollary \autoref{cor:exp} follows from Proposition \autoref{prop:link_geodesics}, while Proposition \autoref{prop:link_geodesics} just exploits compatibility with $J$ and the formula in Proposition \autoref{prop:helpful_formula} as well as the fact that $\nabla^{h_R}$ induces $\nabla^h$. These results are also valid for associated connections.
\end{remark}

At last, we arrive at the following lemma:

\begin{lemma}\label{lem:normal_of_h_and_h_R}
 Let $X$ be a complex manifold with holomorphic metric\linebreak $h = h_R + ih_I$. Let $p\in X$ be any point. Then, holomorphic normal coordinates\linebreak $(z_1 = x_1 + iy_1,\ldots, z_n = x_n + iy_n)$ of $\nabla^h$ near $p$ give rise to normal coordinates $(x_1,\ldots, y_n)$ of $\nabla^{h_R}$ near $p$.
\end{lemma}

\begin{proof}
 Combine all previous results from this section, in particular use Corollary \autoref{cor:exp}. Note that the same result is only true for $h_I$ after applying a linear transformation, since
 \begin{gather*}
  h_{R,p} = \sum^n_{j = 1} dx^2_{j,p} - dy^2_{j,p}
 \end{gather*}
 is in standard form at $p$ in the coordinates $(x_1,\ldots, y_n)$, while the same is not true for
 \begin{gather*}
  h_{I,p} = \sum^n_{j=1} dx_{j,p}\otimes dy_{j,p} + dy_{j,p}\otimes dx_{j,p}.
 \end{gather*}
\end{proof}

 \chapter[Complex Lie Groups, Algebras, and Coadjoint Orbits]{Complex Lie Groups, Algebras, and Coadjoint Orbits\chaptermark{Complex Lie Groups}}
 \chaptermark{Complex Lie Groups}
 \label{app:complex_lie_groups}
 This part shall serve as an introduction to complex Lie groups. Since most notions of \autoref{sec:lie_groups}, which offers a general introduction to Lie groups, simply transfer to the complex case (conjugation, (co)adjoint action, exponential map, and so on), we only focus on those concepts that have no real analogue or differ in one way or the other from the real case. Particularly, we are concerned with the complex and real structures on Lie groups, the various definitions of the complex Lie algebra, and how the complex and real structures of the group descend to the (co)adjoint orbits.\\
We start by recalling the definition of a complex Lie group:

\begin{definition}\label{def:comp_group}
 $G$ is called a \textbf{complex Lie group} if $G$ is a smooth Lie group equipped with an almost complex structure $I$ such that the left and right multiplication are $I$-holomorphic.
\end{definition}

\begin{remark}
 It directly follows from the definition that the adjoint action commutes with $I_e$, the almost complex structure at the neutral element $e\in G$. This implies that the Lie algebra $\mathfrak{g}$ of $G$ is a complex Lie algebra:
 \begin{gather*}
  [I_e u, v] = [u, I_e v] = I_e [u,v]\quad\forall u,v\in\mathfrak{g}.
 \end{gather*}
 One can deduce from the last equation that $I$ is integrable, i.e., a complex structure (cf. \cite{Salamon2019}). Moreover, it follows from the holomorphic implicit function theorem that the inversion $G\to G, g\mapsto g^{-1}$ is holomorphic.
\end{remark}

The complex Lie groups we are interested also carry a real structure:

\begin{definition}\label{def:group_real_str}
 Let $G$ be a complex Lie group. A \textbf{real structure} on $G$ is an antiholomorphic involution $\sigma:G\to G$ which is also a group homomorphism. Its fixed point set $G_\R\coloneqq\Fix\sigma$ is called \textbf{real form}. We say $G_\R$ is nice if $G_\R$ meets every connected component of $G$. A complex Lie group is called \textbf{reductive} if it admits a nice compact real form.
\end{definition}

As one would expect, $G_\R\subset G$ is a (real) Lie subgroup of half its (real) dimension:

\begin{proposition}[$G_\R\subset G$ Lie subgroup]\label{prop:lie_real_form}
 Let $G$ be a complex Lie group with real structure $\sigma$. Then, its real form $G_\R\subset G$ is an (real) embedded Lie subgroup with $\dim_\R G = 2\dim_\R G_\R$. In particular, $G_\R$ is not empty.
\end{proposition}

\begin{proof}
 We first prove that $G_\R\subset G$ is an embedded Lie subgroup. By Cartan's subgroup theorem, it suffices to show that $G_\R\subset G$ is a closed subgroup. $G_\R\subset G$ is a subgroup, since it is the fixed point set of a group homomorphism. $G_\R\subset G$ is also closed, as it is the preimage of the closed set $\{e\}$ under the continuous map $g\mapsto \sigma (g)g^{-1}$.\\
 To conclude the proof, we note that $G_\R$ is not empty, because $G_\R$ is a group and contains at least the neutral element $e$. Thus, it is also a real form in the sense of Definition \autoref{def:real_form}, hence, $\dim_\R G = 2\dim_\R G_\R$.
\end{proof}

One example of a complex Lie group is given by $\GL (n,\C)$. It possesses several real forms, for instance $\GL(n,\R)$ and $\U(n)$. The real structures corresponding to the real forms $\GL(n,\R)$ and $\U(n)$ are $\GL (n,\C)\to\GL (n,\C)$, $A\mapsto\bar{A}$ and $\GL (n,\C)\to\GL (n,\C)$, $A\mapsto (A^{-1})^\ast = \overline{(A^{-1})}^t$, respectively. We see that $\GL (n,\C)$ admits the compact real form $\U (n)$, even though $\GL (n,\C)$ itself is not compact. This is true for most reductive groups. Indeed, if a connected complex Lie group was compact, its adjoint representation would be trivial due to the maximum principle implying that the complex group in question would be Abelian. In this thesis, we are not interested in Abelian groups, as their (co)adjoint orbits are just points.

\begin{remark}
 From now on, we always denote a complex Lie group by $G$, its neutral element by $e$, its complex structure by $I$, a real structure on $G$ by $\sigma$, and the corresponding real form by $G_\R$.
\end{remark}

Next, let us consider the Lie algebras of complex Lie groups in more detail. For a (real) Lie group, we think of its Lie algebra as the tangent space at the neutral element which we usually identify with the space of left-invariant\footnote{Less commonly, with the space of right-invariant vector fields which leads to a change of sign in the Lie bracket.} vector fields to equip the Lie algebra with a Lie bracket. For a complex Lie group $G$, however, there are multiple ways to define its (complex) Lie algebra: As usual, we can take the tangent space $T_eG$ and identify it with the space of left-invariant vector fields, where the complex structure of the Lie algebra is $I_e$. Alternatively, we might think of the Lie algebra as the holomorphic tangent space $T^{(1,0)}_e G$ and identify it with the space of left-invariant, holomorphic vector fields, where the complex structure is simply given by multiplication with the imaginary unit $i$. Additionally, if $G$ admits a real form $G_\R$, we may complexify the Lie algebra of $G_\R$ and take it to be the complex Lie algebra of $G$. The following proposition states that all these models are isomorphic:

\begin{proposition}[Equivalent models of complex Lie algebras]\label{prop:models_of_Lie_algebras}
 Let $G$ be a complex Lie group with real form $G_\R\subset G$. Denote the Lie algebras by\linebreak $\mathfrak{g}_\R\coloneqq T_e G_\R$, $\mathfrak{g}_\C\coloneqq\mathfrak{g}_\R\otimes_\R\C$, $\mathfrak{g}\coloneqq T_e G$, and $\mathfrak{g}^{(1,0)}\coloneqq T^{(1,0)}_e G$. Then, the maps
 \begin{gather*}
  h_1:\mathfrak{g}_\C\to\mathfrak{g},\, u+iv\mapsto u+I_e v\quad \text{and}\quad h_2:\mathfrak{g}\to\mathfrak{g}^{(1,0)},\, w\mapsto \frac{1}{2}\left(w-iI_e w\right)
 \end{gather*}
 are complex Lie algebra isomorphisms.
\end{proposition}

\begin{proof}
 It is widely known and easy to check that $h_1:(\mathfrak{g}_\C, i)\to (\mathfrak{g}, I_e)$ and $h_2:(\mathfrak{g}, I_e)\to (\mathfrak{g}^{(1,0)}, I_e = i)$ are $\C$-linear bijections. Hence, we only need to check that $h_1$ and $h_2$ respect the Lie bracket. We compute:
 \begin{align*}
  [h_1 (u_1 + iv_1), h_1 (u_2 + iv_2)] &= [u_1 + I_e v_1, u_2 + I_e v_2]\\
  &= [u_1, u_2] - [v_1, v_2] + I_e\left([u_1, v_2] + [v_1, u_2]\right)\\
  &= h_1\left([u_1, u_2] - [v_1, v_2] + i\left([u_1, v_2] + [v_1, u_2]\right)\right)\\
  &= h_1 \left([u_1 + iv_1, u_2 + iv_2]\right)\quad\forall u_1,u_2,v_1,v_2\in\mathfrak{g}_\R
 \end{align*}
 and:
 \begin{align*}
  [h_2 (w_1), h_2 (w_2)] &= \frac{1}{4}[w_1 - iI_e w_1, w_2 - iI_e w_2]\\
  &= \frac{1}{4}\left([w_1, w_2] - [I_e w_1, I_e w_2] -i\left([w_1, I_e w_2] + [I_e w_1, w_2]\right)\right)\\
  &= \frac{1}{2}\left([w_1, w_2] - iI_e [w_1, w_2]\right)\\
  &= h_2 ([w_1, w_2])\quad\forall w_1, w_2\in\mathfrak{g},
 \end{align*}
 where we used $[I_e u, v] = [u, I_e v] = I_e [u,v]$ and $[i u, v] = [u, i v] = i [u,v]$.
\end{proof}

\begin{remark}
 Often, we mean $T_e G$ when we talk about the Lie algebra of $G$ and denote it simply by $\mathfrak{g}$. Sometimes, however, we implicitly use a different model for the complex Lie algebra and do not change the symbol, so caution is advised.
\end{remark}

In \autoref{sec:duality} and especially in \autoref{app:hyp_orb}, we consider semisimple complex Lie groups:

\begin{definition}[Semisimple]\label{def:semisimple}
 A real Lie algebra $(\mathfrak{g},[\cdot,\cdot])$ is called \textbf{semisimple} if it contains no proper non-zero Abelian ideals. We say a Lie group $G$ is semisimple if its Lie algebra is semisimple.
\end{definition}

\begin{remark}\ 
 \begin{enumerate}
  \item It is a famous result from Lie group theory, known as Cartan's criterion, that a Lie algebra is semisimple if and only if its Killing form is non-degenerate (cf. \autoref{sec:lie_groups}).
  \item One easily verifies with the help of Cartan's criterion that a complex Lie group $G$ with real form $G_\R$ is semisimple if and only if its real form $G_\R$ is semisimple.
 \end{enumerate}
\end{remark}

We now turn our attention to (co)adjoint orbits of complex Lie groups. We only consider adjoint orbits in the upcoming discussions, since adjoint and coadjoint orbits are isomorphic in the relevant cases\footnote{By this, we mean the case of complex reductive groups. The adjoint and coadjoint orbits of such groups are isomorphic due to the existence of an $\Ad$-invariant, non-degenerate scalar product.}. Recall that there are two ways to think of orbits of group actions: We can view them as group quotients or as immersed submanifolds. First, let us investigate the quotient model for complex adjoint orbits more closely. To do that, we consider complex stabilizer groups:

\begin{proposition}\label{prop:comp_stab}
 Let $G$ be a complex Lie group with real structure $\sigma$\linebreak and real form $G_\R\subset G$. Further, let $w_0\in\mathfrak{g}_\R\subset\mathfrak{g}$ be a point and\linebreak $G_{w_0}\coloneqq\{g\in G\mid \Ad (g)w_0 = w_0\}$ be the stabilizer of $w_0$. Then, $G_{w_0}\subset G$ is an embedded complex Lie subgroup with real structure $\sigma\vert_{G_{w_0}}$ and real form
 \begin{gather*}
  \Fix \sigma\vert_{G_{w_0}} = G_{w_0}\cap G_\R = (G_\R)_{w_0}\coloneqq\{g\in G_\R\mid \Ad (g) w_0 = w_0\}.
 \end{gather*}
\end{proposition}

\begin{proof}
 As in the proof of Lemma \autoref{lem:immersed_orbit}, we find that $G_{w_0}\subset G$ is an embedded Lie subgroup and its Lie algebra is $\ker d\ev_{w_0,e} = \ker\ad_{w_0}$. Since $\ad_{w_0}$ commutes with $I_e$, $\ker\ad_{w_0}$ is closed under $I_e$. Hence, the tangent spaces of $G_{w_0}$ are complex subspaces of tangent spaces of $G$ turning $G_{w_0}$ into a complex Lie subgroup.\\
 To show that $\sigma\vert_{G_{w_0}}$ is a real structure on $G_{w_0}$, we just need to prove that $G_{w_0}$ is $\sigma$-invariant, i.e., $\sigma (G_{w_0}) = G_{w_0}$. This is true, since:
 \begin{gather*}
  \Ad (\sigma (g))w_0 = d\sigma_e\left(\Ad (g) (d\sigma_e w_0)\right) = d\sigma_e\left(\Ad (g) w_0\right) = d\sigma_e w_0 = w_0\quad\forall g\in G_{w_0},
 \end{gather*}
 where we used $w_0\in\mathfrak{g}_\R = \Fix d\sigma_e$. The last equation in Proposition \autoref{prop:comp_stab} immediately follows from $\Fix\sigma\vert_{G_{w_0}} = G_{w_0}\cap\Fix\sigma$.
\end{proof}

The adjoint orbit through $w_0$ is now the quotient $G/G_{w_0}$. Since both $G$ and $G_{w_0}$ carry a complex and real structure, the quotient does so as well:

\begin{proposition}\label{prop:comp_quotient}
 Let $G$ be a complex Lie group with real structure $\sigma$\linebreak and real form $G_\R\subset G$. Further, let $w_0\in\mathfrak{g}_\R\subset\mathfrak{g}$ be a point and\linebreak $G_{w_0}$ be the stabilizer of $w_0$. Then, $G/G_{w_0}$ is a complex manifold with holomorphic submersion $\pi:G\to G/G_{w_0}$. Moreover, $\sigma$ descends to the antiholomorphic involution $[\sigma]$ on $G/G_{w_0}$.
\end{proposition}

\begin{proof}
 $G_{w_0}\subset G$ is a smooth embedded Lie subgroup, hence, the quotient $G/G_{w_0}$ is a smooth manifold and $\pi:G\to G/G_{w_0}$ is a smooth submersion (cf. Remark \autoref{rem:godement}). Since $G_{w_0}\subset G$ is also a complex submanifold, there exists a well-defined, unique, and smooth section $\hat I\in\Gamma \End (T (G/G_{w_0}))$ such that the diagram
 \begin{center}
 \begin{tikzcd}
  TG \arrow[r, "I"] \arrow[d,"d\pi"]
  & TG \arrow[d, "d\pi"] \\
  T(G/G_{w_0}) \arrow[r, "\hat I"]
  & T(G/G_{w_0})
 \end{tikzcd}
 \end{center}
 commutes. $I^2 = -\mathds{1}$ now implies $\hat I^2 = -\mathds{1}$, thus, $\hat I$ is an almost complex structure on $G/G_{w_0}$. Furthermore, $\hat I$ is integrable, as $I$ is integrable. This shows the first part of Proposition \autoref{prop:comp_quotient}.\\
 Let us now consider the real structure $\sigma$. By Proposition \autoref{prop:comp_stab}, $G_{w_0}$ is\linebreak $\sigma$-invariant, hence, the map $[\sigma]:G/G_{w_0}\to G/G_{w_0}$, $[g]\mapsto [\sigma (g)]$ is well-defined. Next, observe that the diagram
 \begin{center}
 \begin{tikzcd}
  G \arrow[r, "\sigma"] \arrow[d,"\pi"]
  & G \arrow[d, "\pi"] \\
  G/G_{w_0} \arrow[r, "{[\sigma]}"]
  & G/G_{w_0}
 \end{tikzcd}
 \end{center}
 commutes. As the vertical maps are submersions, we can infer the smoothness of $[\sigma]$ from the smoothness of $\sigma$. Likewise, the fact that $\sigma$ is an antiholomorphic involution implies that $[\sigma]$ is an antiholomorphic involution.
\end{proof}

We now examine how the complex and real structures of the adjoint orbits look like if we take them to be subsets of $\mathfrak{g}$:

\begin{proposition}\label{prop:holo_immersed}
 Let $G$ be a complex Lie group with real structure $\sigma$\linebreak and real form $G_\R\subset G$. Further, let $w_0\in\mathfrak{g}_\R\subset\mathfrak{g}$ be a point and\linebreak $G_{w_0}$ be the stabilizer of $w_0$. Then, the map $f_{w_0}:G/G_{w_0}\to\mathfrak{g}$, $[g]\mapsto\Ad (g)w_0$ is a well-defined, injective immersion whose image is the adjoint orbit $\mathcal{O}$ through $w_0$. Moreover, $f_{w_0}$ is $\hat I$-$I_e$-holomorphic and $[\sigma]$-$d\sigma_e$-equivariant.
\end{proposition}

\begin{proof}
 The first part follows from Lemma \autoref{lem:immersed_orbit}. Thus, we only need to show that $f_{w_0}$ intertwines the complex and real structures on $G/G_{w_0}$ and $\mathfrak{g}$, where the complex and real structure on $\mathfrak{g}$ are given by $I_e$ and $d\sigma_e$, respectively. To do that, we consider the evaluation of the adjoint action at $w_0$, i.e. $\ev_{w_0}:G\to\mathfrak{g}$, $g\mapsto \Ad (g)w_0$. The relation between $\ev_{w_0}$ and $f_{w_0}$ is encoded in the following commuting diagram:
 \begin{center}
 \begin{tikzcd}
  G \arrow[r, "\ev_{w_0}"] \arrow[d,"\pi"]
  & \mathfrak{g}\\
  G/G_{w_0} \arrow[ur, "f_{w_0}", dashed]
 \end{tikzcd}
 \end{center}
 $\pi$ is a submersion intertwining complex and real structures, hence, if $\ev_{w_0}$ intertwines complex and real structures, $f_{w_0}$ does so as well.\\
 Due to the $G$-equivariance of $\ev_{w_0}$ with respect to the left $G$-actions on $G$ and $\mathfrak{g}$, it suffices to check the Cauchy-Riemann equation at $e$ in order to show that $\ev_{w_0}$ is holomorphic. One finds $d\ev_{w_0,e} = -\ad_{w_0}$. Since $\ad_{w_0}$ commutes with $I_e$, $d\ev_{w_0,e}$ does so as well proving that $\ev_{w_0}$ is holomorphic.\\
 Lastly, we show the $\sigma$-$d\sigma_e$-equivariance of $\ev_{w_0}$. It follows from a straightforward calculation:
 \begin{align*}
  \ev_{w_0}\circ\sigma (g) &= \Ad (\sigma (g)) w_0 = d\sigma_e \left(\Ad (g) (d\sigma_e w_0)\right) = d\sigma_e \left(\Ad (g) w_0\right)\\
  &= d\sigma_e\circ\ev_{w_0} (g)\quad\forall g\in G.
 \end{align*}
\end{proof}

To conclude this subsection, we consider the real forms of complex adjoint orbits. We wish to establish a relation between them and adjoint orbits of $G_\R$. Intuitively, one might expect that the real form of the adjoint orbit of $G$ through $w_0$ is the adjoint orbit of $G_\R$ through $w_0$. In the quotient model, this would mean that the real form of $(G/G_{w_0}, [\sigma])$ is $G_\R/ (G_\R)_{w_0}$. In general, however, its real form might be strictly larger, since the sets
$B_1(G_{w_0},\sigma)$ and $B_2(G_{w_0},\sigma)$ (definition below) might not be equal:

\pagebreak

\begin{proposition}\label{prop:real_form_adjoint}
 Let $G$ be a complex Lie group with real\linebreak structure $\sigma$ and real form $G_\R\subset G$. Take a point $w_0\in\mathfrak{g}_\R\subset\mathfrak{g}$ and\linebreak let $G_{w_0}\subset G$ and $(G_\R)_{w_0}\subset G_\R$ be the stabilizers of $w_0$. Then, the map\linebreak $i:G_\R/(G_\R)_{w_0}\to G/G_{w_0}$, $[g]_{(G_\R)_{w_0}}\mapsto [g]_{G_{w_0}}$ is an embedding whose image is an open subset of $\Fix [\sigma]$. Further, one has:
 \begin{gather*}
  \im i = \Fix [\sigma] \quad\Leftrightarrow\quad B_1(G_{w_0},\sigma) = B_2(G_{w_0},\sigma),
 \end{gather*}
 where $B_1(G_{w_0},\sigma)\subset B_2(G_{w_0},\sigma)\subset G_{w_0}$ are defined as follows:
 \begin{gather*}
  B_1(G_{w_0},\sigma)\coloneqq\{\hat h^{-1}\sigma (\hat h)\mid \hat h\in G_{w_0}\},\quad B_2(G_{w_0},\sigma)\coloneqq G_{w_0}\cap \{ g^{-1}\sigma (g)\mid g\in G\}.
 \end{gather*}
\end{proposition}

\begin{remark}\label{rem:real_form_adjoint}
  For an adjoint orbit $\mathcal{O}\subset\mathfrak{g}$, its real form is $\mathcal{O}\cap\mathfrak{g}_\R$. Even though the real form of $\mathcal{O}$ is not closed under the adjoint action of the full group $G$ anymore, it is still closed under the adjoint action of $G_\R$. Thus, the real form is a disjoint union of adjoint orbits of $G_\R$ which, by Proposition \autoref{prop:real_form_adjoint}, are open subsets of $\mathcal{O}\cap \mathfrak{g}_\R$. In particular, if $G_\R$ is connected, the connected components of $\mathcal{O}\cap\mathfrak{g}_\R$ are just adjoint orbits of $G_\R$.
\end{remark}

\begin{proof}
 Consider the following commuting diagram:
 \begin{center}
 \begin{tikzcd}
  G_\R/(G_\R)_{w_0} \arrow[r, "i"] \arrow[d]
  & G/G_{w_0} \arrow[d, "f_{w_0}"] \\
  \mathfrak{g}_\R \arrow[r]
  & \mathfrak{g}
 \end{tikzcd}
 \end{center}
 The map on the left is given by $[g]_{(G_\R)_{w_0}}\mapsto \Ad (g)w_0$ and the bottom map is the natural inclusion. We know that the vertical maps as well as the map at the bottom are injective immersions. Thus, the map $i$ is also an injective immersion. We also know that $\im i\subset\Fix [\sigma]$, since:
 \begin{gather*}
  [\sigma]\circ i \left([g]_{(G_\R)_{w_0}}\right) = [\sigma (g)]_{G_{w_0}} = [g]_{G_{w_0}} = i \left([g]_{(G_\R)_{w_0}}\right)\quad\forall g\in G_\R = \Fix\sigma.
 \end{gather*}
 By Definition \autoref{def:real_form}, $\Fix [\sigma]$ is an embedded submanifold of $G/G_{w_0}$ whose dimension is given by:
 \begin{align*}
  \dim \Fix [\sigma] &= \frac{1}{2}\dim G/G_{w_0} = \frac{1}{2}\dim G - \frac{1}{2}\dim G_{w_0} = \dim G_\R - \dim (G_\R)_{w_0}\\
  &= \dim G_\R /(G_\R)_{w_0}.
 \end{align*}
 Hence, $i$ viewed as a map $G_\R/(G_\R)_{w_0}\to \Fix [\sigma]$ is an injective immersion between manifolds of the same dimension and, therefore, an embedding. This proves that $i:G_\R/(G_\R)_{w_0}\to G/G_{w_0}$ is the composition of two embeddings and, thus, itself an embedding. From this, we can also infer that $\im i$ is an open subset of $\Fix [\sigma]$.\\
 Lastly, we need to show the equivalence in Proposition \autoref{prop:real_form_adjoint}. It is a direct consequence of the following lemma.
\end{proof}
 
 \begin{lemma}\label{lem:group_involution_homology}
  Let $G$ be a group with involutive group homomorphism\linebreak $\sigma:G\to G$. Further, let $H\subset G$ be a $\sigma$-invariant subgroup. Then, $\sigma$ descends to an involution $[\sigma]$ on $G/H$ and the map $i: G_\sigma/H_\sigma\to G/H$, $[g]_{H_\sigma}\mapsto [g]_H$ is injective, where $G_\sigma\coloneqq\Fix\sigma$ and $H_\sigma\coloneqq\Fix\sigma\vert_H$. Further, one has:
  \begin{gather*}
   \im i = \Fix [\sigma] \quad\Leftrightarrow\quad B_1(H,\sigma) = B_2(H,\sigma),
  \end{gather*}
  where $B_1(H,\sigma)\subset B_2(H,\sigma)\subset H$ are defined as follows:
 \begin{gather*}
  B_1(H,\sigma)\coloneqq\{\hat h^{-1}\sigma (\hat h)\mid \hat h\in H\},\quad B_2(H,\sigma)\coloneqq H\cap \{g^{-1}\sigma (g)\mid g\in G\}.
 \end{gather*}
 \end{lemma}
 
 \begin{proof}
  Clearly, the maps $[\sigma]:G/H\to G/H$, $[\sigma]([g]_H) \coloneqq [\sigma (g)]_H$ and\linebreak $i: G_\sigma/H_\sigma\to G/H$ are well-defined. Moreover, it is obvious that $[\sigma]$ is an involution and the inclusion $\im i\subset\Fix [\sigma]$ holds. Next, we note that $i$ is injective: If we have $g_1, g_2\in G_\sigma$ such that $i ([g_1]_{H_\sigma}) = i ([g_2]_{H_\sigma})$, there exists an $h\in H$ such that $g_2 = g_1 h$. Hence, we obtain:
  \begin{gather*}
   h = g^{-1}_1 g_2 = \sigma (g_1)^{-1}\sigma (g_2) = \sigma (h)
  \end{gather*}
  This shows $h\in H_\sigma$ and, therefore, $[g_1]_{H_\sigma} = [g_2]_{H_\sigma}$.\\
  To prove the equivalence, we first consider the fixed points of $[\sigma]$ ($g\in G$):
  \begin{gather*}
   [g]_H \in \Fix [\sigma] \quad\Leftrightarrow\quad [\sigma (g)]_H = [g]_H \quad\Leftrightarrow\quad \sigma (g) = gh\text{ for }h\in H.
  \end{gather*}
  $\sigma (g) = gh$ implies $h = g^{-1}\sigma (g)\in B_2(H,\sigma)$. Thus, we obtain:
  \begin{gather}\label{eq:equivalence_1}
   [g]_H \in \Fix [\sigma] \quad\Leftrightarrow\quad \sigma (g) = gh\text{ for }h\in B_2(H,\sigma).
  \end{gather}
  Next, we investigate points in the image of $i$ ($g\in G$):
  \begin{gather*}
   [g]_H\in\im i \quad\Leftrightarrow\quad g = s\hat h\text{ for }s\in G_\sigma\text{ and }\hat h\in H.
  \end{gather*}
  To express this equivalence in a different way, we make the following considerations: If $s\in G_\sigma$ and $\hat h\in H$, then we find for $g\coloneqq s\hat h$:
  \begin{gather*}
   \sigma (g) = s\sigma (\hat h) = s\hat h\hat h^{-1}\sigma (\hat h) = g\hat h^{-1}\sigma (\hat h).
  \end{gather*}
  On the other hand, if $g\in G$, $\hat h\in H$, and $\sigma (g) = g\hat h^{-1}\sigma (\hat h)$, we obtain for $s\coloneqq g\hat h^{-1}$:
  \begin{gather*}
   \sigma (s) = \sigma (g)\sigma (\hat h)^{-1} = g\hat h^{-1} = s,
  \end{gather*}
  i.e., $g = s\hat h$ for $s\in G_\sigma$ and $\hat h \in H$. In total, this gives us the following equivalences:
  \begin{alignat}{3}
   \phantom{\Rightarrow}\quad &[g]_H\in\im i \quad && \Leftrightarrow\quad \sigma (g) = g\hat h^{-1}\sigma (\hat h) && \text{ for }\hat h\in H\notag\\
   \Rightarrow\quad &[g]_H\in\im i \quad && \Leftrightarrow\quad \sigma (g) = gh && \text{ for }h\in B_1(H,\sigma)\label{eq:equivalence_2},
  \end{alignat}
  where we set $h\coloneqq \hat h^{-1}\sigma (\hat h)\in B_1(H,\sigma)$. Combining Equivalence \eqref{eq:equivalence_1} and \eqref{eq:equivalence_2} now shows
  \begin{gather*}
   \im i = \Fix [\sigma] \quad\Leftrightarrow\quad B_1(H,\sigma) = B_2(H,\sigma)
  \end{gather*}
  concluding the proof.
 \end{proof}

 \chapter[Hyperkähler Structure of Coadjoint Orbits]{Hyperkähler Structure of Coadjoint Orbits\chaptermark{Hyperkähler Orbits}}
 \chaptermark{Hyperkähler Orbits}
 \label{app:hyp_orb}
 In this part, we sketch how coadjoint orbits of semisimple complex reductive\footnote{Throughout \autoref{app:hyp_orb}, a complex reductive group $G$ denotes the universal complexification of its compact real form $G_\R$.} groups obtain a Hyperkähler structure following \cite{Kronheimer1990} and \cite{Kovalev1996}. This task is carried out in two steps: We first explain how to equip the moduli space of instantons with a Hyperkähler structure. Afterwards, we put this structure on orbits by identifying them with the space of instantons. During this discussion, we also see that cotangent bundles of complex reductive groups obtain a Hyperkähler structure in very similar fashion (cf. \cite{Kronheimer2004}).\\
We begin by choosing a semisimple complex reductive group $G$ with real form $G_\R$. Next, we pick a Riemannian four-fold $(M,g)$. In \cite{Kronheimer1990}, $\R^4\backslash\{0\}\cong\C^2\backslash\{0\}$ plays the role of $M$, while $M$ is given by $\R_{\geq 0}\times T^3$ in \cite{Kovalev1996} (both spaces are equipped with the standard Euclidean metric). Let $P\coloneqq G_\R\times M$ be the trivial $G_\R$-principal bundle over $M$ and $A\in\Omega^1(P,\mathfrak{g}_\R)$ be a connection $1$-form on $P$. As $P$ is trivial, we can view $A$ as a $\mathfrak{g}_\R$-valued $1$-form on $M$ by fixing a gauge, i.e., by pulling $A$ back to $M$ via a global section, for instance $M\to P, p\mapsto (e,p)$. Now consider the curvature $2$-form $F\coloneqq dA + \frac{1}{2}[A\wedge A]$ associated with $A$. In our setup, we want $A$ to describe an instanton\footnote{In physics, instantons are special vacuum solutions of the Yang-Mills equations. For Abelian gauge theories, the Yang-Mills equations in matter (for instance the inhomogeneous Maxwell equations) take the form $\delta F = m$, where $\delta$ is the adjoint operator of $d$ and $m$ is a term involving matter fields. $m$ vanishes in the absence of matter, i.e., in vacuum. In this case, the closed $2$-form $F$ automatically satisfies the Yang-Mills equations if it is self-dual or anti-self-dual, i.e., $\ast F = \pm F$. These solutions describe particles known as instantons.}. To be precise, we want $F$ to satisfy the following anti-self-duality equation:
\begin{gather*}
 \ast F = -F,
\end{gather*}
where $\ast$ is the Hodge star associated with $(M,g)$. After choosing coordinates, the equation $\ast F = -F$ becomes a set of three PDEs depending on four variables. Solving these PDEs can be quite tricky. To simplify the problem, we want to transform the PDEs into ODEs which means we have to eliminate three of the four variables. We can achieve this by requiring $A$ to be invariant under a three-dimensional symmetry group $S$ acting freely on $M$. In \cite{Kronheimer1990}, the symmetry group is $\SU (2)$, while it is $T^3$ in \cite{Kovalev1996}. The anti-self-duality equation for $S$-invariant connections $A$ now takes the following form in suitable gauge and coordinates:
\begin{align*}
 \dot A_1 + [A_0,A_1] + [A_2,A_3] &= 0,\\
 \dot A_2 + [A_0,A_2] + [A_3,A_1] &= 0,\\
 \dot A_3 + [A_0,A_3] + [A_1,A_2] &= 0,
\end{align*}
where $A_j:\R\to\mathfrak{g}_\R$ ($j\in\{0,1,2,3\}$) are the $\mathfrak{g}_\R$-valued coefficients of the connection $A$. These equations are known as \textbf{Nahm's equations}\footnote{Confer \cite{Mayrand2019} for an in-depth analysis of Nahm's equations and their relation to Hyperkähler structures.}.\\
We can now formulate the strategy for equipping coadjoint orbits of $G$ with a Hyperkähler structure:
\begin{enumerate}
 \item First, one has to show that the moduli space of $S$-invariant instantons $A$, i.e., the space of solutions of Nahm's equations modulo gauge transformations, is a smooth manifold admitting a Hyperkähler structure.
 \item Afterwards, one uses objects like $2\alpha\coloneqq A_0 + iA_1$ and $2\beta\coloneqq A_2 + iA_3$ to identify this space with coadjoint orbits.
\end{enumerate}
The execution of this strategy is rather cumbersome and intransparent in \cite{Kronheimer1990} and \cite{Kovalev1996} due to some tedious details (boundary conditions, analysis of stable manifolds, involvement of the group structure, and so on). That is why we present the simpler approach taken in \cite{Kronheimer2004} here. In this paper, Kronheimer shows how to construct a Hyperkähler structure on the cotangent bundle of a complex reductive group with essentially the same strategy and without the tedious details.\\

We begin by considering the Banach space $\Omega\coloneqq C^1 ([0,1],\mathfrak{g}^4_\R)$, where\linebreak $\mathfrak{g}^4_\R = \mathfrak{g}_\R\times\mathfrak{g}_\R\times\mathfrak{g}_\R\times\mathfrak{g}_\R$. It carries a flat Hyperkähler structure induced by the linear maps
\begin{align*}
 I(A_0,A_1,A_2,A_3) &\coloneqq (-A_1, A_0, -A_3, A_2),\\
 J(A_0,A_1,A_2,A_3) &\coloneqq (-A_2, A_3, A_0, -A_1),\\
 K(A_0,A_1,A_2,A_3) &\coloneqq (-A_3, -A_2, A_1, A_0)
\end{align*}
and the $L^2$-metric\footnote{The $L^2$-metric $\skcdot_{L^2}$ is not the Banach metric on $\Omega$.} $\skcdot_{L^2}$ on $\Omega$ defined via an $\Ad$-invariant, positive definite scalar product $\skcdot_\R$ on $\mathfrak{g}_\R$:
\begin{gather*}
 \sk{A}{B}_{L^2} = \sum^3_{j= 0}\int\limits^1_0 \sk{A_j(t)}{B_j(t)}_{\R} dt.
\end{gather*}
Now consider the Banach space $\Gamma\coloneqq C([0,1],\mathfrak{g}^3_\R)$ and the map $\mu:\Omega\to\Gamma$ which maps $(A_0,A_1,A_2,A_3)$ to the left-hand side of Nahm's equations. We want to show that the space $N\coloneqq\mu^{-1}(0)$ of solutions of Nahm's equations is a smooth Banach submanifold of $\Omega$. By the infinite-dimensional regular value theorem (cf. Appendix A in \cite{Mcduff2012}), it suffices to show that for every $A\in N$ the differential $d\mu_A:\Omega\to\Gamma$ has a bounded linear right-inverse. Simply put, for every $A = (A_0,A_1,A_2,A_3)\in N$ and $\gamma = (\gamma_1,\gamma_2,\gamma_3)\in\Gamma$ we have to find a solution $a = (a_0,a_1,a_2,a_3)\in\Omega$ of the equation $d\mu_A (a) = \gamma$. Writing out the components of $d\mu_A (a) = \gamma$, we obtain:
\begin{align*}
 \dot a_1 - [A, Ia] &= \gamma_1,\\
 \dot a_2 - [A, Ja] &= \gamma_2,\\
 \dot a_3 - [A, Ka] &= \gamma_3,
\end{align*}
where we used the notation $[A,B] \coloneqq \sum^3_{j=0}[A_j,B_j]$ for maps $A:[0,1]\to\mathfrak{g}_\R^4$ and  $B:[0,1]\to\mathfrak{g}_\R^4$. These equations are just ODEs which possess plenty of solutions. We can single out a unique solution $a$ by requiring $a_0\equiv 0$ and $a_1(0) = a_2 (0) = a_3 (0) = 0$. Denote the map which assigns every $\gamma$ the specified solution $a$ by $d\mu^{-1}_A:\Gamma\to\Omega$. It is clear from the definition of $d\mu^{-1}_A$ that $d\mu^{-1}_A$ is a linear right-inverse of $d\mu_A$. To check that $d\mu^{-1}_A$ is bounded, we note that the image $V\coloneqq\im d\mu^{-1}_A$ is a closed subspace of $\Omega$ and, therefore, a Banach space. Thus, $d\mu_A\vert_V:V\to\Gamma$ is a bijective bounded operator between Banach spaces. Hence, $d\mu^{-1}_A$ as its inverse is also bounded due to the bounded inverse\linebreak theorem.\\
Even though $N\subset\Omega$ is a Banach submanifold, it does not inherit a Hyperkähler structure from $\Omega$. Indeed, the tangent space $T_A N = \ker d\mu_A$ is not invariant under the action of $I$, $J$, and $K$. For example, $\dot a_1 - [A,Ia] = 0$ turns into $\dot a^\prime_0 + [A,a^\prime] = 0$ under the transformation $a\mapsto a^\prime \coloneqq Ia$. However, this is the only obstruction, since $I$ transforms $\dot a_2 - [A, Ja] = 0$ into $\dot a^\prime_3 - [A, Ka^\prime] = 0$ and vice versa (a similar behavior occurs if we exchange $I$ for $J$ or $K$). Thus, the maximal $I$-$J$-$K$-invariant subspace
\begin{gather*}
 H_A\coloneqq T_AN\cap I(T_AN)\cap J(T_AN)\cap K(T_AN)\subset T_AN
\end{gather*}
is just the space of all $a\in\ker d\mu_A$ satisfying $\dot a_0 + [A,a] = 0$. Explicitly, this means that $a\in\Omega$ is contained in $H_A$ if and only if $a$ fulfills:
\begin{alignat*}{2}
 \dot a_0 + [A,a] &= 0,\quad \dot a_1 - [A, Ia] &&= 0,\\
 \dot a_2 - [A, Ja] &= 0,\quad \dot a_3 - [A, Ka] &&= 0.
\end{alignat*}
We see that the elements of $H_A$ are fully described by a set of four first-order ODEs implying that $H_A$ is isomorphic to $\mathfrak{g}^4_\R$ via $H_A\to\mathfrak{g}^4_\R$, $a\mapsto a(0)$. In particular, $H_A$ is a finite-dimensional subspace of $T_A N$ with dimension $4\dim_\R \mathfrak{g}_\R$, while $T_A N$ is infinite-dimensional. The $L^2$-metric on $\Omega$ naturally descends to a scalar product on $H_A$ turning $H_A$ into a Hyperkähler vector space.\\
Our next task is to find a manifold whose tangent spaces can be identified\linebreak with $H_A$. To do that, consider the group of gauge transformations\linebreak $\mathcal{G}_\R\coloneqq C^2([0,1], G_\R)$ and its normal subgroup $\mathcal{G}^0_\R\coloneqq\{g\in\mathcal{G}_\R\mid g(0) = g(1) = e\}$. The group $\mathcal{G}_\R$ acts as follows on $\Omega$:\footnote{In physics literature, the term $\Ad (g)w_g$ is sometimes incorrectly called Maurer-Cartan form and sloppily denoted by $gdg^{-1}$. For any Lie group $G$, the actual Maurer-Cartan form $\omega_{\text{MC}}\in\Omega^1 (G,\mathfrak{g})$ is the canonical (up to sign and $L$-$R$ conventions) $\mathfrak{g}$-valued one-form on $G$ defined as follows:
\begin{gather*}
 \omega_{\text{MC},g} (v)\coloneqq -dR_{g^{-1},g} (v)\quad\forall v\in T_gG\, \forall g\in G.
\end{gather*}
The term physicists denote by $gdg^{-1}$ is actually $g^\ast\omega_{\text{MC}}$, the pullback of the Maurer-Cartan form with respect to the gauge transformation $g:M\to G$.}
\begin{gather*}
 (gA_0) (t)\coloneqq \Ad (g(t))\left(A_0 (t) + w_g (t)\right),\quad (gA_j) (t)\coloneqq \Ad (g(t)) A_j (t)
\end{gather*}
where $g\in\mathcal{G}_\R$, $A = (A_0,A_1,A_2,A_3)\in\Omega$, $t\in[0,1]$, $j\in\{1,2,3\}$, and
\begin{gather*}
 w_g (t)\coloneqq \left.\frac{d}{ds}\right\vert_{s = t} g(s)^{-1}g(t)\in\mathfrak{g}_\R.
\end{gather*}
It is straightforward to check that gauge transformations map a solution of Nahm's equations to another solution of Nahm's equations. Hence, $\mathcal{G}_\R$ and $\mathcal{G}^0_\R$ also act on the Banach submanifold $N\subset\Omega$. We will not go into detail here why the quotient $N/\mathcal{G}^0_\R$ is again a Banach manifold\footnote{The full explanation can be found in \cite{Kronheimer2004}.}. The basic idea is to find a ``slice'' -- in our case $\{A+a\mid \dot a_0 + [A,a] = 0\}$ -- that hits every $\mathcal{G}^0_\R$-orbit exactly once and does so transversely. It is more important for us to realize that the tangent space of $N$ and the one of this slice intersect in $H_A$ at the point $A\in N$. Therefore, we can identify the tangent spaces of $N/\mathcal{G}^0_\R$ with the spaces $H_A$ which carry a linear Hyperkähler structure. This equips the quotient $N/\mathcal{G}^0_\R$ with a pre-Hyperkähler structure. The pre-Hyperkähler manifold $N/\mathcal{G}^0_\R$ is finite-dimensional, as its tangent space $H_A$ is finite-dimensional. At this point, it is not clear why this structure is integrable. The precise reason is a bit more involved, but the gist of it is that $N/\mathcal{G}^0_\R$ is the result of an infinite-dimensional Hyperkähler reduction of the space $\Omega$ with respect to the action $\mathcal{G}^0_\R$ and the Hyperkähler moment map $\mu:\Omega\to\Gamma$ (cf. \cite{Hitchin1987b} and \cite{Hitchin1987}).\\
So far, we have executed part (i) of our strategy, but have not used the complex group $G$. Indeed, $G$ is only important for part (ii), the identification of $N/\mathcal{G}^0_\R$ with the cotangent bundle of $G$. The core idea of part (ii) goes back to Donaldson (cf. \cite{Donaldson1984}) and essentially revolves around complexifying the constructions from part (i): We now consider the Banach space $\Omega_G\coloneqq C^1 ([0,1],\mathfrak{g}^2)$ and the group of complex gauge transformations
\begin{gather*}
 \mathcal{G}^0\coloneqq \{g\in C^2 ([0,1],G)\mid g(0) = g(1) = e\},
\end{gather*}
where $\mathcal{G}^0$ acts as follows on $\Omega_G$:
\begin{gather*}
 (g\alpha) (t)\coloneqq \Ad (g(t))\left(\alpha (t) + \frac{1}{2}w_g (t)\right),\quad (g\beta) (t)\coloneqq \Ad (g(t)) \beta (t)
\end{gather*}
where $g\in\mathcal{G}^0$, $(\alpha,\beta)\in\Omega_G$, $t\in[0,1]$, and
\begin{gather*}
 w_g (t)\coloneqq \left.\frac{d}{ds}\right\vert_{s = t} g(s)^{-1}g(t)\in\mathfrak{g}.
\end{gather*}
Instead of Nahm's equations, we consider their complex counterpart:
\begin{gather}
 \dot \beta + 2[\alpha,\beta] = 0\label{eq:complex_nahm}.
\end{gather}
We denote the space of solutions of \autoref{eq:complex_nahm} by $N_G\subset\Omega_G$. In \cite{Donaldson1984} and \cite{Kronheimer2004}, it is shown that
\begin{gather*}
 \alpha = \frac{1}{2}(A_0 + iA_1),\quad \beta = \frac{1}{2}(A_2 + iA_3)
\end{gather*}
furnishes an isomorphism between the quotients $N/\mathcal{G}^0_\R$ and $N_G/\mathcal{G}^0$. The upshot of this isomorphism is that the space $N_G$ is much easier to describe than $N$, because \autoref{eq:complex_nahm} can be solved directly by integration: For every $(\alpha,\beta)\in N_G$, there exists exactly one element $\eta\in\mathfrak{g}$ and exactly one curve $u\in C^2([0,1],G)$ with $u(0) = e$ such that:
\begin{gather*}
 \alpha(t) = \frac{1}{2}\Ad (u(t)) w_u (t),\quad \beta (t) = \Ad (u(t))\eta.
\end{gather*}
It is straightforward to verify that replacing $(\alpha,\beta)\in N_G$ by $(g\alpha,g\beta)$ is equivalent to replacing $u$ by $gu$. Thus, $(u,\eta)$ and $(u^\prime,\eta^\prime)$ describe the same $\mathcal{G}^0$-orbit in $N_G$ if and only if there exists $g\in \mathcal{G}^0$ such that $u^\prime = gu$ and $\eta^\prime = \eta$. If $u(1) = u^\prime (1)$, this is always the case, as we can simply set $g\coloneqq u^\prime u^{-1}$. Conversely, if $u^\prime = gu$, then we must have $u(1) = u^\prime (1)$, because $g(1) = e$. In total, we have shown that every element in $N/\mathcal{G}^0_\R\cong N_G/\mathcal{G}^0$ is completely described by a tuple $(u(1),\eta)\in G\times\mathfrak{g}$, i.e., $N/\mathcal{G}^0_\R\cong G\times\mathfrak{g}$. The space $G\times\mathfrak{g}$ itself is isomorphic to $TG$ by identifying $T_gG$ with $\mathfrak{g}$ via left or right multiplication (cf. \autoref{sec:lie_groups}). The tangent bundle $TG$ can be, in turn, identified with $T^\ast G$ via the bi-invariant semi-Riemannian metric $\skcdot$ on $G$ obtained from $\skcdot_\R$ (cf. \autoref{sec:lie_groups} and Proposition \autoref{prop:comp_of_metric}). It is in this way that $TG$ and $T^\ast G$ inherit the Hyperkähler structure from $N/\mathcal{G}^0_\R$.\\

At this stage, we should address some differences between the procedure in \cite{Kronheimer2004} and the one in \cite{Kronheimer1990} and \cite{Kovalev1996}. The first difference regards gauge fixing. In \cite{Kronheimer1990} and \cite{Kovalev1996}, the authors do not consider a quotient like $N/\mathcal{G}^0_\R$, but get rid of the gauge freedom by fixing a gauge, i.e., picking one connection from each $\mathcal{G}^0_\R$-orbit. The chosen gauge is in both papers the temporal gauge, i.e., $A_0\equiv 0$.\\
Arguably the most important difference concerns boundary conditions, because they allow us to describe coadjoint orbits instead of cotangent bundles. The approach taken in \cite{Kronheimer2004}, the approach we presented, has no need for boundary conditions, since all maps are defined on the closed interval $[0,1]$. However, the connections examined in \cite{Kronheimer1990} and \cite{Kovalev1996} might exhibit problematic behavior, as their basepoint approaches the origin in the case of \cite{Kronheimer1990}\linebreak ($M = \R^4\backslash\{0\}$) or as the first component tends to $+\infty$ in the case of \cite{Kovalev1996} ($M = \R_{\geq 0}\times T^3$). To account for this, the connections in \cite{Kronheimer1990} and \cite{Kovalev1996} do not only solve Nahm's equations, but also have to satisfy certain boundary conditions. In both cases, we first need to fix a triple\footnote{The same one mentioned in the introduction of \autoref{sec:duality}.} $\tau = (\tau_1,\tau_2,\tau_3)\in\mathfrak{g}^3_\R$.\pagebreak\linebreak
Afterwards, we impose the following conditions:
\begin{itemize}
 \item The connections $A_j (t)$ in \cite{Kronheimer1990} ($t\in\R$ and $j\in\{1,2,3\}$) approach zero, as $t$ tends to $+\infty$ (corresponds to the point at infinity turning $\R^4$ into $S^4$), and behave like $e^{-2t}\Ad (g)\tau_j$ for some $g\in G_\R$, as $t$ goes to $-\infty$ (corresponds to the origin of $\R^4$). The space of these connections, denoted by $M(\tau_1,\tau_2,\tau_3)$ in \cite{Kronheimer1990}, then obtains a Hyperkähler structure in the previously explained manner.
 \item The connections $A_j (t)$ in \cite{Kovalev1996} ($t\in\R_{\geq 0}$ and $j\in\{1,2,3\}$), on the other hand, merely need to coincide with $\Ad (g)\tau_j$ for some $g\in G_\R$ in the limit $t\to + \infty$. We do not need to specify a boundary condition for $t\to -\infty$, as $A_j$ is only defined on $\R_{\geq 0}$. Counterintuitively, the space of these connections, called $M(\tau)$ in \cite{Kovalev1996}, does not possess a Hyperkähler structure. Only the subspace $M_{\exp} (\tau)\subset M(\tau)$ of exponentially fast converging solutions and the subspace $M (\tau;\rho)\subset M(\tau)$ of solutions asymptotic to the solution determined by the homomorphism $\rho$\footnote{Confer \cite{Kovalev1996} for details.} carry a Hyperkähler structure which they obtain in the usual manner.
\end{itemize}
The last difference regards the identification of $M(\tau_1,\tau_2,\tau_3)$, $M_{\exp} (\tau)$, and $M (\tau;\rho)$ with adjoint orbits\footnote{We only consider adjoint orbits here, as coadjoint orbits can be identified with adjoint orbits via an $\Ad$-invariant scalar product on $\mathfrak{g}$.}. In \cite{Kronheimer2004}, the identification of $N/\mathcal{G}^0_\R$ with $TG$ and $T^\ast G$ is rather involved. For adjoint orbits, it is much easier to construct the isomorphism. We simply map the spaces $M(\tau_1,\tau_2,\tau_3)$, $M_{\exp} (\tau)$, and $M (\tau;\rho)$ to an adjoint orbit of $G$ via:\footnote{A change of variables is needed in the case of $M(\tau_1,\tau_2,\tau_3)$ (cf. \cite{Kronheimer1990}).}
\begin{gather*}
 (A_1,A_2,A_3)\mapsto 2\beta (0) = A_2 (0) + iA_3 (0).
\end{gather*}
Which space is mapped to which adjoint orbit, depends on the choice of $\tau$ and $\rho$. Nevertheless, two choices might be mapped to the same adjoint orbit giving us two different Hyperkähler structures on this orbit.\\
Before we conclude \autoref{app:hyp_orb}, we want to point out two similarities between the construction of Hyperkähler and holomorphic Kähler structures on coadjoint orbits. The first similarity concerns the holomorphic symplectic forms underlying Hyperkähler/holomorphic Kähler manifolds. As explained in \autoref{app:kaehler}, every holomorphic Kähler manifold $(X,\omega,J,I)$ also admits a holomorphic symplectic form $\Omega\coloneqq\omega - i\omega(I\cdot,\cdot)$, while Hyperkähler manifolds even possess three such forms, one for each complex structure $I$, $J$, and $K$. We have seen in \autoref{sec:holo_semi-kaehler} that the tensors $\omega$ and $I$ associated with a holomorphic Kähler structure of a coadjoint orbit are just the KKS form and the complex structure coming from the underlying complex group, respectively. Together, they give rise to the canonical form $\Omega_{\text{KKS}}$, known as the holomorphic KKS form. Similarly, it is shown in \cite{Kronheimer1990} that we can choose the Hyperkähler structure on coadjoint orbits such that the holomorphic symplectic form associated with $I$ is also just $\Omega_{\text{KKS}}$.\\
The second similarity is related to the fact that both structures come in families, to be specific, that both structures depend on the choice of an $\Ad$-invariant, positive definite scalar product $\skcdot_\R$ on $\mathfrak{g}_\R$. We need $\skcdot_\R$ to identify adjoint with coadjoint orbits in the holomorphic Kähler case, while, in the Hyperkähler case, we have to choose a scalar product $\skcdot_\R$ in order to define $\skcdot_{L^2}$.

 \chapter[Reduction: From Cotangent Bundles to Coadjoint Orbits]{Reduction: From Cotangent Bundles to Coadjoint Orbits\chaptermark{Reduction}}
 \chaptermark{Reduction}
 \label{app:reduction}
 In this part, we explain the process of reduction which allows us to transfer a symplectic, Kähler, holomorphic Kähler, or Hyperkähler structure\footnote{Confer \autoref{app:kaehler} for the definition of the various Kähler structures.} from a manifold to a quotient. Our goal is to illustrate how one can use reduction to relate the various Kähler structures on cotangent bundles to similar structures on coadjoint orbits.\\
Consider a symplectic manifold $(M,\omega)$ and a Lie group $G$ acting on $M$. We want to divide $M$ by $G$ in such a way that the quotient still carries a symplectic structure. In order for this to be possible, $G$ has to be compatible with $(M,\omega)$ in the following sense:

\begin{definition}[Hamiltonian action and moment map]
 Let $(M,\omega)$ be a symplectic manifold and $G$ be a Lie group acting smoothly on $M$. We call the $G$-action on $M$ \textbf{symplectic} if for every $g\in G$ the diffeomorphism $g:M\to M$, $p\mapsto gp$ preserves $\omega$, i.e., $g^\ast \omega = \omega$. The action is said to be \textbf{Hamiltonian} if it is symplectic and a smooth \textbf{moment map} $\mu:M\to\mathfrak{g}^\ast$ exists such that:
 \begin{enumerate}
  \item $d\mu_v = \iota_{X_v}\omega$\footnote{In \autoref{chap:PHHS}, we say $X_H$ is the Hamiltonian vector field of the Hamiltonian system $(M,\omega,H)$ if $\iota_{X_H}\omega = -dH$ (mind the minus sign). In this sense, $X_v$ is the Hamiltonian vector field of $(M,\omega,\mu_{-v})$.}, where $v\in\mathfrak{g}$, $\mu_v\in C^\infty (M)$ is the function on $M$ defined by $\mu_v (p)\coloneqq \mu (p)(v)$, and $X_v$ is the fundamental vector field (cf. \autoref{sec:semi-kaehler}) associated with $v$, i.e., $X_v (p)\coloneqq \left.\frac{d}{dt}\right\vert_{t = 0}\exp (tv)p$.
  \item $\mu (gp) = \Ad^\ast (g)\mu (p)$ for every $g\in G$ and $p\in M$, where $\Ad^\ast$ denotes the coadjoint action (cf. \autoref{sec:lie_groups}).
 \end{enumerate}
\end{definition}

\begin{remark}\ 
 \begin{enumerate}[label = (\arabic*)]
  \item There is a deep connection between Condition (i) and an action being symplectic: For every symplectic action and point $p\in M$, there exists an open neighborhood $U\subset M$ of $p$ and a map $\mu:U\to\mathfrak{g}^\ast$ satisfying Condition (i). This is due to the fact that the form $\iota_{X_v}\omega$ is closed if the $G$-action is symplectic. However, $\iota_{X_v}\omega$ does not need to be exact, i.e., Condition (i) does not need to be fulfilled globally.\\
  In the converse direction, a similar problem occurs: If $G$ is connected and the $G$-action satisfies Condition (i), then the action is automatically symplectic. This statement fails for Lie groups $G$ with multiple connected components, for instance $G = \Z_2$ acting on $M$ by an antisymplectic involution.
  \item There are several ways to interpret Condition (ii): Condition (ii) in the given form states that $\mu$ is $G$-equivariant with respect to the action on $M$ and the coadjoint action. If we assume Condition (i) and that $G$ is connected, then Condition (ii) is equivalent to the statement that\linebreak $\mu:(M,\{\cdot,\cdot\}_M)\to (\mathfrak{g}^\ast,\{\cdot,\cdot\}_{\mathfrak{g}^\ast})$ is a Poisson map, i.e.
  \begin{gather*}
   \{f\circ\mu, g\circ\mu\}_M = \{f,g\}_{\mathfrak{g}^\ast}\circ\mu\quad\forall f,g\in C^\infty (\mathfrak{g}^\ast),
  \end{gather*}
  where $\{\cdot,\cdot\}_M$ is the Poisson structure on $M$ induced by $\omega$ and $\{\cdot,\cdot\}_{\mathfrak{g}^\ast}$ denotes the canonical Poisson structure on $\mathfrak{g}^\ast$ (cf. Proposition \autoref{prop:dual_lie}).\linebreak $\mu$ being a Poisson map is, in turn, equivalent to the statement that the map $(\mathfrak{g}, [\cdot,\cdot])\to(C^\infty (M),\{\cdot,\cdot\}_M)$, $v\mapsto \mu_v$ is a Lie algebra homomorphism.
 \end{enumerate}
\end{remark}

Before we turn our attention to the reduction theorems, it is instructive to discuss two examples of Hamiltonian actions: The coadjoint action and the $G$-action on $T^\ast G$. These actions will play a vital role in the course of \autoref{app:reduction}.

\begin{example}[Coadjoint action]\label{ex:coadj}
 Choose any Lie group $G$ and take $M$ to be a coadjoint orbit $\mathcal{O}^\ast\subset\mathfrak{g}^\ast$ of $G$. For the symplectic form $\omega$, we choose the KKS form $\kks$ (cf. \autoref{sec:semi-kaehler}). The moment map $\mu:\mathcal{O}^\ast\to\mathfrak{g}^\ast$ is just the inclusion, i.e., $\mu(\alpha)\coloneqq\alpha$. We verify that the coadjoint action constitutes a Hamiltonian action with these choices: The coadjoint action is symplectic, because $\kks$ is $G$-invariant, as we have checked in \autoref{sec:semi-kaehler}. A quick calculation further reveals that Condition (i) is also satisfied:
\begin{gather*}
 \omega_{\text{KKS},\alpha} (X^\ast_v (\alpha), X^\ast_w (\alpha)) = \alpha ([v,w]) = X^\ast_w (\alpha)(v) = d\mu_{v,\alpha} (X^\ast_w (\alpha)),
\end{gather*}
where $\alpha\in\mathcal{O}^\ast$, $v,w\in\mathfrak{g}$, $X^\ast_w (\alpha)\coloneqq -\alpha\circ\ad_w$ is the fundamental vector field associated with the coadjoint action, and we identify $T_\alpha\mathcal{O}^\ast$ with the span of the fundamental vector fields. Condition (ii) is trivially fulfilled.
\end{example}

\begin{example}[$G$-action on $T^\ast G$]\label{ex:cotangent}
 We again choose any Lie group $G$, but take $M$ to be the cotangent bundle $T^\ast G$ this time. $G$ acts on $T^\ast G$ by right multiplication, i.e.:
 \begin{gather*}
  G\times T^\ast G\to T^\ast G,\quad (g,\alpha\in T^\ast_h G)\mapsto \alpha\circ dR_g\in T^\ast_{hg^{-1}}G.
 \end{gather*}
 For the symplectic form $\omega$, we choose the canonical two-form $\omega_{\can}$ on cotangent bundles. To construct the moment map $\mu:T^\ast G\to\mathfrak{g}^\ast$, we first identify $T^\ast G$ with the trivial bundle $G\times\mathfrak{g}^\ast$ using left multiplication and afterwards project onto the second component, i.e.:
 \begin{gather*}
  \mu (\alpha)\coloneqq \alpha\circ dL_h\quad\forall \alpha\in T^\ast_h G\ \forall h\in G.
 \end{gather*}
 Let us check that the action in question is indeed Hamiltonian. First, we employ a general fact to show that the action is symplectic: If $f:Q\to Q$\linebreak is a diffeomorphism of a manifold $Q$, then $f$ gives rise to the following symplectomorphism of $(T^\ast Q,\omega_{\text{can}})$:
 \begin{gather*}
  f_{T^\ast Q}:T^\ast Q\to T^\ast Q,\quad \alpha\in T^\ast_p Q\mapsto \alpha\circ df^{-1}\in T^\ast_{f(p)}Q.
 \end{gather*}
 We now realize that for any $g\in G$ the transformation $\alpha\mapsto g\alpha = \alpha\circ dR_g$ is just the symplectomorphism induced by the diffeomorphism $h\mapsto hg^{-1}$. Thus, the $G$-action on $T^\ast G$ is symplectic.\\
 To verify Condition (ii), we compute:
 \begin{align*}
  \mu (g\alpha) &= \mu (\alpha\circ dR_g) = \alpha\circ dR_g\circ dL_{hg^{-1}}\\
  &= \alpha\circ dR_g\circ dL_h\circ dL_{g^{-1}}\\
  &= \alpha\circ dL_h\circ dR_g\circ dL_{g^{-1}}\\
  &= \mu (\alpha)\circ \Ad (g^{-1}) = \Ad^\ast (g)\mu (\alpha),
 \end{align*}
 where $\alpha\in T^\ast_h G$, $g,h\in G$, we exploited the fact that left and right multiplication commute\footnote{By associativity, we have $a(bc) = (ab)c$, so it does not matter whether we first multiple by $a$ from the left or by $c$ from the right.}, and used the definition of the adjoint action:
 \begin{gather*}
  \Ad (g^{-1}) = dR_g\circ dL_{g^{-1}}.
 \end{gather*}
 Lastly, we check Condition (i). To keep the calculation simple, we identify $T^\ast G$ with $T^\ast_L G\coloneqq G\times\mathfrak{g}^\ast$:
 \begin{gather*}
  T^\ast G\to T^\ast_L G,\quad \alpha\in T^\ast_g G\mapsto (g,\alpha\circ dL_g).
 \end{gather*}
 $G$ now acts as follows on $T^\ast_L G$:
 \begin{gather*}
  G\times T^\ast_L G\to T^\ast_L G,\quad (g, (h,\alpha))\mapsto (hg^{-1}, \Ad^\ast (g)\alpha).
 \end{gather*}
 The moment map $\mu^L:T^\ast_L G\to\mathfrak{g}^\ast$ is just the projection onto the second component. Under the given identification, $\omega_{\can}$ becomes the two-form $\omega^L$ on $T^\ast_L G$:
 \begin{gather*}
  \omega^L_{(g,\alpha)}\left((v_1,\beta_1), (v_2,\beta_2)\right) = \beta_1 (v_2) - \beta_2 (v_1) - \alpha\left([v_1,v_2]\right),
 \end{gather*}
 where $(g,\alpha)\in T^\ast_L G$ and $(v_1,\beta_1), (v_2,\beta_2)\in\mathfrak{g}\times\mathfrak{g}^\ast$. Here, we take the tangent spaces of $T^\ast_L G = G\times\mathfrak{g}^\ast$ to be $\mathfrak{g}\times\mathfrak{g}^\ast$ by identifying $T_gG$ with $\mathfrak{g}$ via $dL_{g^{-1}}$. Using the same identification, we can express the fundamental vector fields $X_v$ of the $G$-action on $T^\ast_L G$ as follows:
 \begin{gather*}
  X_v (g,\alpha) = - (v,\alpha\circ\ad_v),
 \end{gather*}
 where $v\in\mathfrak{g}$ and $(g,\alpha)\in T^\ast_L G$. We now calculate:
 \begin{align*}
  \left(\iota_{X_v}\omega^L\right)_{(g,\alpha)} (w,\beta) &= \omega^L_{(g,\alpha)} \left(X_v (g,\alpha), (w,\beta)\right)\\
  &= -\alpha ([v,w]) + \beta (v) + \alpha ([v,w])\\
  &= \beta (v) = d\mu^L_{v, (g,\alpha)} (w,\beta).
 \end{align*}
 This shows $\iota_{X_v}\omega^L = d\mu^L_v$ for all $v\in\mathfrak{g}$ and proves that Condition (i) is satisfied.
\end{example}

Most books and papers discussing reduction assume that the group $G$ is compact to simplify proofs. In our case, however, this is not possible, since we also want to consider reduction for complex groups $G$ which are usually not compact. Dropping compactness introduces some technical difficulties. The following proposition takes care of these issues:

\begin{proposition}\label{prop:prep}
 Let $(M,\omega)$ be a symplectic manifold, $G$ a Lie group, and $G\times M\to M$ a Hamiltonian action with moment map $\mu:M\to\mathfrak{g}^\ast$. Further, let $\mathcal{O}^\ast\subset\mathfrak{g}^\ast$ be a coadjoint orbit of $G$.
 \begin{enumerate}[label = (\alph*)]
  \item If $(N_1,\omega_1)$ and $(N_2,\omega_2)$ are symplectic manifolds equipped\linebreak with Hamiltonian $G$-actions and moment maps $\mu_j:N_j\to\mathfrak{g}^\ast$, then\linebreak $(N_1\times N_2,\pr_1^\ast\omega_1 + \pr_2^\ast\omega_2)$ is also a symplectic manifold with Hamiltonian $G$-action and moment map $\mu(p,q)\coloneqq \mu_1 (p) + \mu_2 (q)$. In particular, the $\mathcal{O}^\ast$-extended space
  \begin{gather*}
   (\tilde M\coloneqq M\times\mathcal{O}^\ast, \pr_1^\ast\omega - \pr_2^\ast\kks)
  \end{gather*}
  carries a Hamiltonian $G$-action with moment map $\tilde\mu (p,\eta)\coloneqq \mu (p)-\eta$.
  \item Let $G_\eta\subset G$ be the stabilizer of $\eta\in\mathcal{O}^\ast$. The following statements are equivalent:
        \begin{enumerate}[label = (\roman*)]
         \item $G$ acts freely on $\mu^{-1}(\mathcal{O}^\ast)$.
         \item $G_\eta$ acts freely on $\mu^{-1}(\eta)$ for every $\eta\in\mathcal{O}^\ast$.
         \item $G_\eta$ acts freely on $\mu^{-1}(\eta)$ for one $\eta\in\mathcal{O}^\ast$.
        \end{enumerate}
        If (i), (ii), or (iii) holds, then $\mu$ is a submersion on an open neighborhood of $\mu^{-1}(\mathcal{O}^\ast)$. In particular, $\mu^{-1}(\eta)\subset M$ ($\eta\in\mathcal{O}^\ast$) and $\tilde\mu^{-1}(0)\subset\tilde M$ are embedded submanifolds. $\mu^{-1}(\mathcal{O}^\ast)\subset M$ is an immersed submanifold and embedded iff $\mathcal{O}^\ast\subset\mathfrak{g}^\ast$ is an embedded submanifold. $\mu^{-1}(\mathcal{O}^\ast)$ and $\tilde\mu^{-1}(0)$ are isomorphic via the diffeomorphism $\mu^{-1}(\mathcal{O}^\ast)\to\tilde\mu^{-1}(0)$, $p\mapsto (p,\mu (p))$.
  \item Let $N_1$ and $N_2$ be two manifolds with smooth $G$-actions on them. If the action on $N_1$ is proper, then the action on $N_1\times N_2$ is proper as well. In particular, the $G$-action on $\tilde M$ is proper if the action on $M$ is proper. If $G$ acts freely on $\mu^{-1}(\mathcal{O}^\ast)$ and the action on $M$ is proper, then $\mu^{-1}(\eta)/G_\eta$ ($\eta\in\mathcal{O}^\ast$), $\mu^{-1}(\mathcal{O}^\ast)/G$, and $\tilde\mu^{-1}(0)/G$ are manifolds which are diffeomorphic via the maps:
  \begin{alignat*}{2}
   \mu^{-1}(\eta)/G_\eta\to &\mu^{-1}(\mathcal{O}^\ast)/G &&\to \tilde\mu^{-1}(0)/G,\\
   [p]_{G_\eta}\mapsto &\qquad [p]_G &&\mapsto [p,\mu (p)]_G.
  \end{alignat*}
 \end{enumerate}
\end{proposition}

\begin{proof}
 It is straightforward to check the first part of (a). The second part of (a) follows from the first by setting $N_1 = M$ and $N_2 = \mathcal{O}^\ast$ equipped with the (negative) KKS structure from Example \autoref{ex:coadj}.\\
 We begin the proof of (b) by showing the equivalence. The implications ``(i)$\Rightarrow$(ii)'' and ``(ii)$\Rightarrow$(iii)'' are trivial, so we only prove ``(iii)$\Rightarrow$(i)''. For this, it suffices to show that $G_p\subset G_\eta$ for every $(p,\eta)\in M\times\mathfrak{g}^\ast$ with $\mu (p) = \eta$. $G_p\subset G_\eta$ immediately follows from the definition of a moment map:
 \begin{gather*}
  \Ad^\ast (g)\eta = \Ad^\ast (g)\mu (p) = \mu (gp) = \mu (p) = \eta\quad\forall g\in G_p.
 \end{gather*}
 For the next assertion, we assume that $G$ acts freely on $\mu^{-1}(\mathcal{O}^\ast)$. It is shown in \cite{Mayrand2016} (cf. Lemma 4.1 in this paper) that, under this condition, every $\eta\in\mathcal{O}^\ast$ is a regular value\footnote{In \cite{Mayrand2016}, this result is only shown for compact groups. However, the same proof still applies in our case.} of $\mu$. Thus, $d\mu_p$ has full rank for every $p\in\mu^{-1}(\mathcal{O}^\ast)$. Having full rank is an open property, hence, $\mu$ is a submersion on a neighborhood of $\mu^{-1}(\mathcal{O}^\ast)$. This also implies that $\tilde\mu$ is a submersion on a neighborhood\linebreak of $\tilde\mu^{-1}(0)$.\\
 To show that the sets $\mu^{-1}(\eta)$, $\tilde\mu^{-1}(0)$, and $\mu^{-1}(\mathcal{O}^\ast)$ are manifolds, we employ a standard result from differential geometry: Let $N_1$ and $N_2$ be two manifolds, $f:N_1\to N_2$ a smooth submersion, and $S\subset\im f$ a set. Then, $f^{-1}(S)\subset N_1$ is an immersed submanifold if $S\subset N_2$ is an immersed submanifold. $f^{-1}(S)\subset N_1$ is embedded if and only if $S\subset N_2$ is embedded.\\
 As $\{\eta\}\subset\mathfrak{g}^\ast$ and $\{0\}\subset\mathfrak{g}^\ast$ are clearly embedded submanifolds, $\mu^{-1}(\eta)\subset M$ and $\tilde\mu^{-1}(0)\subset\tilde M$ are embedded submanifolds\footnote{This is just the regular value theorem.} as well. The coadjoint orbit $\mathcal{O}^\ast\subset\mathfrak{g}^\ast$, on the other hand, is usually just an immersed submanifold (cf. Lemma \autoref{lem:immersed_orbit}). Thus, $\mu^{-1}(\mathcal{O}^\ast)\subset M$ is an immersed submanifold and embedded iff $\mathcal{O}^\ast\subset\mathfrak{g}^\ast$ is embedded.\\
 To prove the last statement in (b), we note that $\mu^{-1}(\mathcal{O}^\ast)\to\tilde\mu^{-1}(0)$,\linebreak $p\mapsto (p,\mu (p))$ is a bijective immersion. By dimensional reasoning, it must also be a submersion concluding the proof of (b).\\
 Lastly, we show (c). For the first assertion, denote the map
 \begin{gather*}
  G\times N_1\times N_2\to N^2_1\times N^2_2,\quad (g,n_1,n_2)\mapsto (gn_1,n_1,gn_2,n_2)
 \end{gather*}
 by $f$ and consider a compact subset $K\subset N^2_1\times N^2_2$. We need to show that the preimage $f^{-1}(K)$ is compact. $f^{-1}(K)$ is closed, because $f$ is continuous and $K$ is closed. Thus, it suffices to show that $f^{-1}(K)$ is contained in a compact set. Let $\pr_1:N^2_1\times N^2_2\to N^2_1$ be the projection onto the first two components and $\pr_2:N^2_1\times N^2_2\to N_2$ the projection onto the last component. If we denote the map $G\times N_1\to N^2_1$, $(g,n_1)\mapsto (gn_1, n_1)$ by $f_1$, we find:
 \begin{gather*}
  f^{-1}(K)\subset f_1^{-1}(\pr_1 (K))\times\pr_2 (K).
 \end{gather*}
 The compactness of $K$ implies that $\pr_1 (K)$ and $\pr_2 (K)$ are compact. Moreover, we know that the action on $N_1$ is proper meaning that $f_1$ is proper. Thus, $f_1^{-1}(\pr_1 (K))$ is also compact. This shows that $f^{-1}(K)$ is contained in the compact set $f_1^{-1}(\pr_1 (K))\times\pr_2 (K)$.\\
 The second statement in (c) follows from the first one by applying it to the $\mathcal{O}^\ast$-extended space $\tilde M$. To prove the last assertion, we observe that, in the given setup, $G$ acts freely and properly on $\mu^{-1}(\mathcal{O}^\ast)$ and $\tilde\mu^{-1}(0)$, while $G_\eta$ does so on $\mu^{-1}(\eta)$. Hence, the quotients inherit a manifold structure by Theorem \autoref{thm:godement}. The maps in (c) are clearly bijections and, due to the universal property of quotient spaces, also immersions. Dimensional reasoning now shows that they are diffeomorphisms concluding the proof.
\end{proof}

We are now ready to describe the process of symplectic reduction\footnote{Symplectic (or Marsden-Weinstein) reduction is a fundamental tool in symplectic geometry and can be found in most introductory books, for instance \cite{Abraham1978} and \cite{mcduff2017}. Even though Kähler and Hyperkähler reduction are not as standard, they are still known to most symplectic geometers. \cite{Mayrand2016} provides an overview over the different kinds of reduction.}:

\pagebreak

\begin{theorem}[Marsden-Weinstein reduction]\label{thm:marsden-weinstein}
 Let $(M,\omega)$ be a symplectic manifold, $G$ a Lie group, and $G\times M\to M$ a proper Hamiltonian action with moment map $\mu:M\to\mathfrak{g}^\ast$. Assume that $\mathcal{O}^\ast\subset\mathfrak{g}^\ast$ is a coadjoint orbit such that $G$ acts freely on $\mu^{-1}(\mathcal{O}^\ast)$. Then, there exists a unique symplectic form $\bar{\omega}$ on $M//G (\mathcal{O}^\ast)\coloneqq \mu^{-1}(\mathcal{O}^\ast)/G$ such that $\pi^\ast\bar{\omega} = i^\ast\omega - (\mu\circ i)^\ast\kks$, where $\pi:\mu^{-1}(\mathcal{O}^\ast)\to\mu^{-1}(\mathcal{O}^\ast)/G$ is the natural projection, $i:\mu^{-1}(\mathcal{O}^\ast)\hookrightarrow M$ denotes the inclusion of $\mu^{-1}(\mathcal{O}^\ast)$ into $M$, and $\kks$ is the KKS form on $\mathcal{O}^\ast$.
\end{theorem}

\begin{proof}
 There are two methods to construct the form $\bar{\omega}$: The idea of the first one is to identify $M//G (\mathcal{O}^\ast)$ with $\tilde\mu^{-1}(0)/G$ (cf. Proposition \autoref{prop:prep}) and equip it with a symplectic structure via coisotropic reduction. For the second method, we pick an element $\eta\in\mathcal{O}^\ast$ and identify $M//G (\mathcal{O}^\ast)$ with $\mu^{-1}(\eta)/G_\eta$ as in Proposition \autoref{prop:prep}. We then use a slightly generalized reduction procedure to construct the form $\bar{\omega}$ on $\mu^{-1}(\eta)/G_\eta$. The first method is well suited to show the properties of $\bar{\omega}$, while the second method is preferable if we want to construct $\bar{\omega}$ explicitly. We present both methods here, but focus on the first one.\\
 
 \textbf{Method 1:} $M//G (\mathcal{O}^\ast)\cong\tilde\mu^{-1}(0)/G$\\

 Method 1 works in two steps: First, we show Theorem \autoref{thm:marsden-weinstein} for the zero orbit $\mathcal{O}^\ast = \{0\}$. In the second step, we apply the \textbf{shifting trick} (cf. \cite{Mayrand2016}), i.e., we identify $\mu^{-1}(\mathcal{O}^\ast)$ with $\tilde\mu^{-1}(0)$ of the modified moment map $\tilde\mu:M\times\mathcal{O}^\ast\to\mathfrak{g}^\ast$, $\tilde\mu(p,\eta) \coloneqq \mu (p)- \eta$.\\
 
 \textit{Step 1:} $\mathcal{O}^\ast = \{0\}$\\
 
 The idea of Step 1 is based on a simple construction from linear algebra:
 
 \begin{proposition}[Coisotropic reduction]\label{prop:coiso_red}
  Let $(V,\omega)$ be a symplectic vector space and $W\subset V$ a coisotropic subspace, i.e.:
   \begin{gather*}
    W^{\perp \omega}\coloneqq\{v\in V\mid \omega (v,w) = 0\ \forall w\in W\}\subset W.
   \end{gather*}
   Then, $(W/W^{\perp\omega}, \bar{\omega})$ is also a symplectic vector space, where $\bar{\omega}$ is defined as follows:
   \begin{gather*}
    \bar{\omega} ([v],[w])\coloneqq \omega (v,w)\quad\forall v,w\in W.
   \end{gather*}
 \end{proposition}
 
 Our goal is to construct the form $\bar{\omega}_{[p]}$ at each point $[p]\in\mu^{-1}(0)/G$ by applying Proposition \autoref{prop:coiso_red} to $W = T_p(\mu^{-1}(0))$ and $W^{\perp\omega} = T_p (G\cdot p)$, where $p\in\mu^{-1}(0)$, $G\cdot p$ is the orbit of $G$ through $p$, and we identify $T_{[p]} (\mu^{-1}(0)/G)$ with\linebreak $W/W^{\perp\omega} = T_p(\mu^{-1}(0))/T_p (G\cdot p)$. For this to be possible, we need to show:
 \begin{gather}
  \left(T_p\mu^{-1}(0)\right)^{\perp\omega} = T_p (G\cdot p)\subset T_p(\mu^{-1}(0)).\label{eq:sym_red_eq}
 \end{gather}
 $T_p (G\cdot p)\subset T_p(\mu^{-1}(0))$ immediately follows from the fact that $G\cdot p\subset \mu^{-1}(0)$ is an immersed submanifold (cf. Lemma \autoref{lem:immersed_orbit}). To prove $\left(T_p\mu^{-1}(0)\right)^{\perp\omega} = T_p (G\cdot p)$, we compute:
 \begin{gather*}
  \omega_p (X_v (p), w) = d\mu_{v,p} (w) = d\mu_p (w)(v) = 0\ \forall v\in\mathfrak{g}\ \forall w\in T_p(\mu^{-1}(0)) = \ker d\mu_p.
 \end{gather*}
 This implies $\left(T_p\mu^{-1}(0)\right)^{\perp\omega}\supset T_p (G\cdot p)$. Equality now follows from dimensional reasoning.\\
 A priori, we do not know whether $\bar{\omega}$ constructed this way is well-defined, since multiple points $p\in\mu^{-1}(0)$ are mapped to the same point $[p]\in\mu^{-1}(0)/G$. However, the $G$-action on $M$ is symplectic meaning that $\omega$ is $G$-invariant which takes care of this problem. To conclude the proof of Step 1, we note that, by \autoref{eq:sym_red_eq}, the $G$-invariant and closed form $i^\ast\omega$ is also horizontal with respect to the $G$-principal\footnote{$\mu^{-1}(0)\stackrel{\pi}{\to}\mu^{-1}(0)/G$ is a $G$-principal bundle, as $G$ acts freely and properly on $\mu^{-1}(0)$ (cf. Remark \autoref{rem:godement}).} bundle $\mu^{-1}(0)\stackrel{\pi}{\to}\mu^{-1}(0)/G$. This allows us to apply the following proposition\footnote{Even though Proposition \autoref{prop:hor_forms} incorporates Proposition \autoref{prop:coiso_red}, we have listed Proposition \autoref{prop:coiso_red} separately for the sake of clarity.} showing that $\bar{\omega}$ has the desired properties:
 
 \begin{proposition}\label{prop:hor_forms}
  Let $P\stackrel{\pi}{\to}B$ be a $G$-principal bundle and let $\omega\in\Omega^k (P)$ be $G$-invariant and horizontal, i.e. $\iota_v\omega = 0\ \forall v\in\ker d\pi$. Then, there exists a unique $k$-form $\bar{\omega}\in\Omega^k (B)$ with $\pi^\ast\bar{\omega} = \omega$. Furthermore, $\bar{\omega}$ is closed if and only if $\omega$ is closed. For $k=2$, the two-form $\bar{\omega}$ is non-degenerate if and only if $\ker\omega = \ker d\pi$.
 \end{proposition}
 
 \begin{proof}[Proof of Proposition \autoref{prop:hor_forms}]
  As $\pi$ is a surjective submersion and $\bar{\omega}$ needs to satisfy $\pi^\ast\bar{\omega} = \omega$, $\bar{\omega}$ is uniquely defined by:
  \begin{gather*}
   \bar{\omega} (v_1,\ldots, v_k)\coloneqq \omega (w_1,\ldots, w_k),
  \end{gather*}
  where $d\pi (w_j) = v_j$. $\bar{\omega}$ is well-defined, since $\omega$ is $G$-invariant and horizontal. We can see that $\bar{\omega}$ is smooth by choosing local trivializations of $P\stackrel{\pi}{\to}B$. To prove the first equivalence, we note that $d\bar{\omega}\in\Omega^{k+1} (B)$ is the unique form satisfying $\pi^\ast d\bar{\omega} = d\omega$. If $\omega$ is closed, then $0\in\Omega^{k+1} (B)$ also fulfills the equation $\pi^\ast 0 = d\omega$, thus, $d\bar{\omega} = 0$ by uniqueness. The converse direction as well as the last equivalence are trivial.
 \end{proof}
 
 \textit{Step 2:} Shifting trick\\
 
 By Proposition \autoref{prop:prep}, the $G$-action on $\tilde M$ is proper and $\mu^{-1}(\mathcal{O}^\ast)$ is isomorphic to $\tilde\mu^{-1}(0)$. Since $G$ acts freely on $\mu^{-1}(\mathcal{O}^\ast)$, $G$ also acts freely on $\mu^{-1}(0)$. This allows us to apply Step 1 to the $\mathcal{O}^\ast$-extended space $\tilde M$ with modified moment map $\tilde\mu$ yielding a symplectic structure on $\tilde\mu^{-1}(0)/G$. We can use Proposition \autoref{prop:prep} again to identify $\tilde\mu^{-1}(0)/G$ with $\mu^{-1}(\mathcal{O}^\ast)/G$ giving us a symplectic form $\bar{\omega}$ on $M//G (\mathcal{O}^\ast)$. One easily checks that $\bar{\omega}$ possesses all properties stipulated in Theorem \autoref{thm:marsden-weinstein}.\\
 
 \textbf{Method 2:} $M//G (\mathcal{O}^\ast)\cong \mu^{-1}(\eta)/G_\eta$\\
 
 Ideally, we would like to repeat Step 1 of Method 1 to equip $\mu^{-1}(\eta)/G_\eta$ with a symplectic structure. Unfortunately, this is not possible, as $T_p (\mu^{-1}(\eta))$ is not a coisotropic subspace of $T_pM$ anymore. To rectify this, we have to improve the linear reduction process:
 
 \pagebreak
 
 \begin{proposition}[Generalized linear reduction]\label{prop:gen_coiso_red}
  Let $(V,\omega)$ be a symplectic vector space and $W\subset V$ a subspace. Denote the $\omega$-orthogonal complement of $W$ by $W^{\perp\omega}$. Then, $(W/(W^{\perp\omega}\cap W), \bar{\omega})$ is also a symplectic vector space, where $\bar{\omega}$ is defined as follows:
  \begin{gather*}
   \bar{\omega} ([v],[w])\coloneqq \omega (v,w)\quad\forall v,w\in W.
  \end{gather*}
 \end{proposition}
 
 To apply Proposition \autoref{prop:gen_coiso_red}, we need to show:
 \begin{align}
  \left(T_p\mu^{-1}(\eta)\right)^{\perp\omega} &= T_p (G\cdot p),\label{eq:sym_red_eq_alt}\\
  T_p (G\cdot p)\cap \left(T_p\mu^{-1}(\eta)\right) &= T_p (G_\eta\cdot p).\label{eq:sym_red_eq_2}
 \end{align}
 \autoref{eq:sym_red_eq_alt} is proven in the same way as \autoref{eq:sym_red_eq}. One easily verifies \autoref{eq:sym_red_eq_2}. Thus, we obtain a non-degenerate two-form on $\mu^{-1}(\eta)/G_\eta$. Under the isomorphism $\mu^{-1}(\eta)/G_\eta\cong\mu^{-1}(\mathcal{O}^\ast)/G$ from Proposition \autoref{prop:prep}, this form coincides with the one from Method 1, as one can check.
\end{proof}

Let us apply symplectic reduction to cotangent bundles (cf. Example \autoref{ex:cotangent}). First, we check that the conditions of Theorem \autoref{thm:marsden-weinstein} are fulfilled. In Example \autoref{ex:cotangent}, we have seen that the $G$-action on $T^\ast G$ is isomorphic to the following action on $T^\ast_L G = G\times \mathfrak{g}^\ast$:
\begin{gather*}
 G\times T^\ast_L G\to T^\ast_L G,\quad (g, (h,\alpha))\mapsto (hg^{-1}, \Ad^\ast (g)\alpha).
\end{gather*}
This action is just the product of the right action on $G$ and the coadjoint action on $\mathfrak{g}^\ast$. As the right action on $G$ is free and proper, the action on $T^\ast _LG$ is as well.\\
Now consider the reduced space $T^\ast G//G (\mathcal{O}^\ast)$. $\mu^L:T^\ast_L G\to\mathfrak{g}^\ast$ is just the projection onto the second component, so we find for any $\eta\in\mathcal{O}^\ast$:
\begin{gather*}
 T^\ast G//G (\mathcal{O}^\ast)\cong (\mu^L)^{-1}(\eta)/G_\eta\cong G/G_\eta\cong\mathcal{O}^\ast.
\end{gather*}
We see that the reduced space $T^\ast G//G(\mathcal{O}^\ast)$ is just $\mathcal{O}^\ast$ itself. The form $\bar{\omega}$ on $T^\ast G//G (\mathcal{O}^\ast)$ is also quite familiar. To compute $\bar{\omega}$, we use the second method. First, note that $T^\ast G$ admits another $G$-action induced by left instead of right multiplication. On $T^\ast_L G$, this action becomes:
\begin{gather*}
 A_L:G\times T^\ast_L G\to T^\ast_L G,\quad (g, (h,\alpha))\mapsto (gh, \alpha).
\end{gather*}
$A_L$ descends to a transitive action on $T^\ast G//G(\mathcal{O}^\ast)$. In fact, it is just the coadjoint action on $T^\ast G//G(\mathcal{O}^\ast)\cong\mathcal{O}^\ast$. In any case, the upshot of this observation is that it suffices to consider the fundamental vector fields associated with $A_L$ to calculate $\bar{\omega}$. They are given by\footnote{Here, we identify the tangent spaces of $T^\ast_L G$ with $\mathfrak{g}\times\mathfrak{g}^\ast$ as in Example \autoref{ex:cotangent}.}:
\begin{gather*}
 X^{A_L}_v (g,\alpha) = (\Ad (g^{-1})v,0)\quad\forall v\in\mathfrak{g}\ \forall (g,\alpha)\in T^\ast_L G.
\end{gather*}
We fix $\eta\in\mathcal{O}^\ast$ and compute:
\begin{align*}
 \bar{\omega}_\alpha (X^\ast_v (\alpha), X^\ast_w (\alpha)) &= \omega^L_{(g,\eta)} (X^{A_L}_v(g,\eta), X^{A_L}_w(g,\eta))\\
 &= \omega^L_{(g,\eta)}((\Ad (g^{-1})v,0), (\Ad (g^{-1})w,0))\\
 &= -\eta ([\Ad (g^{-1})v,\Ad (g^{-1})w])\\
 &= - (\Ad^\ast (g)\eta) ([v,w])\\
 &= -\alpha ([v,w]),
\end{align*}
where $\alpha\in\mathcal{O}^\ast\cong G/G_\eta\cong T^\ast G//G (\mathcal{O}^\ast)$, $g\in G$ is chosen such that $\alpha = \Ad^\ast (g)\eta$, $v,w\in\mathfrak{g}$, and $X^\ast_v, X^\ast_w$ are fundamental vector fields associated with the coadjoint action. Comparing $\bar{\omega}$ with \autoref{eq:kks_form}, we realize that $\bar{\omega}$ is just the negative KKS form. This may not come as a surprise, since Theorem \autoref{thm:marsden-weinstein} involves the negative KKS form. Nevertheless, this result is indeed quite remarkable, since Method 2 from the proof of Theorem \autoref{thm:marsden-weinstein}, which we employed to compute $\bar{\omega}$, does not use the KKS structure on $\mathcal{O}^\ast$.\\
The Marsden-Weinstein reduction of Example \autoref{ex:cotangent} is, in some sense, a blueprint for reductions of cotangent bundles of Lie groups: If the cotangent bundle $T^\ast G$ and a coadjoint orbit $\mathcal{O^\ast}$ of a Lie group $G$ admit a certain structure (symplectic, Kähler, holomorphic Kähler, Hyperkähler,\dots), we expect that the structure on $T^\ast G$ transfers to a similar structure on the reduced space $T^\ast G//G (\mathcal{O}^\ast)$ via the shifting trick. As the reduced space is diffeomorphic to $\mathcal{O}^\ast$, the structure on $T^\ast G//G (\mathcal{O}^\ast)$ should coincide with the one on $\mathcal{O}^\ast$ used for the reduction. If this is the case, we say that the structures on $T^\ast G$ and $\mathcal{O}^\ast$ are \textbf{compatible}.\\
Let us now consider Kähler reductions. As we will see soon, we cannot apply Method 2 from the proof of Theorem \autoref{thm:marsden-weinstein} in the Kähler case. The reason for this is that the Kähler structure on the reduced space depends on the choice of Kähler structure on the coadjoint orbit $\mathcal{O}^\ast$. While the symplectic structure on $\mathcal{O}^\ast$ is canonical and unique, the same does not hold for the Kähler structures on $\mathcal{O}^\ast$, since they depend on the choice of an $\Ad$-invariant scalar product (cf. \autoref{sec:semi-kaehler}). Method 2 does not ``see'' the Kähler structures on $\mathcal{O}^\ast$, so we are forced to use the shifting trick.\\
Recall that the shifting trick consists of two steps: One first performs the reduction at $0\in\mathfrak{g}^\ast$ and then considers the $\mathcal{O}^\ast$-extended space $\tilde M$. Step 1 is more or less the same for all Kähler structures, while Step 2 differs for the various structures (Kähler, holomorphic Kähler,\dots). We, therefore, separate Step 1 and 2 and list them as their own theorems:

\begin{theorem}[Semi-Kähler reduction at $0\in\mathfrak{g}^\ast$]\label{thm:semi-kaehler-red-zero}
 Let $(M,\omega, J)$ be a semi-Kähler manifold with semi-Riemannian metric $g\coloneqq\omega (\cdot,J\cdot)$, $G$ a Lie group, and $G\times M\to M$ a proper Hamiltonian action preserving\footnote{This means $d\phi_g\circ J_p = J_{gp}\circ d\phi_g$ for every $p\in M$ and $g\in G$, where $\phi_g:M\to M$ is defined by $\phi_g (p) \coloneqq gp$.} $J$ with moment map $\mu:M\to\mathfrak{g}^\ast$. Assume that $G$ acts freely on $\mu^{-1}(0)$ and the following condition is fulfilled:
 \begin{gather}
  T_p (\mu^{-1}(0))\cap (T_p\mu^{-1}(0))^{\perp g} = \{0\}\quad\forall p\in\mu^{-1}(0).\label{eq:weird_con}
 \end{gather}
 Then, the reduced space $M//G(0)\coloneqq \mu^{-1}(0)/G$ admits a unique semi-Kähler structure given by $(M//G (0), \bar{\omega},\bar{J})$, where $\bar{\omega}$ is the symplectic form from Theorem \autoref{thm:marsden-weinstein} and $\bar{J}$ is uniquely determined by:
 \begin{gather}
  \bar{J}_{\pi (p)}\circ d\pi_p (v) = d\pi_p\circ J_p(v)\quad\forall v\in H_p\ \forall p\in\mu^{-1}(0)\label{eq:J_bar_def}
 \end{gather}
 with $H$ being the Ehresmann connection, i.e. $G$-invariant horizontal distribution, of the $G$-principal bundle $\mu^{-1}(0)\stackrel{\pi}{\to}M//G(0)$ defined by:
 \begin{gather}
  H_p\coloneqq T_p(\mu^{-1}(0))\cap J_p (T_p\mu^{-1}(0))\quad\forall p\in\mu^{-1}(0).\label{eq:def_of_H}
 \end{gather}
 If $(M,\omega,J)$ is Kähler, i.e., $g$ is positive definite, then \autoref{eq:weird_con} is automatically satisfied and $(M//G (0), \bar{\omega},\bar{J})$ is also Kähler.
\end{theorem}

\begin{proof}
 We know by Theorem \autoref{thm:marsden-weinstein} that $(M//G (0),\bar{\omega})$ is a symplectic manifold. Thus, we only need to construct $\bar{J}$ and prove that it has the desired properties. The following proposition shows everything except integrability:
 
 \begin{proposition}[Linear semi-Kähler reduction]\label{prop:lin_kaehler_red}
  Let $(V,\omega, J)$ be a semi-Kähler vector space, i.e., $\omega\in\Lambda^2 V^\ast$ is non-degenerate and $J\in\End (V)$ satisfies $J^2 = -\id_V$ as well as $\omega (J\cdot,J\cdot) = \omega$. Further, let $W\subset V$ be a coisotropic subspace. Set $g\coloneqq\omega (\cdot,J\cdot)$ and $H\coloneqq W\cap JW$. Assume that $W$ satisfies $W\cap W^{\perp g} = \{0\}$. Then, we have the decomposition $W = W^{\perp\omega}\oplus H$. In particular, the natural projection $\pi:W\to W/W^{\perp\omega}$ restricts to an isomorphism between $H$ and $W/W^{\perp\omega}$. Moreover, $(W/W^{\perp\omega},\bar{\omega},\bar{J})$ is a semi-Kähler vector space, where $\bar{\omega}$ and $\bar{J}$ are defined by $\pi^\ast\bar{\omega} = \omega\vert_{W\times W}$ and $\bar{J}\circ\pi\vert_H = \pi\vert_H\circ J\vert_H$, respectively. If $(V,\omega,J)$ is Kähler, i.e. $g$ is positive definite, then $W\cap W^{\perp g} = \{0\}$ automatically holds and $(W/W^{\perp\omega},\bar{\omega},\bar{J})$ is also Kähler.
 \end{proposition}
 
 \begin{proof}[Proof of Proposition \autoref{prop:lin_kaehler_red}]
  We first show $W^{\perp\omega}\cap H = \{0\}$. Take a vector $\tilde w\in W^{\perp\omega}\cap H$, then there exists a vector $w\in W$ such that $\tilde w = Jw$ and:
  \begin{gather*}
   \omega (v,Jw) = 0\ \forall v\in W\quad\Rightarrow g(v,w) = 0\ \forall v\in W.
  \end{gather*}
  Thus, $w\in W\cap W^{\perp g} = \{0\}$ and $\tilde w = 0$.\\
  To prove $W = W^{\perp\omega}\oplus H$, it now suffices to show:
  \begin{gather*}
   \dim H = \dim W - \dim W^{\perp\omega}.
  \end{gather*}
  Setting $2n\coloneqq\dim V$ and $k\coloneqq\dim W^{\perp\omega}$, we obtain $\dim W = \dim JW = 2n-k$ and compute:
  \begin{align*}
   \dim H &= \dim W + \dim JW - \dim (W+JW) = 2\dim W -\dim V\\
   &= 4n - 2k - 2n = 2(n-k) = (2n -k) - k\\
   &= \dim W - \dim W^{\perp\omega},
  \end{align*}
  where we used $W\cap W^{\perp g} = \{0\}$, $W^{\perp g} = JW^{\perp\omega}$, and $W^{\perp\omega}\subset W$ to show:
  \begin{gather*}
   V = W\oplus W^{\perp g} = W\oplus JW^{\perp\omega} = W + JW.
  \end{gather*}
  Given the decomposition $W = W^{\perp\omega}\oplus H$, it is clear that $\pi\vert_H:H\to W/W^{\perp\omega}$ is an isomorphism. Thus, $\bar{J}$ is well-defined and uniquely determined by\linebreak $\bar{J}\circ\pi\vert_H = \pi\vert_H\circ J\vert_H$. $\bar{\omega}$ is obviously well-defined. It is now trivial to check that $(W/W^{\perp\omega},\bar{\omega},\bar{J})$ is a semi-Kähler vector space.\\
  If $g$ is positive definite, then the $g$-orthogonal complement $W^{\perp g}$ is an actual complement of $W$ in the sense that $V = W\oplus W^{\perp g}$ implying $W\cap W^{\perp g} = \{0\}$. In this case, $\bar{g} \coloneqq \bar{\omega}(\cdot,\bar{J}\cdot)$ is also positive definite, since $(\pi\vert_H)^\ast\bar{g} = g\vert_{H\times H}$. This turns $(W/W^{\perp\omega},\bar{\omega},\bar{J})$ into a Kähler vector space concluding the proof.
 \end{proof}
 
 Let us return to the proof of Theorem \autoref{thm:semi-kaehler-red-zero}. To construct $\bar{J}$, we want to apply Proposition \autoref{prop:lin_kaehler_red} at each point $[p]\in\mu^{-1}(0)/G$ to $W = T_p(\mu^{-1}(0))$,\linebreak $W^{\perp\omega} = T_p (G\cdot p)$, and $H = H_p$, where we identify $T_{[p]} (\mu^{-1}(0)/G)$ with\linebreak $W/W^{\perp\omega} = T_p(\mu^{-1}(0))/T_p (G\cdot p)$. Before we can do this, we have to check the conditions stipulated in Proposition \autoref{prop:lin_kaehler_red}. We have already checked in the proof of Theorem \autoref{thm:marsden-weinstein} that $T_p(\mu^{-1}(0))$ is a coisotropic subspace of $T_pM$, while $W\cap W^{\perp g} = \{0\}$ is just \autoref{eq:weird_con}.\\
 As in the proof of Theorem \autoref{thm:marsden-weinstein}, we do not know yet whether $\bar{J}$ constructed this way is well-defined. This is the case if $H$ is $G$-invariant which we can infer directly: Since $J$ and $\mu^{-1}(0)$ are $G$-invariant, it immediately follows from \autoref{eq:def_of_H} that $H$ is also $G$-invariant. At this point, we can see that $H$ is even an Ehresmann connection. By Proposition \autoref{prop:lin_kaehler_red}, $H_p$ is a horizontal, i.e., a complement of the vertical space $T_p (G\cdot p) = (T_p\mu^{-1}(0))^{\perp\omega}$. This also implies that the dimension of $H_p$ is independent of $p$, turning $H$ into a distribution. The distribution $H$ is smooth, because it is the $g$-orthogonal complement of the vertical bundle of $\mu^{-1}(0)\to\mu^{-1}(0)/G$. Indeed, it is straightforward to check that the space $H$ from Proposition \autoref{prop:lin_kaehler_red} satisfies $H = (W^{\perp\omega})^{\perp g\vert_{W\times W}}$.\\
 The smoothness of $H$ infers that $\bar{J}$ is also smooth which becomes clear by going into trivializations of $\mu^{-1}(0)\stackrel{\pi}{\to}\mu^{-1}(0)/G$. The algebraic properties of $\bar{\omega}$ and $\bar{J}$ follow from Proposition \autoref{prop:lin_kaehler_red}. It remains to be shown that $\bar{J}$ is integrable. $\bar{J}$ being integrable is equivalent to:
 \begin{gather*}
  [T^{(0,1)}(\mu^{-1}(0)/G), T^{(0,1)}(\mu^{-1}(0)/G)] \subset T^{(0,1)}(\mu^{-1}(0)/G),
 \end{gather*}
 i.e., for every vector field $\bar{X}$ and $\bar{Y}$ on $\mu^{-1}(0)/G$ there exists a vector field $\bar{Z}$ such that:
 \begin{gather}
  [\bar{X} + i\bar{J}\bar{X}, \bar{Y} + i\bar{J}\bar{Y}] = \bar{Z} + i\bar{J}\bar{Z}.\label{eq:aux_hor_J}
 \end{gather}
 As $H$ is an Ehresmann connection, we can find unique horizontal lifts $X$ and $Y$ of the vector fields $\bar{X}$ and $\bar{Y}$. $J$ is integrable, hence, there is a vector field $Z$ such that:
 \begin{gather*}
  [X + iJX, Y + iJY] = Z + iJZ
 \end{gather*}
 In fact, $Z$ is just $[X,Y] - [JX,JY]$. $H$ is $J$-invariant and the commutator preserves horizontal vector fields, thus, $Z$ is also a horizontal vector field. This allows us to identify $Z$ via $\pi$ with a vector field $\bar{Z}$ on $\mu^{-1}(0)/G$ satisfying \autoref{eq:aux_hor_J} and concluding the proof.
\end{proof}

\begin{remark}\ 
 \begin{enumerate}
  \item At first glance, \autoref{eq:weird_con} seems to be an unnatural condition and one wonders why this condition is not always satisfied. Even though \autoref{eq:weird_con} holds for all Riemannian metrics $g$, the same is not true for semi-Riemannian metrics $g$. Consider, for instance, the standard Lorentzian metric $g = dx^2_1- dx^2_2$ on $\R^2$ and the subspace $W$ spanned by $v = \begin{pmatrix*}1 & 1\end{pmatrix*}$. In this case, one has $W^{\perp g} = W$. There is different, but equivalent formulation of \autoref{eq:weird_con} which is much more natural: One can show that \autoref{eq:weird_con} is satisfied if and only if $g$ restricts to a semi-Riemannian metric on $\mu^{-1}(0)$.
  \item We can see now why Method 2 is not applicable for Kähler reductions: In the symplectic case, we could generalize Proposition \autoref{prop:coiso_red} to Proposition \autoref{prop:gen_coiso_red} by dropping the condition that $W\subset V$ is a coisotropic subspace. However, Proposition \autoref{prop:lin_kaehler_red} completely fails if $W\subset V$ is not coisotropic.
 \end{enumerate}
\end{remark}

To execute Step 2, we have to find examples where the $\mathcal{O}^\ast$-extended space $\tilde M$ naturally fulfills \autoref{eq:weird_con} in Theorem \autoref{thm:semi-kaehler-red-zero}. Kähler manifolds with compact group actions and holomorphic Kähler manifolds with complex reductive group actions are among those examples:

\begin{theorem}[Kähler reduction]\label{thm:kaehler_red}
 Let $(M,\omega, J)$ be a Kähler manifold, $G$ a compact Lie group, and $G\times M\to M$ a Hamiltonian action preserving $J$ with moment map $\mu:M\to\mathfrak{g}^\ast$. Assume that $\mathcal{O}^\ast\subset\mathfrak{g}^\ast$ is a coadjoint orbit of $G$ such that $G$ acts freely on $\mu^{-1}(\mathcal{O}^\ast)$. Then, the reduced space $M//G(\mathcal{O}^\ast)$ admits a complex structure $\bar{J}_{\skcdot}$ for every $\Ad$-invariant positive definite scalar product $\skcdot$ on $\mathfrak{g}$ such that $(M//G (0), \bar{\omega},\bar{J}_{\skcdot})$ is a Kähler manifold, where $\bar{\omega}$ is the symplectic form from Theorem \autoref{thm:marsden-weinstein}.
\end{theorem}

\begin{remark}
 By Proposition \autoref{prop:averaging}, a compact group $G$ admits a plethora of $\Ad$-invariant scalar products $\skcdot$.
\end{remark}

\begin{proof}
 Equip $\mathcal{O}^\ast$ with the (negative\footnote{For the shifting trick, the symplectic form of this Kähler structure should be $-\kks$, hence, the complex structure picks up an additional sign as well.}) Kähler structure determined by $\skcdot$ (cf. Theorem \autoref{thm:semi-kaehler}). This turns $\tilde M = M\times\mathcal{O}^\ast$ into a Kähler manifold. We\linebreak already know that the $G$-action on $\tilde M$ is Hamiltonian with moment\linebreak map $\tilde\mu (p,\eta) = \mu (p)-\eta$. As $G$ acts freely on $\mu^{-1}(\mathcal{O}^\ast)$, it also acts freely on $\tilde\mu^{-1}(0)\cong\mu^{-1}(\mathcal{O}^\ast)$. The action is proper, because every compact group action is proper. The complex structure of $\tilde M$ is preserved by the action, since the $G$-actions on $M$ and $\mathcal{O}^\ast$ preserve their respective complex structures. \autoref{eq:weird_con} is trivially satisfied, as the Kähler metric of $\tilde M$ is positive definite. Hence, all conditions of Theorem \autoref{thm:semi-kaehler-red-zero} are fulfilled allowing us to apply it to $\tilde M$ and $\tilde\mu$ which concludes the proof.
\end{proof}

As in the symplectic case, we wish to apply Theorem \autoref{thm:kaehler_red} to Example \autoref{ex:cotangent}. To do this, we first have to construct a Kähler structure on $T^\ast G$. The (various) Kähler structures on (double) cotangent bundles are studied thoroughly in \autoref{sec:duality}. Here, a quick summary shall suffice: Take a compact Lie group $G$ and a positive definite $\Ad$-invariant scalar product $\skcdot$ on $\mathfrak{g}$. To define a complex structure on $T^\ast G$, we identify $T^\ast G$ with $TG$ using the bi-invariant Riemannian metric induced by $\skcdot$ (cf. \autoref{sec:lie_groups}). $TG$ is isomorphic to the universal complexification $G_\C$ of $G$ and, as such, inherits the complex structure from $G_\C$. The diffeomorphism between $TG$ and $G_\C$ is given by left trivialization and subsequent polar decomposition\footnote{If one is unsure why polar decomposition gives rise to a diffeomorphism, one can confer Lemma \autoref{lem:kaehler_reductive_groups}.}:
\begin{alignat*}{2}
 TG &\to\quad T_LG\coloneqq G\times\mathfrak{g} &&\to G_\C,\\
 v\in T_gG &\mapsto\, (g, w = dL_{g^{-1},g}v) &&\mapsto g\exp (iw).
\end{alignat*}
Using the identification $TG\cong T^\ast G$, the complex structure on $TG$ becomes a complex structure on $T^\ast G$. One can check\footnote{A similar computation is done, for instance, in \cite{Bremigan2000}. Confer the last subsection of \autoref{sec:duality} for details.} that this complex structure is compatible with the (negative\footnote{In \cite{Bremigan2000}, they use the convention $g\coloneqq \omega (J\cdot,\cdot)$, while we use the convention $g\coloneqq\omega (\cdot,J\cdot)$.}) canonical symplectic form giving us a Kähler structure on $T^\ast G$.\\
We expect that the Kähler structures on $T^\ast G$ and $\mathcal{O}^\ast$ determined by $\skcdot$ are compatible in the following sense: If we apply Theorem \autoref{thm:kaehler_red} to the Kähler structure on $T^\ast G$ determined by $\skcdot$, then the reduction procedure using $\skcdot$ should yield the (negative) Kähler structure on $\mathcal{O}^\ast$ determined by $\skcdot$. Unfortunately, this computation is quite tedious, as the explicit formula for the complex structure on $T^\ast G$ is rather long and cumbersome (cf. \cite{Bremigan2000}). Therefore, we stop the discussion of Kähler reduction here and move on to the holomorphic Kähler case instead.\\
To describe the reduction process for holomorphic Kähler manifolds, we need to impose additional conditions on the $G$-action. To be precise, we require that the action is compatible with the underlying real structures (cf. \autoref{app:real_structures}):

\begin{definition}[Action compatible with real structures]\label{def:action_comp_with_real_str}
 Let $X$ be a holomorphic symplectic manifold with real structure $\rho$ and $G$ a complex Lie group with real structure $\sigma$. A Hamiltonian $G$-action on $X$ is called \textbf{holomorphic} if the map $G\times X\to X$, $(g,p)\mapsto gp$ and the moment map $\mu:X\to\mathfrak{g}^\ast$ are holomorphic. We say a holomorphic Hamiltonian action is \textbf{compatible} with $\rho$ and $\sigma$ if:
 \begin{gather*}
  \rho (gp) = \sigma (g)\rho (p)\ \ \forall g\in G\,\forall p\in X,\quad\text{and}\quad \mu\circ\rho = d\sigma^\ast_e\circ \mu,
 \end{gather*}
 where $d\sigma^\ast_e:\mathfrak{g}^\ast\to\mathfrak{g}^\ast$ is the map dual to $d\sigma_e:\mathfrak{g}\to\mathfrak{g}$.
\end{definition}

We can now formulate the reduction theorem for holomorphic Kähler manifolds:

\begin{theorem}[Holomorphic Kähler reduction]\label{thm:holo_kaehler_red}
 Let $(X,\omega, J, I)$ be a holomorphic Kähler manifold with real form $M$, $G$ a complex reductive group with real form $G_\R$, and $G\times X\to X$ a holomorphic (w.r.t. $I$) Hamiltonian action with moment map $\mu:X\to\mathfrak{g}^\ast$. Assume that the Hamiltonian action is proper, preserves $J$, and is compatible with the real structures on $X$ and $G$. Assume further that $\mathcal{O}^\ast\subset\mathfrak{g}^\ast$ is a coadjoint orbit of $G$ such that $\mu^{-1}(\mathcal{O}^\ast)\cap M\neq \emptyset$ and $G$ acts freely on $\mu^{-1}(\mathcal{O}^\ast)$. Then, $\mu^{-1}(\mathcal{O}^\ast)$ is an immersed complex submanifold of $(X,I)$ with real structure. The reduced space $X//G(\mathcal{O}^\ast) = \mu^{-1}(\mathcal{O}^\ast)/G$ inherits a complex ($\bar{I}$) and real structure from $\mu^{-1}(\mathcal{O}^\ast)$. For every $G_\R$-invariant positive definite scalar product $\skcdot_\R$ on $\mathfrak{g}_\R$, there exists an open neighborhood $U_{\skcdot_\R}\subset X//G(\mathcal{O}^\ast)$ of $M//G_\R(\mathcal{O}^\ast\cap\mathfrak{g}^\ast_\R)$ and a complex structure $\bar{J}_{\skcdot_\R}$ on $U_{\skcdot_\R}$ such that $(U_{\skcdot_\R}, \bar{\omega},\bar{J}_{\skcdot}, \bar{I})$ is a holomorphic Kähler manifold, where $\bar{\omega}$ is the symplectic form from Theorem \autoref{thm:marsden-weinstein}.
\end{theorem}


\begin{proof}
 By Proposition \autoref{prop:prep}, $\mu:X\to\mathfrak{g}^\ast$ is a holomorphic (w.r.t. $I$) submersion on an open neighborhood of $\mu^{-1}(\mathcal{O}^\ast)$. As in the real case, the preimage of an immersed complex submanifold under a holomorphic submersion is again an immersed complex submanifold. Hence, $\mu^{-1}(\mathcal{O}^\ast)$ is an immersed complex submanifold of $(X,I)$. Denote the real structure on $X$ by $\rho$ and the one on $G$ by $\sigma$. We find $\mathcal{O}^\ast\cap \mathfrak{g}^\ast_\R\neq \emptyset$ which follows from $\mu^{-1}(\mathcal{O}^\ast)\cap M\neq\emptyset$ and the fact that $\mu$ intertwines $\rho$ and $d\sigma^\ast_e$. $\mathcal{O}^\ast\cap \mathfrak{g}^\ast_\R\neq \emptyset$ and $d\sigma_e^\ast\circ\Ad^\ast (g) = \Ad^\ast (\sigma(g))\circ d\sigma^\ast_e$ imply $d\sigma^\ast_e (\mathcal{O}^\ast) = \mathcal{O}^\ast$ giving us:
 \begin{gather*}
  \rho(\mu^{-1}(\mathcal{O}^\ast)) = \mu^{-1}(d\sigma^\ast_e (\mathcal{O}^\ast)) = \mu^{-1}(\mathcal{O}^\ast),
 \end{gather*}
 where we used again that $\mu$ intertwines $\rho$ and $d\sigma^\ast_e$. The last equation tells us that $\rho$ restricts to a real structure on $\mu^{-1}(\mathcal{O}^\ast)$. By Godement's theorem (Theorem \autoref{thm:godement}), we know that $X//G (\mathcal{O}^\ast) = \mu^{-1}(\mathcal{O}^\ast)/G$ is a smooth manifold and that the canonical projection $\pi:\mu^{-1}(\mathcal{O}^\ast)\to X//G(\mathcal{O}^\ast)$ is a smooth submersion. One easily checks that $X//G (\mathcal{O}^\ast)$ obtains a well-defined and unique complex structure $\bar{I}$ by requiring $\pi$ to be holomorphic. In similar fashion, the real structure $\rho$ on $\mu^{-1}(\mathcal{O}^\ast)$ descends to the real structure $\bar{\rho}$ on $X//G (\mathcal{O}^\ast)$ defined by\linebreak $\bar{\rho}([p])\coloneqq [\rho (p)]$. $\bar{\rho}$ is well-defined, since the $G$-action is compatible with $\rho$ and $\sigma$, specifically because of $\rho (gp) = \sigma (g)\rho (p)$. Now note that the map
 \begin{gather*}
  M//G_\R(\mathcal{O}^\ast\cap\mathfrak{g}^\ast_\R)\to X//G (\mathcal{O}^\ast),\quad [p]_{G_\R}\mapsto [p]_G
 \end{gather*}
 is an embedding whose image\footnote{One expects that the image of this embedding not only lies in the real form, but is equal to $\Fix\bar{\rho}$. Often, this is indeed the case. In general, however, there might be obstructions (cf. \autoref{app:complex_lie_groups}, particularly the discussion at the end).} lies in the real form $\Fix\bar{\rho}$. This embedding allows us to view $M//G_\R(\mathcal{O}^\ast\cap\mathfrak{g}^\ast_\R)$ as a subset of $X//G (\mathcal{O}^\ast)$ and an open subset of $\Fix\bar{\rho}$. In particular, $\Fix\bar{\rho}$ is not empty.\\
 The next step is to construct the semi-Kähler structure on $X//G (\mathcal{O}^\ast)$. For this, we choose $\skcdot_\R$ as specified in the theorem. By Proposition \autoref{prop:comp_of_metric}, $\skcdot_\R$ gives rise to a unique $G$-invariant non-degenerate scalar product $\skcdot$ on $\mathfrak{g}$. Equip $\mathcal{O}^\ast$ with the (negative\footnote{As in the Kähler case, the symplectic and complex structure each pick up an additional sign.}) semi-Kähler structure determined by $\skcdot$ (cf. Theorem \autoref{thm:semi-kaehler}). This allows us to put a semi-Kähler structure on $\tilde X = X\times\mathcal{O}^\ast$. Note that $\tilde X$ carries, in fact, a holomorphic Kähler structure, since the semi-Kähler structure on $\mathcal{O}^\ast$ comes from a holomorphic Kähler structure (cf. Corollary \autoref{cor:reductive}). In any case, the $G$-action on $\tilde X$ is Hamiltonian with moment map $\tilde\mu (p,\eta) = \mu (p) - \eta$.\linebreak Ideally, we would apply Theorem \autoref{thm:semi-kaehler-red-zero} to $\tilde X$ and $\tilde\mu$. It is trivial to check all assumptions except for \autoref{eq:weird_con}. Unfortunately, \autoref{eq:weird_con} might not be satisfied for all points $p\in\tilde\mu^{-1}(0)$. We can still apply Theorem \autoref{thm:semi-kaehler-red-zero} if we find an open and $G$-invariant subset $V_{\skcdot_\R}\subset\mu^{-1}(\mathcal{O}^\ast)\cong\tilde\mu^{-1}(0)$ such that \autoref{eq:weird_con} is fulfilled for all points in $V_{\skcdot_\R}$. In this case, only the space $U_{\skcdot_\R}\coloneqq \pi (V_{\skcdot_\R})$ carries the semi-Kähler structure. The space $V_{\skcdot_\R}$ we construct shortly contains $\mu^{-1}(\mathcal{O}^\ast)\cap M$ implying that $U_{\skcdot_\R}\subset X//G (\mathcal{O}^\ast)$ is an open neighborhood of $M//G_\R (\mathcal{O}^\ast\cap\mathfrak{g}^\ast_\R)$.\\
 We begin the construction of $V_{\skcdot_\R}$ by observing that all points in the real form $\mu^{-1}(\mathcal{O}^\ast)\cap M$ of $\tilde\mu^{-1}(0)$ satisfy \autoref{eq:weird_con} due to the following proposition:
 
 \begin{proposition}\label{prop:aux_for_holo_kaehler}
  Let $V$ be a real vector space with linear maps\linebreak $I,\rho\in\End (V)$ and non-degenerate scalar product $g:V\times V\to \R$ satisfying:
  \begin{gather*}
   I^2 = -\rho^2 = -\id_V,\quad \rho I = -I\rho,\quad g(I\cdot, I\cdot) = -g ,\quad g(\rho\cdot,\rho\cdot) = g.
  \end{gather*}
  Further, let $W\subset V$ be a subspace fulfilling $I (W) = I$ and $\rho (W) = W$. Then, we have the following equivalence:
  \begin{gather*}
   V_\rho = W_\rho\oplus W^{\perp g_\rho}_\rho\quad\Leftrightarrow\quad V = W\oplus W^{\perp g},
  \end{gather*}
  where $V_\rho\coloneqq \Fix\rho$, $W_\rho\coloneqq W\cap V_\rho$, and $g_\rho\coloneqq g\vert_{V_\rho\times V_\rho}$.
 \end{proposition}
 
 \begin{proof}[Proof of Proposition \autoref{prop:aux_for_holo_kaehler}]
  By assumption, $(V,I)$ is a complex vector space with real structure $\rho$ and real form $\Fix\rho = V_\rho$. Similarly, $(W,I\vert_W)$ is a complex vector space with real structure $\rho\vert_W$ and real form $\Fix\rho\vert_W = W\cap V_\rho = W_\rho$, as $W$ satisfies $I(W) = W$ and $\rho (W) = W$. Now note that $W^{\perp g}$ also satisfies $I(W^{\perp g}) = W^{\perp g}$ due to $g(I\cdot, I\cdot) = -g$, $I^2 = -\id_V$, and $I(W) = W$. Likewise, we have $\rho (W^{\perp g}) = W^{\perp g}$ because of $g(\rho\cdot, \rho\cdot) = g$, $\rho^2 = \id_V$, and $\rho(W) = W$. Thus, $(W^{\perp g},I\vert_{W^{\perp g}})$ is also a complex vector space with real structure $\rho\vert_{W^{\perp g}}$ and real form $\Fix\rho\vert_{W^{\perp g}} = W^{\perp g}\cap V_\rho$. It is easy to show that $W^{\perp g}\cap V_\rho = W^{\perp g_\rho}_\rho$.\\\\
  ``$\Rightarrow$'': $V_\rho$, $W_\rho$, and $W^{\perp g_\rho}_\rho$ are the real forms of $V$, $W$, and $W^{\perp g}$, respectively. Hence, we find:
  \begin{gather*}
   V = V_\rho\oplus I(V_\rho),\quad W = W_\rho\oplus I(W_\rho),\quad W^{\perp g} = W^{\perp g_\rho}_\rho\oplus I(W^{\perp g_\rho}_\rho).
  \end{gather*}
  This yields:
  \begin{align*}
   V &= V_\rho \oplus I(V_\rho) = W_\rho\oplus W^{\perp g_\rho}_\rho\oplus I(W_\rho)\oplus I(W^{\perp g_\rho}_\rho)\\
   &= W_\rho\oplus I(W_\rho)\oplus W^{\perp g_\rho}_\rho\oplus I(W^{\perp g_\rho}_\rho)\\
   &= W\oplus W^{\perp g}.
  \end{align*}
  ``$\Leftarrow$'': Let $v\in V_\rho$, then there are unique vectors $w_1\in W$ and $w_2\in W^{\perp g}$ such that $v = w_1 + w_2$. Applying $\rho$ to this equation yields:
  \begin{gather*}
   \rho (w_1) + \rho (w_2) = \rho (v) = v = w_1 + w_2.
  \end{gather*}
  Since $w_1$ and $w_2$ are unique, the last equation implies $\rho (w_1) = w_1$ and\linebreak $\rho (w_2) = w_2$. This shows $w_1\in W_\rho$, $w_2\in W^{\perp g}\cap V_\rho = W^{\perp g_\rho}_\rho$, and\linebreak $V_\rho = W_\rho\oplus W^{\perp g_\rho}_\rho$ concluding the proof.
 \end{proof}
 
 Let us return to Theorem \autoref{thm:holo_kaehler_red}. For every point $p\in\mu^{-1}(\mathcal{O}^\ast)\cap M\subset \tilde\mu^{-1}(0)$,\linebreak we now apply Proposition \autoref{prop:aux_for_holo_kaehler} to $V = T_p\tilde X$, $W = T_p (\tilde\mu^{-1}(0))$, $I$ being the complex structure at $p$, $\rho$ being the differential of the real structure at $p$, and $g$ being the semi-Kähler metric of $\tilde X$ at $p$. If we pullback the semi-Kähler metric of $\tilde X$ to its real form $M\times (\mathcal{O}^\ast\cap\mathfrak{g}^\ast_\R)$, it becomes a Kähler metric (recall that $\tilde X$ is a holomorphic Kähler manifold). With the notation from above, this means that $g_\rho$ is a positive definite scalar product implying\linebreak $V_\rho = W_\rho\oplus W^{\perp g_\rho}_\rho$. By Proposition \autoref{prop:aux_for_holo_kaehler}, we now find $V = W\oplus W^{\perp g}$ or, equivalently, $W\cap W^{\perp g} = \{0\}$. Hence, we have shown that \autoref{eq:weird_con} is satisfied for every point $p\in\mu^{-1}(\mathcal{O}^\ast)\cap M$.\\
 Next, we observe that \autoref{eq:weird_con} is an open condition, i.e., if \autoref{eq:weird_con} holds for a point $p\in\tilde\mu^{-1}(0)$, then it also holds on an open neighborhood $V_p\subset\tilde\mu^{-1}(0)$ of $p$. Furthermore, the $G$-invariance of the semi-Kähler metric infers that \autoref{eq:weird_con} is also $G$-invariant, i.e., if \autoref{eq:weird_con} holds for a point $p\in\tilde\mu^{-1}(0)$, then it holds for all points in the orbit $G\cdot p$. Combining everything now proves that \autoref{eq:weird_con} is fulfilled for all points in the $G$-invariant and open neighborhood
 \begin{gather*}
  V_{\skcdot_\R}\coloneqq \bigcup_{p\in\mu^{-1}(\mathcal{O}^\ast)\cap M}\bigcup_{g\in G} gV_p\subset\mu^{-1}(\mathcal{O}^\ast)\cong\tilde\mu^{-1}(0)
 \end{gather*}
 of $\mu^{-1}(\mathcal{O}^\ast)\cap M$. Applying Theorem \autoref{thm:semi-kaehler-red-zero} to $V_{\skcdot_\R}$ now gives us a semi-Kähler structure on $U_{\skcdot_\R} = \pi (V_{\skcdot_\R})$. It is straightforward to verify that this semi-Kähler structure forms together with the complex structure $\bar{I}$ and the real structure $\bar{\rho}$ a holomorphic Kähler structure concluding the proof.
\end{proof}

As in the symplectic and Kähler case, we want to apply Theorem \autoref{thm:holo_kaehler_red} to Example \autoref{ex:cotangent}. To do that, we need to construct a holomorphic Kähler structure on $T^\ast G$ for complex reductive groups $G$. A detailed account of this construction can be found in \autoref{sec:duality}. For now, it shall suffice that, given a $G_\R$-invariant positive definite scalar product $\skcdot_\R$ on $\mathfrak{g}_\R$ and the resulting $G$-invariant scalar product $\skcdot$ on $\mathfrak{g}$, $\omega$ and $J$\footnote{To be precise, $J$ also needs to be twisted with a diffeomorphism $\phi:T^\ast G\to T^\ast G$ (cf. Conjecture \autoref{con:commutation_relation}).} on $T^\ast G$ are obtained as in the Kähler case. The complex structure $I$ on $T^\ast G$ is the one induced by the complex structure of $G$. We have already discussed how $\omega$ and $J$ are reduced, so we only need to investigate the reduction of $I$. First, we note that both $\mu^{-1}(\mathcal{O}^\ast)$ and $\mu^{-1}(\eta)$ ($\eta\in\mathcal{O}^\ast$) are complex submanifolds of $(T^\ast G,I)$ and that the submersions $\mu^{-1}(\mathcal{O}^\ast)\to\mu^{-1}(\mathcal{O}^\ast)/G = T^\ast G//G (\mathcal{O}^\ast)$ and $\mu^{-1}(\eta)\to\mu^{-1}(\eta)/G_\eta\cong T^\ast G //G (\mathcal{O}^\ast)$ induce the same complex structure $\bar{I}$. Recall that $\mu^{-1}(\eta)$ becomes $G\times \{\eta\}\cong G$ under the identification $T^\ast G\cong T^\ast_L G = G\times\mathfrak{g}^\ast$. Given this identification, the complex structure of $\mu^{-1}(\eta)$ coincides with the complex structure of $G$. As shown in \autoref{app:complex_lie_groups}, $G/G_\eta\cong \mathcal{O}^\ast$ is a biholomorphism, where the complex structure of $G/G_\eta$ is induced by the submersion $G\to G/G_\eta$. Thus, the complex structure $\bar{I}$ of $T^\ast G//G (\mathcal{O}^\ast)$ coincides with the complex structure $I_e$ of $\mathcal{O}^\ast$ (cf. Theorem \autoref{thm:complex_orbit_real_form}) under the identification $T^\ast G//G (\mathcal{O}^\ast)\cong \mathcal{O}^\ast$. In particular, the holomorphic Kähler structures on $T^\ast G$ and $\mathcal{O}^\ast$ determined by $\skcdot$ are compatible if the Kähler structures of $T^\ast G_\R$ and $\mathcal{O}^\ast\cap\mathfrak{g}^\ast_\R$ determined by $\skcdot_\R$ are compatible (cf. discussion after Theorem \autoref{thm:kaehler_red}).\\
One noteworthy aspect of Example \autoref{ex:cotangent} is that we can choose\linebreak $U_{\skcdot_\R} = T^\ast G//G (\mathcal{O}^\ast)$. To prove this, it suffices to show that \autoref{eq:weird_con} holds for every point $p\in\mu^{-1}(\mathcal{O}^\ast)$. Recall that $T^\ast G\cong T^\ast_L G = G\times\mathfrak{g}^\ast$ admits two $G$-actions:
\begin{gather*}
 A_L (g, (h,\alpha)) = (gh,\alpha),\quad A_R (g, (h,\alpha)) = (hg^{-1}, \Ad^\ast (g)\alpha).
\end{gather*}
$A_R$ is the Hamiltonian action, while the left action $A_L$ conserves the moment map $\mu^L(g,\alpha) = \alpha$. Both actions preserve the holomorphic Kähler structure of $T^\ast G$. Therefore, the $G\times G$-action obtained by combining $A_L$ and $A_R$ also preserves the holomorphic Kähler structure. Note further that the $G\times G$-action is transitive on $(\mu^L)^{-1}(\mathcal{O}^\ast) = G\times\mathcal{O}^\ast$. Thus, if \autoref{eq:weird_con} is fulfilled for one point, it is satisfied for every point $p\in\mu^{-1}(\mathcal{O}^\ast)$. In the proof of Theorem \autoref{thm:holo_kaehler_red}, we have really shown that \autoref{eq:weird_con} is satisfied for at least one point proving $U_{\skcdot_\R} = T^\ast G//G (\mathcal{O}^\ast)$.\\
For the sake of brevity, we will not discuss the Hyperkähler case in detail. The idea is the same as for the other reduction theorems: First, one performs Hyperkähler reduction\footnote{Hyperkähler reduction for the zero orbit was first introduced by Hitchin, Karlhede, Lindström, and Ro\v{c}ek in 1987 \cite{Hitchin1987b}. Confer \cite{Mayrand2016} for a more modern approach.} for the zero orbit $\mathcal{O}^\ast = \{0\}$. Afterwards, one should be able to apply the shifting trick\footnote{The shifting trick for Hyperkähler quotients was already proposed by Kronheimer in his original paper on Hyperkähler structures of coadjoint orbits (cf. Section 4(b) in \cite{Kronheimer1990}). Even though several similar constructions have been developed throughout the years (cf. Section 4 and 5 in \cite{Mayrand2019} for a selection, in particular Proposition 4.9), nobody expanded on Kronheimer's idea to the author's extent of knowledge.}, where the Hyperkähler structures of $\mathcal{O}^\ast$ introduced in \autoref{app:hyp_orb} are used to equip $\tilde M = M\times\mathcal{O}^\ast$ with a Hyperkähler structure. As before, we expect compatibility meaning that the Hyperkähler quotient of $T^\ast G$, where $T^\ast G$ carries the Hyperkähler structure from \autoref{app:hyp_orb}, gives us the Hyperkähler structure of $\mathcal{O}^\ast$.\\
One remarkable circumstance is that Hyperkähler reduction needs three moment maps, one for each symplectic form $\omega_I$, $\omega_J$, $\omega_K$ associated with the Hyperkähler manifold $(M,g,I,J,K)$. At first glance, this may strike the reader as odd, since all our previous reduction theorems only need one moment map. Nevertheless, this is actually in line with the spirit of our reduction theorems: In the symplectic and Kähler case, we only have one symplectic form $\omega$ and, consequently, only one moment map $\mu$. In the holomorphic Kähler case, however, we have two symplectic forms, the real and imaginary part of $\Omega = \omega - i\omega (I\cdot,\cdot)$. Decomposing $\mathfrak{g}^\ast$ into $\mathfrak{g}^\ast_\R\oplus i\mathfrak{g}^\ast_\R$ also gives us two moment maps\footnote{We are a bit sloppy with our wording here: The real and imaginary part of $\mu$ are moment maps with respect to $G_\R$, not $G$. However, the three moment maps used for Hyperkähler reduction are also just moment maps with respect to the compact group $G_\R$ (cf. \cite{Hitchin1987b} and \cite{Mayrand2016}).}, the real and imaginary part of $\mu:X\to\mathfrak{g} = \mathfrak{g}^\ast_\R\oplus i\mathfrak{g}^\ast_\R$. Consequently, Hyperkähler manifolds should have three moment maps, as is the case.

\end{appendix}

\newpage
\phantom{x}
\newpage
\pagenumbering{Roman}
\addcontentsline{toc}{chapter}{Bibliography}
\markboth{}{Bibliography}
\newcommand{\etalchar}[1]{$^{#1}$}

\end{document}